\documentclass[reqno,10pt]{amsart}
\usepackage{amssymb}
\usepackage{amsmath}
\usepackage{amsthm}
\usepackage{amsopn}
\usepackage{mathrsfs}
\usepackage{color}
\usepackage{xcolor}
\usepackage{esint}
\usepackage{stmaryrd}
\usepackage{stmaryrd}
\usepackage[colorlinks,linkcolor=blue,anchorcolor=blue,citecolor=blue]{hyperref}
\usepackage[utf8]{inputenc}
\usepackage{enumitem}
\usepackage{geometry}
\allowdisplaybreaks

\newtheorem{theorem}{Theorem}

\newtheorem{corollary}[theorem]{Corollary}

\newtheorem{lemma}[theorem]{Lemma}
\newtheorem{proposition}[theorem]{Proposition}

\newcommand{\R}{\mathbb{R}}
\renewcommand{\c}{\mathfrak{u} }
\renewcommand{\P}{\mathbf{P}}
\renewcommand{\S}{\mathbf{S}}
\newcommand{\ip}{(\mathbf{I}-\mathbf{P})}
\newcommand{\ipt}{(\mathbf{I}-\tilde{\mathbf{P}})}
\newcommand{\yb}{y_{\mathbf{b}}}
\newcommand{\tb}{t_{\mathbf{b}}}
\newcommand{\vb}{v_{\mathbf{b}}}
\newcommand{\1}{\mathbf{1}}
\newcommand{\dd}{\mathrm{d}}
\let\pt=\partial
\let\d=\delta
\let\a=\alpha
\let\e=\varepsilon

\let\s=\sigma

\let\O=\Omega
\let\G=\Gamma
\let\o=\omega
\let\g=\gamma
\let\th=\theta
\let\l=\lambda

\let\b=\beta

\newcommand\normmm[1]{\left\vert\kern-0.25ex \left\vert\kern-0.25ex \left\vert #1 \right\vert\kern-0.25ex \right\vert\kern-0.25ex \right\vert}
\newcommand\normm[1]{\left\lVert#1\right\rVert}
\newcommand\norm[1]{\left\lvert#1\right\rvert}
\newcommand\inn[1]{\left\langle#1\right\rangle}

\numberwithin{equation}{section}
\numberwithin{theorem}{section}

\allowdisplaybreaks[3]

\begin{document}

\title[Strong diffusive limit of Boltzmann equation with Maxwell boundary]
{Strong diffusive limit of the Boltzmann equation with Maxwell boundary condition}

\author[Y. Guo]{Yan Guo}
\thanks{Y. Guo: Division of Applied Mathematics, Brown University, Providence, RI 02812, U.S.A.; email: yan\_guo@brown.edu}

\author[J. Jung]{Junhwa Jung}
\thanks{J. Jung: Department of Mathematics, The Pennsylvania State University, State college, PA 16801, U.S.A.; email: jbj5730@psu.edu}

\author[F. Zhou]{Fujun Zhou}
\thanks{F. Zhou: School of Mathematics, South China University of Technology, Guangzhou 510640, P.R. China; email: fujunht@scut.edu.cn}

%\footnote{\emph{Mathematics Subject Classification (2010):} }
%\date{\today}

%\thanks{}

%\subjclass{Primary XXXX; Secondary YYYY}

\begin{abstract}

 While weak diffusive limit from the Boltzmann equation to the incompressible Navier-Stokes-Fourier system was established for the Maxwell boundary condition
 within renormalized solutions framework \cite{SaintRaymond2009, jiang2017boundary}, the corresponding strong diffusive limit has remained outstanding except when the accommodation coefficient $\a \sim \e^{1/2}$ \cite{jiang2017boundary}. We establish global in time strong diffusive limit for all accommodation coefficients $\a \in [0, 1]$ within strong solutions framework. The main novelties of our proof include: (1) a $\e$-stretching method for reduction to a single-bounce $L^\infty$ estimate; (2) a dissipation estimate for a carefully constructed rotating Maxwellian in the near-specular regime $\a\ll \e$.
\end{abstract}

\maketitle

\tableofcontents
%\makeindex
\thispagestyle{empty}

%%%%%%%%%%%%%%%%%%%%%%%%%%%%%%%%%%
%%%%%%%%%%%%%%%%%%%%%%%%%%%%%%%%%%
%%%%%%%%%%%%%%%%%%%%%%%%%%%%%%%%%%
%%%%%%%%%%%%%%%%%%%%%%%%%%%%%%%%%%

%\input{intoduction}
\section{Introduction}
\subsection{Problem Formulation}\
\medskip

This paper is devoted to the study of the strong diffusive limit, within the framework of strong solutions,  of the Boltzmann equation to the incompressible Navier-Stokes-Fourier (INSF) system under the renowned Maxwell boundary condition.

In the diffusive scaling, the evolution of a rarefied gas is governed by the following rescaled Boltzmann equation
\begin{equation} \label{F - Boltzmann equation}
\begin{split}
\e \pt_t F + v \cdot\nabla_x F  =  \e^{-1}Q(F,F) \quad
   &\text{in }  \mathbb{R}^{+}\times \O \times \mathbb{R}^3,   \\
%------
F |_{\g_-}  =   (1-\a) \mathscr{R}F + \a \mathscr{P}F  \;\;\;\;
   &\text {on }  \mathbb{R}^{+}\times \pt\O \times \mathbb{R}^3,  \\
%------
F(t, x, v)|_{t=0}  =   F_0(x,v)  \;\;\;\;
&\text{on }   \O \times \mathbb{R}^3.
\end{split}
\end{equation}
Here, $F(t,x,v)$ represents the distribution density of particles at time $t\geq 0$, position $x\in \O$ and velocity $v\in \mathbb{R}^3$. The Boltzmann collision operator for hard-sphere interactions is given by
\begin{equation*}
\begin{split}
Q(F,H)(v):=&\int_{\R^{3}}\int_{\S^{2}}|(v-v_*)\cdot \s |[F(v')H(u') - F(v)H(u)] \dd \s \dd u\\
:=&Q_{+}(F,H)(v) - Q_{-}(F,H)(v),
\end{split}
\end{equation*}
 where $v' = v - [(v-u) \cdot \s]\s$ and $u' = u + [(v-u) \cdot \s]\s$.
Throughout this work, $\O=\{x:\xi(x)<0\}$ denotes a general bounded  (possibly  non-convex) domain in $\R^3$, with $C^3$ boundary $\partial\O=\{x: \xi(x)=0\}$. We assume $\nabla\xi(x)\neq 0$ on $\partial\O$. The outward unit normal at the boundary is
\begin{equation}\label{n-def}
\begin{split}
n=n(x)=\frac{\nabla\xi(x)}{|\nabla\xi(x)|},
\end{split}
\end{equation}
which admits a smooth extension to a neighborhood of $\partial\Omega$.
The boundary phase space $\g:= \pt \Omega \times \R^3$ decomposes into the outgoing, incoming, and grazing sets:
\begin{equation*}
\begin{split}
\g_{+} :=& \{(x,v) \in \pt \Omega \times \R^3 : n(x) \cdot v > 0\},\\
\g_{-} :=& \{(x,v) \in \pt \Omega \times \R^3 : n(x) \cdot v < 0\},\\
\g_{0} :=& \{(x,v) \in \pt \Omega \times \R^3 : n(x) \cdot v = 0\}.
\end{split}
\end{equation*}

The physical boundary condition in \eqref{F - Boltzmann equation}, known as the {\it Maxwell boundary condition}, was introduced by Maxwell \cite{Maxwell} in 1879 to model gas-surface interactions. The dimensionless accommodation coefficient $\a\in [0, 1]$ characterizes boundary roughness: $\a=0$ represents specular reflection for perfectly smooth surface,
\begin{equation} \label{F - Reflective bdd condition}
\mathscr{R}F(x,v) := F(x,R_{x}v) = F(x,  v - 2[n \cdot v]n);
\end{equation}
 while $\a=1$ denotes diffuse reflection for rough surface,
\begin{equation} \label{F - Diffusive bdd condition}
\mathscr{P}F(x,v) := \sqrt{2 \pi} \mu \int_{n \cdot u > 0} F(x,u) [n \cdot u] \dd u.
\end{equation}
Here $R_x v = v - 2[n \cdot v]n$ is the velocity reflection operator,
\begin{equation}
M_{\rho,u,T}:= \frac{\rho}{(2 \pi T)^{3/2}} \exp \Big(-\frac{\norm{v-u}^2}{2T} \Big)
\end{equation}
denotes the local Maxwellian with density $\rho$, bulk velocity $u$ and temperature $T$, and
\begin{equation}\label{def-mu}
\mu=\mu(v) := M_{1,0,1}= \frac{1}{(2 \pi )^{3/2}}  \exp \Big (-\frac{\norm{v}^2}{2} \Big )
\end{equation}
is the global Maxwellian. The Maxwell boundary condition in \eqref{F - Boltzmann equation} ensures zero net mass flux across boundary:
\begin{equation}\label{zero-flux-bdy}
\int_{\R^3} F(x,v) [n \cdot v] \dd v = 0,  \quad \forall x \in \pt \Omega.
\end{equation}

Let $\mathcal{R}(\Omega)$ denote the finite-dimensional space of rigid motions on $\Omega$ (see \cite{Desvillettes2002}):
\begin{equation*}%\label{infinitesimal_rigid}
\mathcal{R}(\Omega) := \left\{ x \mapsto Ax + x_0 : \; A \in \mathfrak{so}(3,\mathbb{R}),\ x_0 \in \mathbb{R}^3 \right\},
\end{equation*}
where
$$
\mathfrak{so}(3,\mathbb{R}): = \left\{ A = (a_{ij}):\; a_{ij} \in \R, \ i,j=1,2,3, \ A + A^T=0\right\}
$$
is the Lie algebra of $3\times3$ real antisymmetric matrices, equipped with the basis
\begin{equation}\label{A1A2A3 - def}
A_1 = \begin{pmatrix}
    0&  0& 0\\
    0&  0&-1\\
    0&  1& 0
\end{pmatrix},\quad
A_2 = \begin{pmatrix}
    0&  0&-1\\
    0&  0& 0\\
    1&  0& 0
\end{pmatrix},\quad
A_3 = \begin{pmatrix}
    0& -1& 0\\
    1&  0& 0\\
    0&  0& 0
\end{pmatrix}.
\end{equation}
The infinitesimal rigid displacement fields preserving $\Omega$ are defined as
\begin{equation}\label{RO}
\mathcal{R}_\Omega := \left\{ R(x) \in \mathcal{R}(\Omega) :\; x_0 = 0,\ R(x) \cdot n(x) = 0\;\;\;\forall x \in \partial\Omega \right\}.
\end{equation}
 For a bounded domain $\Omega \subset \mathbb{R}^3$ with nonempty boundary $\partial\Omega$, $\dim\mathcal{R}_\Omega \in \{0,1,3\}$. More precisely,
\begin{equation}\label{R-space-form}
\mathcal{R}_\Omega=
\begin{cases}
\{0 \} & \text{ if }\dim\mathcal{R}_\Omega = 0, \\
\text{span}\{A (\cdot) \} & \text{ if }\dim\mathcal{R}_\Omega = 1, \\
\text{span}\{A_1  (\cdot), A_2 (\cdot), A_3 (\cdot)\} & \text{ if }\dim\mathcal{R}_\Omega = 3.
\end{cases}
\end{equation}
 When $\text{dim} \mathcal{R}_\Omega =1$ we take $A =A_{3}$ without loss of generality.
 This dimensional classification corresponds to the following geometric types of the domain:
\begin{equation}\label{def:axiss}
\Omega \text{ is called }
\begin{cases}
\text{non-axisymmetric }  & \text{ if }\dim\mathcal{R}_\Omega = 0, \\
\text{axisymmetric }  & \text{ if }\dim\mathcal{R}_\Omega = 1, \\
\text{spherical } & \text{ if }\dim\mathcal{R}_\Omega = 3.
\end{cases}
\end{equation}
  For conciseness, we shall denote a generic basis element of $\mathcal{R}_\Omega$ by $Ax$ or $R(x)$,  for all three geometric types of $\O$.

Without loss of generality, we assume that the initial data $F_0$ satisfies the following conservation laws:
\begin{equation}\label{initial-conservation-F}
\begin{split}
&\iint_{ \O \times \mathbb{R}^3}  F_0  \dd v\dd x = \iint_{ \O \times \mathbb{R}^3}  \mu  \dd v\dd x =\norm{\Omega}, \\
&\iint_{ \O \times \mathbb{R}^3} Ax \cdot v F_0  \dd v \dd x= \iint_{ \O \times \mathbb{R}^3} Ax \cdot v \mu  \dd v\dd x = 0 \quad \text { for all } Ax \in \mathcal{R}_\O, \\
&\iint_{ \O \times \mathbb{R}^3} \norm{v}^2 F_0  \dd v \dd x= \iint_{ \O \times \mathbb{R}^3} \norm{v}^2 \mu  \dd v\dd x = 3 \norm{\Omega}.
\end{split}
\end{equation}

In the hydrodynamic limit $\varepsilon \to 0$, the relative scaling $\a / {\e}$ plays a critical role in the treatment of boundary conditions. We adopt the following
conventions:
\begin{equation}\label{alpha-classification}
\begin{split}
  \e \lesssim \alpha\leq 1:\;\;  & \lim_{\e \to 0}\frac{\a}{\e} \in (0, \infty] \;\; (\alpha \text{ is of lower or the same order as } \e);\\
 0\leq \alpha \ll \e:\;\;  & \lim_{\e \to 0} \frac{\a}{\e} =0 \;\; (\alpha \text{ is of higher order than } \e, \text{ or } \a=0).
\end{split}
\end{equation}
Thus, for  $\e\in (0,1)$,  the full parameter range $[0,1]$ for $\a$ is partitioned as
\begin{equation}\label{alpha-decomposition}
[0,1]=\{\a: \e \lesssim \a \leq 1\} \cup  \{\a: 0\leq \a\ll \e\}.
\end{equation}
\bigskip

%%%%%%%%%%%%%%%%%%%%%%%%%%%%%%%%%%%%%%%%%%%%%%%%%%%%%%%%%%%%
%%%%%%%%%%%%%%%%%%%%%%%%%%%%%%%%%%%%%%%%%%%%%%%%%%%%%%%%%%%%

\subsection{Strong Limit Result for the Case $\e\lesssim \a\leq 1$}\
\medskip

In the regime $\e\lesssim \a \leq 1$, we define the key limiting parameter
\begin{equation}\label{a-e-limit-infty}
\displaystyle  \lambda:=\frac{1}{\sqrt{2\pi}} \lim_{\e\to  0} \frac{\a}{\e} \in(0,\infty].
\end{equation}
We consider fluctuations around the global Maxwellian $\mu$ via the rescaling
\begin{equation}\label{f-def-F}
F= \mu + \e \sqrt{\mu}f, \quad F_0= \mu+ \e \sqrt{\mu} f_0,
\end{equation}
where $f$ and $f_0$ denote the fluctuation fields.  Under this scaling, the Boltzmann equation \eqref{F - Boltzmann equation} transforms into
\begin{equation}\label{f-eq}
\begin{split}
 \e\pt_t f+ v\cdot\nabla_x  f   +  \e^{-1} L f =  \Gamma(f,f) \;\;\;\;   &\text{ in } \mathbb{R}^{+}\times\O\times\mathbb{R}^3,  \\
 f |_{\g_-} = (1-\a) \mathscr{R}f + \a \mathscr{P}_{\g}f \;\;\;\;   &\text{ on } \mathbb{R}^{+}\times\pt\O\times\mathbb{R}^3, \\
 f(t, x, v)|_{t=0}  =   f_0(x,v)  \;\;\;\;     &\text{ on }   \O \times \mathbb{R}^3,
\end{split}
\end{equation}
with the operators $\G$,  $L$ and $\mathscr{P}_\g$ defined by
\begin{equation}\label{L-Definition}
\begin{split}
{\G}(f,g) :=&  \frac{1}{\sqrt{{\mu}}}Q(\sqrt{{\mu}}f,\sqrt{{\mu}}g),\quad
 %------
{L}(f) :=  -{\G}(\sqrt{{\mu}},f) -{\G}(f,\sqrt{{\mu}}), \\
%-------
\mathscr{P}_\g f:=& \sqrt{2\pi\mu} \int_{n\cdot u>0}f(u) \sqrt{\mu(u)} [n\cdot u] \dd u.
\end{split}
\end{equation}

The null space of $L$ is the five-dimensional subspace of $L^2(\mathbb{R}^3)$ given by
\begin{equation}\label{ker-P}
\ker L = \operatorname{span} \left \{1,   v, |v|^{2} \right \}\sqrt{\mu}.
\end{equation}
An orthonormal basis for $\ker L$ is $\{\chi_{i}\}_{i=0}^{4}$, where
\begin{align}\label{base-hat-chi}
&\chi_{0} := \sqrt{{\mu}}, \quad
\chi_{i} := v_i\sqrt{{\mu}} \;\; (i=1,2,3),\quad
\chi_{4} := \frac{\norm{v}^2-3}{\sqrt{6}}\sqrt{{\mu}}.
\end{align}
The orthogonal projection of $f$ onto $\ker L$ is denoted by
\begin{equation}\label{Pf-abc-hat}
  \P f  =  a \chi_{0} + \sum_{i=1}^{3}b_{i} \chi_{i} + c \chi_{4},
\end{equation}
with coefficients
\begin{equation}\label{hat-abc-def}
a:=\inn{\chi_{0}, f },\quad
b:=\inn{\chi_{i}, f} \;\; (i=1,2,3),\quad
c:=\inn{\chi_{4}, f}.
\end{equation}
Let $(\mathbf{I}-\mathbf{P})f$ denote projection onto the orthogonal
complement of $\ker L $.
\medskip

  We introduce the instant energy functional
\begin{equation}
\begin{split}
\mathscr{E}_{1}[f](t) := & \sup_{0\le s\le t} \Big\{ \normm{f(s)}_{L^2_{x,v}}^{2} +\normm{\pt_t f(s)}_{L^2_{x,v}}^{2} \Big\}
 \end{split}
 \end{equation}
and the dissipation functional
\begin{equation} \begin{split}
\mathscr{D}_{1}[f](t) := & \int_{0}^{t} \Big\{  \normm{\P  f(s)}_{L^2_{x,v}}^{2} + \normm{\P \pt_t f(s)}_{L^2_{x,v}}^{2} \Big\} \dd s \\
%--------------------------
&+ \int_{0}^{t}  \Big\{ \frac{1}{\e^2}  \normm{\ip  f(s)}_{L^2_{x,v}(\nu)}^{2} + \frac{1}{\e^2} \normm{\ip \pt_t f(s)}_{L^2_{x,v}(\nu)}^{2} \Big\} \dd s \\
%---------------
& +  \int_{0}^{t}  \Big\{\frac{\a}{\e} \norm{(1-\mathscr{P}_{\g})f(s)}_{L^2_{\g_{+}}}^2
+   \norm{\mathscr{P}_{\g}f(s)}_{L^2_{\g_{+}}}^2 \Big\}\dd s
\\
%---------------
& + \int_{0}^{t} \Big\{\frac{\a}{\e} \norm{(1-\mathscr{P}_{\g})\pt_t f(s)}_{L^2_{\g_{+}}}^2
 +\norm{\mathscr{P}_{\g}\pt_t f(s)}_{L^2_{\g_{+}}}^2 \Big\}\dd s.
 \end{split} \end{equation}
The total energy functional is defined as
\begin{equation}%\label{total-energy}
\begin{split}
 \normmm{f}_{1}(t) :=\,&\mathscr{E}_{1}^{\frac{1}{2}}[f](t)  + \mathscr{D}_{1}^{\frac{1}{2}}[f](t) + \e^{\frac{1}{2}}\sup_{0\leq s\leq t}\normm{\o f(s)}_{L^{\infty}_{x,v}}\\
&+ \e^{\frac{3}{2}}\sup_{0\leq s\leq t} \normm{\o \pt_t f(s)}_{L^{\infty}_{x,v}}
 + \sup_{0\leq s\leq t}\normm{\P f(s)}_{L^6_{x,v}},
\end{split}
\end{equation}
where the weight function is
\begin{equation}\label{weight-w}
\begin{split}
 \o=\o(v) := e^{\b \norm{v}^2} \quad  \text{ with } 0 < \b \ll \frac{1}{4}.
\end{split}
\end{equation}
The corresponding  norm for the initial data is
\begin{equation}\label{initial-data-total-norm}
\begin{split}
\left[\!\left[{f}_{0}\right]\!\right]_{1} :=&
\normmm{f}_{1}(0)  + \e^{-1} \normm{\ip {f}_{0}}_{L^2_{x,v}(\nu)}
 +  \Big(\frac{\a}{\e}\Big)^{\frac{1}{2}}    \norm{(1-\mathscr{P}_{\g}) f_0 }_{L^2_{\g_{+}}}
 \\ &+
\normm{v\cdot \nabla_x {f}_0}_{L^2_{x,v}}+\normm{v\cdot \nabla_x \pt_t{f}_0}_{L^2_{x,v}},
\end{split}\end{equation}
where $\pt_t {f}_0$ is determined from the perturbation equation \eqref{f-eq}.

We now state the first main result for the regime $\e\lesssim \a\leq 1$.
\medskip

\begin{theorem}[Case $\e\lesssim \a\leq 1$]\label{main-th-1} \
Let $F_0 = \mu + \e \sqrt{\mu} f_0 \geq  0$. Then there exists $\e_0>0$ such that for all $0<\e<\e_0$, if the initial fluctuation satisfies
\begin{align}\label{0-initial-data-total-norm}
  \left[\!\left[{f}_{0}\right]\!\right]_{1}
% \equiv \left[\left[\frac{F_0-\mu}{\e \sqrt{\mu}}\right]\right]_{1}
 \leq  \d_0
\end{align}
for some small constant $\delta_0>0$ independent of $\e$, then the Boltzmann equation with Maxwell boundary condition \eqref{F - Boltzmann equation} admits a unique global strong solution $F = \mu + \e \sqrt{\mu} f \geq 0$ satisfying the uniform bound
\begin{align}\label{0-uniform-bound}
\normmm{f}_{1}(\infty) \leq C \left[\!\left[{f}_{0}\right]\!\right]_{1}
\end{align}
for some constant $C>0$ independent of $\e$.

Moreover, suppose there exist fluid initial data  $(\varrho_0, u_0, \vartheta_0)\in \mathbb{H}_{\vartheta}\times\mathbb{H}_{u}\times\mathbb{H}_{\vartheta}$ (see \eqref{initial-space-Navier}) such that
\begin{align}
\begin{split}\label{initial-converge-condition}
&f _0\to {f}^{*}_0=\Big(\varrho_0+u_0\cdot v+\vartheta_0\frac{|v|^2-3}{2}\Big) \sqrt{\mu}
\ \  \text { strongly in  }  L^2(\O \times \mathbb{R}^3) \text{ as } \e \to 0.
\end{split}
\end{align}
Then the following convergence results hold as $\e \to 0$:
\begin{align}
\frac{F -\mu}{\e} \to  \sqrt{\mu}{f}^{*} = \Big({\varrho}+ {u}\cdot v+{\vartheta}\frac{|v|^2-3}{2} \Big ) {\sqrt{\mu}}
\;\;
&\begin{array}{l}
\text{strongly in }L^2\big(\mathbb{R}^{+};L^2(\O\times \mathbb{R}^3)\big), \\
%-------------
\text{weakly}\!-\!* ~\text{in}~ L^\infty\big( \mathbb{R}^+; L^2(\O\times \mathbb{R}^3)\big),
\end{array}\label{tilde-f-strong-convergence} \\
%-------------
%-------------
 \int_{\mathbb{R}^3} \frac{F -\mu}{\e} \Big[1, v,  \frac{|v|^2-3}{2}\Big] \dd v  \to
  \left(\varrho, u, \vartheta\right )
\;\;  &\hbox{ strongly in }L^2\big(\mathbb{R}^{+};L^2(\O)\big), \label{momentums-strong-convergence}
\end{align}
 where $(\varrho, u, \vartheta)\in C(\mathbb{R}_+, L^2(\O))\cap L^2(\mathbb{R}_+, H^1(\O))$ is the unique weak solution of the INSF system
\begin{equation}\label{INSF-unst}
\begin{split}
\pt_t {u} + {u} \cdot\nabla_x {u}  +\nabla_x {p}
= \s \Delta {u}, \ \ \  \nabla_x\cdot {u} = 0 \ \ \ \    & \text{ in } \mathbb{R}^{+}\times \O, \\
\pt_t {\vartheta} +  u\cdot\nabla_x {\vartheta}
= \kappa  \Delta {\vartheta}, \ \ \  \nabla_x ({\varrho} + {\vartheta} ) = 0 \ \ \ \
  &\text{ in } \mathbb{R}^{+}\times \O, \\
{u}|_{t=0}={u}_0, \ \ \   {\vartheta}|_{t=0}={\vartheta}_0  \ \ \ \     &\text{ on }  \O,
\end{split}
\end{equation}
with viscosity $\s$ and heat conductivity $\kappa$ defined in \eqref{sigma-def} and \eqref{kappa-def}, respectively.

Furthermore, if $\displaystyle \lim_{\e\rightarrow 0}\frac{\alpha}{\e}=\infty$, then \eqref{INSF-unst} is  supplemented
by the Dirichlet boundary condition
\begin{equation}\label{Diri-bdy-unst}
\begin{split}
u=0,  \ \ \   \vartheta = 0 \ \ \   \text{on} \ \   \mathbb{R}^{+}\times\pt\O;
\end{split}
\end{equation}
and if $\displaystyle  \l=\frac{1}{\sqrt{2\pi}} \lim_{\e\rightarrow 0}\frac{\a}{\e} \in (0,\infty)$, then
\eqref{INSF-unst} is  supplemented  by the Navier slip boundary condition
\begin{equation}\label{Navier-bdy-unst}
\begin{split}
 \left [ \sigma\big(\nabla_x u + {(\nabla_x u)}^{\mathrm{T}}\big) \cdot n
  +\l u \right ]^{\mathrm{tan}} = 0, \ \ \  u\cdot n =0 \ \ \ \
    &\text{ on }\mathbb{R}^{+}\times\pt\O, \\
 \kappa \partial_n \vartheta + \frac{4}{5}\l \vartheta = 0 \ \ \ \  &\text{ on }\mathbb{R}^{+}\times\pt\O.
\end{split}
\end{equation}
\end{theorem}
\medskip

Proof of Theorem \ref{main-th-1} will be presented in Section \ref{proof-main-1}.
We remark that the initial requirement  \eqref{0-initial-data-total-norm}, which arises primarily from the $L^2$ and $L^6$ estimates, is natural: only the microscopic part
 $(\mathbf{I}-\mathbf{P})f_0$ and the boundary dissipation $(1-\mathscr{P}_\gamma)f_0$ depend explicitly on $\e$. Hence a wide class of admissible fluctuations $f_0$ satisfies \eqref{0-initial-data-total-norm}; for example,  any $f_0$ whose macroscopic projection $\mathbf{P}f_0$ coincides with the fluid initial data $(\varrho_0, u_0, \vartheta_0)$ of the INSF system \eqref{INSF-unst}--\eqref{Navier-bdy-unst} fulfills this condition.

\bigskip

\subsection{Methodology 1: Streaching Method for $L^\infty$ Estimate}\
\medskip

The inherent low regularity of Boltzmann solutions under physical boundary conditions \cite{guo2017regularity} precludes the use of high-order energy methods. Consequently, we adopt the $L^2$-$L^\infty$ framework pioneered by \cite{guo2010decay}.
A standard $L^2$ energy estimate for \eqref{f-eq} yields
 \begin{equation}\label{model-L2}
  \begin{split}
& \| f ( t ) \|_{L^2_{x,v}}^2 + \frac{1}{\e^{2}} \int_0^t \normm{ ( \mathbf{I} - \mathbf{P} ) f }_{L^2_{x,v}(\nu)}^2
   + \frac{\a}{\e} \int_0^t \norm{ ( 1 - \mathscr{P}_\g )  f }_{L^2_{\g_{+}}}^2
 \lesssim  \frac{1}{\e} \int_0^t \normm{\Gamma( f,  f)}_{L^2_{x,v}}^2   +  \cdots,
\end{split}
\end{equation}
which follows from the Maxwell boundary condition in \eqref{F - Boltzmann equation}. To close the energy estimate, it is necessary to control both $\int_0^t \normm{  \mathbf{P}  f }_{L^2_{x,v}}^2$ and $\normm{ \P   f }_{L^{6}_{x,v}}$ (these bounds are established in Section \ref{subsection-Nonlinear-Estimate}):
\begin{proposition} \label{0-macro-L2-L6-estimate} Let $\e \lesssim \a \leq 1$, and let $f$ be a solution of \eqref{f-eq} satisfying mass conservation law
\begin{equation} \label{0-f-a-conservation-law}
\begin{split}
\iint_{\Omega\times\R^3} f(t,x,v) \dd v \dd x =& 0 \;\;  \text{ for all }  t\in [0, T]
\end{split}
\end{equation}
with $0<T\leq \infty$. Then,  for all $0\leq s\leq t \leq T$, the following estimates hold:
\begin{align}
\int_{s}^{t} \normm{ \P f}_{L^{2}_{x,v}}^2\dd \tau \lesssim \; &
 \e \big[ G_{0} (t)-  G_{0} (s)\big]
 + \int_{s}^{t} \norm{  (1-\mathscr{P}_{\g}) f }_{L^2_{\g_{+}}}^2 \dd \tau \nonumber \\
 %-------
 &+  \int_{s}^{t}   \Big [  \normm{\e^{-1}  \ip f}_{L^{2}_{x,v}(\nu)}^{2} +   \normm{
  {\nu}^{-\frac{1}{2}} \Gamma( f, f) }_{L^{2}_{x,v}}^{2} \Big]  \dd \tau, \label{0-P-f-macro-L2}\\
%----------------
%----------------
\normm{ \P   f }_{L^{6}_{x,v}} \lesssim \; & \e\normm{\pt_t  f }_{L^{2}_{x,v}} + \normm{\P f }_{L^{2}_{x,v}}+ \a \norm{  (1-\mathscr{P}_{\g}) f }_{L^2_{\g_{+}}}^{\frac{1}{2}}
\normm{\omega  f }_{L^{\infty}_{x,v}}^{\frac{1}{2}} \nonumber\\
%-----------
&   + \normm{\e^{-1}\ip  f}_{L^{2}_{x,v}(\nu)}  + \normm{\ip  f}_{L^{6}_{x,v}} + \normm{ {\nu}^{-\frac{1}{2}} \Gamma(f,f)}_{L^{2}_{x,v}},\label{0-Pf-L6}
\end{align}
 where $|G_{0}(t)| \lesssim \normm{f(t)}_{L^{2}_{x,v}}^2$.
\end{proposition}
\medskip

To elucidate the core methodology for  obtaining $L^\infty$ estimates with Maxwell boundary condition in general domains, we first consider a simplified model problem with a specular reflection boundary condition:
\begin{equation} \label{model-eq-st}
\begin{split}
\e\pt_t f+ v\cdot\nabla_x  f   +  \e^{-1} \nu_0 f = \e^{-1} \int_{|v^{\prime}|\leq N}
   f(v^{\prime}) \dd v^{\prime}    \ \ \ \ &\text{ in } \mathbb{R}^{+}\times\O\times\mathbb{R}^3,  \\
f |_{\g_-} = \mathscr{R}f  \ \ \ \  &\text{ on } \mathbb{R}^{+}\times\pt\O\times\mathbb{R}^3,\\
  f |_{t=0} = f_0    \ \ \ \  &\text{ on } \O\times\mathbb{R}^3
\end{split}
\end{equation}
where $\nu_0$ denotes a uniform lower bound of the collision frequency $\nu(v)$, and the integral term on the right-hand side arises from a truncation of $Kf$ (see \eqref{K-def-absolute}). Define the back-time cycles
 \begin{equation*}
 \begin{split}
 \displaystyle
 &  X_{\mathbf{cl}}(s;t,x,v):=\sum_{k} \mathbf{1}_{[{t}_{k+1},{t}_{k})} (s)
              X (s; {t}_{k}, x_{k} , v_{k} ), \\
             &  V_{\mathbf{cl}}(s;t,x,v):= \sum_{k} \mathbf{1}_{[{t}_{k+1},{t}_{k})} (s)
             V (s; {t}_{k}, x_{k} , v_{k} ),
\end{split}
\end{equation*}
where $[X(s;t,x,v),\,V(s;t,x,v)]$ denotes the characteristic trajectories, and $(t_k, x_k, v_k)$ marks the $k$-th bounce of the backward trajectory against $\pt\O$.
The solution of \eqref{model-eq-st} admits the Duhamel representation
 \begin{equation}\label{f-Duhamel-1}
\begin{split}
 f(t, x, v) = & \;\frac{1}{\e} \int_{0}^t  e^{- \frac{\nu_0 }{\e}(t-s) }
\int_{|v^{\prime}|\leq N} f\left(s,X_{\mathbf{cl}}(s;t,x,v),v^{\prime}\right) \dd v^{\prime}  \dd s + \cdots,
\end{split}
\end{equation}
which incorporates boundary effects through repeated application of the specular reflection boundary condition in \eqref{model-eq-st}. Substituting \eqref{f-Duhamel-1} into itself yields
\begin{equation}\label{double-duhamel}
\begin{split}
f(t, x,v) =
& \frac{1}{\e^2}  \int_{0}^t  \int_{0}^{s}  e^{- \frac{\nu_0}{\e}(t-\tau)}
   \iint_{|v^{\prime}|\leq N, |v^{\prime\prime}|\leq N}
 f \left(\tau,  X_{\mathbf{cl}}(\tau; s, X_{\mathbf{cl}}(s;t,x,v),v^{\prime}), v^{\prime\prime}  \right )
\dd v^{\prime\prime}  \dd v^{\prime}  \dd \tau \dd s   +   \cdots.
\end{split}
\end{equation}
The central insight of \cite{guo2010decay} --- subsequently employed in \cite{Duan-Liu2021,esposito2013non,Esposito2017,Guo-Jang2010,
Guo-Jang-Jiang2010,guo2017regularity,guo2017incompressible}
--- is to gain $L^p$ control via the change of variables
$$
  \left [ v^{\prime}\mapsto z:= X_{\mathbf{cl}}(\tau; s, X_{\mathbf{cl}}(s;t,x,v),v^{\prime})\right].
$$
A critical requirement for this approach is establishing a uniform lower bound on the associated Jacobian:
\begin{equation}\label{low-bound-Jacobian}
\begin{split}
\norm{ \mathbf{J}}:= \Big | \det  \Big  [ \frac{\pt  X_{\mathbf{cl}}(\tau; s, X_{\mathbf{cl}}(s;t,x,v),v^{\prime})}
  {\pt v^{\prime}} \Big ] \Big  | \gtrsim \delta>0
\end{split}
\end{equation}
 away from a small set of parameters $s$. When \eqref{low-bound-Jacobian} holds,
 the $L^\infty$ norm can be controlled as
$$
\|f(t)\|_{L^{\infty}_{x,v}}  \lesssim  \d^{-\frac{1}{p}}\Big  ( \int_{\O}\int_{|v^{\prime\prime}|\leq N}
    | f (t, z, v^{\prime\prime}  ) |^p  \dd v^{\prime\prime} \dd z  \Big )^{\frac{1}{p}}
   + \cdots.
$$
However, for the specular reflection boundary condition in \eqref{model-eq-st}, there is no apparent inductive way to analyze the back-time cycles $\frac{\pt  X_{\mathbf{cl}}(\tau; s, X_{\mathbf{cl}}(s;t,x,v),v^{\prime})} {\pt v^{\prime}} $ inductively with finite bounces, making \eqref{low-bound-Jacobian} extremely difficult to verify.

For the standard Boltzmann equation ($\e=1$) in convex domains with analytic boundary, Guo \cite{guo2010decay} established an asymptotic Jacobian lower bound
\begin{align*}
\left| \det  \left [ \frac{\partial v_k}{\partial v_1}  \right ] \right| \gtrsim \delta > 0 \quad \text{for near-tangential back-time cycles}.
\end{align*}
Kim-Lee \cite{kim2018boltzmann} later removed the analyticity requirement via triple Duhamel expansions while preserving the core strategy.

For hydrodynamic limit problems ($\e \to 0$), precise quantification of the Jacobian lower bound dependence $\d(\e)$ becomes essential --- a stark contrast to standard Boltzmann theory ($\e=1$) \cite{guo2010decay,kim2018boltzmann} where $\delta > 0$ suffices.  This distinction introduces a fundamental difficulty: after multiple specular reflections, the map $[v' \mapsto X_{\mathbf{cl}}(\tau; s, X_{\mathbf{cl}}(s;t,x,v),v^{\prime})]$ develops pathological dependence on $\e$ that precludes asymptotic control and renders the key estimate \eqref{low-bound-Jacobian} unverifiable. Consequently,  the core techniques of \cite{Esposito2017, guo2010decay, kim2018boltzmann} fail  catastrophically for  hydrodynamic limits  involving  specular reflection component.

To overcome this fundamental difficulty, we introduce the \textit{stretching method}: for sufficiently small $\e \ll 1$, we transform the spatial and temporal domains via
\begin{equation}\label{stretch-x-y}
\begin{split}
\O \to \O_\e : = \e^{-1}\O, \ \ \  x\mapsto y:=&\e^{-1}x,\\
 [0,\infty] \rightarrow [0,\infty],\ \ \ \  t\mapsto \bar{t}:=&\e^{-2}t.
\end{split}
\end{equation}
This stretching method enables us to enforce a single-bounce constraint along characteristic trajectories and leads to a uniform-in-$\e$ $L^\infty$ estimate.
One of our main contributions is the following $L^\infty$ estimate for the linear Boltzmann equation on the stretched domain $[0, T_0] \times \Omega_{\varepsilon} \times \mathbb{R}^3$:

\begin{proposition}\label{lemma-fbar-infty-unst-0} \
Let $T_0\geq 1$ be a sufficiently large constant (to be determined later), and let
$\bar{f}$ satisfy
\begin{equation}\label{eq-unst-stretch-t}
\begin{split}
\pt_{\bar{t}}\bar{f} + v\cdot\nabla_y \bar{f}  +  L \bar{f}
   =\e \bar{g}  \ \   &\text{ in } [0,T_0]\times \O_{\e}\times\mathbb{R}^3, \\
    \bar{f} |_{\g_-} = (1-\a) \mathscr{R}\bar{f} + \a \mathscr{P}_\g \bar{f} \ \   &\text{ on } [0,T_0]\times\pt\O_{\e}\times\mathbb{R}^3,\\
  \bar{f} |_{t=0} = \bar{f}_0    \ \ \ \  &\text{ on } \O_{\e}\times\mathbb{R}^3,
\end{split}
\end{equation}
where the transformed functions are defined via the stretching \eqref{stretch-x-y}:
\begin{equation}\label{change-fun-t}
\begin{split}
& \bar{f}(\bar{t},y,v):=f(t,x,v), \quad \bar{f}_0(y,v):=f_0(x,v),\quad \bar{g}(\bar{t},y,v):=g(t,x,v).
\end{split}
\end{equation}
Then there exists a constant $\e_0\in (0,1)$ such that for any $0<\e\leq \e_0$, the following estimates hold for all $\bar{t}\in [0, T_0]$:
\begin{align}
\| \o \bar{f}(\bar{t})\|_{L^\infty_{y,v}(\O_\e\times \mathbb{R}^3)}
\lesssim \
& e^{-\frac{\nu_0}{2}\bar{t}}   \| \o \bar{f}_0\|_{L^\infty_{y,v}(\O_\e\times \mathbb{R}^3)}
+  o(1)  \sup_{0\leq s\leq T_0}\| \o \bar{f}(s)\|_{L^{\infty}_{y,v}  (\O_\e\times \mathbb{R}^3)}\nonumber \\
%-----------------------
%-----------------------
& +  \sup_{0\leq s\leq T_0}  \| \o ^{-1}\P \bar{f}(s)\|_{L^6_{y,v}  (\O_\e\times \R^3)}\label{fbar-infty-unst}\\
& +  \sup_{0\leq s\leq T_0} \|  \o ^{-1} \ip \bar{f}(s)\|_{L^2_{y,v} (\O_\e\times \R^3)}\nonumber\\
%-----------------------
& + \e \sup_{0\leq s\leq T_0} \|\langle v\rangle^{-1}  \o \bar{g}(s)\|_{L^\infty_{y,v} (\O_\e\times \mathbb{R}^3)},\nonumber \\
%-----------------------
%-----------------------
\|\o\bar{f}(\bar{t})\|_{L^\infty_{y,v}(\O_\e\times \mathbb{R}^3)}   \lesssim \;
& e^{ -\frac{\nu_0} {2} \bar{t}}    \| \o\bar{f}_0 \|_{L^\infty_{y,v}(\O_\e\times \mathbb{R}^3)}
+  o(1)  \sup_{0\leq s\leq T_0}\|\o\bar{f}(s)\|_{L^{\infty}_{y,v}  (\O_\e\times \mathbb{R}^3)} \nonumber\\
%-----------------------
& +  \sup_{0\leq s\leq T_0} \|\bar{f}(s)\|_{L^2_{y,v} (\O_\e\times \mathbb{R}^3)}
 +   \sup_{0\leq s\leq T_0} \|\e \langle v\rangle^{-1} \o\bar{g}(s)\|_{L^\infty_{y,v} (\O_\e\times \mathbb{R}^3)}.\label{fbar-infty-L2-unst}
\end{align}
\end{proposition}

The proof is given in Section 2.1. We note that $T_0>0$ creates desired decay property. This approach yields  the first  uniform $L^\infty$ estimate for ``large stretched'' non-convex domains.
\medskip

Applying the transformation \eqref{change-fun-t}  to the model equation \eqref{model-eq-st} yields the equivalent problem on the stretched domain:
%\begin{equation}\label{change-fun}
%\begin{split}
%& \bar{f}(\bar{t},y,v):=f(t, x,v)
%\end{split}
%\end{equation}
\begin{equation}\label{model-eq-st-stretch}
\begin{split}
  \pt_{\bar{t}}\bar{f} + v\cdot\nabla_y \bar{f}   +  \nu_0 \bar{f}
   = \int_{|v^{\prime}|\leq N}  \bar{f}(t,y,v^{\prime})\dd v^{\prime}  \ \ \ \
   &\text{ in }  \mathbb{R}^{+}\times \O_\e\times\mathbb{R}^3, \\
  \bar{f}|_{\g_-}  =   \mathscr{R} \bar{f}  \ \ \ \
  &\text{ on }  \mathbb{R}^{+}\times\pt\O_\e\times\mathbb{R}^3, \\
  \bar{f} |_{t=0} = \bar{f}_0(y,v)    \ \ \ \  &\text{ on } \O_{\e}\times\mathbb{R}^3.
\end{split}
\end{equation}
Crucially, while $\O_\epsilon$ becomes asymptotically large, the outward unit normal remains invariant under this scaling:
\begin{equation}\label{ny=nx}
n(y)=\frac{\nabla_y[\xi(\e y)]}{|\nabla_y[\xi(\e y)])|}
%=\frac{\nabla_x\xi(\e y)}{|\nabla_x\xi(\e y)|}
   =\frac{\nabla_x\xi(x)}{|\nabla_x\xi(x)|}=n(x)
       \;\text { for } x\in\partial\O,\; y=\e^{-1}x\in\partial\O_\e.
\end{equation}
The characteristic trajectories for \eqref{model-eq-st-stretch} are simply
\begin{equation}\label{model-char}
\begin{split}
  &\left[Y(s;\bar{t},y,v), \, V(s;\bar{t},y,v)\right]=[y+v (s-\bar{t}), \, v].
\end{split}
\end{equation}
Denote the first boundary collision along the backward specular trajectory by
 \begin{equation}\label{def-tb ch1}
\begin{split}
&(t_1, y_1):= (\bar{t}-t_{\bf b}(y, v), Y(t_1; \bar{t}, y, v) ),
\end{split}
\end{equation}
where
\begin{equation}\label{def-tb-yb}
\begin{split}
&\tb(y,v):=\inf \{\bar{t}\geq 0:Y(-\bar{t};0,y,v)\notin \Omega \},\\
& \yb(y,v):=Y(-\tb(y,v);0,y,v), \\
&  \vb(y,v):= V( -\tb (y,v) ;0,y,v).
\end{split}
\end{equation}
 From \eqref{model-char} we obtain the relation
\begin{equation}\label{model-y1-y2}
   |y- y_1| = | v (\bar{t}-t_1)|.
\end{equation}
Now consider $(\bar{t}, y, v)\in [0,T_0]\times \O_\e\times \big\{|v|\leq N, \ \big |v\cdot \frac{\nabla_x\xi(\e y)}  {|\nabla_x\xi(\e y)|}\big  |>\eta \big \}$ for sufficiently large constants $T_0, N>0$ and a small constant $\eta>0$. Due to the stretching \eqref{stretch-x-y}, the left-hand side $|y- y_1|$ in \eqref{model-y1-y2} is of order $O(\frac{1}{\e})$, while the right-hand side $|v (\bar{t}-t_1)|$ in \eqref{model-y1-y2} is bounded by $T_0 N$. This implies that, for sufficient small $\e\ll 1$, the specular backward trajectory starting from $(\bar{t},y,v)$ undergoes  at most a single bounce (see  Lemma \ref{1-bounce-claim-model}). Consequently, we can establish a uniform-in-$\e$ Jacobian lower bound analogous to \eqref{low-bound-Jacobian} along this single-bounce trajectory, which ultimately leads to a $\e$-independent $L^\infty$ estimate.
\bigskip

\subsection{Strong Limit Result for the Case $0\leq \a \ll \e$}\
\medskip

In the regime $0\leq \a \ll \e$, we have
\begin{equation}\label{a-e-limit-0}
\displaystyle \l:= \frac{1}{\sqrt{2\pi}}\lim_{\e\to  0} \frac{\a}{\e} =0.
\end{equation}
  Proposition \ref{0-macro-L2-L6-estimate} fails to provide an uniform estimate for $\P f$, as the boundary dissipation in \eqref{model-L2} becomes nearly negligible. For the pure specular reflection case $\a=0$,  uniform estimate for $\P f$ can still be obtained through conservation laws of mass, angular momentum and energy. However, when $\alpha \neq 0$, the latter two conservation laws no longer hold, precluding the control of $\mathbf{P} f$ via this  method.

 To overcome this essential difficulty, we introduce the following {\em rotating Maxwellian}:
\begin{equation}\label{tilde-mu-def}
\begin{split}
\displaystyle \tilde{\mu}=\tilde{\mu}(t,x,v)  :=\frac{\rho(t,x)}{[2\pi T(t)]^{3/2}} \exp \Big(-\frac{\norm{v-\c (t,x)}^2}{2 T(t)}\Big),
\end{split}
\end{equation}
where the temperature is $T(t) := 1 + \th(t)$, the rigid velocity field $\c$ is defined by
\begin{equation}\label{u-definition}
 \c=\c (t,x) :=\begin{cases}
\displaystyle 0 & \text{ if } \text{dim} \mathcal{R}_\Omega = 0,
\\
\displaystyle w(t) Ax  & \text{ if } \text{dim} \mathcal{R}_\Omega = 1 \; ( Ax  \in\mathcal{R}_\Omega),
\\
\displaystyle \sum_{i=1}^3 w_i(t) A_{i}x  & \text{ if } \text{dim} \mathcal{R}_\Omega =3\; ( A_{i}x \in \mathcal{R}_\Omega, i=1,2,3)
\end{cases}
\end{equation}
(see \eqref{R-space-form} and \eqref{def:axiss}), and the density $\rho$ is given by
\begin{equation}\label{rho-def}
\displaystyle \rho= \rho(t,x):= \displaystyle
\frac { \norm{\Omega}\exp \left (\frac{\norm{\c (t,x)}^2}{2T(t)} \right ) }
 { \int_{\Omega} \exp \left(\frac{\norm{\c (t,x)}^2}{2T(t)}\right) \dd x}.
\end{equation}
Here, $\th(t)$, $w(t)$ and $w_i(t)$ ($i\in \{1,2,3\}$)  are scalar functions (to be determined  in Lemma \ref{tildemu - existence}), subject to the initial conditions
\begin{equation}\label{parameters-initials}
\th(0)=0, \quad w(0)=0, \quad  w_{i}(0)=0\; (i=1,2,3).
\end{equation}
In what follows, a summation of the form $\sum w_i A_{i}x $ without explicit indices will denotes $w Ax$ for an axiymmetric domain or $\sum_{i=1}^3 w_i A_{i}x $ for a spherical domain.

We now define the parallel fluctuation field $\tilde{f}$ by
\begin{equation} \label{tildef - definition ch1}
F= \tilde{\mu}+ \e \sqrt{\tilde{\mu}} \tilde{f}, \quad F_0= \tilde{\mu}+ \e \sqrt{\tilde{\mu}} \tilde{f}_0.
\end{equation}
Consequently, the original equation \eqref{F - Boltzmann equation}  can be rewritten in terms of   $\tilde{f}$ as
\begin{equation}\label{tildef - Boltzmanneq ch1}
\begin{split}
\e \pt_t \tilde{f} + v \cdot \nabla_x \tilde{f} + \e^{-1} \tilde{L} \tilde{f} =\tilde{g}   \;\;\;\;  &\text{ in } \mathbb{R}^{+}\times\O\times\mathbb{R}^3, \\
%---------------
 \tilde{f}|_{\g_-} = (1- \a) \mathscr{R} \tilde{f} + \a \tilde{\mathscr{P}}_{\g} \tilde{f} +\a r  \;\;\;\;   &\text{ in } \mathbb{R}^{+}\times\pt\O\times\mathbb{R}^3,\\
%---------------
 \tilde{f}|_{t=0}  =   \tilde{f}_0(x,v)  \;\;\;\;    &\text{ on }   \O \times \mathbb{R}^3,
\end{split}
\end{equation}
with the operators
\begin{equation}\label{tildeg - definition}
\begin{split}
\tilde{\G}(f,g) :=& \frac{1}{\sqrt{\tilde{\mu}}}Q(\sqrt{\tilde{\mu}}f,\sqrt{\tilde{\mu}}g), \quad
\tilde{L}(f) := -\tilde{\G}(\sqrt{\tilde{\mu}},f) -\tilde{\G}(f,\sqrt{\tilde{\mu}}),\\
%-------------------
%-------------------
\tilde{g} :=& \tilde{\G}(\tilde{f},\tilde{f}) - \frac{\pt_t \tilde{\mu}}{\sqrt{\tilde{\mu}}}  -  \e \frac{ \pt_t \sqrt{\tilde{\mu}}} {\sqrt{\tilde{\mu}}} \tilde{f}, \\
%--------------------------
\tilde{\mathscr{P}}_{\g} f :=
&% \frac{1}{\sqrt{\tilde{\mu}}} \mathscr{P} (\sqrt{\tilde{\mu}} f)=
\sqrt{2 \pi} \frac{\mu}{\sqrt{\tilde{\mu}}} \int_{n \cdot u >0}  f\sqrt{\tilde{\mu}(u)} [n \cdot u]  \dd u,  \quad
%-------------------
 r := \frac{1}{\e \sqrt{\tilde{\mu}}}(\mathscr{P}\tilde{\mu} -\tilde{\mu}).
\end{split}
\end{equation}
For the transport operator $v \cdot \nabla_x \tilde{f} $, we have used the identities (valid for all three geometric types of $\O$)
\begin{equation} \label{tildemu - nabla zero ch1}
v \cdot \nabla_x \tilde{\mu} =\frac{1}{T} w(v \cdot Av)  \tilde{\mu} = 0,\quad
v \cdot \nabla_x \tilde{\mu} =\frac{1}{T}  \sum_{i=1}^3 w_i  (v \cdot A_iv)  \tilde{\mu} = 0.
\end{equation}

The null space of $\tilde{L}$ is a five-dimensional subspace of $L^2(\mathbb{R}^3)$ given by
\begin{equation}\label{ker-tilde-P}
\ker   \tilde{L} = \text{span} \big \{1, v-
\c, \norm{v-\c}^2  \big \} \sqrt{\tilde{\mu}}  = \text{span} \big \{ 1,   v ,    |v|^{2}\big \}\sqrt{\tilde{\mu}},
\end{equation}
equipped with orthonormal basis $\{\bar{\chi}_{i}\}_{i=0}^{4}$:
\begin{equation}
\begin{split} \label{base-bar-chi}
&\bar{\chi}_{0} := {\frac{1}{\sqrt{\rho}}}\sqrt{\tilde{\mu}}, \quad
\bar{\chi}_{i} := \frac{(v_i - \c_{i})}{\sqrt{{\rho}T}}\sqrt{\tilde{\mu}} \;\; (i=1,2,3),\quad
\bar{\chi}_{4} := \frac{\norm{v-\c}^2-3T}{\sqrt{6{\rho}}T}\sqrt{\tilde{\mu}}.
\end{split}
\end{equation}
The orthogonal projection of $\tilde{f}$ onto $\ker   \tilde{L} $ is
\begin{equation}\label{Pf-barabc -ch}
  \tilde{\P} \tilde{f}  =  \bar{a} \bar{\chi}_{0} + \sum_{i=1}^{3}\bar{b}_{i} \bar{\chi}_{i} + \bar{c} \bar{\chi}_{4},
\end{equation}
with coefficients
\begin{equation}\label{barabc-def}
    \bar{a}:= \langle \bar{\chi}_{0} , \tilde{f} \rangle,\quad
    \bar{b}:=\langle \bar{\chi}_{i}, \tilde{f}\rangle(i=1,2,3),\quad
     \bar{c}:=\langle \bar{\chi}_{4} , \tilde{f}\rangle.
\end{equation}
We denote by $(\mathbf{I}-\tilde{\P})\tilde{f}$ the projection on the orthogonal complement of $\ker \tilde{L}$.

A crucial observation is the relationship between $f$ and $\tilde{f}$:
\begin{equation}\label{f-relation-tildef}
\tilde{f} = \frac{\mu - \tilde{\mu}}{\e \sqrt{\tilde{\mu}}} + \frac{\sqrt{\mu}}{\sqrt{\tilde{\mu}}} f.
\end{equation}
Moreover, the initial conditions in \eqref{parameters-initials} imply
\begin{equation}\label{tilde-mu-equals-mu-t0}
\tilde{\mu} = \mu,\;\; \tilde{\P}=\P, \;\;  \tilde{\mathscr{P}}_{\g}=\mathscr{P}_{\g} \;\;  \text{ at } t=0.
\end{equation}
Consequently, the two perturbation equations \eqref {f-eq} and \eqref{tildef - Boltzmanneq ch1} actually satisfy the same initial condition:
\begin{equation}\label{tilde-f0-equals-f0-t0}
\tilde{f}_0(x,v)= f_0(x,v).
\end{equation}

We define the instant energy functional
\begin{equation}
\begin{split}
\mathscr{E}_{2}[\tilde{f}](t) :=  &\sup_{0\le s\le t}  \Big\{
\normm{\tilde{f}(s)}_{L^2_{x,v}}^{2} + \normm{\pt_t \tilde{f}(s)}_{L^2_{x,v}}^{2}
+  \norm{\frac{\th(s)}{\e}}^2 + \norm{ \frac{w(s)}{\e}}^2 \Big\}\\
&+ \sup_{0\le s\le t}  \Big\{\norm{ \frac{\pt_t \th(s)}{\e} }^2+  \norm{ \frac{\pt_t w(s)}{\e} }^2  \Big\}.
\end{split}
\end{equation}
%where $\tilde{\mu}$, $\th(t)$ and $w(t)$ will be determined in Lemma \ref{tildemu - %existence}.
The dissipation functional is defined as
\begin{equation}
\begin{split}
\mathscr{D}_{2}[\tilde{f}](t) := &   \int_{0}^{t}  \Big\{
\normm{\tilde{\P} \tilde{f}(s)}_{L^2_{x,v}(\tilde{\nu})}^{2} + \normm{\tilde{\P} \pt_t  \tilde{f}(s)}_{L^2_{x,v}(\tilde{\nu})}^{2}  \Big\}\dd s          \\
 %-------------
 &+\int_{0}^{t}  \Big\{
 \frac{1}{\e^{2}} \normm{\ipt \tilde{f}(s)}_{L^2_{x,v}(\tilde{\nu})}^{2} + \frac{1}{\e^{2}} \normm{\ipt \pt_t \tilde{f}(s)}_{L^2_{x,v}(\tilde{\nu})}^{2}  \Big\}\dd s \\
 %---------------
& + \int_{0}^{t}  \Big\{   \frac{\a}{\e} \norm{\tilde{f}(s)}_{L^2_{\g_{+}}}^2
 + \frac{\a}{\e} \norm{\pt_t \tilde{f}(s)}_{L^2_{\g_{+}}}^2+ \frac{\a}{\e}  \norm{\frac{\th(s)}{\e}}^2
+\frac{\a}{\e}   \norm{ \frac{w(s)}{\e}}^2 \Big\}\dd s\\
&+\int_{0}^{t}  \Big\{
\frac{\a}{\e}    \norm{ \frac{\pt_t \th(s)}{\e} }^2
+ \frac{\a}{\e}   \norm{ \frac{\pt_t w(s)}{\e} }^2  \Big\}\dd s.
\end{split}
\end{equation}
The total energy functional is defined by
\begin{equation}\label{total-energy}
\begin{split}
\normmm{\tilde{f}}_{2}(t) :=\,& \mathscr{E}_{2}^{\frac{1}{2}}[\tilde{f}](t)  + \mathscr{D}_{2}^{\frac{1}{2}}[\tilde{f}](t) +
\e^{\frac{1}{2}}\sup_{0\leq s\leq t}\normm{
%---------
{ \o f(s) }
%---------
}_{L^{\infty}_{x,v}}\\
&+ \e^{\frac{3}{2}}\sup_{0\leq s\leq t} \normm{
%---------
{ \o \pt_t f(s)}
%----------
}_{L^{\infty}_{x,v}}
 + \sup_{0\leq s\leq t}\normm{\tilde{\P} \tilde{f}(s)}_{L^6_{x,v}}.
\end{split}
\end{equation}
The corresponding norm for the initial data is
\begin{equation}\label{initial-data-total-norm-tilde defn}
\begin{split}
\left[\!\left[\tilde{f}_{0}\right]\!\right]_{2} :=&
\normmm{\tilde{f}}_{2}(0)  + \e^{-1} \normm{\ipt \tilde{f}_{0}}_{L^2_{x,v}(\tilde{\nu})}  +  \Big(\frac{\a}{\e}\Big)^{\frac{1}{2}}  \norm{(1-\tilde{\mathscr{P}}_{\g} )\tilde{f}_{0} }_{L^2_{\g_{+}}} \\
&
+\normm{v\cdot \nabla_x \tilde{f}_0}_{L^2_{x,v}}+\normm{v\cdot \nabla_x \pt_t\tilde{f}_0}_{L^2_{x,v}}\\
= &  \left[\!\left[f_{0}\right]\!\right]_{1},
\end{split}
\end{equation}
 where the final equality follows from \eqref{tilde-mu-equals-mu-t0} and \eqref{tilde-f0-equals-f0-t0}.
\medskip

We now state the second main result for the regime $0\leq \a \ll \e$.
\medskip

\begin{theorem}[Case $0\leq \a \ll \e$] \label{main-th-2} \
Let $F_0 =\mu + \e \sqrt{\mu} f_0 \geq 0 $. Then there exists $\e_0>0$ such that for every $0<\e<\e_0$, if the initial  fluctuation $f_0$ satisfies
\begin{align}\label{initial-data-total-norm-tilde}
  \left[\!\left[{f}_{0}\right]\!\right]_{1}
% \equiv \left[\left[\frac{F_0-\mu}{\e \sqrt{\mu}}\right]\right]_{1}
 \leq  \d_0
\end{align}
for some small constant $\delta_0>0$ independent of $\e$ (the same initial condition as in \eqref{0-initial-data-total-norm}), then the Boltzmann equation \eqref{F - Boltzmann equation} admits a unique global solution $F =  \tilde{\mu} + \e \sqrt{\tilde{\mu}} \tilde{f} \geq 0$ satisfying the uniform bound
\begin{align}\label{uniform-bound-tilde}
\normmm{\tilde{f}}_{2}(\infty) \leq  C
% \left[\!\left[\tilde{f}_{0}\right]\!\right]_{2}
\left[\!\left[f_{0}\right]\!\right]_{1}
\end{align}
for some constant $C>0$ independent of $\e$.

Moreover, if the strong initial convergence \eqref{initial-converge-condition} holds,
%if there exist fluid initial data $(\varrho_0, u_0, \vartheta_0)\in %\mathbb{H}_{\vartheta}\times\mathbb{H}_{u}\times\mathbb{H}_{\vartheta}$
%such that
%\begin{align}
%\begin{split}\label{initial-converge-condition-tilde}
%&\tilde{f}_0\to {f}^{*}_0=\;\Big(\varrho_0+u_0\cdot v+\vartheta_0\frac{|v|^2-3}{2}\Big) %\sqrt{\mu}
%\ \  \text { strongly in  }  L^2(\O \times \mathbb{R}^3) \text{ as } \e \to 0,
%\end{split}
%\end{align}
then the convergence results \eqref{tilde-f-strong-convergence}--\eqref{momentums-strong-convergence} are also valid.
 Here, $(\varrho, u, \vartheta)\in C(\mathbb{R}_+, L^2(\O))\cap L^2(\mathbb{R}_+, H^1(\O))$ is the unique weak solution  of the INSF system \eqref{INSF-unst}, now supplemented with the perfect Navier slip boundary condition:
\begin{equation}\label{Navier-bdy-unst-perfect}
\begin{split}
 \left [ \big(\nabla_x u + {(\nabla_x u)}^{\mathrm{T}}\big) \cdot n
   \right ]^{\mathrm{tan}} = 0, \ \ \  u\cdot n =0 \ \ \ \
    &\text{ on }\mathbb{R}^{+}\times\pt\O, \\
  \partial_n \vartheta  = 0 \ \ \ \  &\text{ on }\mathbb{R}^{+}\times\pt\O.
\end{split}
\end{equation}
\end{theorem}

\medskip

Proof of Theorem \ref{main-th-2} will be presented in Section \ref{proof-main-2}.

\bigskip

\subsection{Methodology 2: Dissipative Decomposition Mechanism}\
\medskip

To overcome the difficulties caused by the nearly negligible boundary dissipation in \eqref{model-L2} and the loss of conservation laws of angular momentum and energy, we uncover a dissipative decomposition mechanism via the construction of a rotating Maxwellian. More precisely, we design the rotating Maxwellian $\tilde{\mu}$ as in \eqref{tilde-mu-def} and reformulate the Boltzmann solution $F$ around $\tilde{\mu}$ via \eqref{tildef - definition ch1}. This decomposition splits the original equation  \eqref{F - Boltzmann equation} into two dissipative subsystems: one for spatially averaged macroscopic variables $(\c,\theta)$, and another for the fluctuation $\tilde{f}$ satisfying the following conservation laws of mass, angular momentum and energy:
\begin{equation} \label{Psi - conservation law}
\begin{split}
&\iint_{\Omega\times\R^3}  \sqrt{\tilde{\mu}}\tilde{f} \dd v \dd x = 0, \\
&\iint_{\Omega\times\R^3} Ax \cdot v \sqrt{\tilde{\mu}}\tilde{f} \dd v \dd x = 0 \;\;  \text{ for all } Ax \in \mathcal{R}_{\Omega},\\
&\iint_{\Omega\times\R^3} \norm{v}^2 \sqrt{\tilde{\mu}}\tilde{f} \dd v \dd x = 0,
\end{split}
\end{equation}
 guaranteed by \eqref{initial-conservation-F}. These conservation laws allow us to control the macroscopic components  $ \int_0^t \| \P f \|_{L^2_{x,v}}^2 $ and $\normm{\tilde{\P} \tilde{f}}_{L^6_{x,v}}$ even with weak boundary dissipation via a test function approach \cite{esposito2013non,Chen2024}.

The velocity field $\c(t,x)$ and the temperature deviation $\theta(t) = T(t)-1$ are determined via the implicit function theorem (with density $\rho$ depending on $\c$ and  $\theta$ through \eqref{rho-def}), from the full conservation laws of the original solution $F$:
\begin{equation} \label{F - consercation law - ch 1}
\begin{split}
&\iint_{\O\times\mathbb{R}^3}  F(t) \dd v\dd x= |\O|, \\
 &\iint_{\O\times\mathbb{R}^3} Ax \cdot v F(t) \dd v\dd x=  \int_{\Omega} \rho Ax \cdot \c \dd x\; \;\text{ for all } Ax\in \mathcal{R}_\Omega,
 \\
 & \iint_{\O\times\mathbb{R}^3} \norm{v}^2 F(t) \dd v\dd x=\int_{\Omega} (3\rho T + \rho \norm{\c}^2)\dd x,
\end{split}
\end{equation}
  as shown in Lemma \ref{tildemu - existence}.  Crucially, $\theta^2$ and $|\c |^2$ satisfy a dissipative ODE system:
\begin{equation}
\begin{split}
&\frac{3}{2} \pt_t \int_{\Omega} \theta^2 \dd x  + \frac{\a}{\e\sqrt{2\pi}} \int_{\pt \Omega} 4\theta^2 \dd S_x + {\a}\iint_{\g_{+}} ( \norm{v}^2 -4) \sqrt{\tilde{\mu}} \tilde{f} \th \dd \g   =\text{higher-order terms}, \\
%----------------------
&\frac{1}{2}\pt_t \int_{\Omega} \norm{\c}^2 \dd x  + \frac{\a}{\e\sqrt{2\pi}} \int_{\pt \Omega} \norm{\c}^2 \dd S_x
+ {\a}\iint_{\g_{+}} (\c \cdot v) \sqrt{\tilde{\mu}} \tilde{f} \dd \g  =\text{higher-order terms},
\end{split}
\end{equation}
derived in Propositions \ref{th tht w wt - ODE}.

 Although $\tilde{\mu}$ and $\tilde{\mathbf{P}}$ do not commute with $\partial_t$ and $\nabla_x$, a careful analysis shows that
$$
 v\cdot \nabla_x\tilde{\mu}=0, \quad [\partial_t, \tilde{\mathbf{P}}] \approx O(\alpha), \quad \partial_t \tilde{\mu}\approx O(\alpha).
$$
Combining these observations with a standard energy estimate yields
\begin{equation}
\begin{split}
&\frac{1}{2} \pt_t \normm{\tilde{f}}_{L^2_{x,v}}^2 + \frac{1}{\e^2} \iint_{\Omega \times \R^3} \tilde{f} \tilde{L} \tilde{f} \dd x \dd v + \frac{3}{2} \pt_t \int_{\Omega} \Big(\frac{|\theta|}{\e}\Big)^2 \dd x + \pt_t \int_{\Omega} \Big(\frac{\norm{\c}}{\e}\Big)^2 \dd x\\
& +  \frac{\a(2-\a)}{\e}  \iint_{\g_{+}} \Big[ \frac{1}{2}\frac{\theta}{\e}( \norm{v}^2 -4) \sqrt{\tilde{\mu}}   + v\cdot \frac{\c}{\e}\sqrt{\tilde{\mu}}  +  [(1- \tilde{\mathscr{P}}_{\g})\tilde{f}]  \Big]^2 \dd \g \\
\le & \frac{1}{\e}\iint_{\Omega \times \R^3} \norm{ \tilde{f} \tilde{g} } \dd x \dd v +\text{higher-order terms}.
\end{split}
\end{equation}
The boundary dissipation in this estimate covers all directions except those parallel to $(|v|^2-4)\sqrt{\tilde{\mu}}$, $v\cdot Ax\sqrt{\tilde{\mu}}$, and $\tilde{\mathscr{P}}_\gamma$. Applying Ukai's trace lemma to these rapidly decaying directions ultimately yields complete boundary dissipation (see Proposition \ref{tildef tildeft - Energy estimate}).

For brevity, we state only the key a priori estimates. Assume that \eqref{tildef - Boltzmanneq ch1} admits a solution $\tilde{f}(t)$ on $[0, T]$ with $0 < T \leq \infty$. To simplify the derivation, we impose the following a priori assumption:  there exists a sufficiently small constant $\delta_1>0$ (to be chosen later), independent of $\e$,  such that
\begin{equation}\label{theta-u-smallness-assumption}
\begin{split}
 \sup_{0\leq t\leq T}\Big(\frac{\norm{\th (t)}}{\e}+\frac{\norm{w (t)}}{\e}+\frac{\norm{\pt_t \th (t)}}{\e}+\frac{\norm{\pt_t w (t)}}{\e}\Big)\leq \delta_1.
\end{split}
\end{equation}

Our main estimate on the macroscopic part $\tilde{\P} \tilde{f}$ in the regime $0\leq \a \ll \e$ is summarized as follows.

\begin{proposition} \label{Psi - L2 and L6 estimate}
Let $\tilde{f}$ be a solution of \eqref{tildef - Boltzmanneq ch1} satisfying
the conservation laws of mass, angular momentum and energy given in \eqref{Psi - conservation law}.
%\begin{equation} \label{Psi - conservation law}
%\begin{split}
%&\iint_{\Omega\times\R^3}  \sqrt{\tilde{\mu}}\tilde{f}\dd v \dd x = 0, \\
%&\iint_{\Omega\times\R^3}Ax \cdot v \sqrt{\tilde{\mu}}\tilde{f} \dd v \dd x = 0, \quad \forall A x\in \mathcal{R}_{\Omega},\\
%&\iint_{\Omega\times\R^3} \norm{v}^2 \sqrt{\tilde{\mu}}\tilde{f} \dd v \dd x = 0.
%\end{split}
%\end{equation}
 Then, under the a priori assumption \eqref{theta-u-smallness-assumption}, the following estimates hold  for all $0\leq s \leq t \leq T$:
\begin{align}
\int_{s}^{t} \normm{\tilde{\P}\tilde{f}}_{L^2_{x,v}}^{2} \dd \tau \lesssim
&\e \big[ \tilde{G}_0(t) - \tilde{G}_0(s) \big ] + \a^2\int_{s}^{t} \Big[ \norm{\tilde{f}}_{L^2_{\g_{+}}}^2 +\norm{r}_{L^2_{\g_{-}}}^2 \Big]\dd \tau
+  \int_{s}^{t}   \normm{\tilde{\nu}^{-\frac{1}{2}} \tilde{g} }_{L^2_{x,v}}^{2} \dd \tau \nonumber \\
&+\a^2\e^2 \int_{s}^{t} \norm{\tilde{f}}_{L^2_{\g_{+}}}^2 \normm{ \tilde{f}}_{L^{2}_{x,v}}^2 \dd \tau  + \int_{s}^{t} \normm{\e^{-1}\ipt \tilde{f} }_{L^2_{x,v}(\tilde{\nu})}^2 \dd \tau,  \label{P-tilde-f-macro-L2}
\\
%-----------------%-----------------
%-----------------%-----------------
\normm{\tilde{\P} \tilde{f} }_{L^6_{x,v}} \lesssim
& \e\normm{ \pt_t \tilde{f}}_{L^{2}_{x,v}} + \a \norm{\tilde{f}}_{L^2_{\g_{+}}} +\a\norm{r}_{L^4_{\g_{-}}}+  \a  \norm{ \tilde{f} }_{L^2_{\g_{+}}}^{\frac{1}{2}}
\normm{\omega^{\frac{1}{2}} \tilde{f} }_{L^{\infty}_{x,v}}^{\frac{1}{2}}+ \e^{\frac{1}{2}}\mathfrak{h}_{1} \normm{\e^{\frac{1}{2}} \omega^{\frac{1}{2}} \tilde{f} }_{L^{\infty}_{x,v}}  \nonumber\\
&
+ \normm{\ipt  \tilde{f} }_{L^{6}_{x,v}}  + \normm{\e^{-1}\ipt  \tilde{f} }_{L^{2}_{x,v}(\tilde{\nu})} +\normm{\tilde{\nu}^{-\frac{1}{2}} \tilde{g} }_{L^2_{x,v}}. \label{P-tilde-f-macro-L6}
\end{align}
where $\tilde{G}_0(t)\lesssim \normm{\tilde{f}(t)}_{2}^2$.
\end{proposition}

Proposition \ref{Psi - L2 and L6 estimate} (proved in Section 4.3) supplies the essential dissipative control on the macroscopic component $\tilde{\mathbf{P}}\tilde{f}$, thereby completing the uniform energy framework for the regime $0 \leq \a \ll \e$.

\bigskip

%%%%%%%%%%%%%%%%%%%%%%%%%%%%%%%%%%%%%%%%%%%%%%%%%%%%%%%
%%%%%%%%%%%%%%%%%%%%%%%%%%%%%%%%%%%%%%%%%%%%%%%%%%%%%%%

\subsection{Background and Progress}\
\medskip

The derivation of fluid dynamical equations from kinetic theory
constitutes a cornerstone of mathematical physics since the pioneering works of Maxwell and Boltzmann. Maxwell \cite{Maxwell} and Boltzmann \cite{BOLTZMANN2003} demonstrated that microscopic particle interactions could explain macroscopic phenomena (e.g., viscosity and thermal conductivity), providing foundational insights into molecular dynamics. Based on these foundations, Hilbert formalized the kinetic-continuum connection through his Sixth Problem \cite{Hilbert1902}. His pioneering work \cite{Hilbert1912} established mathematical links between the Boltzmann equation and hydrodynamic models, thereby inaugurating sustained research into hydrodynamic limits.

Building on Hilbert's foundational vision, rigorous hydrodynamic limits of the Boltzmann equation have been established across several principal scaling regimes:
(1) Compressible Euler limit for classical  and renormalized solutions \cite{Caflisch1980, Guo-Huang-Wang2021, Guo-Jang2010, Guo-Jang-Jiang2010, Jiang-Luo-Tang2024, Nishida1978, Ukai-Asano1983, Yu2005};
(2) Compressible Navier-Stokes approximation via Chapman-Enskog expansion \cite{Duan-Liu2021, Jiang-Luo2025, Kawashima-Matsumura-Nishida1979, Liu-Yang-Zhao2014};
(3) Incompressible Euler limit confirmed for renormalized solutions \cite{de1989incompressible, Saint-Raymond2003, Saint-Raymond2009} and analytic solutions in half-space \cite{Cao-Jang-Kim2023, Jang-Kim2021, Kim-Nguyen2024}.
In contrast, the incompressible Navier-Stokes-Fourier limit --- characterized by diffusive scaling and low Mach asymptotic --- demands specialized analysis due to its physical prevalence and mathematical depth.
As the most extensively studied hydrodynamic limit paradigm, the INSF limit exhibits fundamental methodological divergences dictated by domain topology: whole-space and periodic domains; domains with boundary. We now delineate seminal advances in these settings.

For the whole space or periodic domains, the INSF limit  has attained substantial resolution through three frameworks:\\
\noindent (a) Renormalized solutions framework.
Bardos-Golse-Levermore \cite{Bardos1991, Bardos1993} pioneered the convergence of DiPerna-Lions renormalized solutions \cite{DiPerna1989} to Leray-Hopf weak solutions of INSF, contingent on specific a priori assumptions. Subsequent research \cite{Bardos1998, Bardos2000, Golse2001, DavidLevermore2009, Saint-Raymond2003} progressively weakened these constraints. A foundational breakthrough came with Golse-Saint-Raymond's complete proof for bounded collision kernels \cite{Golse2003}, which catalyzed extensions to more general kernels \cite{Golse2009, DavidLevermore2009}, see also comprehensive surveys in \cite{SaintRaymond2009, villani2002review}.\\
\noindent (b)  Classical solutions framework. DeMasi-Esposito-Lebowitz adapted Caflisch's approach \cite{Caflisch1980} to examine the INSF limit \cite{de1989incompressible}. Guo \cite{guo2006boltzmann} later provided rigorously justification, incorporating higher-order correction for both Boltzmann cutoff potentials and Landau collision kernels. Related developments are documented in \cite{bardos1991classical, Briant2015, Briant-Merino-Aceituno-Mouhot2015}.\\
\noindent (c) Strong solutions framework. Inspired by Guo-Jang's $L^2$-$W^{1,\infty}$ framework \cite{Guo-Jang2010}, Gong-Zhou-Wu studied the diffusive limit in $\mathbb{R}^3$ within the $H^1$-$W^{1,\infty}$ framework \cite{GZW2021} (for the more general Vlasov-Poisson-Boltzmann system), while Liu-Yu investigated the diffusive limit in $\mathbb{T}^3$ within the $L^2$-$L^\infty$ framework \cite{LY2024}.

 For domains with boundary, the analysis of INSF limit presents significantly greater complexity than the whole-space or periodic settings. Boundary interactions inherently degrade the regularity of the Boltzmann solutions \cite{guo2017regularity}, precluding classical solutions in general domains. Consequently, research is confined to two frameworks:\\
\noindent (1)  Renormalized solutions framework. Masmoudi-Saint-Raymond \cite{Masmoudi2003} established hydrodynamic limit of renormalized solution \cite{mischler2010kinetic} to the linear Stokes-Fourier system for  the Maxwell boundary. Then Saint-Raymond  extended to the weak INSF limit for cutoff hard potentials \cite{SaintRaymond2009}. Later on, by constructing boundary layer Jiang-Masmoudi  proved  weak convergence for all $\alpha \in [0,1]$ and strong convergence \emph{only} for $\alpha \sim \e^{1/2}$ \cite{jiang2017boundary}.

\noindent (2) Strong solutions framework.
Pioneered by Guo's $L^2$-$L^\infty$ theory \cite{guo2010decay}, this approach achieved critical advances under diffuse boundary conditions.
For interior domains,  Esposito-Guo-Kim-Marra \cite{Esposito2017} justified the steady/unsteady limit by using an $L^2$-$L^6$-$L^\infty$ approach, while Esposito-Guo-Marra-Wu \cite{Esposito-Guo-Marra-Wu2023} and Wu-Ouyang \cite{Wu-Ouyang2022} conducted detailed boundary-layer analyses. For exterior domains,  progress was made by Esposito-Guo-Marra for steady flows \cite{esposito2018hydrodynamic} and by Jung \cite{jung2023diffusivelimitboltzmannequation} for unsteady flows.

However, extant results on diffusive limit with Maxwell boundary --- including the significant works \cite{jiang2017boundary,SaintRaymond2009} --- remain confined to weak convergence within renormalized solutions framework, with strong convergence established only for $\a \sim \e^{1/ 2 }$.
In this work, we establish strong convergence to the INSF system within strong solutions framework for the full range $\a \in[0,1]$. This result encompasses both the pure specular reflection case ($\alpha=0$) and the challenging near-specular regime ($0 < \alpha \ll \varepsilon$), which had previously resisted analysis.
%Two new core methodological innovations underpin our results:\\
%\noindent (I) \emph{Stretching method for $L^\infty$ estimate}.  Indeed, $L^\infty$ estimate for specular reflection boundary in general domains remains a longstanding open problem. By exploiting the scaling regime $\e \ll 1$, we propose a novel \textit{stretching method} to enforce an \textit{at most single-bounce} constraint along characteristics, yielding the first uniform-in-$\e$ $L^\infty$ estimate for arbitrary domains.\\
%\noindent (II) \emph{Dissipative decoupling mechanism for near-specular regime $0<\a \ll \e $}.  In fact, the critical near-specular regime $0<\a \ll \e $ presents dual challenges: Energy estimates yield only $O(\a)$ boundary dissipation, meanwhile lack of angular momentum and energy conservation laws induces uncontrollable boundary growth in macroscopic estimate.
%We overcome such dual challenges through a canonical decomposition, revealing a hidden dissipative  subsystem governed by conservation laws (which ensure controllable boundary terms) coupled with mean-field ODEs governing spatial averages of macroscopic quantities. This architecture culminates in overcoming the boundary growth.\\
%\noindent Integrating advances (I) and (II), we establish the $\e$-uniform global estimate within the $L^2$-$L^6$-$L^\infty$ framework \cite{Esposito2017}, and eventually prove strong convergence from the Boltzmann equation with Maxwell boundary \eqref{F - Boltzmann equation} to INSF for all $\a\in [0,1]$.
\bigskip

\subsection{Notations}\
\medskip

Throughout this paper we adopt the following asymptotic  conventions:
\begin{itemize}%[topsep=0pt, itemsep=4pt, leftmargin=*]
\item [$\cdot$] $C$ denotes a generic positive constant independent of $\e$ and $\a$;
\item [$\cdot$]  $X \lesssim Y$ indicates $X \leq C Y$ for some constant $C>0$ independent of $\e$ and $\a$;
\item [$\cdot$]  $X \approx Y$ denotes $X \lesssim Y$ and $Y \lesssim X$;
\item [$\cdot$]  $X \lesssim_\beta Y$ denotes dependence on parameter $\beta$;
\item [$\cdot$]   $o(1)$ represents a small constant independent of $\e$ and $\a$;
\item [$\cdot$]  $\ll 1$ signifies a sufficiently small positive bound.
\end{itemize}

For $1 \leq p \leq \infty$, we define
\begin{itemize}
\item [$\cdot$] $\| \cdot \|_{L^p_{x,v}}$,  $\| \cdot \|_{L^p_{x}}$ or $\| \cdot \|_{L^p_{v}}$ denote $L^p(\Omega \times \mathbb{R}^3)$, $L^p(\Omega)$ or $L^p( \mathbb{R}^3)$ norms;
%\item [$\cdot$] $\| \cdot \|_{L^p_y}$ specifies $L^p(Y)$ norm over variable $y$;
\item [$\cdot$] $\| \cdot \|_{L^p_{x} L^q_{v}} := \big\| \| \cdot \|_{L^q_{x}} \big\|_{L^p_{v}}$ for mixed norms;
\item [$\cdot$] $\| \cdot \|_{L^p(m)} := \| m^{1/2} \cdot \|_{L^p}$ with weight $m$;
\item [$\cdot$] $\langle \cdot, \cdot \rangle$ : $L^2(\mathbb{R}^3_v)$ inner product;
\item [$\cdot$] $\langle v \rangle := (1 + |v|^2)^{1/2}$.
\end{itemize}

Boundary measure and integrals are denoted by
\begin{itemize}
 \item [$\cdot$] $d\gamma := |n \cdot v| \dd v \dd S_x $ (surface measure);
 \item [$\cdot$] $|f|_{p,\g_{\pm}} := \big( \int_{\g_{\pm}} |f|^p \dd \gamma \big)^{1/p}$ for $1 \leq p < \infty$;
%    \item [$\cdot$] $|f|_p := \left( \int_{\gamma} |f|^p d\gamma \right)^{1/p}$ for $1 %\leq p < \infty$;
 \item [$\cdot$]$|f|_\infty := \text{esssup}_{(x,v)\in\gamma} |f(x,v)|$;
 \item [$\cdot$]  $| \cdot |_{L^p_{x}}$ denotes $L^p(\pt\Omega)$ boundary norm;
\end{itemize}

For the perfect Navier slip boundary condition $\l=0$ (arises when $0\leq \a \ll \e$), domain symmetry $\Omega$ also affects the uniqueness of solutions to the INSF system. We define the admissible function spaces for initial data:
\begin{equation}\label{initial-space-Navier}
\begin{split}
&\mathbb{H}_{u}:= \begin{cases}
     \big \{u\in L^2(\O): \nabla_x\cdot u=0\big\}
&\begin{array}{l}
\displaystyle \text{if } 0<\l \leq \infty, \text{ or if } \l=0  \\ \text{and }\O \text{ is non-axisymmetric};
\end{array}\\[2mm]
%---------------------
     \big \{u\in L^2(\O): \nabla_x\cdot u=0, \int_{\O}u\cdot R\dd x=0\; \forall R\in \mathcal{R}_\Omega \big\}  &\begin{array}{l}\displaystyle  \text{if }  \l=0
     \text{ and }  \O \text{ is axisymmetric}\\
\text{or spherical},
\end{array}
                  \end{cases} \\[3mm]
&\mathbb{H}_{\vartheta}:=\begin{cases}
       L^2(\O)
&\displaystyle \text{if }  0<\l \leq \infty;\\
     \big \{\vartheta\in L^2(\O): \int_{\O}\vartheta\dd x=0 \big\}
&\displaystyle \text{if } \l=0,\\
                   \end{cases}
\end{split}
\end{equation}
where $R = R(x)$ denotes basis element of $\mathcal{R}_\Omega$ (see \eqref{RO}), which generates non-trivial special solutions to the INSF system under perfect Navier slip boundary $\l=0$.

\bigskip

The remainder of this paper is organized as follows. Section 2 presents $L^\infty$ estimates for the linear Boltzmann equation on the stretched domain.
Section 3 establishes uniform-in-$\e$ global estimates and the strong convergence for the case $\e \lesssim \a \leq 1$. Sections 4 addresses the strong convergence for the case $ 0\leq \a\ll \e$. Technical supporting results are collected in the appendices: Appendix A provides an $L^2L^3$ estimate, Appendix B gives the uniqueness of weak solutions to the INSF system, and Appendix C contains auxiliary facts on Gaussian integration and elliptic estimates.
\bigskip

%%%%%%%%%%%%%%%%%%%%%%%%%%%%%%%%%%%%%%%%%%%%%%%%%%%%%%%
%%%%%%%%%%%%%%%%%%%%%%%%%%%%%%%%%%%%%%%%%%%%%%%%%%%%%%%

\noindent
\section{$L^\infty$ Estimate} \
\medskip

This section establishes the $L^\infty$ estimate for the linear Boltzmann equation \eqref{eq-unst-stretch-t} on the stretched domain $[0, T_0] \times \Omega_{\varepsilon} \times \mathbb{R}^3$.
The main result is Proposition \ref{lemma-fbar-infty-unst-0}, whose proof is presented at the end of the section after several preparatory lemmas.
\medskip

For the linearized Boltzmann operator $L$ defined in \eqref{L-Definition}, it is standard that
${L} {f} = {\nu} {f} -{K} {f}$, where the collision frequency $\nu$ and the compact operator ${K}$ on $L^2(\R^3_{v})$ are given by
\begin{equation}\label{K-def-absolute}
\begin{split}
&\nu={\nu}(v) := \frac{1}{\sqrt{ {\mu}}}Q_{-} (\sqrt{ {\mu}},  {\mu}) = \int_{\R^{3}}\int_{\mathbb{S}^2} \norm{(v-u) \cdot \omega}  {\mu}(u) \dd \omega \dd u,\\
%--------------
&K f := \frac{1}{\sqrt{\mu}} \left[Q_{+} (\mu, \sqrt{\mu} f)+ Q_{+} (\sqrt{\mu} f, \mu) - Q_{-} (\mu, \sqrt{\mu} f)\right]
=\int_{\R^{3}}[k_{1}(v,u)-k_{2}(v,u)] f(u) \dd u.
\end{split}
\end{equation}
For hard sphere cross sections, there exist positive constants $C_0$ and $C_1$ such that
\begin{equation*}
\begin{split}
& \nu_0\leq C_{0} \inn{v} \le \nu(v) \le C_{1} \inn{v},
\end{split}
\end{equation*}
%\begin{equation*}
%\begin{split}
%\big(f,Lg\big)_{L^2_v}=\big(Lf,g\big)_{L^2_v}\; \text{ on }\; D_{L} = & \big\{f \in %L^2(\R^3_{v}) | \; \nu^{1/2} f \in L^{2}(\R^3_{v}) \big\},\\
%\big(f,\tilde{L}g\big)_{L^2_v}=\big(\tilde{L}f,g\big)_{L^2_v} \; \text{ on }\; %D_{\tilde{L}} = & \big\{\tilde{f} \in L^2(\R^3_{v}) |\; \tilde{\nu}^{1/2} \tilde{f} \in %L^{2}(\R^3_{v})\big \}.
%\end{split}
%\end{equation*}
 with the uniform lower bound $\nu_0$. The operator $L$ is symmetric with the spectral inequality:
\begin{equation*}%\label{positive-definiteness}
\begin{split}
\langle f,Lf \rangle_{2} \gtrsim& \normm{\ip f}_{L^2_v(\nu)}^2 \; \text{ for } \; f\in D_{L}=\big\{f \in L^2(\R^3_{v}) | \; \nu^{1/2} f \in L^{2}(\R^3_{v}) \big\}.
\end{split}
\end{equation*}

Multiplying equation \eqref{eq-unst-stretch-t}  by the weight function $\o$ defined in \eqref{weight-w} yields the equivalent formulation
\begin{equation}\label{eq-unst-stretch-t-h-with-K}
\begin{split}
\pt_{\bar{t}}h + v\cdot\nabla_y h     +   \nu(v)h= \o K( \o ^{-1}h)+\e \o  \bar{g}
 \ \  & \text{ in } [0, T_0]\times \O_{\e}\times\mathbb{R}^3, \\
h \big|_{\g_-} = (1-\a) \mathscr{R}h + \a \o \sqrt{\mu} \int_{n(y)\cdot u>0}h \dd \sigma \ \   & \text{ on }[0, T_0]\times \pt\O_{\e}\times\mathbb{R}^3,\\
h |_{t=0} = h_0    \ \   & \text{ on } \O_{\e}\times\mathbb{R}^3.
\end{split}
\end{equation}
Here and in the following, we use the notations
\begin{align}
 &h(\bar{t},y,v):=\o \bar{f}(\bar{t},y,v),\qquad  h_0(y,v):=\o \bar{f}_0(y,v), \label{def-h}\\
%---------------
 &\dd \sigma:=\o ^{-1}\sqrt{2\pi}{\mu^{\frac{1}{2}}} [n(y)\cdot u] \dd u, \qquad
 C_*:= \int_{n(y)\cdot u>0}  \dd \sigma. \label{def-C-star}
\end{align}

Given $(\bar{t},y,v)\in [0,T_0]\times\bar{\O}_\e\times \mathbb{R}^3$, recall the characteristic trajectory \eqref{model-char}.
% Let $(t_1, y_1, v_1)$ be the first bounce time, position and velocity
% of the characteristic equation
%\begin{equation}\label{model-char-eq}
%\begin{split}
%  &\frac{\dd Y(s,\bar{t},y,v)}{\dd s}=v,\quad Y(\bar{t};\bar{t},y,v)=y.
%\end{split}
%\end{equation}
%along the backward trajectory (cf. \eqref{def-tb ch1}), where
%\begin{equation}\label{def-bounce-1-2}
%\begin{split}
%&t_1=\bar{t}-\tb(\bar{t},y,v),\;\;
%y_1=Y(t_1;\bar{t},y,v),\;\;
%v_1=\left\{
%             \begin{array}{ll}
%             R_{y_1}(v), & \hbox{\text{specular reflection};} \\
%             V_{1}^*, & \hbox{\text{diffuse reflection}}
%             \end{array}
%             \right.
%\end{split}
%\end{equation}
Let $(t_{k+1}, y_{k+1}, v_{k+1})$ denote the $(k+1)$-th ($k\in\mathbb{N}$) bounce along the backward trajectory (cf. \eqref{def-tb ch1} and \eqref{def-tb-yb}):
%, be the first bounce time, position and velocity
% of the characteristic equation
\begin{equation}\label{def-bounce-1-2}
\begin{split}
&t_{k+1}=t_k-\tb(t_k,y_k,v_k),\;\;
y_{k+1}=Y(t_{k+1};t_k,y_k,v_k),\;\; %=\yb(t_k,y_k,v_k)=y_k+(t_{k+1}-t_k)v,
v_{k+1}=\left\{        \begin{array}{ll}
                         R_{y_{k+1}}(v_k), & \hbox{\text{specular reflection};} \\
                         v_{k+1}^*, & \hbox{\text{diffuse reflection},}
                         \end{array}
                         \right.
\end{split}
\end{equation}
where we set $(t_0, y_0, v_0):=(\bar{t}, y, v)$. This yields a sequence $t_{k+1} < t_k < \cdots < t_2 < t_1 < t_0=\bar{t} \leq T_0$.

Because $\pt\O\in C^3$ is compact and $\nabla_x \xi\neq 0$ on $\pt\O$, there exist positive constants $0<C_{\xi_1}<C_{\xi_2}$, independent of $\e$,  such that
\begin{equation}\label{xi-low-upp-bd}
\begin{split}
 \|\xi\|_{C^3(\pt\O)}\leq C_{\xi_2}, \quad  |\nabla_x \xi|\geq C_{\xi_1} \;\; \text{ on }\; \pt\O.
\end{split}
\end{equation}
For given $(\bar{t},y,v)\in [0,T_0]\times\overline{\O}_\e\times \mathbb{R}^3$, define
the grazing set
\begin{equation}\label{def-Sy(v)}
\begin{split}
 S_y(v):=\big \{ v\in \mathbb{R}^3 : \  n\big (y_{\mathbf{b}}(y,v)\big) \cdot v=0 \big \},
\end{split}
\end{equation}
By Lemma 17 in \cite{guo2010decay}, the set $S_y(v)$ has zero Lebesgue measure.
\bigskip

\noindent
\subsection{$L^\infty$ Estimate for the Semigroup} \

\medskip

This subsection establishes the $L^\infty$ estimate for the semigroup generated by the linear homogeneous equation of \eqref{eq-unst-stretch-t-h-with-K} without collision $K$.
\medskip

We begin with an estimate for the backward bounce time.

\medskip
%-------------------Begin of Lemma------------------------
\begin{lemma}\label{k-bounce-time-model}  \
Let $(t_k,y_k,v_k)$ be  the $k$-th bounce of the backward trajectory \eqref{model-char}. Then
\begin{equation}\label{t2-t1-low-bd-model}
\begin{split}
& \tb(t_k,y_k,v_k) \geq \frac{C_{\xi_1}|v_k\cdot n(y_k)| } {\e C_{\xi_2} |v_k|^2}.
\end{split}
\end{equation}
\end{lemma}
%------------------end of Lemma------------------------

%----------------proof---------------------
\begin{proof}[\textbf{Proof}.] \
  By Taylor expansion of $\xi(\e y_{k+1})$ about $y_k$, we obtain
\begin{equation*}\label{model-xi-expan}
\begin{split}
  \xi(\e y_{k+1})=&\xi(\e y_k)+\e \nabla_x \xi(\e y_k)\cdot (y_{k+1}-y_k)\\
  &+
  \e^2 (y_{k+1}-y_k)\cdot\big[ \nabla_x^2 \xi(\tilde{\theta} \e  y_k+(1-\tilde{\theta})\e  y_{k+1})\big] \cdot(y_{k+1}-y_k), \;\; \tilde{\theta} \in (0,1).
\end{split}
\end{equation*}
Since $\xi(\e y_{k+1})=0=\xi(\e y_k)$ and $\nabla_x \xi\neq 0$, we have
\begin{equation}\label{y2-y1-normal-ineq-0}
\begin{split}
  \Big|\frac{\nabla_x \xi(\e y_k)}{|\nabla_x \xi(\e y_k)|}\cdot (y_{k+1}-y_k)\Big|
  =\e  \frac{\big| \nabla_x^2 \xi(\tilde{\theta} \e  y_k+(1-\tilde{\theta})\e  y_{k+1})\big| }{|\nabla_x \xi(\e y_k)|}|y_{k+1}-y_k|^2.
\end{split}
\end{equation}
 Using \eqref{xi-low-upp-bd} and \eqref{y2-y1-normal-ineq-0}, we obtain
\begin{equation}\label{y2-y1-normal-ineq}
\begin{split}
|n(y_k)\cdot (y_{k+1}-y_k)|
  \leq \e C_{\xi_2}C^{-1}_{\xi_1} |y_{k+1}-y_k|^2,
\end{split}
\end{equation}
where we have used the fact $ n(\e y_k)=n(y_k)$ derived from \eqref{ny=nx}.
%\begin{equation}
%\begin{split}
% n(\e y_k)=\frac{\nabla_x \xi\left(\e y_k\right)}{|\nabla_x \xi\left(\e y_k\right)|}
% =\frac{\nabla_y [\xi\left(\e y_k\right)] } {|\nabla_y [ \xi\left(\e y_k\right)]|}=n(y_k).
%\end{split}
%\end{equation}
Along the backward trajectory, we have
$y_{k+1}=y_k+ v_k (t_{k+1}-t_k)$,
which implies
\begin{equation}\label{y2-y1-v1-normal}
\begin{split}
  | y_{k+1}-y_k|=|t_{k+1}-t_k| \  |v_k|,\;\;\;\;  (y_{k+1}-y_k)\cdot n(y_k)=(t_{k+1}-t_k) [v_k\cdot n(y_k)].
\end{split}
\end{equation}
Substituting \eqref{y2-y1-v1-normal} into \eqref{y2-y1-normal-ineq} yields \eqref{t2-t1-low-bd-model}.
\end{proof}
\medskip
%-------------------end proof------------------------

The following lemma shows that for small $\e$, a backward specular trajectory in a non-grazing regime undergoes at most one bounce.

%-------------------Begin of Lemma------------------------
\begin{lemma}[Single-bounce for specular trajectory]\label{1-bounce-claim-model} \
Let $(\bar{t}, y, v)\in [0,T_0]\times \O_\e\times \big\{|v|\leq N, \ \big |v\cdot \frac{\nabla_x\xi(\e y)}  {|\nabla_x\xi(\e y)|}\big  |>\eta \big \}$ be given, with sufficient large constants $T_0, N>0$ and a small constant $\eta>0$. Define
\begin{equation}\label{def-epsilon-1}
\begin{split}
\e_1:=\frac{C_{\xi_1}^2\eta } { 2C_{\xi_2}^2 N^2T_0}\in(0,1).
\end{split}
\end{equation}
 If $0<\e\leq \e_1$,  then the backward specular trajectory \eqref{model-char} starting from $ (\bar{t},y,v)$ has at most one bounce.
  \end{lemma}
%------------------end of Lemma------------------------

\begin{proof}[\textbf{Proof.}] \
If $t_1\leq 0$, there is no bounce before reaching the initial plane $\{\bar{t}=0\}$. If $t_1>0$, it suffices to show that the backward time $t_{\bf b}(t_1,y_1,v_1)$ exceeds $T_0$ for sufficiently small $\e$.

Since $0<t_1<\bar{t}\leq T_0$ and $|v|\leq N$, we have
$|(t_1-\bar{t})v|\leq T_0N$. Because $y_1\in \pt\O_{\e}$, we have $\e y_1\in \pt\O$.
From the relation
\begin{equation}\label{y1-y}
\begin{split}
  y_1=y +  (t_1-\bar{t})v,
\end{split}
\end{equation}
 we see that $\e y \in \O$ lies close to the boundary $\partial\O$ for sufficiently small $\e$:
$$
 \e y=\e y_1 - \e (t_1-\bar{t})v= \e y_1 +O(\e)\sim \pt\O.
$$
 Indeed, for bounded velocity $|v|\leq N$, if the backward trajectory hits the boundary $\pt\O_\e$, the distance between the starting point $y$ and the boundary $\pt\O_\e$ must be bounded; consequently $\e y=x$ is near $\pt\O$.

Now observe that $n(y)= n(\e y)$, because $\nabla_x \xi(\e y)\neq 0$ near the boundary.
Expanding $\nabla_x\xi(\e y_1)$ about $y$ gives
%\begin{equation*}\label{ny'1}
%\begin{split}
%&\nabla_x\xi(\e y_1)=\nabla_x\xi(\e y)+ \e\nabla_x^2\xi\left(\bar{\theta} \e  %y_1+(1-\bar{\theta})\e  y\right)\cdot (y_1-y), \;\;\; \bar{\theta}\in (0,1),
%\end{split}
%\end{equation*}
\begin{equation}\label{ny1-ny}
\begin{split}
v\cdot \frac{\nabla_x\xi(\e y_1)}{|\nabla_x\xi(\e y_1)|}&=
v\cdot \frac{\nabla_x\xi(\e y)}{|\nabla_x\xi(\e y)|}
\frac{|\nabla_x\xi(\e y)|}{|\nabla_x\xi(\e y_1)|}+
\frac{v\cdot\e\nabla_x^2\xi\left(\bar{\theta} \e  y_1+(1-\bar{\theta})\e  y_2\right)\cdot (y_1-y)}{|\nabla_x\xi(\e y_1)|},
%&=n(\e y_1)+O(\e)
\end{split}
\end{equation}
where $\bar{\theta}\in (0,1)$.
%By (\ref{xi-low-upp-bd}) and (\ref{y1-y}), we have
%\begin{equation*}\label{ny'1}
%\begin{split}
%&C_{\xi_1} \leq |\nabla_x\xi(\e y)|\leq C_{\xi_2},\qquad C_{\xi_1} \leq |\nabla_x\xi(\e %y_1)|\leq C_{\xi_2},\\
%&\left|v\cdot\e\nabla_x^2\xi\left(\bar{\theta} \e  y_1+(1-\bar{\theta})\e  y_2\right)\cdot %(y_1-y)\right|\leq \e C_{\xi_2} T_0N^2.
%\end{split}
%\end{equation*}
Using \eqref{ny1-ny}, \eqref{y1-y} and \eqref{xi-low-upp-bd}, we obtain
\begin{equation}\label{ny1-lower}
\begin{split}
\big | v\cdot n(y_1)\big|=\Big |v\cdot \frac{\nabla_x\xi(\e y_1)}{|\nabla_x\xi(\e y_1)|}\Big |\geq \frac{C_{\xi_1}}{C_{\xi_2}}\Big| v\cdot \frac{\nabla_x\xi(\e  y)}{|\nabla_x\xi(\e  y)|}\Big|
-\e \frac{C_{\xi_2}}{C_{\xi_1}} T_0N^2.
\end{split}
\end{equation}
Substituting \eqref{ny1-lower} into \eqref{t2-t1-low-bd-model} and using $v_1=R_x v$ for specular reflection, we have
\begin{equation*}%\label{t2-t1-low-bd-ny}
\begin{split}
& \tb(t_1,y_1,v_1) \geq  \frac{C_{\xi_1}|v\cdot n(y_1)| } { \e C_{\xi_2} |v|^2} \geq  \Big[ \frac{ C_{\xi_1}^2\big | v\cdot  \frac{\nabla_x\xi(\e y)}  {|\nabla_x\xi(\e y)|} \big|} { \e C_{\xi_2}^2}- T_0N^2\Big]
 \frac{1} {|v|^2}
 \geq  \Big[ \frac{  C_{\xi_1}^2\eta} {\e C_{\xi_2}^2}- T_0N^2\Big]
 \frac{1} {N^2}\geq T_0,
\end{split}
\end{equation*}
provided $  0<\e\leq \e_1$. Hence, the backward trajectory reaches the initial plane $\{\bar{t}=0\}$ before any further bounce after $(t_1,y_1,v_1)$. The assertion is thus proved.
\end{proof}
%-------------End of proof--------------
\medskip

The following complementary result holds for a backward diffuse trajectory.

\medskip

%-------------------Begin of Lemma------------------------
\begin{lemma}[No further bounce for diffuse trajectory]\label{0-bounce-claim-model} \
 Let $(t_1, y_1, v^*_1)\in [0,T_0]\times \partial\O_\e\times \{|v^*_1|\leq N, \ | n(y_1) \cdot v^*_1 |>\eta \}$ be given, with sufficiently large constants $T_0, N>0$ and a small constant $\eta>0$. Define
\begin{equation}\label{def-epsilon-2}
\begin{split}
 \e_2:= \frac{C_{\xi_1} \eta } {C_{\xi_2} N^2 T_0} \in (0,1).
\end{split}
\end{equation}
If $0<\e\leq \e_2$, then the backward trajectory  \eqref{model-char}  starting from $(t_1,y_1,v^*_1)$ has no further collision.
\end{lemma}
%-------------------end of Lemma------------------------

%----------------proof---------------------

\begin{proof}[\textbf{Proof.}] \
 Following the proof of  \eqref{t2-t1-low-bd-model} in Lemma \ref{k-bounce-time-model},
 we obtain
\begin{equation*}\label{t'-t1-low}
\begin{split}
&   \tb(t_1, y_1, v^*_1)  \geq \frac{C_{\xi_1} | v^*_1 \cdot n(y_1)| } {\e C_{\xi_2} |v^*_1|^2}\geq \frac{C_{\xi_1} \eta } {\e C_{\xi_2} N^2}\geq T_0,
\end{split}
\end{equation*}
provided $0<\e \leq \e_2$. Thus no further collision occurs after leaving $(t_1, y_1, v^*_1)$.
\end{proof}
\medskip

Finally, we state the semigroup estimate for the linear homogeneous Boltzmann equation  without collision $K$ under the Maxwell boundary condition.
\medskip

\begin{lemma}[Semigroup estimate]\label{semigroup-estimate}  \
Let  $h_0\in L^\infty(\O_\e\times\mathbb{R}^3)$, and let $\e_2$ be the constant defined in \eqref{def-epsilon-2}. Then, for every $0<\e\leq \e_2$, the weighted linear problem
%This equation does not contain K and \bar{g}, different with %
\begin{equation}\label{eq-unst-stretch-t-h}
\begin{split}
 \pt_{\bar{t}}h + v\cdot\nabla_y h      +   \nu(v)h=0
 \ \  &\text{ in } \mathbb{R}^+\times \O_{\e}\times\mathbb{R}^3, \\
    h \big|_{\g_-} = (1-\a) \mathscr{R}h + \a \o \sqrt{\mu} \int_{n(y)\cdot u>0}h \dd \sigma \ \   &\text{ on }\mathbb{R}^+\times \pt\O_{\e}\times\mathbb{R}^3,\\
 h |_{t=0} = h_0    \ \   &\text{ on } \O_{\e}\times\mathbb{R}^3
\end{split}
\end{equation}
admits a unique solution $h(\bar{t},y,v) =\left\{G(\bar{t})h_0\right\}(y,v)$ satisfying
 \begin{equation}\label{semigroup-esti}
\begin{split}
&\left\|G(\bar{t}) h_0 \right\|_{L^\infty_{\bar{t},y,v}(\mathbb{R}^+\times\O_{\e}\times \mathbb{R}^3)}
\leq (2C_*+1)e^{-\frac{\nu_0}{2}\bar{t}} \left\| h_0 \right\|_{L^\infty_{y,v}(\O_{\e}\times \mathbb{R}^3)} \quad \text{ for all }\bar{t}>0.
\end{split}
\end{equation}
\end{lemma}
\medskip
%------------------end of Lemma------------------------

\begin{proof}[\textbf{Proof.}]  \ The proof is divided into two steps. In Step 1, we derive the uniform estimate on a bounded time interval. In Step 2, we extend the result to the entire  $\mathbb{R}^+$.
\medskip

\noindent \textbf{Step 1. Uniform estimate on a bounded time interval.}

We claim that for any sufficiently large $T_0>0$ satisfying $(2C_*+1)e^{-\frac{\nu_0}{2}T_0}\leq 1$,  the following estimate holds:
 \begin{equation}\label{semigroup-esti-local}
\begin{split}
&\sup_{0\leq s\leq T_0}\big[ e^{\nu_0 s} \|h(s) \|_{_{L^\infty_{y,v}}} \big]
\leq (2C_*+1)\| h_0 \|_{L^\infty_{y,v}}.
\end{split}
\end{equation}

 To prove this, we construct an iterative sequence $\{h^{n+1}\}_{n=0}^\infty$ via
\begin{equation}\label{model-eq-stretch-t}
\begin{split}
& \pt_{\bar{t}} h^{n+1} + v\cdot\nabla_y  h^{n+1}  +   \nu(v)  h^{n+1}=0\ \  \text{ in } \mathbb{R}^+\times\O_\e\times\mathbb{R}^3,\\
&  h^{n+1} \big|_{\g_-}  = (1-\a) \mathcal{L}h^{n} + \a \o \sqrt{\mu} \int_{n(y)\cdot u>0} h^{n+1} \dd \sigma \ \   \text{ on }  \mathbb{R}^+\times \pt\O_\e\times\mathbb{R}^3,\\
&  h^{n+1}=   h_0 \ \   \text{ on } \O_\e\times\mathbb{R}^3,
\end{split}
\end{equation}
 with the initial iterate
\begin{equation}\label{h^{0}-def}
   h^{0}= h^{0}(\bar{t},y,v):= e^{-\nu_0 \bar{t}} h_0(y,v).
\end{equation}
To establish \eqref{semigroup-esti-local}, it suffices to show that
 \begin{equation}\label{semigroup-esti-local-h^{n+1}}
\begin{split}
&\sup_{0\leq s\leq T_0}\big[ e^{\nu_0 s} \|h^{n+1}(s) \|_{L^\infty_{y,v}}\big]
\leq (2C_*+1)\| h_0 \|_{L^\infty_{y,v}}\;\;  \text{ for all } n=0,1,2,\cdots.
\end{split}
\end{equation}
Indeed, once the uniform estimate \eqref{semigroup-esti-local-h^{n+1}} is verified, there exists a function $h\in L^\infty([0,T]\times \O_\e\times\mathbb{R}^n)$ such that
a subsequence of $\{h^{n+1}\}$ (still denoted by $\{ h^{n+1}\}$) satisfies
$$
h^{n+1}\rightarrow h \;\; \text{ weakly}\!-\!* \text{ in } L^\infty([0,T]\times \O_\e\times\mathbb{R}^n) \; \text{ as } n\rightarrow \infty,
$$
and the limit $h$ satisfies the uniform estimate \eqref{semigroup-esti-local}
and the linear problem \eqref{eq-unst-stretch-t-h} in the weak sense.

We now verify the uniform estimate \eqref{semigroup-esti-local-h^{n+1}} in four sub-steps.
\medskip

\noindent \textbf{Step 1.1. The first bounce.}

For $\e\in (0,1]$, $\a\in [0,1]$, $n\in\mathbb{N}$, $\bar{t}\in [0,T_0]$ and $(y,v)\in \overline{\O}_\e\times\mathbb{R}^3\setminus \gamma_0$ with $v\notin S_y(v)$, using the characteristic trajectory \eqref{model-char} and the equation $\eqref{model-eq-stretch-t}_1$, we obtain
\begin{equation}\label{model-h-formula}
\begin{split}
&\frac{d}{ds } \Big[ e^{-\int^{\bar{t}}_s\nu(v)\dd \tau }  h^{n+1} \big(s, Y(s;\bar{t},y,v), v \big)  \Big]=0
\end{split}
\end{equation}
 for $t_1< s \leq \bar{t}$. Integrating along the backward trajectory yields
\begin{equation}\label{J1}
\begin{split}
  h^{n+1}(\bar{t},y,v) = &\mathbf{1}_{\{t_1\leq 0\}}   e^{-\int^{\bar{t}}_0\nu(v)\dd \tau }
   h_0  \big( Y(0; \bar{t},y,v ),v \big) \\
 &+\mathbf{1}_{\{t_1>0\}}     e^{-\int^{\bar{t}}_{t_1}\nu(v)\dd \tau }  (1-\a) h^{n}(t_1,y_1,\,v_1) \\
%-----------------%---------------
 &+  \mathbf{1}_{\{t_1>0\}}    e^{-\int^{\bar{t}}_{t_1}\nu(v)\dd \tau } \a \o \sqrt{\mu}
 \int_{n(y_1)\cdot v_1^*>0}   h^{n+1}(t_1,y_1,v_1^*)
  \dd\sigma^{*}_1\\
   := & J_0^1+J_{sp}^1+J_{di}^1,
\end{split}\end{equation}
where $\dd \sigma^{*}_1=\o ^{-1}\sqrt{2\pi}\mu^{\frac{1}{2}} [n(y_1)\cdot v_1^*] \dd v_1^*$ similarly as in  \eqref{def-C-star}.
Obviously, $J_0^1(\bar{t},y,v)$ is bounded by
\begin{equation} \label{J1-0-bd}
\begin{split}
& \left| J_0^1(\bar{t},y,v) \right| \leq \mathbf{1}_{\{t_1\leq 0\}} e^{-\nu_0\bar{t}}\|h_0\|_{L^\infty_{y,v}}.
\end{split}
\end{equation}

For the diffuse boundary term $J_{di}^1$, we partition the integration domain:
\begin{align}
 &\{n(y_1) \cdot v^{*}_{1}>0  \}\nonumber\\
 =&\{|v^{*}_{1}|>N,\ n(y_1) \cdot v^{*}_{1}>0 \}\cup \{|v^{*}_{1}|\leq N,\  0<  n(y_1) \cdot v^{*}_{1}  <\eta \}
 \cup \{|v^{*}_{1}|\leq N,\ n(y_1) \cdot v^{*}_{1}  \geq \eta\} \nonumber \\
  :=& A^*_1(v^{*}_{1})\cup A^*_2(v^{*}_{1})\cup M^*_{y_1}(v^{*}_{1}),\label{domain-split-diff}
\end{align}
%\begin{equation}\label{domain-split-diff}
%\begin{split}
%A^*_1(v^{*}_{1}) & := \left\{|v^{*}_{1}|>N,\ n(y_1) \cdot v^{*}_{1}>0 \right \},\\
%A^*_2(v^{*}_{1}) & :=\left\{|v^{*}_{1}|\leq N,\  0<  n(y_1) \cdot v^{*}_{1}  <\eta %\right\},\\
%M^*_{y_1}(v^{*}_{1}) & :=\left\{|v^{*}_{1}|\leq N,\ n(y_1) \cdot v^{*}_{1}  \geq %\eta\right\},
%\end{split}
%\end{equation}
with positive constants $N$ and $\eta$ to be determined later.
On $A^*_1(v^{*}_{1})$,  $J_{di}^1$ is bounded as
\begin{equation}\label{J{di}1-A1}
\begin{split}
 \big| J_{di}^1{\bf 1}_{A^*_1(v^{*}_{1})} \big|
%\leq \; &  \a e^{-\nu_0(\bar{t}-t_1)} \int_{|v^{*}_{1}|\geq N,\ n(y_1) \cdot v^{*}_{1} >0} %\big |h^{n+1}(t_1,y_1,v_1^*)\big|  \dd\sigma^{*}_1\\
 \leq \; & \a  o(1) e^{-\nu_0(\bar{t}-t_1)}|h^{n+1}(t_1)|_{L^\infty_{y,v}(\pt\O_\e\times\mathbb{R}^3)},
\end{split}
\end{equation}
provided $N>0$ is sufficiently large.  For $A^*_2(v^{*}_{1})$ and fixed $N$, we apply the decomposition $v^{*}_{1,\perp}
=v^{*}_{1}- v^{*}_{1,\parallel}$ with $v^{*}_{1,\parallel}=[v^{*}_{1}\cdot n(y_1)] n(y_1)$ for $|v^{*}_{1}\cdot n(y_1)|<\eta$ to obtain
\begin{equation}\label{cut-inpd-e}
\begin{split}
 \big|J_{di}^1{\bf 1}_{A^*_2(v^{*}_{1})}  \big|
% &\leq \a e^{-\nu_0(\bar{t}-t_1)}
% \int_{|v^{*}_{1}|\leq N,\ | n(y_1) \cdot v^{*}_{1} | <\eta} \big %|h^{n+1}(t_1,y_1,v_1^*)\big| \widetilde{w}(v_1^*)  \dd\sigma^{*}_1 \\
%---------------
 &\leq \a e^{-\nu_0(\bar{t}-t_1)} C_N |h^{n+1}(t_1)|_{L^\infty_{y,v}(\pt\O_\e\times\mathbb{R}^3)} \int^{\eta}_{-\eta}\dd v^{*}_{1,\parallel} \int_{|v^{*}_{1,\perp}|\leq N}\dd v^{*}_{1,\perp}\\
%---------------
%&\leq \a e^{-\nu_0(\bar{t}-t_1)}\eta C_N %\|h^{n+1}(t_1)\|_{L^\infty_{y,v}(\pt\O_\e\times\mathbb{R}^3)}\\
%---------------
&\leq \a o(1) e^{-\nu_0\bar{t}}\sup_{0\leq s\leq T_0}\big[ e^{\nu_0 s} |h^{n+1}(s) |_{L^\infty_{y,v}(\pt\O_\e\times\mathbb{R}^3)} \big],
\end{split}
\end{equation}
provide $\eta>0$ is  sufficiently small.
For the bulk $M^*_{y_1}(v^{*}_{1})$, Lemma \ref{0-bounce-claim-model} implies that
for $0<\e\leq \e_2$,  the backward trajectory starting from $(t_1,y_1,v_1^*)$ undergoes
no further collisions. Thus, $J_{di}^1(\bar{t},y,v){\bf 1}_{M^*_{y_1}(v^{*}_{1})}$ traces  back  to the initial plane $\{\bar{t}=0\}$
%\begin{equation*}\label{J1-di-0}
%\begin{split}
% & J_{di}^1(\bar{t},y,v){\bf 1}_{M^*_{y_1}(v^{*}_{1})}\\
% =\;&\mathbf{1}_{\{t_1>0\}} \a  e^{-\int^{\bar{t}}_{t_1}\nu(v)\dd \tau }
%  w\sqrt{\mu}\int_{ {|v_1^*|\leq N},{|n(y_1)\cdot v_1^*| \geq\eta }}
%  e^{-\int^{t_1}_{0}\nu(v_1^*)\dd \tau}
%    h_0 \big(Y(0; t_1,y_1,v_1^*), v_1^*\big)\dd\sigma_1^*,
%\end{split}
%\end{equation*}
and is bounded as:
\begin{equation} \label{J1-di-0-bd}
\begin{split}
& \big| J_{di}^1{\bf 1}_{M^*_{y_1}(v^{*}_{1})}\big|
 \leq \mathbf{1}_{\{t_1>0\}} \a e^{-\nu_0\bar{t}}
 \Big| \int_{n(y)\cdot u>0}  \dd \sigma \Big|
 \|h_0\|_{L^\infty_{y,v}}\leq \mathbf{1}_{\{t_1>0\}} \a C_*  e^{-\nu_0\bar{t}}
  \|h_0\|_{L^\infty_{y,v}}.
\end{split}
\end{equation}
%Combining the above estimates yields
%\begin{equation}\label{J1-di-bd}
%\begin{split}
%  \left| J_{di}^1(\bar{t},y,v) \right|
%  \leq & \mathbf{1}_{\{t_1>0\}}\a  C_* e^{-\nu_0\bar{t}} \|h_0\|_{L^\infty_{y,v}} \\
% &+   \mathbf{1}_{\{t_1>0\}} \a o(1) e^{-\nu_0\bar{t}} \sup_{0\leq s\leq T_0}\left[ %e^{\nu_0 s} \left\|h^{n+1}(s) \right\|_{L^\infty_{y,v}(\pt\O_\e\times\mathbb{R}^3)} %\right].
%%\end{split}
%\end{equation}
Note that the $o(1)$ term depends only on $N>0$ and $\eta>0$, and is independent of $\e$ and $\a$.

Combining estimates \eqref{J1}, \eqref{J1-0-bd} and \eqref{J{di}1-A1}--\eqref{J1-di-0-bd}, we obtain
\begin{equation}\label{J1-bd}
\begin{split}
 |h^{n+1}(\bar{t},y,v)| \leq & \big(\mathbf{1}_{\{t_1\leq 0\}}+
 \mathbf{1}_{\{t_1>0\}}\a{C_*} \big) e^{-\nu_0\bar{t}}\|h_0\|_{L^\infty_{y,v}}\\
 &+  \mathbf{1}_{\{t_1>0\}} \a o(1)  e^{-\nu_0\bar{t}}
 \sup_{0\leq s\leq T_0}\big[ e^{\nu_0 s} |h^{n+1}(s) |_{L^\infty_{y,v}(\pt\O_\e\times\mathbb{R}^3)} \big]\\
 &+ \mathbf{1}_{\{t_1>0\}}   (1-\a)   e^{-\nu_0(\bar{t}-t_1)} \big| h^{n}(t_1,y_1,\,v_1)\big|.
\end{split}
\end{equation}

\noindent\textbf{Step 1.2. The 2nd bounce.}\;

After the first collision at $(t_1,y_1,\,v_1)$, the term $J_{sp}^1$ may continue to undergo reflection along the specular backward trajectory.
Note that the equation of $h^{n}$ shares the same specular backward trajectory as that of $h^{n+1}$. Consequently, we have
\begin{equation}%\label{J2}
\begin{split}
   h^{n}(t_1,y_1,v_1) = &
 \mathbf{1}_{\{t_2\leq 0<t_1\}}  e^{-\int^{t_1}_{0}\nu(v)\dd \tau }
 h_0\left(Y(0;t_1,y_1,v_1),v_1\right) \\
%---------------------------------
& +  \mathbf{1}_{\{t_2>0\}}  (1-\a)   e^{-\int^{t_1}_{t_2}\nu(v)\dd \tau }
 h^{n-1}\left(t_2,y_2,v_2\right)  \\
%-----------------%---------------
 & +  \mathbf{1}_{\{t_2>0\}} \a   e^{-\int^{t_1}_{t_2}\nu(v)\dd \tau }\o\sqrt{\mu}
 \int_{n(y_2)\cdot v_2^*>0}   h^{n}(t_2,y_2,v_2^*) \dd\sigma^{*}_2\\
 := &   J_{0}^2+J_{sp}^2+J_{di}^2.
\end{split}\end{equation}
Similarly to \eqref{J1-0-bd}, $J_0^2$ is bounded by $\mathbf{1}_{\{t_2\leq 0<t_1\}}e^{-\nu_0t_1}\|h_0\|_{L^\infty_{y,v}}$.
Following the same procedure as that of $J_{di}^1$, we partition the integration domain $\{n(y_2) \cdot v^{*}_{2}>0 \}$ and bound $J_{di}^2$ as:
\begin{equation}\label{J2-di-bd}
\begin{split}
 |J_{di}^2|  \leq  &  \mathbf{1}_{\{t_2>0\}} \a {C_*} e^{-\nu_0\bar{t}}\|h_0\|_{L^\infty_{y,v}}
  +  \mathbf{1}_{\{t_2>0\}} \a  o(1) e^{-\nu_0t_1} \sup_{0\leq s\leq T_0}\big[ e^{\nu_0 s} |h^n(s)|_{L^\infty_{y,v}(\pt\O_\e\times\mathbb{R}^3)} \big].
\end{split}
\end{equation}
Thus, $ \big| h^{n}(t_1,y_1,\,v_1)\big|$ satisfies the bound:
\begin{equation}\label{J_{sp}^1-bd}
\begin{split}
  \big| h^{n}(t_1,y_1,\,v_1)\big| \leq\;
  & \big(\mathbf{1}_{\{t_2\leq 0<t_1\}}
 + \mathbf{1}_{\{t_2>0\}} \a {C_*} \big)e^{-\nu_0t_1}\|h_0\|_{L^\infty_{y,v}} \\
 & +  \mathbf{1}_{\{t_2>0\}}  \a  o(1)  e^{-\nu_0t_1}
 \sup_{0\leq s\leq T_0} \big[ e^{\nu_0 s} |h^n(s) |_{L^\infty_{y,v}(\pt\O_\e\times\mathbb{R}^3)} \big]\\
 & + \mathbf{1}_{\{t_2>0\}}  (1-\a)   e^{-\nu_0(t_1-t_2)}
 |  h^{n-1}(t_2,y_2,v_2) |.
\end{split}
\end{equation}

\noindent\textbf{Step 1.3. The $k$-th bounce.}\;

Proceeding inductively, after the $(k-1)$-th collision,  the term $J_{sp}^{k-1}$
 may continue to undergo reflections along the specular backward trajectory, leading to
the $k$-th collision:
\begin{equation}%\label{J2}
\begin{split}
  & h^{n+1-(k-1)}\left(t_{k-1},y_{k-1},v_{k-1}\right)\\
= &  \mathbf{1}_{\{t_{k}\leq 0<t_{k-1}\}}    e^{-\int^{t_{k-1}}_{0}\nu(v_{k-1})\dd \tau }  h_0\left(Y(0;t_{k-1},y_{k-1},\,v_{k-1}),v_{k-1}\right)  \\
%---------------------------------
 & + \mathbf{1}_{\{t_k> 0\}} (1-\a)   e^{-\int^{t_{k-1}}_{t_k}\nu(v_{k-1})\dd \tau }
 h^{n+1-k}\left(t_k,y_k,v_k\right)  \\
%-----------------%---------------
& +\mathbf{1}_{\{t_k>0\}}     \a  e^{-\int^{t_{k-1}}_{t_k}\nu(v_{k-1})\dd \tau }
 \int_{n(y_k)\cdot v_k^*>0}  h^{n+1-(k-1)}(t_k,y_k,v_k^*) \dd\sigma^{*}_k
  \\  :=  & J_{0}^k+J_{sp}^k
  +J_{di}^k.
\end{split}\end{equation}
Analogous to the derivation of  \eqref{J_{sp}^1-bd}, we obtain the bound:
\begin{equation}\label{J_{sp}^{n+1-(k-1)}-bd}
 \begin{split}
  &| h^{n+1-(k-1)}(t_{k-1},y_{k-1},v_{k-1})|\\
  \leq &
 \big(\mathbf{1}_{\{t_{k}\leq 0<t_{k-1}\}} + \mathbf{1}_{\{t_k>0\}} \a  {C_*}
  \big)e^{-\nu_0 t_{k-1} }\|h_0\|_{L^\infty_{y,v}}\\
 &+ \mathbf{1}_{\{t_k>0\}} \a o(1)e^{-\nu_0t_{k-1}}
 \sup_{0\leq s\leq T_0} \big[ e^{\nu_0 s}| h^{n+1-(k-1)}(s) |_{L^\infty_{y,v}(\pt\O_\e\times\mathbb{R}^3)}  \big]\\
 &+\mathbf{1}_{\{t_k> 0\}} (1-\a)   e^{-\nu_0(t_{k-1}-{t_k}) }
|  h^{n+1-k}(t_k,y_k,v_k)|.
\end{split}
\end{equation}

\noindent\textbf{ Step 1.4. Bounce back trajectory starting from $(t_k,y_k,v_k)$.}\;

 After the $k$-th collision at  $(t_{k},y_{k},\,v_{k})$, the term $J_{sp}^{k}(t_{k},y_{k},\,v_{k})$
 may continue to propagate along the specular backward trajectory:
\begin{equation}\label{h^{n+1-k}}
\begin{split}
  &h^{n+1-k}\left(t_k,y_k,v_k\right)\\
= &\mathbf{1}_{\{t_{k+1}\leq 0<t_k\}}  e^{-\int^{t_k}_{0}\nu(v_k)\dd \tau }
 h_0\left(Y(0;t_k,y_k,\,v_k),v_k\right),  \\
%---------------------------------
& +\mathbf{1}_{\{t_{k+1}> 0\}} (1-\a)   e^{-\int^{t_{k}}_{t_{k+1}}\nu(v_k)\dd \tau }
 h^{n-k}\left(t_{k+1},y_{k+1},v_{k+1}\right),  \\
%-----------------%---------------
 & + \mathbf{1}_{\{t_{k+1}>0\}}     \a  e^{-\int^{t_{k}}_{t_{k+1}}\nu(v_k)\dd \tau }
 \int_{n(y_{k+1})\cdot v_{k+1}^*>0}
h^{n+1-k}(t_{k+1},y_{k+1},v_{k+1}^*) \dd\sigma^{*}_{k+1},\\
 := & J_{0}^{k+1}+J_{sp}^{k+1}+J_{di}^{k+1}.
\end{split}\end{equation}

Clearly, $k\leq n$, since the term $h^{n-k}(t_{k},y_{k},v_{k})$ on the right hand side of the expression for $ J_{sp}^{k+1}$ in \eqref{h^{n+1-k}} generates the initial iterate $h^0$ when $k=n$, and no further collision occur for given initial iteration $h^0$. Recall that $t_0=\bar{t}$. For any fixed $n\in\mathbb{N}$, there are two possible cases:
$(1)$ There exists some $k\in \{0,1,2, \cdots, n\}$ such that $t_{k+1}\leq 0< t_k$;
$(2)$ $t_{k+1}>0$ for all $k\in \{0,1,2, \cdots, n\}$.
We now estimate $h^{n+1}(\bar{t},y,v)$ according to these two cases.

\noindent\textbf{Case 1: There exists $k\in \{0,1,2, \cdots, n\}$ such that $t_{k+1}\leq 0< t_k$.}\;

 In this case, for such a $k\in \{0,1,2, \cdots, n\}$, we have
\begin{equation}\label{t-i-index}
\begin{split}
\mathbf{1}_{\{ t_{k+1}\leq 0< t_{k}\}}=1, \;\; \mathbf{1}_{\{ t_{k+1}>0\}}=0;\;\; \mathbf{1}_{\{ t_{i}\leq 0 < t_{i-1}\}}=0,  \;\;  \mathbf{1}_{\{ t_{i}>0\}}=1,\;\;
\forall i\in \{1,2, \cdots,k\}.
\end{split}
\end{equation}
This means that the backward trajectory starting from  $(t_k,y_k,v_k)$
reaches the initial plane $\{\bar{t}=0\}$ with no further collision. Therefore,
\begin{equation*}%\label{h^{n+1-k}-initial}
 \begin{split}
  & | h^{n+1-k}(t_k,y_k,v_k)| =
| J_0^{k+1}| \leq
 \mathbf{1}_{\{t_{k+1}\leq 0<t_k\}}  e^{-\nu_0 t_{k} }\|h_0\|_{L^\infty_{y,v}},
\end{split}
\end{equation*}
and hence
\begin{equation}\label{h^{n+1-k}-initial}
\begin{split}
   &\sup_{0\leq s \leq T_0}\big[   e^{\nu_0s } | h^{n+1-k}(s ) |_{L^\infty_{y,v}(\partial\O\times\mathbb{R}^3)} \big]
    \leq    \|h_0\|_{L^\infty_{y,v}}.
\end{split}
\end{equation}
Substituting \eqref{h^{n+1-k}-initial} into the right-hand side of  \eqref{J_{sp}^{n+1-(k-1)}-bd} and  using \eqref{t-i-index}, we obtain
\begin{equation}\label{h^{n+1-(k-1)}}
\begin{split}
   &[ 1- \a o(1)]\sup_{0\leq s \leq T_0}\big[   e^{\nu_0 s } |h^{n+1-(k-1)}(s ) |_{L^\infty_{y,v}(\pt\O\times{R}^3)} \big]
    \leq
    \a{C_*}\|h_0\|_{L^\infty_{y,v}}  +      (1-\a) \|h_0\|_{L^\infty_{y,v}}.
\end{split}
\end{equation}

Next, substituting \eqref{h^{n+1-(k-1)}} into the estimate of $h^{n+1-(k-2)}(t_{k-2})$ and deducing similarly,
\begin{equation}\label{h^{n+1-(k-2)}-bd-factor}
\begin{split}
   &[ 1- \a o(1)]^2\sup_{0\leq s \leq T_0}\big[   e^{\nu_0s }
|h^{n+1-(k-2)}(s )|_{L^\infty_{y,v}(\pt\O_\e\times{R}^3)} \big] \\
    \leq
    & \{ \a{C_*} [ 1- \a o(1)]
    +\a{C_*}  (1-\a)
  +   (1-\a)^2  \} \|h_0\|_{L^\infty_{y,v}}.
\end{split}
\end{equation}
Repeating this process for $h^{n+1-(k-3)}(t_{k-3} )$, we obtain
\begin{equation}\label{h^{n+1-(k-2)}-bd}
\begin{split}
   &[ 1- \a o(1)]^3\sup_{0\leq s \leq T_0}\big[   e^{\nu_0 s } |h^{n+1-(k-3)}(s ) |_{L^\infty_{y,v}(\pt\O_\e\times{R}^3)} \big]  \\
    \leq
    & \{ \a{C_*}  [ 1- \a o(1)]^2
        + \a{C_*} [ 1- \a o(1)]   ( 1- \a)
    + \a{C_*} ( 1- \a)^2  +(1-\a)^3   \} \|h_0\|_{L^\infty_{y,v}}.
\end{split}
\end{equation}
By induction and  \eqref{t-i-index}, we arrive at
\begin{equation}%\label{h^{n+1}-bd}
\begin{split}
   &[ 1- \a o(1)]^k\sup_{0\leq s\leq T_0}\big[   e^{\nu_0 s } |h^{n+1}(s ) |_{L^\infty_{y,v}(\O_\e\times{R}^3)} \big]
   \\
%    \leq
%    & \; \a{C_*} \left[ 1- \a o(1)\right]^{k-1} \|h_0\|_{L^\infty_{y,v}}\\
%----------------------
%    &+  \a{C_*}  \left[ 1- \a o(1)\right]^{k-2}( 1- \a) \|h_0\|_{L^\infty_{y,v}}\\
%----------------------
%    &+\a{C_*} \left[ 1- \a o(1)\right]^{k-3}( 1- \a)^2 \|h_0\|_{L^\infty_{y,v}}\\
%----------------------
%    &+\cdots\\
%      & + \a{C_*} \left[ 1- \a o(1)\right]^2 ( 1- \a)^{k-3}\|h_0\|_{L^\infty_{y,v}}\\
%----------------------
%       & + \a{C_*}  \left[ 1- \a o(1)\right] ( 1- \a)^{k-2} \|h_0\|_{L^\infty_{y,v}}\\
%----------------------
%   & + \a{C_*}  ( 1- \a)^{k-1}\|h_0\|_{L^\infty_{y,v}}\\
%----------------------
% &
%  +    (1-\a)^k  \|h_0\|_{L^\infty_{y,v}}\\
   \leq  \,
    & \a {C_*} \sum_{i=1}^{k} \left[ 1- \a o(1)\right]^{k-i}( 1- \a)^{i-1}\|h_0\|_{L^\infty_{y,v}}
    +   (1-\a)^k \|h_0\|_{L^\infty_{y,v}}.
 %   = \,
%    &   {C_*} \left[ 1- \a o(1)\right]^{k} \left[ 1- \left(\frac{1-\a}{1-\a %o(1)}\right)^{k}\right] \frac{1}{1-o(1)} \|h_0\|_{L^\infty_{y,v}}
%    +     (1-\a)^k \|h_0\|_{L^\infty_{y,v}},
\end{split}
\end{equation}

Finally, we obtain the following uniform bound for $h^{n+1}(\bar{t})$:
\begin{equation}\label{h^{n+1}-bd-uniform-k}
\begin{split}
   &\sup_{0\leq s\leq T_0}\big[   e^{\nu_0s } \|h^{n+1}(s ) \|_{L^\infty_{y,v}(\O_\e\times{R}^3)} \big]
   \\
    \leq \,
    &  {C_*} \Big[ 1-  \Big(\frac{1-\a}{1-\a o(1)} \Big)^{k} \Big] \frac{1}{1-o(1)} \|h_0\|_{L^\infty_{y,v}}
    +   \Big(\frac{1-\a}{1-\a o(1)} \Big)^{k}    \|h_0\|_{L^\infty_{y,v}}
    \\
    \leq \,
    & (2C_*+1)   \|h_0\|_{L^\infty_{y,v}},
\end{split}
\end{equation}
where the last inequality follows from the bounds
\begin{equation}\label{o(1)-bd}
o(1)\leq \frac{1}{2},\,\,\,\, 1-o(1)\geq \frac{1}{2},\,\,\,\,  1-\a\leq 1-\a o(1).
\end{equation}

\noindent\textbf{Case 2:  $t_{k+1}>0$ for all $k\in \{0,1,2, \cdots, n\}$.}

In this case, after the $n$-th collision at $(t_{n},y_{n},\,v_{n})$, the specular trajectory continues to  propagate and produce an $(n+1)$-th collision. Taking $k=n$ in \eqref{h^{n+1-k}}, we obtain
\begin{equation}%\label{J2}
 \begin{split}
  h^{1}(t_{n},y_{n},v_{n})
 =& \mathbf{1}_{\{t_{n+1}\leq 0<t_{n}\}}   e^{-\int^{t_{n}}_{0}\nu(v_{n})\dd \tau }
  h_0\left(Y(0;t_{n},y_{n},\,v_{n}),v_{n}\right)  \\
%---------------------------------
 &+  \mathbf{1}_{\{t_{n+1}> 0\}} (1-\a)   e^{-\int^{t_{n}}_{t_{n+1}}\nu(v_{n})\dd \tau }
  h^{0}(t_{n+1},y_{n+1},v_{n+1})  \\
%-----------------%---------------
& +  \mathbf{1}_{\{t_{n+1}>0\}}     \a  e^{-\int^{t_{n}}_{t_{n+1}}\nu(v_{n})\dd \tau }
 \int_{n(y_{n+1})\cdot v_{n+1}^*>0}  h^1(t_{n+1},y_{n+1},v_{n+1}^*) \dd\sigma^{*}_{n+1}  \\
 :=& J_{0}^{n+1}+J_{sp}^{n+1} +J_{di}^{n+1}.
\end{split}\end{equation}
Following a similar procedure as in \eqref{J_{sp}^1-bd}, we bound $\left| h^{1}\left(t_{n},y_{n},v_{n}\right)\right|$ as
\begin{equation}\label{J_{sp}^{h1}-bd}
 \begin{split}
 | h^{1}(t_{n},y_{n},v_{n})| \leq \; &
 \left(\mathbf{1}_{\{t_{n+1}\leq 0<t_{n}\}} + \mathbf{1}_{\{t_{n+1}>0\}} \a  {C_*}
 \right)e^{-\nu_0 t_{n} }\|h_0\|_{L^\infty_{y,v}}\\
 &+ \mathbf{1}_{\{t_{n+1}>0\}} \a o(1)e^{-\nu_0t_{n}}
 \sup_{0\leq s\leq T_0}\big[ e^{\nu_0 s} | h^{1}(s) |_{L^\infty_{y,v}(\pt\O\times{R}^3)} \big] \\
 &  + \mathbf{1}_{\{t_{n+1}>0\}} (1-\a)   e^{-\nu_0  t_{n}} \| h_0\|_{L^\infty_{y,v}},
\end{split}
\end{equation}
where the last term has used the initial iterate $h^0\equiv e^{-\nu_0\bar{t}}h_0$ and the bound
$$
  e^{\nu_0t_{n+1}}| h^{0}\left(t_{n+1},y_{n+1},v_{n+1}\right)|
  =| h_{0}\left(y_{n+1},v_{n+1}\right)|
  \leq  \| h_0\|_{L^\infty_{y,v}}.
$$
 Since $t_{n+1}>0$, we have $\mathbf{1}_{\{ t_{i}\leq 0 < t_{i-1}\}}=0$ for all $i\in \{1,2, \cdots,n,n+1\}$. Then, \eqref{J_{sp}^{h1}-bd} implies
\begin{equation}\label{h^{1}(t_{n}-bd}
\begin{split}
   &[ 1- \a o(1)]\sup_{0\leq s \leq T_0}\big[   e^{\nu_0 s }|h^{1}( s ) |_{L^\infty_{y,v}(\pt\O\times{R}^3)} \big]
    \leq
    \a{C_*} \|h_0\|_{L^\infty_{y,v}}+   (1-\a)  \| h_0\|_{L^\infty_{y,v}}.
\end{split}
\end{equation}

Substituting \eqref{h^{1}(t_{n}-bd} into  the estimate for $h^{2}(t_{n-1} )$,  we derive
\begin{equation*} %\label{h^{n+1-(k-1)}-bd-factor}
\begin{split}
   &[ 1- \a o(1)]^2\sup_{0\leq s \leq T_0}\big[   e^{\nu_0 s } |h^{2}(s ) |_{L^\infty_{y,v}(\pt\O\times{R}^3)} \big] \\
    \leq &
     \{ \a{C_*}  [1-\a o(1)] + \a{C_*} (1-\a) +  (1-\a)^2  \} \|h_0\|_{L^\infty_{y,v}}.
\end{split}
\end{equation*}
 Proceeding iteratively as in  case 1, we finally obtain
\begin{equation}%\label{h^{n+1}-bd}
\begin{split}
   &[ 1- \a o(1)]^{n+1}\sup_{0\leq s \leq T_0}
\big[   e^{\nu_0 s } \|h^{n+1}( s ) \|_{L^\infty_{y,v}(\O_
\e\times{R}^3) } \big]\\
%   \\
%    \leq
%    & \; \a{C_*} \left[ 1- \a o(1)\right]^{n} \|h_0\|_{L^\infty_{y,v}}\\
%----------------------
%    &+ \a{C_*}  \left[ 1- \a o(1)\right]^{n-1}( 1- \a) \|h_0\|_{L^\infty_{y,v}}\\
%----------------------
%    &+ \a{C_*}  \left[ 1- \a o(1)\right]^{n-2}( 1- \a)^2 \|h_0\|_{L^\infty_{y,v}}\\
%----------------------
%    &+\cdots\\
% %     & +\a{C_*} \left[ 1- \a o(1)\right]^2 ( 1- \a)^{n-2}\|h_0\|_{L^\infty_{y,v}}\\
%----------------------
%       & + \a{C_*}  \left[ 1- \a o(1)\right] ( 1- \a)^{n-1} \|h_0\|_{L^\infty_{y,v}}\\
%----------------------
%   & + \a{C_*}  ( 1- \a)^{n}\|h_0\|_{L^\infty_{y,v}}\\
%----------------------
 %&
%  +  (1-\a)^{n+1} \|h_0\|_{L^\infty_{y,v}}\\
      \leq  \,
    & \a {C_*}\sum_{i=0}^{n} [ 1- \a o(1)]^{n-i}( 1- \a)^{i} \|h_0\|_{L^\infty_{y,v}}
    +     (1-\a)^{n+1} \|h_0\|_{L^\infty_{y,v}}.
%    \leq  \,
%    &  {C_*} \left[ 1- \a o(1)\right]^{n+1} \left[ 1- \left(\frac{1-\a}{1-\a %o(1)}\right)^{n+1}\right] \frac{1}{1-o(1)}\|h_0\|_{L^\infty_{y,v}}
%   +     (1-\a)^{n+1} \|h_0\|_{L^\infty_{y,v}},
\end{split}
\end{equation}
 This, combined with \eqref{o(1)-bd}, yields the uniform bound for $h^{n+1}(\bar{t})$:
\begin{equation}\label{h^{n+1}-bd-uniform-n}
\begin{split}
   &\sup_{0\leq s \leq T_0}\big[   e^{\nu_0 s }\|h^{n+1}( s ) \|_{L^\infty_{y,v}(\O_
\e\times{R}^3)} \big]
   \\
    \leq \,
    &   {C_*} \Big[ 1- \Big(\frac{1-\a}{1-\a o(1)}\Big)^{n+1}\Big] \frac{1}{1-o(1)} \|h_0\|_{L^\infty_{y,v}}
    +  \Big(\frac{1-\a}{1-\a o(1)}\Big)^{n+1}    \|h_0\|_{L^\infty_{y,v}}
    \\
    \leq \,
    & (2C_*+1)   \|h_0\|_{L^\infty_{y,v}}.
\end{split}
\end{equation}

Combing \eqref{h^{n+1}-bd-uniform-k} in Case 1 and \eqref{h^{n+1}-bd-uniform-n} in Case 2, we verify the claim \eqref{semigroup-esti-local-h^{n+1}}.
%Note that for the pure specular reflection boundary ($\a=0$) and the pure diffusive
%boundary ($\a=1$), sharper estimates are available:
%$$
%\sup_{0\leq s \leq T_0}\left[   e^{\nu_0 s } \left\|h^{n+1}( s ) %\right\|_{L^\infty_{y,v}(\O_
%\e\times{R}^3)} \right]\leq
%\left\{
%  \begin{array}{ll}
%      \|h_0\|_{L^\infty_{y,v}}, & \text{ if }  \a=0; \\[0.3mm]
%     2C_*   \|h_0\|_{L^\infty_{y,v}}, & \text{ if }\a=1,
%  \end{array}
%\right.
%$$
%as can be seen from  \eqref{h^{n+1}-bd-uniform-k} in Case 1 and %\eqref{h^{n+1}-bd-uniform-n} in Case 2.
Note that excluding the zero-measure sets $\g_0$ and $S_y(v)$ does not affect this uniform $L^\infty$ estimate.
\medskip

\noindent\textbf{Step 2. Proof of the uniform estimate \eqref{semigroup-esti}.}\;

 From \eqref{semigroup-esti-local}, we obtain
\begin{equation}\label{semigroup-esti-0}
\begin{split}
&\|h(T_0) \|_{_{L^\infty_{y,v}}}
\leq (2C_*+1)e^{-\nu_0T_0}\| h(0)\|_{L^\infty_{y,v}}
\leq  e^{-\frac{\nu_0}{2}T_0} \| h_0\|_{L^\infty_{y,v}},
\end{split}
\end{equation}
provided $T_0$ is sufficiently large. Then, we apply the estimate \eqref{semigroup-esti-0} iteratively on the intervals $[T_0,2T_0]$, $[2T_0,3T_0]$, $\cdots$, $[(j-1)T_0,jT_0]$ ($j\in \mathbb{Z}_+$), yielding
\begin{equation}\label{semigroup-esti-nT0}
\begin{split}
\|h(jT_0) \|_{_{L^\infty_{y,v}}}
&\leq  e^{-\frac{\nu_0}{2}T_0} \| h\big((j-1)T_0\big) \|_{L^\infty_{y,v}}
%\leq  e^{-2\frac{\nu_0}{2}T_0} \left\| h((j-2)T_0) \right\|_{L^\infty_{y,v}}\\&
\leq \cdots
\leq  e^{-j\frac{\nu_0}{2}T_0} \| h_0\|_{L^\infty_{y,v}}.
\end{split}
\end{equation}
Finally, for an arbitrary $\bar{t}>0$, choose $j\in \mathbb{Z}^+$
such that $jT_0\leq \bar{t}<(j+1)T_0$. Applying \eqref{semigroup-esti-local} on
the interval $[jT_0, \bar{t}]$ and using \eqref{semigroup-esti-nT0},
we obtain
\begin{equation}\label{semigroup-esti-t}
\begin{split}
\|h(\bar{t}) \|_{_{L^\infty_{y,v}}}
&\leq  (2C_*+1)e^{-\nu_0(\bar{t}-jT_0)}\| h(jT_0)\|_{L^\infty_{y,v}}
%\leq  (2C_*+1)e^{-\frac{\nu_0}{2}\bar{t}} e^{-\frac{\nu_0}{2}(\bar{t}-jT_0)}
% \| h_0 \|_{L^\infty_{y,v}}\\&
\leq  (2C_*+1) e^{-\frac{\nu_0}{2}\bar{t}} \| h_0\|_{L^\infty_{y,v}}.
\end{split}
\end{equation}
This completes the proof of Lemma \ref{semigroup-estimate}.
 \end{proof}
\bigskip

\noindent
\subsection{$L^\infty$ Estimate for the Linear Equation} \

\medskip

We establish the $L^\infty$ estimate for the linear equation \eqref{eq-unst-stretch-t} and give the proof of Proposition \ref{lemma-fbar-infty-unst-0}.
\medskip

\begin{proof}[\textbf{Proof of Proposition \ref{lemma-fbar-infty-unst-0}.}] \
We first claim that, for any given $(\bar{t}, y,v)\in [0, T_0]\times \bar{\O}_\e\times \mathbb{R}^3$ with $(y,v)\notin \gamma_0$ or $v\notin S_y(v)$, the following bounds hold:
\begin{align}
|h(\bar{t},y,v)|  \lesssim \;
& e^{-\frac{\nu_0}{2}\bar{t}}  \|  h_0  \|_{L^\infty_{y,v} (\Omega_\e\times \mathbb{R}^3)}
+  o(1)  \sup_{0\leq s\leq T_0}\|h(s)\|_{L^{\infty}_{y,v}  (\Omega_\e\times \mathbb{R}^3)} \nonumber \\
%-----------------------
& +  \sup_{0\leq s\leq T_0}  \|\mathbf{P}\bar{f}(s)\|_{L^6_{y,v}  (\Omega_\e\times \mathbb{R}^3)}
 +  \sup_{0\leq s\leq T_0} \|(\mathbf{I-P})\bar{f}(s)\|_{L^2_{y,v} (\Omega_\e\times \mathbb{R}^3)}\label{|h(bar-t,y,v)|-L6}\\
%-----------------------
&+    \e \sup_{0\leq s\leq T_0} \|\langle v\rangle^{-1} \o \bar{g}(s)\|_{L^\infty_{y,v} (\Omega_\e\times \mathbb{R}^3)},\nonumber \\
%----------------------------
%----------------------------
|h(\bar{t},y,v)|  \lesssim \;
& e^{-\frac{\nu_0}{2}\bar{t}}    \|  h_0  \|_{L^\infty_{y,v} (\Omega_\e\times \mathbb{R}^3)}
+  o(1)  \sup_{0\leq s\leq T_0}\|h(s)\|_{L^{\infty}_{y,v}  (\Omega_\e\times \mathbb{R}^3)}\nonumber \\
%-----------------------
& +  \sup_{0\leq s\leq T_0} \|\bar{f}(s)\|_{L^2_{y,v} (\Omega_\e\times \mathbb{R}^3)}
   +  \e \sup_{0\leq s\leq T_0} \|\langle v\rangle^{-1} \o \bar{g}(s)\|_{L^\infty_{y,v} (\Omega_\e\times \mathbb{R}^3)}.\label{|h(bar-t,y,v)|-L2}
\end{align}
 Once \eqref{|h(bar-t,y,v)|-L6} and \eqref{|h(bar-t,y,v)|-L2} are verified, the main estimates
 \eqref{fbar-infty-unst} and \eqref{fbar-infty-L2-unst} follow by applying \eqref{def-h} and taking the $L^{\infty}_{y,v}$ norm on both sides.
 Note that excluding the zero-measure sets $\g_0$ and $S_y(v)$ in  \eqref{|h(bar-t,y,v)|-L6} and \eqref{|h(bar-t,y,v)|-L2} does not affect the validity of the uniform $L^\infty$ estimate.

 We now establish the estimates \eqref{|h(bar-t,y,v)|-L6} and \eqref{|h(bar-t,y,v)|-L2}.
 From \eqref{eq-unst-stretch-t-h-with-K}, for $t_1< s \leq \bar{t}$,
 we have
\begin{equation}\label{h-formula}
\begin{split}
&\frac{d}{ds } \Big[ e^{-\int^{\bar{t}}_s \nu(v) \dd \tau }  h \big(s, Y(s;\bar{t},y,v), v \big)  \Big] \\
= \; &  e^{-\int^{\bar{t}}_s\nu(v) \dd \tau }\Big[
\int_{\mathbb{R}^3} \mathbf{k}_{\beta}(v,u) \frac{ \o (v)}{ \o (u)}  h \big( s,Y(s;\bar{t},y,v), u \big)\dd u  + \e \big(  \o \bar{g} \big) \big(s, Y(s;\bar{t},y,v), v \big)\Big] .
\end{split}
\end{equation}
Define  the principal set
\begin{equation}\label{main-set-mathcal-M}
\begin{split}
 \mathcal{M}(y,v):=\Big\{ (y,v)\in \overline{\O}_\e\times \mathbb{R}^3 :\; |v|\leq N \text{ and  }  \Big|v\cdot \frac{\nabla_x\xi(\e y)}  {|\nabla_x\xi(\e y)|}\Big| \geq \eta \Big\},
\end{split}
\end{equation}
where $N > 0$ is a large constant and $\eta > 0$ is a small constant, both to be specified later. Let $\e_1$ and $\e_2$ be the small constants defined in \eqref{def-epsilon-1} of Lemma \ref{1-bounce-claim-model} and in \eqref{def-epsilon-2} of Lemma \ref{0-bounce-claim-model}, respectively.
Let $\e$ satisfy the restriction
\begin{equation}\label{epsilon-12}
\begin{split}
 0< \e \leq \e_0:=\min\{\e_1, \e_2 \}.
\end{split}
\end{equation}

The proofs of  \eqref{|h(bar-t,y,v)|-L6} and \eqref{|h(bar-t,y,v)|-L2} are divided into two steps.
\medskip

\noindent \textbf{Step 1. Estimate of $h(\bar{t}, y,v)\mathbf{1}_{ \mathcal{M}(y,v)}$.}

Applying the Duhamel principle along the backward trajectory, we obtain
\begin{equation}\label{h(t,y,v){1}}
\begin{split}
 h(\bar{t}, y,v)\mathbf{1}_{ \mathcal{M}{(y,v)}}= J_0(\bar{t},y,v)  +  J_k(\bar{t},y,v)  +  J_g(\bar{t},y,v)  +  J_{sp}(\bar{t},y,v)  +  J_{di}(\bar{t},y,v),
\end{split}
\end{equation}
where
\begin{equation}\label{J-t}
\begin{split}
 J_0(\bar{t},y,v) & := \mathbf{1}_{\{t_1\leq 0\}}  e^{-\int^{\bar{t}}_0 \nu(v)\dd \tau }
   \big|h  \big(0, Y(0), v \big) \big|,   \\
%----------------------
J_k(\bar{t},y,v) & :=\int_{\max{\{0, t_1\}}}^{\bar{t}} \dd s \ e^{-\int^{\bar{t}}_{s}
   \nu(v) \dd \tau }
  \int_{\mathbb{R}^{3}} \dd u \  \mathbf{k}_{\beta}(v,u) \frac{ \o (v)}{ \o (u)}  h\big(s, Y(s), u \big),    \\
%-----------------%---------------
J_g(\bar{t},y,v) & := \int_{\max{\{0, t_1\}}}^{\bar{t}} \dd s \ e^{-\int^{\bar{t}}_{s} \nu(v)\dd \tau }
   \big| \e \big(  \o \bar{g} \big)\big(s, Y(s),v \big)\big|,    \\
%-----------------%---------------
J_{sp}(\bar{t},y,v) &  := \mathbf{1}_{\{t_1>0\}}     e^{-\int^{\bar{t}}_{t_1} \nu(v)\dd \tau }  (1-\a)
  \big| h(t_1,y_1,\,v_1)\big|,  \\
%-----------------%---------------
J_{di}(\bar{t},y,v)  &  :=  \mathbf{1}_{\{t_1>0\}}    e^{-\int^{\bar{t}}_{t_1} \nu(v) \dd \tau}  \a
 \int_{n(y_1)\cdot v_1^*>0}   \big |h(t_1,y_1,v_1^*)\big|  \dd\sigma^{*}_1.
\end{split}
\end{equation}

Direct estimates yield
\begin{equation} \label{J0g-bd}
\begin{split}
& \left |J_0(\bar{t},y,v)\right | \lesssim e^{-\nu_0\bar{t}}\|h_0\|_{L^\infty_{y,v}},\quad
\left |J_g(\bar{t},y,v)\right | \lesssim \e \sup_{0\leq s\leq T_0}\|{\langle v\rangle}^{-1} w\bar{g}(s)\|_{L^\infty_{y,v}}.
\end{split}
\end{equation}
We now estimate the remaining terms $J_{sp}(\bar{t},y,v)$, $J_{k}(\bar{t},y,v)$ and $J_{di}(\bar{t},y,v)$ in Steps 1.1--1.3.

\noindent \textbf{Step 1.1. Estimate of $J_{sp}(\bar{t},y,v)$.}

  By Lemma \ref{1-bounce-claim-model} and \eqref{epsilon-12},  the specular backward trajectory starting from $(\bar{t},y,v)\in [0,T]\times M(y,v)$ undergoes at most single-bounce against $\pt\O_{\e}$. Thus, after the first collision at $(t_1,{y}_1, {v}_1)$, the term $J_{sp}(\bar{t},y,v)$
 propagates back to the initial plane $\{\bar{t}=0\}$:
\begin{align}
 J_{sp}(\bar{t},y,v)
  = & \mathbf{1}_{\{t_1>0\}}  e^{-\int^{\bar{t}}_{t_1}\nu(v_1)\dd \tau}
   e^{-\int^{t_1}_{0}\nu(v_1)\dd \tau }   h \big(0, Y_1(0), v_1 \big)\nonumber  \\
%-----------------%--------------------
& +\mathbf{1}_{\{t_1>0\}}  e^{-\int^{\bar{t}}_{t_1}\nu(v_1)\dd \tau}
 \int_0^{t_1} \dd s \ e^{-\int^{t_1}_{s}\nu(v_1)\dd \tau}     \e \big( \o \bar{g}\big) \big(s, Y_1(s), v_1 \big)     \label{Jsp-decom}  \\
%-----------------%--------------------
& +\mathbf{1}_{\{t_1>0\}}    e^{-\int^{\bar{t}}_{t_1}\nu(v_1)\dd \tau}
 \int_0^{t_1} \dd s \ e^{-\int^{t_1}_{s}\nu(v_1)\dd \tau}  \int_{\mathbb{R}^{3}} \dd v^{\prime}\
\mathbf{k}_{\beta}(v,u) \frac{ \o (v)}{ \o (u)}
        h \big(s, Y_1(s),v^{\prime} \big)\nonumber  \\
 := & J_{sp,0}+J_{sp,g}+J_{sp,k},\nonumber
\end{align}
where we have used the abbreviation
\begin{equation}\label{Y1s-def}
 Y_1(s):=Y(s;t_1,y_1,v_1).
\end{equation}

The terms $J_{sp,0}$ and $J_{sp,g}$ are bounded similarly to \eqref{J0g-bd}.
%Similar to (\ref{J0g-bd}), we have
%\begin{equation}\label{Jsp-0g-bd}
%\begin{split}
%& \left | J_{sp,0} \right |\lesssim e^{-\nu_0\bar{t}}\|h_0\|_{L^{\infty}_{y,v}},\quad
%  \left | J_{sp,g}\right | \lesssim \e  \sup_{0\leq s\leq %T_0}\|w\bar{g}(s)\|_{L^{\infty}_{y,v}}.
%\end{split}
%\end{equation}
To estimate $J_{sp,k}$, we invoke Lemma 3 from \cite{guo2010decay}, which ensures the existence of $\tilde{\beta}=\tilde{\beta}(\beta, \beta^{\prime})>0$ such that
\begin{equation}\label{k-tilde-beta}
\begin{split}
\mathbf{k}_{\beta}(v,u) \frac{ \o (v)}{ \o (u)} \lesssim \mathbf{k}_{\tilde{\beta}} (v,u).
\end{split}
\end{equation}
Moreover, for any $m\geq 1$, we can choose $N=N(m)\gg 1$ further large so that
\begin{equation*}\label{k-cut}
\begin{split}
 & \mathbf{k}_{N} (V, v^{\prime}) : = \mathbf{1}_{|V-v^{\prime}|\geq \frac{1}{N}}
\mathbf{1}_{|v^{\prime}| \leq N} \mathbf{1}_{|V| \leq N}
\mathbf{k}_{\tilde{\beta} } (V, v^{\prime}), \\
&\sup_{V} \int_{\mathbb{R}^{3}} | \mathbf{k}_{N} (V, v^{\prime})
- \mathbf{k}_{\tilde{\beta} } (V, v^{\prime}) | \dd v^{\prime} \leq \frac{1}{m}.
\end{split}
\end{equation*}
We decompose the kernel as
\begin{equation}\label{k-split}
\mathbf{k}_{\tilde{\beta} }  (V, v^{\prime}) = [ \mathbf{k}_{\tilde{\beta} }  (V, v^{\prime})- \mathbf{k}_{N}  (V, v^{\prime})  ] \
+  \  \mathbf{k}_{N}  (V, v^{\prime}).
\end{equation}
The first term in \eqref{k-split} contributes at most $o(1) \|h\|_{L^\infty_{y,v}}$ for
sufficiently large $m \gg 1$.
For $y'\in \overline{\O}_\e$, define the principal set
\begin{equation}\label{main-set}
\begin{split}
 M_{y'}(v'):=\Big \{ v'\in \mathbb{R}^3 :\; |v'|\leq N \text{ and  }  \Big|v'\cdot \frac{\nabla_x\xi(\e y')}  {|\nabla_x\xi(\e y')|}\Big | \geq \eta
\Big\}.
\end{split}
\end{equation}
The second term in \eqref{k-split} leads to
$$
C_m \int_{|v'|\leq N,\;|v'\cdot \frac{\nabla_x\xi(\e Y_1(s))}  {|\nabla_x\xi(\e Y_1(s))|}|<\eta}   +   C_m\int_{M'_{Y_1(s)}(v')},
$$
which is further bounded by
\begin{equation}\label{Jspk-bd}
\begin{split}
 o(1) \sup_{0\leq s\leq T_0}\|h(s)\|_{L^\infty_{y,v}}+ C_{m}J_{sp,k*},
\end{split}
\end{equation}
where  $\eta>0$ is chosen sufficiently small, and
\begin{equation}\label{Jspk*}
\begin{split}
 J_{sp,k*}:=&\mathbf{1}_{\{t_1>0\}} \int_0^{t_1} \dd s \ e^{-\nu_0(t-s)}
\int_{ M'_{Y_1(s)}(v')} \dd v^{\prime}
\big|{\underbrace{  h\big(s,Y_1(s),v^{\prime}\big) }}\big|.
\end{split}
\end{equation}
Note that the $o(1)$ coefficient in \eqref{Jspk-bd} depends on $N$ and $\eta$ but
is independent of $\e$.

We now apply the Duhamel principle and \eqref{h-formula} to the under braced term in \eqref{Jspk*},
considering the backward trajectory starting from $(s,Y_1(s),v')$:
\begin{equation}\label{Jspk*-duh}
\begin{split}
 J_{sp, k*}= &\; J_{sp, k*}^{0}  +  J_{sp, k*}^{k}  +  J_{sp, k*}^{g} + J_{sp, k*}^{sp} +  J_{sp, k*}^{di},
\end{split}
\end{equation}
where
\begin{align*}
 J_{sp, k*}^{0}=\; &  \mathbf{1}_{\{t_1>0\}}\mathbf{1}_{\{t^\prime_1\leq 0\}} \int_0^{t_1}\dd s   \  e^{-\nu_0(\bar{t}-s)}
\int_{M_{Y_1(s)}(v')} \dd v^{\prime}
  e^{-\int^{s}_{0}\nu(v')\dd \tau}\big|h\big(0,Y(0; s,Y_1(s),v'), v'\big)\big|,\\
%---------------------------------------
 J_{sp, k*}^{k}=\;  &  \mathbf{1}_{\{t_1>0\}} \int_0^{t_1} \dd s \ e^{-\nu_0(\bar{t}-s)}
\int_{M_{Y_1(s)}(v')}   \dd v^{\prime}\int_{\max{\{0, t^\prime_1\}}}^{s} \dd \tau
 e^{-\int^{s}_{\tau}\nu(v')\dd \tau} \\
&\times
\int_{\mathbb{R}^{3}} \dd u \   \mathbf{k}_{\beta}(v^{\prime},u) \frac{ \o (v^{\prime})}{ \o (u)}
 \big| h\big(\tau, Y(\tau; s,Y_1(s),v^{\prime}),u\big)\big|,\\
%----------------------------------------
 J_{sp, k*}^{g} = \; &  \mathbf{1}_{\{t_1>0\}} \int_0^{t_1} \dd s \ e^{-\nu_0(\bar{t}-s)}
\int_{M_{Y_1(s)}(v')}   \dd v^{\prime}\int_{\max{\{0, t^\prime_1\}}}^{s} \dd \tau
e^{-\int^{s}_{\tau} \nu(v')  \dd \tau}\\
&\times
\big| \e \big(   \o \bar{g} \big)   \big(\tau,  Y(\tau;s,Y_1(s),v'), v' \big) \big|,
\\
%------------------------------------
 J_{sp, k*}^{sp}=\; &  \mathbf{1}_{\{t_1>0\}} \mathbf{1}_{\{t^\prime_1>0\}} \int_0^{t_1} \dd s \ e^{-\nu_0(\bar{t}-s)}
\int_{M_{Y_1(s)}(v')}   \dd v^{\prime}
 e^{-\int^{s}_{t^\prime_1}  \nu(v')  \dd \tau}
\big|h\big(t_{1}^\prime, y_{1}^\prime,v_1^\prime\big)\big|,
\\
%------------------------------------
 J_{sp, k*}^{di}= \;&  \mathbf{1}_{\{t_1>0\}} \mathbf{1}_{\{t^\prime_1>0\}}
 \int_0^{t_1} \dd s \ e^{-\nu_0(\bar{t}-s)}
\int_{M_{Y_1(s)}(v')}  \dd v^{\prime}
e^{-\int^{s}_{t^\prime_1}\nu(v^{\prime})\dd \tau}
\int_{n(y'_1) \cdot u_1^{\prime*}>0}
\big|h\big(t'_1,y_{1}^\prime,u^{\prime*}_{1}\big)\big|\dd\sigma^{\prime*}_1.
\end{align*}
%where $v_1^\prime=R_{y_{1}^\prime}(v')$ and $\dd \sigma^{\prime*}_1$ is defined as in %\eqref{def-C-star}.
%\dd \sigma^{\prime*}_1$
%\begin{equation}\label{y'1-def}
%\begin{split}
%&t_{1}^\prime:=s-t_{\mathbf{b}}(s, Y_1(s),v'), \ \ \ \
%y_{1}^\prime:=Y(t^\prime_1;s,Y_1(s),v')= y_{\bf b}(Y_1(s),v'), \\
%&  v_1^\prime:=R_{y_{1}^\prime}(v'), \qquad\qquad\quad\;\;\;
%\dd \sigma^{\prime*}_1:=\sqrt{2\pi \mu(u^{\prime*}_{1})}w^{-1}(u^{\prime*}_{1})  %[n(Y^\prime(t^\prime_1)) \cdot u^{\prime*}_{1} ] \dd u^{\prime*}_{1}.
%\end{split}
%\end{equation}
%---------------------
The terms $J_{sp, k*}^{0}$ and $J_{sp, k*}^{g}$ are bounded similarly to \eqref{J0g-bd}.
%\begin{equation}\label{Jspk*-0-g-bd}
%\begin{split}
%&\big| J_{sp, k*}^{0} \big |
%\lesssim \int_0^{t_1} e^{-\nu_0\bar{t}}\dd s \|h_0\|_{L^\infty_{y,v}}
%\leq \bar{t} e^{-\nu_0\bar{t}}\|h_0\|_{L^\infty_{y,v}}
%\leq e^{-\frac{\nu_0}{2}\bar{t}}\|h_0\|_{L^\infty_{y,v}},\quad
%\big|J_{sp, k*}^{g}\big|
%\lesssim \e\sup_{0\leq s\leq T_0}\|{\langle v\rangle}^{-1} w\bar{g}(s)\|_{L^\infty_{y,v}}.
%\end{split}
%\end{equation}
The remaining terms $J_{sp,k*}^{sp} $, $J_{sp, k*}^{k}$ and $J_{sp,k*}^{di} $ will be estimated in the following Steps 1.1.1--1.1.3.

\noindent \textbf{Step 1.1.1. Estimate of $J_{sp,k*}^{sp}$.}

 For $(s,Y_1(s),v^\prime)$ with $0 \leq t'_1\leq s<t_1<\bar{t}\leq T_0$, $Y_1(s)\in\bar{\Omega}_\e$ and $v'\in M_{Y_1(s)}(v')$,  similarly as \eqref{Jsp-decom}, Lemma \ref{1-bounce-claim-model} ensures that the specular backward trajectory
 starting from $(s,Y_1(s),v^\prime)$ reaches the initial plane \{$\bar{t}=0$\} after the first collision at $({y}_1^\prime, {v}_1^\prime)$. Thus,
\begin{equation}\label{Jspk*-sp}
 J_{sp,k*}^{sp} =  J_{sp,k*}^{sp,0} +   J_{sp,k*}^{sp,g} +  J_{sp,k*}^{sp,k},
\end{equation}
where
\begin{align*}
 J_{sp,k*}^{sp,0}=&\mathbf{1}_{\{t_1>0\}}\mathbf{1}_{\{t^\prime_1>0\}}
  \int_0^{t_1} \dd s \ e^{-\nu_0(\bar{t}-s)}
\int_{M_{Y_1(s)}(v')}
\dd v^{\prime}\\
%------------------
&  \times
e^{-\int^{s}_{t^\prime_1}\nu({v}^\prime)\dd \tau}  e^{-\int^{t^\prime_1}_{0}\nu({v}_1^\prime)\dd \tau}
 \big| h \big( 0, Y(0; t_1^\prime, {y}_1^\prime, {v}_1^\prime), {v}_1^\prime\big)\big|,\\
%-------------------
J_{sp,k*}^{sp,g}=&\mathbf{1}_{\{t_1>0\}}\mathbf{1}_{\{t^\prime_1>0\}}
\int_0^{t_1} \dd s \ e^{-\nu_0(\bar{t}-s)}
\int_{M_{Y_1(s)}(v')}
\dd v^{\prime}
\\
%------------------
&\times
 e^{-\int^{s}_{t^\prime_1}\nu(v^\prime) \dd \tau}
 \int_0^{t^\prime_1} \dd\tau  \  e^{-\int^{t^\prime_1}_{\tau}  \nu({v}_1^\prime)\dd \tau^{\prime}} \big| \e \big(  \o  \bar{g} \big)  \big(\tau,  Y(\tau; t_1^\prime, {y}_1^\prime, {v}_1^\prime),{v}_1^\prime \big) \big|,\\
%--------------------------
J_{sp,k*}^{sp,k}=
&\mathbf{1}_{\{t_1>0\}}\mathbf{1}_{\{t^\prime_1>0\}}
\int_0^{t_1} \dd s \ e^{-\nu_0(\bar{t}-s)}
\int_{M_{Y_1(s)}(v')}   \dd v^{\prime}
e^{-\int^{s}_{t^\prime_1}\nu(v^\prime)\dd \tau}
  \int_0^{t^\prime_1} \dd \tau  \\
%------------------
&\times
e^{-\int^{t^\prime_1}_{\tau}     \nu({v}_1^\prime)      \dd \tau^{\prime}}
  \int_{\mathbb{R}^{3}} \dd \tilde{u}\  \mathbf{k}_{\beta}({v}_1^\prime,\tilde{u}) \frac{ \o ({v}_1^\prime)}{ \o (\tilde{u})}
\ \big|  h\big(\tau,  Y(\tau; t_1^\prime, {y}_1^\prime, {v}_1^\prime), \tilde{u}  \big)  \big|.
\end{align*}

The terms $J_{sp,k*}^{sp,0}$ and $J_{sp,k*}^{sp,g}$ are estimated similarly to  \eqref{J0g-bd}. For  $J_{sp,k*}^{sp,k}$, we proceed as in Step 1.1: bound the kernel by $\mathbf{k}_{\tilde{\beta} }(V^\prime,\tilde{u})$, decompose it as $[ \mathbf{k}_{\tilde{\beta} } (V^\prime,\tilde{u})- \mathbf{k}_{N} (V^\prime,\tilde{u})  ] + \mathbf{k}_{N} (V^\prime,\tilde{u})$, and split the time interval $[0, \; t^\prime_1]=[0,\;  t^\prime_1-\delta] \cup  [t^\prime_1-\delta, \; t^\prime_1 ]$.   This yields
\begin{equation}\label{Jspk*-spk-decom}
\begin{split}
 J_{sp,k*}^{sp,k}\lesssim o(1)\sup_{0\leq s\leq T_0}\|h(s)\|_{L^\infty_{y,v}}  +  J_{sp,k*}^{sp,k*},
\end{split}
\end{equation}
where
\begin{equation}\label{Jspk*-spk*-def}
\begin{split}
 J_{sp,k*}^{sp,k*}:=\;&\mathbf{1}_{\{t_1>0\}}\mathbf{1}_{\{t^\prime_1>0\}} \int_0^{t_1} \dd s \ e^{-\nu_0(t-s)}
\int_{M'_{Y_1(s)}(v')}  \dd v^{\prime}\\[0.1cm]
%------------------
&\times \int_0^{t^\prime_1-\delta} \dd \tau \ e^{-\nu_0(s-\tau)}\int_{|\tilde{u}|\leq N} \dd  \tilde{u}
 \ \big|\underbrace{h\big(\tau, Y(\tau; t_1^\prime, {y}_1^\prime, {v}_1^\prime), \tilde{u}\big)}\big|.
\end{split}
\end{equation}

Now consider the change of variables:
\begin{equation}\label{v'-Y'1}
 v^\prime\mapsto Y_1^\prime(\tau)  :=  Y(\tau; t_1^\prime, {y}_1^\prime, {v}_1^\prime)= y_1^\prime  +  (\tau-t^\prime_1)v_1^\prime.
\end{equation}
Since $0<t^\prime_1<s<t_1<t\leq T_0$ and $|v^\prime|\leq N$,
the relation
\begin{equation}\label{y'1}
\begin{split}
  y_1^\prime=Y(t^\prime_1;s,Y_1(s),v^\prime)
&   =  Y_1(s)  +  (t^\prime_1-s)v^\prime
\end{split}
\end{equation}
implies that $|(t^\prime_1-s)v'|\leq NT_0$. This further indicates that $|y_1^\prime-Y_1(s)|$ must be bounded by $NT_0$. While $y^{\prime}_1\in \pt\O_{\e}$ and $\e y^{\prime}_1\in \pt\O$,
so that $\e Y_1(s)\in \O$ lies near the boundary $\partial\O$, and thus \begin{equation*}%\label{intfor}
\begin{split}
n(\e Y_1(s))  =  \frac{\nabla_x\xi(\e Y_1(s))}  {|\nabla_x\xi(\e Y_1(s))|},
\end{split}
\end{equation*}
where we used the fact that $\nabla_x \xi(\e Y_1(s))\neq 0$ near the boundary.
It follows that
\begin{equation*}\label{ny'1}
\begin{split}
&n(y_1^\prime)=n(\e y_1^\prime)=n\big(\e Y_1(s)+\e (t^\prime_1-s)v^\prime \big) =n(\e Y_1(s))+O(\e).
\end{split}
\end{equation*}
Thus, we derive
\begin{equation}\label{v'1}
\begin{split}
 v_1^\prime  &=R_{y_{1}^\prime}(v^\prime \big )
 =v^\prime-2[n(y_{1}^\prime) \cdot v^\prime]\ n(y_{1}^\prime)
 =v^\prime-2\big[n(\e Y_1(s))\cdot v^\prime\big] n(\e Y_1(s))+O(\e).
\end{split}
\end{equation}
It follows from \eqref{v'-Y'1}, \eqref{y'1} and  \eqref{v'1} that
\begin{equation}%\label{intfor}
\begin{split}
&Y(\tau; t_1^\prime, {y}_1^\prime, {v}_1^\prime)\\
%-------------------------------
%=&y_1^\prime  +  (\tau-t^\prime_1)v_1^\prime \\
%-------------------------------
=&Y_1(s)  +  (t^\prime_1-s)v^\prime  +  \big[(\tau-s)-(t^\prime_1-s)\big]
  \big \{ v^\prime-2\big[n(\e Y_{1}(s))\cdot v^\prime\big]  n(\e Y_{1}(s))\big \}  +  O(\e)\\
%---------------------------
=&Y_1(s)  +  (\tau-s)v^\prime  + 2\big[(t^\prime_1-s)- (\tau-s)\big]  \big[n(\e Y_1(s))\cdot v^\prime\big]\ n(\e Y_1(s))  +  O(\e)\\
%-------------------------------
=&Y_1(s)  +  (\tau-s)v^\prime  +2\big[(t^\prime_1-s)- (\tau-s)\big] \frac{\big[\nabla_x\xi(\e Y_1(s))\cdot v^\prime\big]
  \nabla_x\xi(\e Y_1(s))}  {|\nabla_x\xi(\e Y_1(s))|^2} + O(\e).
\end{split}
\end{equation}
We now compute the Jacobian entries:
\begin{equation}\label{Y'1i-v'j}
\begin{split}
\frac{\partial Y(\tau; t_1^\prime, {y}_1^\prime, {v}_1^\prime)_i}  {\partial v^\prime_j}
&=(\tau-s)\delta_{ij}+ 2\big[(t^\prime_1-s)- (\tau-s)\big]
\frac{\partial_i \xi(\e Y_1(s))\ \partial_j \xi(\e Y_1(s))}{|\nabla_x\xi(\e Y_1(s))|^2} \\
%------------------------
&\ \ \ +\frac{2\big[\nabla_x\xi(\e Y_1(s))\cdot v^\prime\big]}   {|\nabla_x\xi(\e Y_1(s))|^2}
   \partial_i\xi(\e Y_1(s))   \frac{\partial(t^\prime_1-s)}{\partial v^\prime_j}    +O(\e)\\
 &:= (\tau-s)\delta_{ij}  +a_{ij}  +O(\e),\qquad i,j=1,2,3,
\end{split}
\end{equation}
where $\partial_i\xi=\frac{\partial \xi}{\partial x_i}$ denotes the spatial derivative, with  the notations
\begin{equation*}
\begin{split}
& a_{ij}:=b_ic_j,  \ \ \ \  b_i:=\partial_i \xi(\e Y_1(s)), \\
  &c_j:=2\big[(t^\prime_1-s)- (\tau-s)\big]\frac{\partial_j \xi(\e Y_1(s))}{|\nabla_x\xi(\e Y_1(s))|^2}
+\frac{2\big[\nabla_x\xi(\e Y_1(s))\cdot v^\prime\big]}   {|\nabla_x\xi(\e Y_1(s))|^2}
  \frac{\partial(t^\prime_1-s)} {\partial v^\prime_j}.
\end{split}
\end{equation*}
Elementary computations yield:
\begin{equation*}
\begin{split}
& \sum_{k=1}^3 a_{kk}=2\big[(t^\prime_1-s)- (\tau-s)\big]
 +  \frac{2\big[\nabla_x\xi(\e Y_1(s))\cdot v^\prime\big]}   {|\nabla_x\xi(\e Y_1(s))|^2}
 \big[\nabla_v (t^\prime_1-s)\cdot \nabla_x\xi(\e Y_1(s))\big], \\
&\det B_{ij}
:=\det \left(
           \begin{array}{cc}
             a_{ii} & a_{ij} \\
             a_{ji} & a_{jj} \\
           \end{array}
         \right)
         =\det \left(
           \begin{array}{cc}
             b_i c_i & b_ic_j  \\
             b_jc_i  & b_jc_j  \\
           \end{array}
         \right)=0 \ \ \ \ \ \hbox{ for }  i\neq j,     \\
&\det C:=\det
  \left(
  \begin{array}{ccc}
     a_{11} & a_{12} & a_{13}  \\
     a_{21} & a_{22} & a_{23}  \\
     a_{31} & a_{32} & a_{33} \\
           \end{array}
         \right)
= \det \left(
  \begin{array}{ccc}
      b_{1}c_{1} & b_{1}c_{2} & b_{1}c_{3}  \\
      b_{2}c_{1} & b_{2}c_{2} & b_{2}c_{3}  \\
      b_{3}c_{1} & b_{3}c_{2} & b_{3}c_{3} \\
           \end{array}
         \right)=0.
     \end{split}
\end{equation*}
From these relations and \eqref{Y'1i-v'j}, we obtain
\begin{align}
&\det\big[ \nabla_{v^\prime}Y(\tau; t_1^\prime, {y}_1^\prime, {v}_1^\prime) \big] \nonumber\\
=&(\tau-s)^3  +  (\tau-s)^2 \sum_{k=1}^3 a_{kk} + (\tau-s)\sum_{1\leq i<j\leq 3}^3 \det B_{ij}
+ \det C +O(\e), \label{det-Y'1-v'}\\
=&-(\tau-s)^3+2(\tau-s)^2\Big\{(t^\prime_1-s)+\frac{\big[\nabla_x \xi(\e Y_1(s))\cdot v^\prime\big]}
{|\nabla_x\xi(\e Y_1(s))|^2}\big[\nabla_{v^\prime}(t^\prime_1-s)\cdot\nabla_x\xi(\e Y_1(s))\big]\Big\}  + O(\e).\nonumber
\end{align}

Recall that $y^\prime_1\in \partial\O_\e$, $\e y^\prime_1\in \partial\O$ and $\e Y_1(s)$ is near the boundary $\partial\O$. Since $\frac{1}{2}C_{\xi_1} \leq |\nabla_x \xi (x)|\leq 2C_{\xi_2} $ for $x$ near the boundary $\partial\O$, we have $\frac{1}{4}C_{\xi_1} \leq |\nabla_x\xi(\e Y_1(s))|\leq 4C_{\xi_2} $. It follows that
\begin{equation*}
\begin{split}
\big| \big[v^\prime \cdot \nabla_x\xi(\e Y_1(s))\big] \big|  & =  |\nabla_x\xi(\e Y_1(s))| \
\big|  \big[v^\prime \cdot  n(\e Y_1(s)) \big] \big|  =  |\nabla_x\xi(\e Y_1(s))| \
\big|  \big[v^\prime \cdot  n(Y_1(s)) \big] \big| \geq \frac{\eta C_{\xi_1}}{4},
\end{split}
\end{equation*}
where we have used the condition $|v^\prime \cdot  n( Y_1(s))|\geq \eta$. From the expansion
\begin{equation*}
\begin{split}
0=\xi\big(\e y^\prime_1\big)&=\xi\big(\e Y_1(s)+\e(t^\prime_1-s)v^\prime\big)
=\xi(\e Y_1(s))+ \e(t^\prime_1-s) \big[\nabla_x\xi(\e Y_1(s)) \cdot  v^\prime\big]   +   O(\e^2),
\end{split}
\end{equation*}
we take the partial derivative $\partial_{v^\prime_j}$:
\begin{equation*}
\begin{split}
(t^\prime_1-s)\partial_{j}\xi(\e Y_1(s))+ \frac{\partial (t^\prime_1-s)}{\partial v^\prime_j}
\big[\nabla_x\xi(\e Y_1(s)) \cdot  v^\prime\big]+O(\e)=0, \;\; j=1,2,3.
\end{split}
\end{equation*}
Taking inner product with $\nabla_x\xi(\e Y_1(s))$ yields
\begin{equation*}
\begin{split}
 (t^\prime_1-s)|\nabla_x\xi(\e Y_1(s))|^2+
 \big[\nabla_x \xi(\e Y_1(s))\cdot v^\prime\big]   \big[\nabla_{v^\prime}(t^\prime_1-s)\cdot\nabla_x\xi(\e Y_1(s))\big]
=  O(\e).
\end{split}
\end{equation*}
It follows that
\begin{equation}\label{det-Y'1-v'-1}
\begin{split}
 & (t^\prime_1-s)  +  \frac{\big[\nabla_x \xi(\e Y_1(s))\cdot v^\prime\big]}{|\nabla_x\xi(\e Y_1(s))|^2}
 \big[\nabla_{v^\prime}(t^\prime_1-s)\cdot\nabla_x\xi(\e Y_1(s))\big]=O(\e).
\end{split}
\end{equation}

Since $0<\tau\leq t^\prime_1-\delta<t^\prime_1<s<t<T_0$, we have $s-\tau>t^\prime_1-\tau\geq \delta$.
Combining \eqref{det-Y'1-v'} and \eqref{det-Y'1-v'-1}, we obtain the lower bound for the Jacobian:
\begin{equation*}%\label{intfor}
\begin{split}
\big | \det\nabla_{v^\prime}Y(\tau; t_1^\prime, {y}_1^\prime, {v}_1^\prime) \big |  \gtrsim |s-\tau|^3
+ O(\e) \geq \frac{1}{2} \delta^3,
\end{split}
\end{equation*}
for sufficiently small $\e\leq \e_1$. Note that this lower bound is independent of $\e$.

Integrating over time first and using $|h(\tilde{u})| = \o (\tilde{u}) |\bar{f}(\tilde{u})| \lesssim_{N} \o^{-1} |\bar{f}(\tilde{u})|$ for $|\tilde{u}| \leq N$, we have
\begin{equation}\label{Jspk*-k*-decom}
\begin{split}
J_{sp, k*}^{sp,k*}
 &\lesssim \sup_{0 \leq \tau \leq s-\delta < s \leq t_1} \int_{|v^\prime|\leq N} \int_{|\tilde{u}|\leq N}
 \big| h\big(\tau, Y(\tau; t_1^\prime, {y}_1^\prime, {v}_1^\prime),\tilde{u}\big)\big| \dd \tilde{u} \dd v^{\prime}\\
%---------------------------
&\lesssim \sup_{0 \leq \tau \leq s-\delta < s \leq t_1} \int_{|v^\prime|\leq N} \int_{|\tilde{u}|\leq N}
 \big| \o^{-1}  \mathbf{P} \bar{f}\big(\tau, Y(\tau; t_1^\prime, {y}_1^\prime, {v}_1^\prime),\tilde{u}\big)\big| \langle \tilde{u}\rangle^2 \sqrt{\mu(\tilde{u})} \dd \tilde{u} \dd v^{\prime}\\
%---------------------------
&\quad + \sup_{0 \leq \tau \leq s-\delta < s \leq t_1} \int_{|v^\prime|\leq N} \int_{|\tilde{u}|\leq N}
 \big| \o^{-1} (\mathbf{I}- \mathbf{P}) \bar{f}\big(\tau, Y(\tau; t_1^\prime, {y}_1^\prime, {v}_1^\prime),\tilde{u}\big)\big| \dd \tilde{u} \dd v^{\prime}\\
 &:=J_{sp, k*}^{sp,k*,1}+J_{sp, k*}^{sp,k*,2}.
\end{split}
\end{equation}
For $\mathbf{P} \bar{f}$ contribution $J_{sp, k*}^{sp,k*,1}$,
\begin{equation}\label{Jspk*-k*-Pf}
\begin{split}
J_{sp, k*}^{sp,k*,1}
%\sup_{0 \leq \tau \leq s-\delta < s \leq t_1} \int_{|v^\prime|\leq N} \int_{|u|\leq N}
% \big| \mathbf{P} \bar{f}\big(\tau, Y(\tau; s,Y_1(s),v^{\prime}),u\big)\big|\langle %u\rangle^2 \sqrt{\mu(u)}  \dd u \dd v^{\prime}\\
%---------------------------
&\lesssim_{N}  \sup_{0 \leq \tau \leq s-\delta < s \leq t_1}\Big[ \int_{v^{\prime}}\big \| \o^{-1} \mathbf{P}\bar{f} \big(\tau, Y(\tau; t_1^\prime, {y}_1^\prime, {v}_1^\prime) \big) \big\|^{6}_{L^6(\mathbb{R}^3_{\tilde{u}})}
 \dd v^{\prime} \Big]^{1/6}\\
%---------------------------
&\lesssim_{N} \sup_{0\leq \tau\leq T_0}\Big[ \int_{\Omega_\e}
  \big \|  \o^{-1}  \mathbf{P} \bar{f} \big(\tau, y  \big) \big \|^{6}_{L^6(\mathbb{R}^3_{\tilde{u}})}
  \frac{2}{\delta^{3}}
  \dd y\Big]^{1/6} \\
%---------------------------
&   \lesssim_{N}     \sup_{0\leq \tau\leq T_0} \|  \o^{-1}   \mathbf{P} \bar{f}(\tau)  \|_{L^{6}(\Omega_\e \times \mathbb{R}^{3} )} .
\end{split}
\end{equation}
For $(\mathbf{I}- \mathbf{P})\bar{f}$ contribution $J_{sp, k*}^{sp,k*,2}$,
\begin{equation}\label{Jspk*-k*-(I-P)f}
\begin{split}
J_{sp, k*}^{sp,k*,2}
%&\sup_{0 \leq \tau < s-\delta < s \leq t_1} \int_{|v^\prime|\leq N} \int_{|u|\leq N}
% \big|(\mathbf{I}- \mathbf{P}) \bar{f}\big(\tau, Y(\tau; s,Y_1(s),v^{\prime}),u\big)\big| %\dd u \dd v^{\prime}\\
%------------------------------
& \lesssim_{N} \sup_{0 \leq \tau < s-\delta < s \leq t_1} \Big[ \iint  \big| \o^{-1} (\mathbf{I}- \mathbf{P}) \bar{f}\big(\tau, Y(\tau; t_1^\prime, {y}_1^\prime, {v}_1^\prime),\tilde{u}\big)\big|^{2}
 \dd v^{ \prime} \dd \tilde{u} \Big]^{1/2}\\
& \lesssim_{N}  \sup_{0\leq \tau\leq T_0}\Big[ \iint_{\Omega_\e \times \mathbb{R}^{3}}  \big| \o^{-1} (\mathbf{I}- \mathbf{P}) \bar{f}\big(\tau,y,\tilde{u}\big)\big|^{2}
 \frac{2}{\delta^{3}} \dd y \dd \tilde{u} \Big]^{1/2}\\
&  \lesssim_{N}\sup_{0\leq \tau\leq T_0} \| \o^{-1} (\mathbf{I} - \mathbf{P})  \bar{f}(\tau)  \|_{L^{2}(\Omega_\e \times \mathbb{R}^{3} )}.
\end{split}
\end{equation}

%Now integrating over time first and deducing similarly as %(\ref{Jspk*-k*-decom})--(\ref{Jspk*-k*-(I-P)f}), we get
%\begin{equation}\label{Jspk*-spk*-bd}
%\begin{split}
%&  J_{sp,k*}^{sp,k*} \lesssim  \sup_{0\leq s\leq T_0}\|\mathbf{P} \bar{f}(s)\|_{L^6_{y,v}}  %+  \sup_{0\leq s\leq T_0}\|(\mathbf{I-P}) \bar{f}(s)\|_{L^2_{y,v}}.
%\end{split}
%\end{equation}
Collecting \eqref{Jspk*-spk-decom}, \eqref{Jspk*-k*-decom}--\eqref{Jspk*-k*-(I-P)f}, we obtain
\begin{equation}\label{Jspk*-sp-bd}
\begin{split}
J_{sp,k*}^{sp}
\lesssim  \; &  e^{-\frac{\nu_0}{2}\bar{t}}\|h_0\|_{L^{\infty}_{y,v}} + o(1) \sup_{0\leq s\leq T_0}\|h(s)\|_{L^\infty_{y,v}}
+ \sup_{0\leq s\leq T_0} \| \o^{-1} \mathbf{P}\bar{f}(s)\|_{L^6_{y,v}}\\
 &
 +  \sup_{0\leq s\leq T_0}\| \o^{-1} (\mathbf{I-P})\bar{f}(s)\|_{L^2_{y,v}}+\e\sup_{0\leq s\leq T_0}\|{\langle v\rangle}^{-1} \o\bar{g}(s)\|_{L^\infty_{y,v}}.
\end{split}
\end{equation}

\noindent \textbf{Step 1.1.2. Estimation of $J_{sp, k*}^{k}$.}

We decompose the kernel in $J_{sp, k*}^{k}$ similarly to \eqref{k-split},
%$$
%\mathbf{k}_{\tilde{\beta} } (v,u) = [ \mathbf{k}_{\tilde{\beta} } (v,u)- \mathbf{k}_{N} %(v,u)  ] + \mathbf{k}_{N} (v,u),
%$$
where  the first term contributes at most (cf. \eqref{Jspk-bd}) $o(1) \sup_{0\leq s\leq T_0}\|h(s)\|_{L^\infty_{y,v}}$.  For the second term, we split the time integration:
\begin{equation}\label{Jspk*-k-t-decom}
\begin{split}
 \int^{s}_{s-\delta}  +   \int^{s-\delta}_{\max\{0,  t^\prime_1\}}:=J_{sp, k*}^{k,1}+J_{sp, k*}^{k,2},
\end{split}
\end{equation}
 where $J_{sp, k*}^{k,1}$ is bounded by $\delta\sup_{0\leq s\leq T_0}\|h(s)\|_{L^\infty_{y,v}}$ due to the small-time truncation.
%\begin{equation*}
%\begin{split}
%\delta  \sup_{v} \int_{|u| \leq N} \mathbf{k}_{N} (v, u) \dd u \sup_{0\leq s\leq %T_0}\|h(s)\|_{L^\infty_{y,v}}
%\lesssim  ,
%\end{split}
%\end{equation*}
The second term $J_{sp, k*}^{k,2}$ in \eqref{Jspk*-k-t-decom} satisfies
\begin{equation}\label{Jspk*-k*}
\begin{split}
 J_{sp, k*}^{k,2}\lesssim & \int_0^{t_1} \dd s
\int_{M'_{Y_1(s)}(v')}   \dd v^{\prime}
\int_{0}^{s-\delta} \dd \tau \   e^{-\nu_0(t-\tau)}
\int_{|u|\leq N} \dd u \ \big| h\big(\tau, Y(\tau; s,Y_1(s),v^{\prime}),u\big)\big|.
\end{split}
\end{equation}
%------------------%------------------%------------------

Consider the change of variables
\begin{equation*}
 v^\prime\mapsto y:=Y(\tau;\, s, Y_1(s), v^\prime)=Y_1(s)+(\tau-s){v}^\prime.
\end{equation*}
 Since $\tau\geq 0$, the trajectory $Y(\tau;\, s, Y_1(s), v^\prime)$ does not collide with the boundary $\partial\Omega_\e$. Then, for $0 \leq \tau \leq s-\delta <s$, we compute
\begin{eqnarray*}
\big| \det \big[\nabla_{v^{\prime}} Y(\tau; s,Y_1(s),v^\prime)\big] \big|
= |s-\tau|^3  \big|  \det \big(\delta_{ij}   +   O(\e^3)\big) \big|
\geq \frac{1}{2} \delta^3.
 \end{eqnarray*}
Deducing similarly as \eqref{Jspk*-k*-decom}--\eqref{Jspk*-k*-(I-P)f} and
collecting \eqref{Jspk*-k-t-decom} and \eqref{Jspk*-k*}, we obtain
\begin{equation}\label{Jspk*-k-bd}
\begin{split}
& J_{sp,k*}^{k}  \lesssim  o(1) \sup_{0\leq s\leq T_0}\|h(s)\|_{L^\infty_{y,v}}  +  \sup_{0\leq s\leq T_0}\| \o^{-1} \mathbf{P} \bar{f}(s)\|_{L^6_{y,v}}
+  \sup_{0\leq s\leq T_0}\| \o^{-1} (\mathbf{I-P}) \bar{f}(s)\|_{L^2_{y,v}}.
\end{split}
\end{equation}

\noindent \textbf{Step 1.1.3.  Estimate of $J_{sp, k*}^{di}$.}

Similarly to \eqref{domain-split-diff}, we partition the integration domain $\{n(y_1^\prime) \cdot u^{\prime*}_{1}>0 \}= A^{\prime*}_{1}(u^{\prime*}_{1}) \cup A^{\prime*}_{2}(u^{\prime*}_{1}) \cup  M^{\prime*}_{y_1^\prime}(u^{\prime*})$, where $M^{\prime*}_{y_1^\prime}(u^{\prime*}):=\{|u^{\prime*}_{1}|\leq N,\ n(y_1^\prime) \cdot u^{\prime*}_{1}  \geq \eta\}$.
%\begin{equation*}%\label{domain-split-diff}
%\begin{split}
% &\left\{|u^{\prime*}_{1}|>N,\ n(y_1^\prime) \cdot u^{\prime*}_{1}>0 \right \} \cup
%\left\{|u^{\prime*}_{1}|\leq N,\  n(y_1^\prime) \cdot u^{\prime*}_{1}  <\eta \right\}
%\cup \left\{|u^{\prime*}_{1}|\leq N,\ n(y_1^\prime) \cdot u^{\prime*}_{1}  \geq %\eta\right\}\\
%\{n(y_1^\prime) \cdot u^{\prime*}_{1}>0 \}= A^{\prime*}_{1}(u^{\prime*}_{1}) \cup %A^{\prime*}_{2}(u^{\prime*}_{1}) \cup  M^{\prime*}_{y_1^\prime}(u^{\prime*}),
%\end{split}
%\end{equation*}
The set  $A^{\prime*}_{1}(u^{\prime*}_{1})$ and $A^{\prime*}_{2}(u^{\prime*}_{1})$ yield small contribution, similar to \eqref{cut-inpd-e}.  Thus,
\begin{equation}\label{Jspk*-di-0}
\begin{split}
J_{sp, k*}^{di} \lesssim  o(1)  \sup_{0\leq s\leq T_0}\|h(s)\|_{L^\infty_{y,v}} +  J_{sp, k*}^{di*},
\end{split}
\end{equation}
where the bulk $J_{sp, k*}^{di*}$ is given by
\begin{equation*}\label{Jsp-k*di}
\begin{split}
 J_{sp, k*}^{di*}=&\mathbf{1}_{\{t_1>0\}}    \mathbf{1}_{\{t^\prime_1>0\}}   \int_0^{t_1} \dd s \ e^{-\nu_0(t-s)}
\int_{|v^\prime|\leq N}  \dd v^{\prime}    e^{-\nu_0(s-t^\prime_1)}  \int_{M^{\prime*}_{y_1^\prime}(u^{\prime*}_{1})}
\big|h\big(y_{1}^\prime,u^{\prime*}_{1}\big)\big|\dd\sigma^{\prime*}_1.
\end{split}
\end{equation*}

Let $ (y_{1}^\prime,u^{\prime*}_{1})\in \partial \O_{\e}\times  {M^{\prime*}_{y_1^\prime}(u^{\prime*})} $ be given.
 Lemma \ref{0-bounce-claim-model}  ensures that the backward trajectory starting from $(y_{1}^\prime, u^{\prime*}_{1})$ undergoes no further collisions. Thus, $J_{sp,k}^{di*}$ propagates back to the initial plane $\{\bar{t} = 0\}$:
\begin{equation}\label{Jspk*-di*}
 J_{sp,k*}^{di*} =  J_{sp,k*}^{di*,0}   +   J_{sp,k*}^{di*,g} +  J_{sp,k*}^{di*,k},
\end{equation}
where
\begin{equation*}\label{Jspk*-di*0-k}
\begin{split}
 J_{sp,k*}^{di*,0} = \; &  \mathbf{1}_{\{t_1>0\}}    \mathbf{1}_{\{t^\prime_1>0\}}
 \int_0^{t_1} \dd s \ e^{-\nu_0(t-s)}  \int_{|v^\prime|\leq m}  \dd v^{\prime}  e^{-\nu_0(s-t^\prime_1)}
 \int_{M^{\prime*}_{y_1^\prime}(u^{\prime*}_1)}  \dd\sigma^{\prime*}_1 \\
%------------------
&  \times
e^{-\int^{t^\prime_1}_{0}\nu({u}_1^{\prime*})\dd \tau}
  h \big(0, Y(0; t_1^\prime, {y}_1^\prime, u^{\prime*}_{1}), u^{\prime*}_{1}\big),  \\
%------------------%------------------%------------------%------------------
 J_{sp,k*}^{di*,g} = \; &  \mathbf{1}_{\{t_1>0\}}    \mathbf{1}_{\{t^\prime_1>0\}}
 \int_0^{t_1} \dd s \ e^{-\nu_0(t-s)}  \int_{ |v^\prime|\leq m }  \dd v^{\prime}    e^{-\nu_0(s -t^\prime_1)}
 \int_{M^{\prime*}_{y_1^\prime}(u^{\prime*})}  \dd\sigma^{\prime*}_1  \\
%------------------
& \times  \int_0^{t^\prime_1} \dd\tau
e^{-\int^{t^\prime_1}_{\tau}\nu(u^{\prime*}_{1})\dd \tau^\prime}
 \e \big(  \o \bar{g}  \big )
\big(\tau,  Y(\tau; t_1^\prime, {y}_1^\prime, u^{\prime*}_{1} ),  u^{\prime*}_{1}  \big) ,  \\
%------------------%------------------%------------------%------------------
 J_{sp,k*}^{di*,k} = \; &   \mathbf{1}_{\{t_1>0\}}    \mathbf{1}_{\{t^\prime_1>0\}}
 \int_0^{t_1} \dd s \  e^{-\nu_0(t-s)}   \int_{|v^\prime|\leq m}  \dd v^{\prime}    e^{-\nu_0(s-t^\prime_1)}
 \int_{M^{\prime*}_{y_1^\prime}(u^{\prime*})}  \dd\sigma^{\prime*}_1 \\
%------------------
& \times \int_0^{t^\prime_1} \dd \tau   e^{ -\int^{t^\prime_1}_{\tau}  \nu(u^{\prime*}_{1})  \dd \tau^{\prime}}
\int_{\mathbb{R}^{3}} \dd \tilde{u}
\mathbf{k}_{\beta}(u^{\prime*}_{1},\tilde{u}) \frac{\o (u^{\prime*}_{1})}{\o(\tilde{u})}
  h\big(  Y(\tau; t_1^\prime, {y}_1^\prime, u^{\prime*}_{1}), \tilde{u}  \big).
\end{split}
\end{equation*}

The terms $J_{sp,k*}^{di*,0}$ and $J_{sp,k*}^{di*,g}$ are estimated similarly to \eqref{J0g-bd}. For $J_{sp,k*}^{di*,k}$, we follow the approach used for $J_{sp,k*}^{k}$: bound and decompose the kernel by $\mathbf{k}_{\tilde{\beta} }(u^{\prime*}_{1},\tilde{u})$, and split the time interval $[0,  t^{\prime}_{1}]= [0,  t^{\prime}_{1}-\delta] \cup [t^{\prime}_{1}-\delta , t^{\prime}_{1}]$.
To handle the integration on $[0,  t^{\prime}_{1}-\delta]$, consider the change of variables:
\begin{equation*}
\begin{split}
u^{\prime*}_{1} \mapsto Y(\tau; t^{\prime}_{1}, y^{\prime}_{1}, u^{\prime*}_{1})  =   y^{\prime}_{1}  +  (\tau- t^{\prime}_{1}) u^{\prime*}_{1}.
\end{split}
\end{equation*}
For $0 \leq \tau \leq t^{\prime}_{1}-\delta$, we compute
\begin{eqnarray*}
  \big| \det \big[\nabla_{u^{\prime*}_{1}} Y(\tau; t^{\prime}_{1}, y^{\prime}_{1}, u^{\prime*}_{1}) \big]\big |
 = |t^{\prime}_{1}-\tau|^3
 \geq  \delta^3.
 \end{eqnarray*}
 Following the same argument as \eqref{Jspk*-k*-decom}--\eqref{Jspk*-k-bd}, and
%\begin{equation}\label{Jspk*-di*k-bd}
%\begin{split}
%& J_{sp,k*}^{di*,k}\lesssim o(1) \sup_{0\leq s\leq T_0} \|h(s)\|_{L^\infty_{y,v}}  +   %\sup_{0\leq s\leq T_0}\|\mathbf{P}\bar{f}(s)\|_{L^6_{y,v}}
% + \sup_{0\leq s\leq T_0}\|(\mathbf{I-P})\bar{f}(s)\|_{L^2_{y,v}}.
%\end{split}
%\end{equation}
combining \eqref{Jspk*-di-0} and \eqref{Jspk*-di*}, we obtain
\begin{equation}\label{Jspk*-di-bd}
\begin{split}
 J_{sp, k*}^{di} \lesssim \; &  e^{-\frac{\nu_0}{2}\bar{t}}\|h_0\|_{L^{\infty}_{y,v}}+ o(1) \sup_{0\leq s\leq T_0}\|h(s)\|_{L^\infty_{y,v}}  +  \e \sup_{0\leq s\leq T_0}\|{\langle v\rangle}^{-1} \o\bar{g}(s)\|_{L^\infty_{y,v}}
 \\ &+  \sup_{0\leq s\leq T_0}\| \o^{-1} \mathbf{P}\bar{f}(s)\|_{L^6_{y,v}}+\sup_{0\leq s\leq T_0}\| \o^{-1} (\mathbf{I-P})\bar{f}(s)\|_{L^2_{y,v}}.
\end{split}
\end{equation}

 Finally, we collect  \eqref{Jspk-bd},  \eqref{Jspk*-duh}, \eqref{Jspk*-sp-bd},   \eqref{Jspk*-k-bd} and \eqref{Jspk*-di-bd} in Steps 1.1.1--1.1.3 to get
\begin{equation}\label{Jsp-bd}
\begin{split}
 |J_{sp}(\bar{t},y,v)|  \lesssim \;
 &  e^{-\frac{\nu_0}{2}\bar{t}}\|h_0\|_{L^{\infty}_{y,v}}
+ o(1) \sup_{0\leq s\leq T_0}\|h(s)\|_{L^\infty_{y,v}}
+ \sup_{0\leq s\leq T_0}\| \o^{-1} \mathbf{P}\bar{f}(s)\|_{L^6_{y,v}}
\\
&+ \sup_{0\leq s\leq T_0} \| \o^{-1} (\mathbf{I-P})\bar{f}(s)\|_{L^2_{y,v}}
+ \e \sup_{0\leq s\leq T_0}\| {\langle v\rangle}^{-1} \o \bar{g}(s)\|_{L^\infty_{y,v}}.
\end{split}
\end{equation}
\medskip

\noindent \textbf{ Step 1.2.  Estimate of $J_{k}(\bar{t},y,v)$.}

For $J_{k}(\bar{t},y,v)$ in \eqref{J-t}, the backward trajectory does not collide with the boundary $\pt\O_\e$. Following the same approach as in the estimation of $J_{sp,k}(\bar{t},y,v)$ in Step 1.1, we partition the integration domain
$\mathbb{R}^3= A_{1}(u)\cup A_{2}(u) \cup \hat{M}_{Y(s)}(u)$,
where $A_{1}(u)$ and $A_{2}(u)$ yield the small term  $ o(1)  \sup_{0\leq s\leq T_0}\|h(s)\|_{L^\infty_{y,v}}$ and the bulk set is defined by
\begin{equation*}%\label{Jk-bd}
\begin{split}
%A_{1}(u)&:=\left \{ u\in \mathbb{R}^3 :\; |u|\geq N  \right\},\\
%A_{2}(u)&:=\left \{ u\in \mathbb{R}^3 :\; |u|\leq  N  \text{ and  }  \left|u\cdot %\frac{\nabla_x\xi(\e Y(s))}  {|\nabla_x\xi(\e Y(s))|}\right | \leq  \eta \right\},\\
\hat{M}_{Y(s)}(u)&:=\Big \{ u\in \mathbb{R}^3 :\; |u|\leq N \text{ and  }  \Big|u\cdot \frac{\nabla_x\xi(\e Y(s))}  {|\nabla_x\xi(\e Y(s))|}\Big | \geq  \eta \Big\}.
\end{split}
\end{equation*}
 For $\hat{M}_{Y(s)}(u)$, we apply the Duhamel principle to the integrand $h(s,Y(s),u)$ in $J_{k}(t,y,v)$, obtaining an expression similar to \eqref{Jspk*-duh}. Following the same estimation procedure as in \eqref{Jspk*-duh} in Step 1.1, we finally obtain
\begin{equation}\label{Jk-bd}
\begin{split}
 |J_{k}(t,y,v)| \lesssim \;  &  e^{-\frac{\nu_0}{2}\bar{t}}\|h_0\|_{L^{\infty}_{y,v}}
+ o(1) \sup_{0\leq s\leq T_0}\|h(s)\|_{L^\infty_{y,v}}
+ \sup_{0\leq s\leq T_0}\| \o^{-1} \mathbf{P}\bar{f}(s)\|_{L^6_{y,v}}
\\
&
+ \sup_{0\leq s\leq T_0} \| \o^{-1} (\mathbf{I-P})\bar{f}(s)\|_{L^2_{y,v}}
+ \e \sup_{0\leq s\leq T_0}\| {\langle v\rangle}^{-1} \o \bar{g}(s)\|_{L^\infty_{y,v}}.
\end{split}
\end{equation}
\medskip

\noindent \textbf{ Step 1.3.  Estimate of $J_{di}(t,y,v)$.}

 Following the approach used for estimating $J_{sp, k*}^{di}$, we partition the integration domain  $\{n(y_1)\cdot v_1^*>0 \} = A^*_1(v_1^*) \cup A^*_2(v_1^*) \cup M^*_{y_1}(v_1^*)$,
where $M^*_{y_1}(v_1^*) :=\{|v^{*}_{1}|\leq N,\ n(y_1) \cdot v_1^*  \geq \eta\}$.
%\begin{equation*}%\label{domain-split-diff}
%\begin{split}
%A^*_1(v_1^*) & := \left\{v_1^*\in\mathbb{R}^3:\;|v^{*}_{1}|>N,\ n(y_1) \cdot v_1^*>0 \right %\},\\
%A^*_2(v_1^*) & :=\left\{v_1^*\in\mathbb{R}^3:\;|v^{*}_{1}|\leq N,\  0<  n(y_1) \cdot v_1^*  %<\eta \right\},\\
%M^*_{y_1}(v_1^*) & :=\{|v^{*}_{1}|\leq N,\ n(y_1) \cdot v_1^*  \geq \eta\}.
%\end{split}
%\end{equation*}
The contributions from $A^*_1(v_1^*)$ and $A^*_2(v_1^*)$ yield small term. Thus, $J_{di}(t,y,v)$ is bounded by
\begin{equation*}\label{Jdi-0}
\begin{split}
 |J_{di}(\bar{t},y,v)|\lesssim &\; o(1)  \sup_{0\leq s\leq T_0}\|h(s)\|_{L^\infty_{y,v}}+ J_{di*}(\bar{t},y,v),
\end{split}
\end{equation*}
where
\begin{equation*}\label{Jdi-0}
\begin{split}
 J_{di*}(\bar{t},y,v):= \mathbf{1}_{\{t_1>0\}} \ e^{-\nu_0(t-t_1)}
\int_{M^*_{y_1}(v_1^*)}
\big |h(t_1,y_1,v_1^*)\big|\dd v_1^*.
\end{split}
\end{equation*}
For $(y_1,v_1^*) \in  \partial \O_{\e} \times {M^*_{y_1}(v_1^*)}$, Lemma \ref{0-bounce-claim-model} ensures that the backward trajectory starting from $(y_1,v_1^*)\in M_{y_1}(v_1^*)$ propagates back to the initial plane $\{\bar{t}=0\}$. Thus,
\begin{equation*}%\label{Jdi*}
\begin{split}
&J_{di*}(\bar{t},y,v) = J_{di*}^{0} +  J_{di*}^{g} +  J_{di*}^{k},
\end{split}
\end{equation*}
where
\begin{equation*}%\label{Jdi*}
\begin{split}
 &J_{di*}^{0}=   \mathbf{1}_{\{t_1>0\}} e^{-\nu_0(t-t_1)}
  \int_{M^*_{y_1}(v_1^*)} \dd v_1^*  e^{-\int^{t_1}_{0} \nu(v_1^*)  \dd \tau }
 h \big(0,Y(0; t_1,y_1,v_1^*), v_1^* \big),  \\
 %----------------------%----------------------%----------------------
&J_{di*}^{g}=   \mathbf{1}_{\{t_1>0\}} e^{-\nu_0(t-t_1)}
\int_{M^*_{y_1}(v_1^*)}\dd v_1^*
\int_0^{t_1} \dd s
e^{-\int^{t_1}_{0}  \nu(v_1^*)  \dd \tau }
   \e \big( \o \bar{g} \big)  \big(  Y(s;t_1,y_1,v_1), v_1 \big),  \\
%----------------------%----------------------%----------------------
&J_{di*}^{k}=  \mathbf{1}_{\{t_1>0\}}  e^{-\nu_0(t-t_1)}
\int_{M^*_{y_1}(v_1^*)} \dd v_1^* e^{-\int^{t_1}_{0}  \nu(v_1^*)  \dd \tau }
\int_{\mathbb{R}^{3}} \dd v^{\prime}  \
\mathbf{k}_{\beta}(v_1^*,v^{\prime}) \frac{\o(v_1^*)}{\o (v^{\prime})}
  \big|  h(Y(s;t_1,y_1,v_1^*),v^{\prime})   \big|.
\end{split}
\end{equation*}
The terms $J_{di*}^{0}$ and $J_{di*}^{g}$ are estimated similarly to \eqref{J0g-bd}.
The term $J_{di*}^{k}$ is estimated by the change of variable $v_1^*\mapsto Y(s;t_1,y_1,v_1^*)$, similar to the approach used for $J_{sp,k*}^{di*,k}$ in Step 1.1.3. We conclude
\begin{equation}\label{Jdi-bd}
\begin{split}
|J_{di}(t,y,v)| \lesssim \; &  e^{-\frac{\nu_0}{2}\bar{t}}\|h_0\|_{L^{\infty}_{y,v}}+ o(1) \sup_{0\leq s\leq T_0}\|h(s)\|_{L^\infty_{y,v}}
+ \sup_{0\leq s\leq T_0}\| \o^{-1} \mathbf{P}\bar{f}(s)\|_{L^6_{y,v}}
\\
&+ \sup_{0\leq s\leq T_0}\| \o^{-1} (\mathbf{I-P})\bar{f}(s)\|_{L^2_{y,v}}
+ \e \sup_{0\leq s\leq T_0}\|{\langle v\rangle}^{-1} \o\bar{g}(s)\|_{L^\infty_{y,v}} .
\end{split}
\end{equation}

Finally, combining  \eqref{h(t,y,v){1}}, \eqref{J0g-bd},  \eqref{Jsp-bd}, \eqref{Jk-bd} and \eqref{Jdi-bd} in Steps 1.1--1.3, we obtain
the following estimate for $h(\bar{t}, y,v)$ restricted on $ \mathcal{M}{(y,v)}$:
\begin{equation}\label{h(t,y,v){1}-bd}
\begin{split}
 \left \|h(\bar{t}, y,v)\mathbf{1}_{ \mathcal{M}{(y,v)}} \right \|_{L^\infty_{y,v}}   \lesssim \; &
 e^{-\frac{\nu_0}{2}\bar{t}}\left \|h_0 \right \|_{L^{\infty}_{y,v}}
 +  o(1) \sup_{0\leq s\leq T_0}\left \|h(s) \right \|_{L^\infty_{y,v}} + \sup_{0\leq s\leq T_0}\left\| \o^{-1} \mathbf{P}\bar{f}(s) \right \|_{L^6_{y,v}}\\
&+ \sup_{0\leq s\leq T_0}\left \| \o^{-1} (\mathbf{I-P})\bar{f}(s) \right \|_{L^2_{y,v}}
+ \e \sup_{0\leq s\leq T_0}\left \|{\langle v\rangle}^{-1} w\bar{g}(s) \right \|_{L^\infty_{y,v}}.
\end{split}
\end{equation}
\medskip

\noindent \textbf{Step 2.  Estimate of $h(\bar{t}, y,v)$.}

Apply the semigroup representation from Lemma \ref{semigroup-estimate} and the Duhamel principle to \eqref{eq-unst-stretch-t-h-with-K}:
\begin{equation}\label{h(t,y,v)}
\begin{split}
 h(\bar{t}, y,v) =  \;
 &  G(\bar{t})h_0(y,v)    +
   \int^{\bar{t}}_0 G(\bar{t}-s)  \left[ \e \o  \bar{g} \left(s, Y(s), v \right)\right]  \dd s \\
 & + \int^{\bar{t}}_0 G(\bar{t}-s)
 \Big[\int_{\mathbb{R}^3} \mathbf{k}_{\beta}(v,u) \frac{ \o (v)}{ \o (u)}
  h \big(s, Y(s), u \big)\dd u \Big] \dd s.
\end{split}
\end{equation}
 Applying the semigroup estimate
\eqref{semigroup-esti} from Lemma \ref{semigroup-estimate}, we derive:
\begin{equation}\label{h(t,y,v)-semi-esti}
\begin{split}
 \|h(\bar{t})\|_{L^{\infty}_{y,v}} \leq  \; &  (2C_*+1)e^{-\frac{\nu_0}{2}\bar{t}} \|h_0\|_{L^{\infty}_{y,v}}
  +   (2C_*+1)\int^{\bar{t}}_0 e^{-\frac{\nu_0}{2}(\bar{t}-s)} \| \e w\bar{g} (s)\|_{L^{\infty}_{y,v}}\dd s \\
 & + (2C_*+1)\int^{\bar{t}}_0 e^{-\frac{\nu_0}{2}(\bar{t}-s)}
 \sup_{y,v}\Big|\underbrace{\int_{\mathbb{R}^3} \mathbf{k}_{\tilde{\beta}} \big( v ,u \big) \big| h \big(s, Y(s), u \big)\big|\dd u }_{:=I(s;\bar{t},y,v)} \Big| \dd s,
\end{split}
\end{equation}
where we have used the kernel bound \eqref{k-tilde-beta}.

We decompose the integral  $I(s;\bar{t},y,v)$ into two parts:
\begin{equation*}%\label{h(t,y,v)-I}
\begin{split}
 I(s;\bar{t},y,v)=\;&\int_{\mathbb{R}^3} \mathbf{k}_{\tilde{\beta}} \big( v ,u \big)  h \big(s, Y(s), u \big) \left\{1-{\bf 1}_{\mathcal{M}(Y(s),u)} \right\}\dd u\\
 & +\int_{\mathbb{R}^3} \mathbf{k}_{\tilde{\beta}} \big( v ,u \big) h \big(s, Y(s), u \big) {\bf 1}_{\mathcal{M}(Y(s),u)} \dd u\\
 :=\; &I_1(s;\bar{t},y,v) + I_2(s;\bar{t},y,v).
\end{split}
\end{equation*}
By the definition of $\mathcal{M}(Y(s),u)$ in \eqref{main-set-mathcal-M}, the first term $I_1(s;\bar{t},y,v)$ is bounded by
\begin{equation}\label{h(t,y,v)-I1-bd}
\begin{split}
\|I_1(s;\bar{t},y,v)\|_{L^\infty_{y,v}}
\lesssim & \Big(
\int_{|u|>N} |\mathbf{k}_{\tilde{\beta}} \big( v ,u \big) | \dd u
+\int_{|u\cdot \frac{\nabla_x\xi(\e Y(s))}  {|\nabla_x\xi(\e Y(s))|}| < \eta} |\mathbf{k}_{\tilde{\beta}} \big( v ,u \big) | \dd u
\Big) \| h (s)\|_{L^{\infty}_{y,v}}\\
\lesssim &o(1)\sup_{0\leq s\leq T_0} \| h (s)\|_{L^{\infty}_{y,v}},
\end{split}
\end{equation}
where we have used the compact support approximation of $\mathbf{k}_{\tilde{\beta}}$ by $\mathbf{k}_{N}$ as  in \eqref{k-split}.
 For the second term $I_2(s;\bar{t},y,v)$, we apply the estimate \eqref{h(t,y,v){1}-bd}  to obtain
 \begin{equation}\label{h(t,y,v)-I2-bd}
\begin{split}
 \|I_2(s;\bar{t},y,v)\|_{L^\infty_{y,v}}  \lesssim \; &
 e^{-\frac{\nu_0}{2}s}\|h_0\|_{L^{\infty}_{y,v}}
 +  o(1) \sup_{0\leq s\leq T_0}\|h(s)\|_{L^\infty_{y,v}}+ \sup_{0\leq s\leq T_0}\| \o^{-1} \mathbf{P}\bar{f}(s)\|_{L^6_{y,v}} \\
&
+ \sup_{0\leq s\leq T_0}\| \o^{-1} (\mathbf{I-P})\bar{f}(s)\|_{L^2_{y,v}}+ \e \sup_{0\leq s\leq T_0}\|{\langle v\rangle}^{-1} \o \bar{g}(s)\|_{L^\infty_{y,v}}.
\end{split}
\end{equation}
Substituting \eqref{h(t,y,v)-I1-bd} and \eqref{h(t,y,v)-I2-bd} into \eqref{h(t,y,v)-semi-esti}, we thus prove the claim \eqref{|h(bar-t,y,v)|-L6}.

 The estimate \eqref{|h(bar-t,y,v)|-L2} can be proved  in a similar way to \eqref{|h(bar-t,y,v)|-L6}. The main difference lies in the change of variables analogous to  \eqref{Jspk*-k*-decom}--\eqref{Jspk*-k*-(I-P)f}: here we directly use the norm $\|\bar{f}(\bar{s})\|_{L^2_{y,v}
 (\Omega_\e\times \mathbb{R}^3)}$, without splitting $\bar{f}$ into $\mathbf{P}\bar{f}$ and  $\ip \bar{f}$.
We skip the details for brevity. This completes the proof.
\end{proof}
\bigskip

%%%%%%%%%%%%%%%%%%%%%%%%%%%%%%
%%%%%%%%%%%%%%%%%%%%%%%%%%%%%%
%%%%%%%%%%%%%%%%%%%%%%%%%%%%%%
\section{Strong Limit for the Case $\e \lesssim \a \leq 1$}
\medskip

This section studies the perturbation equation \eqref{f-eq} and presents the proof of Theorem~\ref{main-th-1}. The proof relies on  Propositions \ref{0-macro-L2-L6-estimate} and \ref{lemma-fbar-infty-unst-0}, which are established first.

%%%%%%%%%%%%%%%%%%%%%%%%%%%%%%
%%%%%%%%%%%%%%%%%%%%%%%%%%%%%%
\subsection{Energy Estimate}\
\medskip

In this subsection, we derive the basic energy estimates for the fluctuation $f$ and its time derivative $\pt_t f$.

\begin{proposition} \label{f-ft-Energy-estimate} Let $f \in L^2([0,T]\times \Omega \times \R^3)$ be a solution of \eqref{f-eq} with $0<T\leq \infty$. Then the following energy estimates hold for all $t\in [0,T]$:
\begin{align}
& \normm{f(t)}_{L^2_{x,v}}^{2}  + \frac{1}{\e^2}\int_{0}^t\normm{\ip  f}_{L^2_{x,v}(\nu)}^2 \dd s
+  \int_{0}^t  \big(     \norm{\mathscr{P}_{\g} f}^2_{L^2_{\gamma_+}}
+   \frac{\a}{\e} \norm{(1-\mathscr{P}_{\g}) f}^2_{L^2_{\gamma_+}}   \big) \dd s \nonumber \\
%-----------
\lesssim & \normm{f_0 }_{L^2_{x,v}}^{2}  +  \e \int_{0}^t \| \nu^{-\frac{1}{2}}  \Gamma(f,f) \|_{L^2_{x,v}}^2 \dd s +  \eta \int_{0}^t\normm{\P  f}_{L^2_{x,v}}^2 \dd s,\label{f-energy-estimate}\\
%--------------------------
& \normm{\pt_t f(t)}_{L^2_{x,v}}^{2}
+ \frac{1}{\e^2}\int_{0}^t\normm{\ip \pt_t f}_{L^2_{x,v}(\nu)}^2 \dd s
+  \int_{0}^t  \big( \norm{\mathscr{P}_{\g} \pt_tf}^2_{L^2_{\gamma_+}} +   \frac{\a}{\e} \norm{(1-\mathscr{P}_{\g}) \pt_tf}^2_{L^2_{\gamma_+}} \big) \dd s  \nonumber \\
%-----------
\lesssim & \normm{\pt_t f_0 }_{L^2_{x,v}}^{2}   +  \e  \int_{0}^t  \big[\| \nu^{-\frac{1}{2}}   \Gamma(\pt_tf,f) \|_{L^2_{x,v}}^2 + \| \nu^{-\frac{1}{2}}   \Gamma(f, \pt_tf)  \|_{L^2_{x,v}}^2\big]\dd s
+  \eta \int_{0}^t\normm{\P  \pt_t f}_{L^2_{x,v}}^2 \dd s, \label{ft-energy-estimate}
\end{align}
where $0<\eta \ll \min\{ 1, \frac{\l}{4}\}$  is a sufficiently small constant with $\l$ defined in \eqref{a-e-limit-infty}.
\end{proposition}
\medskip

\begin{proof}[\textbf{Proof}.] \  Standard $L^2$ energy estimate for \eqref{f-eq} gives
\begin{align*}
\e \frac{1}{2} \pt_{t} \normm{f(t)}_{L^2_{x,v}}^2 + \frac{1}{2} \iint_{\pt\Omega \times \R^3} f^2 [n \cdot v] \dd v \dd S_x  + \frac{1}{\e} \iint_{\Omega \times \R^3} f L f \dd v \dd x = \iint_{\Omega \times \R^3} \Gamma(f,f)f \dd v \dd x.
\end{align*}
Using the Maxwell boundary condition and the change of variables  $R_x v \mapsto v$, we obtain
\begin{align*}
\iint_{\pt\Omega \times \R^3}  f^2 [n \cdot v] \dd v \dd S_x
%=&\iint_{\g_{+}} f^2 [n \cdot v] \dd v\dd S_x + \iint_{\g_{-}} \big[(1- \a) \mathscr{R} f %+ \a \mathscr{P}_{\g} f \big]^{2} [n \cdot v] \dd S_x \dd v\\
%---------------------
=&\iint_{\g_{+}} f^2 \dd \g -  \iint_{\g_{+}} \big[(1- \a)(1-\mathscr{P}_{\g}) f + \mathscr{P}_{\g} f \big]^2 \dd \g \\
%---------------------
=&  \iint_{\g_{+}} f^2 \dd \g -\iint_{\g_{+}} \big[(1- \a)^2|(1-\mathscr{P}_{\g})f|^2 + |\mathscr{P}_{\g} f|^2   \big] \dd \g\\
=&\; \a(2-\a)\iint_{\g_+} |(1-\mathscr{P}_{\g})f|^2 \dd \g,
\end{align*}
where we have used the orthogonal decomposition
\begin{equation}\label{P-gamma-orthogonal}
 f=(1-\mathscr{P}_{\g})f+ \mathscr{P}_{\g}f \text { on } L^2_{\g_{+}}, \qquad \norm{ f}^2_{L^2_{\gamma_+}} =\norm{\mathscr{P}_{\g} f}^2_{L^2_{\gamma_+}} + \norm{(1-\mathscr{P}_{\g}) f}^2_{L^2_{\gamma_+}}.
\end{equation}
By H\"{o}lder's inequality and the coercivity of $L$, we derive
\begin{equation}\label{f-energy-estimate-1}
\begin{split}
& \normm{f(t)}_{L^2_{x,v}}^{2}  + \frac{1}{\e^2}\int_{0}^t\normm{\ip  f}_{L^2_{x,v}(\nu)}^2 \dd s
+ \frac{\a}{\e} \int_{0}^t  \norm{(1-\mathscr{P}_{\g}) f}^2_{L^2_{\gamma_+}} \dd s  \\
%-----------
\lesssim & \normm{f_0 }_{L^2_{x,v}}^{2}  +  \int_{0}^t \| \nu^{-\frac{1}{2}}  \Gamma(f,f) \|_{L^2_{x,v}}^2 \dd s.
\end{split}
\end{equation}

Define the non-grazing set
\begin{equation}\label{non-grazing-set}
\g_{\pm}^{\d} := \Big\{(x,v) \in \g_{\pm} : \norm{n(x)\cdot v} > \d,\; \d \le \norm{v} \le \frac{1}{\d} \Big\}.
\end{equation}
Note that $\iint_{\g_{+}\backslash\g_{+}^\d} \mu \dd \g\lesssim o(\d)$,
which implies
\begin{align*}
&\int_{\g_{+}\backslash\g_{+}^\d} \norm{\mathscr{P}_{\g} f }^2 \dd \g
\lesssim o(\d)
\iint_{\g_{+}} \norm{\mathscr{P}_{\g} f }^2 \dd \g.
\end{align*}
Applying the trace lemma (cf. Lemma 3.2 in \cite{Esposito2017}) to the non-grazing part, we obtain
\begin{equation} \label{f - Bundary ukai lemma}
\begin{split}
&\int_{0}^t\iint_{\g_{+}} \norm{\mathscr{P}_{\g} f }^2 \dd \g\dd s
%= &\Big\{  \iint_{\g_{+}^\d} + \iint_{\g_{+}\backslash\g_{+}^\d}\Big \}\norm{\mathscr{P}_{\g} f }^2 %\dd \g\\
\lesssim  \; \d \int_{0}^t\iint_{\g_{+}} \norm{(1-\mathscr{P}_{\g}) f }^2 \dd \g\dd s + \int_{0}^t\iint_{\g_{+}^\d} \norm{f}^2 \dd \g\dd s
 \\
%-------------------
\lesssim &  \int_{0}^{t}\iint_{\g_{+}} \norm{(1-\mathscr{P}_{\g}) f}^2 \dd \g \dd s
+ \e \iint_{\Omega \times \R^3} \norm{f_0}^2 \dd v \dd x +  \int_{0}^t\iint_{\Omega \times \R^3} \norm{f(s)}^2 \dd v \dd x \dd s
\\
&+ \int_{0}^t\iint_{\Omega \times \R^3} \norm{(\Gamma(f,f) -\e^{-1} L f)f} \dd v \dd x \dd s\\
%-------------------
\lesssim & \int_{0}^{t}\norm{(1-\mathscr{P}_{\g}) f}_{L^2_{\g_{+}}}^2  \dd s
+ \e \normm{f_0}_{L^2_{x,v}}^2 +   \int_{0}^t \normm{f}_{L^2_{x,v}}^2  \dd s
+   \frac{1}{\e} \int_{0}^t \normm{\ip f}_{L^2_{x,v}(\nu)}^2   \dd s\\
&  + \e \int_{0}^t  \normm{\nu^{-\frac{1}{2}}\Gamma(f,f)}_{L^2_{x,v}}^2 \dd s.
\end{split}
\end{equation}
Multiplying \eqref{f - Bundary ukai lemma} by a sufficiently small constant $0<\eta \ll \min\{ 1, \frac{\l}{4}\}$  and adding to \eqref{f-energy-estimate-1}, we obtain \eqref{f-energy-estimate}.

The estimate \eqref{ft-energy-estimate} follows by applying the same procedure
to the equation for $\pt_t f$. We omit the details for brevity.
\end{proof}
\bigskip

%%%%%%%%%%%%%%%%%%%%%%%%%%%%%%%%%%
%%%%%%%%%%%%%%%%%%%%%%%%%%%%%%%%%%

\subsection{Macroscopic $L^2$ and $L^6$ Estimates}\
\medskip

This subsection establishes macroscopic $L^2$ and $L^6$ estimates  and proves Proposition \ref{0-macro-L2-L6-estimate}.

 By virtue of \eqref{0-f-a-conservation-law}, the coefficient $a(t,x)$ of $\P f$ satisfies the zero-mean condition
\begin{equation} \label{a-consercation-law}
\begin{split}
\int_{\Omega} a(t,x) \dd x = 0, \;\;  \forall  t\in [0, T].
\end{split}
\end{equation}
Note that $b$ and $c$ do not satisfy the zero-mean condition due to the lack of conservation laws of angular momentum and energy for $f$. Define the Burnett functions:
\begin{equation}\label{0-Burnette-function}
\begin{split}
&A_{ij}(v) := \Big(v_{i}v_{j} - \frac{\delta_{ij}}{3}\norm{v}^2 \Big)\sqrt{\mu},
\qquad
B_{i}(v) := v_{i}\frac{\norm{v}^2-5}{\sqrt{10}} \sqrt{\mu},\quad  i,j=1,2,3.
\end{split}
\end{equation}
For each $i,j=1,2,3$,  $A_{ij}(v)$ and $B_{i}(v)$ are orthogonal to every basis element $\chi_{k} $ of $\ker L$:
\begin{equation}\label{Burnette-orthogonal}
\int_{\R^3} \chi_{k}(v) A_{ij}(v) \dd v=0, \quad
\int_{\R^3} \chi_{k}(v) B_{i}(v) \dd v=0, \quad k=0,\cdots, 4.
\end{equation}
\medskip

\begin{proof}[\textbf{Proof of Proposition \ref{0-macro-L2-L6-estimate}}] \
The proof relies on the test function method \cite{esposito2013non, Esposito2017} and elliptic theory.

Multiplying \eqref{f-eq} by a test function $\psi_{p,q}$, we obtain
\begin{equation}\label{0-test-equation-uniform-form}
\begin{split}
&\underbrace{\e \iint_{\Omega \times \R^3}  \psi_{p,q} \pt_t f \dd v \dd x}_{:=\Xi_{p,q}^{1}}
+ \underbrace{\iint_{\g_{+}} \psi_{p,q} f\dd \g
 - \iint_{\g_{-}} \psi_{p,q} f \dd \g}_{:=\Xi_{p,q}^{2}}
%----------------
\underbrace{- \iint_{\Omega \times \R^3}  (v \cdot \nabla_x \psi_{p,q}) f \dd v \dd x}_{:=\Xi_{p,q}^{3}}\\
%----------------
= & \underbrace{\iint_{\Omega \times \R^3} \Big[\e^{-1}\psi_{p,q}Lf + \psi_{p,q}\G(f,f)\Big] \dd v \dd x}_{:=\Xi_{p,q}^{4}}.
\end{split}
\end{equation}
Here the temporary index $p\in \{a,b,c\}$ marks estimates of $a, b$ and $c$, and $q\in \{2,6\}$ indicates the norms $\|\cdot \|_{L^2_{x,v}}$ and $\| \cdot \|_{L^6_{x,v}}$.

To estimate $\P f$, by the representation \eqref{Pf-abc-hat}, it suffices to estimate $a, b$ and $c$.
\medskip

\noindent\textbf{Step 1. Estimates for $a$.}

\noindent\textbf{Step 1.1.  Estimates for $\int_{s}^{t}\|a\|_{L^{2}_{x}}\dd \tau$ and $\|a\|_{L^{6}_{x}}$.}

 In \eqref{0-test-equation-uniform-form}, we consider the test function
  $$
  \psi_{a,q}(t,x,v) := \sum_{i=1}^{3} \pt_{i} \varphi_{a,q}(t,x) \big[\sqrt{10} B_{i}(v) - 5\chi_{i}(v)\big],\quad q\in \{2,6\},
  $$
  where $\varphi_{a,2}(x)$  and $\varphi_{a,6}(x)$ are solutions to the elliptic equations
\begin{align}
- \Delta_x \varphi_{a,2} = a \;\text{ in }\O, \quad& {\pt_n}\varphi_{a,2}= 0  \;\text{ on }\pt\O,  \quad \int_{\Omega} \varphi_{a,2} \dd x =0, \label{0-a-2-elliptic-equation} \\
%---------------
- \Delta \varphi_{a,6} = a^{5} - \frac{1}{\norm{\Omega}}\int_{\Omega}a^{5} \dd x, \;\text{ in }\O, \quad& {\pt_{n}}\varphi_{a,6}= 0 \;\text{ on }\pt\O,  \quad \int_{\Omega} \varphi_{a,6} \dd x =0, \label{0-a-6-elliptic-equation}
\end{align}
respectively. By \eqref{a-consercation-law}, Lemma \ref{Poisson-equation-theory} guarantees unique solutions $\varphi_{a,2}$ and $\varphi_{a,6}$ satisfying
\begin{align}
\normm{\nabla^2 \varphi_{a,2}}_{L^2_x} + \normm{\nabla \varphi_{a,2}}_{L^2_x} + \normm{\varphi_{a,2}}_{L^2_x} &\lesssim \normm{a}_{L^2_x}, \label{0-a-2-elliptic-estimate}\\
%---------------
\normm{\nabla^2 \varphi_{a,6}}_{L^{\frac{6}{5}}_x} + \normm{\nabla \varphi_{a,6}}_{L^2_x} + \normm{\varphi_{a,6}}_{L^6_x}
&\lesssim \normm{a^{5}}_{L^{\frac{6}{5}}_x}=\normm{a}_{L^{6}_x}^5.
\label{0-a-6-elliptic-estimate}
\end{align}
%Moveover, ${\varphi}_{a,q}$ satisfies Poincar\'{e}'s inequality.
%The Poincar\'{e}'s inequality gives
%\begin{equation}\label{0-Poincare-inequality-varphi-a}
%\begin{split}
%\normm{\varphi_{a,2}}_{L^2_x} \lesssim \normm{\nabla_x  \varphi_{a,2}}_{L^2_x}.
%\end{split}\end{equation}

We now estimate each term in \eqref{0-test-equation-uniform-form}.
For $ \Xi_{a,2}^{1}$, integration by parts yields
\begin{equation*}
\begin{split}
\int_{s}^{t}   \Xi_{a,2}^{1}\dd \tau = & \; \e\big[ G_{a}(t)-G_{a}(s) \big] -\e\int_{s}^{t}\iint_{\Omega \times \R^{3}} \pt_t\big(\psi_{a,2}\big) f.
\end{split}
\end{equation*}
By H\"{o}lder's inequality and \eqref{0-a-2-elliptic-equation}, $G_{a} (t)$ is bounded by $\normm{f(t)}_{L^2_{x,v}}^2$. Moreover,
\begin{equation}\label{Theta-a2-01-estimate 1}
\begin{split}
 \iint_{\Omega \times \R^{3}} \pt_t\big(\psi_{a,2}\big) f
%---------------
%=&\sum_{i=1}^{3} \iint_{\Omega \times \R^{3}}  \pt_t \pt_{i} \varphi_{a,2} (\sqrt{10} %B_{i} - 5\chi_{i}) f\\
%---------------
%=& \sum_{i=1}^{3}\int_{\Omega} \pt_t \pt_{i} \varphi_{a,2} \int_{\R^{3}}   v_i(|v|^2 - %10)\sqrt{\mu} \big[\P f+\ip f \big] \\
%---------------
%=& \sum_{i=1}^{3}\int_{\Omega} \pt_t \pt_{i} \varphi_{a,2} \int_{\R^{3}}   v_i(|v|^2 - %10)\sqrt{\mu} \big(a \chi_{0}+ \sum_{i=1}^{3} b_{i} \chi_{i}+ c \chi_{4}+\ip f \big) \\
%---------------
\lesssim &\normm{\pt_t \nabla_x \varphi_{a,2}}_{L^2_x} \big( \normm{b}_{L^2_x} + \normm{\ip f}_{L^2_{x,v}}\big),
\end{split}
\end{equation}
since contribution from $a$ and $c$ vanish due to the fact
\begin{equation*}
\begin{split}
&  \int_{\R^{3}}  \big[\sqrt{10} B_{i}(v) - 5\chi_{i}(v)\big]  \chi_{j}(v)\dd v   =
\int_{\R^{3}}   v_i(|v|^2 - 10)\sqrt{\mu}  \chi_{j}(v)\dd v
%=   \int_{\R^{3}}   v_i v_j (|v|^2 - 10)\tilde{\mu}\dd v
=     -5\delta_{ij}
\end{split}
\end{equation*}
for $i=1,2,3$ and $j=0,\cdots,4$.   Thus, we obtain
\begin{equation} \label{0-Theta1 - a estimate}
\begin{split}
\int_{s}^{t}   \norm{ \Xi_{a,2}^{1}  }
 \leq  &  \e  \left[G_{a} (t)-  G_{a} (s)\right] + \e  \int_{s}^{t}\normm{\pt_t \nabla_x \varphi_{a,2}}_{L^2_x}  \big( \normm{b}_{L^2_x} + \normm{\ip f}_{L^2_{x,v}}\big).
\end{split}
\end{equation}

For $ \Xi_{a,6}^{1} $, using \eqref{0-a-6-elliptic-estimate}, we have
\begin{equation}
\begin{split} \label{0-Theta1 - a L6 estimate}
\norm{ \Xi_{a,6}^{1}  }
%---------------
=&  \e \Big| \sum_{i=1}^{3}  \iint_{\Omega\times \R^3 }  \pt_{i} \varphi_{a,6}   v_i(|v|^2 - 10)\sqrt{\mu} \pt_t f \Big|
%---------------
\lesssim  { \e  \normm{\pt_t f }_{L^2_{x,v}} }
 \normm{a}_{L^6_{x}}^{5}.
%---------------
\end{split}
\end{equation}

For ${\Xi}_{a,q}^{2}$ ($q\in \{2,6\}$), the condition ${\pt_{n}}\varphi_{a,q}|_{\pt\O}= 0$ implies $\mathscr{R}(\psi_{a,q} )= \psi_{a,q}$. Thus, by the Maxwell boundary condition and the change of variables $R_x v \mapsto v$, we obtain
\begin{equation}\label{0-Theta2-boundary-calculation-a}
\begin{split}
\Xi_{a,q}^{2}
%=& \iint_{\g_{+}} \psi_{a,q}f \dd \g
%- \iint_{\g_{-}} \psi_{a,q} \big( (1-\a)  \mathscr{R}f
%+\a  \mathscr{P}
%_{\g} f \big)\dd \g  \\
%---------------------
=& \iint_{\g_{+}} \psi_{a,q}f \dd \g
-\iint_{\g_{+}} \mathscr{R}\big ( \psi_{a,q} \big)\big( (1-\a) f
+\a  \mathscr{P}_{\g}f\big)\dd \g
%-----------------------
%=& \iint_{\g_{+}} \psi_{a,q}f \dd \g
%-\iint_{\g_{+}}  \psi_{a,q}   \big( (1-\a) f
%+\a  \mathscr{P}_{\g}f\big)\dd \g  \\
%-----------------------
= \a \iint_{\g_{+}}  \psi_{a,q}(1-\mathscr{P}_{\g})f  \dd \g,
\end{split}
\end{equation}
where we used \eqref{P-gamma-orthogonal}. For $\Xi_{a,2}^{2}$, applying the trace theorem and \eqref{0-a-2-elliptic-estimate} gives
\begin{equation}\label{0-Theta2 - L2 estimate}
\begin{split}
\norm{\Xi_{a,2}^{2} } \lesssim  \a  \norm{(1-\mathscr{P}_{\g})f}_{L^2_{\g_{+}}} \norm{\varphi_{a,2}}_{L^2(\pt\O)}
\lesssim & \a  \norm{(1-\mathscr{P}_{\g})f}_{L^2_{\g_{+}}} \normm{a}_{L^2_x}.
%-----------------
%\norm{{\Xi}_{a,2t}^{1,2} } \lesssim  & \a  \norm{\mathscr{P}_{\perp} %\pt_tf}_{L^2_{\g_{+}}} \normm{\pt_t {a}}_{L^2_x}.
\end{split}
\end{equation}
For $\Xi_{a,6}^{2}$,  using \eqref{0-a-6-elliptic-estimate} and interpolation, we derive
\begin{equation}\label{0-Theta2 - L2 estimate-L6-a}
\begin{split}
\norm{\Xi_{a,6}^{2}}
\lesssim & \a  \norm{(1-\mathscr{P}_{\g})f}_{L^4_{\g_{+}}}  \norm{\nabla_x  \varphi_{a,6}}_{L^{\frac{4}{3}}(\pt\O)}
 %--------------------
%  \lesssim & \a   \big( \e^{-\frac{1}{2}}\norm{ %(1-\mathscr{P}_{\g})f}_{L^2_{\g_{+}}}\big)^{\frac{1}{2}} \big( %\e^{\frac{1}{2}}\norm{(1-\mathscr{P}_{\g})f}_{L^\infty_{x,v}} \big )^{\frac{1}{2}}  %\normm{\varphi_{a,6}}_{W^{2,\frac{6}{5}}_x(\O)}\\
 %--------------------
 \lesssim  \a  \norm{(1-\mathscr{P}_{\g})  f}_{L^2_{\g_{+}}}^{\frac{1}{2}}
\normm{\omega f}_{L^{\infty}_{x,v}}^{\frac{1}{2}} \normm{ a}_{L^6_x}^5.
\end{split}
\end{equation}
where we used the Soboelv embedding
$\normm{\phi}_{L^{\frac{4}{3}}(\pt \Omega)} \lesssim \normm{\phi}_{W^{1,\frac{6}{5}}(\Omega)}$ (cf. \cite{leoni2024first}).

For $\Xi_{a,q}^{3}$ ($q\in \{2,6\}$), direct computation gives
\begin{equation}\label{0-Theta3 - a q estimate}
\begin{split}
\Xi_{a,q}^{3}%=&-\iint_{\Omega \times \R^3} v\cdot \nabla_x  \psi_{a,2} f \dd v \dd x \\
%---------------
=& -\sum_{i,j=1}^{3}\int_{\Omega} \pt_i\pt_j \varphi_{a,q}\int_{\R^3} v_iv_j(|v|^2-10) \sqrt{\mu} \big[\P f+\ip f \big]
%---------------
=  \int_{\Omega} 5\Delta_x \varphi_{a,q} a   + E_{a,q},
\end{split}
\end{equation}
where the remainder $E_{a,q}$ arises from the $\ip f$ contribution, and we have used
\begin{equation}\label{v_iv_j(|v|^2 - 10)-orthogonal}
\begin{split}
&\displaystyle\int_{\R^{3}}   v_iv_j(|v|^2 - 10)\sqrt{\mu}  \chi_{k}(v)
=    0,   \quad
\int_{\R^{3}}   v_iv_j(|v|^2 - 10)\sqrt{\mu}  \chi_{0}(v)
=     -5\delta_{ij}
\end{split}
\end{equation}
for $i,j=1,2,3$ and $k=1,2,3,4$.
Using \eqref{0-a-2-elliptic-equation} and  \eqref{0-a-6-elliptic-equation}, we have
\begin{align}
&\Xi_{a,2}^{3}%=&-\iint_{\Omega \times \R^3} v\cdot \nabla_x  \psi_{a,2} f \dd v \dd x \\
%---------------
=  \int_{\Omega} 5\Delta_x \varphi_{a,2} a \dd x   + E_{a,2}
= -5 \normm{ a }^2_{L^{2}_{x}}  + E_{a,2}, \label{0-Theta3 - a estimate}\\
%---------------
&\Xi_{a,6}^{3}
%---------------
=  \int_{\Omega} 5\Delta_x \varphi_{a,6} a\dd x
 + E_{a,6}
= -5 \normm{ a }^6_{L^{6}_{x}}  + E_{a,6}. \label{0-Theta3 - a L6 ettimate}
\end{align}
 The remainders $E_{a,2}$ and $E_{a,6}$ are controlled via \eqref{0-a-2-elliptic-estimate} and \eqref{0-a-6-elliptic-estimate}:
\begin{equation}\label{0-E-a-2-estimate}
\begin{split}
\norm{E_{a,2}} \lesssim \normm{a}_{L^{2}_{x}} \normm{\ip f}_{L^{2}_{x,v}}, \quad
\norm{E_{a,6}} \lesssim  \normm{a}_{L^{6}_{x}}^{5} \normm{\ip f}_{L^{6}_{x,v}}.
\end{split}
\end{equation}

For $\Xi_{a,q}^{4}$, by H\"{o}lder's inequality and \eqref{0-a-2-elliptic-estimate} and \eqref{0-a-6-elliptic-estimate}, we obtain
\begin{equation} \label{0-Theta4-a2-estimate}
\begin{split}
\norm{\Xi_{a,2}^{4}}
\lesssim  & \Big( \e^{-1}\normm{\ip  f}_{L^{2}_{x,v}(\nu)} +\normm{ \nu^{-\frac{1}{2}} \G(f,f)}_{L^{2}_{x,v}} \Big)  \normm{a}_{L^{2}_{x}}, \\
\norm{\Xi_{a,6}^{4}}
\lesssim & \Big( \e^{-1}\normm{\ip  f}_{L^{2}_{x,v}(\nu)} +\normm{ \nu^{-\frac{1}{2}} \G(f,f)}_{L^{2}_{x,v}} \Big)  \normm{a}_{L^{6}_{x}}^5.
\end{split}
\end{equation}

Integrating \eqref{0-test-equation-uniform-form} and combining  \eqref{0-Theta1 - a estimate}, \eqref{0-Theta2 - L2 estimate},  \eqref{0-Theta3 - a estimate}, \eqref{0-E-a-2-estimate} and \eqref{0-Theta4-a2-estimate}, we derive
\begin{equation} \label{0-tildea - L2 estimate step1}
\begin{split}
\int_{s}^{t} \normm{a}_{L^2_{x}}^2  \leq \;
&
 \e  \big[ G_{a} (t)-  G_{a} (s)\big]
 + \a^2\int_{s}^{t}  \norm{(1-\mathscr{P}_{\g}) f }_{L^2_{\g_{+}}}^2
 +  \int_{s}^{t}  \normm{\e^{-1}  \ip f}_{L^{2}_{x,v}({\nu})}^{2}\\
 %-------
 & +   \normm{ {\nu}^{-\frac{1}{2}}  \G(f,f)  }_{L^{2}_{x,v}}^{2} + \e  \int_{s}^{t} \normm{\pt_t \nabla \varphi_{a,2}}_{L^2_x} \Big( \normm{b}_{L^2_x} + \normm{\ip f}_{L^2_{x,v}} \Big).
\end{split}
\end{equation}

Combining  \eqref{0-test-equation-uniform-form}, \eqref{0-Theta1 - a L6 estimate}, \eqref{0-Theta2 - L2 estimate-L6-a},   \eqref{0-Theta3 - a L6 ettimate}, \eqref{0-E-a-2-estimate} and \eqref{0-Theta4-a2-estimate},   we obtain
\begin{equation} \label{0-tildea - L6 estimate}
\begin{split}
\normm{a}_{L^{6}_{x,v}} \lesssim \; &  \e\normm{\pt_t f}_{L^{2}_{x,v}} + \a  \norm{(1-\mathscr{P}_{\g})  f}_{L^2_{\g_{+}}}^{\frac{1}{2}}
\normm{\omega f}_{L^{\infty}_{x,v}}^{\frac{1}{2}} \\
%-----------
&  + \normm{\e^{-1}\ip  f}_{L^{2}_{x,v}(\nu)}  + \normm{\ip  f}_{L^{6}_{x,v}}+ \normm{ {\nu}^{-\frac{1}{2}} \Gamma(f,f)}_{L^{2}_{x,v}}.
\end{split}
\end{equation}

\noindent\textbf{Step 1.2. Estimate for $  \normm{\pt_t \nabla_x \varphi_{a,2}}_{L^2_x}$.}

In \eqref{0-test-equation-uniform-form} we choose the test function
$\psi_{a,2} = \pt_t\varphi_{a,2} \sqrt{\mu}$ and estimate each term.
Clearly, $\Xi_{a,2}^{4} = 0$. By \eqref{0-a-2-elliptic-equation}, we have
\begin{equation}\label{a-theta1-estimate-1}
\begin{split}
   \Xi_{a,2}^{1}
=\e \int_{\Omega} \pt_t \varphi_{a,2} \pt_t a
=\e \int_{\Omega} - \pt_t \varphi _{a,2}  \Delta_x( \pt_t\varphi_{a,2})
=\e \normm{\nabla_x \pt_t \varphi_{a,2}}_{L^2_x}^2.
\end{split}
\end{equation}
Noting $\mathscr{R}(\psi_{a,2})=\psi_{a,2}$, by the change of variables
as in \eqref{0-Theta2-boundary-calculation-a}, we have
\begin{equation}\label{0-Theta2 - L2 estimate-1}
\begin{split}
\norm{\Xi_{a,2}^{2} } \lesssim & \, \a  \norm{(1-\mathscr{P}_{\g})f}_{L^2_{\g_{+}}} \norm{\pt_t \varphi_{a,2}}_{L^2(\pt\O)}\lesssim  \a  \norm{(1-\mathscr{P}_{\g})f}_{L^2_{\g_{+}}} \normm{\nabla_x\pt_t \varphi_{a,2}}_{L^2_x},
\end{split}
\end{equation}
 where we used the trace theorem and Poincar\'{e}'s inequality. By H\"{o}lder's inequality,
\begin{equation} \label{a-theta3-estimate-1}
\begin{split}
\norm {\Xi_{a,2}^{3} } = &\Big| \iint_{\Omega \times \R^3}  v\cdot \nabla_x  \pt_t\varphi_{a,2} \sqrt{\mu} f \Big|
=  \Big|\int_{\Omega} \nabla_x  \pt_t\varphi_{a,2} \cdot b\dd x \Big|
\lesssim  \normm{b}_{L^2_x} \normm{\nabla_x \pt_t\varphi_{a,2}}_{L^2_x}.
\end{split}
\end{equation}
 Combining \eqref{0-test-equation-uniform-form} with \eqref{a-theta1-estimate-1}--\eqref{a-theta3-estimate-1} gives
\begin{align}
&\e \normm{\nabla_x \pt_t \varphi_{a,2} }_{L^2_x} \lesssim  \normm{b}_{L^2_x}+  \a  \norm{(1-\mathscr{P}_{\g})f}_{L^2_{\g_{+}}} .  \label{0-varphia - pt t estimate}
\end{align}

Finally, substituting \eqref{0-varphia - pt t estimate} into \eqref{0-tildea - L2 estimate step1}, we obtain
\begin{equation} \label{0-tildea - l2 estimate final}
\begin{split}
\int_{s}^{t} \normm{a}_{L^2_{x}}^2  \leq
& C_{a} \Big \{
 \e  G_{a} (t)-  \e G_{a} (s)
 + \a^2\int_{s}^{t}  \norm{(1-\mathscr{P}_{\g}) f }_{L^2_{\g_{+}}}^2 \\
 %-------
 &\qquad  +  \int_{s}^{t}  \Big ( \normm{b}_{L^2_x} + \normm{\e^{-1}  \ip f}_{L^{2}_{x,v}({\nu})}^{2} +   \normm{ {\nu}^{-\frac{1}{2}}  \G(f,f) }_{L^{2}_{x,v}}^{2} \Big)
 \Big \}.
\end{split}
\end{equation}
\medskip

\noindent\textbf{ Step 2.  Estimate for $b$.}

Because the estimates for $\int_{s}^{t}\|b\|_{L^{2}_{x}}\dd \tau$ and  $\| b\|_{L^{6}_{x}}$ require different test functions, we treat them separately.

\noindent\textbf{ Step 2.1.  Estimate for $\int_{s}^{t}\|b\|_{L^{2}_{x}}\dd \tau$.}

In \eqref{0-test-equation-uniform-form}, we choose the test function
\begin{equation}\label{0-psi-b-2-definition}
\begin{split}
 \psi_{b,2}(t,x,v) := & \sum_{i,j=1}^{3} \pt_{j} \varphi_{b,2,i} A_{ij}(v) + \sum_{i=1}^{3} \pt_{i}\varphi_{b,2,i}\chi_{4}(v) \frac{\sqrt{6}}{6}\\
 =&\sum_{i,j=1}^{3} \pt_{j} \varphi_{b,2,i} v_iv_j\sqrt{\mu}
 -\sum_{i=1}^{3} \pt_{i} \varphi_{b,2,i}\frac{|v|^2-1}{2}\sqrt{\mu}.
\end{split}
\end{equation}
Here the vector-valued function $\varphi_{b,2}$ satisfies the elliptic system
\begin{equation}\label{0-b-2-elliptic-equation}
\begin{split}
- \Delta_x \varphi_{b,2} = b \;\;  &\text{in } \Omega, \quad
 \varphi_{b,2} =0 \;\; \text{on } \pt \Omega.
\end{split}
\end{equation}
Standard elliptic theory \cite{gilbarg1977elliptic} ensures that \eqref{0-b-2-elliptic-equation} admits a unique solution satisfying
\begin{equation}\label{0-b-2-elliptic-estimate}
\begin{split}
\normm{\nabla^2_x \varphi_{b,2} }_{L^{2}_{x}} + \normm{\nabla_x \varphi_{b,2} }_{L^{2}_{x}} + \normm{\varphi_{b,2} }_{L^{2}_{x}} & \lesssim \normm{b}_{L^{2}_{x}}.
\end{split}
\end{equation}

We now estimate each term in \eqref{0-test-equation-uniform-form}.
For $ \Xi_{b,2}^{1}$, integration by parts yields
\begin{equation*}
\begin{split}
\int_{s}^{t}  \Xi_{b,2}^{1}\dd \tau = &\; \e\big[ G_{b} (t) - G_{b} (s)\big] -\e \int_{s}^{t}\iint_{\Omega \times \R^{3}} \pt_t\psi_{b,2}f,
\end{split}
\end{equation*}
where $G_{b} (t)$ is bounded by $\normm{f(t)}_{L^2_{x,v}}^2$. The contributions from $a$ and $b$ varnish due to \eqref{Burnette-orthogonal} and the identity $\int_{\R^{3}}\chi_{4}f\dd v=c$.
% and \eqref{0-b-2-elliptic-estimate}, the second term equals
%\begin{equation}\label{Theta-b2-01-estimate}
%\begin{split}
%&-\e\int_{s}^{t}\int_{\Omega} \Big[\sum_{i,j=1}^{3} \pt_t \pt_{j} \varphi_{b,2,i}
%\Big(\int_{\R^{3}} A_{ij}(v) f \dd v\Big)  + \sum_{i=1}^{3} \pt_t \pt_{i}\varphi_{b,2,i} %\frac{\sqrt{6}}{6} \Big(\int_{\R^{3}}\chi_{4}(v)f\dd v \Big)  \Big]\\
%---------------
%\lesssim &\; \e \int_{s}^{t}\normm{\pt_t \nabla_x \varphi_{b,2}}_{L^2_x} \big( %\normm{c}_{L^2_x} + \normm{\ip f}_{L^2_{x,v}}\big).
%\end{split}
%\end{equation}
Thus, by \eqref{0-b-2-elliptic-estimate}, we obtain
\begin{equation} \label{0-Theta1 - b estimate}
\begin{split}
 \norm{ \int_{s}^{t}  \Xi_{b,2}^{1} }
 \leq  &    \e\big[ G_{b} (t) - G_{b}(s) \big] + \e  \int_{s}^{t}\normm{\pt_t \nabla_x \varphi_{b,2}}_{L^2_x} \big( \normm{c}_{L^2_x} + \normm{\ip f}_{L^2_{x,v}} \big).
%---------------
\end{split}
\end{equation}

For $\Xi_{b,2}^{2}$, noting  $\mathscr{R} ( \psi_{p,2})  \neq \psi_{b,2}$,
we apply the change of variables $R_x v \mapsto v$ to obtain
\begin{equation}\label{0-Theta2-boundary-calculation-b-1}
\begin{split}
\Xi_{b,2}^{2}
%=& \iint_{\g_{+}} \psi_{b,2}f \dd \g
%- \iint_{\g_{-}} \psi_{b,2}
%\big[ (1-\a)  \mathscr{R}f
%+\a  \mathscr{P}_\g f \big]\dd \g  \\
%--------------
 =& \iint_{\g_{+}} \psi_{p,2} f \dd \g
-\iint_{\g_{+}} \mathscr{R} ( \psi_{p,2}) [ (1-\a) f
+\a  \mathscr{P}_\g f ]\dd \g \\
%-----------------------
%  =& \iint_{\g_{+}} \psi_{p,2}
% \big[ (1-\mathscr{P}_{\g})  f
% +\mathscr{P}_{\g} f\big]
%  \dd \g\\
%------------
%&-\iint_{\g_{+}} \mathscr{R} \big( \psi_{p,2}\big) \big[ (1-\a) (1-\mathscr{P}_{\g}) f
%+  \mathscr{P}_\g f \big]\dd \g \\
%-----------------------
 =& \iint_{\g_{+}}
  [\psi_{p,2}-\mathscr{R} ( \psi_{p,2})  ]\mathscr{P}_{\g}f
  \dd \g + \iint_{\g_{+}}
 [\psi_{p,2} - (1-\a)\mathscr{R} ( \psi_{p,2}) ]
  (1-\mathscr{P}_{\g}) f
  \dd \g \\
:= &I_{1}+I_{2},
\end{split}
\end{equation}
where we used \eqref{P-gamma-orthogonal}. For $I_{1}$, using the change of variables $R_x v \mapsto v$ and \eqref{Burnette-orthogonal}, we have
\begin{equation}\label{0-Theta2-boundary-calculation-b-2}
\begin{split}
I_{1}
%=& \iint_{\g_{+}}
%  \psi_{p,2}  \mathscr{P}_{\g} f
%  \dd \g - \iint_{\g_{+}} \mathscr{R} \big(\psi_{p,2}\big)  \mathscr{P}_{\g} f\dd \g\\
%-----------------------
% =& \iint_{\g_{+}} \psi_{p,2}  \sqrt{\mu} Z \dd \g
% -\iint_{\g_{+}} \mathscr{R} \big( \psi_{p,2}\big)  \sqrt{\mu} Z \dd \g\\
%-----------------------
 =& \iint_{\g_{+}} \psi_{p,2}  \sqrt{\mu} z \dd \g
-\iint_{\g_{-}}  \psi_{p,2} \sqrt{\mu} z\dd \g
%-----------------------
 = \iint_{\pt\O\times \R^3} \psi_{p,2}  \sqrt{\mu} z [n\cdot v] \dd v \dd S_x=0,
%-----------------------
%=&  \int_{\pt\O} z \Big[
%\sum_{i,j,l=1}^{3} \pt_{j} \varphi_{b,2,i} n_l
 %\big(\int_{\R^3}  \chi_{l} A_{ij} \dd v \big)+
% \sum_{i,l=1}^{3} \pt_{i}\varphi_{b,2,i} \frac{\sqrt{6}}{6}n_l
% \big(\int_{\R^3}  \chi_{l} \chi_{4}  \dd v  \big) \Big]\dd S_x=0.
\end{split}
\end{equation}
where we used the notation $z=z(t,x):=\sqrt{2\pi}\int_{n\cdot v>0} f [n\cdot v]\dd v$ and $\mathscr{P}_{\g} f = \sqrt{\mu} z$. For $I_{2}$, we apply the trace theorem and $\eqref{0-b-2-elliptic-estimate}$. Thus, we obtain
\begin{equation}\label{0-Thetbc2 - L2 estimate}
\begin{split}
%-----------------
\norm{\Xi_{b,2}^{2} }
\lesssim  \norm{ (1-\mathscr{P}_{\g}) f}_{L^2_{\g_{+}}} \norm{\nabla_x \varphi_{b,2}}_{L^2(\pt\O)}
\lesssim \norm{ (1-\mathscr{P}_{\g}) f}_{L^2_{\g_{+}}} \normm{b}_{L^2_x}.
\end{split}
\end{equation}

 To estimate $\Xi_{b,2}^{3}$, we use the expression in the second line of \eqref{0-psi-b-2-definition} and split
\begin{equation}\label{0-Theta3 - b estimate-psi-b-2}
\begin{split}
- v\cdot \nabla_x  \psi_{b,2}
%---------------
%---------------
%=&- \sum_{i,j,k=1}^{3}  \pt_{j} \pt_{k}\varphi_{b,2,i}
%v_i v_jv_k \sqrt{\mu}
%+ \sum_{i,k=1}^{3} \pt_{i} \pt_k  \varphi_{b,2,i}  v_k \frac{|v|^2-1}{2} \sqrt{\mu}\\
%---------------
%---------------
=&- \sum_{i,j,k=1}^{3}  \pt_{j} \pt_{k} \varphi_{b,2,i}
\P \left (v_i v_jv_k \sqrt{\mu} \right)
 +  \sum_{i,k=1}^{3} \pt_{i} \pt_k  \varphi_{b,2,i}  v_k \frac{|v|^2-1}{2} \sqrt{\mu}\\
%---------------
%---------------
& - \sum_{i,j,k=1}^{3}  \pt_{j} \pt_{k} \varphi_{b,2,i}
\ip \left (v_i v_jv_k \sqrt{\mu} \right):={K}_1+{K}_2+{K}_3.
\end{split}
\end{equation}
Direct calculation yields
\begin{equation}\label{b-estimate-K1}
\begin{split}
{K}_1
%=&- \sum_{i,j,k=1}^{3}  \pt_{j} \pt_{k} {\varphi}_{b,2,i} \P
%\big( v_i v_jv_k \sqrt{\mu}\big) \\
%---------------
=&- \sum_{i,j,k=1}^{3}  \pt_{j} \pt_{k} {\varphi}_{b,2,i} \sum_{l=1}^{3} v_l \sqrt{\mu} \big(\int_{\R^{3}}  v_i v_j v_k v_l \mu\dd v  \big)\\
= & -\sum_{l=1}^3 v_l\sqrt{\mu} \Big( \sum_{i=j,k=l} + \sum_{i\neq j, i=k,j=l}+ \sum_{i\neq j, i=l,k=j}  \Big)
  \pt_{j} \pt_{k}{\varphi}_{b,2,i} \int_{\mathbb{R}^3}v_i v_j v_k v_l \mu \dd v,
\end{split}
\end{equation}
where in the first equality we have used the identities
\begin{equation}\label{b2-estimate-Gauss-3}
\begin{split}
%&\int_{\R^{3}}
% v_i v_jv_k \mu\dd v = 0, \;\;\;\;
% \int_{\R^{3}}
% v_i v_jv_k (|v|^2-3) \mu\dd v = 0, \;\; i,j,k=1,2,3,\\
%---------------
&\int_{\mathbb{R}^3}v_i^2 v_j^2\mu \dd v=
 \left\{
   \begin{array}{ll}
     3, & \hbox{if $i=j$,} \\
     1, & \hbox{if $i\neq j$.}
   \end{array}
 \right.
\end{split}
\end{equation}
For each fixed $l\in \{1,2,3\}$, the inner sums in \eqref{b-estimate-K1} are computed as: \begin{equation}
\begin{split}
\sum_{i=j,k=l}
= & \Big( \sum_{i=l}+ \sum_{i\neq l}\Big) \partial_{i}\partial_{l} \varphi_{b,2,i}  \int_{\mathbb{R}^3}v_i^2 v_l^2 \mu \dd v
= 3\sum_{i=l} \pt_{i} \pt_{l}   \varphi_{b,2,i} + \sum_{i\neq l} \pt_{i} \pt_{l} \phi^{b}_i, \\
%-----------------
\sum_{i\neq j, i=k,j=l}
= &\sum_{i\neq l} \pt_i \pt_l  \varphi_{b,2,i} \int_{\mathbb{R}^3}v_i^2 v_l^2\mu \dd v
 = \sum_{i\neq l} \pt_i \pt_l \varphi_{b,2,i},\\
%-------------------------
\sum_{i\neq j, i=l,k=j}
= &\sum_{i\neq l} \pt_i \pt_i \varphi_{b,2,l} \int_{\mathbb{R}^3}v_i^2 v_l^2\mu \dd v
 = \sum_{i\neq l} \pt_i \pt_i \varphi_{b,2,l},
\end{split}
\end{equation}
where we have used \eqref{b2-estimate-Gauss-3} again. Substituting these into \eqref{b-estimate-K1} yields
\begin{equation}\label{b-estimate-K1-final-form}
\begin{split}
{K}_1= & -\sum_{l=1}^3 v_l\sqrt{\mu} \Big(3\sum_{i=l} \pt_{i} \pt_{l}{\varphi}_{b,2,i} + 2\sum_{i\neq l} \pt_{i} \pt_{l}{\varphi}_{b,2,i}
+ \sum_{j\neq l} \pt_{j} \pt_{j}{\varphi}_{b,2,l}\Big).
\end{split}
\end{equation}
This further leads to
\begin{equation}\label{0-pf-b-left1}
\begin{split}
 &\iint_{\Omega\times\mathbb{R}^3}   {K}_1  \mathbf{P}f
=-\sum_{l=1}^3 \int_{\O} b_l\Big(3 \pt_{l}\pt_{l}{\varphi}_{b,2,l} + \sum_{i\neq l} \pt_{i}\pt_{l}{\varphi}_{b,2,i}+ \sum_{i\neq l}
\pt_{i}\pt_{l}{\varphi}_{b,2,i}+\sum_{j\neq l} \pt_{j}\pt_{j}{\varphi}_{b,2,l}\Big).
\end{split}
\end{equation}
Moreover, direct calculation implies
\begin{equation}\label{0-pf-b-left2}
\begin{split}
 \iint_{\Omega\times\mathbb{R}^3} {K}_2  \mathbf{P}f
% =&\iint_{\Omega\times\mathbb{R}^3}  \sum_{i,k=1}^{3}\pt_{i}\pt_{k}{\varphi}_{b,2,i} v_k %\frac{|v|^2-1}{2}\sqrt{\mu}  \mathbf{P}f \dd v \dd x\\
= & \int_{\O} \Big(2\sum_{l=1}^3{b}_l\pt_{l}\pt_{l}{\varphi}_{b,2,l}+ 2 \sum_{i=1}^3 \sum_{k\neq i}{b}_k\pt_{k}\pt_{i}{\varphi}_{b,2,i} \Big),
\end{split}
\end{equation}
where we have used
$$
 \int_{\mathbb{R}^3}v_i^2\frac{|v|^2-1}{2} \mu \dd v=2,\;\; \;\; i=1, 2, 3.
 $$
 Combining \eqref{0-Theta3 - b estimate-psi-b-2}, \eqref{0-pf-b-left1} and \eqref{0-pf-b-left2}, we obtain
\begin{equation}\label{0-Theta3 - b estimate}
\begin{split}
{\Xi}_{b,2}^{3}
%=&-\iint_{\Omega \times \R^3} v\cdot \nabla_x {\psi}_{b,2} f \dd v \dd x \\
%------------
%=&\iint_{\Omega\times\mathbb{R}^3} ({K}_1 + {K}_2 + {K}_R) [\mathbf{P}f+ \ip f]\dd v \dd x\\
%------------
=&\iint_{\Omega\times\mathbb{R}^3} ({K}_1 + {K}_2)  \mathbf{P}f \dd v \dd x+ E_{b,2}\\
=&-\sum_{l=1}^3\int_{\O}{b}_l\big(\pt_{l}\pt_{l}{\varphi}_{b,2,l}+ \sum_{i\neq l}\pt_{i}\pt_{i}{\varphi}_{b,2,l}\big)\dd x+ E_{b,2}\\
=&-\sum_{l=1}^3\int_{\O}{b}_l\Delta_x{\varphi}_{b,2,l} \dd x+ E_{b,2}  =\|{b}\|_{L^2_x}^2+ E_{b,2},
\end{split}
\end{equation}
 where we used \eqref{0-b-2-elliptic-equation} and the orthogonality of $\mathbf{P}f$ and ${K}_R$. By \eqref{0-b-2-elliptic-estimate},
\begin{equation}\label{0-E-b-2-estimate-L6}
\begin{split}
 |E_{b,2}|=&
 \Big | \iint_{\Omega\times\mathbb{R}^3} ({K}_1 + {K}_2 + {K}_R) \ip f \dd v \dd x\Big |
%  \lesssim &  \normm{\ip f}_{L^{2}_{x,v}} \normm{\nabla_x^2 {\varphi}_{b,2} }_{L^{2}_{x}}\\
 \lesssim   \normm{\ip f}_{L^{2}_{x,v}} \normm{{b}}_{L^{2}_{x}}.
 \end{split}
\end{equation}
The term $\Xi_{b,2}^{4}$ is estimated as \eqref{0-Theta4-a2-estimate}.

Integrating \eqref{0-test-equation-uniform-form} and combining \eqref{0-Theta1 - b estimate},  \eqref{0-Thetbc2 - L2 estimate}, \eqref{0-Theta3 - b estimate} and \eqref{0-E-b-2-estimate-L6}, we obtain
\begin{equation} \label{0-tildeb - L2 estimate step1}
\begin{split}
\int_{s}^{t} \normm{{b}}_{L^2_{x}}^2 \lesssim \;&
 \e   \big [{G}_{b} (t)-  {G}_{b} (s) \big ]
 + \int_{s}^{t}  \norm{(1-\mathscr{P}_{\g})f}_{L^2_{\g_{+}}}^2  + \int_{s}^{t}   \normm{\e^{-1}  \ip  f }_{L^{2}_{x,v}( \nu )}^{2}   \\
 %-------
 &+  \int_{s}^{t}    \normm{  {\nu}^{-\frac{1}{2}} \G(f,f)  }_{L^{2}_{x,v}}^{2} + \e  \int_{s}^{t}\normm{\pt_t \nabla_x {\varphi}_{b,2}}_{L^2_x}  \big ( \normm{{c}}_{L^2_x} + \normm{\ip f }_{L^2_{x,v}}  \big ).
\end{split}
\end{equation}

\noindent\textbf{Step 2.2. }  Estimate for $  \normm{\pt_t \nabla_x
{\varphi}_{b,2}}_{L^2_x} $.

In \eqref{0-test-equation-uniform-form}, we choose the test function
${\psi}_{b,2} = \pt_t{\varphi}_{b,2} \cdot v \sqrt{ \mu }$ and estimate each term.
Clearly, ${\Xi}_{b,2}^{4}=0$. By \eqref{0-b-2-elliptic-equation}, we obtain
\begin{equation} \label{0-b-t-theta1-estimate}
\begin{split}
{\Xi}_{b,2}^{1}
% =&  \e \iint_{\Omega \times \R^{3}} \pt_t{\varphi}_{b,2} \cdot v  \sqrt{ \mu  } \pt_t f\dd %v \dd x
=  \e \int_{\Omega} \pt_t {\varphi}_{b,2}\cdot  \pt_t {b}  \dd x
=  -\e \int_{\Omega} \pt_t {\varphi}_{b,2}\cdot  \Delta_x \pt_t {\varphi}_{b,2}  \dd x
=  \e \normm{\nabla_x \pt_t {\varphi}_{b,2}}_{L^2_{x}}^2.
\end{split}
\end{equation}
Similarly to \eqref{0-Theta2-boundary-calculation-b-1}--\eqref{0-Thetbc2 - L2 estimate},  using the trace theorem and Poincar\'{e}'s inequality, we obtain
\begin{equation}\label{0-Thetb2 - L2 estimate}
\begin{split}
%-----------------
\norm{{\Xi}_{b,2}^{2} } \lesssim   \norm{(1-\mathscr{P}_{\g}) f}_{L^2_{\g_{+}}} \norm{ \pt_t {\varphi}_{b,2}}_{L^2(\pt\O)}
 \lesssim   \norm{(1-\mathscr{P}_{\g}) f}_{L^2_{\g_{+}}} \normm{\nabla_x \pt_t {\varphi}_{b,2}}_{L^2_x}.
\end{split}
\end{equation}
 Elementary computation and Poincar\'{e}'s inequality yield
\begin{equation} \label{0-b-t-theta3-estimate}
\begin{split}
\norm{ {\Xi}_{b,2}^{3}}
%= & \Big| \iint_{\Omega \times \R^3}  \pt_i \pt_t\varphi_{b,2,j}v_iv_j \sqrt{\mu} \big[\P %f+ \ip f \big]\Big|  \\
%---------------
 \lesssim & \normm{\nabla_x \pt_t{\varphi}_{b,2}}_{L^2_{x}} \big(\normm{{a}}_{L^2_x} + \normm{{c}}_{L^2_x} + \normm{\ip f }_{L^2_{x,v}}  \big).
\end{split}
\end{equation}
Collecting \eqref{0-test-equation-uniform-form} and \eqref{0-b-t-theta1-estimate}--\eqref{0-b-t-theta3-estimate} yields
\begin{equation}
\label{0-varphib - pt t estimate}
\begin{split}
\e \normm{\nabla_x \pt_t {\varphi}_{b,2} }_{L^2_{x}} \lesssim& \normm{{a}}_{L^2_{x}} + \normm{{c}}_{L^2_{x}}   + \normm{\ip f }_{L^2_{x,v}}+ \norm{(1-\mathscr{P}_{\g}) f}_{L^2_{\g_{+}}}.
\end{split}
\end{equation}

Finally, substituting \eqref{0-varphib - pt t estimate} into \eqref{0-tildeb - L2 estimate step1}, we obtain
\begin{equation} \label{0-tildeb - l2 estimate final}
\begin{split}
\int_{s}^{t} \normm{{b}}_{L^2_{x}}^2 \leq C_{b} \Big \{ &
 \e \big[ G_{b} (t)-   G_{b} (s)\big]
 + \int_{s}^{t}\norm{(1-\mathscr{P}_{\g}) f}_{L^2_{\g_{+}}}^2
 +  \delta_b \int_{s}^{t} \normm{{a}}_{L^2_x}^2\\
 %---------
 &+  \int_{s}^{t}   \Big ( \normm{{c}}_{L^2_x}^2+ \normm{\e^{-1}  \ip f}_{L^{2}_{x,v}(\nu)}^{2} +   \normm{ {\nu}^{-\frac{1}{2}} \G(f,f)   }_{L^{2}_{x,v}}^{2} \Big)  \Big \},
\end{split}
\end{equation}
where the small constant $\delta_b>0$ arises from H\"{o}lder's inequality.
\medskip

\noindent\textbf{Step 2.3.  Estimate for $\| b\|_{L^{6}_{x}}$.}

Note that the estimate for $\normm{b}_{L^6_{x,v}}$ cannot be established simultaneously with $\int_{s}^{t}\normm{b}_{L^2_{x,v}}$, since $\e^{-\frac{1}{2}}\norm{  (1-\mathscr{P}_{\g}) f }_{L^2_{\g_{+}}}$ (as in \eqref{0-Theta2 - L2 estimate-L6-a}) exceeds the boundary dissipation $\a^{\frac{1}{2}} \e^{-\frac{1}{2}} \norm{  (1-\mathscr{P}_{\g}) f }_{L^2_{\g_{+}}}$ in Proposition \ref{f-ft-Energy-estimate} when $\e \leq  \a <1$.  To overcome this, we estimate $\normm{\P f}_{L^6_{x,v}}$ separately by choosing a new test function.

In \eqref{0-test-equation-uniform-form}, we choose the test function
\begin{equation}\label{0-psi-b-n-definition-L6}
\begin{split}
 {\psi}_{b,6}(t,x,v) := & \sum_{i,j=1}^{3} \pt_{j} \varphi_{b,6,i} {A}_{ij}(v) + \sum_{i=1}^{3} \pt_{i}{\varphi}_{b,6,i} \chi_{4}(v) \frac{\sqrt{6}}{3}\\
 =&\sum_{i,j=1}^{3} \pt_{j} {\varphi}_{b,6,i} v_iv_j\sqrt{ \mu }
 -\sum_{i=1}^{3} \pt_{i} \varphi_{b,6,i}\sqrt{ \mu }.
\end{split}
\end{equation}
Here ${\varphi}_{b,6}(t,x)=({\varphi}_{b,6,1}(t,x), {\varphi}_{b,6,2}(t,x), {\varphi}_{b,6,3}(t,x))$ satisfies the elliptic system
\begin{equation}\label{0-b-6-elliptic-equation-L6}
\begin{split}
- \text{div} (\nabla^{\text{s}}_{x} {\varphi}_{b,6}) = {b}^5 - \sum \frac{\int_{\Omega}A_ix \cdot {b}^5 \dd x}{\int_{\Omega}\norm{A_i x}^2\dd x}A_ix \;\; &\text{in } \Omega,\\
{\varphi}_{b,6} \cdot n =0 \;\; &\text{on } \pt \Omega,\\
(\nabla^{\text{s}}_{x} {\varphi}_{b,6}) n = (\nabla^{\text{s}}_{x}{\varphi}_{b,6}: n \otimes n)n  \;\;  &\text{on } \pt \Omega,
\end{split}
\end{equation}
where $A_ix\in\mathcal{R}_{\O}$ defined in \eqref{def:axiss}, and ${b}^5=({b}^5_1, {b}^5_2, {b}^5_3)$. For a vector field
$M = (m_{i})_{i=1,2,3} : \Omega \to \R^{3}$, we define the gradient $\nabla_x M$, the symmetric gradient $\nabla^{\text{s}}_x M$ and  the antisymmetric gradient $\nabla^{\text{a}}_x M$ by
\begin{equation}\label{def-nabla-s}
\begin{split}
&(\nabla_x M)_{ij} := \frac{\pt m_{i}}{\pt x_{j}}, \quad
(\nabla^{\text{s}}_x M)_{ij} := \frac{1}{2}  \Big( \frac{\pt m_{i}}{\pt x_{j}} + \frac{\pt m_{j}}{\pt x_{i}}\Big), \quad
(\nabla^{\text{a}}_x M)_{ij} := (\nabla_x M)_{ij} - (\nabla^{\text{s}}_x M)_{ij}.
\end{split}
\end{equation}
The inner product of two matrixes $P=(p_{ij})_{i,j=1,2,3}$ and $Q=(q_{ij})_{i,j=1,2,3}$ is defined by $P : Q = \sum_{i,j=1}^3 p_{ij}q_{ij}$.

For each $j=1,2,3$, direct computation gives
\begin{align}\label{0-tilde-b5-compatible}
\int_{\Omega} A_{j}x \cdot \Big({b}^5 - \sum \frac{\int_{\Omega}A_{i}x \cdot {b}^5 \dd x}{\int_{\Omega}\norm{A_{i}x}^2\dd x}A_{i}x\Big) \dd x = \int_{\Omega}A_{j}x \cdot {b}^5 \dd x - \int_{\Omega}A_{j}x \cdot {b}^5 \dd x = 0,
\end{align}
which verifies the compatibility condition \eqref{elliptic-system-compatible-condition} for the elliptic system \eqref{0-b-6-elliptic-equation-L6} in all non-axisymmetric, axisymmetric, and spherical domains. Thus, by Lemma \ref{elliptic-system-theory} and \eqref{0-tilde-b5-compatible}, the elliptic system \eqref{0-b-6-elliptic-equation-L6} admits a unique strong solution satisfying
\begin{align}\label{0-b-6-elliptic-estimate-L6}
\normm{\nabla^2_x {\varphi}_{b,6}}_{L^{\frac{6}{5}}_{x}} + \normm{\nabla_x {\varphi}_{b,6}}_{L^{2}_{x}} + \normm{{\varphi}_{b,6}}_{L^{6}_{x}} &\lesssim \normm{{b}^5}_{L^{\frac{6}{5}}_{x}}=\normm{{b}}_{L^{6}_x}^5.
\end{align}

For ${\Xi}_{b,6}^{1}$ and ${\Xi}_{b,6}^{4}$, applying H\"{o}lder's inequality and  \eqref{0-b-6-elliptic-estimate-L6} directly yields
\begin{equation} \label{0-Theta1 - b estimate-L6}
\begin{split}
%-----------------
\norm{ {\Xi}_{b,6}^{1}  } \lesssim   \e  \normm{\pt_t {f}}_{L^2_{x,v}} \normm{{b}}_{L^6_{x}}^{5},\;\;
%---------------
\norm{{\Xi}_{b,6}^{4}}
\lesssim & \big( \e^{-1}\normm{\ip f}_{L^{2}_{x,v}( \nu )} +\normm{  \nu ^{-\frac{1}{2}} \Gamma(f,f) }_{L^{2}_{x,v}} \big)\normm{{b}}_{L^6_{x}}^{5}.
\end{split}
\end{equation}

For $ \Xi_{b,6}^{2}$, the boundary condition $(\nabla^{\text{s}}_{x} {\varphi}_{b,6}) n = (\nabla^{\text{s}}_{x}{\varphi}_{b,6}: n \otimes n)n$   on $\pt \Omega$ implies
\begin{equation}\label{reflection-invariant-condition}
\begin{split}
&\mathscr{R}\big(\psi_{b,6}(t,x,v) \big)-\psi_{b,6}(t,x,v)\\
%--------------
= & \sum_{i,j=1}^{3} \pt_{j} \varphi_{b,6,i} \big[ A_{ij}(R_x v)-{A}_{ij}(v) \big] + \frac{\sqrt{6}}{3}\sum_{i=1}^{3} \pt_{i}\varphi_{b,6,i} \big[\chi_{4}(R_x v)-\chi_{4}(v)  \big]\\
%--------------
 = & -2(v\cdot n) \sum_{i,j=1}^{3} \pt_{j} \varphi_{b,6,i}  \big[
 v_i n_j+ n_iv_j -2(v\cdot n)n_i n_j \big]\\
 %--------------
  = & -2(v\cdot n)\Big[  \sum_{k=1}^{3} v_k \Big( \sum_{j=1}^{3} \pt_{j} \varphi_{b,6,k} n_j+\sum_{j=1}^{3} \pt_{k} \varphi_{b,6,j}n_j - 2 \sum_{i,j=1}^{3} \pt_{j} \varphi_{b,6,i} n_i n_j n_k\Big) \Big]
\\
%--------------
  = & -4(v\cdot n)\Big[   v\cdot  \Big( (\nabla^{\text{s}}_{x}\varphi_{b,6}) n -  (\nabla^{\text{s}}_{x}\varphi_{b,6}: n \otimes n) n\Big) \Big]
  = 0.
\end{split}\end{equation}
Thus, similar to \eqref{0-Theta2-boundary-calculation-a} and \eqref{0-Theta2 - L2 estimate-L6-a}, we derive
\begin{equation}\label{0-Theta2 - b estimate-L6-b}
\begin{split}
 \norm{{\Xi}_{b,6}^{2}}
  \lesssim\; &  \a  \norm{(1-\mathscr{P}_{\g})  f}_{L^2_{\g_{+}}}^{\frac{1}{2}}
\normm{\omega f}_{L^{\infty}_{x,v}}^{\frac{1}{2}} \normm{{b}}_{L^6_{x}}^{5}.
\end{split}
\end{equation}

For ${\Xi}_{b,6}^{3}$, using the expression in the second line of \eqref{0-psi-b-n-definition-L6}, we have
\begin{equation}\label{0-Theta3 - b estimate-psi-b}
\begin{split}
- v\cdot \nabla_x {\psi}_{b,6}
%---------------
%---------------
%=&- \sum_{i,j,k=1}^{3}  \pt_{j} \pt_{k} {\varphi}_{b,6,i}
%v_i v_jv_k \sqrt{{\mu}}
%+ \sum_{i,l=1}^{3} \pt_{i} \pt_l  {\varphi}_{b,6,i}  v_l \sqrt{{\mu}}\\
%---------------
%---------------
=&- \sum_{i,j,k=1}^{3}  \pt_{j} \pt_{k} {\varphi}_{b,6,i}
{\P} \left (v_i v_jv_k \sqrt{{\mu}} \right)
+  \sum_{i,l=1}^{3} \pt_{i} \pt_l  {\varphi}_{b,6,i}  v_l \sqrt{{\mu}} \\
%-----------------
& - \sum_{i,j,k=1}^{3}  \pt_{j} \pt_{k} {\varphi}_{b,6,i}
\ip \left (v_i v_jv_k \sqrt{{\mu}} \right) \\
 :=&\hat{K}_1  + \hat{K}_2 + \hat{K}_{R}.
\end{split}
\end{equation}
For $\hat{K}_1$, calculations similarly to \eqref{b-estimate-K1}--\eqref{b-estimate-K1-final-form} yield
\begin{equation}\label{0-b-macro-pb-1-K1}
\begin{split}
 \hat{K}_{1}=& -\sum_{l=1}^3 v_l\sqrt{{\mu}} \Big(3\sum_{i=l} \pt_{i} \pt_{l}{\varphi}_{b,6,i} + 2\sum_{i\neq l} \pt_{i} \pt_{l}{\varphi}_{b,6,i}
+ \sum_{j\neq l} \pt_{j} \pt_{j}{\varphi}_{b,6,l}\Big).
\end{split}
\end{equation}
Substituting \eqref{0-b-macro-pb-1-K1} into \eqref{0-Theta3 - b estimate-psi-b} gives
\begin{equation}\label{0-Theta3 - b estimate-psi-b2}
\begin{split}
&- v\cdot \nabla_x  {\psi}_{b,6}
%---------------
=\hat{K}_1+ \hat{K}_2+ \hat{K}_3  \\
%------
= & \sum_{l=1}^3 v_l\sqrt{\mu} \Big[\sum_{i=1}^3  \pt_{i}\pt_{l}{\varphi}_{b,6,i} -  \Big(3\sum_{i=l} \pt_{i}\pt_{l}{\varphi}_{b,6,i}
+ 2\sum_{i\neq l} \pt_{i}\pt_{l} {\varphi}_{b,6,i}  + \sum_{i\neq l}  \pt_{i}\pt_{l} {\varphi}_{b,6,l}\Big) \Big]+\hat{K}_{R}\\
%------
= & \sum_{l=1}^3 v_l\sqrt{\mu} \Big[- 2\sum_{i=l}  \pt_{i}\pt_{l}{\varphi}_{b,6,i}-
\sum_{i\neq l}  \pt_{i}\pt_{l}{\varphi}_{b,6,i} - \sum_{i\neq l}  \pt_{i}\pt_{l}{\varphi}_{b,6,l} \Big]+\hat{K}_{R}\\
%------
= & \sum_{l=1}^3 v_l\sqrt{\mu} \Big[-\Big(\sum_{i=l}  \pt_{i}\pt_{l}{\varphi}_{b,6,i}+\sum_{i\neq l} \pt_{i}\pt_{l}{\varphi}_{b,6,i}\Big)
-\Big( \sum_{i=l}  \pt_{i}\pt_{l}{\varphi}_{b,6,i}+\sum_{i\neq l}  \pt_{i}\pt_{i}{\varphi}_{b,6,l} \Big)\Big]+\hat{K}_{R}\\
%------
= & \sum_{l=1}^3 v_l\sqrt{\mu} \big[- \partial_{l}(\text{div}{\varphi}_{b,6})-  \Delta_x {\varphi}_{b,6,l}\big]
+\hat{K}_{R}\\
%------
= & - \sqrt { \mu } v \cdot \text {div} \big( \nabla^s_x {\varphi}_{b,6} \big )+\hat{K}_{R}.
\end{split}
\end{equation}
Thus, using \eqref{0-b-6-elliptic-equation-L6}, we have
\begin{equation}\label{0-Theta3 - b L6 estimate}
\begin{split}
{\Xi}_{b,6}^{3}
%=&-\iint_{\Omega \times \R^3} v\cdot \nabla_x {\psi}_{b,6} f \dd v \dd x\\
%---------------
=&   \iint_{\Omega \times \R^3}\big[ -\sqrt {\mu } v \cdot \text {div} \big( \nabla^s_x {\varphi}_{b,6} \big ) +\hat{K}_{R}\big]\big[\P f+\ip f \big]\dd v \dd x\\
%---------------
=& -\int_{\O}{ b}\cdot \text {div} \big( \nabla^s_x {\varphi}_{b,6} \big ) \dd x + {E}_{b,6}
%------------
= \int_{\O} {b}\cdot  \Big( {b}^5 - \sum \frac{\int_{\Omega}A_ix \cdot {b}^5 \dd x}{\int_{\Omega}\norm{A_ix}^2\dd x}A_i x \Big)\dd x + {E}_{b,6} \\
=&  \normm{{b} }^6_{L^{6}_{x}}  + {E}_{b,6}+  {F}_{b,6}.
\end{split}
\end{equation}
The terms ${E}_{b,6}$ and ${F}_{b,6}$ are bounded via
%\begin{equation*}
%\begin{split}
% {E}_{b,6}:=& \iint_{\Omega\times\R^{3}}  -\text{div} \big( \nabla^s_x {\varphi}_{b,6} \big %)\cdot  v
%   \sqrt {\mu } \ip f \dd v  \dd x
%  +
%  \iint_{\Omega \times \R^3} {K}_{R} \ip f \dd v \dd x,\\
%------------
% {F}_{b,6}:=&  \sum_{i=1}^{3} \Big(\int_{\O} {b}\cdot   A_i x \dd x \Big)
%   \frac{\int_{\Omega}A_ix \cdot {b}^5 \dd x}{\int_{\Omega}\norm{A_ix}^2\dd x}.
%\end{split}
%\end{equation*}
\eqref{0-b-6-elliptic-estimate-L6}:
\begin{equation}\label{0-EF-b-2t-estimate-L6}
\begin{split}
\norm{{E}_{b,6}} \lesssim  \normm{\ip  f}_{L^{6}_{x,v}} \normm{{b}}_{L^{6}_{x}}^{5},
\quad
%--------------
\norm{{F}_{b,6}} \lesssim   \normm{{b}}_{L^{2}_{x}}\normm{{b}}_{L^{6}_{x}}^{5}.
\end{split}
\end{equation}

Combining \eqref{0-test-equation-uniform-form}, \eqref{0-Theta1 - b estimate-L6}, \eqref{0-Theta2 - b estimate-L6-b},  \eqref{0-Theta3 - b L6 estimate} and \eqref{0-EF-b-2t-estimate-L6}, we obtain
\begin{equation} \label{0-tildeb - L6 estimate}
\begin{split}
\normm{{b}}_{L^{6}_{x,v}} \lesssim &  \e\normm{\pt_t f}_{L^{2}_{x,v}}+ \normm{{b}}_{L^{2}_{x}} + \a  \norm{(1-\mathscr{P}_{\g})  f}_{L^2_{\g_{+}}}^{\frac{1}{2}}
\normm{\omega f}_{L^{\infty}_{x,v}}^{\frac{1}{2}} \\
%-----------
&  + \normm{\ip  f}_{L^{6}_{x,v}} + \normm{\e^{-1}\ip  f}_{L^{2}_{x,v}(\nu)}  + \normm{ {\nu}^{-\frac{1}{2}} \Gamma(f,f)}_{L^{2}_{x,v}}.
\end{split}
\end{equation}

\noindent{\bf Step 3.  Estimate for ${c}$.}

\noindent{\bf Step 3.1.  Estimate for $\int_{s}^{t}\|{c}\|_{L^{2}_{x}}\dd \tau$ and $\|{c}\|_{L^{6}_{x}}$.}

In \eqref{0-test-equation-uniform-form}, we choose the test function
\begin{equation}\label{0-psi-b-n-definition}
\begin{split}
 {\psi}_{c,q}(t,x,v) := &  \sum_{i=1}^{3} \pt_{i} {\varphi}_{c,q}(t,x) \sqrt{10} {B}_{i}(v),\quad q\in\{2, 6\},
\end{split}
\end{equation}
where ${\varphi}_{c,2}(x)$ and ${\varphi}_{c,6}(x)$ satisfy the elliptic equations
\begin{align}
&- \Delta_x {\varphi}_{c,2} = {c} \;\;
\text{in } \Omega, \quad  {\varphi}_{c,2}= 0\;\; \text{on } \pt \Omega,\label{0-c-2-elliptic-equation} \\
%-------------
&- \Delta_x {\varphi}_{c,6} = {c}^{5} -\frac{1}{\norm{\Omega}}\int_{\Omega}{c}^{5} \dd x\;\; \text{in } \Omega, \quad  \pt_n {\varphi}_{c,6}= 0\;\; \text{on } \pt \Omega, \quad
\int_{\Omega} {\varphi}_{c,6} \dd x =0, \label{0-c-6-elliptic-equation}
\end{align}
respectively. ${\varphi}_{c,2}$ and ${\varphi}_{c,6}$ satisfy elliptic estimates like analogous to those in \eqref{0-a-2-elliptic-estimate} and \eqref{0-a-6-elliptic-estimate}.

We now estimate each term in \eqref{0-test-equation-uniform-form}.
For $ \Xi_{c,2}^{1}$, integration by parts shows that the contribution from ${\P}f$ vanishes due to \eqref{Burnette-orthogonal}. Thus, similarly to \eqref{0-Theta1 - a estimate}, we obtain
%\begin{equation}\label{Theta-c2-01-estimate}
%\begin{split}
%\norm{ \iint_{\Omega \times \R^{3}} \pt_t{\psi}_{c,2} f }
%---------------
%=&\sum_{i=1}^{3} \iint_{\Omega \times \R^{3}}  \pt_t \pt_{i} {\varphi}_{c,2}  \sqrt{10} %{B}_{i} f\\
%---------------
%=& \sum_{i=1}^{3}\int_{\Omega} \pt_t \pt_{i} \tilde{\varphi}_{c,2} \int_{\R^{3}}   v_i(|v|^2 - %5)\sqrt{\tilde{\mu}} \Big(\tilde{\P} \frac{\tilde{\Phi}}{\sqrt{\tilde{\mu}}}+\ipt\frac{\tilde{\Phi}}{\sqrt{\tilde{\mu}}} \Big) \\
%---------------
%=& \norm{ \sum_{i=1}^{3}\int_{\Omega} \pt_t \pt_{i} {\varphi}_{c,2} \int_{\R^{3}}   \sqrt{10} {B}_{i}(v) \big[{\P}f +\ip f \big] }\\
%---------------
%\lesssim &\normm{\pt_t \nabla_x {\varphi}_{c,2}}_{L^2_x} \normm{\ip f}_{L^2_{x,v}}.
%\end{split}
%\end{equation}
\begin{equation} \label{0-Theta1 - c estimate}
\begin{split}
\Big |\int_{s}^{t}    \Xi_{c,2}^{1}  \dd \tau \Big|
\lesssim  &\;  \e  \big[G_{c} (t)-  G_{c} (s)\big] + \e  \int_{s}^{t}\normm{\pt_t \nabla_x \varphi_{c,2}}_{L^2_x}  \normm{\ip f}_{L^2_{x,v}}.
\end{split}
\end{equation}

For $\Xi_{c,6}^{1}$, the elliptic estimate for $\varphi_{c,6}$ yields
\begin{equation}
\begin{split} \label{0-Theta1 - c L6 estimate}
\norm{ \Xi_{c,6}^{1}  }
%---------------
=&  \e  \Big | \iint_{\Omega\times\R^{3} } \sum_{i=1}^{3}  \pt_{i} {\varphi}_{c,6}   \sqrt{10} {B}_{i}(v)  \pt_t f \Big |
\lesssim   \e  \normm{\pt_t {f}}_{L^2_{x,v}} \normm{{c}}_{L^6_{x}}^{5}.
\end{split}
\end{equation}

For $ \Xi_{c,2}^{2}$, noting that ${\psi}_{c,2}$ is not specular reflection invariant, we use the change of variables to obtain
\begin{equation}\label{Theta2-boundary-calculation 1}
\begin{split}
\Xi_{c,2}^{2}
%=& \iint_{\g_{+}} \frac{\tilde{\psi}_{p,q}}{\sqrt{\tilde{\mu}}} \Tilde{\Phi}_{} \dd \g
%- \iint_{\g_{-}} \frac{\tilde{\psi}_{p,q}}{\sqrt{\tilde{\mu}}}  \Tilde{\Phi}_{} \dd \g \\
%-----------------------
%=& \iint_{\g_{+}} {\psi}_{c,2} {f} \dd \g
%- \iint_{\g_{-}} {\psi}_{c,2} \Big( (1-\a)  \mathscr{R}{f}
%+\a  {\mathscr{P}}_{\g}{f}\Big)\dd \g \\
%----------------------
%=& \iint_{\g_{+}} {\psi}_{c,2}  \big[(1-\mathscr{P}_{\g}) {f}
%+  {\mathscr{P}}_{\g}{f}  \big]\dd \g
%- \iint_{\g_{+}}  \mathscr{R}({\psi}_{c,2}) \big[ (1-\a) (1-\mathscr{P}_{\g}) {f}
%+  {\mathscr{P}}_{\g}{f}\big]\dd \g \\
%-----------------------
%----------------------
=& \iint_{\g_{+}} \big[{\psi}_{c,2}-(1-\a)\mathscr{R}({\psi}_{c,2})
\big] (1-\mathscr{P}_{\g}) {f}\dd \g
-  \iint_{\g_{+}}  \mathscr{R}({\psi}_{c,2}) {\mathscr{P}}_{\g}{f}\dd \g,
\end{split}\end{equation}
where we used \eqref{P-gamma-orthogonal}. The term involving ${\mathscr{P}}_{\g}{f}$ vanishes due to the identities
\begin{equation*}
\begin{split}
(n\cdot v)^2=\sum_{i,j=1}^3 v_iv_jn_in_j, \ \ \  \ \ \  \int_{n\cdot v>0}(|v|^2-5)v_k^2\mu \dd v=0, \ \ \ k=1,2,3.
\end{split}
\end{equation*}
Thus, by the trace theorem and the elliptic estimate of $\varphi_{c,2}$, we obtain
\begin{equation}\label{0-Theta-c2-2-L2-estimate}
\begin{split}
%-----------------
\norm{\Xi_{c,2}^{2} } \lesssim   \norm{ (1-\mathscr{P}_{\g}) f}_{L^2_{\g_{+}}} \norm{\nabla_x \varphi_{c,2}}_{L^2(\pt\O)}\lesssim  \norm{ (1-\mathscr{P}_{\g}) f}_{L^2_{\g_{+}}} \normm{c}_{L^2_x}.
\end{split}
\end{equation}

For ${\Xi}_{c,6}^{2}$, the condition ${\pt_n} {\varphi}_{c,6}|_{\pt
\O}= 0$ implies that $\mathscr{R}({\psi}_{c,6}= {\psi}_{c,6}$.
Consequently, ${\Xi}_{c,6}^{2}$ can be treated
similarly to \eqref{0-Theta2-boundary-calculation-a} and  \eqref{0-Theta2 - L2 estimate-L6-a}:
\begin{equation}\label{0-Theta2 - L2 estimate-L6-c}
\begin{split}
\norm{{\Xi}_{c,6}^{2}}
%=& \Big| \a \iint_{\g_{+}} \sum_{i=1}^{3} \pt_{i} {\varphi}_{c,6} \sqrt{10} {B}_{i}(v) %(1-\mathscr{P}_{\g})  f \dd \g\Big| \\
 \lesssim &\a  \norm{(1-\mathscr{P}_{\g})  f}_{L^2_{\g_{+}}}^{\frac{1}{2}}
\normm{\omega f}_{L^{\infty}_{x,v}}^{\frac{1}{2}} \normm{ {c}}_{L^6_x}^5.
\end{split}
\end{equation}

For $\Xi_{c,q}^{3}$ ($q\in \{2,6\}$), direct computation gives
\begin{equation}\label{0-Theta3 - c p estimate}
\begin{split}
\Xi_{c,q}^{3}
%=&-\iint_{\Omega \times \R^3} (v\cdot \nabla_x {\psi}_{c,2}) f  \\
%---------------
%---------------
=& -\sum_{i,j=1}^{3}\int_{\Omega} \pt_i \pt_{j} {\varphi}_{c,q} \int_{\R^{3}}   v_iv_j(|v|^2 - 5)\sqrt{\mu} f
%---------------
=  -\frac{10}{\sqrt{6}}\int_{\Omega} c \Delta_x \varphi_{c,q}  + E_{c,q},
\end{split}
\end{equation}
where $E_{c,q}$  arises from the contribution $\ip f$, and for the $\P f$ contribution we have used
\begin{equation}\label{0-v_iv_j(|v|^2 - 5)-orthogonal}
\begin{split}
&\int_{\R^{3}}   v_iv_j(|v|^2 - 5)\sqrt{\mu}  \chi_{k}(v)\dd v
=  0,   \; \;
\int_{\R^{3}}   v_iv_j(|v|^2 - 5)\sqrt{\mu}  \chi_{4}(v)\dd v
=    \frac{10}{\sqrt{6}}\d_{ij}
%&\int_{\R^{3}}   v_iv_j(|v|^2 - 5)\sqrt{\mu} \chi_{k}(v)\dd v
%=   \int_{\R^{3}}   v_iv_jv_k  (|v|^2 - 5)\tilde{\mu}
%=   0,   \quad    i,j,k=1,2,3.
\end{split}\end{equation}
for $i,j=1,2,3$ and  $k=0,1,2,3$. Using \eqref{0-c-2-elliptic-equation} and \eqref{0-c-6-elliptic-equation}, we have
\begin{equation}\label{0-Theta3 - c estimate}
\begin{split}
\Xi_{c,2}^{3}
%=&-\iint_{\Omega \times \R^3} (v\cdot \nabla_x {\psi}_{c,2}) f  \\
%---------------
%---------------
%=& -\sum_{i,j=1}^{3}\int_{\Omega} \pt_i \pt_{j} {\varphi}_{c,2} \int_{\R^{3}}   %v_iv_j(|v|^2 - 5)\sqrt{\mu} \big[\P f+\ip f \big] \\
%---------------
= & -\frac{10}{\sqrt{6}}\int_{\Omega} c \Delta_x \varphi_{c,2}  + E_{c,2}
= \frac{10}{\sqrt{6}} \normm{ c }^2_{L^{2}_{x}}  + E_{c,2},\\
%--------------
\Xi_{c,6}^{3}
%= & -\iint_{\Omega \times \R^3} (v\cdot \nabla_x \psi_{c,6}) f \dd v \dd x\\
%---------------
% = & -\int_{\Omega}  \sum_{i=1}^{3} \pt_{i}\pt_j {\varphi}_{c,6}  \int_{\R^3} v_i v_j %(|v|^2-5)\sqrt{\mu}  \big[ \P {f} + \ip f\big] \dd v \dd x\\
%---------------
= & -\frac{10}{\sqrt{6}} \int_{\Omega}{c}  \Delta_x {\varphi}_{a,6} \dd x
 + {E}_{c,6}
 =  \frac{10}{\sqrt{6}} \int_{\Omega}  {c} \big(  {c}^5 -\frac{1}{|\O|} \int_{\O}   {c}^5\dd x \big)\dd x
 + {E}_{c,6}\\
  = & \frac{10}{\sqrt{6}}\normm{ c}^6_{L^{6}_{x}} + {F}_{c,6}+ {E}_{c,6}.
\end{split}
\end{equation}
 The remainders $E_{c,2}$, $E_{c,6}$  and $F_{c,6}$ are controlled via elliptic estimates as in \eqref{0-a-2-elliptic-estimate} and \eqref{0-a-6-elliptic-estimate}:
\begin{equation}\label{0-E-c-2-estimate}
\begin{split}
\norm{E_{a,2}} \lesssim \normm{c}_{L^{2}_{x}} \normm{\ip f}_{L^{2}_{x,v}}, \;\;
\norm{E_{a,6}} \lesssim  \normm{c}_{L^{6}_{x}}^{5} \normm{\ip f}_{L^{6}_{x,v}}\;\;
\norm{F_{a,6}} \lesssim  \normm{{c}}_{L^{2}_{x}}\normm{{c}}_{L^{6}_{x}}^{5}.
\end{split}
\end{equation}

The terms $\Xi_{c,2}^{4}$ and $\Xi_{c,6}^{4}$ are estimated similarly to \eqref{0-Theta4-a2-estimate}.
%\begin{equation} \label{0-Theta1 - c estimate-L6}
%\begin{split}
%\norm{\Xi_{c,2}^{4}}
%\lesssim  & \Big( \e^{-1}\normm{\ip  f}_{L^{2}_{x,v}(\nu)} +\normm{ \nu^{-\frac{1}{2}} %\G(f,f)}_{L^{2}_{x,v}} \Big)  \normm{c}_{L^{2}_{x}},\\
%-----------------
%\norm{{\Xi}_{c,6}^{4}}
%\lesssim & \Big( \e^{-1}\normm{\ip f}_{L^{2}_{x,v}( \nu )} +\normm{  \nu ^{-\frac{1}{2}} %\Gamma(f,f) }_{L^{2}_{x,v}} \Big)\normm{{c}}_{L^6_{x}}^{5}.
%\end{split}
%\end{equation}

Integrating \eqref{0-test-equation-uniform-form} and combining \eqref{0-Theta1 - c estimate}, \eqref{0-Theta-c2-2-L2-estimate}, \eqref{0-Theta3 - c estimate} and \eqref{0-E-c-2-estimate}, we have
\begin{equation} \label{0-tildec - L2 estimate step1}
\begin{split}
\int_{s}^{t} \normm{c}_{L^2_{x}}^2  \lesssim \; &
 \e  \big[G_{c} (t)- G_{c} (s)\big]
 + \int_{s}^{t} \norm{(1-\mathscr{P}_{\g})f}_{L^2_{\g_{+}}}^2 +\int_{s}^{t}  \normm{\e^{-1}  \ip f}_{L^{2}_{x,v}(\nu)}^{2}\\
 %-------
 &+  \int_{s}^{t}     \normm{ \nu^{-\frac{1}{2}} \G(f,f)   }_{L^{2}_{x,v}}^{2} + \e  \int_{s}^{t}\normm{\pt_t \nabla \varphi_{c,2}}_{L^2_x} \normm{\ip f}_{L^2_{x,v}}.
\end{split}
\end{equation}

Combining \eqref{0-test-equation-uniform-form},
\eqref{0-Theta1 - c L6 estimate}, \eqref{0-Theta2 - L2 estimate-L6-c}, \eqref{0-Theta3 - c estimate} and \eqref{0-E-c-2-estimate},  we obtain
\begin{equation} \label{0-tildec - L6 estimate}
\begin{split}
\normm{{c}}_{L^{6}_{x,v}} \lesssim \; &  \e\normm{\pt_t f}_{L^{2}_{x,v}} +\normm{{c}}_{L^{2}_{x}}+ \a  \norm{(1-\mathscr{P}_{\g})  f}_{L^2_{\g_{+}}}^{\frac{1}{2}}
\normm{\omega f}_{L^{\infty}_{x,v}}^{\frac{1}{2}} \\
%-----------
&  + \normm{\e^{-1}\ip  f}_{L^{2}_{x,v}(\nu)}  + \normm{\ip  f}_{L^{6}_{x,v}} + \normm{ {\nu}^{-\frac{1}{2}} \Gamma(f,f)}_{L^{2}_{x,v}}.
\end{split}
\end{equation}

\noindent\textbf{Step 3.2. Estimate for $  \normm{\pt_t \nabla_x \varphi_{c,2}}_{L^2_x} $.}

In \eqref{0-test-equation-uniform-form}, we choose the test function as
$\psi_{c,2} = \pt_t\varphi_{c,2}  \chi_4(v)$ and estimate each term.
Clearly, $\Xi_{c,2}^{4}=0$. Using \eqref{0-c-2-elliptic-equation}, we obtain
\begin{equation}\label{0-c-theta1-estimate}
\begin{split}
\Xi_{c,2}^{1}
=& \e \iint_{\Omega \times \R^{3}}  \pt_t \varphi_{c,2} \int_{\R^3} \chi_4  \pt_t f\dd v
%= \e \int_{\Omega} \pt_t \varphi_{c,2} \pt_t c \\
%-------------
= \e \int_{\Omega} - \pt_t\varphi_{c,2} \Delta_x \pt_t\varphi_{c,2}
= \e \normm{\nabla_x \pt_t \varphi_{c,2}}_{L^2_x}^2.
\end{split}
\end{equation}
Noting $\mathscr{R}(\psi_{c,2})=\psi_{c,2}$, we deduce similarly to \eqref{0-Theta2 - L2 estimate-1} that
\begin{equation}\label{0-Thetc2 - L2 estimate-1}
\begin{split}
\norm{\Xi_{c,2}^{2} } \lesssim & \a  \norm{(1-\mathscr{P}_{\g})f}_{L^2_{\g_{+}}} \normm{\pt_t \varphi_{c,2}}_{L^2_x}\lesssim  \a  \norm{(1-\mathscr{P}_{\g})f}_{L^2_{\g_{+}}} \normm{\nabla_x\pt_t \varphi_{c,2}}_{L^2_x}.
\end{split}
\end{equation}
By oddness of the integrand involving $a$ and $c$ contributions, we have
\begin{equation} \label{0-c-theta3-estimate}
\begin{split}
\norm {  \Xi_{c,2}^{3} }
%= & \norm{\iint_{\Omega \times \R^3}  (v\cdot \nabla_x  \pt_t\varphi_{c,2}) \chi_4(v) f %\dd v \dd x}\\
%------------------
=& \Big|\sum_{i=1}^{3} \int_{\Omega }  \pt_i \pt_t \varphi_{c,2}\int_{ \R^3}v_i \chi_4(v) f \Big|
%------------------
\lesssim   \normm{\nabla_x \pt_t \varphi_{c,2}}_{L^2_x}  \big( \normm{b}_{L^2_x}
+  \normm{\ip f}_{L^{2}_{x,v}}
 \big).
\end{split}
\end{equation}
Combining \eqref{0-test-equation-uniform-form} with \eqref{0-c-theta1-estimate}--\eqref{0-c-theta3-estimate} gives
\begin{equation}\label{0-varphic - pt t estimate}
\begin{split}
\e \normm{\nabla_x \pt_t \varphi_{c,2} }_{L^2_x} \lesssim
& \normm{b}_{L^2_x}+ \a \norm{(1-\mathscr{P}_{\g})f}_{L^2_{\g_{-}}}
+ \normm{\ip f}_{L^{2}_{x,v}}.
\end{split}
\end{equation}

Finally, substituting \eqref{0-varphic - pt t estimate} into \eqref{0-tildec - L2 estimate step1} yields
\begin{equation} \label{0-tildec - l2 estimate final}
\begin{split}
\int_{s}^{t} \normm{c}_{L^2_{x}}^2  \leq  C_c\Big \{
 &
 \e G_{c} (t)-  \e G_{c} (s)
 + \int_{s}^{t}\norm{(1-\mathscr{P}_{\g})f}_{L^2_{\g_{-}}} +
  \d_c\int_{s}^{t}
  \normm{b}_{L^2_x}^2 \\
 %-------
 &+  \int_{s}^{t}   \Big ( \normm{\e^{-1}  \ip f}_{L^{2}_{x,v}(\nu)}^{2}  +  \normm{\nu^{-\frac{1}{2}}  \G(f,f)  }_{L^{2}_{x,v}}^{2} \Big)   \Big \},
\end{split}
\end{equation}
where the small constant $\delta_c>0$   arises from Young's inequality.

\noindent\textbf{Step 4.  Combination of the estimates for $a$,  $b$ and  $c$.}

Choose $\delta_{b} = (2^8 4 C_{b} C_{c}^2)^{-1}$ and $\delta_{c} = (4 C_{c})^{-1}$.
A direct computation of
$$
(2^{8}C_{b}C_{c}^2)^{-1} \times \eqref{0-tildea - l2 estimate final} + (2^{11/2}C_{b}C_{c}^2)^{-1}  \times\eqref{0-tildeb - l2 estimate final} +\eqref{0-tildec - l2 estimate final}
$$
yields \eqref{0-P-f-macro-L2}. Furthermore, combining \eqref{0-tildea - L6 estimate}, \eqref{0-tildeb - L6 estimate} and \eqref{0-tildec - L6 estimate},
we obtain \eqref{0-Pf-L6}. This completes the proof of Proposition \ref{0-macro-L2-L6-estimate}.
\end{proof}

\medskip

The equation for $\pt_t f$ shares the same linear structure as equation \eqref{f-eq} for $f$, differing only in the source term. Moreover, $\pt_t f$ also satisfies the mass conservation law $\iint_{\Omega\times\R^3}  \sqrt{\mu}\pt_t f\dd v \dd x = 0$. Therefore, Proposition \ref{0-macro-L2-L6-estimate} applies to  $\pt_t f$ and yields the following result:

\begin{corollary} \label{0-macro-L2-L6-estimate-pt}
Under the same assumption as in Proposition \ref{0-macro-L2-L6-estimate}, there holds
\begin{equation}  \label{0-t-P-f-macro-L2-pt}
\begin{split}
\int_{s}^{t} \normm{ \P \pt_t f }_{L^{2}_{x,v}}^2 \lesssim \;
& \e  G_{1}(t) - \e G_{1}(s)+ \int_{s}^{t}  \norm{ (1-\mathscr{P}_{\g}) \pt_t f  }_{L^2_{\g_{+}}}^2 \\
%-------------
+ & \int_{s}^{t}  \big[    \normm{\e^{-1} \ip \pt_t f }_{L^{2}_{x,v}( \nu)}^{2}
+\normm{ {\nu}^{-\frac{1}{2}} [\Gamma( \pt_t f, f)+\Gamma( f, \pt_t f)]}_{L^{2}_{x,v}}^2  \big],
\end{split}
\end{equation}
where $|G_1(t)| \lesssim \normm{f(t)}_{L^{2}_{x,v}}^2 + \normm{\pt_t f(t)}_{L^{2}_{x,v}}^2$.
\end{corollary}
\bigskip

%%%%%%%%%%%%%%%%%%%%%%%%%%%%%%%%%%
%%%%%%%%%%%%%%%%%%%%%%%%%%%%%%%%%%

\subsection{Nonlinear Estimates} \label{subsection-Nonlinear-Estimate}\

This subsection establishes an $L^\infty$ estimate for the linear equation
 and provides nonlinear estimates for the collision operator $\G(f,f)$.

\begin{proposition}\label{L-infty-bd-unst} \
Let $0<\e\leq \e_0$, where $\e_0$ is the constant determined in Proposition \ref{lemma-fbar-infty-unst-0}. Assume $g, \pt_t g \in L^{\infty} ([0, T] \times \Omega \times \mathbb{R}^{3})$ and $f_0, \pt_t f_0 \in L^{\infty} (\Omega \times \mathbb{R}^{3})$ with $0 < T \leq \infty$. Let $f$ be the solution to the linear Boltzmann equation  on $[0, T]$:
\begin{equation}\label{unst-lin-orig}
\begin{split}
 \e\pt_t f+ v\cdot\nabla_x  f   +  \e^{-1} L f = g \;\;\;\;   &\text{ in } [0, T]\times\O\times\mathbb{R}^3,  \\
 f |_{\g_-} = (1-\a) \mathscr{R}f + \a \mathscr{P}_{\g}f \;\;\;\;   &\text{ on } [0, T]\times\pt\O\times\mathbb{R}^3, \\
 f(t, x, v)|_{t=0}  =   f_0(x,v)  \;\;\;\;     &\text{ on }   \O \times \mathbb{R}^3.
\end{split}
\end{equation}
Then the following estimates hold  for all $t\in [0,T]$:
\begin{align}
\| \o f(t)\|_{L^\infty_{x,v}} \lesssim \; &  \| \o  f_0\|_{L^\infty_{x,v}} + \e^{-\frac{1}{2}}  \sup_{0\leq s\leq t}  \|{\P}{f}(s)\|_{L^6_{x,v}}
+  \e^{-\frac{3}{2}}  \sup_{0\leq s\leq t}  \|\ip {f}(s)\|_{L^2_{x,v}} \nonumber \\
&  +
\e   \sup_{0\leq s\leq t}  \|\langle v\rangle^{-1}\o  g(s)\|_{L^\infty_{x,v}}, \label{Lifnty-bd-unst}\\
%-----------------------
%-------------------------
\| \o f(t)\|_{L^\infty_{x,v}} \lesssim \;
&  \| \o f_0\|_{L^\infty_{x,v}}+  \e^{-\frac{3}{2}}  \sup_{0\leq s\leq t}  \| {f}(s)\|_{L^2_{x,v}} +  \e   \sup_{0\leq s\leq t}  \|\langle v\rangle^{-1}\o   g(s)\|_{L^\infty_{x,v}}.
\label{Lifnty-bd-unst-2}
\end{align}
\end{proposition}

\begin{proof}[\textbf {Proof}] \
The proof relies on Proposition \ref{lemma-fbar-infty-unst-0}.  Recall the scaling transformations \eqref{stretch-x-y} and \eqref{change-fun-t} for the domain $\O\subset \mathbb{R}^3$.
For $0\leq t\leq \e^2 T_0$, we have
\begin{equation}\label{f-barf-relat}
\begin{split}
% & \sup_{0\leq \bar{s}\leq \bar{t}}  \|wf\|_{L^\infty_{x,v}(\Omega\times \mathbb{R}^3)}
% =  \sup_{0\leq \bar{s}\leq \bar{t}}   \|w \hat{f}\|_{L^\infty_{y,v}(\Omega_\e\times % \mathbb{R}^3)}, \\
& \sup_{0\leq t\leq \e^2 T_0}  \|\mathbf{P}f(t)\|_{L^6_{x,v}(\Omega\times \mathbb{R}^3)}
   =  \sup_{0\leq \bar{t}\leq T_0}  \e^{1\over 2}\|\mathbf{P}\bar{f}(\bar{t})\|_{L^6_{y,v}(\Omega_\e\times \mathbb{R}^3)}, \\
%-----------------------
&  \sup_{0\leq t\leq \e^2 T_0}  \|(\mathbf{I-P})f(t)\|_{L^2_{x,v}(\Omega\times \mathbb{R}^3)}
   =    \sup_{0\leq \bar{t}\leq T_0}   \e^{\frac{3}{2}} \|(\mathbf{I-P})\bar{f}(\bar{t})\|_{L^2_{y,v}(\Omega_\e\times \mathbb{R}^3)},
\end{split}
\end{equation}
where $\bar{t}=\e^{-2}t\in [0,T_0]$ from \eqref{stretch-x-y}.
Applying Proposition \ref{lemma-fbar-infty-unst-0} and these relations, we obtain for $0\leq t\leq \e^2 T_0$:
\begin{align}
\| \o f(t)\|_{L^\infty_{x,v}}   \lesssim
\; &  e^{-\frac{\nu_0 }{2\e^2}t} \|  \o f_0\|_{L^\infty_{x,v}}
   + o(1) \sup_{0\leq s\leq \e^2 T_0}  \|   \o f(s)  \|_{L^\infty_{x,v}}+  \e^{-\frac{1}{2}} \sup_{0\leq s\leq \e^2 T_0} \| \mathbf{P}f(s)  \|_{L^6_{x,v} } \nonumber\\
&
   + \e^{-\frac{3}{2}} \sup_{0\leq s\leq \e^2 T_0} \|  (\mathbf{I-P})f(s)  \|_{L^2_{x,v}}
   +  \e  \sup_{0\leq s\leq \e^2 T_0}  \| \langle v\rangle^{-1}  \o g(s) \|_{L^\infty_{x,v}},\label{Lifnty-bd-T0}
\\
%--------------------
\|  \o f(t)\|_{L^\infty_{x,v}}   \lesssim \
& e^{-\frac{\nu_0 }{2\e^2}t} \|   \o f_0\|_{L^\infty_{x,v}}
   + o(1) \sup_{0\leq s\leq \e^2 T_0}  \|   \o f(s)  \|_{L^\infty_{x,v}}+  \e^{-\frac{3}{2}} \sup_{0\leq s\leq \e^2 T_0} \| f(s)  \|_{L^2_{x,v}} \nonumber \\
&
         +  \e  \sup_{0\leq s\leq \e^2 T_0}  \| \langle v\rangle^{-1}  \o g(s) \|_{L^\infty_{x,v} }. \label{Lifnty-bd-T0-2}
\end{align}

Define
\begin{eqnarray*}
\begin{split}
D(s) := & \  o(1)\| w f(s)  \|_{L^\infty_{x,v}}
  + \e^{-\frac{1}{2}} \| \mathbf{P}f(s) \|_{L^6_{x,v}}  +    \e^{-\frac{3}{2}} \|  (\mathbf{I} - \mathbf{P})f(s) \|_{L^2_{x,v}}
  +  \e \| \langle v\rangle^{-1} \o  g(s)\|_{L^\infty_{x,v}}.
\end{split}\end{eqnarray*}
Then \eqref{Lifnty-bd-T0} becomes
\begin{equation}\label{Lifnty-bd-T0-1}
\begin{split}
\| wf(t)\|_{L^\infty_{x,v}}   \lesssim
\; &  e^{-\frac{\nu_0 }{2\e^2}t} \| wf_0\|_{L^\infty_{x,v}}
   +  \sup_{0\leq s\leq \e^2 T_0} D(s),\qquad 0\leq t\leq \e^2 T_0.
\end{split}
\end{equation}

Applying \eqref{Lifnty-bd-T0-1} iteratively yields
\begin{equation}\label{wf-n}
\begin{split}
  \|  \o f( n \e^2 T_{0})\|_{L^\infty_{x,v}}
\leq \; &  e^{- \frac{\nu_0 }{2}T_{0}} \|  \o f((n-1) \e^2 T_{0}) \|_{L^\infty_{x,v}}
+  \sup_{(n-1) \e^2 T_{0} \leq s \leq n \e^2 T_{0} }  D(s) \\
\leq \;&  e^{- \frac{2\nu_0 }{2} T_{0}} \|  \o f((n-2) \e^2 T_{0}) \|_{L^\infty_{x,v}} + \sum_{j=0}^{1}
  e^{- \frac{j\nu_0 }{2} T_{0}}
 \sup_{(n-2) \e^2 T_{0} \leq s \leq n \e^2 T_{0} } D(s)
\\
& \ \ \ \vdots  \\
\leq \;&  e^{- \frac{n\nu_0 }{2} T_{0}}\|  \o f_{0}  \|_{L^\infty_{x,v}}  +  \sum_{j=0}^{n-1}
  e^{- \frac{j\nu_0 }{2} T_{0}}
  \sup_{0 \leq s \leq  n \e^2 T_{0}}  D(s) \\
\leq \;& C_1 \|  \o f_{0}  \|_{L^\infty_{x,v}} +  C_1\sup_{0 \leq s \leq  n \e^2 T_{0}}  D(s)
\end{split}
\end{equation}
for some constant $C_1>0$, provided $ T_0>0$ is sufficiently large.
% where in the last step we have used
%$$
% e^{- \frac{n\nu_0 }{2}T_{0}}<\infty, \ \ \ \  \sum_{j=0}^{n-1} e^{- \frac{j\nu_0 %}{2}T_{0}} < \infty
%\ \ \ \ \hbox{ for large }  T_0>0.
%$$

For arbitrary $t>0$, choose $n\in \mathbb{N}$ such that $t \in [n \e^2 T_{0}, (n+1) \e^2 T_{0}]$. Combining the estimate \eqref{wf-n} with \eqref{Lifnty-bd-T0-1}, we obtain
\begin{equation*}
\begin{split}
\|  \o f(t)\|_{L^\infty_{x,v}}
&\leq  e^{- \frac{\nu_0 }{2\e^2}(t-n\e^2 T_{0})} \|  \o f(n \e^2 T_{0}) \|_{L^\infty_{x,v}}
   +  \sup_{n\e^2 T_{0} \leq s \leq t}  D(s)
\leq  C\|  \o f_0 \|_{L^\infty_{x,v}} + C\sup_{0 \leq s \leq t}  D(s)
\end{split}
\end{equation*}
for some constant $C>0$. Absorbing the small term $ C  o(1) \sup_{0\leq s\leq t} \|  \o  f(s)  \|_{L^\infty_{x,v}}$, we proves \eqref{Lifnty-bd-unst}. The estimate \eqref{Lifnty-bd-unst-2} follows similarly using \eqref{Lifnty-bd-T0-2}.
\end{proof}
\medskip

We now derive estimates for the nonlinear collision operator $\G(f,f)$.

\begin{lemma} \label{Gamma-L2-Linfty-estimate-1}
Recall the definition of $\G$ in \eqref{L-Definition}. For $\omega = e^{ {\b} \norm{v}^2}$ with $0 <  {\b} \ll \frac{1}{4}$, we have
\begin{align}
&\normm{ \nu^{-\frac{1}{2}} {\G}({f}, {g})}_{L^2_{x,v}} \lesssim   \normm{\omega g}_{L^\infty_{x,v}} \normm{f}_{L^2_{x,v}(\nu)}, \label{G-L2-estimate-1}\\
%------------------------
&\normm{ \nu^{-\frac{1}{2}} {\G}({f}, {g})}_{L^2_{x,v}} \lesssim  \normm{\omega f}_{L^\infty_{x,v}} \normm{ g}_{L^2_{x,v}(\nu)},\label{G-L2-estimate-2}\\
%---------------
&\normm{\nu^{-1}\omega {\G}({f}, {g})}_{L^\infty_{x,v}} \lesssim  \normm{\omega f}_{L^\infty_{x,v}} \normm{\omega g}_{L^\infty_{x,v}}, \label{weighted-G-Linfty-estimate-1}\\
%--------------------------
&\normm{  {\nu}^{-\frac{1}{2}}   {\G}( {\P}f, {\P}g)}_{L^2_{x,v}} \lesssim \normm{ {\P}f  {\P}g}_{L^2_{x,v}}. \label{0-G-Pf}
\end{align}
\end{lemma}

\begin{proof}[\textbf{Proof.}] \
We first note that
\begin{align*}
&\normm{  {\nu}^{-\frac{1}{2}}  {\G}({f}, {g})  }_{L^2_{x,v}}
\lesssim \normm{ \omega g }_{L^\infty_{x,v}}
\normm{  {\nu}^{-\frac{1}{2}} {\G}(f, \omega^{-1}) }_{L^2_{x,v}}.
%------------------
%&\normm{ {\nu}^{-\frac{1}{2}} {\G}({f}, {g})  }_{L^2_{x,v}}
%\lesssim \normm{ \omega f }_{L^\infty_{x,v}} \normm{ {\nu}^{-\frac{1}{2}}  %{\G}(\omega^{-1}, {g}) }_{L^2_{x,v}}.
\end{align*}
%Using the conservation laws $\norm{v}^2+\norm{u}^2 = \norm{v'}^2+\norm{u'}^2$ and $v +u = %v'+u'$, we have
%$$
% {\nu}^{-\frac{1}{2}}  \norm{(v-u)\cdot w} {\mu}^{\frac{1}{2}}(u) \lesssim  %{\nu}^{-\frac{1}{2}} (\norm{v}+ \norm{u}) {\mu}^{\frac{1}{2}}(u) \lesssim (1 + \norm{v} + %\norm{u})^{\frac{1}{2}} {\mu}^{\frac{1}{2}-\eta}(u)
%$$
%for some $\eta \in (0, \frac{1}{2})$. Hence,
%\begin{align*}
%&\int_{\R^3}  \norm{ {\nu}^{-\frac{1}{2}}   {\G}(f,\omega^{-1})(v) }^2 \dd v \\
%\lesssim&\iint_{\R^3 \times \R^3} \left[1+\norm{v'}+\norm{u'}\right]  \norm{ f(v')}^2 %\omega^{-2}(u')  \dd u \dd v \\
%&+\iint_{\R^3 \times \R^3 } \left[1+\norm{v'}+\norm{u'}\right]\norm{f(u')}^2\omega^{-2}(v')  %\dd u \dd v \\
%&+\iint_{\R^3 \times \R^3 } \left[1+\norm{v}+\norm{u}\right]\norm{f(v)}^2\omega^{-2}(u)  %\dd u \dd v \\
%&+\iint_{\R^3 \times \R^3 } \left[1+\norm{v}+\norm{u}\right]\norm{f(u)}^2\omega^{-2}(v)  %\dd u \dd v.
%\end{align*}
%Applying the change of variables $(v,u) \leftrightarrow (v',u')$ to the first term, $(v,u) %\leftrightarrow (u',v')$ to the second term, and $(v,u) \leftrightarrow (u,v)$ to the last term, we derive
Following Lemma 2.13 in \cite{Esposito2017}, we obtain
\begin{align*}
\int_{\R^3} \norm{   {\nu}^{-\frac{1}{2}} {\G}(f,\omega^{-1})(v)}^2 \dd v
\lesssim&\iint_{\R^3 \times \R^3 } (1+\norm{v}+\norm{u})f^2(v)\omega^{-2}(u)  \dd u \dd v
\lesssim \int_{\R^3}  {\nu} \norm{f(v)}^2 \dd v,
\end{align*}
which proves \eqref{G-L2-estimate-1}. The estimate \eqref{G-L2-estimate-2}  follows similarly.

Next, \eqref{weighted-G-Linfty-estimate-1} follows from the bound $\normm{\nu^{-1}\omega  {\G}(\omega^{-1}, \omega^{-1})}_{L^\infty_{x,v}} \lesssim  1$, due to the exponential decay of $ {\mu}$.

Finally, for $0< \d \ll 1$, we have
$$
  \normm{ {\mu}^{-\d} \norm{ {\P}f}}_{L^\infty_{v}} \lesssim \normm{ {\P}f}_{L^p_{v}} \;\;\; \text{for any } \; 1 \le p \le \infty.
$$
It follows that
\begin{align*}
\normm{ {\nu}^{-\frac{1}{2}}  {\G}( {\P}f, {\P}g)}_{L^2_{x,v}} \lesssim& \normm{ {\nu}^{-\frac{1}{2}} {\G}( {\mu}^{\d},  {\mu}^{\d})}_{L^2_{v}}\normm{ {\P}f  {\P}g}_{L^2_{x,v}}
\lesssim \normm{ {\P}f  {\P}g}_{L^2_{x,v}},
\end{align*}
which complete the proof of \eqref{0-G-Pf}.
\end{proof}
\medskip

\begin{corollary} \label{G - L2 Linfty L3 L6 estimate}
Let $f,g \in L^2([0,T]\times \Omega \times \R^3)$ with $0<T\leq \infty$, and let $S_{j}f, S_{j}g \ge 0$ $(j=1,2)$ be defined as in  Proposition \ref{f - L2L3 estimate}. Suppose that for $t\in [0,T]$,
\begin{equation*}
\begin{split}
\norm{ {a}(h)}+\sum_{i=1}^{3}\norm{ {b}_{i}(h)} +\norm{ {c}(h)}  \le S_{1}h(t,x) + S_{2}h(t,x)\quad \text{ for } h\in \{f,g\},
\end{split}
\end{equation*}
 where ${a}(h),  {b}_{i}(h)$ and ${c}(h)$ are coefficients of $ \P h$ with respect to the basis $\{ \chi_i \}$.
Then for $\omega = e^{ {\b} \norm{v}^2}$ with $0 <  {\b} \ll \frac{1}{4}$, the following estimate holds:
\begin{equation}\label{0-Gamma-f-g-estimate}
\begin{split}
&\normm{  \nu^{-\frac{1}{2}} \G(f,g)}_{L^2_{t,x,v}}+
\normm{  \nu^{-\frac{1}{2}}  \G(g,f)}_{L^2_{t,x,v}}\\
%----------------
\lesssim \; &\e^{\frac{1}{2}}\Big [ {\e}^{-1} \normm{\ip f}_{L^2_{t,x,v}(\nu)} +  {\e}^{-1} \normm{S_{2} f}_{L^2_{t,x}} \Big] \Big[\e^{\frac{1}{2}} \normm{\o g}_{L^\infty_{t,x,v}}\Big] \\
%----------------
& +   \normm{S_{1}f}_{L^2_{t}L^3_{x}} \Big[ \e^{\frac{1}{2}}
 \normm{ \o  g}_{L^\infty_{t,x,v}}\Big]^{\frac{2}{3}} \Big[ {\e}^{-1} \normm{\ip  g}_{L^\infty_{t}L^2_{x,v}( {\nu})}\Big]^{\frac{1}{3}}
 %----------------
+\normm{S_{1}f}_{L^2_{t}L^3_{x}}\normm{\P g}_{L^\infty_{t} L^6_{x,v}}.
\end{split}
\end{equation}

\end{corollary}

\begin{proof}[\textbf{Proof.}] \
To estimate $ {\G}(f,g)$, we decompose
\begin{equation*}
\begin{split}
\norm{    {\G}(f,g)} \le
&\norm{  {\G}( {\P}f, {\P}g)} + \norm{   {\G}( {\P}f,\ip g)}
 + \norm{   {\G}(\ip f,g)}.
\end{split}\end{equation*}
By Lemma \ref{Gamma-L2-Linfty-estimate-1}, we obtain
\begin{equation}\label{Gamma-estimate -1}
\begin{split}
\normm{   {\nu}^{-\frac{1}{2}}   {\G}( {\P}f, {\P}g)}_{L^2_{t,x,v}} \lesssim & \normm{S_{1}f}_{L^2_{t}L^3_{x}}\normm{ {\P} g}_{L^\infty_{t} L^6_{x,v}} + \normm{S_{2}f}_{L^2_{t,x}}\normm{ \o  g}_{L^\infty_{t,x,v}},\\
 %----------------
\normm{   {\nu}^{-\frac{1}{2}}   {\G}(\ip f,g)}_{L^2_{t,x,v}} \lesssim &
\e^{\frac{1}{2}}\Big [ {\e}^{-1} \normm{\ip f}_{L^2_{t,x,v}( {\nu})}  \Big] \Big[\e^{\frac{1}{2}} \normm{ \o  g}_{L^\infty_{t,x,v}}\Big],
\\
\normm{   {\nu}^{-\frac{1}{2}}   {\G}( {\P}f,\ip g)}_{L^2_{t,x,v}}
\lesssim & \normm{S_{1}f}_{L^2_{t}L^3_{x}}\normm{\ip  g}_{L^\infty_{t} L^6_{x,v}} + \normm{S_{2}f}_{L^2_{t,x}}\normm{ \o  g}_{L^\infty_{t,x,v}}\\
\lesssim &  \normm{S_{1}f}_{L^2_{t}L^3_{x}} \Big[ \e^{\frac{1}{2}}
 \normm{ \o  g}_{L^\infty_{t,x,v}}\Big]^{\frac{2}{3}} \Big[\frac{1}{\e}\normm{\ip  g}_{L^\infty_{t}L^2_{x,v}( {\nu})}\Big]^{\frac{1}{3}}\\
 %----------------
&+\e^{\frac{1}{2}} \Big[\frac{1}{\e} \normm{S_{2} f}_{L^2_{t,x}} \Big] \Big[\e^{\frac{1}{2}} \normm{ \o  g}_{L^\infty_{t,x,v}}\Big],
\end{split}
\end{equation}
where the last inequality uses interpolation. This establishes
\eqref{0-Gamma-f-g-estimate} for $ {\G}(f,g)$.

For the term $ {\G}(g,f)$,  we decompose it similarly:
\begin{equation*}
\begin{split}
\norm{    {\G}(g,f)} \le
&\norm{    {\G}( {\P}g, {\P}f)} + \norm{     {\G}(g,\ip f)}
 + \norm{    {\G}(g, {\P} f)}.
\end{split}\end{equation*}
The first two terms can be bounded  in the same way as \eqref{Gamma-estimate -1}.
%Applying Lemma \ref{Gamma-L2-Linfty-estimate-1} yields
%\begin{equation*}
%\begin{split}
%&\normm{   {\nu}^{-\frac{1}{2}}   {\G}( {\P}g, {\P}f)}_{L^2_{t,x,v}} \lesssim %\normm{S_{1}f}_{L^2_{t}L^3_{x}}\normm{ {\P} g}_{L^\infty_{t} L^6_{x,v}} + %\normm{S_{2}f}_{L^2_{t,x}}\normm{wg}_{L^\infty_{t,x,v}},\\
%----------------
%&\normm{   {\nu}^{-\frac{1}{2}}   {\G}(g,\ip f)}_{L^2_{t,x,v}} \lesssim %\e^{\frac{1}{2}}\Big [ \e^{\frac{1}{2}} \normm{\o g }_{L^\infty_{t,x,v}}  \Big]
%\Big [ \frac{1}{\e} \normm{\ip f}_{L^2_{t,x,v}} \Big],
%\end{split}
%\end{equation*}
For the last term, we first use
\begin{equation*}
\begin{split}
\normm{   {\nu}^{-\frac{1}{2}}   {\G}(g, {\P} f)}_{L^2_{t,x,v}}
\lesssim & \normm{S_{1}f}_{L^2_{t}L^3_{x}}\normm{g}_{L^\infty_{t} L^6_{x,v}} + \normm{S_{2}f}_{L^2_{t,x}}\normm{w g}_{L^\infty_{t,x,v}},\\
%\lesssim & \normm{S_{1}f}_{L^2_{t}L^3_{x}}
% \Big[ \e^{\frac{1}{2}}\normm{w g}_{L^\infty_{t,x,v}}\Big]^{\frac{2}{3}} %\Big[\frac{1}{\e}\normm{\ip  g}_{L^\infty_{t}L^2_{x,v}( {\nu})}\Big]^{\frac{1}{3}} \\
% %----------------
%&+\normm{S_{1}f}_{L^2_{t}L^3_{x}}\normm{ {\P} g}_{L^\infty_{t} L^6_{x,v}}
%+\normm{S_{2} f}_{L^2_{t,x}}  \normm{w g}_{L^\infty_{t,x,v}},
\end{split}
\end{equation*}
%where interpolation is used in the last inequality. This completes the proof of \eqref{0-Gamma-f-g-estimate} for $\Gamma(g,f)$.
and then handle it analogously to \eqref{Gamma-estimate -1}.
\end{proof}
\medskip

\begin{corollary} \label{f-L6-estimate-final} \
%Let $\e \in (0, \e_0]$, where $\e_0\in (0,1)$ is the constant from Proposition %\ref{L-infty-bd-unst}.
Let $f$ be the solution to \eqref{f-eq} on $ [0,T]$ with $0<T\leq \infty$. Then, for any $t\in [0,T]$,
\begin{equation} \label{0-pf6 - bound3-final}
\begin{split}
\normm{\P f}_{L^6_{x,v}}^2
\lesssim &\left[\!\left[{f}_{0}\right]\!\right]_{1}^2 + \mathscr{E}_{1}[f](t) + \mathscr{D}_{1}[f](t)
 + \d\e\normm{\omega f}_{L^{\infty}_{x,v}}^2\\
 & + \left[\!\left[{f}_{0}\right]\!\right]_{1}^4 +\mathscr{E}_{1}^3[f](t) + \mathscr{D}_{1}^2[f](t) +\e^2\normm{\omega f}_{L^{\infty}_{x,v}}^4,
\end{split}
\end{equation}
where $\d>0$ is a sufficiently small constant and $\omega = e^{ {\b} \norm{v}^2}$ with $0 <  {\b} \ll \frac{1}{4}$.
\end{corollary}

\begin{proof}[\textbf{Proof.}] \
 We start from the estimate \eqref{0-Pf-L6}. Both $ \e\normm{\pt_t  f }_{L^{2}_{x,v}}$ and
$\normm{\P f }_{L^{2}_{x,v}}$ are bounded by $\mathscr{E}_{1}[f](t)$.
 For the boundary term in \eqref{0-Pf-L6},  Young's inequality yields
\begin{equation*}%\label{L-t-infty-boundary-term}
\begin{split}
  \a  \norm{(1-\mathscr{P}_{\g})  f}_{L^2_{\g_{+}}}^{\frac{1}{2}}
\normm{\omega f}_{L^{\infty}_{x,v}}^{\frac{1}{2}}
\lesssim &
\a^\frac{3}{2}\left(\frac{\a}{\e}\right)^{\frac{1}{2}} \norm{  (1-\mathscr{P}_{\g}) f }_{L^2_{\g_{+}}} + \d  \e^{\frac{1}{2}} \normm{\omega  f }_{L^{\infty}_{x,v}} \\
 \lesssim
&  \left[\!\left[f_{0}\right]\!\right]_1+ \mathscr{D}_{1}^{\frac{1}{2}}[f](t)  +\d  \e^{\frac{1}{2}}\normm{\omega  f }_{L^{\infty}_{x,v}},
\end{split}
\end{equation*}
where $\d>0$ is sufficiently small, and we used the estimate
\begin{equation}\label{L-t-infty-boundary-term}
\begin{split}
\frac{\a}{\e} \norm{(1-\mathscr{P}_{\g})f(t)}_{L^2_{\g_{+}}}^2
%=& \frac{\a^4}{\e} \iint_{\g_{+}} |(1-\mathscr{P}_{\g})f(t)|^2(t) \dd \g \\
=& \frac{\a}{\e}  \iint_{\g_{+}} |(1-\mathscr{P}_{\g})f_0|^2
+ \frac{\a}{\e}  \int_{0}^{t} \iint_{\g_{+}} \frac{\dd\big[(1-\mathscr{P}_{\g})f\big]^2(s)}{\dd s}
%=& \frac{\a^4}{\e} \norm{(1-\mathscr{P}_{\g}) f_0 }_{\g_{+}}^2 + \frac{2\a^4}{\e} %\int_{0}^{t} \iint_{\g_{+}}(1-\mathscr{P}_{\g})f \pt_t (1-\mathscr{P}_{\g})f  \\
\lesssim  \left[\!\left[f_{0}\right]\!\right]_1^2  +\mathscr{D}_{1}(t).
\end{split}
\end{equation}
 Meanwhile,  the term $\normm{\e^{-1}\ip  f }_{L^{2}_{x,v}(\nu)} $ satisfies
\begin{equation}\label{L-t-infty-L2-ipt-f}
\begin{split}
&\e^{-2} \normm{ \ip f(t)}_{L^2_{x,v}(\nu)}^2
%=& \e^{-2}\iint_{\O\times\R^3} [\ip f]^2(t)\nu  \dd v
%\dd x \\
%=& \e^{-2} \iint_{\O\times\R^3} [\ip f]^2(0) \nu\dd v
%\dd x
%+ \e^{-2} \int_{0}^{t} \iint_{\O\times\R^3}  \frac{\dd [\ip f]^2(s)}{\dd s} \nu \dd v
%\dd x \dd s  \\
\lesssim   \e^{-2} \iint_{\O\times\R^3} [\ip f_0]^2\nu   + \mathscr{D}_{1}[f](t)
\lesssim  \left[\!\left[f_{0}\right]\!\right]_1^2 +\mathscr{D}_{1}[f](t).
\end{split}
\end{equation}
For $\normm{\ip f }_{L^{6}_{x,v}} $,  interpolation combined with Young's inequality and \eqref{L-t-infty-L2-ipt-f} gives
\begin{equation}
\begin{split}
 \normm{\ip f }_{L^{6}_{x,v}}\leq &
 \big[ \e^{\frac{1}{2}} \normm{ \o^{\frac{1}{2}}  f }_{L^{\infty}_{x,v}} \big]^{\frac{2}{3}}
  \big[ \e^{-1} \normm{ \ip f }_{L^2_{x,v}} \big]^{\frac{1}{3}}
% \leq &\d \e^{\frac{1}{2}} \normm{ \o^{\frac{1}{2}}  f }_{L^{\infty}_{x,v}}
%  +C_\d \e^{-1} \normm{ \ip f}_{L^2_{x,v}}\\
  \leq
 \d \e^{\frac{1}{2}} \normm{ \o f}_{L^{\infty}_{x,v}}
  + \left[\!\left[f_{0}\right]\!\right]_1 +\mathscr{D}_{1}^{\frac{1}{2}}[f](t).
\end{split}
\end{equation}
Moreover, by Lemma \ref{Gamma-L2-Linfty-estimate-1}, interpolation ($L^4\subseteq L^2\cap L^6$) and \eqref{L-t-infty-L2-ipt-f}, we obtain
\begin{equation} \label{0-Gamma-L-t-infty-L2}
\begin{split}
\normm{ \nu^{-\frac{1}{2}} \Gamma(f,f) }_{L^{2}_{x,v}}
% \lesssim   &
%\normm{ \nu^{-\frac{1}{2}} \Gamma(f,\ip f) }_{L^{2}_{x,v}}
%+\normm{ \nu^{-\frac{1}{2}} \Gamma(\P f, \P f) }_{L^{2}_{x,v}} \\
 \lesssim   &
 \big( \e \normm{\omega^{\frac{1}{2}} f }_{L^\infty_{x,v}} \big ) \big(\e^{-1}\normm{\ip f}_{L^2_{x,v}(\nu)}\big) +
 \normm{ \P f }_{L^{4}_{x,v}}^{2}  \\
% \lesssim   &
%    \e^2 \normm{\omega^{\frac{1}{2}} f}_{L^\infty_{x,v}}^2
%+\e^{-2}\normm{\ipt  \tilde{f}}_{L^2_{x,v}(\nu)}^2  +
%  \normm{ \P f }_{L^{2}_{x,v}}^{\frac{3}{2}}
%  \normm{ \P f }_{L^{6}_{x,v}}^{\frac{1}{2}}\\
 %-------------
 \lesssim  &
  \e^2 \normm{\omega f}_{L^\infty_{x,v}}^2 +\left[\!\left[f_{0}\right]\!\right]_1^2 +\mathscr{D}_{1}[f](t)
  +
 \mathscr{E}_{1}^{\frac{3}{2}}[f](t)+
 \d \normm{ \P f  }_{L^{6}_{x,v}},
\end{split}
\end{equation}
where $\d>0$ is a sufficiently small constant.

Combining all the estimates with \eqref{0-Pf-L6} and absorbing the small term
$\d \normm{ \P f }_{L^{6}_{x,v}}$ from \eqref{0-Gamma-L-t-infty-L2}, we arrive at \eqref{0-pf6 - bound3-final}.
\end{proof}
\bigskip

\subsection{Proof of Main Result for the Case $\e \lesssim \a \leq 1$} \label{proof-main-1} \
\medskip

This subsection presents the proof of  Theorem \ref{main-th-1}.
\medskip

\begin{proof}[\textbf{Proof of Theorem \ref{main-th-1}}]  \
 We work with the perturbation formulation \eqref{f-eq} around the global Maxwellian $\mu$.  The proof proceeds in three main steps.
 \medskip

%  To track the $\e$-dependence of functions and equations explicitly in the limit $\e \to %0$, we denote all relevant quantities with a superscript $\e$.
% Thus, let $F^\e=\mu+\e \sqrt{\mu}f^{\e}$ be the solution to \eqref{F - Boltzmann %equation}, where the fluctuation  $f^{\e}$ satisfies
% \begin{equation}\label{f-epsilon-equation}
%\begin{split}
% \e \pt_t f^{\e} + v \cdot \nabla_x f^{\e} + \e^{-1}L f^{\e} = \G(f^{\e},f^{\e})       %\;\;\;\;   &\text{ in } \mathbb{R}^{+}\times\O\times\mathbb{R}^3,\\
% f^{\e}\big|_{\g_{-}}
%= (1-\a)\mathscr{ R}f^{\e} + \a \mathscr{P}_{\g}f^{\e} \;\;\;\;   &\text{ on } %\mathbb{R}^{+}\times\pt\O\times\mathbb{R}^3, \\
% f^{\e} \big|_{t=0}  =   f^{\e}_0(x,v) \;\;\;\;  &\text{ on } \O\times\mathbb{R}^3.
%\end{split}
%\end{equation}
%  The proof proceeds in three main steps: (1) establishing the uniform bound % % % \eqref{0-uniform-bound}; (2) verifying the strong convergence % % % \eqref{tilde-f-strong-convergence}--\eqref{momentums-strong-convergence} and deriving the % limiting INSF system \eqref{INSF-unst}; (3) justifying the boundary conditions % \eqref{Diri-bdy-unst} and \eqref{Navier-bdy-unst}.

\noindent\textbf{Step 1. Global existence and uniform $\e$-independent estimates.}

 We first establish the global a priori estimate \eqref{0-uniform-bound} under the initial condition \eqref{initial-data-total-norm}. Assume that a solution $f$ to \eqref{f-eq} exists on $[0, T]$ for some $0 < T \leq \infty$.

First, applying  Corollary \ref{G - L2 Linfty L3 L6 estimate} and
Proposition \ref{f - L2L3 estimate} in Appendix A with source terms $g=-\e^{-1}L {{f}}+\Gamma(f,f)$ (for  $S_1 f$) and
$g=-\e^{-1} L \pt_t {f}+\Gamma(f,\pt_t f)+\Gamma(\pt_t f, f)$ (for $S_1\pt_t {f}$), we obtain
\begin{equation}\label{0-Gamma-L2-estimate}
\begin{split}
 &\normm{ {\nu}^{-\frac{1}{2}} {\G}({f}, {f})}^2_{L^2_{t,x,v}}
+\normm{ {\nu}^{-\frac{1}{2}} {\G}({{f}},\pt_t {{f}})}^2_{L^2_{t,x,v}}
+\normm{  {\nu}^{-\frac{1}{2}} {\G}(\pt_t {{f}},{{f}})}^2_{L^2_{t,x,v}}
%------------------
\lesssim  \left[\!\left[{f}_{0}\right]\!\right]_{1}^2\normmm{{{f}}}_{1}^{2}(t) + \normmm{{{f}}}_{1}^{4}(t).
\end{split}
\end{equation}

Second, multiplying the estimate \eqref{0-P-f-macro-L2} from  Proposition \ref{0-macro-L2-L6-estimate} and \eqref{0-t-P-f-macro-L2-pt} from Corollary \ref{0-macro-L2-L6-estimate-pt} by a small coefficient $ \eta_1 $ satisfying $0<\eta  \ll \eta_1 \ll \min\{ 1, \frac{\l}{4}\}$ (cf. the definition of $\l_1$ in \eqref{a-e-limit-infty}), and adding the result to the estimates \eqref{f-energy-estimate} and \eqref{ft-energy-estimate} in Proposition \ref{f-ft-Energy-estimate}, we obtain
\begin{equation}\label{0-EfDf - bound1}
\begin{split}
\mathscr{E}_{1}[{f}](t)  + \mathscr{D}_{1}[{f}](t)
\lesssim \left[\!\left[{f}_{0}\right]\!\right]_{1}^2
+ \left[\!\left[{f}_{0}\right]\!\right]_{1}^2\normmm{{f}}_{1}^{2}(t)
+ \normmm{{f}}_{1}^{3}(t)
+ \normmm{{f}}_{1}^{4}(t).
\end{split}\end{equation}

Third, combining Proposition \ref{L-infty-bd-unst} and Lemma \ref{Gamma-L2-Linfty-estimate-1} gives
\begin{equation} \label{0-wfinfty - bound2}
\begin{split}
\e\normm{\omega {f}}_{L^{\infty}_{t,x,v}}^2 + \e^{3}\normm{\omega \pt_t {f}}_{L^{\infty}_{t,x,v}}^{2}
\lesssim & \left[\!\left[{f}_{0}\right]\!\right]_{1}^2
+ \mathscr{E}_{1}[{f}](t) + \mathscr{D}_{1}[{f}](t)+ \normmm{{f}(t)}_{1}^{4}
 + \normm{\P {f}}_{L^{\infty}_{t}L^6_{x,v}}^2.
\end{split}\end{equation}
Applying Corollary \ref{f-L6-estimate-final} yields
\begin{equation} \label{0-pf6 - bound3}
\begin{split}
\normm{ \P {f} }_{L^{\infty}_{t}L^6_{x,v}}^2
\lesssim & \left[\!\left[{f}_{0}\right]\!\right]_{1}^2 + \left[\!\left[{f}_{0}\right]\!\right]_{1}^4+ \mathscr{E}_{1}[{f}](t) + \mathscr{D}_{1}[{f}](t) + \normmm{{f}}_{1}^{4}(t) + \normmm{{f}}_{1}^{6}(t)
 + \d\e\normm{\omega {f}}_{L^{\infty}_{t,x,v}}^2,
\end{split}
\end{equation}
where $\d>0$ is a sufficiently small constant.
Combining \eqref{0-wfinfty - bound2} and \eqref{0-pf6 - bound3} and absorbing $\d\e\normm{\omega {f}}_{L^{\infty}_{t,x,v}}^2$ on the right-hand side of \eqref{0-pf6 - bound3} and $\normm{\P {f}}_{L^{\infty}_{t}L^6_{x,v}}^2$ on the right-hand side of \eqref{0-wfinfty - bound2}, we obtain
\begin{equation}\label{0-infty-macro-bound}
\begin{split}
&\e\normm{\omega {f}}_{L^{\infty}_{t,x,v}}^2 + \e^{3}\normm{\omega \pt_t {f}}_{L^{\infty}_{t,x,v}}^{2}  + \normm{ \P  {f}}_{L^{\infty}_{t}L^6_{x,v}}^2 \\
\le &
\left[\!\left[{f}_{0}\right]\!\right]_{1}^2 + \left[\!\left[{f}_{0}\right]\!\right]_{1}^4+ \mathscr{E}_{1}[{f}](t) + \mathscr{D}_{1}[{f}](t) + \normmm{{f}}_{1}^{4}(t) + \normmm{{f}}_{1}^{6}(t).
\end{split}
\end{equation}

Finally, multiplying \eqref{0-infty-macro-bound} by a small constant, adding it to \eqref{0-EfDf - bound1} and  absorbing small terms, we obtain
\begin{equation}
\normmm{{f}}_{1}^2(t) \lesssim \left[\!\left[{f}_{0}\right]\!\right]_{1}^2 + \normmm{{f}}_{1}^{3}(t) + \normmm{{f}}_{1}^{4}(t)+ \normmm{{f}}_{1}^{6}(t)
\end{equation}
 for any $0 \leq t \leq T$, provided $\left[\!\left[{f}_{0}\right]\!\right]_{1}^2 \leq \d_0$ is sufficiently small.
This establishes the global a priori estimate \eqref{0-uniform-bound}.

The existence of a global solution ${f}$ on $[0,\infty]$ then follows from a standard continuity argument (see, e.g. \cite{Guo2003}); the routine local existence theory is omitted for brevity.
\medskip

\noindent \textbf{Step 2. Derivation of strong convergence \eqref{tilde-f-strong-convergence}--\eqref{momentums-strong-convergence}  and INSF system \eqref{INSF-unst}.}

The uniform bound on  $\normmm{{f}}_{1}(\infty)$ given by \eqref{0-uniform-bound} implies:
\begin{align}
&\sup_{0\leq s\leq \infty}\Big(
 \normm{{f}(s)}_{L^2_{x,v}}
+\normm{\pt_t {f}(s)}_{L^2_{x,v}}
+\normm{\P {f}(s)}_{L^6_{x,v}}
\Big)\leq C\d_0, \label{0-tilde-f-uniform-bound} \\
%----------------------
& \int_{0}^{\infty} \Big( \normm{\P {f}(s)}_{L^2_{x,v}}^{2}
+  \int_{0}^{t} \normm{\pt_t \P {f}(s)}_{L^2_{x,v}}^{2} \Big)\dd s
\leq C\d_0,\label{0-tilde-f-L2-bound} \\
%----------------
&\int_{0}^{\infty}
\Big(
\normm{\ip {f}(s)}_{L^2_{x,v}(\nu)}^{2}
+\normm{\ip \pt_t {f}(s)}_{L^2_{x,v}(\nu)}^{2}
\Big) \dd s  \to  0\;\;\; \hbox{ as } \; \e\to  0.\label{0-ip-to-0}
\end{align}
Hence, there exists $f^{*}\in L^\infty\left( \mathbb{R}^+; L^2(\O\times \mathbb{R}^3)\right)$ such that, up to a subsequence,
\begin{align}
f \to  f^{*} \;\; &\text{~~weakly}\!-\!* ~\text{in}~ L^\infty\left( \mathbb{R}^+; L^2(\O\times \mathbb{R}^3)\right) \; \hbox{ as } \; \e\to  0.\label{0-f-es-L2-lim}
\end{align}
%This further implies
%\begin{equation*}
% {L} {f}\to Lf^{*} \;\; \hbox{ in the sense of distributions}  \; \hbox{ as } \; %\e\to  0.
%\end{equation*}
On the other hand, \eqref{0-ip-to-0} gives
\begin{equation*}
L {f} \to 0 \ \hbox{ strongly in } L^2\left( \mathbb{R}^+\times \O\times \mathbb{R}^3\right) \; \hbox{ as } \; \e\to  0.
\end{equation*}
By the uniqueness of distribution limits, we conclude $Lf^{*}=0$. Hence, there exist functions $\varrho_{f^{*}},u_{f^{*}},\vartheta_{f^{*}}\in L^\infty\left( \mathbb{R}^+; L^2( \O)\right)$ such that
\begin{equation} \label{0-Ltildef and tildef}
 f^{*} = \Big(\varrho_{f^{*}} + u_{f^{*}} \cdot v + \vartheta_{f^{*}} \frac{\norm{v}^2-3}{2}\Big)\sqrt{\mu}.
\end{equation}

Furthermore, the uniform boundedness of $\normmm{f}_{1}(\infty)$ together with \eqref{0-Gamma-L2-estimate} implies
\begin{equation}\label{0-L2-bound-average}
\begin{split}
 \pt_t {f}, \; \e^{-1} \nu^{-\frac{1}{2}} L{f}, \;  \nu^{-\frac{1}{2}} \G({f}, {f}) \in L^2(\mathbb{R}^{+}\times\O\times \mathbb{R}^3).
\end{split}
\end{equation}
Consequently, equation \eqref{f-eq} indicates that
$ \nu^{-\frac{1}{2}} v \cdot \nabla_x {f}\in L^2\left( \mathbb{R}^+; L^2(\O\times \mathbb{R}^3)\right)$ and hence admits a weak limit.
On the other hand, \eqref{0-tilde-f-uniform-bound}  implies
\begin{equation*}
\nu^{-\frac{1}{2}} v \cdot \nabla_x {f} \to \nu^{-\frac{1}{2}} v \cdot \nabla_x f^{*}\;\; \hbox{ in the sense of distributions} \; \hbox{ as } \; \e\to  0.
\end{equation*}
By the uniqueness of distribution limits, we obtain
\begin{equation}\label{0-nabla-f-convergence}
 \nu^{-\frac{1}{2}} v \cdot \nabla_x {f} \to \nu^{-\frac{1}{2}}  v \cdot \nabla_x f^{*} \;\; \text{ weakly} ~\text{in}\; L^2\left( \mathbb{R}^+; L^2(\O\times \mathbb{R}^3)\right) \; \hbox{ as } \; \e\to  0.
\end{equation}
Using the linear independence of $\nu^{-\frac{1}{2}}v \big \{ 1,v,v\otimes v,|v|^{2}, v|v|^2 \big  \} \sqrt{\mu}$ and \eqref{0-nabla-f-convergence}, we conclude that
$$
 \varrho_{f^{*}},u_{f^{*}},\vartheta_{f^{*}} \in L^2\left( \mathbb{R}^+;  H^{1}(\O)\right).
$$

We now prove the strong convergence stated in  \eqref{tilde-f-strong-convergence}--\eqref{momentums-strong-convergence}.
First, we claim that
\begin{equation}\label{tilde-f-epsilon-strong-convergence}
{f} \to  f^* \;\; \hbox{ strongly in } L^2\big(\mathbb{R}^{+};L^2(\O\times \mathbb{R}^3)\big) \; \hbox{ as } \; \e\to  0.
\end{equation}
To prove this claim, by virtue of \eqref{0-L2-bound-average}, we truncate ${f}$ as in \eqref{fdelta- Defintion} to obtain ${f}_{\d}$. Then we  apply the extension Lemma 3.6 from \cite{Esposito2017} to define $\overline{{f}_{\d}}$ on $\mathbb{R}\times \mathbb{R}^3\times \mathbb{R}^3$, and invoke the $L^2$ averaging lemma (cf. Proposition 3.3.1 in \cite{SaintRaymond2009}) to obtain
\begin{equation}\label{f-delta-L2-H1half-bound}
\begin{split}
\normm{\int_{\mathbb{R}^3} \nu^{-\frac{1}{2}} \overline{{f}_{\d}} \psi \dd v }_{L^2_t\big(\mathbb{R};H^{\frac{1}{2}}_x(\mathbb{R}^3)\big)}\leq C,
\end{split}
\end{equation}
where $\psi\in L^\infty(\mathbb{R}^3)$ represents any compactly supported test function, and the constant $C$
%depends only on $\psi$, the upper bounds of $\normm{{f}_{0}}_{L^2_{x,v}}$,  %$\e^{-1}\normm{\nu^{-\frac{1}{2}}L{f}}_{L^2_{x,v}}$ and $\normm{\nu^{-\frac{1}{2}} %\G({f}, {f})}_{L^2_{x,v}}$, and
is independent of $\e$.
By compact embedding, up to a subsequence, we have
\begin{equation}\label{f-delta-L2-strong-converge}
\begin{split}
\int_{\mathbb{R}^3} \nu^{-\frac{1}{2}} {f}_{\d} \psi \dd v  \;\; & \hbox{ converges strongly in }L^2_{\text{loc}}\big(\mathbb{R}^{+};L^2_x(\O)\big) \; \hbox{ as } \; \e\to  0.
\end{split}
\end{equation}
 Using \eqref{f-delta-L2-strong-converge} and a decomposition similar to \eqref{int-f-delta-decomposition}, for each  $i=0,1, \cdots, 4$, we deduce
\begin{equation}
\begin{split}
\int_{\mathbb{R}^{3}} {f}_{\d}  \chi_{i} (v) \mathrm{d} v
%=& \mathbf{1}_{t\geq 0} \Big\{\overline{a^{\e}_{i}}(t,x)  + O(\delta) \sum_{j=0}^{4} %|\overline{a^{\e}_{j}}
%(t,x)| + O_{\delta}(1) \int_{\mathbb{R}^{3}} | \ipt
%\tilde{f}^{\e} | \bar{\chi}_{i} (v) \mathrm{d} v\Big\} \\
%&+ \mathbf{1} _{ t \leq 0} \chi(t) \int_{\mathbb{R}^{3}} \tilde{f}^{\e}_{0}
%\bar{\chi}_{i} (v) \mathrm{d} v \quad \text{ for each } i=0,1, \cdots, 4\\
=& \mathbf{1}_{t\geq 0} \Big\{ {
a_{i}}  + O(\delta) \sum_{j=0}^{4} | {a_{j}}
| \Big\} \\
& +\mathbf{1}_{t\geq 0}  \int_{\mathbb{R}^{3}} \Big[ 1- \chi( \frac{n(x) \cdot v}{\delta} ) \chi(
\frac{\xi(x)}{\delta}) \Big]
%--------
  \Big[ 1- \chi \big(\frac{|v|}{2\d} \big)  \Big]
 %--------
 \chi( \delta |v|)  \ip
f  \chi_{i}(v) \mathrm{d} v \\
%--------------------------------------
&+ \mathbf{1} _{ t \leq 0} \chi(t) \int_{\mathbb{R}^{3}} \Big[ 1- \chi( \frac{n(x) \cdot v}{\delta} ) \chi(
\frac{\xi(x)}{\delta}) \Big]
%--------
 \Big[ 1- \chi \big(\frac{|v|}{2\d} \big)  \Big]
 %--------
 \chi( \delta |v|)  f_{0}
\chi_{i}(v) \mathrm{d} v.
\end{split}
\end{equation}
Here and in what follows,  we use the temporary notations
$$
{a}_{0}={a},\;  {a}_{i}={b}_{i}\, (i=1,2,3),\; {a}_{4}={c};\;\; \;
{a}_{0}^*={\varrho_{f^*}},\;\;  {a}_{i}^*={u_{f^*}}\, (i=1,2,3),\; \;{a}_{4}^*={\vartheta_{f^*}}.
$$
From \eqref{initial-converge-condition} and \eqref{0-ip-to-0}, we obtain
 for each $i=0,1, \cdots, 4$,
\begin{equation}\label{a-i-strong-converges 1}
\begin{split}
a_{i}  + O(\delta) \sum_{j=0}^{4} |a_{j}
|   \text{ converges strongly in } L^2_{\text{loc}}\big(\mathbb{R}^{+};L^2(\O)\big)
 \; \hbox{ as } \; \e\to  0.
\end{split}
\end{equation}
  Combining this with the weak convergence \eqref{0-f-es-L2-lim}, we obtain
  for each $i=0,1, \cdots, 4$:
\begin{equation*}\label{a-i-plus-strong-converges}
\begin{split}
a_{i}  + O(\delta) \sum_{j=0}^{4} |a_{j}
| \to a^{*}_{i}  + O(\delta) \sum_{j=0}^{4} |a^{*}_{j}
|  \; & \hbox{ strongly in }L^2_{\text{loc}}\big(\mathbb{R}^{+};L^2(\O)\big)
 \; \hbox{ as } \; \e\to  0.
\end{split}
\end{equation*}
 Consequently,
\begin{equation*}\label{a-i-strong-converges 2}
\begin{split}
\big(1-5 O(\d)\big) \sum_{i=0}^{4} \big  \|   a_{i}-a^{*}_{i}   \big  \|_{L^2_{t,x}}
\leq  \sum_{i=0}^{4} \Big \| a_{i}-a^{*}_{i}+ O(\d)\sum_{j=0}^{4}\big ( |a_{j}| -|a^{*}_{j}| \big) \Big \|_{L^2_{t,x}}.
\end{split}
\end{equation*}
Since $\d >0$ is sufficiently small, we conclude that for each $i=0,1, \cdots, 4$,
\begin{equation}\label{a-i-strong-converges 3}
\begin{split}
a_{i} \to  a^{*}_{i}  \;\; & \hbox{ strongly in }L^2_{\text{loc}}\big(\mathbb{R}^{+};L^2(\O)\big) \; \hbox{ as } \; \e\to  0.
\end{split}
\end{equation}
This indicates
\begin{equation}\label{Pf-strong-converges}
\begin{split}
\P {f} \to  \P f^{*}   \;\; & \hbox{ strongly in }L^2_{\text{loc}}\big(\mathbb{R}^{+};L^2(\O\times \R^3)\big) \; \hbox{ as } \; \e\to  0.
\end{split}
\end{equation}
Together with the strong convergence of $\ip {f}$ and the integrability of $\P f$ from \eqref{0-uniform-bound} as well as the integrability of $\P f^*$ from Lemma \ref{lem:global-estimate-NS}, this proves the claim \eqref{tilde-f-epsilon-strong-convergence}.
Together with  \eqref{0-ip-to-0}, this yields the claim \eqref{tilde-f-epsilon-strong-convergence}. Moreover, \eqref{Pf-strong-converges} gives
\begin{equation}\label{Pf-momentum-strong-convergence}
\begin{split}
&\int_{\mathbb{R}^3} {f} \sqrt{\mu}\Big[1,v, \frac{|v|^2-3}{2}\Big] \dd v  \to  \big( \varrho_{f^{*}},  u_{f^{*}}, \vartheta_{f^{*}}\big)
\;\;  \hbox{ strongly in }L^2\big(\mathbb{R}^{+};L^2(\O)\big) \text { as } \e\to  0.
\end{split}
\end{equation}
The strong convergence properties \eqref{tilde-f-strong-convergence}--\eqref{momentums-strong-convergence} now follow
readily.

Using \eqref{tilde-f-epsilon-strong-convergence}, we take the weak limit of equation \eqref{f-eq} in $L^2\left( \mathbb{R}^+; L^2(\O\times \mathbb{R}^3)\right)$ to obtain
\begin{equation} \label{0-tildeL tildef - epsilon goes 0}
\lim_{\e \to 0} \e^{-1} \nu^{-\frac{1}{2}} L {f} = \nu^{-\frac{1}{2}}\G(f^{*},f^{*}) - \nu^{-\frac{1}{2}} (v \cdot \nabla_x f^{*}) \;\;\text{in the weak sense}.
\end{equation}
Multiplying \eqref{f-eq} by $\sqrt{\mu}$ and $v\sqrt{\mu}$, and integrating over $\R^3$, we have
\begin{equation}
\nabla_{x} \cdot u_{f^{*}}  = 0, \quad \nabla_x (\rho_{f^{*}} + \vartheta_{f^{*}})  = 0.
\end{equation}
Multiplying  \eqref{f-eq}  by $\e^{-1} \frac{\norm{v}^2-5}{2} \sqrt{\mu}$, integrating over $\R^3$ and following the procedure in \cite{Bardos1991}, we obtain
\begin{equation} \label{ener-lim-st}
\begin{split}
-\pt_t \th_{f^{*}}  =
&\; - \lim_{\e \to 0} \int_{\R^3}
 \frac{\norm{v}^2-5}{2}\sqrt{\mu} \pt_t {f} \dd v
= \lim_{\e \to 0}  \frac{1}{\e}  \int_{\R^3} \frac{\norm{v}^2-5}{2}\sqrt{\mu}
( v \cdot \nabla_x {f} ) \dd v \\
= & \; \lim_{\e \to 0} \frac{1}{\e} \nabla_x \cdot \int_{\R^3} L^{-1}\big(\frac{\norm{v}^2-5}{2} v\sqrt{\mu}\big)  L {f}\dd v \\
=& \; \nabla_{x} \cdot \int_{\R^{3}} {L}^{-1}\big(\frac{\norm{v}^2-5}{2} v\sqrt{{\mu}}\big) \big(\G(f^{*},f^{*})
- (v \cdot \nabla_x f^{*})\big)\dd v \\
= & \; \nabla_{x} \cdot \big (  \frac{5}{2} \kappa \nabla_x  \vartheta_{f^{*}} - \frac{5}{2}u_{f^{*}} \vartheta_{f^{*}}   \big ),
\end{split}
\end{equation}
where we have used \eqref{0-tildeL tildef - epsilon goes 0} and the decay  property of ${L}^{-1}\big(\frac{\norm{v}^2-5}{2} v\sqrt{{\mu}}\big) $. Here the thermal conductivity
is defined as
\begin{equation}\label{kappa-def}
\kappa:=\frac{2}{5}\int_{\R^{3}}\big(\frac{\norm{v}^2-5}{2} v\sqrt{{\mu}}\big) {L}^{-1}\big(\frac{\norm{v}^2-5}{2} v\sqrt{{\mu}}\big)\dd v.
\end{equation}
Similarly, multiplying  \eqref{f-eq}  by $\e^{-1} v\sqrt{{\mu}}$ and integrating over $\R^3$, we obtain
\begin{equation} \label{mome-lim-st}
\begin{split}
-\pt_t u_{f^{*}} =
&\; - \lim_{\e \to 0} \int_{\R^3}
v \sqrt{{\mu}} \pt_t {f}\dd v
%------------
=  \lim_{\e \to 0}  \frac{1}{\e}  \int_{\R^3} v \sqrt{{\mu}} \big(v \cdot \nabla_x {f}\big) \dd v \\
%------------
= &\; \lim_{\e \to 0}  \frac{1}{\e}  \nabla_x \cdot \int_{\R^3} \big[ L^{-1}\big(v\big(v \otimes v- \frac{\norm{v}^2}{3}\mathbb{I}\big)\sqrt{{\mu}}\big) L {f} + \frac{\norm{v}^2}{3}\sqrt{{\mu}}  {f}\
\big] \dd v \\
%------------
= &\;\nabla_{x} \cdot \int_{\R^{3}} {L}^{-1}\big(\big(v \otimes v- \frac{\norm{v}^2}{3}\mathbb{I}\big)\sqrt{{\mu}}\big)  \big(\G(f^{*},f^{*}) - v \cdot \nabla_x f^{*}\big) \dd v + \nabla_{x} p_{f^{*}} \\
%------------
= &\; \nabla_x\cdot \big [ 2u_{f^{*}}\otimes u_{f^{*}}  -  \frac{2}{3}\norm{u_{f^{*}}}^2\mathbb{I}
  -  \s \big(  \nabla_x u_{f^{*}}+{(\nabla_x u_{f^{*}})}^{\mathrm{T}}\big)\big ] +\nabla_x p_{f^{*}}.
\end{split}
\end{equation}
%where we have used \eqref{0-tildeL tildef - epsilon goes 0} and the decay property of ${L}^{-1}\big((v \otimes v- \frac{\norm{v}^2}{3}\mathbb{I})\sqrt{{\mu}}\big)$.
Here the viscosity is defined as
\begin{equation}\label{sigma-def}
\s:=\frac{1}{10}\int_{\R^{3}}\big[\big(v \otimes v- \frac{\norm{v}^2}{3}\mathbb{I}\big)\sqrt{{\mu}}\big] : {L}^{-1}\big[\big(v \otimes v- \frac{\norm{v}^2}{3}\mathbb{I}\big)\sqrt{{\mu}}\big]  \dd v,
\end{equation}
and we have used the notation
$$
 p_{f^{*}} : = \lim_{\e \to 0} \frac{1}{\e} \int_{\R^3} \frac{\norm{v}^2}{3}\sqrt{\mu} {f}\dd v.
$$
Hence, $(\rho_{f^{*}}, u_{f^{*}}, \th_{f^{*}})$ satisfies the INSF system \eqref{INSF-unst} in the weak sense.
\medskip

%\begin{equation}\label{u-vartheta-equation}
%\begin{split}
%&\pt_t {u} + {u} \cdot\nabla_x {u} + u \cdot\nabla_x {u} +\nabla_x {p}
%= \s \Delta {u}, \ \ \  \nabla_x\cdot {u} = 0 \ \ \ \    \text{ in } \mathbb{R}^{+}\times %\O, \\
%&\pt_t { \vartheta} + {u} \cdot\nabla_x { \vartheta}+ u\cdot\nabla_x { \vartheta}
%= \kappa  \Delta { \vartheta}, \ \ \  \nabla_x ({\rho} + { \vartheta} ) = 0 \ \ \ \
%\text{ in } \mathbb{R}^{+}\times \O.
%\end{split}
%\end{equation}

\noindent\textbf{Step 3. Derivation of the boundary conditions   \eqref{Diri-bdy-unst} and \eqref{Navier-bdy-unst}.}

  Consider the identity
\begin{equation*}
\begin{split}
\iint_{\partial\Omega\times\mathbb{R}^3}\nu^{-\frac{1}{2}} \phi {f}|_{\partial\Omega}[n\cdot v]\dd v \dd S_x
&=\iint_{\Omega\times\mathbb{R}^3}\nu^{-\frac{1}{2}}(v\cdot\nabla_x\phi) {f}
+\iint_{\Omega\times\mathbb{R}^3}\nu^{-\frac{1}{2}}(v\cdot\nabla_x {f})\phi,
\end{split}
\end{equation*}
where $\phi(x, v)$ is test function satisfying $\phi(\cdot, v)\in C^\infty(\bar{\Omega})$ and $\phi(x, \cdot)\in
C^\infty_0(\mathbb{R}^3)$.  Using the weak convergence of  ${f}$ and $v\cdot\nabla_x {f}$, we obtain
\begin{equation}\label{0-g-gamma-distr-lim}
\begin{split}
\nu^{-\frac{1}{2}}{f}|_{\partial\Omega} \to \nu^{-\frac{1}{2}} f^{*} |_{\partial\Omega} \quad  \hbox{ in the sense of distributions} \; \hbox{ as } \; \e\to  0.
\end{split}
\end{equation}
The uniform bound on $\normmm{{f}}_{1}(\infty)$ in \eqref{0-uniform-bound} implies
\begin{equation}\label{0-f-bdy-norm-bd}
\begin{split}
\left(\frac{\a}{\e}\right)^{\frac{1}{2}} \norm{(1-\mathscr{P}_{\g}){f}}_{L^2_tL^2_{\gamma_+}}
+ \norm{\mathscr{P}_\g {f}}_{L^2_tL^2_{\gamma_+}}
 \   \hbox{ is uniformly bounded}.
\end{split}
\end{equation}
On the other hand, by \eqref{a-e-limit-infty},  the quantity $\norm{f^\e}_{L^2_tL^2_{\gamma_+}}$
is uniformly bounded, and hence, up to a subsequence, has a weak limit in
$ L^2(\mathbb{R}_+\times\dd \gamma)$. From \eqref{0-g-gamma-distr-lim} and the uniqueness of distribution limits, we conclude that $\nu^{-\frac{1}{2}}f^{*}\big|_{\partial\Omega}\in L^2(\mathbb{R}_+\times\dd \gamma)$ and
\begin{equation}\label{f-epsilon-bd-con}
\begin{split}
\;\;
 \nu^{-\frac{1}{2}} {f}\big|_{\partial\Omega}    \to
 \nu^{-\frac{1}{2}}f^{*}\big|_{\partial\Omega} \;\;  \hbox{weakly in } L^2(\mathbb{R}^+\times\dd \gamma)  \text{ as } \e\to 0.
\end{split}
\end{equation}

We now define
\begin{equation*}
  \begin{split}
          \langle g\rangle_{\partial\Omega}:=\sqrt{2\pi}\int_{v\cdot n>0}g\big|_{\partial\Omega}\sqrt{\mu}[n\cdot v]\dd v.
  \end{split}
\end{equation*}
From \eqref{f-epsilon-bd-con} and the fact that $\langle {f}\rangle_{\partial\Omega}$ is independent of $v$, we have
\begin{equation}\label{0-ges-bd-con}
\begin{split}
\langle {f}\rangle_{\partial\Omega}\to  \langle f^{*}\rangle_{\pt\O}\;\;
\hbox{ weakly in } L^2(\mathbb{R}^+\times\dd \g)  \text{ as } \e\to 0.
\end{split}
\end{equation}
Combining this with \eqref{f-epsilon-bd-con} gives
\begin{equation}\label{0-f-epsilon-<f-epsilon>-converge}
\begin{split}
 \nu^{-\frac{1}{2}} \big (  {f}|_{\partial\Omega}-\sqrt{\mu}\langle {f}\rangle_{\partial\Omega}\big)
  \to
 \nu^{-\frac{1}{2}}\big (  f^{*}|_{\partial\Omega}-\sqrt{\mu}\langle f^{*}\rangle_{\partial\Omega} \big)\;
 \hbox{ weakly in }  L^2(\mathbb{R}^+\times\dd \g) \; \text{ as } \e\to 0.
\end{split}
\end{equation}

We now derive the boundary conditions \eqref{Diri-bdy-unst} and \eqref{Navier-bdy-unst} according to the limit value $\lambda$ defined in \eqref{a-e-limit-infty}.
\medskip

\noindent\textbf{Step 3.1. Dirichlet boundary condition \eqref{Diri-bdy-unst} for $\l=\infty$.}

 In this case, we can take the limit in the Maxwell boundary condition directly and show strong convergence. The uniform  boundedness \eqref{0-f-bdy-norm-bd} implies
\begin{equation}\label{f-epsilon-conv-on-bdy}
\begin{split}
{f}|_{\pt\O}-\sqrt{\mu}\langle {f}\rangle_{\pt\O}=(1-\mathscr{P}_{\g}){f}\to  0 \ \  \text{ strongly in }
L^2(\mathbb{R}^{+}\times\dd \g)\; \text{ as } \e\to 0.
\end{split}
\end{equation}
Combining \eqref{0-f-epsilon-<f-epsilon>-converge} and \eqref{f-epsilon-conv-on-bdy}, we obtain
$$
 \nu^{-\frac{1}{2}}\big( f^{*}|_{\partial\Omega}-\sqrt{\mu}\langle f^{*}\rangle_{\partial\Omega} \big)=0,
$$
 which, together with \eqref{0-Ltildef and tildef}, yields the Dirichlet boundary condition \eqref{Diri-bdy-unst}:
\begin{equation}\label{Diri-bd-g}
\begin{split}
u_{f^{*}}|_{\partial\Omega}=0,\quad \theta_{f^{*}}|_{\partial\Omega}=0.
\end{split}
\end{equation}

\noindent\textbf{Step 3.2. Navier boundary condition \eqref{Navier-bdy-unst}  for $ \l \in ( 0,+\infty )$.}

  By \eqref{0-f-bdy-norm-bd}, we take the weak limit in the Maxwell boundary condition in \eqref{f-eq} to obtain
$$
 f^{*}|_{\g^-}=\mathscr{R}(f^{*}|_{\g^+}).
$$
 This, together with \eqref{0-Ltildef and tildef}, implies the zero mass flux condition
$$
   n\cdot u_{f^{*}}\big|_{\partial\Omega}=0.
$$

To verify the Navier boundary condition, we pass to the limit in the weak formulation of  \eqref{f-eq} and show that the moments $u_{f^{*}}$ and $\theta_{f^{*}}$ satisfy the weak form of the INSF system.  To this  end, we take a test
function $\phi\in C^{\infty}(\bar{\Omega})$ and a divergence-free test vector field $\vec{\omega}\in C^{\infty}(\bar{\Omega})$ with $n \cdot \vec{\omega}|_{\partial\Omega}=0$. Multiplying \eqref{f-eq} by
$\e^{-1}\frac{|v|^2-5}{2}\sqrt{\mu} \phi$ and $\e^{-1}(v\cdot\vec{\omega}) \sqrt{\mu} $, respectively, integrating over $[t_{1},t_{2}] \times \Omega \times \R^{3}$
and passing to the weak limit in $L^2(\O\times \mathbb{R}^3)$, we obtain
\begin{align}
& \lim_{\e\to 0} \int_{t_1}^{t_2} \iint_{\O \times \R^{3}} \pt_t {f}  \frac{|v|^2-5}{2} \sqrt{\mu} \phi
-\lim_{\e\to 0}  \frac{1}{\e} \int_{t_1}^{t_2} \int_{\O}  \big \langle v\frac{|v|^2-5}{2}\sqrt{\mu} ,
\; {f} \big \rangle \cdot \nabla _{x} \phi  \nonumber  \\
&  + \lim_{\e\to 0}   \frac{1}{\e} \int_{t_1}^{t_2}  \iint_{\pt \O \times \R^{3}} {f} \frac{|v|^2-5}{2}
\sqrt{\mu}  \phi [n\cdot v] = 0, \label{mome-lim-unst} \\
%--------------------
&
\lim_{\e\to 0} \int_{t_1}^{t_2}\iint_{\O \times \R^{3}} \pt_t {f}(v\cdot \vec{\omega})  \sqrt{\mu}
- \lim_{\e\to 0} \frac{1}{\e}  \int_{t_1}^{t_2} \int_{\O} \big \langle ( v\otimes v-\frac{|v|^2}3 \mathbb{I} ) \sqrt{\mu} ,
\, {f}\big \rangle :\nabla_x \vec{\omega} \nonumber  \\
 &  + \lim_{\e\to 0} \frac{1}{\e} \int_{t_1}^{t_2}  \iint_{\pt\O\times\mathbb{R}^3}( v\cdot \vec{\omega} ) \sqrt{\mu} {f}[n\cdot v] = 0. \label{ener-lim-unst}
\end{align}
It follows from \eqref{ener-lim-st}  and \eqref{mome-lim-st}  that
\begin{equation}\label{mome-ener-lim-st}
\begin{split}
&\lim_{\e\to  0}\frac{1}{\e}\big \langle v\frac{|v|^2-5}{2}\sqrt{\mu} ,\; {f}\big \rangle
= \frac{5}{2}\kappa\nabla_x\theta_{f^{*}}-\frac{5}{2}u_{f^{*}}\theta_{f^{*}},\\
&\lim_{\e\to  0}\frac{1}{\e}\big \langle (v\otimes v-\frac{|v|^2}3 \mathbb{I})\sqrt{\mu} ,\; {f}\big \rangle=
2u_{f^{*}}\otimes u_{f^{*}}-\frac{2}{3}|u_{f^{*}}|^2\mathbb{I}-\nu \big[  \nabla_x u_{f^{*}} + ({\nabla_x u_{f^{*}}})^\mathrm{T}  \big]
\end{split}
\end{equation}
in the weak sense.

For the boundary term in \eqref{mome-lim-unst}, using \eqref{f-eq} , \eqref{0-Ltildef and tildef} and \eqref{0-f-epsilon-<f-epsilon>-converge} and the change of variables $v\mapsto R_x v$ on $\g_{-}$, we obtain
\begin{equation}\label{theta-bdy-lim-st}
\begin{split}
\lim_{\e\to  0} \frac{1}{\e} \iint_{\partial\Omega \times \R^{3}} {f}\frac{|v|^2-5}{2}\sqrt{\mu}\phi [n\cdot v]\dd v \dd S_x
%--------------------------------------
=&\lim_{\e\to  0} \frac{\alpha}{\e}    \int_{\gamma^+}   \frac{|v|^2-5}{2}\sqrt{\mu}\phi
\big[{f}|_{\partial\Omega}  -  \sqrt{\mu}\langle {f}\rangle_{\partial\Omega} \big]  \dd \g\\
%----------------------
=& \l  \sqrt{2\pi}  \int_{\gamma^+}    \frac{|v|^2-5}{2}\sqrt{\mu}\phi   \big[  {f^{*}}|_{\partial\Omega}
   -  \sqrt{\mu}\langle {f^{*}}\rangle_{\partial\Omega}\big] \dd \g  \\
%----------------------
=& 2\l \int_{\partial\Omega}  \theta_{f^{*}}\phi \dd S_x.
\end{split}
\end{equation}
For the boundary term in \eqref{ener-lim-unst}, using $n\cdot \vec{\omega}|_{\pt\O}=0$ and a similar computation,  we obtain
\begin{equation}\label{u-bdy-lim-st}
\begin{split}
 \lim_{\e\to 0} \frac{1}{\e}   \iint_{\pt \O \times\mathbb{R}^3}  (v\cdot \vec{\omega}) \sqrt{\mu}  {f} [n\cdot v]  \dd v \dd S_x
%--------------------------------------
=&\lim_{\e\to  0} \frac{\a}{\e}  \int_{\g^+}(v\cdot \vec{\omega}) \sqrt{\mu}
\big[   {f}|_{\pt\O}- \sqrt{\mu}\langle {f}\rangle_{\pt\O}   \big] \dd \g\\
%---------------------------------------
=&\l\sqrt{2\pi}\int_{\g^+}  (v\cdot \vec{\omega})\sqrt{\mu}   \big[  {f^{*}}|_{\pt\O}
  -  \sqrt{\mu}\langle {f^{*}}\rangle_{\pt\O} \big] \dd \g \\
=&\l\int_{\pt\O}\vec{\omega}\cdot u_{f^{*}} \dd S_x.
\end{split}
\end{equation}
Thus, \eqref{mome-lim-unst} and \eqref{ener-lim-unst} become
\begin{align}
&   \frac{5}{2} \int_{\O}  \left[ \theta_{f^{*}}(t_2) - \theta_{f^{*}}(t_1)\right  ] \phi\dd x
+  2 \l \int_{t_1}^{t_2} \int_{\pt \O} \theta_{f^{*}}  \phi \dd x\dd s \nonumber
\\
&  -  \int_{t_1}^{t_2} \int_{\O}
 \left ( u_{f^{*}} \theta_{f^{*}} - \kappa \nabla_x \theta_{f^{*}} \right ) \cdot \nabla_x \phi\dd x\dd s = 0, \label{theta-star-weak-equation} \\
%------------------
 & \int_{\O} \left[ u_{f^{*}}(t_2) - u_{f^{*}}(t_1)  \right ] \cdot \vec{\omega}\dd x
+\l \int_{t_1}^{t_2}\int_{\pt \O} u_{f^{*}}\cdot \vec{\omega}\dd x\dd s \nonumber
\\
 &  - \int_{t_1}^{t_2} \int_{\O} \big [ u_{f^{*}} \otimes u_{f^{*}}
     - \s \big ( \nabla_x u_{f^{*}} + {(\nabla_x u_{f^{*}})}^{\mathrm{T}}  \big  ) \big  ]: \nabla_x \vec{\omega} \dd x\dd s = 0. \label{u-star-weak-equation}
\end{align}
 The equations \eqref{theta-star-weak-equation} and \eqref{u-star-weak-equation} constitute  the weak formulation of the INSF system with Navier boundary condition \eqref{Navier-bdy-unst}, satisfied by $\rho_{f^{*}}, u_{f^{*}}$ and $\th_{f^{*}}$.

 Finally, Lemma \ref{A-lemma-NS} in Appendix B guarantees the uniqueness of weak solutions to the INSF system \eqref{INSF-unst} with either Dirichlet boundary condition \eqref{Diri-bdy-unst} or the Navier boundary condition \eqref{Navier-bdy-unst} in the setting of Theorem \ref{main-th-1}. Consequently, all weak limits points coincide with the unique solution to the INSF system.

 This completes the proof of Theorem \ref{main-th-1}.
\end{proof}
\bigskip

\section{Strong Limit for the Case $0 \leq \a \ll \e$}
\medskip

This section investigates the perturbation equation \eqref{tildef - Boltzmanneq ch1} and gives the proof of Theorem~\ref{main-th-2}. The proof relies on Proposition \ref{Psi - L2 and L6 estimate}, which is established first.
\medskip

For clarity and to maintain correspondence with the respective  unknown functions $f$ and $\tilde{f}$, we keep the distinct notations $f_0$ and $\tilde{f}_0$ throughout, although they are equal at the initial time (see \eqref{tilde-f0-equals-f0-t0}).

\subsection{Construction of the Rotating Maxwellian}\
\medskip

In this subsection, we construct the rotating Maxwellian $\tilde{\mu}$ introduced in \eqref{tilde-mu-def} by deriving the ordinary differential equations that govern its component functions $\c$ and $\theta$.
\medskip

We begin with the following Taylor expansion with remainder.

\begin{lemma} \label{approximate - lemma}
Let $h(v,\th,\c) : \R^3 \times [-\d,\d] \times [-\d,\d]^3 \to \R$ be a $C^{\infty}$ function. Define the $n$th-order Taylor expansion of $h$ with respect to $\th$ and $\c$ by
\begin{equation*}
h^{n}(v,\th,\c) := \sum_{\a+\norm{\b} \le n} \frac{1}{\a ! \b !}\th^{\a} \c^{\b}
\frac{\pt^{\a + \norm{\b}}}{
(\pt {\th})^{\a} (\pt {\c})^{\b}} h(v,0,0).
\end{equation*}
Then, the following estimate holds:
\begin{equation*}
\norm{h(v,\th,\c) - h^{n}(v,\th,\c)} \lesssim \e^{n+1} \Big(\frac{\norm{\th}}{\e} + \frac{\norm{\c}}{\e} \Big)^{n+1} \sup_{(\xi,\omega) \in [-\d,\d]^{4}} \norm{\nabla_{\th,\c}^{n+1}h(v,\xi,\omega)}.
\end{equation*}
\end{lemma}

\begin{proof}[\textbf{Proof.}] \
This follows directly from Taylor's theorem with remainder.
\end{proof}
\medskip

For each $n \in \mathbb{N}$, we define the sets of higher-order terms as
\begin{equation} \label{High order set - Definition}
\begin{split}
\mathfrak{H}_{n} :=&\Big\{ \mathfrak{h}_{n}(s) \in \R : \norm{\mathfrak{h}_{n}(s)} \lesssim \Big[\frac{\norm{\th (s)}}{\e} + \frac{\norm{w (s)}}{\e}\Big]^{n}\Big\}, \\
\mathfrak{H}_{n,t} :=&\Big\{ \mathfrak{h}_{n}(s) \in \R : \norm{\mathfrak{h}_{n}(s)} \lesssim \Big[\frac{\norm{\th (s)}}{\e} + \frac{\norm{w (s)}}{\e} +\frac{\norm{\pt_t \th (s)}}{\e} + \frac{\norm{\pt_t w (s)}}{\e}  + \Big(\frac{\a}{\e}\Big)^{\frac{1}{2}} \norm{\tilde{f}(s)}_{L^2_{\g_{+}}}\Big]^{n}\Big\},
\end{split}
\end{equation}
where $\th(s) = T(s) -1$ and $w(s)$ will be determined in Lemma \ref{tildemu - existence}.
Under the a priori assumption \eqref{theta-u-smallness-assumption}, we have
\begin{equation}\label{hn-hm}
\begin{split}
  \mathfrak{H}_{m} \subseteq \mathfrak{H}_{n}\;  \text{ and} \;  \mathfrak{H}_{m,t} \subseteq \mathfrak{H}_{n,t}\quad \text{ for } n \le m \text{ and } m,n\in\mathbb{N}.
\end{split}
\end{equation}

By Lemma \ref{approximate - lemma}, if $\sup_{(\xi,\omega) \in [-\d,\d]^{4}} \big| \nabla_{\th,\c}^{n+1}h(v,\xi,\omega)\big|$ is uniformly bounded and decays sufficiently fast as $v\to \infty $, then the $L^p_{v}$ norm of the difference is bounded by $\e^{n+1}\mathfrak{h}_{n+1}$.
\medskip

 The next lemma quantifies the error between the Maxwellians $\mu$ and $\tilde{\mu}$.
\begin{lemma} \label{approximate p - lemma}
Let $\norm{\d} <1$ and $p>0$ be given. For $x, y \in \R^3$ with $\norm{y}=1$, there exists a constant $c_{p}>0$ such that
\begin{equation*}
\begin{split}
&\Big | \exp\big(-\frac{\norm{x}^2}{p}\big) - (1+ \d) \exp\big(-\frac{\norm{x}^2}{p}\big)\Big |  \le c_{p} \norm{\d}, \\
&\Big | \exp\big(-\frac{\norm{x}^2}{p}\big) -  \exp\big(-\frac{\norm{x+ \d y}^2}{p}\big)\Big | \le c_{p} \norm{\d}, \\
&\Big | \exp\big(-\frac{\norm{x}^2}{p}\big) - \exp\big(-\frac{\norm{x}^2}{p(1+\d)}\big)\Big | \le c_{p} \norm{\d}.
\end{split}
\end{equation*}
\end{lemma}

\begin{proof}[\textbf{Proof.}] \
This follows directly from Lemma \ref{approximate - lemma}.
\end{proof}
\medskip

The following lemma estimates the error between $f$ and $\tilde{f}$ in weighted $L^p$ norms.

\begin{lemma}\label{tilde-f-bounded-by-f}
Let $w^\beta = e^{\beta\norm{v}^2}$ be a weight function with $0 \le \b < \b' < \frac{1}{4}$. Under the a priori assumption \eqref{theta-u-smallness-assumption}, for any $1 \le p \le \infty$,
\begin{equation*}
\begin{split}
\normm{w^{\b} \tilde{f}}_{L^p_v} \lesssim & \normm{w^{\b'} f}_{L^p_v} + \frac{\norm{\th}}{\e} + \frac{\norm{\c}}{\e}, \\
\normm{w^{\b} \pt_t \tilde{f}}_{L^p_v} \lesssim & \normm{w^{\b'} \pt_t f}_{L^p_v} + \normm{w^{\b'} f}_{L^p_v} + \frac{\norm{\th}}{\e} + \frac{\norm{\c}}{\e} + \frac{\norm{\pt_t \th}}{\e} + \frac{\norm{\pt_t \c}}{\e}.
\end{split}
\end{equation*}
\end{lemma}

\begin{proof}[\textbf{Proof.}] \
From the definition of $\tilde{f}$,
\begin{equation*}
\tilde{f} - f
%= \frac{\mu - \tilde{\mu}}{\e \sqrt{\tilde{\mu}}} + \frac{\sqrt{\mu}}{\sqrt{\tilde{\mu}}} %f - f
= \frac{\mu - \tilde{\mu}}{\e \sqrt{\tilde{\mu}}} + \Big(\frac{\sqrt{\mu}}{\sqrt{\tilde{\mu}}}-1\Big) f,
\end{equation*}
and similarly for the time derivative,
\begin{equation*}
\pt_t \tilde{f} - \pt_t f = \pt_t \frac{\mu - \tilde{\mu}}{\e \sqrt{\tilde{\mu}}} + \pt_t \Big(\frac{\sqrt{\mu}}{\sqrt{\tilde{\mu}}}-1\Big) f + \Big(\frac{\sqrt{\mu}}{\sqrt{\tilde{\mu}}}-1\Big) \pt_t f.
\end{equation*}
Using the structure of $\tilde{\mu}$ and Lemma \ref{approximate p - lemma}, we obtain for any $\b'' > 0$:
\begin{equation*}
\begin{split}
&\normm{\frac{\mu - \tilde{\mu}}{\e \sqrt{\tilde{\mu}}}}_{L^p_v}, \;\; \normm{w^{\b''}\Big(\frac{\sqrt{\mu}}{\sqrt{\tilde{\mu}}}-1\Big)}_{L^p_v}\lesssim \frac{\norm{\th}}{\e} + \frac{\norm{\c}}{\e}, \\
%------------------------
&\normm{\pt_t \Big(\frac{\mu - \tilde{\mu}}{\e \sqrt{\tilde{\mu}}}\Big)}_{L^p_v},  \;\;
\normm{w^{\b''}\pt_t \Big(\frac{\sqrt{\mu}}{\sqrt{\tilde{\mu}}}-1\Big)}_{L^p_v}\lesssim \frac{\norm{\th}}{\e} + \frac{\norm{\c}}{\e} +  \frac{\norm{\pt_t \th}}{\e} + \frac{\norm{\pt_t \c}}{\e}.
\end{split}
\end{equation*}
The desired estimates follow by the triangle inequality, absorbing the weight shift from $\b$ to $\b'$ where necessary.
\end{proof}
\medskip

We introduce an alternative, non-orthogonal basis $\{\tilde{\chi}_{i}\}_{i=0}^{4}$ for $\ker \tilde{L}$:
\begin{equation}
\begin{split}\label{base-tilde-chi}
&\tilde{\chi}_{0} := \sqrt{\tilde{\mu}}, \quad
\tilde{\chi}_{i} := v_i\sqrt{\tilde{\mu}} \;\; (i=1,2,3),\quad
\tilde{\chi}_{4} := \frac{\norm{v}^2-3}{\sqrt{6}}\sqrt{\tilde{\mu}}.
\end{split}
\end{equation}
The relation between the two bases $\{\bar{\chi}_{i}\}_{i=0}^{4}$ and $\{\tilde{\chi}_{i}\}_{i=0}^{4}$ is described in the following lemma.

\begin{lemma}\label{chi - almost same}\
The sets $\{{\bar{\chi}_{i}}\}_{i=0}^{4}$ defined in \eqref{base-bar-chi} and   $\{\tilde{\chi}_{i}\}_{i=0}^{4}$ defined in \eqref{base-tilde-chi} are both bases of $\ker \tilde{L}$, with $\{{\bar{\chi}_{i}}\}_{i=0}^{4}$ being orthonormal. Moreover,
for every $p \in [1,\infty]$,
\begin{equation*}
\normm{\tilde{\chi}_{i} - {\bar{\chi}}_{i}}_{L^{p}_{v}} \lesssim \e \Big(\frac{\norm{\theta}}{\e}+ \frac{\norm{\c}}{\e}\Big), \quad  i = 0, \dots, 4.
\end{equation*}
\end{lemma}

\begin{proof}[\textbf{Proof.}] \
By \eqref{ker-tilde-P},
%\begin{equation*}
%\ker \tilde{L} =\text{span} \left\{1, v, \norm{v}^2 \right\} \sqrt{\tilde{\mu}} = %\text{span} \left\{1, v-\c, \norm{v-\c}^2 \right\} \sqrt{\tilde{\mu}},
%\end{equation*}
 both $\{\tilde{\chi}_{i}\}_{i=0}^{4}$ and {$\{{\bar{\chi}_{i}}\}_{i=0}^{4}$} are bases of $\ker \tilde{L}$. The orthogonality of $\{{\bar{\chi}_{i}}\}_{i=0}^{4}$ follows from a direct computation. Furthermore, for each $i = 0, \dots, 4$, we have
\begin{align*}
\norm{\tilde{\chi}_{i} - \bar{\chi}_{i}} \lesssim \tilde{\mu}^{\frac{1}{4}}\e \Big(\frac{\norm{\theta}}{\e}+ \frac{\norm{\c}}{\e}\Big),
\end{align*}
which implies the desired estimate in $L^p_v(\R^3)$.
\end{proof}
\medskip

%
% The following result is not used directly, so is deleted.
%
%Since the error between the two bases $\{{\bar{\chi}_{i}}\}_{i=0}^{4}$ and %$%\{\tilde{\chi}_{i}\}_{i=0}^{4}$ is of order $\e \left(\frac{\norm{\theta}}{\e}+ %\frac{\norm{\c}}{\e}\right)$, we obtain the following corollary.

%\begin{corollary} \label{tildechi - barchi - similar}
%Let the coordinate transformation matrix $\mathscr{M}$ be defined by
%\begin{eqnarray}
%\tilde{\chi}_{i} = \sum_{j=0}^{4} \mathscr{M}_{ij} \bar{\chi}_{j}.
%\end{eqnarray}
%Then
%\begin{equation}
%\begin{split}
%\norm{\mathscr{M} - I_{5}}_{\infty} \lesssim \e \left(\frac{\norm{\theta}}{\e}+ %\frac{\norm{\c}}{\e}\right), \quad
%\norm{\mathscr{M}^{-1} - I_{5}}_{\infty} \lesssim \e \left(\frac{\norm{\theta}}{\e}+ %\frac{\norm{\c}}{\e}\right).
%\end{split}
%\end{equation}
%\end{corollary}

Recall the expansion \eqref{Pf-barabc -ch} of $\tilde{\P} \tilde{f} $ with coefficients \eqref{barabc-def}. Analogously, we define the coefficients of $\tilde{\P} \tilde{f} $ with respect to the basis $\{\tilde{\chi}_{i}\}_{i=0}^{4}$:
\begin{equation}\label{tilde-abc-def}
  \tilde{a}(t,x):=\langle  \tilde{\chi}_{0} , \tilde{f} \rangle, \quad
  \tilde{b}_i(t,x):=\langle \tilde{\chi}_{i} , \tilde{f}\rangle  \;\; (i=1,2,3),\quad
  \tilde{c}(t,x):=\langle \tilde{\chi}_{4} , \tilde{f}\rangle.
\end{equation}
The relationship between these two sets of coefficients is characterized by the following lemma.

\begin{lemma} \label{Pf - abc similar estimate}
Assume that the a priori assumption \eqref{theta-u-smallness-assumption} holds. Then for any $1 \le p, q \le \infty$, the following norm equivalence holds:
\begin{equation*}
\begin{split}
\normm{\tilde{\P}\tilde{f}}_{L^p_{x} L^q_{v}}
&\approx \normm{\bar{a}}_{L^p_{x}} +\sum_{i=1}^{3} \normm{\bar{b}_{i}}_{L^p_{x}}+\normm{\bar{c}}_{L^p_{x}}
\approx \normm{\tilde{a}}_{L^p_{x}} +\sum_{i=1}^{3} \normm{\tilde{b}_{i}}_{L^p_{x}}+\normm{\tilde{c}}_{L^p_{x}}.
\end{split}
\end{equation*}
\end{lemma}

\begin{proof}[\textbf{Proof.}] \
From the definition of $\tilde{\P} \tilde{f}$ and the expansion \eqref{Pf-barabc -ch}, we have
\begin{equation*}
\begin{split}
\normm{\tilde{\P}\tilde{f}(t)}_{L^p_{x} L^q_{v}}
&\approx \normm{\bar{a}(t)}_{L^p_{x}} +\sum_{i=1}^{3} \normm{\bar{b}_{i}(t)}_{L^p_{x}}+\normm{\bar{c}(t)}_{L^p_{x}}.
%\approx \normm{\tilde{a}(t)}_{L^p_{x}} +\sum_{i=1}^{3} %\normm{\tilde{b}_{i}(t)}_{L^p_{x}}+\normm{\tilde{c}(t)}_{L^p_{x}}.
\end{split}
\end{equation*}
%\begin{equation*}
%\begin{split}
%\normm{\tilde{\P}\tilde{f}(t)}_{L^p_{x} L^q_{v}}
%\lesssim  &\normm{\bar{a} \bar{\chi}_{0}}_{L^p_{x} L^q_{v}}
%+  \sum_{i=1}^{3} \normm{\bar{b}_{i} \bar{\chi}_{i}}_{L^p_{x} L^q_{v}}
%+ \normm{\bar{c} \bar{\chi}_{4}}_{L^p_{x} L^q_{v}} \\
%\lesssim & \normm{\bar{a} \normm{\bar{\chi}_{0}}_{L^q_{v}}}_{L^p_{x}} +
%\sum_{i=0}^{4} \normm{\bar{b}_{i} \normm{\bar{\chi}_{i}}_{L^q_{v}}}_{L^p_{x}}+
%\normm{\bar{c} \normm{\bar{\chi}_{4}}_{L^q_{v}}}_{L^p_{x}} \\
% \lesssim & \normm{\bar{a}}_{L^p_{x}}+ \sum_{i=1}^{3} \normm{\bar{b}_{i}}_{L^p_{x}} + %\normm{\bar{c}}_{L^p_{x}}.
%\end{split}\end{equation*}
%Conversely, by the orthonormality of $\{\bar{\chi}_i\}$,
%\begin{equation*}
%\normm{\bar{a}}_{L^p_{x}}+ \sum_{i=1}^{3} \normm{\bar{b}_{i}}_{L^p_{x}} + %\normm{\bar{c}}_{L^p_{x}}=
%\sum_{i=0}^{4} \normm{\inn{\bar{\chi}_{i},\tilde{\P} \tilde{f}}}_{L^p_{x}}
%\lesssim \normm{\tilde{\P}\tilde{f}(t)}_{L^p_{x} L^q_{v}}.
%\end{equation*}
For the coefficients $\tilde{a}, \tilde{b}_i$ and $\tilde{c}$ associated with the basis $\{\tilde{\chi}_i\}$ defined in \eqref{tilde-abc-def},
\begin{equation*}
\normm{\tilde{a}}_{L^p_{x}}+ \sum_{i=1}^{3} \normm{\tilde{b}_{i}}_{L^p_{x}} + \normm{\tilde{c}}_{L^p_{x}}=
\sum_{i=0}^{4} \normm{\langle \tilde{\chi}_{i},\tilde{\P} \tilde{f}\rangle}_{L^p_{x}}
\lesssim \normm{\tilde{\P}\tilde{f}}_{L^p_{x} L^q_{v}}.
\end{equation*}
Finally, comparing the two sets of coefficients, we obtain
\begin{equation*}
\begin{split}
 \normm{\bar{a}}_{L^p_{x}}+ \sum_{i=1}^{3} \normm{\bar{b}_{i}}_{L^p_{x}} + \normm{\bar{c}}_{L^p_{x}}
% = &\sum_{i=0}^{4} \normm{\inn{\chi_{i},\tilde{\P} \tilde{f}}}_{L^p_{x}} \\
%--------------
\leq & \sum_{i=0}^{4}\normm{\langle \tilde{\chi}_{i},\tilde{\P} \tilde{f} \rangle }_{L^p_{x}} + \sum_{i=0}^{4}\normm{\langle \bar{\chi}_{i}-\tilde{\chi}_{i},\tilde{\P} \tilde{f}\rangle }_{L^p_{x}} \\
%----------------------------
\lesssim& \normm{\tilde{a}}_{L^p_{x}}+ \sum_{i=1}^{3} \normm{\tilde{b}_{i}}_{L^p_{x}} + \normm{\tilde{c}}_{L^p_{x}}+ \e \Big(\frac{\norm{\theta}}{\e}+ \frac{\norm{\c}}{\e}\Big) \normm{\tilde{\P}\tilde{f}}_{L^p_{x}L^{q}_{v}}.
\end{split}
\end{equation*}
Under the smallness assumption on $\frac{\norm{\theta}}{\e} + \frac{\norm{\c}}{\e}$ from \eqref{theta-u-smallness-assumption}, the last term can be absorbed. Combining the estimates above yields the desired norm equivalence.
\end{proof}
\medskip

The next lemma provides a commutator estimate between $\pt_t$ and $\tilde{\P}$.
\begin{lemma}
The following commutator estimate holds:
\begin{equation*}
\normm{\pt_t(\tilde{\P} \tilde{f}) - \tilde{\P} \pt_t \tilde{f}}_{L^p_{x} L^q_{v}} \lesssim (\norm{\pt_t \th} + \norm{\pt_t \c}) \normm{\tilde{f}}_{L^p_{x} L^q_{v}}.
\end{equation*}
\end{lemma}

\begin{proof}[\textbf{Proof.}] \
Using the definition of $\tilde{\P}$ and the basis $\{\bar{\chi}_i\}$, we compute
\begin{equation*}
\begin{split}
\normm{\pt_t(\tilde{\P} \tilde{f}) - \tilde{\P} \pt_t \tilde{f}}_{L^p_{x} L^q_{v}} =& \normm{\sum_{i=0}^{4} \langle \pt_t \bar{\chi}_{i},\tilde{f}\rangle \bar{\chi}_{i} + \sum_{i=0}^{4} \langle \bar{\chi}_{i},\tilde{f}\rangle  \pt_t \bar{\chi}_{i}}_{L^p_{x} L^q_{v}}
\lesssim  (\norm{\pt_t \th} + \norm{\pt_t \c}) \normm{\tilde{f}}_{L^p_{x} L^q_{v}},
\end{split}
\end{equation*}
where we used  the estimate
$\normm{\pt_t \bar{\chi}_{i}}_{L^q_{v}} \lesssim \norm{\pt_t \th} + \norm{\pt_t \c}$ from the structure of $\bar{\chi}_i$.
\end{proof}
\medskip

The following lemma quantifies the approximation error when expressing the projection $\tilde{\mathbf{P}} g$ in the non-orthogonal basis $\{\tilde{\chi}_i\}$.

\begin{lemma}\label{expansion-in-tilde-chi}\
 Under the a priori assumption \eqref{theta-u-smallness-assumption}, for any $p \in [1,\infty]$ and $g\in L^p_v(\R^3)$,
\begin{equation*}
\Big| \tilde{\P} g - \sum_{i=0}^{4} \inn{  g, \tilde{\chi}_{i} } \tilde{\chi}_{i} \Big|
\lesssim \tilde{\mu}^{\frac{1}{4}} \e  \Big(\frac{\norm{\theta}}{\e}+ \frac{\norm{\c}}{\e}\Big) \normm{ \tilde{\P}  g}_{L^p_v}.
\end{equation*}
\end{lemma}

\begin{proof}[\textbf{Proof.}] \
Since $\{\bar{\chi}_{i}\}_{i=0}^{4}$ is an orthonormal basis  of $\ker \tilde{\P}$,
we write
\begin{equation*}
\begin{split}
g= &\sum_{i=0}^{4} \inn{  g, \bar{\chi}_{i} } \bar{\chi}_{i}  =
\sum_{i=0}^{4} \inn{  g, \tilde {\chi}_{i} } \tilde{\chi}_{i}
 -  \sum_{i=0}^{4}  \inn{  g, \tilde {\chi}_{i} } \left(\tilde{\chi}_{i} -\bar{\chi}_{i}  \right )
 -  \sum_{i=0}^{4}   \inn{  g, \tilde {\chi}_{i}-\bar {\chi}_{i} } \bar{\chi}_{i}.
\end{split}\end{equation*}
Observe that $\ipt g$ is orthogonal to both $\tilde {\chi}_{i}$ and $\tilde {\chi}_{i}-\bar {\chi}_{i} $. Applying H\"{o}lder's inequality and Lemma \ref{chi - almost same}, we bound the two error terms in above equality as
\begin{equation*}\begin{split}
& \Big|  \sum_{i=0}^{4}  \inn{  g, \tilde {\chi}_{i} } \left(\tilde{\chi}_{i} -\bar{\chi}_{i}  \right ) \Big|
 \lesssim     \normm{ \tilde{\P}  g}_{L^p_v} \normm{ \tilde{\chi}_{i} }_{L^q_v} \norm{\tilde{\chi}_{i} - {\bar{\chi}}_{i}}
 \lesssim    \tilde{\mu}^{\frac{1}{4}} \e  \Big(\frac{\norm{\theta}}{\e}+ \frac{\norm{\c}}{\e}\Big) \normm{ \tilde{\P}  g}_{L^p_v},\\
 %-------------------
 & \Big|   \sum_{i=0}^{4}   \inn{  g, \tilde {\chi}_{i}-\bar {\chi}_{i} } \bar{\chi}_{i}  \Big|
 \lesssim     \normm{ \tilde{\P}  g}_{L^p_v} \normm{\tilde{\chi}_{i} - {\bar{\chi}}_{i}}_{L^{q}_{v}}   \norm{\bar{\chi}_{i}}
 \lesssim    \tilde{\mu}^{\frac{1}{4}} \e  \Big(\frac{\norm{\theta}}{\e}+ \frac{\norm{\c}}{\e}\Big) \normm{ \tilde{\P}  g}_{L^p_v},
\end{split}\end{equation*}
where $\frac{1}{p}+\frac{1}{q}=1$. This completes the proof.
\end{proof}
\medskip

Next, we construct the functions $\rho$, $\c$, and $T$ in the definition of the rotating Maxwellian $\tilde{\mu}$.

\begin{lemma} \label{tildemu - existence}
Suppose the following conditions hold for $0 < \d \ll 1$:
\begin{equation}\label{F - smallness assumption}
\begin{split}
&\iint_{\O\times\mathbb{R}^3}  F(t) \dd v\dd x = \norm{\Omega}, \quad
\Big |  \iint_{\O\times\mathbb{R}^3}  Ax \cdot v F(t) \dd v\dd x \Big |  < \d, \quad
\Big | \iint_{\O\times\mathbb{R}^3} \norm{v}^2 F(t) \dd v\dd x - 3 \norm{\Omega} \Big |  < \d.
\end{split}
\end{equation}
Then there exist functions $\rho=\rho(t)$, $\c=\c(t,x) = \sum w_i(t) A_ix $  and $T=T(t)$ satisfying the following conservation laws:
\begin{equation} \label{tildemu-conservation-law}
\begin{split}
&\iint_{\O\times\mathbb{R}^3}  F(t) \dd v\dd x=\iint_{\O\times\mathbb{R}^3} \tilde{\mu} \dd v\dd x = |\O|, \\
& \iint_{\O\times\mathbb{R}^3} Ax \cdot v F(t) \dd v\dd x=\iint_{\O\times\mathbb{R}^3} Ax \cdot v \tilde{\mu}\dd v\dd x = \int_{\Omega} \rho Ax \cdot \c \dd x \;\; \text{ for all } A x \in \mathcal{R}_{\Omega},\\
&  \iint_{\O\times\mathbb{R}^3} \norm{v}^2 F(t) \dd v\dd x=\iint_{\O\times\mathbb{R}^3} \norm{v}^2 \tilde{\mu} \dd v\dd x= \int_{\Omega} (3\rho T + \rho \norm{\c}^2)\dd x.
\end{split}
\end{equation}
 Moreover, the perturbation $\tilde{f}$ satisfies:
\begin{equation}\label{abc - average zero}
\begin{split}
 &\iint_{\Omega\times\R^3} \sqrt{\tilde{\mu}}\tilde{f} \dd v \dd x = 0, \;\;
 \iint_{\Omega\times\R^3} Ax \cdot  v \sqrt{\tilde{\mu}}\tilde{f} \dd v \dd x = 0  \; \text{ for all } A x \in \mathcal{R}_{\Omega},\;\;
 \iint_{\Omega\times\R^3}  \norm{v}^2 \sqrt{\tilde{\mu}}\tilde{f} \dd v \dd x = 0.
\end{split}
\end{equation}
\end{lemma}

\begin{proof}[\textbf{Proof.}] \
Conditions in \eqref{F - smallness assumption} can guarantee the existence of a triple
$(\rho,\c,T)$ near $(1,0, 1)$. Using Lemma \ref{wholeinteg} and the definition of $\rho$ in \eqref{rho-def} and that of $\tilde{\mu}$ in \eqref{tilde-mu-def}, we have
\begin{equation*}%\label{tildemu-conservation-integral}
\begin{split}
&\iint_{\O\times\mathbb{R}^3}  \tilde{\mu} \dd x \dd v =  \norm{\O}, \\
&\iint_{\O\times\mathbb{R}^3}  Ax \cdot v \tilde{\mu} \dd v\dd x = \int_{\Omega} \rho Ax \cdot \c \dd x \;\text{ for all } A x \in \mathcal{R}_{\Omega},\;\;\\
&\iint_{\O\times\mathbb{R}^3}  \norm{v}^2 \tilde{\mu} \dd v\dd x = \int_{\Omega} (3\rho T + \rho \norm{\c}^2)\dd x.
\end{split}
\end{equation*}
This establishes the second equality in each line of \eqref{tildemu-conservation-law}.

We now treat the three geometric types of $\Omega$ separately.
\medskip

\noindent\textbf{Case 1. Non-axisymmetric domains.}

In this case, $\mathcal{R}_{\O}=\{0\}$. We define
\begin{equation}\label{non-axisymmetric-condition}
 \rho =1, \quad \c = 0 \quad \text{ and }\;\; T(t) = \frac{1}{3}\iint_{\Omega \times \R^3} \norm{v}^2 F(t) \dd v \dd x.
\end{equation}
Then the first equality in each line of \eqref{tildemu-conservation-law} follows directly.
\medskip

\noindent\textbf{Case 2. Axisymmetric domains.}

 In this case, $\c=  w Ax $. We seek functions $(\rho, w, T)$ satisfying
\begin{equation}\label{rho-w-theta-equation-axismmetric}
\begin{split}
&\rho - \frac{\norm{\Omega}\exp (\frac{\norm{\c (t,x)}^2}{2T(t)}) }{\int_{\Omega} \exp (\frac{\norm{\c (t,x)}^2}{2T(t)}) \dd x} = 0,\\
&w \int_{\Omega} \rho \norm{Ax}^2 \dd x  - \iint_{\Omega \times \R^3} Ax \cdot v F(t) \dd v \dd x = 0,\\
&3 T \int_{\Omega} \rho\dd x + w\iint_{\Omega \times \R^3} Ax \cdot v F(t) \dd v \dd x - \iint_{\Omega \times \R^3} \norm{v}^2 F(t) \dd v \dd x = 0.
\end{split}
\end{equation}
The Jacobian matrix of the system of $(\rho,w,T)$ at $(1,0,1)$ is
\begin{equation*}
\displaystyle
\begin{pmatrix}
 1& 0 & 0\\
 0 & \displaystyle\int_{\Omega} \norm{Ax}^2 \dd x & 0\\
 3\norm{\O} &\displaystyle \iint_{\Omega \times \R^3} Ax \cdot v F(t) \dd x\dd v & 3 \norm{\Omega}
\end{pmatrix},
\end{equation*}
which is invertible. By the implicit function theorem, a solution $(\rho,w,T)$ exists near
$(1,0,1)$.
\medskip

\noindent\textbf{Case 3. Spherical domains.}

For a spherical domain $\O$, $\c=\sum_{i=1}^3 w_i A_ix$. We have the orthogonality relations
\begin{equation}\label{A-x-integral-orthogonal}
\begin{split}
&\int_{\Omega} A_{i} x \cdot A_{j} x \dd x=0, \quad 
%{\int_{\Omega} \rho A_{i} x \cdot A_{j} x \dd x}= 0
%\quad \text{for}\;\; i \neq j,\\
 \int_{\pt \Omega} A_{i} x \cdot A_{j} x \dd S_x =0
% = {\int_{\pt \Omega} \rho A_{i} x \cdot A_{j} x \dd S_x}= 0 
 \quad \text{for}\;\; i \neq j,
\end{split}
\end{equation}
where we used the elementary identities
\begin{equation}\label{A-x-orthogonal}
\begin{split}
A_{1} x \cdot A_{2} x = x_{1} x_{2}, \quad
A_{2} x \cdot A_{3} x = x_{2} x_{3},\quad
A_{3} x \cdot A_{1} x = -x_{3} x_{1}.
\end{split}
\end{equation}
 We seek functions $(\rho, w_{1}, w_{2},  w_{3}, T)$ satisfying
\begin{equation}\label{rho-w-theta-equation-spherical}
\begin{split}
&\rho - \frac{\norm{\Omega}\exp (\frac{\norm{\c (t,x)}^2}{2T(t)}) }{\int_{\Omega} \exp (\frac{\norm{\c (t,x)}^2}{2T(t)}) \dd x} = 0,\\
%-------------------
&\sum_{j=1}^{3} w_{j} \int_{\Omega} \rho (A_{i}x\cdot A_jx) \dd x  
- \iint_{\Omega \times \R^3}  A_{i}x \cdot v F(t) \dd v \dd x = 0\quad \text{ for } \;\; i=1,2,3,\\
%-------------------
&3 T \int_{\Omega} \rho\dd x + \iint_{\Omega \times \R^3}  \sum_{i=1}^{3}w_{i} A_{i}x \cdot v F(t) \dd v \dd x - \iint_{\Omega \times \R^3}  \norm{v}^2 F(t) \dd v \dd x = 0.
\end{split}
\end{equation}
The Jacobian of this system at $(1,0,0,0,1)$ is
\begin{equation*}
\begin{pmatrix}
1 & 0 & 0 & 0& 0\\
0 & \displaystyle\int_{\Omega} \norm{A_1 x}^2 \dd x & 0 & 0 & 0\\
0 & 0 & \displaystyle\int_{\Omega} \norm{A_2 x}^2 \dd x & 0 & 0\\
0 & 0 & 0 & \displaystyle\int_{\Omega} \norm{A_3 x}^2 \dd x & 0\\
3\norm{\O}
& \displaystyle\iint A_{1}x \cdot v F(t) \dd v\dd x
& \displaystyle\iint A_{2}x \cdot v F(t) \dd v\dd x
& \displaystyle\iint A_{3}x \cdot v F(t) \dd v\dd x & 3 \norm{\Omega}
\end{pmatrix},
\end{equation*}
which is invertible. Hence a solution $(\rho,w_{1},w_{2},w_{3},T)$ exists near $(1,0,0,0,1)$.

Finally, \eqref{abc - average zero} follows from the relation $\sqrt{\tilde{\mu}}\tilde{f} = \frac{1}{\e}(F-\tilde{\mu})$ and the conservation laws \eqref{tildemu-conservation-law}.
\end{proof}
\medskip

Next, we derive the ordinary differential equations governing the evolution of $\rho$, $\c$ and $T$. The main result is summarized in the following proposition.

%{\bf \color{purple} (Remark: For the statement of Proposition 4.10 and Proposition 4.11 %(in Subsection 4.2), one may directly use the original equation $\e \pt_t F + v %\cdot\nabla_x F  = \e^{-1}Q(F, F)$ (since $Q$ satisfies the assumptions), rather than %using the linearized form $\e \pt_t F + v \cdot\nabla_x F  = G$ (which needs assumptions %on $G$).)}

\begin{proposition} \label{th tht w wt - ODE}\
Let $F$ be a solution of the Boltzmann equation \eqref{F - Boltzmann equation}, and let $\tilde{\mu}$ be the rotating Maxwellian defined in \eqref{tilde-mu-def}  with parameters $\rho$, $\c$ and $T = 1+ \th$. Let $\tilde{f}= \frac{1}{\e \sqrt{\tilde{\mu}}}(F -\tilde{\mu})$ be the fluctuation defined in \eqref{tildef - definition ch1}.
%linearized Boltzmann equation
%\begin{equation} \label{F - linearlized Boltzmann equation}
%\begin{split}
%\e \pt_t F + v \cdot\nabla_x F  = G \quad&
%   \text{in }  \mathbb{R}_{+}\times \O \times \mathbb{R}^3, \\
%F |_{\g_-}  =   (1-\a) \mathscr{R}F + \a \mathscr{P}F \quad &    \text {on }  %\mathbb{R}_{+}\times \pt\O \times \mathbb{R}^3, \\
%F(t, x, v)|_{t=0}  =   F_0(x,v) \quad & \text{on }   \O \times \mathbb{R}^3
%\end{split}
%\end{equation}
Then, under the a priori assumption \eqref{theta-u-smallness-assumption}, the following estimates hold:
\begin{align}
& \Big| \frac{3}{2} \pt_t \int_{\Omega} \theta^2 \dd x  + \frac{\a}{\sqrt{2 \pi}\e} \int_{\pt \Omega} 4\theta^2 \dd S_x + {\a}\iint_{\g_{+}} \big( \norm{v}^2 -4\big)  \sqrt{\tilde{\mu}} \tilde{f} \th  \dd \g   \Big|
\le  \a \e^{2}\mathfrak{h}_{3} + \a \e^{2} \mathfrak{h}_{2}\norm{\tilde{f}}_{L^2_{\g_{+}}}, \label{th - 2 estimate}\\
%--------------------
&  \Big|  \frac{1}{2}\pt_t \int_{\Omega} \norm{\c}^2 \dd x  + \frac{\a}{\sqrt{2 \pi}\e} \int_{\pt \Omega} \norm{\c}^2 \dd S_x
+ {\a}\iint_{\g_{+}} (\c \cdot v) \sqrt{\tilde{\mu}} \tilde{f} \dd \g    \Big|
\le  \a \e^2\mathfrak{h}_{3} + \a \e^{2} \mathfrak{h}_{2}\norm{\tilde{f}}_{L^2_{\g_{+}}},
\label{w - 2 estimate}\\
%----------------------
&  \Big| \frac{3}{2} \pt_t \int_{\Omega} (\pt_t \theta)^2 \dd x  + \frac{\a}{\sqrt{2 \pi}\e} \int_{\pt \Omega} 4(\pt_t \theta)^2 \dd S_x + {\a}\iint_{\g_{+}} \big(\norm{v}^2 -4\big) \sqrt{\tilde{\mu}} \pt_t \tilde{f} \pt_t \th  \dd \g   \Big|\nonumber \\
\le & \a \e^{2}\mathfrak{h}_{3,t} + \a \e^{2} \mathfrak{h}_{2,t}\norm{\pt_t \tilde{f}}_{L^2_{\g_{+}}}, \label{tht - 2 estimate}\\
%----------------------
&  \Big|  \frac{1}{2}\pt_t \int_{\Omega} \norm{\pt_t \c}^2 \dd x  + \frac{\a}{\sqrt{2 \pi}\e} \int_{\pt \Omega} \norm{\pt_t \c}^2 \dd S_x
+ {\a}\iint_{\g_{+}} (\pt_t \c \cdot v) \sqrt{\tilde{\mu}} \pt_t \tilde{f}  \dd \g  \Big|\nonumber\\
\le & \a \e^2\mathfrak{h}_{3,t} + \a \e^{2} \mathfrak{h}_{2,t}\norm{\pt_t\tilde{f}}_{L^2_{\g_{+}}}, \label{wt - 2 estimate}
\end{align}
where $\mathfrak{h}_{n} \in \mathfrak{H}_{n}$ and $\mathfrak{h}_{n,t} \in \mathfrak{H}_{n,t}$ are defined in \eqref{High order set - Definition}.

Additionally, the following bounds hold:
\begin{equation} \label{rho th u - smallness}
\begin{split}
&\norm{\th} \le \e \mathfrak{h}_{1}, \quad \norm{\c} \le \e \mathfrak{h}_{1}, \quad \norm{\rho - 1} \le \e^{2} \mathfrak{h}_{2}, \\
%-----------
&\norm{\pt_t \th}+\norm{\pt_t \c} \le \a \mathfrak{h}_{1} + \a \norm{\tilde{f}}_{L^2_{\g_{+}}}, \quad \norm{\pt_t \rho} \le \a \e \mathfrak{h}_{2} + \a\e  \mathfrak{h}_{1} \norm{\tilde{f}}_{L^2_{\g_{+}}}, \\
%-----------
&
\norm{\pt_t \pt_t\th}+\norm{\pt_t \pt_t\c} \le \a \mathfrak{h}_{1,t} + \a \norm{\pt_t\tilde{f}}_{L^2_{\g_{+}}}, \quad \norm{\pt_t \pt_t \rho} \le \a \e \mathfrak{h}_{2,t} + \a\e  \mathfrak{h}_{1,t} \norm{\pt_t \tilde{f}}_{L^2_{\g_{+}}}.
\end{split}
\end{equation}
\end{proposition}

\begin{proof}[\textbf{Proof.}] \
 Clearly, the Boltzmann collision operator $Q(F,F)$ satisfies the orthogonal condition
\begin{equation}\label{G-solvability-condition}
\begin{split}
\iint_{\Omega \times \R^3}
\big[1, Ax \cdot v, \norm{v}^2\big] Q(F,F) \dd v \dd x = 0 \;\; \text{ for all } A x \in \mathcal{R}_{\Omega}.
\end{split}
\end{equation}
Therefore, using \eqref{F - Boltzmann equation} and \eqref{G-solvability-condition}, a direct computation shows
%\begin{equation*}
%\begin{split}
%0 =& \iint_{\Omega \times \R^3} G \dd v\dd x  \\
%=& \iint_{\Omega \times \R^3} (\e \pt_t F + v \cdot \nabla F) \dd v\dd x \\
%=& \e \iint_{\Omega \times \R^3} \pt_t F \dd x \dd v + \iint_{\Omega \times \R^3}v \cdot %\nabla F \dd x \dd v \\
%=& \e \iint_{\Omega \times \R^3} \pt_t F \dd v\dd x + \iint_{\pt \Omega \times \R^3}
%F [n \cdot v] \dd v\dd x \\
%=& \e \iint_{\Omega \times \R^3} \pt_t F \dd x \dd v + \iint_{\g_{+}} F \dd \g - %\iint_{\g_{-}} F \dd \g \\
%=& \e \iint_{\Omega \times \R^3} \pt_t F \dd v\dd x + \iint_{\g_{+}} F \dd \g - %\iint_{\g_{-}} [(1-\a) \mathscr{R}F + \a \mathscr{P}F] \dd \g \\
%=& \e \iint_{\Omega \times \R^3} \pt_t F \dd v\dd x,
%\end{split}
%\end{equation*}
%where the boundary condition and conservation property were used. Then we obtain
$$
  \pt_t \iint_{\Omega \times \R^3} F \dd v\dd x
  = -\frac{1}{\e} \iint_{\pt \Omega \times \R^3} F[n\cdot v] \dd v\dd S_x
  + \frac{1}{\e} \iint_{\Omega \times \R^3} Q(F,F) \dd v\dd x=0.
$$
Combined with \eqref{initial-conservation-F}, this implies
$$
 \iint_{\Omega \times \R^3} F(t)  \dd v\dd x  = \iint_{\Omega \times \R^3} F_0  \dd v\dd x  = \norm{\Omega}\;\; \text{ for all }  t\geq 0.
$$
Lemma \ref{tildemu - existence} then guarantees the existence of $\tilde{\mu}$ satisfying \eqref{tildemu-conservation-law}.

We proceed by a case analysis based on the geometry of $\Omega$.
\medskip

\noindent\textbf{Case 1. Non-axisymmetric domains.}

In this case, $\rho = 1$ and $\c = 0$ by \eqref{non-axisymmetric-condition}. Multiplying \eqref{F - Boltzmann equation} by $\norm{v}^2$ and integrating over $\Omega \times \R^3$ yields
\begin{equation}\label{T-equation-energy}
\begin{split}
\e \int_{\Omega} 3 \pt_t  T \dd x
= & \e \pt_t \iint_{\Omega \times \R^3} \norm{v}^2 F \dd v\dd x
%----------------------
%= &- \iint_{\Omega \times \R^3} \norm{v}^2 v \cdot \nabla_x F \dd v\dd x
= - \iint_{\pt\Omega \times \R^3} \norm{v}^2 F [n\cdot v] \dd v \dd S_x  \\
%= &- \iint_{\g_+} \norm{v}^2 F [n\cdot v] \dd v \dd S_x -  \iint_{\g_{-}} \norm{v}^2 F %[n\cdot v]  \dd v \dd S_x \\
=& - \iint_{\g_+} \norm{v}^2 F  \dd \g  + \iint_{\g_{-}} \norm{v}^2 ((1- \a) \mathscr{R} F + \a \mathscr{P} F) \dd \g  \\
=& - \iint_{\g_+} \norm{v}^2 F  \dd \g   +  (1- \a) \iint_{\g_{+}} \norm{R_xv}^2 F \dd \g  + \a \iint_{\g_{+}} \norm{v}^2 \mathscr{P} F  \dd \g  \\
%------------
%=& - \a \iint_{\g_{+}}   \norm{v}^2 (1-\mathscr{P}) F  \dd \g  \\
=&- \a \Big[\iint_{\g_{+}} \norm{v}^2 F  \dd \g - \sqrt{2\pi}\iint_{\g_{+}} \norm{v}^2 \mu  \Big(\int_{n \cdot u >0} F [n \cdot u] \dd u\Big)  \dd \g \Big]\\
%------------------------
=& -\a \Big[\iint_{\g_{+}} \norm{v}^2\tilde{\mu}   \dd \g  - \sqrt{2\pi}\iint_{\g_{+}} \norm{v}^2 \mu  \Big(\int_{n \cdot u >0} \tilde{\mu} [n \cdot u] \dd u\Big)  \dd \g  \Big]\\
&-\a\e \Big[\iint_{\g_{+}} \norm{v}^2  \sqrt{\tilde{\mu}} \tilde{f} \dd \g  - \sqrt{2\pi}\iint_{\g_{+}} \norm{v}^2  \Big(\int_{n\cdot u >0} \sqrt{\tilde{\mu}} \tilde{f} [n \cdot u]  \dd u \Big)  \dd \g \Big]\\
%------------------------
=& -\a \Big[\int_{\pt \Omega} \frac{4 T^{\frac{3}{2}}}{\sqrt{2 \pi}}\dd S_x - \sqrt{2\pi}\int_{\pt \Omega} \frac{4}{\sqrt{2 \pi}} \frac{ T^{\frac{1}{2}}}{\sqrt{2 \pi}} \dd S_x \Big] \\
&-\a\e \Big[\iint_{\g_{+}} \norm{v}^2  \sqrt{\tilde{\mu}} \tilde{f}  \dd \g  - \sqrt{2\pi}\int_{\pt \Omega} \frac{4}{\sqrt{2 \pi}} \int_{n\cdot v >0} \sqrt{\tilde{\mu}} \tilde{f}  \dd \g  \Big] \\
%------------------------
=& -\a \Big[\frac{1}{\sqrt{2 \pi}}
\int_{\pt \Omega} 4(T-1) T^{\frac{1}{2}} \dd S_x + \e \iint_{\g_{+}} ( \norm{v}^2 -4) \sqrt{\tilde{\mu}} \tilde{f}  \dd \g  \Big],
\end{split}
\end{equation}
where we have used \eqref{non-axisymmetric-condition}, \eqref{G-solvability-condition} and Lemma \ref{wholeinteg}. Writing $\th = T-1$, we have
\begin{align}
& \Big|  3 \pt_t \int_{\Omega} \theta \dd x  + \frac{\a}{\sqrt{2 \pi}\e} \int_{\pt \Omega} 4\theta \dd S_x + {\a}\iint_{\g_{+}} \big( \norm{v}^2 -4\big) \sqrt{\tilde{\mu}} \tilde{f}  \dd \g   \Big|
\le \a \e\mathfrak{h}_{2}, \label{th - nonaxi 1 estimate}\\
%------------------
& \Big| \frac{3}{2} \pt_t \int_{\Omega} \theta^2 \dd x  + \frac{\a}{\sqrt{2 \pi}\e} \int_{\pt \Omega} 4\theta^2 \dd S_x + {\a}\iint_{\g_{+}} \big( \norm{v}^2 -4\big) \sqrt{\tilde{\mu}} \tilde{f} \th  \dd \g   \Big|
\le \a \e^{2}\mathfrak{h}_{3},  \label{th - nonaxi 2 estimate}
\end{align}
where $\mathfrak{h}_{n} \in \mathfrak{H}_{n}$. Consequently,
\begin{equation*}
\norm{\pt_t \th} \lesssim \a \mathfrak{h}_{1} + \a \norm{\tilde{f}}_{L^2_{\g_{+}}}.
\end{equation*}
Furthermore, differentiating \eqref{T-equation-energy} with respect to $t$ yields the equation for $\pt_t \th$:
\begin{equation*}
\begin{split}
&3 \pt_t \int_{\Omega} \pt_t \theta \dd x  + \frac{\a}{\sqrt{2 \pi}\e} \int_{\pt \Omega} 4\pt_t \theta(1+\th)^{\frac{1}{2}} \dd S_x + {\a}\iint_{\g_{+}} \big( \norm{v}^2 -4\big) \sqrt{\tilde{\mu}} \pt_t \tilde{f}  \dd \g  \\
&+ \frac{\a}{\sqrt{2 \pi}\e} \int_{\pt \Omega} 2 \theta (1+\th)^{-{\frac{1}{2}}} \pt_t \th \dd S_x + {\a}\iint_{\g_{+}} \big( \norm{v}^2 -4\big) \pt_t \sqrt{\tilde{\mu}}  \tilde{f}  \dd \g   = 0.
\end{split}
\end{equation*}
This leads to the estimates
\begin{align*}
&  \Big| 3 \pt_t \int_{\Omega} \pt_t\theta \dd x  + \frac{4\a}{\sqrt{2 \pi}\e} \int_{\pt \Omega} \pt_t \theta \dd S_x + {\a}\iint_{\g_{+}} \big( \norm{v}^2 -4\big) \sqrt{\tilde{\mu}} \pt_t \tilde{f}  \dd \g   \Big|
\le \a \e\mathfrak{h}_{2,t},  \\ %\label{tht - nonaxi 1 estimate}
%----------------
& \Big| \frac{3}{2} \pt_t \int_{\Omega} (\pt_t \theta)^2 \dd x  + \frac{4\a}{\sqrt{2 \pi}\e} \int_{\pt \Omega} (\pt_t \theta)^2 \dd S_x + {\a}\iint_{\g_{+}} \big( \norm{v}^2 -4\big) \sqrt{\tilde{\mu}} \pt_t \tilde{f} \pt_t \th  \dd \g \Big|
\le \a \e^{2}\mathfrak{h}_{3,t}, %\label{tht - nonaxi 2 estimate}
\end{align*}
where $\mathfrak{h}_{n,t} \in \mathfrak{H}_{n,t}$. Hence,
\begin{equation*}
\norm{\pt_t \pt_t \th} \lesssim \a  \mathfrak{h}_{1,t} + \a \norm{\pt_t \tilde{f}}_{L^2_{\g_{+}}}.
\end{equation*}

\noindent\textbf{Case 2. Axisymmetric domains.}

In this case, it follows from \eqref{rho-def} that $\rho(t,x)= 1+ O(\norm{\c}^2)$.
% =& \frac{\norm{\Omega}\exp (\frac{\norm{\c(t,x)}^2}{2T(x)}) }{\int_{\Omega} \exp (\frac{\norm{\c(t,x)}^2}{2T(x)}) \dd x}
Elementary calculation shows
\begin{align} \label{rho - pt t}
\pt_t \rho =& \bigg[ \frac{\c \cdot \pt_t \c}{T} - \frac{\norm{\c}^2}{2T^2}\pt_t T -\frac{\int_{\Omega}\big( \frac{\c \cdot \pt_t \c}{T} - \frac{\norm{\c}^2}{2T^2}\pt_t T\big) \exp (\frac{\norm{\c(t,x)}^2}{2T(x)}) \dd x}{\int_{\Omega} \exp (\frac{\norm{\c(t,x)}^2}{2T(x)}) \dd x}\bigg]\rho := P [\pt_t w, \pt_t \th]^{t},
\end{align}
where every entry of the matrix $P$ is of order $O(w,\th)$. Similarly,
$\pt_t \pt_t \rho = P [\pt_t \pt_t w, \pt_t \pt_t w]^{t} + Q$,
where $Q$ is bounded by $\norm{\pt_t w}^2+ \norm{\pt_{t} \th}^2$.

Multiplying  \eqref{F - Boltzmann equation} by $\norm{v}^2$ and integrating
over $\Omega \times \R^3$, we obtain
\begin{equation*}
\begin{split}
&\e \int_{\Omega} \pt_t (3{\rho} T  +  {\rho} w^2 \norm{Ax}^2 ) \dd x\\
= &-\a \Big(\frac{1}{\sqrt{2 \pi}}
\int_{\pt \Omega} [4(T-1)+ w^2 \norm{Ax}^2] {\rho}T^{\frac{1}{2}}  \dd S_x + \e \iint_{\g_{+}} \big( \norm{v}^2 -4\big) \sqrt{\tilde{\mu}} \tilde{f}  \dd \g  \Big),
\end{split}
\end{equation*}
where we have used \eqref{rho-w-theta-equation-axismmetric}, \eqref{G-solvability-condition} and Lemma \ref{wholeinteg}. Similarly, multiplying \eqref{F - Boltzmann equation} by $Ax \cdot v$ and integrating
over $\Omega \times \R^3$ yield
\begin{align*}
\e \int_{\Omega} \norm{Ax}^2\pt_t {(\rho w)} \dd x
= &  \e \pt_t \iint_{\Omega \times \R^3} Ax \cdot v F \dd v \dd x
%=& - \iint_{\Omega \times \R^3} (Ax \cdot v) v \cdot \nabla_x F \dd v \dd x
= - \iint_{\pt\Omega \times \R^3} (Ax \cdot v)F [n\cdot v] \dd v \dd S_x  \\
%=& - \iint_{\g_+} (Ax \cdot v) F [n\cdot v] \dd v \dd S_x   -  \iint_{\g_{-}} (Ax \cdot v) %F [n\cdot v] \dd v \dd S_x  \\
=& - \iint_{\g_+} (Ax \cdot v) F\dd \g    +  \iint_{\g_{-}} (Ax \cdot v) ((1- \a) R F + \a \mathscr{P} F) \dd \g   \\
=& - \iint_{\g_+} (Ax \cdot v) F  \dd \g  +  (1- \a) \iint_{\g_{+}} (Ax \cdot Rv) F \dd \g + \a \iint_{\g_{+}} (Ax \cdot v) \mathscr{P} F  \dd \g \\
%------------------------
=& - \a \iint_{\g_{+}}   (Ax \cdot v) (1-\mathscr{P}) F  \dd \g  \quad\\
%------------------------
=&- \a \Big[\iint_{\g_{+}} (Ax \cdot v) F  \dd \g - \sqrt{2\pi}\iint_{\g_{+}} (Ax \cdot v) \tilde{\mu}  \Big(\int_{n \cdot u >0} F [n \cdot u] \dd u\Big)  \dd \g \Big]\\
%------------------------
=& -\a \Big[\iint_{\g_{+}} (Ax \cdot v) \tilde{\mu}   \dd \g
- \sqrt{2\pi}\iint_{\g_{+}} (Ax \cdot v) \mu  \Big(\int_{n \cdot u >0} \tilde{\mu} [n \cdot u] \dd u\Big)   \dd \g \Big]\\
%------------------------
&-\a\e \Big[\iint_{\g_{+}} (Ax \cdot v)  \sqrt{\tilde{\mu}} \tilde{f}  \dd \g
- \sqrt{2\pi}\iint_{\g_{+}} (Ax \cdot v) \mu  \Big(\int_{n \cdot u >0} \sqrt{\tilde{\mu}} \tilde{f} [n \cdot v] \dd u\Big)  \dd \g \Big] \\
%------------------------
=& -\a \Big[\int_{\pt \Omega} \frac{ {\rho} T^{\frac{1}{2}} w \norm{Ax}^2}{\sqrt{2 \pi}} \dd S_x
+ \e \iint_{\g_{+}} (Ax \cdot v)  \sqrt{\tilde{\mu}} \tilde{f}  \dd \g  \Big].
\end{align*}
The equations for $w$ and $\theta$ are
\begin{align*}
&\bigg[
\begin{pmatrix}
\int_{\Omega} \norm{Ax}^2 \dd x & 0\\
0 & \int_{\Omega} 3 \dd x
\end{pmatrix} + P_{1} \bigg]
\begin{pmatrix}
\pt_t w \\
\pt_t \th
\end{pmatrix}
+ \a
\begin{pmatrix}
\frac{1}{\sqrt{2 \pi} \e} \int_{\pt \Omega} w \norm{Ax} \dd S_{x}    + \iint_{\g_{+}} (Ax \cdot v) \sqrt{\tilde{\mu}} \tilde{f} \dd \g \\
\frac{1}{\sqrt{2 \pi} \e} \int_{\pt \Omega} 4 \th \dd S_{x} + \iint_{\g_{+}} \big( \norm{v}^2 -4\big) \sqrt{\tilde{\mu}} \tilde{f} \dd \g
\end{pmatrix} = 0,
\end{align*}
where $P_{1}$ is of order $O(\norm{w},\th)$.
%and $P_{2}$ is of order $O(\norm{w}^2,\th^2)$.
Multiplying this by
$$\begin{pmatrix}
\int_{\Omega} \norm{Ax}^2 \dd x & 0\\
0 & \int_{\Omega} 3 \dd x
\end{pmatrix}\left[
\begin{pmatrix}
\int_{\Omega} \norm{Ax}^2 \dd x & 0\\
0 & \int_{\Omega} 3 \dd x
\end{pmatrix} + P_{1} \right]^{-1}
$$
yields
\begin{align*}
&\pt_t \int_{\Omega} w \norm{Ax}^2 \dd x  + \frac{\a}{\sqrt{2 \pi}\e} \int_{\pt \Omega} w \norm{Ax}^2 \dd S_x
+ {\a}\iint_{\g_{+}} (Ax \cdot v) \sqrt{\tilde{\mu}} \tilde{f}  \dd \g   = {h}_{1}, \\
&\pt_t \int_{\Omega} 3\theta \dd x  + \frac{\a}{\sqrt{2 \pi}\e} \int_{\pt \Omega} 4\theta \dd S_x
+ {\a}\iint_{\g_{+}} \big( \norm{v}^2 -4\big)  \sqrt{\tilde{\mu}} \tilde{f}  \dd \g  = {h}_{2},
\end{align*}
where $h_{1}$ and $h_{2}$ are bounded by
\begin{align*}
\a (\norm{w} + \norm{\th})\Big(\norm{w} + \norm{\th} + \norm{\tilde{f}}_{L^2_{\g_{+}}}\Big).
\end{align*}
Thus, we have
\begin{align}
&  \Big|  3 \pt_t \int_{\Omega} \theta \dd x  + \frac{\a}{\sqrt{2 \pi}\e} \int_{\pt \Omega} 4\theta \dd S_x + {\a}\iint_{\g_{+}} \big( \norm{v}^2 -4\big) \sqrt{\tilde{\mu}} \tilde{f} \dd \g  \Big|
\le \a \e\mathfrak{h}_{2} + \a \e\mathfrak{h}_{1}\norm{\tilde{f}}_{L^2_{\g_{+}}}  \label{th - axi 1 estimate}\\
%------------------
& \Big| \pt_t \int_{\Omega} w \norm{Ax}^2 \dd x  + \frac{\a}{\sqrt{2 \pi}\e} \int_{\pt \Omega} w \norm{Ax}^2 \dd S_x
+ {\a}\iint_{\g_{+}} (Ax \cdot v) \sqrt{\tilde{\mu}} \tilde{f} \dd \g  \Big|
\le \a \e\mathfrak{h}_{2} + \a \e\mathfrak{h}_{1}\norm{\tilde{f}}_{L^2_{\g_{+}}}, \label{w - axi 1 estimate}
\end{align}
where $\mathfrak{h}_{n} \in \mathfrak{H}_{n}$. Moveover,
\begin{align}
&  \Big| \frac{3}{2} \pt_t \int_{\Omega} \theta^2 \dd x  + \frac{\a}{\sqrt{2 \pi}\e} \int_{\pt \Omega} 4\theta^2 \dd S_x + {\a}\iint_{\g_{+}} \big( \norm{v}^2 -4\big) \sqrt{\tilde{\mu}} \tilde{f} \th \dd \g  \Big|
\le \a \e^{2}\mathfrak{h}_{3} + \a \e^{2} \mathfrak{h}_{2}\norm{\tilde{f}}_{L^2_{\g_{+}}}, \label{th - axi 2 estimate}\\
%-------------
& \Big|  \frac{1}{2}\pt_t \int_{\Omega} w^2 \norm{Ax}^2 \dd x  + \frac{\a}{\sqrt{2 \pi}\e} \int_{\pt \Omega} w^2 \norm{Ax}^2 \dd S_x
+ {\a}\iint_{\g_{+}} (w Ax \cdot v) \sqrt{\tilde{\mu}} \tilde{f} \dd \g \Big|
\le \a \e^2\mathfrak{h}_{3} + \a \e^{2} \mathfrak{h}_{2}\norm{\tilde{f}}_{L^2_{\g_{+}}}. \label{w - axi 2 estimate}
\end{align}
It follows that
\begin{equation*}
\begin{split}
\norm{\pt_t \th} \lesssim \a \mathfrak{h}_{1} + \a \norm{\tilde{f}}_{L^2_{\g_{+}}}, \quad
\norm{\pt_t w} \lesssim \a \mathfrak{h}_{1} + \a \norm{\tilde{f}}_{L^2_{\g_{+}}}.
\end{split}
\end{equation*}

Proceeding as in Case 1, we also obtain estimates for $\pt_t \th$ and $\pt_t w$:
\begin{align*}
& \Big| 3 \pt_t \int_{\Omega} \pt_t \theta \dd x  + \frac{\a}{\sqrt{2 \pi}\e} \int_{\pt \Omega} 4 \pt_t \theta \dd S_x + {\a}\iint_{\g_{+}} \big( \norm{v}^2 -4\big)  \sqrt{\tilde{\mu}} \pt_t \tilde{f}  \dd \g \Big|
\le \a \e\mathfrak{h}_{2,t} + \a \e\mathfrak{h}_{1,t}\norm{\pt_t \tilde{f}}_{L^2_{\g_{+}}}, \\ %\label{tht - axi 1 estimate}
%---------------
& \Big|  \pt_t \int_{\Omega} \pt_t w \norm{Ax}^2 \dd x  + \frac{\a}{\sqrt{2 \pi}\e} \int_{\pt \Omega} \pt_t w \norm{Ax}^2 \dd S_x
+ {\a}\iint_{\g_{+}} (Ax \cdot v) \sqrt{\tilde{\mu}} \pt_t \tilde{f} \dd \g \Big|\\
\le & \a \e\mathfrak{h}_{2,t} + \a \e\mathfrak{h}_{1,t}\norm{\pt_t \tilde{f}}_{L^2_{\g_{+}}},\\ %\label{wt - axi 1 estimate}
%------------------
& \Big| \frac{3}{2} \pt_t \int_{\Omega} (\pt_t \theta)^2 \dd x  + \frac{\a}{\sqrt{2 \pi}\e} \int_{\pt \Omega} 4(\pt_t \theta)^2 \dd S_x + {\a}\iint_{\g_{+}} \big( \norm{v}^2 -4\big)  \sqrt{\tilde{\mu}} \pt_t \tilde{f} \pt_t \th  \dd \g \Big|\\
\le &\a \e^{2}\mathfrak{h}_{3,t} + \a \e^{2} \mathfrak{h}_{2,t}\norm{\pt_t \tilde{f}}_{L^2_{\g_{+}}},\\ %\label{tht - axi 2 estimate}
%-------------
& \Big| \frac{1}{2}\pt_t \int_{\Omega} (\pt_t w)^2 \norm{Ax}^2 \dd x  + \frac{\a}{\sqrt{2 \pi}\e} \int_{\pt \Omega} (\pt_t w)^2 \norm{Ax}^2 \dd S_x
+ {\a}\iint_{\g_{+}} (\pt_t w Ax \cdot v) \sqrt{\tilde{\mu}} \pt_t \tilde{f} \dd \g \Big| \\
\le & \a \e^2\mathfrak{h}_{3,t} + \a \e^{2} \mathfrak{h}_{2,t}\norm{\pt_t\tilde{f}}_{L^2_{\g_{+}}},%\label{wt - axi 2 estimate}
\end{align*}
where $\mathfrak{h}_{n,t} \in \mathfrak{H}_{n,t}$. Therefore,
\begin{equation*}
\begin{split}
\norm{\pt_t \pt_t \th} \lesssim \a  \mathfrak{h}_{1,t} + \a \norm{\pt_t \tilde{f}}_{L^2_{\g_{+}}}, \quad
\norm{\pt_t \pt_t w} \lesssim \a  \mathfrak{h}_{1,t} + \a \norm{\pt_t \tilde{f}}_{L^2_{\g_{+}}}.
\end{split}
\end{equation*}

\noindent\textbf{Case 3. Spherical domains.}

Similar to Case 2,  multiplying   \eqref{F - Boltzmann equation} by $\norm{v}^2$ and integrating over $\Omega \times \R^3$ gives
\begin{equation*}
\begin{split}
&\e \int_{\Omega} \pt_t \big[3{\rho} T  +  {\rho}\sum_{i,j=1}^{3} w_{i}w_{j} (A_{i}x\cdot A_{j}x) \big] \dd x \\
=& -\a \Big(\frac{1}{\sqrt{2 \pi}}
\int_{\pt \Omega}  {\rho} T^{\frac{1}{2}} \big[4(T-1)+ |\c|^2\big] \dd S_x + \e \iint_{\g_{+}} \big( \norm{v}^2 -4\big) \sqrt{\tilde{\mu}} \tilde{f} \dd \g \Big).
\end{split}
\end{equation*}
Multiplying   \eqref{F - Boltzmann equation} by $Ax \cdot v$ and integrating
over $\Omega \times \R^3$ yield
\begin{align*}
\e \pt_t  \int_{\Omega} \rho ( A_{i}x\cdot \c) \dd x
+ \a \Big(\frac{1}{\sqrt{2 \pi}} \int_{\pt \Omega} {\rho} T^{\frac{1}{2}} (A_i x \cdot \c) \dd S_x
+ \e \iint_{\g_{+}} (A_{i}x \cdot v)  \sqrt{\tilde{\mu}} \tilde{f} \dd \g\Big)
=0
\end{align*}
for each $i=1,2,3$. Here we have used \eqref{A-x-integral-orthogonal}--\eqref{rho-w-theta-equation-spherical}, \eqref{G-solvability-condition} and Lemma \ref{wholeinteg}. The formulas are identical to those for axisymmetric domains in Case 2. Therefore, the same conclusions follow directly. This completes the proof.
\end{proof}
\bigskip

%%%%%%%%%%%%%%%%%%%%%%%%%%%%%%
\subsection{Energy Estimate}\
\medskip

In this subsection, we establish energy estimates for the fluctuation $\tilde{f}$ and its time derivative $\pt_t\tilde{f}$.

Differentiating equation \eqref{tildef - Boltzmanneq ch1} gives the equation for $\pt_t \tilde{f}$:
\begin{equation} \label{tildef t - Boltzmann with g}
\begin{split}
\e \pt_t (\pt_t \tilde{f}) + v \cdot \nabla_x (\pt_t \tilde{f}) + \e^{-1} \tilde{L}(\pt_t \tilde{f}) = \tilde{g}^{t}   \;\;\;\;  &\text{ in } \mathbb{R}^{+}\times\O\times\mathbb{R}^3, \\
%---------------
\pt_t \tilde{f} |_{\g_{-}} = (1- \a) \mathscr{R} (\pt_t \tilde{f}) + \a \tilde{\mathscr{P}}_{\g} (\pt_t \tilde{f}) +
\a \pt_t r + \a s  \;\;\;\;   &\text{ in } \mathbb{R}^{+}\times\pt\O\times\mathbb{R}^3,\\
%---------------
 \pt_t\tilde{f}|_{t=0}  =   \pt_t\tilde{f}_0  \;\;\;\;    &\text{ on }   \O \times \mathbb{R}^3,
\end{split}
\end{equation}
where  $ \pt_t\tilde{f}_0$ is determined through \eqref{tildef - Boltzmanneq ch1}, the boundary term $r$ is defined in \eqref{tildeg - definition}, and
\begin{equation} \label{tildegt - definition}
\begin{split}
{\tilde{g}^{t}}  :=  & \tilde{\G}(\pt_t \tilde{f}, \tilde{f}) + \tilde{\G}(\tilde{f}, \pt_t \tilde{f}) +\pt_t \Big(\frac{1}{\sqrt{\tilde{\mu}}}\Big) \sqrt{\tilde{\mu}} \tilde{\G}(\tilde{f},\tilde{f})
+ \tilde{\G}^{t}(\tilde{f},\tilde{f}) - \e^{-1}\pt_t
%------------
\Big(\frac{1}{\sqrt{\tilde{\mu}}}\Big) \sqrt{\tilde{\mu}} \tilde{L}\tilde{f}  \\
& - \e^{-1}\tilde{L}^{t} \tilde{f}  - \pt_t \Big( \frac{\pt_t \tilde{\mu}}{\sqrt{\tilde{\mu}}}\Big)
- \e \pt_t \Big(\frac{\pt_t \sqrt{\tilde{\mu}}}{\sqrt{\tilde{\mu}}}\Big) \tilde{f}
- \e \Big(\frac{\pt_t \sqrt{\tilde{\mu}}}{\sqrt{\tilde{\mu}}}\Big)  \pt_t \tilde{f},\\
%------------
{\tilde{\G}^{t}}(f,g)
%=& \frac{1}{\sqrt{\tilde{\mu}}}\Big[Q\big(f \pt_t \sqrt{\tilde{\mu}}, %\sqrt{\tilde{\mu}}g\big) + Q\big( \sqrt{\tilde{\mu}}f, g \pt_t \sqrt{\tilde{\mu}} %\big)\Big]
:=& \tilde{\G}\Big(\frac{\pt_t \sqrt{\tilde{\mu}}}{\sqrt{\tilde{\mu}}}f, g\Big) + \tilde{\G}\Big(f, \frac{\pt_t \sqrt{\tilde{\mu}}}{\sqrt{\tilde{\mu}}}g\Big),\\
%----------
{ \tilde{L}^{t}{f} }
%:=& - \frac{1}{\sqrt{\tilde{\mu}}}\Big[Q\big(\pt_{t} \tilde{\mu}, %\sqrt{\tilde{\mu}}{f}\big) + Q\big(\sqrt{\tilde{\mu}}{f}, \pt_{t} \tilde{\mu}\big) + %Q\big( \tilde{\mu}, f\pt_{t}\sqrt{\tilde{\mu}}\big) + Q\big(f \pt_{t}\sqrt{\tilde{\mu}}, %\tilde{\mu}\big) \Big], \\
:=& \tilde{\G}\Big(\frac{\pt_t \tilde{\mu}}{\sqrt{\tilde{\mu}}}, f\Big) + \tilde{\G}\Big(f, \frac{\pt_t \tilde{\mu}}{\sqrt{\tilde{\mu}}}\Big)
+\tilde{\G}\Big(\sqrt{\tilde{\mu}}, \frac{\pt_t \sqrt{\tilde{\mu}}}{\sqrt{\tilde{\mu}}}f\Big) + \tilde{\G}\Big(\frac{\pt_t \sqrt{\tilde{\mu}}}{\sqrt{\tilde{\mu}}}f, \sqrt{\tilde{\mu}}\Big),\\
%-----------
s := &\sqrt{2 \pi} \pt_t\Big(\frac{\mu}{\sqrt{\tilde{\mu}}}\Big) \int_{n \cdot v >0} \tilde{f} \sqrt{\tilde{\mu}}  [n \cdot v]  \dd v + \sqrt{2 \pi} \frac{\mu}{\sqrt{\tilde{\mu}}} \int_{n \cdot v >0} \tilde{f} \pt_t\sqrt{\tilde{\mu}}
  [n \cdot v]  \dd v.
\end{split}
\end{equation}
\medskip

%{\bf \color{purple} (Remark: For the statement of Proposition 4.11, one may directly use %the original equation for $\tilde{f}$ and $\pt_t\tilde{f}$, with the source terms $ %\tilde{g}$ and $\tilde{g}^{t} $ given (to be estimated in Subsection 4.4).)}

The main result of this subsection is the following energy estimate.
\medskip

\begin{proposition} \label{tildef tildeft - Energy estimate}\
 Let $\tilde{f} \in L^2(\R^+\times\Omega \times \R^3)$ be a solution of the perturbation equation \eqref{tildef - Boltzmanneq ch1} with given source $\tilde{g}$,  and let $\pt_t \tilde{f}\in L^2(\R^+\times\Omega \times \R^3)$ be a solution of \eqref{tildef t - Boltzmann with g} with given source $\tilde{g}^{t}$. Suppose the a priori assumption \eqref{theta-u-smallness-assumption} holds. Then the following estimates hold:
\begin{equation}
\begin{split} \label{tilde-f-energy-estimate}
& \normm{\tilde{f}(t)}_{L^2_{x,v}}^{2}  + \frac{\theta^2(t)}{\e^2} + \sum \frac{w_{i}^2(t)}{\e^2}+ \frac{1}{\e^2}\int_{0}^t\normm{\ipt \tilde{f}}_{L^2_{x,v}(\tilde{\nu})}^2 \dd \tau \\
& + \frac{\a}{\e} \int_{0}^t \iint_{\g_{+}} \norm{\tilde{f}}^2 \dd \g \dd \tau  + \frac{\a}{\e} \int_{0}^t \frac{\theta^2}{\e^2} \dd \tau  + \frac{\a}{\e} \int_{0}^t \sum \frac{w_{i}^2}{\e^2} \dd \tau \\
\lesssim & \normm{\tilde{f} (0)}_{L^2_{x,v}}^{2}  + {\a} \int_{0}^t \normm{\tilde{f}}_{2}^2 \dd \tau  + \int_{0}^t  \tilde{g} \tilde{f} \dd \tau +  \a \int_{0}^{t} \Big( \mathfrak{h}_{3} + \mathfrak{h}_{1}\norm{\tilde{f}}_{L^2_{\g_{+}}}^2 \Big)\dd \tau,
\end{split}
\end{equation}
and
\begin{equation} \label{tilde-ft-energy-estimate}
\begin{split}
& \normm{\pt_t \tilde{f}(t)}_{L^2_{x,v}}^{2}  + \frac{[\pt_t \theta(t)]^2}{\e^2} + \sum \frac{[\pt_t w_{i}(t)]^2}{\e^2} + \frac{1}{\e^2}\int_{0}^t \normm{\ipt \pt_t \tilde{f}}_{L^2_{x,v}(\tilde{\nu})}^2 \dd \tau \\
&  + \frac{\a}{\e} \int_{0}^t \iint_{\g_{+}} \norm{\pt_t \tilde{f}}^2 \dd \g \dd \tau  + \frac{\a}{\e} \int_{0}^t \frac{(\pt_t \theta)^2}{\e^2} \dd \tau  + \frac{\a}{\e} \int_{0}^t \sum \frac{(\pt_t w_{i})^2}{\e^2} \dd \tau \\
\lesssim & \normm{\pt_t \tilde{f} (0)}_{L^2_{x,v}}^{2} + {\a} \int_{0}^t \normm{\pt_t \tilde{f}}_{2}^2 \dd \tau  + \int_{0}^t  \tilde{g}^{t}  \pt_t \tilde{f} \dd \tau + \a \int_{0}^{t}\Big( \mathfrak{h}_{3,t} + \mathfrak{h}_{1,t}\norm{\pt_t \tilde{f}}_{L^2_{\g_{+}}}^2 \Big)\dd \tau,
\end{split}
\end{equation}
where $\mathfrak{h}_{n}\in \mathfrak{H}_{n}$ and  $\mathfrak{h}_{n,t}\in \mathfrak{H}_{n,t}$ for $n\in\mathbb{N}$.
\end{proposition}
\medskip

The estimate for the source terms $ \tilde{g}$ and $\tilde{g}^{t} $ on the right-hand side of \eqref{tilde-f-energy-estimate} and \eqref{tilde-ft-energy-estimate}  will be given in Subsection 4.4. Before giving the proof of Proposition \ref{tildef tildeft - Energy estimate}, we need some preparatory lemmas.

\medskip

Recall the linearized Boltzmann operator $\tilde{L}$ defined in \eqref{tildeg - definition} and its null space $\ker \tilde{L} $ defined in \eqref{ker-tilde-P}. It is standard that $\tilde{L} \tilde{f} = \tilde{\nu} \tilde{f} -\tilde{K} \tilde{f}$ (see e.g. \cite{cercignani1994, Guo2003}), where the collision frequency $\tilde{\nu}$ and
the compact operator $\tilde{K}$ on $L^2(\R^3_{v})$ are
\begin{equation}\label{collsion-frequency-def}
\begin{split}
&\tilde{\nu}=\tilde{\nu}(v) := \frac{1}{\sqrt{\tilde{\mu}}}Q_{-} (\sqrt{\tilde{\mu}}, \tilde{\mu}) = \int_{\R^{3}}\int_{\mathbb{S}^2} \norm{(v-u) \cdot \omega} \tilde{\mu}(u) \dd \omega \dd u,\\
%-------------
&\tilde{K} \tilde{f} = \frac{1}{\sqrt{\tilde{\mu}}} \left[Q_{+} (\tilde{\mu}, \sqrt{\tilde{\mu}} \tilde{f})+ Q_{+} (\sqrt{\tilde{\mu}} \tilde{f}, \tilde{\mu}) - Q_{-} (\tilde{\mu}, \sqrt{\tilde{\mu}} \tilde{f})\right]
=\int_{\R^{3}}[\tilde{k}_{1}(v,u)-\tilde{k}_{2}(v,u)] \tilde{f}(u) \dd u.
\end{split}
\end{equation}
For hard sphere cross sections, there exist positive constants $C_0$ and $C_1$ such that
\begin{equation*}
\begin{split}
\rho\sqrt{T} C_{0} \inn{v} \le \tilde{\nu}(v) \le \rho\sqrt{T} C_{1} \inn{v}.
\end{split}
\end{equation*}
If $\rho$, $\c$ and $T$ are bounded above and below, then
\begin{equation*}
\begin{split}
C_{0} \inn{v} \le \tilde{\nu} (v) \le C_{1} \inn{v},
\end{split}
\end{equation*}
so that $\nu(v) \approx\tilde{\nu}(v)$. Moreover, the operator $\tilde{L}$ is symmetric with
%\begin{equation*}
%\begin{split}
%\big(f,Lg\big)_{L^2_v}=\big(Lf,g\big)_{L^2_v}\; \text{ on }\; D_{L} = & \big\{f \in %L^2(\R^3_{v}) | \; \nu^{1/2} f \in L^{2}(\R^3_{v}) \big\},\\
%\big(f,\tilde{L}g\big)_{L^2_v}=\big(\tilde{L}f,g\big)_{L^2_v} \; \text{ on }\; %D_{\tilde{L}} = & \big\{\tilde{f} \in L^2(\R^3_{v}) |\; \tilde{\nu}^{1/2} \tilde{f} \in %L^{2}(\R^3_{v})\big \}.
%\end{split}
%\end{equation*}
 spectral inequality
\begin{equation}\label{positive-definiteness}
\begin{split}
\langle \tilde{f},\tilde{L}\tilde{f} \rangle_{2} \gtrsim& \normm{ \ipt \tilde{f}}_{L^2_v(\tilde{\nu})}^2\; \text{ for } \; \tilde{f}\in D_{\tilde{L}}=\big\{\tilde{f} \in L^2(\R^3_{v}) |\; \tilde{\nu}^{1/2} \tilde{f} \in L^{2}(\R^3_{v})\big \}.
\end{split}
\end{equation}
Using the relation
$$
\tilde{\mu}(v) = T^{-\frac{3}{2}}\rho \mu\Big(\frac{v-\c}{\sqrt{T}}\Big),
$$
 the implicit constant in  ``$\gtrsim$'' in \eqref{positive-definiteness} is uniform, provided $\rho$, $c$, and $T$ are bounded above and below.

\medskip

\begin{lemma} \label{r ptr s - smallness}
Let $r$ be defined as in \eqref{tildeg - definition} and $s$ as in \eqref{tildegt - definition}.  Under the a priori assumption \eqref{theta-u-smallness-assumption}, the following estimates hold:
\begin{align}
&\Big| r - \frac{\sqrt{\tilde{\mu}}}{\e }\Big(\big(2-\frac{\norm{v}^2}{2}\big)\theta - v\cdot \c\Big)\Big| _{L^2_{\g_{-}}} \lesssim \e \mathfrak{h}_{2}, \label{r-1order-estimate}\\
%---------------
&\norm{r}_{L^2_{\g_{-}}} \lesssim \mathfrak{h}_{1}+\e \mathfrak{h}_{2},\qquad \norm{r}_{L^\infty_{\g_{-}}} \lesssim \mathfrak{h}_{1}+\e \mathfrak{h}_{2}, \label{r-estimate}\\
%---------
&\Big|\pt_t r - \frac{\sqrt{\tilde{\mu}}}{\e }\Big(\big(2-\frac{\norm{v}^2}{2}\big)\pt_t \theta - v\cdot \pt_t \c\Big)\Big|_{L^2_{\g_{-}}} \lesssim \e \mathfrak{h}_{2,t},\label{rt-1order-estimate}\\
%---------------
&\norm{\pt_t r}_{L^2_{\g_{-}}} \lesssim  \a \mathfrak{h}_{1} + \a \norm{\tilde{f}}_{L^2_{\g_{+}}}+\e \mathfrak{h}_{2,t},\label{rt-estimate}\\
%--------------------------
&\norm{s}_{L^2_{\g_{-}}} \lesssim \a \Big(\mathfrak{h}_{1} + \norm{\tilde{f}}_{L^2_{\g_{+}}}\Big)\norm{\tilde{f}}_{L^2_{\g_{+}}} \lesssim \e \mathfrak{h}_{2,t}, \label{s-estimate}
%---------------
\end{align}
where $\mathfrak{h}_{n}\in \mathfrak{H}_{n}$ and $\mathfrak{h}_{2,t}\in \mathfrak{H}_{2,t}$.
\end{lemma}

\begin{proof}[\textbf{Proof.}] \
By direct calculation,
\begin{equation}\label{r-expression}
\begin{split}
r %=& \frac{1}{\e \sqrt{\tilde{\mu}}}\big(\rho T^{\frac{1}{2}} \mu - \tilde{\mu}\big)
%---------
=&\frac{\sqrt{\tilde{\mu}}}{\e }\big(\rho T^{\frac{1}{2}}\frac{\mu}{{\tilde{\mu}}} - 1\big)
%---------
=\frac{\sqrt{\tilde{\mu}}}{\e}\Big[(1+\theta)^2\exp \Big(-\frac{\norm{v}^2}{2(1+\theta)}\theta + \frac{-2v \cdot \c + \norm{\c}^2}{2(1+\theta)}\Big) - 1\Big]  \\
%---------
%= &\frac{\sqrt{\tilde{\mu}}}{\e %}\Big[\Big((1+2\theta)\big(1-\frac{\norm{v}^2}{2}\theta\big) - v \cdot \c\Big) - 1 + %O(\theta^2,\norm{\c}^2)\Big] \\
%---------
= &\frac{\sqrt{\tilde{\mu}}}{\e }\Big[\Big(2-\frac{\norm{v}^2}{2}\Big)\theta - v\cdot \c + O(\theta^2, \norm{\c}^2) p(v)\Big],
\end{split}
\end{equation}
where $p(v)$ is a polynomial in $v$ and $O(\theta^2,\norm{\c}^2)$ denotes terms bounded by $\theta^2+\norm{\c}^2$. Using the exponential decay of $\sqrt{\tilde{\mu}}$, we obtain
\begin{equation*}
\norm{\frac{\sqrt{\tilde{\mu}}}{\e }O(\theta^2, \norm{\c}^2)p(v)} \lesssim \e \mathfrak{h}_{2},
\end{equation*}
which yields  \eqref{r-1order-estimate} and \eqref{r-estimate}.

For the time derivative of $r$, we have
\begin{equation*}
\begin{split}
\pt_t r = \frac{\sqrt{\tilde{\mu}}}{\e }\Big[\Big(2-\frac{\norm{v}^2}{2}\Big)\pt_t \theta - v\cdot \pt_t \c + \pt_t O(\theta^2, \norm{\c}^2)\Big] +  \frac{\pt_t\sqrt{\tilde{\mu}}}{\e} O(\th,\norm{\c}),
\end{split}
\end{equation*}
where
\begin{equation}\label{pt-sqrt-mu}
\begin{split}
\pt_t \sqrt{\tilde{\mu}} = &\frac{\norm{v-\c}^2-3T}{2}\frac{\pt_t T}{2T^2} \sqrt{\tilde{\mu}}+ \frac{(v-\c)\cdot \pt_t \c}{2T} \sqrt{\tilde{\mu}} + \frac{\pt_t \rho}{2\rho} \sqrt{\tilde{\mu}}.
%= \sqrt{\tilde{\mu}} O(\pt_t \th, \pt_t \c, \pt_t \rho).
\end{split}\end{equation}
Using the bounds from \eqref{rho th u - smallness} and the exponential decay of $\sqrt{\tilde{\mu}}$, we obtain \eqref{rt-1order-estimate} and \eqref{rt-estimate}.

Now consider $s$ as defined in \eqref{tildegt - definition}. From \eqref{pt-sqrt-mu} and
%$$
%    \pt_t \Big(\frac{1}{\sqrt{\tilde{\mu}}}\Big) = \frac{1}{\sqrt{\tilde{\mu}}} O(\pt_t %\th, \pt_t \c, \pt_t \rho),
%  $$
\eqref{rho th u - smallness}, we obtain
\begin{equation*}
\begin{split}
&\Big|\pt_t\Big(\frac{\mu}{\sqrt{\tilde{\mu}}}\Big) \int_{n \cdot v >0} \sqrt{\tilde{\mu}} \tilde{f} [n \cdot v]  \dd v\Big|_{L^2_{\g_{-}}}
\lesssim \a \Big(\mathfrak{h}_{1} + \norm{\tilde{f}}_{L^2_{\g_{+}}}\Big)\norm{\tilde{f}}_{L^2_{\g_{+}}}, \\
%-----------------
&\Big|\frac{\mu}{\sqrt{\tilde{\mu}}} \int_{n \cdot v >0}\pt_t(\sqrt{\tilde{\mu}}) \tilde{f} [n \cdot v]  \dd v\Big|_{L^2_{\g_{-}}}
\lesssim \a \Big(\mathfrak{h}_{1} + \norm{\tilde{f}}_{L^2_{\g_{+}}}\Big)\norm{\tilde{f}}_{L^2_{\g_{+}}}.
\end{split}
\end{equation*}
Combining these estimates yields  \eqref{s-estimate}.
\end{proof}
\medskip

The following near-orthogonality properties hold for $\tilde{\mathscr{P}}_{\g} \tilde{f}$,
$(1-\tilde{\mathscr{P}}_{\g}\big) \tilde{f}$, $( \norm{v}^2-4) \sqrt{\tilde{\mu}}$ and
$(v \cdot Ax) \sqrt{\tilde{\mu}}$.

\begin{lemma} \label{tildef - perpendicular}
Let $\tilde{f} \in L^2(\gamma)$ with $\tilde{\mathscr{P}}{\gamma}$ defined as in \eqref{tildeg - definition}.  Under the a priori assumption \eqref{theta-u-smallness-assumption}, the following estimates hold:
\begin{align}\label{}
&\Big| \iint_{\g_{+}}\big[\tilde{\mathscr{P}}_{\g} \tilde{f}\big] \big[(1-\tilde{\mathscr{P}}_{\g}
) \tilde{f}\big] \dd \g \Big|
 \lesssim \e \mathfrak{h}_{1} \norm{\tilde{f}}_{L^2_{\g_{+}}}^2, \label{Pg-almost-perpendicular-1minus-Pg}\\
%--------------
&\int_{n \cdot v >0} \big[\tilde{\mathscr{P}}_{\g} \tilde{f}\big]  \big[( \norm{v}^2-4) \sqrt{\tilde{\mu}}\big] [n \cdot v] \dd v =0, \label{Pg-perpendicular-v2-4}\\
%--------------
&\int_{n \cdot v >0} \big[\tilde{\mathscr{P}}_{\g} \tilde{f}\big]  \big[(v \cdot Ax) \sqrt{\tilde{\mu}}\big] [n \cdot v] \dd v =0,\label{Pg-perpendicular-v}\\
%--------------
&\Big| \int_{n \cdot v >0}  \big[( \norm{v}^2 -4) \sqrt{\tilde{\mu}}\big] \big[(v \cdot Ax) \sqrt{\tilde{\mu}}\big] [n \cdot v] \dd v \Big|  \lesssim \e \mathfrak{h}_{1},\label{v2-4-almost-perpendicular-v}
\end{align}
where $Ax \in \mathcal{R}_\Omega$ and $\mathfrak{h}_{1}\in \mathfrak{H}_{1}$.
\end{lemma}

\begin{proof}[\textbf{Proof.}] \
From the definition of $\mathscr{P}_{\g}$ and a direct computation,
\begin{align*}
 &\iint_{\g_{+}} \big[\tilde{\mathscr{P}}_{\g}\tilde{f}\big] \big[(1- \tilde{\mathscr{P}}_{\g})\tilde{f}\big]  \dd \g\\
%-----------------
%=& \int_{\pt\Omega}\Big(
%\int_{n \cdot v>0} \tilde{f} \big[\tilde{\mathscr{P}}_{\g}\tilde{f}\big] [n \cdot v] \dd %v -
%\int_{n \cdot v>0} [\tilde{\mathscr{P}}_{\g}\tilde{f}]^2 [n \cdot v] \dd v \Big) \dd S_x %\\
%-----------------
%=&\sqrt{2 \pi}  \iint_{\g_+} \sqrt{\tilde{\mu}} \tilde{f} \bigg[
%\int_{n\cdot v >0}  \frac{\tilde{f} \mu}{\sqrt{\tilde{\mu}}}  [n \cdot v] \dd v    - %\sqrt{2 \pi}\int_{n \cdot v >0} \sqrt{\tilde{\mu}} \tilde{f} [n \cdot v]  \dd v
%\int_{n \cdot v >0} \frac{\mu^2}{\tilde{\mu}} [n \cdot v] \dd v \bigg] \dd \g \\
%-----------------
%=& \sqrt{2 \pi}  \iint_{\g_{+}} \sqrt{\tilde{\mu}} \tilde{f}  \bigg[
%%\int_{n\cdot v >0} \sqrt{\tilde{\mu}}\tilde{f}\big[1+ %\big(\frac{\mu}{\tilde{\mu}}-1\big)\big]  [n \cdot v] \dd v \\
%& \qquad\qquad\qquad\;\;  - \sqrt{2 \pi}\int_{n \cdot v >0} \sqrt{\tilde{\mu}} \tilde{f} %[n \cdot v]  \dd v
%\int_{n \cdot v >0} \mu\big[1+ \big(\frac{\mu}{\tilde{\mu}}-1\big)\big] [n \cdot v] \dd v %\bigg] \dd \g \\
%-----------------
= &\sqrt{2 \pi}   \iint_{\g_{+}} \sqrt{\tilde{\mu}} \tilde{f}  \Big[
\int_{n\cdot v >0} \sqrt{\tilde{\mu}}\tilde{f} \big(\frac{\mu}{\tilde{\mu}}-1\big)  [n \cdot v] \dd v \Big] \dd \g\\
&  -\sqrt{2 \pi}   \iint_{\g_{+}} \sqrt{\tilde{\mu}} \tilde{f} \Big[ \sqrt{2 \pi}\int_{n \cdot v >0} \sqrt{\tilde{\mu}} \tilde{f} [n \cdot v]  \dd v
\int_{n \cdot v >0} \mu \big(\frac{\mu}{\tilde{\mu}}-1\big) [n \cdot v] \dd v \Big] \dd \g.
\end{align*}
Since
$\frac{\mu}{\tilde{\mu}}=1+O(|\th|,|\c|)$, we have
\begin{equation*}
\begin{split}
&\Big| \int_{n \cdot v >0} \mu \big(\frac{\mu}{\tilde{\mu}}-1\big) [n \cdot v] \dd v \Big| \lesssim \e \mathfrak{h}_{1}, \quad
\Big(\int_{n \cdot v >0} \tilde{\mu} \big(\frac{\mu}{\tilde{\mu}}-1\big)^2 [n \cdot v] \dd v \Big)^{\frac{1}{2}} \lesssim \e\mathfrak{h}_{1},
\end{split}
\end{equation*}
with $\mathfrak{h}_{1} \in \mathfrak{H}_{1}$.
This proves \eqref{Pg-almost-perpendicular-1minus-Pg}.

Next, using Lemma \ref{boundaryinteg} and the fact that $Ax\cdot n|_{\pt\O}=0$,
\begin{equation*}
\begin{split}
&\int_{n \cdot v >0} \big[ \tilde{\mathscr{P}}_{\g} \tilde{f}\big] \big[( \norm{v}^2-4) \sqrt{\tilde{\mu}}\big] [n \cdot v] \dd v
= \sqrt{2 \pi} \int_{n \cdot v >0} \sqrt{\tilde{\mu}} \tilde{f} [n \cdot v]  \dd v
\int_{n \cdot v >0} ( \norm{v}^2-4) {\mu} [n \cdot v]   \dd v = 0,\\
%------------------
&\int_{n \cdot v >0} \big[ \tilde{\mathscr{P}}_{\g} \tilde{f}\big] \big[(v \cdot Ax)  \sqrt{\tilde{\mu}}\big] [n \cdot v] \dd v
= \sqrt{2 \pi}  \int_{n \cdot v >0} \sqrt{\tilde{\mu}} \tilde{f} [n \cdot v]  \dd v
\int_{n \cdot v >0} (v \cdot Ax) {\mu}  [n \cdot v] \dd v = 0,\\
%------------------
&\Big|\int_{n \cdot v >0}  \big[ ( \norm{v}^2 -4) \sqrt{\tilde{\mu}}\big] \big[ (v \cdot Ax) \sqrt{\tilde{\mu}}\big] [n \cdot v] \dd v \Big|
= \Big| \frac{6  (Ax\cdot \c) {\rho}T^{\frac{3}{2}}}{\sqrt{2 \pi}} + \frac{(Ax\cdot \c)\norm{\c}^2{\rho}T^{\frac{1}{2}}}{\sqrt{2 \pi}}\Big| \leq   \e\mathfrak{h}_{1}.
\end{split}
\end{equation*}
These identities give \eqref{Pg-perpendicular-v2-4}, \eqref{Pg-perpendicular-v} and \eqref{v2-4-almost-perpendicular-v}.
\end{proof}
\medskip

We now prove Proposition \ref{tildef tildeft - Energy estimate}.
\medskip

\noindent
\begin{proof}[\textbf{Proof of Proposition \ref{tildef tildeft - Energy estimate}.}] \
 The proof is divided into three steps. Steps 1 and 2 establish the energy estimates for $\tilde{f}$ and $\pt_t\tilde{f}$, respectively.  Step 3 completes the energy estimates by incorporating the trace lemma.
 \medskip

\noindent\textbf{Step 1.  Energy estimate for $\tilde{f}$.} \;

We first derive the following estimate for  $\tilde{f}$:
\begin{align}
&\frac{1}{2} \pt_t \normm{\tilde{f}}_{L^2_{x,v}}^2 + \frac{1}{\e^2} \iint_{\Omega \times \R^3} \tilde{f} \tilde{L} \tilde{f}  \dd v \dd x+ \frac{3}{2} \pt_t \int_{\Omega} \frac{\theta^2}{\e^2} \dd x + \pt_t \int_{\Omega} \frac{\norm{\c}^2}{\e^2} \dd x  \nonumber \\
%--------------
& +  \frac{\a(2-\a)}{\e}  \iint_{\g_{+}} \Big( \frac{1}{2}\frac{\theta}{\e}( \norm{v}^2 -4) \sqrt{\tilde{\mu}}   + \frac{\c}{\e}\cdot v \sqrt{\tilde{\mu}}  +  [(1- \tilde{\mathscr{P}}_{\g})\tilde{f}]  \Big)^2 \dd \g \label{tildef - Energy estimtae positive form}\\
%--------------
\le & \frac{1}{\e} \Big| \iint_{\Omega \times \R^3} \tilde{f}  \tilde{g} \dd x \dd v  \Big| + \a \mathfrak{h}_{3} + \a \mathfrak{h}_{1}\norm{\tilde{f}}_{L^2_{\g_{+}}}^2,\nonumber
%--------------------------------------------
\end{align}
where $\mathfrak{h}_{n}\in \mathfrak{H}_{n}$ for $n\in\mathbb{N}$.

Standard $L^2$ energy estimate for \eqref{tildef - Boltzmanneq ch1} yields
\begin{align*}
\e \frac{1}{2} \pt_{t} \normm{\tilde{f}}_{L^2_{x,v}}^2 + \frac{1}{2} \iint_{\g} \tilde{f}^2 [n \cdot v] \dd v \dd S_x  + \e^{-1} \iint_{\Omega \times \R^3} \tilde{f} \tilde{L} \tilde{f} \dd v \dd x = \iint_{\Omega \times \R^3} \tilde{f} \tilde{g} \dd v \dd x.
\end{align*}
%Using the boundary condition, we have
%\begin{align*}
%&\iint_{\g} \tilde{f}^2 [n \cdot v] \dd v\dd S_x
%=&  \iint_{\g_{+}} \tilde{f}^2 [n \cdot v] \dd v \dd S_x + \iint_{\g_{-}} \tilde{f}^2 [n %\cdot v] \dd v\dd S_x \\
%=\iint_{\g_{+}} \tilde{f}^2 \dd \g - \iint_{\g_{-}} \big[(1- \a) \mathscr{R} \tilde{f} + %\a \tilde{\mathscr{P}}_{\g} \tilde{f} + \a r\big]^{2}  \dd \g.
%\end{align*}
Using the boundary condition and the change of variables $R_x v \mapsto v$,
\begin{align*}
\iint_{\g_{-}} \tilde{f}^2  \dd \g
%=& -\iint_{\g_{-}} \big[(1- \a) \mathscr{R}  \tilde{f} + \a \tilde{\mathscr{P}}_{\g} %\tilde{f} +\a r\big]^2  \dd \g \\
%---------------------
=&  \iint_{\g_{+}} \big[(1- \a)(1-\tilde{\mathscr{P}}_{\g}) \tilde{f} + \tilde{\mathscr{P}}_{\g} \tilde{f} +\a r\big]^2 \dd \g  \\
%---------------------
=&  \iint_{\g_{+}} \Big\{(1- \a)^2\big[(1-\tilde{\mathscr{P}}_{\g}) \tilde{f}\big]^2 + \big[\tilde{\mathscr{P}}_{\g} \tilde{f}\big]^2 +\a^2r^2 + 2\a r\tilde{\mathscr{P}}_{\g} \tilde{f} \\
%---------------------
& \qquad \quad  + 2(1-\a)\big[(1-\tilde{\mathscr{P}}_{\g}) \tilde{f}\big] \big[\tilde{\mathscr{P}}_{\g} \tilde{f}\big] + 2\a(1- \a)r\big[(1-\tilde{\mathscr{P}}_{\g}) \tilde{f}\big] \Big\}  \dd \g.
\end{align*}
Applying Lemma \ref{r ptr s - smallness} and Lemma \ref{tildef - perpendicular}, we  obtain the intermediate estimate
\begin{equation}\label{tildef - energy mid step}
\begin{split}
&\frac{1}{2} \pt_t \normm{\tilde{f}}_{L^2_{x,v}}^2 + \frac{1}{\e^2}   \iint_{\Omega \times \R^3} \tilde{f} \tilde{L} \tilde{f} \dd v \dd x + \frac{\a(2-\a)}{\e}\iint_{\g_+} \big[(1- \tilde{\mathscr{P}}_{\g})\tilde{f}\big]^2  \dd \g   \\
%-------------------------
& +\frac{\a (1-\a)}{\e} \iint_{\g_+} \Big[ ( \norm{v}^2 -4)\sqrt{\tilde{\mu}} \frac{\theta}{\e} (1- \tilde{\mathscr{P}}_{\g})\tilde{f}   +  2\frac{\c}{\e}\cdot v \sqrt{\tilde{\mu}} (1- \tilde{\mathscr{P}}_{\g})\tilde{f}\Big]    \dd \g \\
%-------------------------
& - \frac{\a^2}{4\e} \iint_{\g_+} \Big[ ( {\norm{v}^2} -4)^2 \tilde{\mu}\frac{\theta^2}{\e^2} + \big( 2\frac{\c}{\e}\cdot v\big)^2 \tilde{\mu} \Big] \dd \g  \\
%-------------------------
\le  &\frac{1}{\e}\Big| \iint_{\Omega \times \R^3} \tilde{f} \tilde{g} \dd v \dd x \Big| + \a \mathfrak{h}_{3} + 2\a \mathfrak{h}_{1}\norm{\tilde{f}}^2_{L^2_{\g_{+}}}.
\end{split}
\end{equation}

From Lemma \ref{boundaryinteg} we compute
\begin{equation}\label{th-u-boundary-calculation}
\begin{split}
&\frac{4}{\sqrt{2 \pi}} \int_{\pt \Omega} \theta^2  \dd S_x = \frac{1}{2}\iint_{\g_{+}} \theta^2 ( \norm{v}^2 -4)^2 {\tilde{\mu}}  \dd \g + \e^{3} \mathfrak{h}_{3}, \\
%-------------------------
&\frac{4}{\sqrt{2 \pi}} \int_{\pt \Omega} (\pt_t \theta)^2  \dd S_x = \frac{1}{2}\iint_{\g_{+}} (\pt_t \theta)^2 ( \norm{v}^2 -4)^2 {\tilde{\mu}}  \dd \g+ \e^{3} \mathfrak{h}_{3,t},\\
%-------------------------
&\frac{4}{\sqrt{2 \pi}} \int_{\pt \Omega} \norm{\c}^2  \dd S_x = \iint_{\g_{+}} (\c \cdot v)^2 {\tilde{\mu}}  \dd \g + \e^{3} \mathfrak{h}_{3}, \\
%-------------------------
&\frac{4}{\sqrt{2 \pi}}\int_{\pt \Omega} \norm{\pt_t \c}^2  \dd S_x = \iint_{\g_{+}} (\pt_t \c \cdot v)^2 {\tilde{\mu}}  \dd \g + \e^{3} \mathfrak{h}_{3,t}.
\end{split}
\end{equation}
Applying \eqref{th-u-boundary-calculation} to \eqref{th - 2 estimate} and \eqref{w - 2 estimate} in Proposition \ref{th tht w wt - ODE} gives
\begin{align}
& \Big|  \frac{3}{2} \pt_t \int_{\Omega} \theta^2 \dd x  + \frac{\a}{2\e}\iint_{\g_{+}} \theta^2 ( \norm{v}^2 -4)^2 {\tilde{\mu}}  \dd \g + {\a}\iint_{\g_{+}} ( \norm{v}^2 -4) \sqrt{\tilde{\mu}} \tilde{f} \th  \dd \g  \Big|
\le \a \e^{2}\mathfrak{h}_{3} + \a \e^{2} \mathfrak{h}_{2}\norm{\tilde{f}}_{L^2_{\g_{+}}},\label{th - 2 new estimate}\\
%----------------------------
&\bigg|\frac{1}{2}\pt_t \int_{\Omega} \norm{\c}^2 \dd x  + \frac{\a}{\e} \iint_{\g_{+}} (\c \cdot v)^2 {\tilde{\mu}}  \dd \g
+ {\a}\iint_{\g_{+}} (\c \cdot v) \sqrt{\tilde{\mu}} \tilde{f}  \dd \g \bigg|
\le  \a \e^2\mathfrak{h}_{3} + \a \e^{2} \mathfrak{h}_{2}\norm{\tilde{f}}_{L^2_{\g_{+}}}.\label{w - 2 new estimate}
\end{align}

Now consider the combination \eqref{tildef - energy mid step} + $\frac{1}{\e^{2}}$ \eqref{th - 2 new estimate} + $\frac{2}{\e^{2}}$\eqref{w - 2 new estimate}:
\begin{align*}
 &\text{Left-hand side of } \Big( \eqref{tildef - energy mid step} + \frac{1}{\e^{2}} \eqref{th - 2 new estimate} + \frac{2}{\e^{2}}\eqref{w - 2 new estimate}\Big) \\
%&\frac{1}{2} \pt_t \normm{\tilde{f}}_{L^2_{x,v}}^2 + \frac{1}{\e^2}   \iint_{\Omega %\times \R^3} \tilde{f} \tilde{L} \tilde{f} \dd v \dd x + \frac{\a(2-\a)}{\e}\iint_{\g_+} %\big[(1- \tilde{\mathscr{P}}_{\g})\tilde{f}\big]^2  \dd \g   \\
%-------------------------
%& +\frac{\a (1-\a)}{\e} \iint_{\g_+} \Big[ ( \norm{v}^2 -4)\sqrt{\tilde{\mu}} %\frac{\theta}{\e} (1- \tilde{\mathscr{P}}_{\g})\tilde{f}   +  2\frac{\c}{\e}\cdot v %\sqrt{\tilde{\mu}} (1- \tilde{\mathscr{P}}_{\g})\tilde{f}\Big]   \dd \g \\
%-------------------------
%& - \frac{\a^2}{4\e} \iint_{\g_+} \Big[ ( {\norm{v}^2} -4)^2 %\tilde{\mu}\frac{\theta^2}{\e^2} + \big( 2\frac{\c}{\e}\cdot v\big)^2 \tilde{\mu} \Big] %\dd \g \\
%-------------------------
%&+ \frac{1}{\e^2} \Big[ \frac{\a}{2\e}\iint_{\g_{+}} \theta^2 ( \norm{v}^2 -4)^2 %{\tilde{\mu}} \dd \g + {\a}\iint_{\g_{+}} ( \norm{v}^2 -4) \sqrt{\tilde{\mu}} \tilde{f} %\th  \dd \g \Big] \\
%--------------------
%&+\frac{2}{\e^{2}} \Big[ \frac{\a}{\e} \iint_{\g_{+}} (\c \cdot v)^2 {\tilde{\mu}} \dd \g
%+ {\a}\iint_{\g_{+}} (\c \cdot v) \sqrt{\tilde{\mu}} \tilde{f}  \dd \g\Big] \\
%--------------------
=&\frac{1}{2} \pt_t \normm{\tilde{f}}_{L^2_{x,v}}^2 + \frac{1}{\e^2}   \iint_{\Omega \times \R^3} \tilde{f} \tilde{L} \tilde{f} \dd v \dd x + \frac{\a(2-\a)}{\e}\iint_{\g_+} \big[(1- \tilde{\mathscr{P}}_{\g})\tilde{f}\big]^2  \dd \g  \\
%--------------------
%--------------------
& +\frac{\a(2-\a)}{4\e}   \iint_{\g_+} \Big[ ( {\norm{v}^2} -4)^2 {\tilde{\mu}} \frac{\theta^2}{\e^2} + \big( 2\frac{\c}{\e}\cdot v \big)^2 {\tilde{\mu}} \Big] \dd \g\\
& +\frac{\a(2-\a)}{\e}  \iint_{\g_+} \Big[ ( \norm{v}^2 -4) [(1- \tilde{\mathscr{P}}_{\g}) \tilde{f}] \sqrt{\tilde{\mu}} \frac{\theta}{\e}   + 2\frac{\c}{\e}\cdot v  [(1- \tilde{\mathscr{P}}_{\g}) \tilde{f}] \sqrt{\tilde{\mu}} \Big] \dd \g \\
& + \frac{\a (2-\a)}{4\e} \iint_{\g_+} \Big[ 2\frac{\c}{\e}\cdot v  ( \norm{v}^2 -4) \frac{\theta}{\e} {\tilde{\mu}} -  2\frac{\c}{\e}\cdot v  ( \norm{v}^2 -4) \frac{\theta}{\e} {\tilde{\mu}} \Big]  \dd \g \\
%--------------------
=&\frac{1}{2} \pt_t \normm{\tilde{f}}_{L^2_{x,v}}^2 + \frac{1}{\e^2}   \iint_{\Omega \times \R^3} \tilde{f} \tilde{L} \tilde{f} \dd v \dd x  \\
%-------------------------
& + \frac{\a(2-\a)}{\e}  \iint_{\g_{+}} \Big[ \frac{1}{2}( \norm{v}^2 -4) \sqrt{\tilde{\mu}}\frac{\theta}{\e}   + \frac{\c}{\e}\cdot v \sqrt{\tilde{\mu}}  +  [(1- \tilde{\mathscr{P}}_{\g})\tilde{f}]  \Big]^2  \dd \g \\
%--------------------
& - \frac{\a (2-\a)}{4\e} \iint_{\g_+} 2\frac{\c}{\e} \cdot v  (\norm{v}^2 -4) \frac{\theta}{\e} {\tilde{\mu}}   \dd \g.
\end{align*}
The last term is bounded by  $\a \mathfrak{h}_{3}$ thanks to Lemma \ref{tildef - perpendicular}. This establishes \eqref{tildef - Energy estimtae positive form}.
\medskip

\noindent\textbf{Step 2.  Energy estimate for $\pt_t\tilde{f}$.} \;

In an analogous way we obtain the corresponding estimate for $\pt_t\tilde{f}$:
\begin{align}
%--------------------------------------------
&\frac{1}{2} \pt_t \normm{\pt_t \tilde{f}}_{L^2_{x,v}}^2 + \frac{1}{\e^2} \iint_{\Omega \times \R^3} \pt_t \tilde{f} \tilde{L} (\pt_t \tilde{f}) \dd x \dd v + \frac{3}{2} \pt_t \int_{\Omega} \frac{(\pt_t \theta)^2}{\e^2} \dd x + \pt_t \int_{\Omega} \frac{\norm{\pt_t \c}^2}{\e^2} \dd x \nonumber\\
%--------------
& +  \frac{\a(2-\a)}{\e}  \iint_{\g_{+}} \Big( \frac{1}{2}\frac{\pt_t \theta}{\e}( \norm{v}^2 -4) \sqrt{\tilde{\mu}}   + \frac{\pt_t \c}{\e}\cdot v \sqrt{\tilde{\mu}}  +  [(1- \tilde{\mathscr{P}}_{\g})\pt_t \tilde{f}]  \Big)^2 \dd \g  \label{tildeft - Energy estimtae positive form}\\
%-----------
\le & \frac{1}{\e} \Big| \iint_{\Omega \times \R^3} \pt_t\tilde{f} \tilde{g}^{t} \dd x \dd v  \Big| + \a \mathfrak{h}_{3,t} + \a \mathfrak{h}_{1,t}\norm{\pt_t \tilde{f}}_{L^2_{\g_{+}}}^2,\nonumber
\end{align}
where $\mathfrak{h}_{n,t}\in \mathfrak{H}_{n,t}$ for $n\in\mathbb{N}$.

The derivation of \eqref{tildeft - Energy estimtae positive form} follows exactly the same pattern as Step 1, using \eqref{tildef t - Boltzmann with g}, Lemma \ref{r ptr s - smallness} and Lemma \ref{tildef - perpendicular} applied to $\partial_t\tilde{f}$, together with estimates \eqref{tht - 2 estimate} and \eqref{wt - 2 estimate}. We omit the repetitive details.
\medskip

\noindent\textbf{Step 3.  Completion of the energy estimates.} \;

Up to now the boundary dissipation has been controlled except for the directions of by $(|v|^2 -4)\sqrt{\tilde{\mu}}$, $Ax \cdot v \sqrt{\tilde{\mu}}$, and $\tilde{\mathscr{P}}_{\gamma}$. The remaining directions are handled via the trace lemma.

%We require the following trace Lemma 3.2  from \cite{Esposito2017}.
%\begin{lemma} \label{f - ukai lemma}
%For $f \in L^1([0,T] \times \Omega \times \R^3)$,
%\begin{align*}
%\int_{0}^T \iint_{\g_{+}^\d} \norm{f(t,x,v)} \dd \g \dd t
%\lesssim \; & \e \iint_{\Omega \times \R^3} \norm{f(0,x,v)} \dd v \dd x +   %\int_{0}^T\iint_{\Omega \times \R^3} \norm{f(t,x,v)} \dd v \dd x \dd t\\
%&+ \int_{0}^T\iint_{\Omega \times \R^3} \norm{\e \pt_t f + v \cdot \nabla_x f }(t,x,v) %\dd v \dd x \dd t,
%\end{align*}
%\end{lemma}

For this purpose, we decompose $\tilde{f}|{\gamma_{+}}$ according to the domain geometry:
\begin{equation}\label{split-f-on-gamma+}
\tilde{f}\big|_{\g_{+}} =\left\{
  \begin{array}{ll}
\displaystyle \tilde{\mathscr{P}}_{\g} \tilde{f} + \tilde{\mathscr{P}}_{\norm{v}^2-4} \tilde{f} + \tilde{\mathscr{P}}_{\perp} \tilde{f} & \text{ for non-axisymmetric domains},\\[0.1cm]
%-------------
\displaystyle\tilde{\mathscr{P}}_{\g} \tilde{f} + \tilde{\mathscr{P}}_{\norm{v}^2-4} \tilde{f} + \tilde{\mathscr{P}}_{v_{Ax}} \tilde{f} + \tilde{\mathscr{P}}_{\perp} \tilde{f} & \text{ for axisymmetric domains},\\[0.1cm]
%---------------
\displaystyle\tilde{\mathscr{P}}_{\g} \tilde{f} + \tilde{\mathscr{P}}_{\norm{v}^2-4} \tilde{f} + \sum_{i=1}^{3}\tilde{\mathscr{P}}_{v_{A_{i}x}} \tilde{f} + \tilde{\mathscr{P}}_{\perp} \tilde{f}  &\text{ for spherical domains}.
  \end{array}
\right.
\end{equation}
Here, $\tilde{\mathscr{P}}_{\g} \tilde{f}$ is defined in \eqref{tildeg - definition} and  the other projections are
\begin{align*}
\tilde{\mathscr{P}}_{\norm{v}^2-4} \tilde{f} :=& C_{\norm{v}^2-4} ( \norm{v}^2 -4)\sqrt{\tilde{\mu}}\int_{n \cdot v >0} \tilde{f}  ( \norm{v}^2 -4)\sqrt{\tilde{\mu}} [n \cdot v] \dd v, \\
\tilde{\mathscr{P}}_{v_{Ax}} \tilde{f} :=& C_{v_{Ax}} \left( v \cdot Ax \right)\sqrt{\tilde{\mu}}\int_{n \cdot v >0} \tilde{f}  v \cdot Ax \sqrt{\tilde{\mu}} [n \cdot v] \dd v,\\
\tilde{\mathscr{P}}_{v_{A_{i}x}} \tilde{f} :=& C_{v_{A_{i}x}} \left( v \cdot A_{i}x \right)\sqrt{\tilde{\mu}}\int_{n \cdot v >0} \tilde{f}  v \cdot A_{i}x \sqrt{\tilde{\mu}} [n \cdot v] \dd v, \quad i=1,2,3,
\end{align*}
with suitable normalization constants $C_{\norm{v}^2-4}, C_{v_{Ax}}$ and $C_{v_{A_{i}x}}$.

By Lemma \ref{tildef - perpendicular}, the terms in \eqref{split-f-on-gamma+} are nearly orthogonal:
\begin{equation}\label{almost-perpendicular-term}
\bigg|\int_{\g_{+}} \tilde{\mathscr{P}}_{X} \tilde{f} \tilde{\mathscr{P}}_{Y} \tilde{f} \dd \g\bigg| \lesssim \e \mathfrak{h}_{1}\norm{\tilde{f}}_{L^2_{\g_{+}}}^2\quad \text{ for } X\neq Y,
\end{equation}
where $X,Y \in \big\{\g, \norm{v}^2-4, v_{Ax}, v_{A_{i}x}, \perp \big\}$.
Using
\begin{equation*}
\begin{split}
\iint_{\g_{+}\backslash\g_{+}^\d} ( \norm{v}^2 -4)^2 {\tilde{\mu}} \dd \g
 +
\iint_{\g_{+}\backslash\g_{+}^\d} \left( v \cdot Ax\right)^2 {\tilde{\mu}} \dd \g
 +
\iint_{\g_{+}\backslash\g_{+}^\d} \frac{\mu^2}{\tilde{\mu}} \dd \g &\lesssim o(\d)\quad
\text{ for } Ax \in \mathcal{R}_\Omega,
\end{split}
\end{equation*}
 the near-grazing part is controlled by
\begin{align*}
&\int_{\g_{+}\backslash\g_{+}^\d} \norm{\tilde{\mathscr{P}}_{\g} \tilde{f} + \tilde{\mathscr{P}}_{\norm{v}^2-4} \tilde{f}
+ \sum \tilde{\mathscr{P}}_{v_{\mathcal{A}x}} \tilde{f}}^2 \dd \g
\lesssim o(\d)
\iint_{\g_{+}} \norm{\tilde{\mathscr{P}}_{\g} \tilde{f} + \tilde{\mathscr{P}}_{\norm{v}^2-4} \tilde{f}+ \sum \tilde{\mathscr{P}}_{v_{\mathcal{A}x}} \tilde{f}}^2 \dd \g.
\end{align*}
where $\sum \tilde{\mathscr{P}}_{v_{\mathcal{A}x}}$   ($\mathcal{A}=A$ or $A_i$) denotes the sum over the relevant axial directions.  Consequently,
\begin{align*}
&\iint_{\g_{+}} \norm{\tilde{\mathscr{P}}_{\g} \tilde{f} + \tilde{\mathscr{P}}_{\norm{v}^2-4} \tilde{f}+ \sum \tilde{\mathscr{P}}_{v_{\mathcal{A}x}} \tilde{f}}^2 \dd \g\\
= &\Big\{  \iint_{\g_{+}^\d} + \iint_{\g_{+}\backslash\g_{+}^\d}\Big \}\norm{\tilde{\mathscr{P}}_{\g} \tilde{f} + \tilde{\mathscr{P}}_{\norm{v}^2-4} \tilde{f} + \sum \tilde{\mathscr{P}}_{v_{\mathcal{A}x}} \tilde{f}}^2 \dd \g\\
\lesssim & \iint_{\g_{+}^\d} \norm{\tilde{\mathscr{P}}_{\g} \tilde{f} + \tilde{\mathscr{P}}_{\norm{v}^2-4} \tilde{f} + \sum \tilde{\mathscr{P}}_{v_{\mathcal{A}x}} \tilde{f}}^2 \dd \g\\
%---------------------
\le &\iint_{\g_{+}^\d} \norm{\tilde{f}}^2 \dd \g
+ 2 \iint_{\g_{+}^\d} \Big({\tilde{\mathscr{P}}_{\g} \tilde{f}
+ \tilde{\mathscr{P}}_{\norm{v}^2-4} \tilde{f}
+ \sum \tilde{\mathscr{P}}_{v_{\mathcal{A}x}}\tilde{f}}\Big){\tilde{\mathscr{P}}_{\perp} \tilde{f}} \dd \g
+ 2\iint_{\g_{+}^\d} \norm{\tilde{\mathscr{P}}_{\perp} \tilde{f}}^2 \dd \g \\
\le &\iint_{\g_{+}^\d} \norm{\tilde{f}}^2 \dd \g + 2\iint_{\g_{+}} \norm{\tilde{\mathscr{P}}_{\perp} \tilde{f}}^2 \dd \g + \e \mathfrak{h}_{1}\norm{\tilde{f}}^2_{L^2_{\g_{+}}},
\end{align*}
where we used \eqref{almost-perpendicular-term} in the last inequality.

Applying the trace lemma (Lemma 3.2 in \cite{Esposito2017}) to the non-grazing part, we obtain
\begin{equation} \label{tildef - Bundary ukai lemma}
\begin{split}
&\int_{0}^t \iint_{\g_{+}} \norm{\tilde{\mathscr{P}}_{\g} \tilde{f} + \tilde{\mathscr{P}}_{\norm{v}^2-4} \tilde{f}+ \sum \tilde{\mathscr{P}}_{v_{\mathcal{A}x}} \tilde{f}}^2 \dd \g \dd s\\
%-------------------
\lesssim & \e \iint_{\Omega \times \R^3} \norm{\tilde{f}(0)}^2 \dd v \dd x +  \int_{0}^t\iint_{\Omega \times \R^3} \norm{\tilde{f}(s)}^2 \dd v \dd x \dd s
+ {\e \int_{0}^{t}\mathfrak{h}_{1}\norm{\tilde{f}}^2_{L^2_{\g_{+}}} \dd s}\\
&+ \int_{0}^t\iint_{\Omega \times \R^3} \norm{(\tilde{g} -\e^{-1}\tilde{L}\tilde{f})\tilde{f} } \dd v \dd x \dd s + \int_{0}^{t}\iint_{\g_{+}} \norm{\tilde{\mathscr{P}}_{\perp} \tilde{f}}^2 \dd \g \dd s.\\
\end{split}
\end{equation}
Note that by \eqref{split-f-on-gamma+},
\begin{align*}
&\iint_{\g_{+}} \Big[ \frac{1}{2}\frac{\theta}{\e}( \norm{v}^2 -4) \sqrt{\tilde{\mu}}   + \frac{\c}{\e}\cdot v \sqrt{\tilde{\mu}}  +  [(1- \tilde{\mathscr{P}}_{\g})\tilde{f}]  \Big]^2  \dd \g \\
=  &\int_{\g_{+}}\Big[ \frac{1}{2}\frac{\theta}{\e}( \norm{v}^2 -4) \sqrt{\tilde{\mu}}  + \tilde{\mathscr{P}}_{( \norm{v}^2 -4) }\tilde{f}  + \frac{\c}{\e}\cdot v \sqrt{\tilde{\mu}}  + \sum \tilde{\mathscr{P}}_{v_{\mathcal{A}x} }\tilde{f} + \tilde{\mathscr{P}}_{\perp }\tilde{f}  \Big]^2  \dd \g.
\end{align*}
Combining \eqref{tildef - Energy estimtae positive form} with $\d \frac{\a}{\e}\times \eqref{tildef - Bundary ukai lemma}$ for a sufficiently small $\d>0$, we derive
the desired estimate \eqref{tilde-f-energy-estimate}.

The energy estimate \eqref{tilde-ft-energy-estimate} for $\partial_t \tilde{f}$ follows in the same way from \eqref{tildeft - Energy estimtae positive form}.
\end{proof}
%\begin{equation} \label{tildeft - Bundary ukai lemma}
%\begin{split}
%&\int_{0}^t \iint_{\g_{+}} \norm{\tilde{\mathscr{P}}_{\g} \pt_t \tilde{f} + %\tilde{\mathscr{P}}_{\norm{v}^2-4} \pt_t \tilde{f}+ \sum %\tilde{\mathscr{P}}_{v_{\mathcal{A}x}} \pt_t \tilde{f}}^2 \dd \g \dd s \\
%\lesssim  & \e \iint_{\Omega \times \R^3} \norm{\pt_t \tilde{f}(0)}^2 \dd v \dd x +  %\int_{0}^t\iint_{\Omega \times \R^3} \norm{\pt_t \tilde{f}(s)}^2 \dd v \dd x \dd s
% + {\e \int_{0}^{t}\mathfrak{h}_{1}\norm{\pt_t \tilde{f}}^2_{L^2_{\g_{+}}} \dd s}\\
%&+ \int_{0}^t\iint_{\Omega \times \R^3} \norm{(\tilde{g}^{t} -\e^{-1}\tilde{L}\pt_t %\tilde{f})\pt_t \tilde{f} } \dd v \dd x \dd t + \int_{0}^{t}\iint_{\g_{+}} %\norm{\tilde{\mathscr{P}}_{\perp} \pt_t \tilde{f}}^2 \dd \g \dd s.
%\end{split}
%\end{equation}
%Combining  \eqref{tildeft - Energy estimtae positive form} with $\d %\frac{\a}{\e}\times\eqref{tildeft - Bundary ukai lemma}$ and proceeding as above yields  %\eqref{tilde-ft-energy-estimate}. This completes the proof of Proposition \ref{tildef %tildeft - Energy estimate}.
\bigskip

%%%%%%%%%%%%%%%%%%%%%%%%%%%%%%%%%%
%%%%%%%%%%%%%%%%%%%%%%%%%%%%%%%%%%
\subsection{Macroscopic $L^2$ and $L^6$ Estimates}\
\medskip

In this subsection, we derive the macroscopic $L^2$ and $L^6$ estimates for the perturbation equation \eqref{tildef - Boltzmanneq ch1} and give the proof of Proposition \ref{Psi - L2 and L6 estimate}.

Recall the non-orthogonal basis $\{\tilde{\chi}_{i}\}_{i=0}^{4}$ of $\ker \tilde{L}$ defined in \eqref{base-tilde-chi} and the coefficients $\tilde{a}, \tilde{b}, \tilde{c}$ defined in \eqref{tilde-abc-def}.
By \eqref{abc - average zero}, the following compatibility conditions hold:
\begin{equation} \label{abc - consercation law}
\begin{split}
\int_{\Omega} \tilde{a}(t,x) \dd x = 0, \;\;
\int_{\Omega} A_jx \cdot \tilde{b}(t,x) \dd x= 0\; (A_jx\in\mathcal{R}_{\O}), \;\;
\int_{\Omega} \tilde{c}(t,x)  \dd x = 0 \;\;\;  \forall  t\geq 0.
\end{split}
\end{equation}
Define the Burnett functions
\begin{equation}\label{Burnett-function}
\begin{split}
&\tilde{A}_{ij}(v) := \Big(v_{i}v_{j} - \frac{\delta_{ij}}{3}\norm{v}^2 \Big)\sqrt{\tilde{\mu}},\;\;\;  \tilde{B}_{i}(v) := v_{i}\frac{\norm{v}^2-5}{\sqrt{10}} \sqrt{\tilde{\mu}},
\quad  i,j=1,2,3.
\end{split}
\end{equation}
By Lemma \ref{wholeinteg}, for every $i,j=1,2,3$,  the following almost orthogonality hold:
\begin{equation}\label{Burnett-orthogonaltilde}
\int_{\R^3} \tilde{\chi}_{k}(v) \tilde{A}_{ij}(v) \dd v= O(\norm{\theta}+ \norm{\c}), \;\;\; \int_{\R^3} \tilde{\chi}_{k}(v) \tilde{B}_i(v) \dd v
= O(\norm{\theta}+ \norm{\c}),
\quad  k=0,\cdots, 4.
\end{equation}

We now give the proof of Proposition \ref{Psi - L2 and L6 estimate}.
\medskip

\begin{proof}[\textbf{Proof of Proposition \ref{Psi - L2 and L6 estimate}.}] \
The proof follows a strategy similar to that of Proposition \ref{0-macro-L2-L6-estimate}, but here we work with the rotating Maxwellian $\tilde{\mu}$ and use the conservation laws of angular momentum and energy provided by \eqref{abc - average zero}. Moveover, the non-orthogonality of the basis $\{\tilde{\chi}_{i}\}_{i=0}^{4}$ introduces additional computational complexity.

We first multiply the  equation \eqref{tildef - Boltzmanneq ch1} by a test function $\tilde{\psi}_{p,q}$:
\begin{equation}\label{test-equation-uniform-form}
\begin{split}
&\underbrace { \e \iint_{\Omega \times \R^3}  \tilde{\psi}_{p,q}  \pt_t \tilde{f} \dd v \dd x }_{ :=\tilde{\Xi}_{p,q}^{1} }
+ \underbrace{\iint_{\g_{+}}  \tilde{\psi}_{p,q}  \tilde{f} \dd \g
 - \iint_{\g_{-}} \tilde{\psi}_{p,q} \tilde{f}  \dd \g }_{ :=\tilde{\Xi}_{p,q}^{2} }
  \underbrace{ - \iint_{\Omega \times \R^3}  \big( v \cdot \nabla_x  \tilde{\psi}_{p,q} \big) \tilde{f} \dd v \dd x}_{:=\tilde{\Xi}_{p,q}^{3}}\\
%----------------
= & \underbrace{\iint_{\Omega \times \R^3} \Big[\e^{-1} \tilde{\psi}_{p,q} \tilde{L}(\tilde{f},\tilde{f}) +  \tilde{\psi}_{p,q}  \tilde{g}  \Big] \dd v \dd x}_{:=\tilde{\Xi}_{p,q}^{4}},
\end{split}
\end{equation}
where we have used \eqref{tildemu - nabla zero ch1} to obtain $\tilde{\Xi}_{p,q}^{3}$.
The test function $\tilde{\psi}_{p,q}$ will be constructed  in the form of $\tilde{\psi}_{p,q}(t,x,v) = h(v){\phi}_{p,q}(t,x)\sqrt{\tilde{\mu}}$,
where $h(v)$ is a polynomial in $v$ and ${\phi}_{p,q}$ satisfies a suitable elliptic boundary value problem.

Note that $\mathscr{R}(\sqrt{\tilde{\mu}} )= \sqrt{\tilde{\mu}}$ by \eqref{tilde-mu-def} and \eqref{u-definition}. If $\tilde{\psi}_{p,q}$ also satisfies $\mathscr{R}(\tilde{\psi}_{p,q} )= \tilde{\psi}_{p,q}$,
 then the boundary term $\Xi_{p,q}^{2}$ in \eqref{test-equation-uniform-form} can be treated
 similarly to \eqref{0-Theta2-boundary-calculation-a} via the Maxwell boundary condition in \eqref{tildef - Boltzmanneq ch1} and the change of variables $v \mapsto R_x v$:
\begin{equation}\label{Theta2-boundary-calculation}
\begin{split}
\tilde{\Xi}_{p,q}^{2}
%=& \iint_{\g_{+}} \frac{\tilde{\psi}_{p,q}}{\sqrt{\tilde{\mu}}} \Tilde{\Phi}_{} \dd \g
%- \iint_{\g_{-}} \frac{\tilde{\psi}_{p,q}}{\sqrt{\tilde{\mu}}}  \Tilde{\Phi}_{} \dd \g \\
%-----------------------
%=& \iint_{\g_{+}} \tilde{\psi}_{p,q} \tilde{f} \dd \g
%- \iint_{\g_{-}} \tilde{\psi}_{p,q} \Big( (1-\a)  \mathscr{R}\tilde{f}
%+\a  \tilde{\mathscr{P}}_{\g} \tilde{f}\Big)\dd \g - \a \iint_{\g_{-}} \tilde{\psi}_{p,q} r %\dd \g \\
%-----------------------
%=& \iint_{\g_{+}} \tilde{\psi}_{p,q} \tilde{f} \dd \g
%-\iint_{\g_{+}} \mathscr{R}\big(\tilde{\psi}_{p,q} \big) \big( (1-\a) \tilde{f}
%+\a  \tilde{\mathscr{P}}_{\g}\tilde{f}  \big) \dd \g
%- \a \iint_{\g_{-}} \tilde{\psi}_{p,q} r \dd \g \\
%-----------------------
%=& \iint_{\g_{+}} \tilde{\psi}_{p,q} \tilde{f} \dd \g
%-\iint_{\g_{+}}  \tilde{\psi}_{p,q}  \big( (1-\a) \tilde{f}
%+\a  \tilde{\mathscr{P}}_{\g}\tilde{f}  \big)  \dd \g
%- \a \iint_{\g_{-}} \tilde{\psi}_{p,q} r \dd \g \\
%-----------------------
=& \a \iint_{\g_{+}} \tilde{\psi}_{p,q} \tilde{f} \dd \g
-\a \iint_{\g_{+}} \tilde{\psi}_{p,q} \tilde{\mathscr{P}}_{\g}\tilde{f} \dd \g
-\a \iint_{\g_{-}} \tilde{\psi}_{p,q} r  \dd \g, \quad  p\in \{a,b,c\}, q\in \{2, 6\}.
\end{split}
\end{equation}
%By H\"{o}lder's inequality,
% \begin{equation*}
%\begin{split}
%&\norm{\iint_{\g_{+}} \tilde{\psi}_{p,q} \tilde{f}  \dd \g} \lesssim \norm{\tilde{f} %}_{L^2_{\g_{+}}} \norm{{\phi}_{p,2}}_{L^2(\pt\O)}, \\
%------------------
%&\norm{\iint_{\g_{+}} \tilde{\psi}_{p,q} \tilde{\mathscr{P}}_{\g}\tilde{f} \dd \g} \lesssim  %\norm{ \tilde{f}  }_{L^2_{\g_{+}}} \norm{{\phi}_{p,2}}_{L^2(\pt\O)}.
%\end{split}
%\end{equation*}
For $\tilde{\Xi}^{2}_{p,2}$ ($p\in \{a,b,c\}$), the trace theorem gives
\begin{equation} \label{Theta2 - L2 estimate}
\begin{split}
\norm{\tilde{\Xi}^{2}_{p,2}}
\lesssim & \a \Big( \norm{\tilde{f} }_{L^2_{\g_{+}}} +\norm{r}_{L^2_{\g_{-}}} \Big) \norm{{\phi}_{p,2}}_{L^2{(\pt \Omega)}}
\lesssim   \a \Big( \norm{\tilde{f} }_{L^2_{\g_{+}}} +\norm{r}_{L^2_{\g_{-}}} \Big) \normm{{\phi}_{p,2}}_{H^1(\Omega)}.
\end{split}
\end{equation}
For $\tilde{\Xi}_{p,6}^{2}$ ($p\in \{a,b,c\}$),  using \eqref{Theta2-boundary-calculation} and deducing as in \eqref{0-Theta2 - L2 estimate-L6-a}, we obtain
\begin{equation}
\begin{split} \label{Theta2 - L6 estimate}
\norm{\tilde{\Xi}_{p,6}^{2}}
%\lesssim &\a \Big( \norm{\tilde{f}}_{L^4_{\g_{+}}} +\norm{r}_{L^4_{\g_{-}}} \Big) \normm{ %{\phi}_{p,6}}_{W^{1,\frac{6}{5}}(\Omega)}\\
%-----------
 \lesssim & \a  \Big( \norm{\tilde{f}}_{L^2_{\g_{+}}}^{\frac{1}{2}} \normm{\omega^{\frac{1}{2}} \tilde{f} }_{L^\infty_{x,v}}^{\frac{1}{2}} +\norm{r}_{L^4_{\g_{-}}} \Big) \normm{ {\phi}_{p,6}}_{W^{1,\frac{6}{5}}(\Omega)}.
\end{split}
\end{equation}
For $\tilde{\Xi}_{p,q}^{4}$ ($ p\in \{a,b,c\}$, $q\in \{2, 6\}$), H\"{o}lder's inequality directly yields
\begin{equation} \label{Theta4 - estimate}
\begin{split}
\norm{\tilde{\Xi}_{p,q}^{4}}
\lesssim & \Big(
  \e^{-1}\normm{\ipt \tilde{f} }_{L^{2}_{x,v}(\tilde{\nu})}
 +\normm{ \tilde{\nu}^{-\frac{1}{2}} \tilde{g}}_{L^{2}_{x,v}} \Big)  \normm{{\phi}_{p,q}}_{L^{2}_{x}}.
\end{split}
\end{equation}

To estimate  $\tilde{\P}\tilde{f}$, by Lemma \ref{Pf - abc similar estimate}, it suffices to control $\tilde{a}, \tilde{b}$ and $\tilde{c}$.
%The proof is divided into fours steps:  Step 1  estimates $\tilde{a}$-related terms %containing $\int_{s}^{t}\|\tilde{a}\|_{L^{2}_{x}}\dd \tau$ and $\|\tilde{a}\|_{L^{6}_{x}}$; %Step 2 and Step 3 focus on estimate of $\tilde{b}$-related and $c$-related terms; Step 4 %combines estimates from Step 1--3 to derive \eqref{P-tilde-f-macro-L2} and %\eqref{P-tilde-f-macro-L6}.
\medskip

\noindent\textbf{Step 1. Estimate for $\tilde{a}$.}

\noindent\textbf{Step 1.1.  Estimate for $\int_{s}^{t}\|\tilde{a}\|_{L^{2}_{x}}\dd \tau$ and $\|\tilde{a}\|_{L^{6}_{x}}$.}\\

In the weak formulation \eqref{test-equation-uniform-form}, we choose the test function
$$
  \tilde{\psi}_{a,q}(t,x,v) := \sum_{i=1}^{3} \pt_{i} \tilde{\varphi}_{a,q}(t,x) \big[\sqrt{10} \tilde{B}_{i}(v) - 5\tilde{\chi}_{i}(v)\big],  \quad q\in\{2, 6\}.
$$
 Here, by Lemma \ref{Poisson-equation-theory} and the  compatibility condition \eqref{abc - consercation law}, $\tilde{\varphi}_{a,2}(x)$ and $\tilde{\varphi}_{a,6}(x)$  are the unique solutions to the elliptic equations
  \begin{align}
- \Delta_x \tilde{\varphi}_{a,2} = \tilde{a} \;\text{ in }\O, \quad& {\pt_n} \tilde{\varphi}_{a,2}= 0  \;\text{ on }\pt\O,  \quad \int_{\Omega} \tilde{\varphi}_{a,2} \dd x =0, \label{a-2-elliptic-equation} \\
%---------------
- \Delta \tilde{\varphi}_{a,6} = \tilde{a}^{5} -\frac{1}{\norm{\Omega}}\int_{\Omega}\tilde{a}^{5} \dd x, \;\text{ in }\O, \quad& {\pt_{n}} \tilde{\varphi}_{a,6}= 0, \;\text{ on }\pt\O,  \quad \int_{\Omega} \tilde{\varphi}_{a,6} \dd x =0, \label{a-6-elliptic-equation}
\end{align}
 with the elliptic estimates
\begin{align}
\normm{\nabla^2 \tilde{\varphi}_{a,2}}_{L^2_x} + \normm{\nabla \tilde{\varphi}_{a,2}}_{L^2_x} + \normm{\tilde{\varphi}_{a,2}}_{L^2_x} &\lesssim \normm{\tilde{a}}_{L^2_x}, \label{a-2-elliptic-estimate}\\
%---------------
\normm{\nabla^2 \tilde{\varphi}_{a,6}}_{L^{\frac{6}{5}}_x} + \normm{\nabla \tilde{\varphi}_{a,6}}_{L^2_x} + \normm{\tilde{\varphi}_{a,6}}_{L^6_x}
&\lesssim \normm{\tilde{a}^{5}}_{L^{\frac{6}{5}}_x}=\normm{\tilde{a}}_{L^{6}_x}^5.
\label{a-6-elliptic-estimate}
\end{align}
%and Poincar\'{e} inequality
%\begin{equation}\label{Poincare-inequality-varphi-a}
%\begin{split}
%\normm{\tilde{\varphi}_{a,2}}_{L^2_x} \lesssim \normm{\nabla_x  %\tilde{\varphi}_{a,2}}_{L^2_x}.
%\end{split}\end{equation}

We now estimate each term in \eqref{test-equation-uniform-form}.
For $ \Tilde{\Xi}_{a,2}^{1}$, integration by parts yields
\begin{equation}\label{Theta-a2-01-estimate-form}
\begin{split}
\int_{s}^{t}   \Tilde{\Xi}_{a,2}^{1}
 = & \e \big[  \tilde{G}_{a} (t)- \tilde{G}_{a} (s)\big]
 - \e \int_{s}^{t}\iint_{\Omega \times \R^{3}} \sum_{i=1}^{3}\big [ \pt_t \pt_{i} \tilde{\varphi}_{a,2} + \pt_t\sqrt{\tilde{\mu}}\big]  v_i(|v|^2 - 10)\tilde{f}\\
 :=& \e \big[ \tilde{G}_{a} (t)- \tilde{G}_{a} (s)\big]- \Tilde{H}_{a,1} - \Tilde{H}_{a,2},
\end{split}
\end{equation}
By \eqref{a-2-elliptic-estimate} and Lemma \ref{Pf - abc similar estimate},
$\tilde{G}_{a} (t)$ is bounded by $\normm{\tilde{f}(t)}_{L^2_{x,v}}^2$.
For $ \Tilde{H}_{a,1}$,  we decompose $\tilde{\P} \tilde{f}$ as
\begin{equation}\label{tilde-P-tilde-f-tilde-base}
\begin{split}
\tilde{\P} \tilde{f}= \sum_{k=0}^{4} \langle \tilde{f}, \tilde{\chi}_{k}\rangle \tilde{\chi}_{k}
+ \Big(
\tilde{\P} \tilde{f} -\sum_{k=0}^{4} \langle \tilde{f}, \tilde{\chi}_{k} \rangle \tilde{\chi}_{k}\Big),
\end{split}
\end{equation}
and obtain
\begin{equation}\label{Theta-a2-01-estimate}
\begin{split}
  \int_{\R^{3}}   v_i(|v|^2 - 10)\sqrt{\tilde{\mu}} \tilde{\P} \tilde{f} \dd v
%---------------
%---------------
= &  \sum_{k=0}^{4} \langle   \tilde{f}, \tilde{\chi}_{k}  \rangle \int_{\R^{3}}   v_i(|v|^2 - 10)\sqrt{\tilde{\mu}}\tilde{\chi}_{k}\dd v  + \tilde{K}_{a,1} \\
%---------------
%---------------
= &  -5\tilde{b}_i + O(|\c|+|\th|)(\tilde{a}+ |\tilde{b}|+\tilde{c})  + \tilde{K}_{a,1},
\end{split}
\end{equation}
where we used \eqref{vi(|v|^2-10)-orthogonality}  for the velocity integral.  The remainder $\tilde{K}_{a,1}$ is bounded via Lemma \ref{expansion-in-tilde-chi}:
\begin{equation}\label{tilde-K-a-1}
\begin{split}
\norm{ \tilde{K}_{a,1} } = \Big|  \int_{\R^{3}}   v_i(|v|^2 - 10)\sqrt{\tilde{\mu}}\Big(
\tilde{\P} \tilde{f} -\sum_{k=0}^{4}  \langle  \tilde{f}, \tilde{\chi}_{k}  \rangle  \tilde{\chi}_{k}\Big)\dd v \Big|
\lesssim  \e  \mathfrak{h}_{1} \normm{ \tilde{\P}  \tilde{f}}_{L^2_v}.
\end{split}
\end{equation}
Then, using Lemma \ref{Pf - abc similar estimate},
\begin{equation}\label{Theta-a2-01-estimate-form-1}
\begin{split}
  \norm{  \Tilde{H}_{a,1} }
%---------------
\lesssim  \;\e  \int_{s}^{t}\normm{\pt_t \nabla_x \tilde{\varphi}_{a,2}}_{L^2_x} \Big[ \normm{\tilde{b}}_{L^2_x}  + \e   \mathfrak{h}_{1} \normm{\tilde{\P} \tilde{f} }_{L^2_{x,v}}+ \normm{\ipt \tilde{f} }_{L^2_{x,v}}\Big].
\end{split}
\end{equation}
The term $\Tilde{H}_{a,2}$ is controlled using \eqref{pt-sqrt-mu}, \eqref{rho th u - smallness} and \eqref{a-2-elliptic-estimate}:
\begin{equation}\label{Theta-a2-01-estimate-form-2}
\begin{split}
    \norm{  \Tilde{H}_{a,2}  }
\lesssim &
 \e\int_{s}^{t}( \norm{ \pt_t \theta } + \norm{ \pt_t w }  ) \normm{ \tilde{f}}_{L^{2}_{x,v}}  \normm{\nabla_x \tilde{\varphi}_{a,2}}_{L^{2}_{x}}
 \lesssim   \e \a \int_{s}^{t}\Big(  \mathfrak{h}_{1} + \norm{\tilde{f}}_{L^2_{\g_{+}}}\Big) \normm{ \tilde{f}}_{L^{2}_{x,v}}\normm{\tilde{a}}_{L^{2}_{x}}.
\end{split}
\end{equation}
Combining \eqref{Theta-a2-01-estimate-form}, \eqref{Theta-a2-01-estimate-form-1}  and \eqref{Theta-a2-01-estimate-form-2} yields
\begin{equation} \label{Theta1 - a estimate}
\begin{split}
\norm{\int_{s}^{t}    \Tilde{\Xi}_{a,2}^{1}   }
 \leq  & \, \e  \big[\tilde{G}_{a} (t)-  \tilde{G}_{a} (s)\big]
 +\e \a  \int_{s}^{t} \Big(  \mathfrak{h}_{1} +  \norm{\tilde{f}}_{L^2_{\g_{+}}}\Big) \normm{ \tilde{f}}_{L^{2}_{x,v}}\normm{\tilde{a}}_{L^{2}_{x}} \\
 %------------------
 &+ \e  \int_{s}^{t}\normm{\pt_t \nabla_x \tilde{\varphi}_{a,2}}_{L^2_x} \Big( \normm{\tilde{b}}_{L^2_x}+\e   \mathfrak{h}_{1} \normm{\tilde{\P} \tilde{f} }_{L^2_{x,v}} + \normm{\ipt \tilde{f} }_{L^2_{x,v}} \Big).
\end{split}
\end{equation}

For $\Tilde{\Xi}_{a,6}^1$, the elliptic estimate \eqref{a-6-elliptic-estimate}  directly yields
\begin{equation}
\begin{split} \label{Theta1 - a L6 estimate}
\norm{ \Tilde{\Xi}_{a,6}^{1}  }
%= & \e \left| \iint_{\Omega \times \R^{3}} \frac{\tilde{\psi}_{a,6}}{\sqrt{\tilde{\mu}}} %\pt_t \tilde{\Phi} \right|\\
%---------------
%=&  \e \Big| \int_{\Omega } \sum_{i=1}^{3}  \pt_{i} \tilde{\varphi}_{a,6} \pt_t %\int_{\R^{3}}   v_i(|v|^2 - 10)\sqrt{\tilde{\mu}} \frac{\tilde{\Phi}}{\sqrt{\tilde{\mu}}} %\Big|\\
%---------------
= &  \e \Big| \sum_{i=1}^{3}   \iint_{\Omega\times\R^3 } \pt_{i} \tilde{\varphi}_{a,6} v_i(|v|^2 - 10)  \sqrt{\tilde{\mu}} \pt_t \tilde{f}\Big|
%---------------
%=&  \e \Big| \int_{\Omega } \sum_{i=1}^{3}  \pt_{i} \tilde{\varphi}_{a,6}  \int_{\R^{3}}   %v_i(|v|^2 - 10)\sqrt{\tilde{\mu}} \frac{\pt_t\tilde{\Phi}}{\sqrt{\tilde{\mu}}} \Big|\\
%---------------------
\lesssim  \e
 \normm{\tilde{a}}_{L^6_{x}}^{5}\normm{\pt_t  \tilde{f}}_{L^2_{x,v}}.
%---------------
\end{split}
\end{equation}

For $ \tilde{\Xi}_{a,q}^{2}$ ($q\in\{2, 6\}$), the homogeneous Neumann boundary condition  ${\pt_n} \tilde{\varphi}_{a,q}\big|_{\pt\O}= 0$ implies $\mathscr{R}(\tilde{\psi}_{a,q})=\tilde{\psi}_{a,q}$.
%$\tilde{\psi}_{a,q}(t,x,v)$ is specular reflection invariant:
%\begin{equation*}%\label{reflection-invariant-condition-a}
%\begin{split}
%\mathscr{R}\big(\tilde{\psi}_{a,q}(t,x,v) \big)
% & \sum_{i=1}^{3} \pt_{i} \tilde{\varphi}_{a,q}(t,x) \big(\sqrt{10} \tilde{B}_{i}(R_x v) %- 5\tilde{\chi}_{i}(R_x v)\big)\\
%= & \sum_{i=1}^{3} \pt_{i} \tilde{\varphi}_{a,q}(t,x) \big(\sqrt{10} \tilde{B}_{i}(v) - %5\tilde{\chi}_{i}(v)\big) \\
%----
%& -\sum_{i=1}^{3} 2(n\cdot v) [\pt_{i} \tilde{\varphi}_{a,q}(t,x)\cdot n]  \big(\sqrt{10} %\tilde{B}_{i}(v) - 5\tilde{\chi}_{i}(v)\big)\\
%= &\tilde{\psi}_{a,q}(t,x,v),   \quad q\in\{2, 6\}.
%\end{split}\end{equation*}
Thus, using \eqref{Theta2 - L2 estimate},  \eqref{Theta2 - L6 estimate}, \eqref{a-2-elliptic-estimate} and \eqref{a-6-elliptic-estimate}, we obtain
\begin{align} \label{Theta2 - L2 estimate-L2}
\norm{\tilde{\Xi}^{2}_{a,2}}
\lesssim &   \a \Big( \norm{\tilde{f} }_{L^2_{\g_{+}}} +\norm{r}_{L^2_{\g_{-}}} \Big) \normm{a}_{L^2_x}, \\
%-----------
%-----------
\norm{\tilde{\Xi}_{a,6}^{2}}
  \lesssim & \a  \Big( \norm{\tilde{f}}_{L^2_{\g_{+}}}^{\frac{1}{2}} \normm{\omega^{\frac{1}{2}} \tilde{f} }_{L^\infty_{x,v}}^{\frac{1}{2}} +\norm{r}_{L^4_{\g_{-}}} \Big) \normm{\tilde{a}}_{L^6_{x}}^{5}. \label{Theta2 - L6 estimate-L6}
\end{align}

For $\Tilde{\Xi}_{a,q}^{3}$ ($q\in\{2, 6\}$),  we compute
\begin{equation}\label{Theta3 - a estimate-form}
\begin{split}
\Tilde{\Xi}_{a,q}^{3}
%=&-\iint_{\Omega \times \R^3} v\cdot \nabla_x\Big(  %\frac{\tilde{\psi}_{a,2}}{\sqrt{\tilde{\mu}}} \Big) \sqrt{\tilde{\mu}} \tilde{f} \\
%---------------
%=& -\sum_{i,j=1}^{3}\int_{\Omega} \pt_i\pt_j\tilde{\varphi}_{a,2}\int_{\R^3} %v_iv_j(|v|^2-10)\sqrt{\tilde{\mu}} \tilde{f} \dd v \dd x\\
%---------------
=& -\sum_{i,j=1}^{3}\int_{\Omega} \pt_i \pt_{j} \tilde{\varphi}_{a,q} \int_{\R^{3}}  v_iv_j(|v|^2 - 10)\sqrt{\tilde{\mu}} \big[\tilde{\P} \tilde{f} + \ipt \tilde{f} \big].
\end{split}
\end{equation}
Using the decomposition \eqref{tilde-P-tilde-f-tilde-base},
\begin{equation}\label{Theta-a3-01-estimate-Pf}
\begin{split}
  \int_{\R^{3}}   v_iv_j(|v|^2 - 10)\sqrt{\tilde{\mu}} \tilde{\P} \tilde{f} \dd v
%---------------
%---------------
= &  \sum_{k=0}^{4} \langle  \tilde{f}, \tilde{\chi}_{k} \rangle  \int_{\R^{3}}   v_iv_j(|v|^2 - 10)\sqrt{\tilde{\mu}}\tilde{\chi}_{k}(v)\dd v  + \tilde{K}_{a,q} \\
%---------------
%---------------
= &  -5\delta_{ij} a + O(|\c|+|\th|)(\tilde{a}+ |\tilde{b}|+\tilde{c})  + \tilde{K}_{a,q},
\end{split}
\end{equation}
where we used \eqref{v_iv_j(|v|^2 - 10)-orthogonal} for the above velocity integral.  Similarly to \eqref{tilde-K-a-1}, the remainder $\tilde{K}_{a,q}$ is bounded by $\e  \mathfrak{h}_{1}\normm{ \tilde{\P}  \tilde{f}}_{L^q_v}$.
%\begin{equation}\label{tilde-K-a-2}
%\begin{split}
%\norm{ \tilde{K}_{a,2} } := & \norm{ \int_{\R^{3}}   v_iv_j(|v|^2 - %10)\sqrt{\tilde{\mu}}\Big(
%\tilde{\P} \tilde{f} -\sum_{k=0}^{4} \inn{  \tilde{f}, \tilde{\chi}_{k} }  %\tilde{\chi}_{k}\Big)\dd v }
%\lesssim  \e  \mathfrak{h}_{1}\normm{ \tilde{\P}  \tilde{f}}_{L^2_v}.
%\end{split}
%\end{equation}
Substituting \eqref{Theta-a3-01-estimate-Pf} into \eqref{Theta3 - a estimate-form} gives
\begin{equation}\label{Theta3 - a estimate-uniform}
\begin{split}
\Tilde{\Xi}_{a,q}^{3}
= & \int_{\Omega} 5\Delta_x \tilde{\varphi}_{a,q} \tilde{a} + \tilde{E}_{a,q},   \quad q\in\{2, 6\},
\end{split}
\end{equation}
where
\begin{equation*}%\label{E-a-2-estimate}
\begin{split}
\tilde{E}_{a,q}
=& \sum_{i,j=1}^{3}\int_{\Omega} \pt_i \pt_{j} \tilde{\varphi}_{a,q}
 \Big[   O(|\c|+|\th|)(\tilde{a}+ |\tilde{b}|+\tilde{c})  + \tilde{K}_{a,q}
 +\int_{\R^{3}}   v_iv_j(|v|^2 - 10)\sqrt{\tilde{\mu}} \ipt \tilde{f}\Big ].
\end{split}
\end{equation*}
Applying \eqref{Theta3 - a estimate-uniform} with the elliptic equations \eqref{a-2-elliptic-equation} and \eqref{a-6-elliptic-equation} yields
\begin{align}
\Tilde{\Xi}_{a,2}^{3}
= & \int_{\Omega} 5\Delta_x \tilde{\varphi}_{a,2} \tilde{a}\dd x + \tilde{E}_{a,2}
= -5 \normm{ \tilde{a} }^2_{L^{2}_{x}}  + \tilde{E}_{a,2},\label{Theta3 - a estimate}\\
%--------------
\Tilde{\Xi}_{a,6}^{3}
=&  \int_{\Omega} 5\Delta_x \tilde{\varphi}_{a,6} \tilde{a}\dd x
 + \tilde{E}_{a,6}
= -5 \normm{ \tilde{a} }^6_{L^{6}_{x}}  + \tilde{E}_{a,6}, \label{Theta3 - a L6 ettimate}
\end{align}
where $\tilde{E}_{a,2}$ and $\tilde{E}_{a,6}$ are bounded via \eqref{a-2-elliptic-estimate}, \eqref{a-6-elliptic-estimate} and Lemma \ref{Pf - abc similar estimate}:
\begin{align}
\norm{ \tilde{E}_{a,q} }
%\leq  & \int_{\Omega} \norm{ \nabla_x^2 \tilde{\varphi}_{a,q} }
% \Big( O(|\c|+|\theta|)\normm{\tilde{\P} \tilde{f} }_{L^{2}_{v}} + \normm{ \ipt \tilde{f} %}_{L^{2}_{v}}   \Big)\\
 \lesssim &
\normm{\tilde{a}}_{L^{2}_{x}} \Big[ \e   \mathfrak{h}_{1}  \normm{\tilde{\P} \tilde{f} }_{L^{2}_{x,v}} + \normm{\ipt  \tilde{f} }_{L^{2}_{x,v}}\Big],\label{E-a-2-estimate}\\
%---------------
\norm{\tilde{E}_{a,6}}
%\lesssim &
%\int_{\Omega} \norm{ \nabla_x^2 \tilde{\varphi}_{a,6} }
% \Big[\e O\Big(\frac{|\c|}{\e}+\frac{|\theta|}{\e}\Big)\normm{\tilde{\P} \tilde{f} %}_{L^{6}_{v}} + \normm{ \ipt \tilde{f} }_{L^{6}_{v}}   \Big]\nonumber\\
 %----------
% \lesssim & \normm{\nabla_x^2 \tilde{\varphi}_{a,6}}_{L^{\frac{6}{5}}_{x}}
%\Big(\e \mathfrak{h}_{1} \normm{\tilde{\P} \tilde{f} }_{L^{6}_{x,v}}  + \normm{\ipt \ %\tilde{f} }_{L^{6}_{x,v}}\Big) \\
 %----------
% \lesssim & \normm{\tilde{a}}_{L^{6}_{x}}^{5}
%\Big( \e   \mathfrak{h}_{1}  \normm{\tilde{\P} \tilde{f} }_{L^{6}_{x,v}}  + \normm{\ipt %\tilde{f} }_{L^{6}_{x,v}}\Big) \\
 %----------
\lesssim&  \normm{\tilde{a}}_{L^{6}_{x}}^{5}
\Big[  \e^{\frac{1}{2}}  \mathfrak{h}_{1}  \normm{\e^{\frac{1}{2}} \omega^{\frac{1}{2}} \tilde{f} }_{L^{\infty}_{x,v}}  + \normm{\ipt \tilde{f} }_{L^{6}_{x,v}}\Big]. \label{E-a-6-estimate}
\end{align}

The bounds for $\tilde{\Xi}_{a,2}^{4}$ and  $\tilde{\Xi}_{a,6}^{4}$ follow directly from
\eqref{Theta4 - estimate} and elliptic estimates \eqref{a-2-elliptic-estimate} and \eqref{a-6-elliptic-estimate}.
 %\begin{align}
 %\norm{\tilde{\Xi}_{a,2}^{4}}
 %\lesssim & \Big(
 %  \e^{-1}\normm{\ipt \tilde{f} }_{L^{2}_{x,v}(\tilde{\nu})}
 % +\normm{ \tilde{\nu}^{-\frac{1}{2}} \tilde{g}}_{L^{2}_{x,v}} \Big)  % % \normm{\tilde{a}}_{L^2_{x}}, \label{Theta4 - estimate-L2} \\
 %------------------
 % \norm{\tilde{\Xi}_{a,6}^{4}}
 %\lesssim & \Big(
 %  \e^{-1}\normm{\ipt \tilde{f} }_{L^{2}_{x,v}(\tilde{\nu})}
  %+\normm{ \tilde{\nu}^{-\frac{1}{2}} \tilde{g}}_{L^{2}_{x,v}}
 % \Big)  \normm{ \tilde{a} }^5_{L^{6}_{x}}. \label{Theta4 - estimate-L6}
 %\end{align}

Integrating \eqref{test-equation-uniform-form} and combining \eqref{Theta1 - a estimate}, \eqref{Theta2 - L2 estimate-L2},  \eqref{Theta3 - a estimate} and \eqref{E-a-2-estimate} yields
\begin{equation} \label{tildea - L2 estimate step1}
\begin{split}
\int_{s}^{t} \normm{\tilde{a}}_{L^2_{x}}^2
\lesssim\; &
 \e \big[ \tilde{G}_{a} (t)- \tilde{G}_{a}(s)\big]
 + \a^2\int_{s}^{t}\Big[ \norm{ \tilde{f} }_{L^2_{\g_{+}}}^2   + \norm{r}_{L^2_{\g_{-}}}^2  + \e^2 \Big(  \mathfrak{h}_{2} +  \norm{\tilde{f}}_{L^2_{\g_{+}}}^2 \Big) \normm{ \tilde{f}}_{L^{2}_{x,v}}^2   \Big]\\
 %-------
 &+ \e  \int_{s}^{t}\normm{\pt_t \nabla \tilde{\varphi}_{a,2}}_{L^2_x} \Big( \normm{\tilde{b}}_{L^2_x} + \e  \mathfrak{h}_{1}   \normm{\tilde{\P} \tilde{f} }_{L^2_{x,v}}+ \normm{\ipt \tilde{f} }_{L^2_{x,v}} \Big)\\
 %----------------
  &+  \int_{s}^{t}  \Big ( \e^2 \mathfrak{h}_{2} \normm{
 \tilde{\P} \tilde{f} }_{L^{2}_{x,v}}^2
 + \normm{\e^{-1}  \ipt \tilde{f} }_{L^{2}_{x,v}(\tilde{\nu})}^{2}
 + \normm{\tilde{\nu}^{-\frac{1}{2}} \tilde{g} }_{L^2_{x,v}}^2
\Big).
 %---------
\end{split}
\end{equation}
Similarly, combining  \eqref{test-equation-uniform-form}, \eqref{Theta1 - a L6 estimate}, \eqref{Theta2 - L6 estimate-L6},   \eqref{Theta3 - a L6 ettimate} and \eqref{E-a-6-estimate} gives
\begin{equation} \label{tildea - L6 estimate}
\begin{split}
\normm{\tilde{a}}_{L^{6}_{x,v}} \lesssim
& \; \e\normm{ \pt_t \tilde{f} }_{L^{2}_{x,v}} + \a \norm{\tilde{f}}_{L^2_{\g_{+}}} + \a\norm{r}_{L^4_{\g_{-}}}+ \a \norm{\tilde{f}}_{L^2_{\g_{+}}}^{\frac{1}{2}} \normm{\omega^{\frac{1}{2}} \tilde{f} }_{L^\infty_{x,v}}^{\frac{1}{2}}  + \e^{\frac{1}{2}}\mathfrak{h}_{1} \normm{\e^{\frac{1}{2}} \omega^{\frac{1}{2}} \tilde{f} }_{L^{\infty}_{x,v}}\\
 &
+ \normm{\ipt  \tilde{f} }_{L^{6}_{x,v}}  + \normm{\e^{-1}\ipt  \tilde{f} }_{L^{2}_{x,v}(\tilde{\nu})} +\normm{\tilde{\nu}^{-\frac{1}{2}} \tilde{g} }_{L^2_{x,v}}.
\end{split}
\end{equation}

\noindent\textbf{Step 1.2. Estimate for $  \normm{\pt_t \nabla_x \tilde{\varphi}_{a,2}}_{L^2_x}$.}

In  \eqref{test-equation-uniform-form}, we now choose the test function
$\tilde{\psi}_{a,2} = \pt_t\tilde{\varphi}_{a,2} \sqrt{\tilde{\mu}}$ and estimate each term.

 For $\Tilde{\Xi}_{a,2}^{1}$, we write
\begin{equation}\label{a-theta1-estimate-t}
\begin{split}
   \Tilde{\Xi}_{a,2}^{1}
%=& \e \iint_{\Omega \times \R^{3}} \pt_t\tilde{\varphi}_{a,2} \sqrt{\tilde{\mu}} \pt_t %\tilde{f}
= \e \iint_{\Omega \times \R^{3}}  \pt_t\tilde{\varphi}_{a,2}\big[ \pt_t(\sqrt{\tilde{\mu}} \tilde{f})- \pt_t(\sqrt{\tilde{\mu}})  \tilde{f}\big].
\end{split}
\end{equation}
The first term in \eqref{a-theta1-estimate-t} is treated using the elliptic equation \eqref{a-2-elliptic-equation}:
\begin{equation}\label{a-theta1-estimate-t-1}
\begin{split}
 &  \e \iint_{\Omega \times \R^{3}} \pt_t\tilde{\varphi}_{a,2} \pt_t(\sqrt{\tilde{\mu}} \tilde{f})
=   \e \int_{\Omega} \pt_t \tilde{\varphi}_{a,2} \pt_t \tilde{a}
=  \e  \int_{\Omega} - \pt_t\tilde{\varphi}_{a,2} \Delta_x \pt_t\tilde{\varphi}_{a,2}
=  \e  \normm{\nabla_x \pt_t \tilde{\varphi}_{a,2}}_{L^2_x}^2.
\end{split}
\end{equation}
The second term in \eqref{a-theta1-estimate-t} is bounded similar to \eqref{Theta-a2-01-estimate-form-2}:
\begin{equation}\label{a-theta1-estimate-t-2}
\begin{split}
   \norm{ \e\iint_{\Omega \times \R^{3}} \pt_t\tilde{\varphi}_{a,2}  \pt_t(\sqrt{\tilde{\mu}})  \tilde{f} }
%\lesssim &
% \e\big( \norm{ \pt_t \theta } + \norm{ \pt_t w }+\norm{ \pt_t \rho } \big) \normm{ %\tilde{f}}_{L^{2}_{x,v}}  \normm{\pt_t \tilde{\varphi}_{a,2}}_{L^{2}_{x}}\\
 \lesssim &  \e \a \Big(  \mathfrak{h}_{1} + \norm{\tilde{f}}_{L^2_{\g_{+}}}\Big) \normm{ \tilde{f}}_{L^{2}_{x,v}}\normm{\pt_t \tilde{\varphi}_{a,2}}_{L^{2}_{x}}.
\end{split}
\end{equation}
Since $\mathscr{R}(\tilde{\psi}_{a,2})=\tilde{\psi}_{a,2}$, the estimate \eqref{Theta2 - L2 estimate} applies to $ \tilde{\Xi}_{a,2}^{2}$:
\begin{equation} \label{Theta2 - L2 estimate-t}
\begin{split}
\norm{\tilde{\Xi}^{2}_{a,2}}
\lesssim &  \a \Big( \norm{\tilde{f} }_{L^2_{\g_{+}}} +\norm{r}_{L^2_{\g_{-}}} \Big) \normm{\pt_t \tilde{\varphi}_{a,2}}_{H^1(\Omega)}.
\end{split}
\end{equation}
 By \eqref{tilde-abc-def}, direct computation implies
\begin{equation} \label{a-theta3-estimate-t}
\begin{split}
\norm {\Tilde{\Xi}_{a,2}^{3} } = &\norm{\iint_{\Omega \times \R^3}  v\cdot \nabla_x  \pt_t\tilde{\varphi}_{a,2} \sqrt{\tilde{\mu}} \tilde{f} }
=  \norm{\int_{\Omega} \nabla_x  \pt_t\tilde{\varphi}_{a,2} \cdot \tilde{b} }
\lesssim  \normm{\tilde{b}}_{L^2_x} \normm{\nabla_x \pt_t\tilde{\varphi}_{a,2}}_{L^2_x}.
\end{split}
\end{equation}
Since the contribution of $\tilde{L} \tilde{f}$ vanishes, $\Tilde{\Xi}_{a,2}^{4}$ is bounded by $\normm{\tilde{\nu}^{-\frac{1}{2}} \tilde{g} }_{L^2_{x,v}}\normm{ \pt_t \tilde{\varphi}_{a,2}}_{L^2_x}$.
%\begin{equation} \label{a-theta4-estimate-t}
%\begin{split}
%\norm {  \Tilde{\Xi}_{a,2}^{4} } = &\norm {  \iint_{\Omega \times \R^3}  \tilde{\psi}_{a,2} %\tilde{g}  }
%\lesssim  \normm{\tilde{\nu}^{-\frac{1}{2}} \tilde{g} }_{L^2_{x,v}}\normm{ \pt_t %\tilde{\varphi}_{a,2}}_{L^2_x}.
%-----------------------
%\lesssim & (\norm{\pt_t \rho}+\norm{\pt_t \th} +\norm{\pt_t \c}) \normm{\nabla_x \pt_t %\tilde{\varphi}_{a,2}}_{L^2_x}\\
%\lesssim &  \a \Big( \mathfrak{h}_{1}+  \norm{\tilde{f}}_{L^2_{\g_{+}}}\Big) %\normm{\nabla_x \pt_t \tilde{\varphi}_{a,2}}_{L^2_x}.
%\end{split}
%\end{equation}

Combining \eqref{test-equation-uniform-form} with these estimates and using Poincar\'{e}'s inequality, we have
\begin{align}
&\e \normm{\nabla_x \pt_t \tilde{\varphi}_{a,2} }_{L^2_x} \lesssim  \normm{\tilde{b}}_{L^2_x}+ \a\Big(  \norm{\tilde{f}}_{L^2_{\g_{+}}}
+  \norm{r}_{L^2_{\g_{-}}}\Big)+\e \a \Big(  \mathfrak{h}_{1} + \norm{\tilde{f}}_{L^2_{\g_{+}}}\Big) \normm{ \tilde{f}}_{L^{2}_{x,v}}+\normm{\tilde{\nu}^{-\frac{1}{2}} \tilde{g} }_{L^2_{x,v}}.  \label{varphia - pt t estimate}
\end{align}

Finally, substituting \eqref{varphia - pt t estimate} into \eqref{tildea - L2 estimate step1}, we arrive at
\begin{equation} \label{tildea - l2 estimate final}
\begin{split}
\int_{s}^{t} \normm{\tilde{a}}_{L^2_{x}}^2
\leq  &
 \e \big[ \tilde{G}_{a} (t)- \tilde{G}_{a} (s)\big]
 + \a^2\int_{s}^{t}\Big[ \norm{ \tilde{f} }_{L^2_{\g_{+}}}^2   + \norm{r}_{L^2_{\g_{-}}}^2  + \e^2 \Big(  \mathfrak{h}_{2} +  \norm{\tilde{f}}_{L^2_{\g_{+}}}^2 \Big) \normm{ \tilde{f}}_{L^{2}_{x,v}}^2   \Big]\\
 & +  \int_{s}^{t}  \Big (  \normm{\tilde{b}}_{L^2_x}^2+ \e^2 \mathfrak{h}_{2}  \normm{ \tilde{\P} \tilde{f}  }_{L^{2}_{x,v}}^2
 + \normm{\e^{-1}  \ipt \tilde{f} }_{L^{2}_{x,v}(\tilde{\nu})}^{2}+ \normm{\tilde{\nu}^{-\frac{1}{2}} \tilde{g} }_{L^2_{x,v}}^2  \Big ).
\end{split}
\end{equation}
\medskip

\noindent\textbf{Step 2. Estimate for $\tilde{b}$.}

\noindent\textbf{\bf Step 2.1.  Estimates for $\int_{s}^{t}\|\tilde{b}\|_{L^{2}_{x}}\dd \tau$ and $\|\tilde{b}\|_{L^{6}_{x}}$.}

In \eqref{test-equation-uniform-form}, we choose the test function
\begin{equation}\label{psi-b-n-definition}
\begin{split}
 \tilde{\psi}_{b,q}(t,x,v) := & \sum_{i,j=1}^{3} \pt_{j} \tilde{\varphi}_{b,q,i} \tilde{A}_{ij}(v) +  \frac{\sqrt{6}}{3}\sum_{i=1}^{3} \pt_{i}\tilde{\varphi}_{b,q,i} \tilde{\chi}_{4}(v) \\
 =&\sum_{i,j=1}^{3} \pt_{j} \tilde{\varphi}_{b,q,i} v_iv_j\sqrt{\tilde{\mu}}
 -\sum_{i=1}^{3} \pt_{i} \tilde{\varphi}_{b,q,i}\sqrt{\tilde{\mu}},  \quad
 q\in\{2, 6\},
\end{split}
\end{equation}
where the vector-valued functions $\tilde{\varphi}_{b,2}(t,x)$ and $\tilde{\varphi}_{b,6}(t,x)$ are solutions to the elliptic systems
\begin{equation}\label{b-2-elliptic-equation}
\begin{split}
- \text{div} (\nabla^{\text{s}}_{x} \tilde{\varphi}_{b,2}) = \tilde{b} \;\;  &\text{in } \Omega,\\
 \tilde{\varphi}_{b,2} \cdot n =0 \;\; &\text{on } \pt \Omega, \\
(\nabla^{\text{s}}_{x} \tilde{\varphi}_{b,2}) n = (\nabla^{\text{s}}_{x}\tilde{\varphi}_{b,2}: n \otimes n)n \;\;  &\text{on } \pt \Omega,
\end{split}
\end{equation}
and
\begin{equation}\label{b-6-elliptic-equation}
\begin{split}
- \text{div} (\nabla^{\text{s}}_{x} \tilde{\varphi}_{b,6}) = \tilde{b}^5 - \sum \frac{\int_{\Omega}A_i x \cdot \tilde{b}^5 \dd x}{\int_{\Omega}\norm{A_i x}^2\dd x}A_ix \;\; &\text{in } \Omega,\\
\tilde{\varphi}_{b,6} \cdot n =0 \;\; &\text{on } \pt \Omega,\\
(\nabla^{\text{s}}_{x} \tilde{\varphi}_{b,6}) n = (\nabla^{\text{s}}_{x}\tilde{\varphi}_{b,6}: n \otimes n)n  \;\;  &\text{on } \pt \Omega,
\end{split}
\end{equation}
respectively. Note that \eqref{b-6-elliptic-equation} has the same structure as \eqref{0-b-6-elliptic-equation-L6}, differing in the source term  and angular momentum conservation law here satisfied by $\tilde{b}$.

By the angular momentum conservation law in \eqref{abc - consercation law}, the system \eqref{b-2-elliptic-equation} satisfies  the compatible
condition \eqref{elliptic-system-compatible-condition} for all non-axisymmetric, axisymmetric and spherical domains. Moreover, for each $j=1,2,3$, a computation analogous to  \eqref{0-tilde-b5-compatible} shows that
%\begin{align}\label{tilde-b5-compatible}
%\int_{\Omega} A_{j}x \cdot \Big(\tilde{b}^5 - \sum_{i=1}^{3} \frac{\int_{\Omega}A_{i}x %\cdot \tilde{b}^5 \dd x} {\int_{\Omega} \norm{A_{i}x}^2\dd x}A_{i}x\Big) \dd x = %\int_{\Omega}A_{j}x \cdot \tilde{b}^5 \dd x - \int_{\Omega}A_{j}x \cdot \tilde{b}^5 \dd x = %0,
%\end{align}
 the system \eqref{b-6-elliptic-equation} also satisfies \eqref{elliptic-system-compatible-condition}.
Therefore, by Lemma \ref{elliptic-system-theory}, the
elliptic systems \eqref{b-2-elliptic-equation} and  \eqref{b-6-elliptic-equation} admit unique solutions satisfying
\begin{align}
\normm{\nabla^2_x \tilde{\varphi}_{b,2} }_{L^{2}_{x}} + \normm{\nabla_x \tilde{\varphi}_{b,2} }_{L^{2}_{x}} + \normm{\tilde{\varphi}_{b,2} }_{L^{2}_{x}} & \lesssim \normm{\tilde{b}}_{L^{2}_{x}}, \label{b-2-elliptic-estimate}\\
\normm{\nabla^2_x \tilde{\varphi}_{b,6}}_{L^{\frac{6}{5}}_{x}} + \normm{\nabla_x \tilde{\varphi}_{b,6}}_{L^{2}_{x}} + \normm{\tilde{\varphi}_{b,6}}_{L^{6}_{x}} &\lesssim \normm{\tilde{b}^5}_{L^{\frac{6}{5}}_{x}}=\normm{\tilde{b}}_{L^{6}_x}^5
\label{b-6-elliptic-estimate}
\end{align}
and
\begin{equation}\label{Korn-inequality-condition}
\begin{split}
 P_{\Omega}\Big(\int_{\Omega} \nabla^{\text{a}}_x \tilde{\varphi}_{b,q} \dd x\Big) = 0, \quad  q\in \{2, 6\},
\end{split}
\end{equation}
where  $P_{\O}$ denotes the orthogonal projection onto the set
$A_\O:=\big \{ A\in \mathfrak{so}(3,\mathbb{R}):  Ax\in \mathcal{R}_{\O} \big \}$.
Moreover, by \cite{Desvillettes2002} and \eqref{Korn-inequality-condition}, the following Korn-type inequality holds:
\begin{equation}\label{Korn-inequality}
\begin{split}
\|  \tilde{\varphi}_{b,q}  \|_{H^1_x}^2
\lesssim & \|\nabla^s_x \tilde{\varphi}_{b,q} \|_{L^2_x}^2 + P_{\Omega}\Big(\int_{\Omega} \nabla^{\text{a}}_x \tilde{\varphi}_{b,q} \dd x\Big)
=  \|\nabla^s_x \tilde{\varphi}_{b,q} \|_{L^2_x}^2,
% \\ \lesssim  &\normm{\nabla_x \tilde{\varphi}_{b,q} }_{L^{2}_{x}}
\;\;  q\in \{2, 6\}.
\end{split}
\end{equation}

We now estimate each term in \eqref{test-equation-uniform-form}.
For $ \Tilde{\Xi}_{b,2}^{1}$, integration by parts gives
\begin{equation}\label{Theta1 - b estimate-form}
\begin{split}
\int_{s}^{t}  \Tilde{\Xi}_{b,2}^{1}
%& \e\iint_{\Omega \times \R^{3}}
%\big[ \tilde{\psi}_{b,2} \tilde{f} (t) - \tilde{\psi}_{b,2} \tilde{f}(s)\big]
% -\e \int_{s}^{t}\iint_{\Omega \times \R^{3}} %\pt_t\Big(\frac{\tilde{\psi}_{b,2}}{\sqrt{\tilde{\mu}}}\Big)  \sqrt{\tilde{\mu}}\tilde{f} %\dd v \dd x \dd \tau  \\
= &  \e  \big[  \tilde{G}_{b} (t)- \tilde{G}_{b}  (s) \big]
-\e \int_{s}^{t}\iint_{\Omega \times \R^{3}} \Big[ \sum_{i,j=1}^{3} \pt_t \pt_{j} \tilde{\varphi}_{b,2,i}
 \tilde{A}_{ij} +\frac{\sqrt{6}}{3}\sum_{i=1}^{3} \pt_t \pt_{i}\tilde{\varphi}_{b,2,i} \tilde{\chi}_{4}\Big]  \tilde{f}\\
%----------------
&  -\e \int_{s}^{t}\iint_{\Omega \times \R^{3}} \Big[ \sum_{i,j=1}^{3}  \pt_{j} \tilde{\varphi}_{b,2,i}
\pt_t \tilde{A}_{ij}+\frac{\sqrt{6}}{3}\sum_{i=1}^{3} \pt_{i}\tilde{\varphi}_{b,2,i} \pt_t\tilde{\chi}_{4}\Big]  \tilde{f}.\\
:= & \big[ \e \tilde{G}_{b}  (t)- \e \tilde{G}_{b}  (s) \big]
- \Tilde{H}_{b,1} - \Tilde{H}_{b,2},
\end{split}
\end{equation}
Clearly, $\tilde{G}_{b}(t)$ is bounded by $\normm{\tilde{f} (t)}_{L^2_{x,v}}^2$.
For $\Tilde{H}_{b,1}$, we use
the decomposition \eqref{tilde-P-tilde-f-tilde-base}:
\begin{equation*}%\label{Theta-a3-01-estimate-Pf-compute}
\begin{split}
  \int_{\R^{3}} \tilde{A}_{ij} \tilde{\P} \tilde{f} \dd v
%---------------
%---------------
= &  \sum_{k=0}^{4} \langle  \tilde{f}, \tilde{\chi}_{k} \rangle  \int_{\R^{3}}   \tilde{A}_{ij}\tilde{\chi}_{k}\dd v  + \tilde{K}_{b,1}
%---------------
%---------------
=   O(|\c|+|\th|)(\tilde{a}+ |\tilde{b}|+\tilde{c})  + \tilde{K}_{b,1},\\
%---------------
   \int_{\R^{3}} \tilde{\chi}_{4} \tilde{\P} \tilde{f} \dd v
%---------------
%---------------
= &  \sum_{k=0}^{4} \langle  \tilde{f}, \tilde{\chi}_{k} \rangle  \int_{\R^{3}}   \tilde{\chi}_{4}\tilde{\chi}_{k}\dd v  + \tilde{K}_{b,2}
%---------------
%---------------
=  \tilde{c}+  O(|\c|+|\th|)(\tilde{a}+ |\tilde{b}|+\tilde{c})  + \tilde{K}_{b,2},
\end{split}
\end{equation*}
where we have used \eqref{Burnett-orthogonaltilde} and \eqref{tilde-chi-4-orthogonal} for the velocity integrals. The remainders $\tilde{K}_{b,1}$ and $\tilde{K}_{b,2}$ can be bounded similarly to \eqref{tilde-K-a-1}. Applying \eqref{b-2-elliptic-estimate} and Lemma \ref{Pf - abc similar estimate}, $\Tilde{H}_{b,1}$ can be estimated analogously to \eqref{Theta-a2-01-estimate-form-1}.
%\begin{equation} \label{Theta1 - b estimate-form-1}
%\begin{split}
%& \norm{ \Tilde{H}_{b,1}  }
% \leq     \e  \int_{s}^{t}\normm{\pt_t \nabla_x \tilde{\varphi}_{b,2}}_{L^2_x} \Big( %\normm{\tilde{c}}_{L^2_x} + \normm{\ipt   \tilde{f} }_{L^2_{x,v}} + \e \mathfrak{h}_{1} %\normm{\tilde{\P}  \tilde{f} }_{L^2_{x,v}}\Big)\dd \tau.
%\end{split}
%\end{equation}
The term $\Tilde{H}_{b,2}$ is bounded similarly to \eqref{Theta-a2-01-estimate-form-2}.
%\begin{equation}\label{Theta1 - b estimate-form-2}
%\begin{split}
% \norm{ \Tilde{H}_{b,2}  }
%\lesssim &
% \e\int_{s}^{t}\big( \norm{ \pt_t \theta } + \norm{ \pt_t w }+\norm{ \pt_t \rho } \big) %\normm{ \tilde{f}}_{L^{2}_{x,v}}  \normm{\nabla_x \tilde{\varphi}_{b,2}}_{L^{2}_{x}}\\
% \lesssim &  \e \a \int_{s}^{t}\Big(  \mathfrak{h}_{1} + %\norm{\tilde{f}}_{L^2_{\g_{+}}}\Big) \normm{ %\tilde{f}}_{L^{2}_{x,v}}\normm{\tilde{b}}_{L^{2}_{x}}.
%\end{split}
%\end{equation}
Combining \eqref{Theta1 - b estimate-form} with these estimates and using Lemma \ref{Pf - abc similar estimate}, we obtain
\begin{equation} \label{Theta1 - b estimate}
\begin{split}
 \norm{ \int_{s}^{t}  \Tilde{\Xi}_{b,2}^{1} }
 \leq  &\;  \e  \big[\tilde{G}_{b} (t)-  \tilde{G}_{b}  (s)\big]
 +\e \a \int_{s}^{t}\Big(  \mathfrak{h}_{1} + \norm{\tilde{f}}_{L^2_{\g_{+}}}\Big) \normm{ \tilde{f}}_{L^{2}_{x,v}}\normm{\tilde{b}}_{L^{2}_{x}}\\
 %----------------
 &+ \e  \int_{s}^{t}\normm{\pt_t \nabla_x \tilde{\varphi}_{b,2}}_{L^2_x} \Big( \normm{\tilde{c}}_{L^2_x}  + \e \mathfrak{h}_{1} \normm{\tilde{\P}  \tilde{f} }_{L^2_{x,v}} + \normm{\ipt   \tilde{f} }_{L^2_{x,v}}\Big).
\end{split}
\end{equation}

For $\Tilde{\Xi}_{b,6}^{1}$, the elliptic estimate \eqref{b-6-elliptic-estimate}  yields directly
\begin{equation}
\begin{split} \label{Theta1 - b L6 estimate}
\norm{ \Tilde{\Xi}_{b,6}^{1}  }\lesssim
 & \e
\normm{\nabla_x \tilde{\varphi}_{b,6}}_{L^{2}_{x}}
\normm{ \pt_t \tilde{f}  }_{L^2_{x,v}}
%----------
\lesssim   \e  \normm{\tilde{b}}_{L^6_{x}}^{5} \normm{\pt_t \tilde{f} }_{L^2_{x,v}}.
%---------------
\end{split}
\end{equation}

For $ \Tilde{\Xi}_{b,q}^{2}$ ($q\in\{2, 6\}$), similarly to \eqref{reflection-invariant-condition}, the boundary condition $(\nabla^{\text{s}}_{x} \tilde{\varphi}_{b,q}) n = (\nabla^{\text{s}}_{x}\tilde{\varphi}_{b,q}: n \otimes n)n$   on $\pt \Omega$ implies
$\mathscr{R}(\tilde{\psi}_{b,q})=\tilde{\psi}_{b,q}$. Therefore, the estimates \eqref{Theta2 - L2 estimate} and \eqref{Theta2 - L6 estimate} apply to $ \Tilde{\Xi}_{b,2}^{2}$ and $ \Tilde{\Xi}_{b,6}^{2}$, which combining with the elliptic estimates \eqref{a-2-elliptic-estimate} and \eqref{a-6-elliptic-estimate} yields
\begin{align} \label{Thetb2 - L2 estimate-L2}
\norm{\tilde{\Xi}^{2}_{b,2}}
\lesssim &   \a \Big( \norm{\tilde{f} }_{L^2_{\g_{+}}} +\norm{r}_{L^2_{\g_{-}}} \Big) \normm{b}_{L^2_x}, \\
%-----------
%-----------
\norm{\tilde{\Xi}_{b,6}^{2}}
  \lesssim & \a  \Big(  \norm{\tilde{f}}_{L^2_{\g_{+}}}^{\frac{1}{2}} \normm{\omega^{\frac{1}{2}} \tilde{f} }_{L^\infty_{x,v}}^{\frac{1}{2}} +\norm{r}_{L^4_{\g_{-}}} \Big) \normm{\tilde{b}}_{L^6_{x}}^{5}. \label{Thetb2 - L6 estimate-L6}
\end{align}

To compute $\Tilde{\Xi}_{b,q}^{3}$ ($q\in \{2, 6\}$),  we employ the treatment as in \eqref{0-Theta3 - b estimate-psi-b}:
\begin{equation}\label{Theta3 - b estimate-psi-b}
\begin{split}
- v\cdot \nabla_x\  \tilde{\psi}_{b,q}
%=&- \sum_{i,j,k=1}^{3}  \pt_{j} \pt_{k} \tilde{\varphi}_{b,q,i}
%v_i v_jv_k \sqrt{\tilde{\mu}}
%+ \sum_{i,l=1}^{3} \pt_{i} \pt_l  \tilde{\varphi}_{b,q,i}  v_l \sqrt{\tilde{\mu}}\\
%---------------
%---------------
=&- \sum_{i,j,k=1}^{3}  \pt_{j} \pt_{k} \tilde{\varphi}_{b,q,i}
\tilde{\P} \big (v_i v_jv_k \sqrt{\tilde{\mu}} \big)
+  \sum_{i,l=1}^{3} \pt_{i} \pt_l  \tilde{\varphi}_{b,q,i}  v_l \sqrt{\tilde{\mu}}
 \\
%-------------
& - \sum_{i,j,k=1}^{3}  \pt_{j} \pt_{k} \tilde{\varphi}_{b,q,i}
\ipt \big (v_i v_jv_k \sqrt{\tilde{\mu}} \big)\\
:=&\tilde{K}_1  + \tilde{K}_2 + \tilde{K}_3.
%---------------
\end{split}
\end{equation}
Noting the basis $\{\tilde{\chi}_k\}_{k=0}^4$ is non-orthogonal, similarly to \eqref{tilde-P-tilde-f-tilde-base},
%$$
%\tilde{\P}  \big( v_i v_jv_k \sqrt{\tilde{\mu}} \big)  = \sum_{l=0}^{4} \langle   v_i %v_jv_k \sqrt{\tilde{\mu}} , \tilde{\chi}_{l} \rangle  \tilde{\chi}_{l}
%+ \big[
%\tilde{\P}  \big( v_i v_jv_k \sqrt{\tilde{\mu}} \big)   -\sum_{l=0}^{4}  \langle    v_i %v_jv_k \sqrt{\tilde{\mu}} , \tilde{\chi}_{l}  \rangle  \tilde{\chi}_{l}\big]
%$$
 we decompose $\tilde{K}_1$ as:
\begin{equation*}\label{K1-decomp}
\begin{split}
\tilde{K}_1
%=&- \sum_{i,j,k=1}^{3}  \pt_{j} \pt_{k} \tilde{\varphi}_{b,q,i} \tilde{\P}
%\big( v_i v_jv_k \sqrt{\tilde{\mu}} \big) \\
%---------------
%---------------
=& - \sum_{i,j,k=1}^{3}  \pt_{j} \pt_{k} \tilde{\varphi}_{b,q,i} \sum_{l=1}^{3} \tilde{\chi}_{l} \int_{\R^{3}}  v_i v_j v_k \sqrt{\tilde{\mu}} \tilde{\chi}_{l}\dd v
%---------------
- \sum_{i,j,k=1}^{3}  \pt_{j} \pt_{k} \tilde{\varphi}_{b,q,i} \sum_{l=0,4} \tilde{\chi}_{l} \int_{\R^{3}}  v_i v_jv_k \sqrt{\tilde{\mu}} \tilde{\chi}_{l} \dd v  \\
 %---------
 &  - \sum_{i,j,k=1}^{3}  \pt_{j} \pt_{k} \tilde{\varphi}_{b,q,i} \Big(
\tilde{\P}  \big( v_i v_jv_k \sqrt{\tilde{\mu}} \big)   -\sum_{l=0}^{4} \langle   v_i v_jv_k \sqrt{\tilde{\mu}} , \tilde{\chi}_{l} \rangle \tilde{\chi}_{l}\Big)\\
:=&\tilde{K}_{11} +\tilde{K}_{12}+ \tilde{K}_{13}.
%---------------
\end{split}
\end{equation*}
The computation of the bulk $\tilde{K}_{11}$ is similar to \eqref{0-b-macro-pb-1-K1}, yielding
\begin{equation}\label{b-macro-pb-1-K1}
\begin{split}
 \tilde{K}_{1}=& -\sum_{l=1}^3 \tilde{\chi}_{l} \Big(3\sum_{i=l} \pt_{i} \pt_{l}\tilde{\varphi}_{b,q,i} + 2\sum_{i\neq l} \pt_{i} \pt_{l}\tilde{\varphi}_{b,q,i}
+ \sum_{j\neq l} \pt_{j} \pt_{j}\tilde{\varphi}_{b,q,l}\Big)\\
& + O(|\c|+|\th|) \sum_{l=0}^{4}\tilde{\chi}_{l}
 \sum_{i,j,k=1}^{3} \pt_j \pt_k \tilde{\varphi}_{b,q,i}+ \tilde{K}_{13},
\end{split}
\end{equation}
where the $O(|\c|+|\th|)$ term arises from the computation of $\tilde{K}_{11}$ and $\tilde{K}_{12}$ via \eqref{b2-estimate-Gauss-theta3}, analogous to \eqref{Theta-a2-01-estimate} and \eqref{Theta-a3-01-estimate-Pf}.
Substituting \eqref{b-macro-pb-1-K1} into \eqref{Theta3 - b estimate-psi-b} and  proceeding as in \eqref{0-Theta3 - b estimate-psi-b2},  we obtain
\begin{equation}\label{Theta3 - b estimate-psi-b2}
\begin{split}
&- v\cdot \nabla_x \tilde{\psi}_{b,q}=\tilde{K}_1  + \tilde{K}_2 + \tilde{K}_3
%------
%= & \sum_{l=1}^3  \tilde{\chi}_{l} \Big[\sum_{i=1}^3  %\pt_{i}\pt_{l}\tilde{\varphi}_{b,q,i} -  \Big(3\sum_{i=l} %\pt_{i}\pt_{l}\tilde{\varphi}_{b,q,i}
%+ 2\sum_{i\neq l} \pt_{i}\pt_{l} \tilde{\varphi}_{b,q,i}  + \sum_{i\neq l}  %\pt_{i}\pt_{l} \tilde{\varphi}_{b,q,l}\Big) \Big]+\tilde{K}_{R,q}\\
%------
%= & \sum_{l=1}^3  \tilde{\chi}_{l} \Big[- 2\sum_{i=l}  %\pt_{i}\pt_{l}\tilde{\varphi}_{b,q,i}-
%\sum_{i\neq l}  \pt_{i}\pt_{l}\tilde{\varphi}_{b,q,i} - \sum_{i\neq l}  %\pt_{i}\pt_{l}\tilde{\varphi}_{b,q,l} \Big]+\tilde{K}_{R,q}\\
%------
%= & \sum_{l=1}^3  \tilde{\chi}_{l} \Big[-\Big(\sum_{i=l}  %\pt_{i}\pt_{l}\tilde{\varphi}_{b,q,i}+\sum_{i\neq l} %\pt_{i}\pt_{l}\tilde{\varphi}_{b,q,i}\Big)
%-\Big( \sum_{i=l}  \pt_{i}\pt_{l}\tilde{\varphi}_{b,q,i}+\sum_{i\neq l}  %\pt_{i}\pt_{i}\tilde{\varphi}_{b,q,l} \Big)\Big]+\tilde{K}_{R,q}\\
%------
%= & \sum_{l=1}^3  \tilde{\chi}_{l} \Big[- \partial_{l}(\text{div}\tilde{\varphi}_{b,q})-  %\Delta_x \tilde{\varphi}_{b,q,l}\Big]
%+\tilde{K}_{R,q}\\
%------
=  - \sqrt{ \tilde{\mu}}v \cdot \text {div} \big( \nabla^s_x \tilde{\varphi}_{b,q} \big )+\tilde{K}_{R,q},
\end{split}
\end{equation}
where
\begin{equation*}%\label{b-macro-pb-1-KR}
\begin{split}
 \tilde{K}_{R,q}:= &  O(|\c|+|\th|) \sum_{l=1}^{4} \tilde{\chi}_{l}
 \sum_{i,j,k=1}^{3} \pt_j \pt_k \tilde{\varphi}_{b,q,i}+ \tilde{K}_{13}+\tilde{K}_{3}.
\end{split}
\end{equation*}
Inserting \eqref{Theta3 - b estimate-psi-b2} into the expression of $\Tilde{\Xi}_{b,q}^{3}$ and using the decomposition \eqref{tilde-P-tilde-f-tilde-base} gives
\begin{equation}\label{Theta3 - b estimate-form-final}
\begin{split}
\Tilde{\Xi}_{b,q}^{3}
%=&-\iint_{\Omega \times \R^3} v\cdot \nabla_x \tilde{\psi}_{b,q} \tilde{f} \\
%---------------
=&   \iint_{\Omega \times \R^3}\big[ -\sqrt {\tilde{\mu} } v \cdot \text {div} \big( \nabla^s_x \tilde{\varphi}_{b,q} \big )+\tilde{K}_{R,q}\big]\big[\tilde{\P} \tilde{f}+ \ipt  \tilde{f}\big]\\
%---------------
%=&   -\iint_{\Omega \times \R^3}\sqrt {\tilde{\mu} } v \cdot \text {div} \big( \nabla^s_x %\tilde{\varphi}_{b,q} \big )\tilde{\P} \tilde{f}
%-\iint_{\Omega \times \R^3}\sqrt {\tilde{\mu} } v \cdot \text {div} \big( \nabla^s_x %\tilde{\varphi}_{b,q} \big )\ipt  \tilde{f}\\
%&+\iint_{\Omega \times \R^3}\tilde{K}_{R,q}\tilde{f} \\
%---------------
=&   -\iint_{\Omega \times \R^3}\sqrt {\tilde{\mu} } v \cdot \text {div} \big( \nabla^s_x \tilde{\varphi}_{b,q} \big )\sum_{k=0}^{4} \langle  \tilde{f}, \tilde{\chi}_{k} \rangle \tilde{\chi}_{k}+\tilde{P}_{R,q}\\
%---------------
=& -\int_{\O}\tilde{ b}\cdot \text {div} \big( \nabla^s_x \tilde{\varphi}_{b,q} \big )  +\tilde{E}_{b,q}, \;\;\; q\in \{2, 6\},
\end{split}
\end{equation}
where in the last identity we used the almost orthogonality of $ \{\tilde{\chi}_{k}\}_{k=0}^4$, and
\begin{align*}
\tilde{P}_{R,q}:=&-\iint_{\Omega \times \R^3}\Big\{ \sqrt {\tilde{\mu} } v \cdot \text {div} \big( \nabla^s_x \tilde{\varphi}_{b,q} \big )\Big[ \Big(
\tilde{\P} \tilde{f} -\sum_{k=0}^{4} \langle   \tilde{f}, \tilde{\chi}_{k} \rangle \tilde{\chi}_{k}\Big)+\ipt  \tilde{f}\Big] -\tilde{K}_{R,q}\tilde{f}\Big\}, \nonumber \\
%-----------------
\tilde{E}_{b,q}:= & \tilde{P}_{R,q}  +\int_{\O}O(|\c|+|\th|)(\tilde{a}+ | \tilde{b}|+\tilde{c})
 \text {div} \big( \nabla^s_x \tilde{\varphi}_{b,q} \big ), \;\;\;\;  q\in \{2, 6\}.%\label{Theta3 - b estimate}
\end{align*}
Combining \eqref{Theta3 - b estimate-form-final} with the elliptic systems  \eqref{b-2-elliptic-equation} and \eqref{b-6-elliptic-equation} gives
\begin{align}
\Tilde{\Xi}_{b,2}^{3}=&-\iint_{\Omega \times \R^3} v\cdot \nabla_x \tilde{\psi}_{b,2} \tilde{f} \dd v \dd x
%---------------
= -\int_{\O}\tilde{ b}\cdot \text {div} \big( \nabla^s_x \tilde{\varphi}^{b}_{b,2} \big )   +\tilde{E}_{b,2}
%---------------
=  \normm{ \tilde{b} }^2_{L^{2}_{x}}  + \tilde{E}_{b,2},\label{Theta3 - b estimate}\\
%-------------------
\Tilde{\Xi}_{b,6}^{3}=&-\iint_{\Omega \times \R^3} v\cdot \nabla_x \tilde{\psi}_{b,6} \tilde{f} \dd v \dd x
%---------------
= -\int_{\O}\tilde{ b}\cdot \text {div} \big( \nabla^s_x \tilde{\varphi}^{b}_{b,6} \big )   +\tilde{E}_{b,6}
%---------------
=  \normm{ \tilde{b} }^6_{L^{6}_{x}}  + \tilde{E}_{b,6}. \label{Theta3 - b L6 estimate}
\end{align}
%Using the almost orthogonality of the basis $ \{\tilde{\chi}_{k}\}_{k=0}^4$, we compute
%\begin{equation}\label{Theta3 - b estimate-main}
%\begin{split}
%&   -\iint_{\Omega \times \R^3}\sqrt {\tilde{\mu} } v \cdot \text {div} \big( \nabla^s_x %\tilde{\varphi}_{b,q}\big )\sum_{k=0}^{4} \langle   \tilde{f}, \tilde{\chi}_{k} \rangle %\tilde{\chi}_{k}\\
%---------------
%=& -\int_{\O}\tilde{ b}\cdot \text {div} \big( \nabla^s_x \tilde{\varphi}_{b,q} \big )   %+\int_{\O}O(|\c|+|\th|)(\tilde{a}+ |\tilde{b}+\tilde{c})
% \text {div} \big( \nabla^s_x \tilde{\varphi}_{b,q} \big ).
%\end{split}
%\end{equation}
% Substituting \eqref{Theta3 - b estimate-main} into \eqref{Theta3 - b estimate-form}, we %arrive at the final expression of $\Tilde{\Xi}_{b,q}^{3}$:
%\begin{equation}\label{Theta3 - b estimate-form-final}
%\begin{split}
%\Tilde{\Xi}_{b,q}^{3}=&-\iint_{\Omega \times \R^3} v\cdot \nabla_x \tilde{\psi}_{b,q} %\tilde{f}
%---------------
%= -\int_{\O}\tilde{ b}\cdot \text {div} \big( \nabla^s_x \tilde{\varphi}_{b,q} \big )   %+\tilde{E}_{b,q}, \;\;\;\;  q\in \{2, 6\},
%\end{split}
%\end{equation}
%where
Note that $\tilde{E}_{b,q}$ consists of two types of terms: the first type involves  $\tilde{\P} \tilde{f}$ with small coefficient $O(|\c|+|\th|)$,
% (the differences $
% \tilde{\P} \tilde{f} -\sum_{k=0}^{4}  \langle \tilde{f}, \tilde{\chi}_{k}  \rangle % \tilde{\chi}_{k}$ and $\tilde{\P}  \big( v_i v_jv_k \sqrt{\tilde{\mu}} \big)   % -\sum_{l=0}^{4} \langle  v_i v_jv_k \sqrt{\tilde{\mu}} , \tilde{\chi}_{l} \rangle % \tilde{\chi}_{l}$ are also of order $O(|\c|+|\th|)$ by virtue of Lemma \ref{Pf - abc similar % estimate}),
 arising from the almost orthogonality of the basis $ \{\tilde{\chi}_{k}\}_{k=0}^4$; the second type includes the microscopic component $\ipt \tilde{f}$ in $\tilde{P}_{R,q}$. Both types can be estimated
%\begin{equation}\label{E-b-q-estimate}
%\begin{split}
%\norm{\tilde{E}_{b,q}}
%\leq  &
% \int_{\Omega} \norm{ \nabla_x^2 \tilde{\varphi}_{b,q} }
% \Big[ \e O\Big(\frac{|\c|}{\e}+\frac{|\theta|}{\e}\Big)\normm{\tilde{\P} \tilde{f} %}_{L^{q}_{v}} + \normm{ \ipt \tilde{f} }_{L^{2}_{v}}   \Big]\dd x, \;\;\;  q\in \{2, 6\}.
%\end{split}
%\end{equation}
 similarly to \eqref{E-a-2-estimate} and \eqref{E-a-6-estimate} by using  \eqref{b-2-elliptic-estimate} and \eqref{b-6-elliptic-estimate}:
\begin{align}
\norm{\tilde{E}_{b,2}}
 \lesssim &
\Big(\e \mathfrak{h}_{1} \normm{\tilde{\P}  \tilde{f} }_{L^{2}_{x,v}} + \normm{\ipt \tilde{f}  }_{L^{2}_{x,v}}\Big) \normm{\tilde{b}}_{L^{2}_{x}},\label{E-b-2-estimate}\\
%----------------------
\norm{\tilde{E}_{b,6}}
%  \lesssim &\Big (\e \mathfrak{h}_{1} \normm{\tilde{\P} \tilde{f}}_{L^{6}_{x,v}} + %\normm{\ipt \tilde{f} }_{L^{6}_{x,v}} \Big ) \normm{\tilde{b}}_{L^{6}_{x}}^{5}\\
%------------------------
\lesssim&  \Big (\e^{\frac{1}{2}}\mathfrak{h}_{1} \normm{\e^{\frac{1}{2}} \omega^{\frac{1}{2}}\tilde{f}}_{L^{\infty}_{x,v}} + \normm{\ipt \tilde{f} }_{L^{6}_{x,v}} \Big ) \normm{\tilde{b}}_{L^{6}_{x}}^{5}.\label{E-b-estimate-L6}
\end{align}

The estimate of $\tilde{\Xi}_{b,2}^{4}$ and  $\tilde{\Xi}_{b,6}^{4}$ follow from
\eqref{Theta4 - estimate} and the elliptic estimates \eqref{b-2-elliptic-estimate} and \eqref{b-6-elliptic-estimate}.
%\begin{align}
%\norm{\tilde{\Xi}_{b,2}^{4}}
%\lesssim & \Big(
%  \e^{-1}\normm{\ipt \tilde{f} }_{L^{2}_{x,v}(\tilde{\nu})}
% +\normm{ \tilde{\nu}^{-\frac{1}{2}} \tilde{g}}_{L^{2}_{x,v}}
%  \Big)  \normm{\tilde{b}}_{L^2_{x}}, \label{Thetb4 - estimate-L2} \\
% %------------------
% \norm{\tilde{\Xi}_{b,6}^{4}}
%\lesssim & \Big(
%  \e^{-1}\normm{\ipt \tilde{f} }_{L^{2}_{x,v}(\tilde{\nu})}
% +\normm{ \tilde{\nu}^{-\frac{1}{2}} \tilde{g}}_{L^{2}_{x,v}}
%  \Big)  \normm{ \tilde{b} }^5_{L^{6}_{x}}. \label{Thetb4 - estimate-L6}
%\end{align}

Integrating \eqref{test-equation-uniform-form} and combining \eqref{Theta1 - b estimate}, \eqref{Thetb2 - L2 estimate-L2},  \eqref{Theta3 - b estimate} and   \eqref{E-b-2-estimate} gives
\begin{equation} \label{tildeb - L2 estimate step1}
\begin{split}
\int_{s}^{t} \normm{\tilde{b}}_{L^2_{x}}^2   \lesssim &
 \e   \big [\tilde{G}_{b}  (t)-  \tilde{G}_{b}  (s) \big ]
 + \a^2\int_{s}^{t}\Big[ \norm{ \tilde{f} }_{L^2_{\g_{+}}}^2   + \norm{r}_{L^2_{\g_{-}}}^2  + \e^2 \Big(  \mathfrak{h}_{2} +  \norm{\tilde{f}}_{L^2_{\g_{+}}}^2 \Big) \normm{ \tilde{f}}_{L^{2}_{x,v}}^2   \Big]\\
  %---------.
  &+ \e  \int_{s}^{t}\normm{\pt_t \nabla_x \tilde{\varphi}_{b,2}}_{L^2_x}  \Big ( \normm{\tilde{c}}_{L^2_x}
 + \e \mathfrak{h}_{1}  \normm{\tilde{\P} \tilde{f} }_{L^2_{x,v}}
 + \normm{\ipt  \tilde{f}  }_{L^2_{x,v}}\Big )\\
 %-------
 &+  \int_{s}^{t}    \Big  (
 \e^2 \mathfrak{h}_{2}\normm{\tilde{\P} \tilde{f} }_{L^{2}_{x,v}}^2
 + \normm{\e^{-1}  \ipt \tilde{f} }_{L^{2}_{x,v}(\tilde{\nu})}^{2}
 + \normm{  \tilde{\nu}^{-\frac{1}{2}} \tilde{g}}_{L^{2}_{x,v}}^{2}  \Big ).
\end{split}
\end{equation}
Combining  \eqref{test-equation-uniform-form} with \eqref{Theta1 - b L6 estimate}, \eqref{Thetb2 - L6 estimate-L6}, \eqref{Theta3 - b L6 estimate} and  \eqref{E-b-estimate-L6} yields
\begin{equation} \label{tildeb - L6 estimate}
\begin{split}
\normm{\tilde{b}}_{L^{6}_{x,v}} \lesssim \; &  \e\normm{\pt_t \tilde{f}}_{L^{2}_{x,v}} + \a \norm{\tilde{f}}_{L^2_{\g_{+}}} +\a\norm{r}_{L^4_{\g_{-}}} + \a   \norm{\tilde{f}}_{L^2_{\g_{+}}}^{\frac{1}{2}} \normm{\omega^{\frac{1}{2}} \tilde{f} }_{L^\infty_{x,v}}^{\frac{1}{2}} +\e^{\frac{1}{2}}\mathfrak{h}_{1} \normm{\e^{\frac{1}{2}} \omega^{\frac{1}{2}} \tilde{f}}_{L^{\infty}_{x,v}}  \\
%-----------
 &   + \normm{\ipt \tilde{f} }_{L^{6}_{x,v}}+ \normm{\e^{-1}\ipt \tilde{f} }_{L^{2}_{x,v}(\tilde{\nu})} + \normm{ \tilde{\nu}^{-\frac{1}{2}} \tilde{g} }_{L^{2}_{x,v}}.
\end{split}
\end{equation}

\noindent\textbf{ Step 2.2. Estimate for $  \normm{\pt_t \nabla_x \tilde{\varphi}_{b,2}}_{L^2_x} $.}

In \eqref{test-equation-uniform-form}, we now choose the test function
$\tilde{\psi}_{b,2} = \pt_t\tilde{\varphi}_{b,2} \cdot v \sqrt{\tilde{\mu}}$ and
estimate each term.

For $\Tilde{\Xi}_{b,2}^{1}$, we decompose $\sqrt{\tilde{\mu}}  \pt_t \tilde{f}= \pt_t (\sqrt{\tilde{\mu}} \tilde{f})-\pt_t (\sqrt{\tilde{\mu}})  \tilde{f}$.
%\begin{equation} \label{b-t-theta1-estimate-t}
%\begin{split}
%\Tilde{\Xi}_{b,2}^{1} =
%& \; \e \iint_{\Omega \times \R^{3}}
% \pt_t \tilde{\varphi}_{b,2} \cdot v \pt_t (\sqrt{\tilde{\mu}} \tilde{f})
% -\e \iint_{\Omega \times \R^{3}}
% \pt_t \tilde{\varphi}_{b,2} \cdot v \pt_t (\sqrt{\tilde{\mu}})  \tilde{f}
%\end{split}
%\end{equation}
Noticing $\pt_t \tilde{\varphi}_{b,2} \in \mathscr{H}(O)$, the variational formulation of \eqref{b-2-elliptic-equation} (cf.  \eqref{elliptic-system-variational-formulation} and \eqref{elliptic-system-variational-formulation-space})  yields, for the first part
\begin{equation} \label{b-t-theta1-estimate-t-1}
\begin{split}
&\e \iint_{\Omega \times \R^{3}}
 \pt_t \tilde{\varphi}_{b,2} \cdot v \pt_t (\sqrt{\tilde{\mu}} \tilde{f})  =
  \e \int_{\Omega} \pt_t \tilde{\varphi}_{b,2}\cdot \pt_t \tilde{b}  \\
=&  \e \int_{\Omega} \big(\nabla^{\text{s}}_x \pt_t \tilde{\varphi}_{b,2}\big):\big(\nabla^{\text{s}}_x \pt_t \tilde{\varphi}_{b,2}\big)
=  \e \normm{\nabla^{\text{s}}_x \pt_t \tilde{\varphi}_{b,2}}_{L^2_{x}}^2.
\end{split}
\end{equation}
The second part is bounded analogously to \eqref{a-theta1-estimate-t-2}:
\begin{equation}\label{b-theta1-estimate-t-2}
\begin{split}
   \norm{ \e\iint_{\Omega \times \R^{3}} \pt_t\tilde{\varphi}_{b,2}\cdot v \pt_t(\sqrt{\tilde{\mu}})  \tilde{f} }
 \lesssim &  \e \a \Big(  \mathfrak{h}_{1} + \norm{\tilde{f}}_{L^2_{\g_{+}}}\Big) \normm{ \tilde{f}}_{L^{2}_{x,v}}\normm{\pt_t \tilde{\varphi}_{b,2}}_{L^{2}_{x}}.
\end{split}
\end{equation}
For $ \Tilde{\Xi}_{b,2}^{2}$, the boundary condition
  $\tilde{\varphi}_{b,2} \cdot n =0$   on $\pt \Omega$ implies
$\mathscr{R}(\tilde{\psi}_{b,2})=\tilde{\psi}_{b,2}$. Therefore, the estimate \eqref{Theta2 - L2 estimate} applies to $\tilde{\Xi}^{2}_{b,2}$:
\begin{align} \label{Thetb2 - L2 estimate-L2-t}
\norm{\tilde{\Xi}^{2}_{b,2}}
\lesssim &   \a \Big( \norm{\tilde{f} }_{L^2_{\g_{+}}} +\norm{r}_{L^2_{\g_{-}}} \Big) \normm{\pt_t \tilde{\varphi}_{b,2}}_{H^1_x}.
\end{align}
 For $\Tilde{\Xi}_{b,2}^{3}$, we
%\begin{equation*}
%\begin{split}
% \Tilde{\Xi}_{b,2}^{3}
%=& \iint_{\Omega \times \R^3} \sum_{i,j=1}^{3}
%\pt_{i}\pt_t\tilde{\varphi}_{b,2,j} v_{i}v_{j}\sqrt{\tilde{\mu}}
%\big[\tilde{\P} \tilde{f} + \ipt \tilde{f} \big].
%\end{split}
%\end{equation*}
use the decomposition \eqref{tilde-P-tilde-f-tilde-base}:
\begin{equation*}\label{Theta-b2-01-estimate}
\begin{split}
  \int_{\R^{3}}   v_iv_j\sqrt{\tilde{\mu}} \tilde{\P} \tilde{f} \dd v
%---------------
%---------------
= &  \sum_{k=0}^{4} \langle  \tilde{f}, \tilde{\chi}_{k} \rangle \int_{\R^{3}}   v_iv_j\sqrt{\tilde{\mu}}\tilde{\chi}_{k}\dd v  + \bar{K}_{b,1} \\
%---------------
%---------------
= & a \d_{ij}+ c \frac{2}{\sqrt{6}} \d_{ij} + O(|\c|+|\th|)(\tilde{a}+ |\tilde{b}|+\tilde{c})  + \bar{K}_{b,1},
\end{split}
\end{equation*}
where we used \eqref{vivj-chi-k-orthogonal} for the velocity integral. The remainder
 $\bar{K}_{b,1}$ can be bounded as in \eqref{tilde-K-a-1}. Following the argument similar to \eqref{Theta-a2-01-estimate-form-1}, we derive
\begin{equation} \label{b-t-theta3-estimate}
\begin{split}
\norm{ \Tilde{\Xi}_{b,2}^{3}}
\lesssim& \Big(\normm{\tilde{a}}_{L^2_x} + \normm{\tilde{c}}_{L^2_x} + \e \mathfrak{h}_{1} \normm{\tilde{\P} \tilde{f}}_{L^2_{x,v}}  + \normm{\ipt \tilde{f} }_{L^2_{x,v}}  \Big)\normm{\nabla_x \pt_t\tilde{\varphi}_{b,2}}_{L^2_{x}}.
\end{split}
\end{equation}
 By the property of $\tilde{L}$, $ \Tilde{\Xi}_{b,2}^{4}$ is bounded directly by $\normm{\tilde{\nu}^{-\frac{1}{2}} \tilde{g} }_{L^2_{x,v}}\normm{ \pt_t \tilde{\varphi}_{a,2}}_{L^2_x}$.
%\begin{equation} \label{b-theta4-estimate-t}
%\begin{split}
%\norm {  \Tilde{\Xi}_{b,2}^{4} } = &\norm {  \iint_{\Omega \times \R^3}  \tilde{\psi}_{b,2} %\tilde{g}  }
%\lesssim  \normm{\tilde{\nu}^{-\frac{1}{2}} \tilde{g} }_{L^2_{x,v}}\normm{ \pt_t %\tilde{\varphi}_{b,2}}_{L^2_x}.
%-----------------------
%\lesssim & (\norm{\pt_t \rho}+\norm{\pt_t \th} +\norm{\pt_t \c}) \normm{\nabla_x \pt_t %\tilde{\varphi}_{a,2}}_{L^2_x}\\
%\lesssim &  \a \Big( \mathfrak{h}_{1}+  \norm{\tilde{f}}_{L^2_{\g_{+}}}\Big) %\normm{\nabla_x \pt_t \tilde{\varphi}_{a,2}}_{L^2_x}.
%\end{split}
%\end{equation}

Combining \eqref{test-equation-uniform-form} with above estimates and using Korn's inequality \eqref{Korn-inequality}, we obtain
\begin{equation}
\label{varphib - pt t estimate}
\begin{split}
\e \normm{\nabla^{\text{s}}_x \pt_t \tilde{\varphi}_{b,2} }_{L^2_{x}} \lesssim& \normm{\tilde{a}}_{L^2_{x}} + \normm{\tilde{c}}_{L^2_{x}} + \e \mathfrak{h}_{1} \normm{\tilde{\P}\tilde{f} }_{L^2_{x,v}}  + \normm{\ipt \tilde{f} }_{L^2_{x,v}}  + \a \norm{\tilde{f} }_{L^2_{\g_{+}}} \\
&  + \a \norm{r}_{L^2_{\g_{-}}} +\e \a \Big(  \mathfrak{h}_{1} + \norm{\tilde{f}}_{L^2_{\g_{+}}}\Big) \normm{ \tilde{f}}_{L^{2}_{x,v}}+\normm{ \tilde{\nu}^{-\frac{1}{2}}  \tilde{g} }_{L^{2}_{x,v}} \\
\end{split}
\end{equation}

Finally, substituting \eqref{varphib - pt t estimate} into \eqref{tildeb - L2 estimate step1} and again using Korn's inequality \eqref{Korn-inequality}, we obtain
\begin{equation} \label{tildeb - l2 estimate final}
\begin{split}
\int_{s}^{t} \normm{\tilde{b}}_{L^2_{x}}^2 \lesssim  &
 \e \tilde{G}_{b} (t)-  \e \tilde{G}_{b}  (s)
 +  \a^2\int_{s}^{t}\Big[ \norm{ \tilde{f} }_{L^2_{\g_{+}}}^2   + \norm{r}_{L^2_{\g_{-}}}^2  + \e^2 \Big(  \mathfrak{h}_{2} +  \norm{\tilde{f}}_{L^2_{\g_{+}}}^2 \Big) \normm{ \tilde{f}}_{L^{2}_{x,v}}^2   \Big]\\
%-------
 &+  \int_{s}^{t}  \Big ( C_{\delta_b}\normm{\tilde{c}}_{L^2_x}^2+ \delta_b \normm{\tilde{a}}_{L^2_x}^2+ \e^2 \mathfrak{h}_{2}\normm{\tilde{\P} \tilde{f} }_{L^{2}_{x,v}}^2  \Big) \\
 %---------
 &+  \int_{s}^{t}   \Big (  \normm{\e^{-1}  \ipt \tilde{f} }_{L^{2}_{x,v}(\tilde{\nu})}^{2} +   \normm{ \tilde{\nu}^{-\frac{1}{2}}  \tilde{g} }_{L^{2}_{x,v}}^{2} \Big),
\end{split}
\end{equation}
where the small constant $\delta_b>0$ arises from Young's inequality.
\medskip

\noindent\textbf{Step 3. Estimate for $\tilde{c}$.}

\noindent\textbf{Step 3.1. Estimates for $\int_{s}^{t}\|\tilde{c}\|_{L^{2}_{x}}\dd \tau$ and $\|\tilde{c}\|_{L^{6}_{x}}$.}

In the weak formulation \eqref{test-equation-uniform-form}, define the test function
\begin{equation}\label{psi-c-n-definition}
\begin{split}
 \tilde{\psi}_{c,q}(t,x,v) := &  \sum_{i=1}^{3} \pt_{i} \tilde{\varphi}_{c,q} \sqrt{10} \tilde{B}_{i}(v),\quad q\in\{2, 6\},
\end{split}
\end{equation}
where $\tilde{\varphi}_{c,2}(x)$  and $\tilde{\varphi}_{c,6}(x)$ are solutions to the elliptic equations
\begin{align}
- \Delta_x \tilde{\varphi}_{c,2} = \tilde{c} \;\;
&\text{in } \Omega, \quad  \pt_n \tilde{\varphi}_{c,2}= 0\;\; \text{on } \pt \Omega, \quad \int_{\Omega}\tilde{\varphi}_{c,2} \dd x =0,\label{c-2-elliptic-equation} \\
%-------------
- \Delta_x \tilde{\varphi}_{c,6} = \tilde{c}^{5} -\frac{1}{\norm{\Omega}}\int_{\Omega}\tilde{c}^{5} \dd x\;\; &\text{in } \Omega, \quad  \pt_n \tilde{\varphi}_{c,6}= 0\;\; \text{on } \pt \Omega, \quad
\int_{\Omega} \tilde{\varphi}_{c,6} \dd x =0,\label{c-6-elliptic-equation}
\end{align}
respectively. Under the compatible conditions in \eqref{abc - consercation law}, Lemma \ref{Poisson-equation-theory} guarantees that the equations \eqref{c-2-elliptic-equation} and  \eqref{c-6-elliptic-equation} admit unique solutions satisfying
\begin{align}
\normm{\nabla^2_x \tilde{\varphi}_{c,2}}_{L^2_x} + \normm{\nabla_x \tilde{\varphi}_{c,2}}_{L^2_x} + \normm{\tilde{\varphi}_{c,2}}_{L^2_x} &\lesssim \normm{\tilde{c}}_{L^2_x}, \label{c-2-elliptic-estimate}\\
%---------------
\normm{\nabla^2_x \tilde{\varphi}_{c,6}}_{L^{\frac{6}{5}}_x} + \normm{\nabla_x \tilde{\varphi}_{c,6}}_{L^2_x} + \normm{\tilde{\varphi}_{c,6}}_{L^6_x}
&\lesssim \normm{\tilde{c}^{5}}_{L^{\frac{6}{5}}_x}=\normm{\tilde{c}}_{L^{6}_x}^5.
\label{c-6-elliptic-estimate}
\end{align}
%and Poincar\'{e}'s inequality
%\begin{equation}\label{Poincare-inequality-varphi-c}
%\begin{split}
%\normm{\tilde{\varphi}_{c,2}}_{L^2_x} \lesssim \normm{\nabla_x  %\tilde{\varphi}_{c,2}}_{L^2_x}.
%\end{split}\end{equation}

We now estimate each term in \eqref{test-equation-uniform-form}.
For $ \Tilde{\Xi}_{c,2}^{1}$, integration by parts yields
\begin{equation}\label{Theta1 - c estimate-from}
\begin{split}
\int_{s}^{t}   \Tilde{\Xi}_{c,2}^{1} =
  &  \e \big[ \tilde{G}_{c} (t)- \tilde{G}_{c} (s)\big] -\int_{s}^{t}  \iint_{\Omega \times \R^{3}}  \sum_{i=1}^{3}\big( \pt_t \pt_{i} \tilde{\varphi}_{c,2}\tilde{B}_{i}+ \pt_{i} \tilde{\varphi}_{c,2}\pt_t \tilde{B}_{i} \big)\tilde{f}\\
:=& \e \big[ \tilde{G}_{c} (t)- \tilde{G}_{c} (s)\big] - \Tilde{H}_{c,1}-\Tilde{H}_{c,2}.
\end{split}
\end{equation}
Clearly, $\tilde{G}_{c} (t)$ is bounded by $\normm{\tilde{f}(t)}_{L^2_{x,v}}^2$.
For $\Tilde{H}_{c,1}$, using the decomposition \eqref{tilde-P-tilde-f-tilde-base} gives
\begin{equation*}
\begin{split}
  \int_{\R^{3}}  \tilde{B}_{i} \tilde{\P} \tilde{f} \dd v
%---------------
%---------------
= &  \sum_{k=0}^{4} \langle  \tilde{f}, \tilde{\chi}_{k} \rangle \int_{\R^{3}}   \tilde{B}_{i}(v)\tilde{\chi}_{k}(v)\dd v  + \tilde{K}_{c,1}
%---------------
%---------------
=  O(|\c|+|\th|)(\tilde{a}+ |\tilde{b}|+\tilde{c})  + \tilde{K}_{c,1},
\end{split}
\end{equation*}
where we used \eqref{Burnett-orthogonaltilde} and  the remainder $\tilde{K}_{c,q}$ is bounded as in \eqref{tilde-K-a-1}.
Then $\Tilde{H}_{c,1}$ and  $\Tilde{H}_{c,2}$ can be estimated analogously to \eqref{Theta-a2-01-estimate-form-1} and \eqref{Theta-a2-01-estimate-form-2}.
%\begin{equation} \label{Theta1 - c estimate-from-1}
%\begin{split}
%\norm{\Tilde{H}_{c,1}}
% \leq  & \e  \int_{s}^{t}\normm{\pt_t \nabla_x \tilde{\varphi}_{c,2}}_{L^2_x} \Big(  \e % \mathfrak{h}_{1} \normm{\tilde{\P} \tilde{f} }_{L^2_{x,v}}+\normm{\ipt \tilde{f} % }_{L^2_{x,v}} \Big)\dd \tau.
%\end{split}
%\end{equation}
We conclude
\begin{equation} \label{Theta1 - c estimate}
\begin{split}
\int_{s}^{t}   \norm{ \Tilde{\Xi}_{c,2}^{1}  }
 \leq  &  \; \e  \big[\tilde{G}_{c} (t)-  \tilde{G}_{c} (s)\big] + \e \a \int_{s}^{t}\Big(  \mathfrak{h}_{1} + \norm{\tilde{f}}_{L^2_{\g_{+}}}\Big) \normm{ \tilde{f}}_{L^{2}_{x,v}}\normm{\tilde{c}}_{L^{2}_{x}}
 \\
 %-----------
 &+ \e  \int_{s}^{t}\normm{\pt_t \nabla_x \tilde{\varphi}_{c,2}}_{L^2_x} \Big(  \e \mathfrak{h}_{1} \normm{\tilde{\P} \tilde{f} }_{L^2_{x,v}}+\normm{\ipt \tilde{f} }_{L^2_{x,v}} \Big).
\end{split}
\end{equation}

For  $\Tilde{\Xi}_{c,6}^{1}$, the elliptic estimate \eqref{c-6-elliptic-estimate} yields directly
\begin{equation}
\begin{split} \label{Theta1 - c L6 estimate}
\norm{ \Tilde{\Xi}_{c,6}^{1}  }\lesssim
 & \e
\normm{\nabla_x \tilde{\varphi}_{c,6}}_{L^{2}_{x}}
\normm{ \pt_t \tilde{f}  }_{L^2_{x,v}}
%----------
\lesssim   \e  \normm{\tilde{c}}_{L^6_{x}}^{5} \normm{  \pt_t \tilde{f} }_{L^2_{x,v}}.
%---------------
\end{split}
\end{equation}

For $ \tilde{\Xi}_{c,q}^{2}$ ($q\in\{2, 6\}$), the Neumann condition  ${\pt_n} \tilde{\varphi}_{c,q}\big|_{\pt\O}= 0$ implies $\mathscr{R}(\tilde{\psi}_{c,q})= \tilde{\psi}_{c,q}$.
% & \sum_{i=1}^{3} \pt_{i} \tilde{\varphi}_{a,q}(t,x) \big(\sqrt{10} \tilde{B}_{i}(R_x v) %- 5\tilde{\chi}_{i}(R_x v)\big)\\
%= & \sum_{i=1}^{3} \pt_{i} \tilde{\varphi}_{a,q}(t,x) \big(\sqrt{10} \tilde{B}_{i}(v) - %5\tilde{\chi}_{i}(v)\big) \\
%----
%& -\sum_{i=1}^{3} 2(n\cdot v) [\pt_{i} \tilde{\varphi}_{a,q}(t,x)\cdot n]  \big(\sqrt{10} %\tilde{B}_{i}(v) - 5\tilde{\chi}_{i}(v)\big)\\
Thus, the estimates \eqref{Theta2 - L2 estimate} and \eqref{Theta2 - L6 estimate}
and the elliptic estimates \eqref{c-2-elliptic-estimate} and \eqref{c-6-elliptic-estimate} apply to $ \tilde{\Xi}_{c,2}^{2}$ and $ \tilde{\Xi}_{c,6}^{2}$:
\begin{align} \label{Thetc2 - L2 estimate-L2}
\norm{\tilde{\Xi}^{2}_{c,2}}
\lesssim &   \a \Big( \norm{\tilde{f} }_{L^2_{\g_{+}}} +\norm{r}_{L^2_{\g_{-}}} \Big) \normm{\tilde{c}}_{L^2_x}, \\
%-----------
%-----------
\norm{\tilde{\Xi}_{c,6}^{2}}
  \lesssim & \a   \Big(  \norm{\tilde{f}}_{L^2_{\g_{+}}}^{\frac{1}{2}} \normm{\omega^{\frac{1}{2}} \tilde{f} }_{L^\infty_{x,v}}^{\frac{1}{2}}  +\norm{r}_{L^4_{\g_{-}}} \Big) \normm{\tilde{c}}_{L^6_{x}}^{5}. \label{Thetc2 - L6 estimate-L6}
\end{align}

For $\Tilde{\Xi}_{c,q}^{3}$ ($q\in\{2, 6\}$),
%\begin{equation}\label{Thetc3 - c estimate-form}
%\begin{split}
%\Tilde{\Xi}_{c,q}^{3}
%=&-\iint_{\Omega \times \R^3} v\cdot \nabla_x\Big(  %\frac{\tilde{\psi}_{c,2}}{\sqrt{\tilde{\mu}}} \Big) \tilde{\Phi} \dd v \dd x \\
%---------------
%=& -\sum_{i,j=1}^{3}\int_{\Omega} \pt_i\pt_j\phi_{c,2}\int_{\R^3} v_iv_j(|v|^2-5) %\tilde{\Phi} \dd v \dd x\\
%---------------
%=& -\sum_{i,j=1}^{3}\int_{\Omega} \pt_i \pt_{j} \tilde{\varphi}_{c,2} \int_{\R^{3}}  v_iv_j(|v|^2 - %5)\sqrt{\tilde{\mu}} \Big(\tilde{\P} \frac{\tilde{\Phi}}{\sqrt{\tilde{\mu}}}+\ipt\frac{\tilde{\Phi}}{\sqrt{\tilde{\mu}}} \Big) \\
%---------------
%=& -\sum_{i,j=1}^{3}\int_{\Omega} \pt_i \pt_{j} \tilde{\varphi}_{c,q} \int_{\R^{3}}   %v_iv_j(|v|^2 - 5)\sqrt{\tilde{\mu}} \big[\tilde{\P} \tilde{f} + \ipt \tilde{f} \big],   %\quad q\in\{2, 6\}.
%\end{split}
%\end{equation}
applying the decomposition \eqref{tilde-P-tilde-f-tilde-base} yields
\begin{equation}\label{Theta-c3-01-estimate-Pf}
\begin{split}
  \int_{\R^{3}}   v_iv_j(|v|^2 -5)\sqrt{\tilde{\mu}} \tilde{\P} \tilde{f} \dd v
%---------------
%---------------
= &  \sum_{k=0}^{4} \langle  \tilde{f}, \tilde{\chi}_{k} \rangle  \int_{\R^{3}}   v_iv_j(|v|^2 - 5)\sqrt{\tilde{\mu}}\tilde{\chi}_{k}(v)\dd v  + \tilde{K}_{c,q} \\
%---------------
%---------------
= &  - c \frac{10}{\sqrt{6}}\delta_{ij} + O(|\c|+|\th|)(\tilde{a}+ |\tilde{b}|+\tilde{c})  + \tilde{K}_{c,q},
\end{split}
\end{equation}
where we used \eqref{v_iv_j(|v|^2 - 5)-orthogonal} and $\tilde{K}_{c,q}$ is bounded as in \eqref{tilde-K-a-1}.
Substituting \eqref{Theta-c3-01-estimate-Pf} into the expression of $\Tilde{\Xi}_{c,q}^{3}$ yields
\begin{equation}\label{Theta3 - c estimate-final-form}
\begin{split}
\Tilde{\Xi}_{c,q}^{3}
= & -\frac{10}{\sqrt{6}}\int_{\Omega} \Delta_x \tilde{\varphi}_{c,q} \tilde{c} + \tilde{E}_{c,q},   \quad q\in\{2, 6\},
\end{split}
\end{equation}
where
\begin{equation*}%\label{E-a-2-estimate}
\begin{split}
\tilde{E}_{c,q}
=& \sum_{i,j=1}^{3}\int_{\Omega} \pt_i \pt_{j} \tilde{\varphi}_{c,q}
 \Big[   O(|\c|+|\th|)(\tilde{a}+ |\tilde{b}|+\tilde{c})  + \tilde{K}_{c,q}-v_iv_j(|v|^2 - 5)\sqrt{\tilde{\mu}} \ipt \tilde{f}\Big ].
\end{split}
\end{equation*}
Combining \eqref{Theta3 - c estimate-final-form} with the elliptic equations \eqref{c-2-elliptic-equation} and \eqref{c-6-elliptic-equation} yields
\begin{align}
\Tilde{\Xi}_{c,2}^{3}
= & -\frac{10}{\sqrt{6}}\int_{\Omega} \Delta_x \tilde{\varphi}_{c,2} \tilde{c} + \tilde{E}_{c,2}
= \frac{10}{\sqrt{6}} \normm{ \tilde{c} }^2_{L^{2}_{x}}  + \tilde{E}_{c,2},\label{Theta3 - c estimate}\\
%-------------------
\Tilde{\Xi}_{c,6}^{3}
= & -\frac{10}{\sqrt{6}}\int_{\Omega} \Delta_x \tilde{\varphi}_{c,6} \tilde{c} + \tilde{E}_{c,6}
= \frac{10}{\sqrt{6}}\normm{ \tilde{a} }^6_{L^{6}_{x}} + \tilde{E}_{c,6}.\label{Theta3 - c estimate-L6}
\end{align}
The remainders $\tilde{E}_{c,2}$ and $\tilde{E}_{c,6}$ are estimated similarly to \eqref{E-a-2-estimate} and \eqref{E-a-6-estimate}:
\begin{align}
\norm{ \tilde{E}_{c,2} }
%\leq  & \int_{\Omega} \norm{ \nabla_x^2 \tilde{\varphi}_{a,q} }
% \Big( O(|\c|+|\theta|)\normm{\tilde{\P} \tilde{f} }_{L^{2}_{v}} + \normm{ \ipt \tilde{f} %}_{L^{2}_{v}}   \Big)\\
 \lesssim &
\normm{\tilde{c}}_{L^{2}_{x}} \Big[ \e   \mathfrak{h}_{1}  \normm{\tilde{\P} \tilde{f} }_{L^{2}_{x,v}} + \normm{\ipt  \tilde{f} }_{L^{2}_{x,v}}\Big],\label{E-c-2-estimate}\\
%---------------
\norm{ \tilde{E}_{c,6} }
%\lesssim &  \normm{\nabla_x^2 \tilde{\varphi}_{c,6}}_{L^{\frac{6}{5}}_{x}}
%\Big(\e \mathfrak{h}_{1} \normm{\tilde{\P} \tilde{f} }_{L^{6}_{x,v}}
%+ \normm{\ipt \tilde{f} }_{L^{6}_{x,v}}\Big),  \\
 %----------
\lesssim&  \normm{\tilde{c}}_{L^{6}_{x}}^{5}
\Big(  \e^{\frac{1}{2}}\mathfrak{h}_{1} \normm{\e^{\frac{1}{2}} \omega^{\frac{1}{2}} \tilde{f} }_{L^{\infty}_{x,v}}  + \normm{\ipt \tilde{f} }_{L^{6}_{x,v}}\Big),\label{E-c-6-estimate}
\end{align}
where \eqref{c-2-elliptic-estimate}, \eqref{c-6-elliptic-estimate} and Lemma \ref{Pf - abc similar estimate} have been used.

The estimates for $\tilde{\Xi}_{c,2}^{4}$ and  $\tilde{\Xi}_{c,6}^{4}$ follow directly from \eqref{Theta4 - estimate}, \eqref{c-2-elliptic-estimate} and \eqref{c-6-elliptic-estimate}.
%\begin{align}
%\norm{\tilde{\Xi}_{c,2}^{4}}
%\lesssim & \Big(
%  \e^{-1}\normm{\ipt \tilde{f} }_{L^{2}_{x,v}(\tilde{\nu})}
% +\normm{ \tilde{\nu}^{-\frac{1}{2}} \tilde{g}}_{L^{2}_{x,v}}
%  \Big)  \normm{\tilde{c}}_{L^2_{x}}, \label{Thetc4 - estimate-L2} \\
 %------------------
% \norm{\tilde{\Xi}_{c,6}^{4}}
%\lesssim & \Big(
%  \e^{-1}\normm{\ipt \tilde{f} }_{L^{2}_{x,v}(\tilde{\nu})}
% +\normm{ \tilde{\nu}^{-\frac{1}{2}} \tilde{g}}_{L^{2}_{x,v}}
% \Big)  \normm{ \tilde{c} }^5_{L^{6}_{x}}. \label{Thetc4 - estimate-L6}
%\end{align}

Integrating \eqref{test-equation-uniform-form}  and combining \eqref{Theta1 - c estimate}, \eqref{Thetc2 - L2 estimate-L2}, \eqref{Theta3 - c estimate} and \eqref{E-c-2-estimate}, we obtain
\begin{equation} \label{tildec - L2 estimate step1}
\begin{split}
\int_{s}^{t} \normm{\tilde{c}}_{L^2_{x}}^2  \lesssim & \;
 \e  \big[\tilde{G}_{c} (t)-  \tilde{G}_{c} (s)\big]
 + \a^2\int_{s}^{t}\Big[ \norm{ \tilde{f} }_{L^2_{\g_{+}}}^2   + \norm{r}_{L^2_{\g_{-}}}^2  + \e^2 \Big(  \mathfrak{h}_{2} +  \norm{\tilde{f}}_{L^2_{\g_{+}}}^2 \Big) \normm{ \tilde{f}}_{L^{2}_{x,v}}^2   \Big]\\
  &+ \e  \int_{s}^{t}\normm{\pt_t \nabla \tilde{\varphi}_{c,2}}_{L^2_x} \Big(   \e \mathfrak{h}_{1} \normm{\tilde{\P} \tilde{f} }_{L^2_{x,v}} + \normm{\ipt \tilde{f}  }_{L^2_{x,v}}\Big)\\
 %-------
 &+  \int_{s}^{t}   \Big ( \e^2 \mathfrak{h}_{2} \normm{ \tilde{\P} \tilde{f} }_{L^{2}_{x,v}}^2 + \normm{\e^{-1}  \ipt \tilde{f} }_{L^{2}_{x,v}(\tilde{\nu})}^{2} +   \normm{  \tilde{\nu}^{-\frac{1}{2}} \tilde{g} }_{L^{2}_{x,v}}^{2} \Big).
\end{split}
\end{equation}
Combining  \eqref{test-equation-uniform-form}, \eqref{Theta1 - c L6 estimate}, \eqref{Thetc2 - L6 estimate-L6}, \eqref{Theta3 - c estimate-L6} and \eqref{E-c-6-estimate},   we derive
\begin{equation} \label{tildec - L6 estimate}
\begin{split}
\normm{\tilde{c}}_{L^{6}_{x,v}} \lesssim & \; \e\normm{\pt_t\tilde{f} }_{L^{2}_{x,v}}+ \a \norm{\tilde{f}}_{L^2_{\g_{+}}}+ \a \norm{r}_{L^4_{\g_{-}}}
 + \a   \norm{\tilde{f}}_{L^2_{\g_{+}}}^{\frac{1}{2}} \normm{\omega^{\frac{1}{2}} \tilde{f} }_{L^\infty_{x,v}}^{\frac{1}{2}} + \e^{\frac{1}{2}}\mathfrak{h}_{1} \normm{\e^{\frac{1}{2}} \omega^{\frac{1}{2}} \tilde{f} }_{L^{\infty}_{x,v}} \\
%---------
 &  + \normm{\ipt \tilde{f} }_{L^{6}_{x,v}}+ \normm{\e^{-1}\ipt \tilde{f} }_{L^{2}_{x,v}(\tilde{\nu})} + \normm{ \tilde{\nu}^{-\frac{1}{2}} \tilde{g} }_{L^{2}_{x,v}}.
\end{split}
\end{equation}

\noindent\textbf{Step 3.2. Estimate for $  \normm{\pt_t \nabla_x \tilde{\varphi}_{c,2}}_{L^2_x} $.}

In \eqref{test-equation-uniform-form}, we now choose the test function
$\tilde{\psi}_{c,2} = \pt_t\tilde{\varphi}_{c,2} \tilde{\chi}_4(v)$ and estimate each term.

For $\Tilde{\Xi}_{c,2}^{1}$, we write $\tilde{\chi}_4\pt_t\tilde{f} =\pt_t(\tilde{\chi}_4\tilde{f})-\pt_t\tilde{\chi}_4\tilde{f}$.
Using definition of $\tilde{c}$ in \eqref{tilde-abc-def} and the elliptic equation  \eqref{c-2-elliptic-equation}, the first term becomes
\begin{equation} \label{c-t-theta1-estimate-t}
\begin{split}
\e \iint_{\Omega \times \R^{3}} \pt_t \tilde{\varphi}_{c,2} \pt_t\big(\tilde{\chi}_4\tilde{f}\big)
=  \e \int_{\Omega}  \pt_t \tilde{\varphi}_{c,2}   \pt_t \tilde{c}
=  -\e \int_{\Omega}\pt_t \tilde{\varphi}_{c,2}  \Delta_x \pt_t\tilde{\varphi}_{c,2}
=  \e \normm{\nabla_x \pt_t \tilde{\varphi}_{c,2}}_{L^2_{x}}^2.
\end{split}
\end{equation}
The second part is bounded analogously to \eqref{a-theta1-estimate-t-2}:
\begin{equation}\label{c-theta1-estimate-t-2}
\begin{split}
   \Big| \e \iint_{\Omega \times \R^{3}} \pt_t \tilde{\varphi}_{c,2} \pt_t\tilde{\chi}_4\tilde{f} \Big|
 \lesssim & \; \e \a \Big(  \mathfrak{h}_{1} + \norm{\tilde{f}}_{L^2_{\g_{+}}}\Big) \normm{ \tilde{f}}_{L^{2}_{x,v}}\normm{\pt_t \tilde{\varphi}_{c,2}}_{L^{2}_{x}}.
\end{split}
\end{equation}
For $ \Tilde{\Xi}_{c,2}^{2}$, since $\mathscr{R}(\tilde{\psi}_{c,2})=\tilde{\psi}_{c,2}$,
the estimate \eqref{Theta2 - L2 estimate} applies:
\begin{align} \label{Thetc2 - L2 estimate-L2-t}
\norm{\tilde{\Xi}^{2}_{c,2}}
\lesssim &   \a \Big( \norm{\tilde{f} }_{L^2_{\g_{+}}} +\norm{r}_{L^2_{\g_{-}}} \Big) \normm{\pt_t \tilde{\varphi}_{c,2}}_{H^1_x}.
\end{align}
 For $\Tilde{\Xi}_{c,2}^{3}$, we use the decomposition \eqref{tilde-P-tilde-f-tilde-base}. For each $i\in\{1, 2, 3\}$,
\begin{equation*}\label{Theta-c2-01-estimate}
\begin{split}
  \int_{\R^{3}}   v_i \tilde{\chi}_4 \tilde{\P} \tilde{f} \dd v
%---------------
%---------------
= &  \sum_{k=0}^{4} \langle \tilde{f}, \tilde{\chi}_{k} \rangle  \int_{\R^{3}}   v_i\tilde{\chi}_4\tilde{\chi}_{k}\dd v  + \bar{K}_{c,1}
%---------------
%---------------
=  b_i\frac{2}{\sqrt{6}} \d_{ik} + O(|\c|+|\th|)(\tilde{a}+ |\tilde{b}|+\tilde{c})  + \bar{K}_{c,1},
\end{split}
\end{equation*}
where we used \eqref{vi-tilde-chi4-orthogonal} and $\bar{K}_{c,1}$ is bounded as in \eqref{tilde-K-a-1}. Similarly to \eqref{b-t-theta3-estimate}, we derive
\begin{equation} \label{c-t-theta3-estimate}
\begin{split}
\norm{ \Tilde{\Xi}_{c,2}^{3}}
\lesssim& \Big(\normm{\tilde{b}}_{L^2_x} + \e \mathfrak{h}_{1} \normm{\tilde{\P} \tilde{f}}_{L^2_{x,v}}  + \normm{\ipt \tilde{f} }_{L^2_{x,v}}  \Big)\normm{\nabla_x \pt_t\tilde{\varphi}_{c,2}}_{L^2_{x}}.
\end{split}
\end{equation}
Finally, $\Tilde{\Xi}_{c,2}^{4}$ is bounded by $\normm{\tilde{\nu}^{-\frac{1}{2}} \tilde{g} }_{L^2_{x,v}}\normm{ \pt_t \tilde{\varphi}_{c,2}}_{L^2_x}$.
%\begin{equation} \label{c-theta4-estimate-t}
%\begin{split}
%\norm {  \Tilde{\Xi}_{c,2}^{4} } = &\norm {  \iint_{\Omega \times \R^3}  \tilde{\psi}_{c,2} %\tilde{g}  }
%\lesssim  \normm{\tilde{\nu}^{-\frac{1}{2}} \tilde{g} }_{L^2_{x,v}}\normm{ \pt_t %\tilde{\varphi}_{c,2}}_{L^2_x}.
%-----------------------
%\lesssim & (\norm{\pt_t \rho}+\norm{\pt_t \th} +\norm{\pt_t \c}) \normm{\nabla_x \pt_t %\tilde{\varphi}_{a,2}}_{L^2_x}\\
%\lesssim &  \a \Big( \mathfrak{h}_{1}+  \norm{\tilde{f}}_{L^2_{\g_{+}}}\Big) %\normm{\nabla_x \pt_t \tilde{\varphi}_{a,2}}_{L^2_x}.
%\end{split}
%\end{equation}

Combining \eqref{test-equation-uniform-form} with \eqref{c-t-theta1-estimate-t}--\eqref{c-t-theta3-estimate} and using Poincar\'{e}'s inequality, we conclude
\begin{equation}\label{varphic - pt t estimate}
\begin{split}
\e \normm{\nabla_x \pt_t \tilde{\varphi}_{c,2} }_{L^2_x} \lesssim
& \normm{\tilde{b}}_{L^2_x}+\e \mathfrak{h}_{1} \normm{\tilde{\P} \tilde{f} }_{L^{2}_{x,v}} + \normm{\ipt \tilde{f} }_{L^{2}_{x,v}} +\a\norm{ \tilde{f} }_{L^2_{\g_{+}}}  + \a \norm{r}_{L^2_{\g_{-}}}
\\
&+ \e \a \Big(  \mathfrak{h}_{1} + \norm{\tilde{f}}_{L^2_{\g_{+}}}\Big) \normm{ \tilde{f}}_{L^{2}_{x,v}} + \normm{\tilde{\nu}^{-\frac{1}{2}} \tilde{g} }_{L^2_{x,v}}.
\end{split}
\end{equation}

Finally, substituting \eqref{varphic - pt t estimate} into \eqref{tildec - L2 estimate step1} yields
\begin{equation} \label{tildec - l2 estimate final}
\begin{split}
\int_{s}^{t} \normm{\tilde{c}}_{L^2_{x}}^2
\leq  &\,
  \e \big[\tilde{G}_{c} (t)-  \tilde{G}_{c} (s)\big]
 +\a^2\int_{s}^{t}\Big[ \norm{ \tilde{f} }_{L^2_{\g_{+}}}^2   + \norm{r}_{L^2_{\g_{-}}}^2  + \e^2 \Big(  \mathfrak{h}_{2} +  \norm{\tilde{f}}_{L^2_{\g_{+}}}^2 \Big) \normm{ \tilde{f}}_{L^{2}_{x,v}}^2   \Big] \\
%---------------------
& +  \int_{s}^{t} \Big [  \d_c \normm{\tilde{b}}_{L^2_x}^2
  +   \e^2 \mathfrak{h}_{2}  \normm{ \tilde{\P}  \tilde{f} }_{L^{2}_{x,v}}^2 + \normm{\e^{-1}  \ipt \tilde{f} }_{L^{2}_{x,v}(\tilde{\nu})}^{2}+ \normm{ \tilde{\nu}^{-\frac{1}{2}} \tilde{g}  }_{L^{2}_{x,v}}^{2} \Big],
\end{split}
\end{equation}
where the small constant $\delta_c>0$ arises from Young's inequality.
\medskip

\noindent\textbf{Step 4.  Combination of the estimates for $\tilde{a}$,  $\tilde{b}$ and  $\tilde{c}$.}

%Choose
%\begin{align*}
%\delta_{b} = (2^8 4 C_{b} C_{c}^2)^{-1}, \quad  \delta_{c} = (4 C_{c})^{-1}.
%\end{align*}
%A direct computation of
%$$(2^{8}C_{b}C_{c}^2)^{-1} \times \eqref{tildea - l2 estimate final} + %(2^{11/2}C_{b}C_{c}^2)^{-1}  \times\eqref{tildeb - l2 estimate final} +\eqref{tildec - l2 %estimate final}
%$$

Following the same pattern as in  Step 4 of the proof of Proposition \ref{0-macro-L2-L6-estimate}, we combine \eqref{tildea - l2 estimate final}, \eqref{tildeb - l2 estimate final} and \eqref{tildec - l2 estimate final} and use Lemma \ref{Pf - abc similar estimate} to obtain
\begin{equation} \label{tildea - L2 estimate final}
\begin{split}
\int_{s}^{t} \normm{\tilde{\P} \tilde{f} }_{2}^{2} \lesssim  & \,
 \e \big[ \tilde{G}_0 (t)-  \tilde{G}_0 (s)\big]
 +\a^2\int_{s}^{t}\Big[ \norm{ \tilde{f} }_{L^2_{\g_{+}}}^2   + \norm{r}_{L^2_{\g_{-}}}^2  + \e^2 \Big(  \mathfrak{h}_{2} +  \norm{\tilde{f}}_{L^2_{\g_{+}}}^2 \Big) \normm{ \tilde{f}}_{L^{2}_{x,v}}^2   \Big] \\
 %-------
 &+  \int_{s}^{t}   \Big [
{ \e^2 \mathfrak{h}_{2}  \normm{\tilde{\P}  \tilde{f} }_{L^{2}_{x,v}}^2}
 + \normm{\e^{-1}  \ipt  \tilde{f} }_{L^{2}_{x,v}(\tilde{\nu})}^{2} +   \normm{ \tilde{\nu}^{-\frac{1}{2}}  \tilde{g} }_{L^{2}_{x,v}}^{2} \Big].
\end{split}
\end{equation}
Using the smallness of $\e$ and $\mathfrak{h}_{1}$ (see definition in \eqref{High order set - Definition}) and writing $\tilde{f}=\tilde{\P} \tilde{f}+\ipt \tilde{f}$, we absorb the terms $\int_{s}^{t}\e^2 \mathfrak{h}_{2} \normm{ \tilde{\P} \tilde{f} }_{L^{2}_{x,v}}^2$ and $\a^2\e^2 \int_{s}^{t}\mathfrak{h}_{2}\normm{ \tilde{\P}\tilde{f} }_{L^{2}_{x,v}}^2$ into the left-hand side of \eqref{tildea - L2 estimate final}. This proves \eqref{P-tilde-f-macro-L2}.

Combining the estimates \eqref{tildea - L6 estimate}, \eqref{tildeb - L6 estimate} and \eqref{tildec - L6 estimate}, we obtain \eqref{P-tilde-f-macro-L6}.
 This completes the proof of Proposition \ref{Psi - L2 and L6 estimate}.
\end{proof}
\medskip

For the derivative $\pt_t\tilde{f}$, we obtain the following consequence of Proposition \ref{Psi - L2 and L6 estimate}.

\begin{corollary} \label{Psi - L2 and L6 estimate-pt}
Under the same assumptions as in Proposition \ref{Psi - L2 and L6 estimate}, we have
\begin{equation} \label{P-pt-tilde-f-macro-L2}
\begin{split}
\int_{s}^{t} \normm{\tilde{\P}\big(\pt_t\tilde{f}\big)}_{L^2_{x,v}}^{2} \lesssim
& \, \e \big[ \tilde{G}_1(t) - \tilde{G}_1(s) \big ] + \int_{s}^{t}\Big( \normm{\e^{-1}\ipt \pt_t\tilde{f} }_{L^2_{x,v}(\tilde{\nu})}^2  +    \normm{\tilde{\nu}^{-\frac{1}{2}} \tilde{g}^t }_{L^2_{x,v}}^{2} \Big)\\
&+ \a^2\int_{s}^{t} \Big( \norm{\pt_t\tilde{f}}_{L^2_{\g_{+}}}^2 +\norm{\pt_t r +s}_{L^2_{\g_{-}}}^2 +\e^2 \norm{\pt_t\tilde{f}}_{L^2_{\g_{+}}}^2 \normm{ \pt_t\tilde{f}}_{L^{2}_{x,v}}^2\Big), \end{split}
\end{equation}
where $\big| \tilde{G}_1(t)\big|\lesssim \normm{\tilde{f}(t)}_{2}^2+\normm{\pt_t \tilde{f}(t)}_{2}^2$  and $\d>0$ is a sufficiently small constant.
\end{corollary}

\begin{proof}[\textbf{Proof}]
The equation \eqref{tildef t - Boltzmann with g} for $\pt_t\tilde{f}$ has exactly the same linear structure as the equation \eqref{tildef - Boltzmanneq ch1} for $\tilde{f}$, differing only in the source term and boundary remainder.  Moreover, $\pt_t(\sqrt{\tilde{\mu}}\tilde{f})$ also satisfies the same conservation laws of mass, angular momentum and energy as \eqref{Psi - conservation law}.
Therefore, Proposition \ref{Psi - L2 and L6 estimate} applied to \eqref{tildef t - Boltzmann with g} directly yields \eqref{P-pt-tilde-f-macro-L2}. The details are omitted for brevity.
\end{proof}
\bigskip

%%%%%%%%%%%%%%%%%%%%%%%%%%%%%%%%%%
%%%%%%%%%%%%%%%%%%%%%%%%%%%%%%%%%%
%%%%%%%%%%%%%%%%%%%%%%%%%%%%%%%%%%
\subsection{Nonlinear Estimates}\
\medskip

This subsection establishes the nonlinear estimates for the source terms $\tilde{g}$ and $\tilde{g}^{t}$, which are used in the energy estimate of Proposition \ref{tildef tildeft - Energy estimate}. The main result is the following.

\begin{proposition} \label{tildeg tildegt - L2 estimate}
Let $\tilde{g}$ and $\tilde{g}^{t}$ be defined as in \eqref{tildeg - definition} and \eqref{tildegt - definition}, respectively. Under the a priori assumption \eqref{theta-u-smallness-assumption}, the following estimates hold:
\begin{align}
 \Big |  \int_{0}^{t} \iint_{\Omega \times \R^3} \tilde{g} \tilde{f} \Big|
\lesssim & \e\int_{0}^{t} \normm{ \tilde{\nu}^{-\frac{1}{2}}  \tilde{\G}(\tilde{f},\tilde{f})}_{L^2_{x,v}}^2\dd s +   \e \normmm{\tilde{f}(t)}_{2}^{2}\Big(1+ \e \normmm{\tilde{f}(t)}_{2}\Big), \label{g-f-nonlinear-estimate}\\
%---------------
%---------------
\Big |  \int_{0}^{t} \iint_{\Omega \times \R^3} \tilde{g}^{t} \pt_t \tilde{f} \Big |  \lesssim
&\e^{\frac{1}{2}}  \normmm{\tilde{f}(t)}_{2}^{2}\Big(1+\normmm{\tilde{f}(t)}_{2}
+ \normmm{\tilde{f}(t)}_{2}^{2} + \left[\!\left[\tilde{f}_{0}\right]\!\right]_2^2  \Big) \label{gt-ft-nonlinear-estimate} \\
%--------
&+ \e\int_{0}^{t} \Big(
  \normm{  \tilde{\nu}^{-\frac{1}{2}}  \tilde{\G}(\pt_t \tilde{f},\tilde{f})}_{L^2_{x,v}}^2
+ \normm{\tilde{\nu}^{-\frac{1}{2}}  \tilde{\G}(   \tilde{f},\pt_t \tilde{f})}_{L^2_{x,v}}^2+ \normm{ \tilde{\nu}^{-\frac{1}{2}}  \tilde{\G}(\tilde{f},\tilde{f})}_{L^2_{x,v}}^2\Big)\dd s.\nonumber
\end{align}
Furthermore, for $\omega = e^{\b \norm{v}^2}$ with $0 < \b \ll \frac{1}{4}$, there hold:
\begin{align}
\int_{0}^{t} \normm{\tilde{g}\omega^{-1}}_{L^2_{x,v}}^2  \lesssim & \int_{0}^{t} \normm{ \tilde{\nu}^{-\frac{1}{2}} \tilde{\G}(\tilde{f},\tilde{f})}_{L^2_{x,v}}^2\dd s + \e^2 \Big(\normmm{\tilde{f}(t)}_{2}^2+\normmm{\tilde{f}(t)}_{2}^4\Big), \label{g-1over-w-estimate}\\
%------------
\int_{0}^{t} \normm{\tilde{g}^{t}\omega^{-1}}_{L^2_{x,v}}^2
\lesssim &  \normmm{\tilde{f}(t)}_{2}^2 \Big(\e^2 + \normmm{\tilde{f}(t)}_{2}^2 + \left[\!\left[\tilde{f}_{0}\right]\!\right]_2^2\Big)\label{gt-1over-w-estimate}\\
&+ \int_{0}^{t} \Big(
\normm{  \tilde{\nu}^{-\frac{1}{2}} \tilde{\G}(\pt_t \tilde{f},\tilde{f})}_{L^2_{x,v}}^2
+ \normm{  \tilde{\nu}^{-\frac{1}{2}}  \tilde{\G}(\tilde{f},\pt_t \tilde{f})}_{L^2_{x,v}}^2+ \normm{ \tilde{\nu}^{-\frac{1}{2}} \tilde{\G}(\tilde{f},\tilde{f})}_{L^2_{x,v}}^2\Big)\dd s.\nonumber
\end{align}
\end{proposition}
\medskip

The proof of Proposition \ref{tildeg tildegt - L2 estimate} will given at the end of this subsection, after several auxiliary lemmas.

\medskip

Recall the relation \eqref{f-relation-tildef}. We have the following $L^\infty$ estimate.

\begin{proposition}\label{w-tildef-tildeft-infty-estimate} \
Let $g, \pt_t g \in L^{\infty} (\mathbb{R}^{+} \times \Omega \times \mathbb{R}^{3})$ and $f_0, \pt_t f_0 \in L^{\infty} (\Omega \times \mathbb{R}^{3})$. Let $f$ be a solution of the linear Boltzmann equation \eqref{unst-lin-orig} on  $[0,T]$ with $0 < T \leq \infty$. For $0<\e\leq \e_0$, if  the a priori assumption  \eqref{theta-u-smallness-assumption} holds, then for all $t \in [0,T]$, we have
\begin{align}
\| \o f(t)   \|_{L^\infty_{x,v}} \lesssim  &  \| \o  f_0\|_{L^\infty_{x,v}} + \e^{-\frac{1}{2}}  \sup_{0\leq s\leq t}  \|\tilde{\P}\tilde{f}(s)\|_{L^6_{x,v}}
+  \e^{-\frac{3}{2}}  \sup_{0\leq s\leq t}  \|(\ipt \tilde{f}(s)\|_{L^2_{x,v}} \nonumber\\
& + \sup_{0 \le s \le t}\Big(\frac{\norm{\th(s)}}{\e} + \frac{\norm{\c(s)}}{\e}\Big) +
\e   \sup_{0\leq s\leq t}  \|\langle v\rangle^{-1}\o  g(s)\|_{L^\infty_{x,v}},\label{Lifnty-tilde-bd-unst}\\
%----------------------
\| \o \pt_t f(t) \|_{L^\infty_{x,v}} \lesssim
&  \| \o f_0\|_{L^\infty_{x,v}}+  \| \o  \pt_t f_0\|_{L^\infty_{x,v}} + \e^{-\frac{3}{2}}  \sup_{0\leq s\leq t}  \|\pt_t \tilde{f}(s)\|_{L^2_{x,v}} +  \sup_{0\leq s\leq t}  \| \tilde{f}(s)\|_{L^\infty_{x,v}}
\nonumber\\
&+ \sup_{0 \le s \le t}\Big(\frac{\norm{\th(s)}}{\e} + \frac{\norm{\c(s)}}{\e}\Big) + \sup_{0 \le s \le t} \Big( \frac{\norm{\pt_t \th (s)}}{\e} + \frac{\norm{\pt_t \c (s)}}{\e}\Big) \label{Lifnty-tilde-bd-unst-2}\\
& +  \e   \sup_{0\leq s\leq t}  \|\langle v\rangle^{-1}\o  \pt_t g(s)\|_{L^\infty_{x,v}},\nonumber
\end{align}
where $\omega = e^{\b \norm{v}^2}$ with $0 < \b \ll \frac{1}{4}$.
\end{proposition}

\begin{proof}[\textbf{Proof.}] \
\textbf{Step 1. \ Proof of \eqref{Lifnty-tilde-bd-unst}.}\

The argument  follows the same strategy as that of Proposition  \ref{L-infty-bd-unst}.

First, in the proof of Proposition \ref{lemma-fbar-infty-unst-0}, when performing the change of variables as in \eqref{Jspk*-k*-decom}--\eqref{Jspk*-k*-(I-P)f}, we adopt a new decomposition  of $f$:
\begin{equation*}
f = \tilde{\P}\tilde{f} + \ipt \tilde{f} + (f-\tilde{f}).
\end{equation*}
 Correspondingly, define
\begin{equation*}
\begin{split}
A_{1}\bar{f}(\bar{t},y,v) &:= \tilde{\P}\tilde{f}(t,x,v),\quad
A_{2}\bar{f}(\bar{t},y,v) := \ipt \tilde{f}(t,x,v),\quad
A_{3}\bar{f}(\bar{t},y,v) := (f-\tilde{f})(t,x,v).
\end{split}
\end{equation*}
Proceeding as before, we obtain an estimate analogous to \eqref{fbar-infty-unst}:
\begin{align}
\|\o\bar{f}(\bar{t})\|_{L^\infty_{y,v}(\O_\e\times \mathbb{R}^3)}   \lesssim
& e^{-\frac{\nu_0}{2}\bar{t}}   \|\o\bar{f}_0\|_{L^\infty_{y,v}(\O_\e\times \mathbb{R}^3)}
+  o(1)  \sup_{0\leq s\leq T_0}\|\o\bar{f}(s)\|_{L^{\infty}_{y,v}  (\O_\e\times \mathbb{R}^3)} \nonumber  \\
%-----------------------
%-----------------------
&+  \sup_{0\leq s\leq T_0}  \|A_{1}\bar{f}(s)\|_{L^6_{y,v}  (\O_\e\times \R^3)}
 +  \sup_{0\leq s\leq T_0} \|A_{2}\bar{f}(s)\|_{L^2_{y,v} (\O_\e\times \R^3)}\label{fbar-infty-unst-tilde}\\
& +  \sup_{0\leq s\leq T_0} \|\o^{-1}A_{3}\bar{f}(s)\|_{L^\infty_{y,v} (\O_\e\times \R^3)}
+   \sup_{0\leq s\leq T_0} \|\e \langle v\rangle^{-1} \o\bar{g}(s)\|_{L^\infty_{y,v} (\O_\e\times \mathbb{R}^3)}.\nonumber
\end{align}

Second, returning to the original time scale $0\leq t\leq \e^2 T_0$ via  \eqref{f-barf-relat},  we have
\begin{equation}\label{Lifnty-bd-T0-tilde}
\begin{split}
\| \o  f(t)\|_{L^\infty_{x,v}}   \lesssim
&  e^{-\frac{\nu_0 }{2\e^2}t} \|  \o f_0\|_{L^\infty_{x,v}}
   + o(1) \sup_{0\leq s\leq \e^2 T_0}  \|   \o f(s)  \|_{L^\infty_{x,v}}\\
&
   + \e^{-\frac{3}{2}} \sup_{0\leq s\leq \e^2 T_0} \|  \ipt \tilde{f}(s)  \|_{L^2_{x,v}}
+  \e^{-\frac{1}{2}} \sup_{0\leq s\leq \e^2 T_0} \| \tilde{\P} \tilde{f}(s)  \|_{L^6_{x,v} }\\
&
+ \sup_{0\leq t\leq \e^2 T_0}  \| \o ^{-1}(f- \tilde{f})(t)\|_{L^\infty_{x,v}(\Omega\times \mathbb{R}^3)} +\e  \sup_{0\leq s\leq \e^2 T_0}  \| \langle v\rangle^{-1}  \o  g(s) \|_{L^\infty_{x,v}}.
\end{split}
\end{equation}
Define
\begin{equation*}
\begin{split}
D(s) := & \  o(1)\|  \o  f(s)  \|_{L^\infty_{x,v}}
+ \e^{-\frac{1}{2}} \| \tilde{\P}\tilde{f}(s) \|_{L^6_{x,v}}  +    \e^{-\frac{3}{2}} \|  \ipt \tilde{f}(s) \|_{L^2_{x,v}}\\
&  + \Big(\frac{\norm{\th(s)}}{\e} + \frac{\norm{\c(s)}}{\e}\Big) +  \e \| \langle v\rangle^{-1} \o  g(s)\|_{L^\infty_{x,v}}.
\end{split}
\end{equation*}
 Applying the previous inequality iteratively and using Lemma \ref{tilde-f-bounded-by-f} yields \eqref{Lifnty-tilde-bd-unst}.

\medskip
\noindent\textbf{Step 2.  Proof of \eqref{Lifnty-tilde-bd-unst-2}.}

The proof is similar. We start with
$\pt_t f = \pt_t \tilde{f} + (\pt_t f - \pt_t \tilde{f})$,
and set
\begin{equation*}
\begin{split}
A_{2}\bar{f}(\bar{t},y,v) &:=  \pt_t \tilde{f}(t,x,v),\quad
A_{3}\bar{f}(\bar{t},y,v) := \pt_t (f-\tilde{f})(t,x,v).
\end{split}
\end{equation*}
Following the same pattern as in Step 1, we derive
\begin{equation*}
\begin{split}
\| \o  \pt_t f(t)\|_{L^\infty_{x,v}}   \lesssim \
& e^{-\frac{\nu_0 }{2\e^2}t} \| \o  \pt_t f_0\|_{L^\infty_{x,v}}
+ o(1) \sup_{0\leq s\leq \e^2 T_0}  \|  \o  \pt_t f(s)  \|_{L^\infty_{x,v}}
+  \e^{-\frac{3}{2}} \sup_{0\leq s\leq \e^2 T_0} \| \pt_t \tilde{f}(s)  \|_{L^2_{x,v}}
\\
&+ \sup_{0\leq t\leq \e^2 T_0}  \|\o ^{-1}\pt_t (f- \tilde{f})(t)\|_{L^\infty_{x,v}(\Omega\times \mathbb{R}^3)}
+  \e  \sup_{0\leq s\leq \e^2 T_0}  \| \langle v\rangle^{-1} \o \pt_t g(s) \|_{L^\infty_{x,v} }.
\end{split}
\end{equation*}
Combining this with Lemma \ref{tilde-f-bounded-by-f} yields \eqref{Lifnty-tilde-bd-unst-2}.
\end{proof}
\medskip

The following lemma controls derivatives of auxiliary functions with algebraic growth in $v$.

\begin{lemma} \label{v polynomial - L2 estimate}
Let $X \in \{\th, \c, \rho\}$, $g \in L^2(\Omega \times \R^3)$, and $p\ge 0$ be an integer. Then for $\omega_1 = e^{\b_1 \norm{v}^2}$ with $0 < \b_1 \ll \frac{1}{4}$, there holds
\begin{align*}
& \normm{\pt_t X \norm{v}^p {g}}_{L^2_{x,v}}
\lesssim   \e  \Big(\frac{\a}{\e}\mathfrak{h}_{1} + \frac{\a}{\e} \norm{\tilde{f}}_{L^2_{\g_{+}}} \Big)  \Big( \frac{1}{\e} \normm{\ipt {g}}_{L^2_{x,v}} +\e^{\frac{3}{2}} \normm{ \omega_1 {g} }_{L^\infty_{x,v}} + \normm{\tilde{\P} {g}}_{L^2_{x,v}} \Big).
\end{align*}
\end{lemma}

\begin{proof}[\textbf{Proof.}]
From \eqref{rho th u - smallness}, we have
$$
 \norm{\pt_t X} \lesssim \e  \Big(\frac{\a}{\e}\mathfrak{h}_{1} + \frac{\a}{\e} \norm{\tilde{f}}_{L^2_{\g_{+}}} \Big).
$$
Decompose $\norm{v}^p g$  as
\begin{align*}
\norm{v}^p g =  \norm{v}^p \tilde{\P}  g + \1_{\norm{v}^p\le \e^{-1}}\norm{v}^p \ipt g + \1_{\norm{v}^p>  \e^{-1}}\norm{v}^p \ipt g.
\end{align*}
The first two terms satisfy
\begin{align*}
&\normm{\norm{v}^p \tilde{\P} g}_{L^2_{x,v}} \lesssim \normm{\tilde{\P} g}_{L^2_{x,v}},\quad
\normm{\1_{\norm{v}^p\le \e^{-1}}\norm{v}^p \ipt g}_{L^2_{x,v}} \lesssim \normm{\e^{-1}\ipt g}_{L^2_{x,v}}.
\end{align*}
For the last term, note that $\norm{v}^{4p} \lesssim \o^{\frac{1}{4}}_1$ for any $p$. Hence,
\begin{align*}
\normm{\1_{\norm{v}^p >  \e^{-1}}\norm{v}^p \ipt  g}_{L^2_{x,v}} \lesssim& \e^{\frac{3}{2}} \normm{\omega_1 \ipt g}_{L^\infty_{x,v}}  \normm{\o^{-\frac{1}{4}}_1}_{L^2_{x,v}}
\lesssim \e^{\frac{3}{2}}  \normm{\omega_1 g}_{L^\infty_{x,v}} .
\end{align*}
Combining these estimates completes the proof.
\end{proof}
\medskip

The next two results provide estimates for the nonlinear collision operator.

\begin{lemma} \label{tildeGamma - L2 and Linfty estimate}
Recall the definition of $\tilde{\G}$ in \eqref{tildeg - definition}. For $\omega_1 = e^{\b_1 \norm{v}^2}$ with $0 < \b_1 \ll \frac{1}{4}$, the following bounds hold:
\begin{align}
&\normm{ \tilde{\nu}^{-\frac{1}{2}} \tilde{\G}({f}, {g})}_{L^2_{x,v}} \lesssim  \normm{\omega_1 g}_{L^\infty_{x,v}} \normm{\tilde{\nu}^{\frac{1}{2}} f}_{L^2_{x,v}},\label{tilde-G-L2-estimate-1}\\
%-----------
&\normm{ \tilde{\nu}^{-\frac{1}{2}}  \tilde{\G}({f}, {g})}_{L^2_{x,v}} \lesssim  \normm{\omega_1 f}_{L^\infty_{x,v}} \normm{\tilde{\nu}^{\frac{1}{2}} g}_{L^2_{x,v}},\label{tilde-G-L2-estimate-2}\\
%-----------
&\normm{\omega_1 \tilde{\G}({f}, {g})}_{L^\infty_{x,v}} \lesssim  \normm{\omega_1 f}_{L^\infty_{x,v}} \normm{\omega_1 g}_{L^\infty_{x,v}}, \label{weighted-tilde-G-Linfty-estimate}\\
%----------------
&\normm{ \tilde{\nu}^{-\frac{1}{2}}  \tilde{\G}(\tilde{\P}f,\tilde{\P}g)}_{L^2_{x,v}} \lesssim \normm{\tilde{\P}f \tilde{\P}g}_{L^2_{x,v}}. \label{tilde-G-Pf}
\end{align}
\end{lemma}

\begin{proof}[\textbf{Proof.}] \
The estimates follow by the same arguments as in the proof of Lemma  \ref{Gamma-L2-Linfty-estimate-1}, using the properties of the collision frequency \eqref{collsion-frequency-def}.  We omit the details for brevity.
\end{proof}
\medskip

\begin{corollary} \label{tildeG - L2 Linfty L3 L6 estimate}
Let $f,g \in L^2([0,T]\times \Omega \times \R^3)$, and let $S_{j}f, S_{j}g \ge 0$ $(j=1,2)$ be defined as in  Proposition \ref{f - L2L3 estimate}. Assume that for $t\in [0,T]$,
\begin{equation*}
\begin{split}
\norm{ \tilde{a}(h)}+\sum_{i=1}^{3}\norm{ \tilde{b}_{i}(h)} +\norm{ \tilde{c}(h)}  \le S_{1}h(t,x) + S_{2}h(t,x)\quad \text{ for } h\in \{f,g\},
\end{split}
\end{equation*}
 where $\tilde{a}(h),  \tilde{b}_{i}(h)$ and $\tilde{c}(h)$ are the coefficients of $ \P h$. Then
\begin{equation}\label{Gamma-f-g-estimate}
\begin{split}
&\normm{  \tilde{\nu}^{-\frac{1}{2}} \tilde{\G}(f,g)}_{L^2_{t,x,v}}+
\normm{  \tilde{\nu}^{-\frac{1}{2}}  \tilde{\G}(g,f)}_{L^2_{t,x,v}}\\
%----------------
\lesssim  &\e^{\frac{1}{2}}\Big [  {\e}^{-1} \normm{\ipt f}_{L^2_{t,x,v}(\tilde{\nu})} +  {\e}^{-1}  \normm{S_{2} f}_{L^2_{t,x}} \Big] \Big[\e^{\frac{1}{2}} \normm{\o_1 g}_{L^\infty_{t,x,v}}\Big] \\
%----------------
 & + \normm{S_{1}f}_{L^2_{t}L^3_{x}} \Big[ \e^{\frac{1}{2}}
 \normm{\o_1 g}_{L^\infty_{t,x,v}}\Big]^{\frac{2}{3}} \Big[ {\e}^{-1} \normm{\ipt  g}_{L^\infty_{t}L^2_{x,v}(\tilde{\nu})}\Big]^{\frac{1}{3}}
 %----------------
+\normm{S_{1}f}_{L^2_{t}L^3_{x}}\normm{\tilde{\P} g}_{L^\infty_{t} L^6_{x,v}},
\end{split}
\end{equation}
where $\omega_1 = e^{\b_1 \norm{v}^2}$ with $0 < \b_1 \ll \frac{1}{4}$.
\end{corollary}

\begin{proof}[\textbf{Proof.}] \
Write
$|f| = |\tilde{\P} f| + |\ipt f|$ and $|g| = |\tilde{\P} g| + |\ipt g|$. The proof then proceeds exactly as that of Corollary \ref{G - L2 Linfty L3 L6 estimate}.
\end{proof}
\medskip

\begin{corollary} \label{tilde-f-L6-estimate-final} \
%Let $\e \in (0, \e_0]$, where $\e_0\in (0,1)$ is the constant from Proposition \ref{L-infty-bd-unst}.
Let $\tilde{f}$ be the solution of \eqref{tildef - Boltzmanneq ch1} on $[0,T]$ with $0<T\leq \infty$. Under the a priori assumption \eqref{theta-u-smallness-assumption}, the following estimates hold for all $t\in [0,T]$:
\begin{equation} \label{pf6 - bound3-final}
\begin{split}
\normm{\tilde{\P} \tilde{f}}_{L^6_{x,v}}^2
\lesssim &\left[\!\left[\tilde{f}_{0}\right]\!\right]_2^2 + \mathscr{E}_{2}[\tilde{f}](t) + \mathscr{D}_{2}[\tilde{f}](t)
 + \d\e\normm{\omega f}_{L^{\infty}_{x,v}}^2  + \left[\!\left[\tilde{f}_{0}\right]\!\right]_2^4+ \mathscr{E}_{2}^2[\tilde{f}](t)\\
 &+\mathscr{E}_{2}^3[\tilde{f}](t) + \mathscr{D}_{2}^2[\tilde{f}](t) +\d\e^2\normm{\omega f}_{L^{\infty}_{x,v}}^4,
\end{split}
\end{equation}
where $\d>0$ is a sufficiently small constant.
\end{corollary}

\begin{proof}[\textbf{Proof.}] \
We start from the estimate \eqref{P-tilde-f-macro-L6}. Both $ \e\normm{\pt_t  \tilde{f} }_{L^{2}_{x,v}}$ and
$\normm{\P \tilde{f} }_{L^{2}_{x,v}}$ are bounded by $\mathscr{E}_{2}[\tilde{f}](t) $. For the boundary term in \eqref{P-tilde-f-macro-L6},  we argue similarly to \eqref{L-t-infty-boundary-term} to obtain
\begin{equation}\label{sup-f-gamma-L2}
\begin{split}
&\frac{\a}{\e} \norm{ \tilde{f} }_{L^{2}_{\g_{+}}}^2
%= \frac{\a^4}{\e}\int_{\g_{+}} \tilde{f}^2(t)  \dd \g \\
%=& \frac{\a^4}{\e}\int_{\g_{+}} \tilde{f}^2(0)  \dd \g
%+ \frac{\a^4}{\e}\int_{0}^{t} \int_{\g_{+}}   \frac{\dd [\tilde{f}^2(t)]^2(s)}{\dd s}  \dd %\g \dd s  \\
=  \frac{\a}{\e}\int_{\g_{+}} \tilde{f}^2_0  \dd \g
+ 2\frac{\a}{\e}\int_{0}^{t} \int_{\g_{+}}   \tilde{f}(s)\pt_{t} \tilde{f}(s) \dd \g \dd s   \lesssim  \left[\!\left[\tilde{f}_{0}\right]\!\right]_2^2 +\mathscr{D}_{2}[\tilde{f}](t).
\end{split}
\end{equation}
In view of \eqref{tilde-mu-equals-mu-t0},  \eqref{tilde-f0-equals-f0-t0} and the orthogonal decomposition \eqref{P-gamma-orthogonal}, the term $\norm{ f_0 }_{L^{2}_{\g_{+}}}^2$ can be controlled via trace lemma similar to \eqref{f - Bundary ukai lemma}  and \eqref{f0-non-grazing}:
\begin{equation} \label{f0-dissipation-part}
\begin{split}
\norm{ f_0 }_{L^{2}_{\g_{+}}}^2= &\norm{(1-\mathscr{P}_{\g}) f_0 }_{L^2_{\g_{+}}}^2 +\norm{\mathscr{P}_{\g} f_0 }_{L^2_{\g_{+}}}^2
%= &\Big\{  \iint_{\g_{+}^\d} + \iint_{\g_{+}\backslash\g_{+}^\d}\Big \}\norm{\mathscr{P}_{\g} f }^2 %\dd \g\\
\lesssim  (1+\d)\norm{(1-\mathscr{P}_{\g}) f_0 }_{L^2_{\g_{+}}}^2 +|f_{0}{\bf 1}_{\g^{\d}_{\pm}}|_{L^{2}_\gamma}^2\\
\lesssim &  (1+\d)\norm{(1-\mathscr{P}_{\g}) f_0 }_{L^2_{\g_{+}}}^2 + \|f_0 \|_{L^{2}_{x,v}}^2  + \|  v \cdot \nabla_x  f_0\|_{L^{2}_{x,v}}^2.
\end{split}
\end{equation}
Consequently, only the contribution  $\norm{(1-\mathscr{P}_{\g}) f_0 }_{L^2_{\g_{+}}}^2$ is required in the definition \eqref{initial-data-total-norm-tilde defn}.

The term $\a\norm{r}_{L^4_{\g_{-}}}$ can be bounded by $\a \mathscr{E}_{2}^{\frac{1}{2}}[\tilde{f}](t)$ via \eqref{r-estimate}.
%\begin{equation*}
%\begin{split}
%\a\norm{r}_{L^4_{\g_{-}}}
%\lesssim \a \norm{\mathfrak{h}_{1}} \leq .
%\end{split}
%\end{equation*}
Moreover, as in \eqref{L-t-infty-L2-ipt-f}, we have
\begin{equation}\label{L-t-infty-L2-ipt-tilde-f}
\begin{split}
&\e^{-2} \normm{ \ipt \tilde{f}(t)}_{L^2_{x,v}(\tilde{\nu})}^2
%=& \e^{-2}\iint_{\O\times\R^3} [\ipt\tilde{f}]^2(t)\tilde{\nu}  \dd v
%\dd x \\
%=& \e^{-2} \iint_{\O\times\R^3} [\ipt\tilde{f}]^2(0)\tilde{\nu}\dd v
%\dd x
%+ \e^{-2} \int_{0}^{t} \iint_{\O\times\R^3}  \frac{\dd [\ipt\tilde{f}]^2(s)}{\dd %s}\tilde{\nu} \dd v
%\dd x \dd s  \\
%=&  \e^{-2} \iint_{\O\times\R^3} [\ipt\tilde{f}]^2(0)\tilde{\nu}\dd v
%\dd x   + \e^{-2}  \int_{0}^{t} \iint_{\O\times\R^3} \ipt\tilde{f} \pt_t %\ipt\tilde{f}\tilde{\nu}\dd v
%\dd x \dd s   \\
\lesssim  \left[\!\left[\tilde{f}_{0}\right]\!\right]_2^2 +\mathscr{D}_{2}[\tilde{f}](t).
\end{split}
\end{equation}
To estimate $\normm{\ipt  \tilde{f} }_{L^{6}_{x,v}} $, we apply interpolation, \eqref{L-t-infty-L2-ipt-tilde-f} and Lemma \ref{tilde-f-bounded-by-f}:
\begin{equation}
\begin{split}
 \normm{\ipt \tilde{f}}_{L^{6}_{x,v}}\leq &
 \Big( \e^{\frac{1}{2}} \normm{ \o^{\frac{1}{2}} \tilde{f}}_{L^{\infty}_{x,v}} \Big)^{\frac{2}{3}}
  \Big( \e^{-1} \normm{ \ipt \tilde{f}}_{L^2_{x,v}} \Big)^{\frac{1}{3}} \\
  \leq
 &\d \e^{\frac{1}{2}} \normm{ \o^{\frac{1}{2}} \tilde{f}}_{L^{\infty}_{x,v}}
  +C_\d \e^{-1} \normm{ \ipt \tilde{f}}_{L^2_{x,v}}\\
  \leq
 &\d \e^{\frac{1}{2}} \normm{ \o f}_{L^{\infty}_{x,v}} + \d \e^{\frac{1}{2}}\mathscr{E}_{2}^{\frac{1}{2}}[\tilde{f}](t)
  + \left[\!\left[\tilde{f}_{0}\right]\!\right]_2 +\mathscr{D}_{2}^{\frac{1}{2}}[\tilde{f}](t),
\end{split}
\end{equation}
where $\d>0$ is a sufficiently small constant. Using the smallness of $\e$ and $\mathfrak{h}_{1}$ (see \eqref{High order set - Definition}),  we absorb the term $\e^{\frac{1}{2}}\mathfrak{h}_{1} \normm{\e^{\frac{1}{2}} \omega^{\frac{1}{2}} \tilde{f} }_{L^{\infty}_{x,v}}$ from \eqref{P-tilde-f-macro-L6} into $\d \e^{\frac{1}{2}} \normm{ \omega f }_{L^\infty_{x,v}}$.

For $\normm{\tilde{\nu}^{-\frac{1}{2}} \tilde{g} }_{L^2_{x,v}}$, recall the definition of $\tilde{g}$ in \eqref{tildeg - definition}. By \eqref{pt-sqrt-mu}, Lemma \ref{v polynomial - L2 estimate}, Lemma \ref{tilde-f-bounded-by-f} and the assumption \eqref{theta-u-smallness-assumption},
\begin{align}
\e \normm{\tilde{\nu}^{-\frac{1}{2}}\frac{ \pt_t \sqrt{\tilde{\mu}} } {\sqrt{\mu}} \tilde{f}
}_{L^2_{x,v}}
%-----------
\lesssim & \e  \Big(\frac{\a}{\e}\mathfrak{h}_{1} + \frac{\a}{\e} \norm{\tilde{f}}_{L^2_{\g_{+}}} \Big)  \Big( \frac{1}{\e} \normm{\ipt \tilde{f}}_{L^2_{x,v}} +\e^{\frac{3}{2}} \normm{ \o^{\frac{1}{2}} \tilde{f} }_{L^\infty_{x,v}} + \normm{\tilde{\P} \tilde{f}}_{L^2_{x,v}} \Big)\nonumber \\
\lesssim & \left[\!\left[\tilde{f}_{0}\right]\!\right]_2 + \mathscr{E}_{2}[\tilde{f}](t) +  \mathscr{D}_{2}[\tilde{f}](t)+\d \e \normm{ \o f}^2_{L^{\infty}_{x,v}}+ \e \mathscr{E}^{\frac{1}{2}}_{2}[\tilde{f}](t). \label{pt-Phi-over-sqrt-mu}
\end{align}
A direct computation shows
\begin{equation}\label{pt-mu}
\begin{split}
\pt_t \tilde{\mu} = &\frac{\norm{v-\c}^2-3T}{2}\frac{\pt_t T}{T^2}\tilde{\mu}+ \frac{(v-\c)\cdot \pt_t \c}{T} \tilde{\mu} + \frac{\pt_t \rho}{\rho} \tilde{\mu}.
\end{split}\end{equation}
From \eqref{rho th u - smallness} and the exponential decay of $\tilde{\mu}$, we obtain
\begin{equation}\label{pt-mu-over-sqrt-mu}
\begin{split}
\normm{ \tilde{\nu}^{-\frac{1}{2}} \frac{\pt_t \tilde{\mu}}{\sqrt{\tilde{\mu}}} }_{L^{2}_{x,v}}
 \lesssim &  \a \mathfrak{h}_{1} + \a \norm{\tilde{f}}_{L^2_{\g_{+}}}
 \lesssim  \a \mathscr{E}_{2}^{\frac{1}{2}}(t)+ \a^{\frac{1}{2}} \e^{\frac{1}{2}}    \Big( \left[\!\left[\tilde{f}_{0}\right]\!\right]_2+  \mathscr{D}_{2}^{\frac{1}{2}}(t) \Big).
\end{split}
\end{equation}
Moreover, by \eqref{tilde-G-Pf} and \eqref{L-t-infty-L2-ipt-tilde-f},
\begin{equation} \label{Gamma-L-t-infty-L2}
\begin{split}
\normm{ \tilde{\nu}^{-\frac{1}{2}} \Gamma(\tilde{f},\tilde{f}) }_{L^{2}_{x,v}}
 \lesssim   &
\normm{ \tilde{\nu}^{-\frac{1}{2}} \Gamma(\tilde{f},\ipt \tilde{f}) }_{L^{2}_{x,v}}
+\normm{ \tilde{\nu}^{-\frac{1}{2}} \Gamma(\tilde{\P} \tilde{f},\tilde{\P} \tilde{f}) }_{L^{2}_{x,v}}
 \\
 \lesssim   &
  \normm{\omega^{\frac{1}{2}} \tilde{f}}_{L^\infty_{x,v}} \normm{\ipt \tilde{f}}_{L^2_{x,v}(\tilde{\nu})}+
 \normm{ \tilde{\P} \tilde{f} }_{L^{4}_{x,v}}^{2}
 %-------------
  \\
 \lesssim   &
  \d \e^2 \normm{\omega^{\frac{1}{2}} \tilde{f}}_{L^\infty_{x,v}}^2
  +\e^{-2}\normm{\ipt \tilde{f}}_{L^2_{x,v}(\tilde{\nu})}^2  +
 \normm{ \tilde{\P} \tilde{f} }_{L^{2}_{x,v}}^{\frac{3}{2}}
 \normm{ \tilde{\P} \tilde{f} }_{L^{6}_{x,v}}^{\frac{1}{2}}\\
 %-------------
 \lesssim  &
 \d \e^2 \normm{\omega f}_{L^\infty_{x,v}}^2 + \e^2\mathscr{E}_{2}[\tilde{f}](t)+\left[\!\left[\tilde{f}_{0}\right]\!\right]_2^2 +\mathscr{D}_{2}[\tilde{f}](t)
  +
 \mathscr{E}_{2}^{\frac{3}{2}}[\tilde{f}](t)+
 \d \normm{ \tilde{\P} \tilde{f} }_{L^{6}_{x,v}}
\end{split}
\end{equation}
hold for a sufficiently small constant $\d>0$ from Young's inequality.

Combining all the estimates above with \eqref{P-tilde-f-macro-L6} and absorbing the small term $\d \normm{ \tilde{\P} \tilde{f} }_{L^{6}_{x,v}}$ from \eqref{Gamma-L-t-infty-L2}, we arrive at \eqref{pf6 - bound3-final}.
\end{proof}
\medskip

Recall the definitions of $\tilde{L}$, $\tilde{\G}^{t}$ and $\tilde{L}^{t}$ in \eqref{tildeg - definition} and \eqref{tildegt - definition}. We have the following estimates.

\begin{corollary} \label{tildeL tildeGammat tildeLt - L2 estimate}
% Let $\e \in (0, \e_0]$, where $\e_0\in (0,1)$ is the constant from Proposition %\ref{L-infty-bd-unst}.
 Let $f$ be a solution of \eqref{tildef - Boltzmanneq ch1} on $[0,T]$ with $0<T\leq \infty$. Under the a priori assumption \eqref{theta-u-smallness-assumption}, the following estimates hold for all $t\in [0,T]$:
\begin{align}
&\normm{\tilde{L}f}_{L^2_{x,v}} \lesssim  \normm{\ipt {f}}_{L^2_{x,v}(\tilde{\nu})},\label{L-L2-estimate}\\
%--------------
&\normm{  \tilde{\nu}^{-\frac{1}{2}}  \tilde{\G}^{t}(f,g)}_{L^2_{x,v}}
\lesssim
 \e  \Big(\frac{\a}{\e}\mathfrak{h}_{1} + \frac{\a}{\e} \norm{\tilde{f}}_{L^2_{\g_{+}}} \Big) \nonumber\\
  & \qquad\qquad\qquad\qquad\quad  \times
\bigg[\Big( \frac{1}{\e} \normm{\ipt {g}}_{L^2_{x,v}} + \e^{\frac{3}{2}}\normm{ \o_1  {g}}_{L^\infty_{x,v}} + \normm{\tilde{\P} {g}}_{L^2_{x,v}} \Big) \normm{ \o_1   f}_{L^\infty_{x,v}} \label{Gt-L2-estimate}\\
&\qquad\qquad\qquad \;\; \qquad \qquad  +\Big( \frac{1}{\e} \normm{\ipt {f}}_{L^2_{x,v}} +\e^{\frac{3}{2}} \normm{  \o_1   {f} }_{L^\infty_{x,v}} + \normm{\tilde{\P} {f}}_{L^2_{x,v}} \Big) \normm{ \o_1   g}_{L^\infty_{x,v}}\bigg],\nonumber\\
%-------------
&\normm{  \tilde{\nu}^{-\frac{1}{2}}  \tilde{L}^{t}f}_{L^2_{x,v}}
\lesssim  \e  \Big(\frac{\a}{\e}\mathfrak{h}_{1} + \frac{\a}{\e} \norm{\tilde{f}}_{L^2_{\g_{+}}} \Big) \Big( \frac{1}{\e} \normm{\ipt {f}}_{L^2_{x,v}} + \normm{  \o_1    {f}}_{L^\infty_{x,v}} + \normm{\tilde{\P} {f}}_{L^2_{x,v}}  \Big), \label{Lt-L2-estimate}
\end{align}
where $\omega_1 = e^{\b_1 \norm{v}^2}$ with $0 < \b_1 \ll \frac{1}{4}$.
\end{corollary}

\begin{proof}[\textbf{Proof.}] \
By the property of $\tilde{L}$ and Lemma \ref{tildeGamma - L2 and Linfty estimate},
\begin{align*}
\normm{\tilde{L}f}_{L^2_{x,v}} %=& \normm{\tilde{L}(\ipt f)}_{L^2_{x,v}} \\
= & \normm{\tilde{\G}(\sqrt{\tilde{\mu}},\ipt f) +\tilde{\G}(\ipt f,\sqrt{\tilde{\mu}})}_{L^2_{x,v}}\\
\lesssim& \normm{ \ipt f}_{L^2_{x,v}(\tilde{\nu})} \normm{ \o_1  \sqrt{\tilde{\mu}}}_{L^\infty_{x,v}}
\lesssim \normm{\ipt f}_{L^2_{x,v}(\tilde{\nu})}.
\end{align*}
%From definition \eqref{tildegt - definition}, we write
%\begin{align*}
%\tilde{\G}^{t}(f,g)
%=& \frac{1}{\sqrt{\tilde{\mu}}}\Big[Q\big(f \pt_t \sqrt{\tilde{\mu}}, %g\sqrt{\tilde{\mu}}\big) + Q\big( f \sqrt{\tilde{\mu}}, g\pt_t %\sqrt{\tilde{\mu}}\big)\Big]
%= \tilde{\G}\Big(\frac{\pt_t \sqrt{\tilde{\mu}}}{\sqrt{\tilde{\mu}}}f, g\Big) + %\tilde{\G}\Big(f, \frac{\pt_t \sqrt{\tilde{\mu}}}{\sqrt{\tilde{\mu}}}g\Big).
%\end{align*}
%

 For $\tilde{\G}^{t}(f,g)$, Lemma \ref{tildeGamma - L2 and Linfty estimate} yields
\begin{align*}
\normm{  \tilde{\nu}^{-\frac{1}{2}}  \tilde{\G}^{t}(f,g)}_{L^2_{x,v}} \lesssim&
\normm{\tilde{\nu}^{\frac{1}{2}}\frac{\pt_t \sqrt{\tilde{\mu}}}{\sqrt{\tilde{\mu}}}f}_{L^2_{x,v}} \normm{ \o_1   g}_{L^\infty_{x,v}} + \normm{\tilde{\nu}^{\frac{1}{2}}\frac{\pt_t \sqrt{\tilde{\mu}}}{\sqrt{\tilde{\mu}}}g}_{L^2_{x,v}} \normm{ \o_1   f}_{L^\infty_{x,v}}.
\end{align*}
Combining this with Lemma \ref{v polynomial - L2 estimate} establishes \eqref{Gt-L2-estimate}.

For $\tilde{L}^{t}f $,
%\begin{align*}
%\tilde{L}^{t}f  = & -2 \big[\tilde{\G}^{t}(\sqrt{\tilde{\mu}},f) + \tilde{\G}^{t}(f, %\sqrt{\tilde{\mu}})\big]\\
%=& -2\Big[
% \tilde{\G}\big(\pt_t \sqrt{\tilde{\mu}}, f\big)
%+\tilde{\G}\big(f, \pt_t \sqrt{\tilde{\mu}}\big)
%+\tilde{\G}\Big(\sqrt{\tilde{\mu}}, \frac{\pt_t %\sqrt{\tilde{\mu}}}{\sqrt{\tilde{\mu}}}f\Big)
%+\tilde{\G}\Big(\frac{\pt_t \sqrt{\tilde{\mu}}}{\sqrt{\tilde{\mu}}}f, %\sqrt{\tilde{\mu}}\Big) \Big].
%\end{align*}
 Lemma \ref{tildeGamma - L2 and Linfty estimate} gives
\begin{align*}
\normm{  \tilde{\nu}^{-\frac{1}{2}}  \tilde{L}^{t}f}_{L^2_{x,v}} \lesssim&
\normm{\tilde{\nu}^{\frac{1}{2}}\frac{\pt_t \sqrt{\tilde{\mu}}}{\sqrt{\tilde{\mu}}}f}_{L^2_{x,v}} \normm{ \o_1   \sqrt{\tilde{\mu}}}_{L^\infty_{x,v}} + \normm{\tilde{\nu}^{\frac{1}{2}}\pt_t \sqrt{\tilde{\mu}}}_{L^2_{x,v}} \normm{ \o_1   f}_{L^\infty_{x,v}}.
\end{align*}
Together with \eqref{rho th u - smallness}, \eqref{pt-sqrt-mu} and Lemma
\ref{v polynomial - L2 estimate}, this proves \eqref{Lt-L2-estimate}.
\end{proof}
\medskip

We now prove Proposition \ref{tildeg tildegt - L2 estimate}.

\begin{proof}[\textbf{Proof of Proposition \ref{tildeg tildegt - L2 estimate}.}]\
The argument proceeds in three steps.
\medskip

\noindent\textbf{Step 1. Estimate for  \eqref{g-f-nonlinear-estimate}.}

Recall the definition of $\tilde{g}$ in \eqref{tildeg - definition}. We decompose
\begin{align*}
\int_{0}^{t} \iint_{\Omega \times \R^3} \tilde{g} \tilde{f}
=&\int_{0}^{t} \iint_{\Omega \times \R^3} \tilde{\G}(\tilde{f},\tilde{f})\ipt \tilde{f}+
 \int_{0}^{t} \iint_{\Omega \times \R^3}  \frac{\pt_t \tilde{\mu}}{\sqrt{\tilde{\mu}}}  \tilde{f} +
 \e \int_{0}^{t} \iint_{\Omega \times \R^3} \frac{ \pt_t \sqrt{\tilde{\mu}}} {\sqrt{\tilde{\mu}}} \tilde{f}^2\\
 := &I_{1} +I_{2} +I_{3}.
\end{align*}
%where
%\begin{align*}
%I_{1} :=&  \Big| \int_{0}^{t} \iint_{\Omega \times \R^3} \tilde{\G}(\tilde{f},\tilde{f}) %\tilde{f} \dd v \dd x  \dd s \Big|,\\
%I_{2} :=& \Big| \int_{0}^{t} \iint_{\Omega \times \R^3}  \frac{\pt_t %\tilde{\mu}}{\sqrt{\tilde{\mu}}}  \tilde{f} \dd v \dd x  \dd s\Big|,\\
%I_{3} :=& \e\Big| \int_{0}^{t} \iint_{\Omega \times \R^3} \frac{ \pt_t %\sqrt{\tilde{\mu}}} {\sqrt{\tilde{\mu}}} \tilde{f} \dd v \dd x  \dd s\Big|.
%\end{align*}
%Since the collision operator is orthogonal to $\tilde{\P}\tilde{f}$, and

For $I_1$, since the collision operator is orthogonal to $\tilde{\P}\tilde{f}$, Lemma \ref{tildeGamma - L2 and Linfty estimate} yields
\begin{align*}
|I_{1}| \lesssim & \e\int_{0}^{t} \normm{ \tilde{\nu}^{-\frac{1}{2}} \tilde{\G}(\tilde{f},\tilde{f})}_{L^2_{x,v}}^2 + \frac{1}{\e} \int_{0}^{t}\normm{\ipt \tilde{f}}_{L^2_{x,v}(\tilde{\nu})}^2.
\end{align*}
To estimate $I_2$, note \eqref{pt-mu} and using Proposition \ref{th tht w wt - ODE} that
\begin{equation*}
\normm{ \frac{\pt_t \tilde{\mu}}{\sqrt{\tilde{\mu}}} - \big[\pt_t \rho \tilde{\chi}_{0}  + \pt_t \c \cdot (\tilde{\chi}_{1}, \tilde{\chi}_{2}, \tilde{\chi}_{3}) + \pt_t \th \tilde{\chi}_{4}\big] }_{L^2_{x,v}} \lesssim \e \mathfrak{h}_{1}(\norm{\pt_t \rho} + \norm{\pt_t \th} + \norm{\pt_t w}).
\end{equation*}
Using the conservation laws \eqref{abc - average zero} and the estimate \eqref{rho th u - smallness}, we obtain
\begin{align*}
|I_{2}| \lesssim & \Big| \int_{0}^{t} \iint_{\Omega \times \R^{3}} \big[\pt_t \rho \tilde{\chi}_{0}  + \pt_t \c \cdot (\tilde{\chi}_{1}, \tilde{\chi}_{2}, \tilde{\chi}_{3}) + \pt_t \th \tilde{\chi}_{4}\big] \tilde{f} \Big|
%----------------------
 + \e \int_{0}^{t} \mathfrak{h}_{1}\Big(\a \mathfrak{h}_{1} + \a \norm{\tilde{f}}_{L^2_{\g_{+}}}\Big)\normm{\tilde{f}}_{L^2_{x,v}}\\
%----------------------
%=& \e \int_{0}^{t} \mathfrak{h}_{1}\Big(\a \mathfrak{h}_{1} + \a %\norm{\tilde{f}}_{L^2_{\g_{+}}}\Big)\normm{\tilde{f}}_{L^2_{x,v}} \dd s\\
%----------------------
\lesssim& \e^{2}\sup_{0 \le s \le t} \normm{\tilde{f}(s)}_{L^2_{x,v}}\frac{\a}{\e} \int_{0}^{t} \mathfrak{h}_{1}\Big(\a \mathfrak{h}_{1} + \a \norm{\tilde{f}}_{L^2_{\g_{+}}}\Big)
\lesssim \e^{2} \normmm{\tilde{f}(t)}_{2}^{3}.
\end{align*}
For $I_3$, using \eqref{pt-sqrt-mu}, \eqref{rho th u - smallness},
%\begin{align*}
%\norm{  \frac{ \pt_t \sqrt{\tilde{\mu}}} {\sqrt{\tilde{\mu}}} }
% \lesssim \big(1+\norm{v}+\norm{v}^2\big) \big( \norm{\pt_t \rho} + \norm{\pt_t \c} + %\norm{\pt_t \th}\big).
%\end{align*}
 Lemma \ref{v polynomial - L2 estimate} (with $\o_1=\o^\frac{1}{2}$) and Lemma \ref{tilde-f-bounded-by-f}, we have
\begin{align*}
|I_{3}| \lesssim & \e  \Big|\int_{0}^{t} \iint_{\Omega \times \R^3} \langle v \rangle^2 \big( \norm{\pt_t \rho} + \norm{\pt_t \c} + \norm{\pt_t \th}\big) \tilde{f}^{2} \Big|\\
\lesssim & \e^2 \int_{0}^{t} \Big(\frac{\a}{\e}\mathfrak{h}_{1} + \frac{\a}{\e} \norm{\tilde{f}}_{L^2_{\g_{+}}} \Big)
 \Big( \frac{1}{\e} \normm{\ipt \tilde{f}}_{L^2_{x,v}} +\e^{\frac{3}{2}} \normm{ \o_1 \tilde{f} }_{L^\infty_{x,v}} + \normm{\tilde{\P} \tilde{f}}_{L^2_{x,v}} \Big) \normm{\tilde{f}}_{L^2_{x,v}}\\
\lesssim& \e^2 \normmm{\tilde{f}(t)}_{2}^{3}.
\end{align*}

Combining the estimates for $I_1$,  $I_2$ and $I_3$ establishes \eqref{g-f-nonlinear-estimate}.
\medskip

\noindent\textbf{Step 2.  Estimate for \eqref{gt-ft-nonlinear-estimate}.}

Recall the definition of $\tilde{g}^{t}$ in \eqref{tildegt - definition}. We decompose
\begin{align*}
 \int_{0}^{t} \iint_{\Omega \times \R^3} \tilde{g}^{t} \pt_t \tilde{f}
=&  \int_{0}^{t} \iint_{\Omega \times \R^3}
\Big[\tilde{\G}(\pt_t \tilde{f}, \tilde{f})
+ \tilde{\G}(\tilde{f}, \pt_t \tilde{f})
+ \pt_t \Big(\frac{1}{\sqrt{\tilde{\mu}}} \Big)
\sqrt{\tilde{\mu}} \tilde{\G}(\tilde{f},\tilde{f}) \Big] \pt_t \tilde{f}\\
%----------
&  +  \int_{0}^{t} \iint_{\Omega \times \R^3} \Big[\e^{-1}\pt_t\Big(\frac{1}{\sqrt{\tilde{\mu}}} \Big) \sqrt{\tilde{\mu}} \tilde{L}\tilde{f} +\e^{-1}\tilde{L}^{t} \tilde{f}  + \tilde{\G}^{t}(\tilde{f},\tilde{f}) \Big] \pt_t \tilde{f}\\
%----------
&   +  \int_{0}^{t} \iint_{\Omega \times \R^3}\Big[ \pt_t \Big( \frac{\pt_t \tilde{\mu}}{\sqrt{\tilde{\mu}}}\Big) - \e \pt_t \Big(\frac{\pt_t \sqrt{\tilde{\mu}}}{\sqrt{\tilde{\mu}}}\Big) \tilde{f} - \e \Big(\frac{\pt_t \sqrt{\tilde{\mu}}}{\sqrt{\tilde{\mu}}}\Big)  \pt_t \tilde{f} \Big] \pt_t \tilde{f}\\
:=& II_{1}+II_{2}+II_{3}.
\end{align*}

For $II_1$, since the collision operator is orthogonal to $\tilde{\P}$, the first two terms are bounded by
\begin{align*}
%&\left| \int_{0}^{t} \iint_{\Omega \times \R^3}  \left[  \tilde{\G}(\pt_t \tilde{f}, %\tilde{f})
%+ \tilde{\G}(\tilde{f}, \pt_t \tilde{f}) \right] \pt_t \tilde{f} \dd v \dd x \dd t %\right|\\
%----------
%\le
%& \e \int_{0}^{t} \Big( \normm{ \tilde{\nu}^{-\frac{1}{2}} \tilde{\G}(\pt_t \tilde{f}, %\tilde{f})}_{L^2_{x,v}} + \normm{ \tilde{\nu}^{-\frac{1}{2}} \tilde{\G}(\tilde{f}, \pt_t %\tilde{f})}_{L^2_{x,v}} \Big) \frac{1}{\e}\normm{\ipt \pt_t %\tilde{f}}_{L^2_{x,v}(\tilde{\nu})} \dd s\\
%----------
%\lesssim
& \e\int_{0}^{t} \Big( \normm{ \tilde{\nu}^{-\frac{1}{2}} \tilde{\G}(\pt_t \tilde{f},\tilde{f})}_{L^2_{x,v}}^2 + \normm{ \tilde{\nu}^{-\frac{1}{2}} \tilde{\G}(\tilde{f},\pt_t \tilde{f})}_{L^2_{x,v}}^2  \Big) + \e \normmm{\tilde{f}(t)}_{2}^2.
\end{align*}
By Lemma \ref{v polynomial - L2 estimate} (with $\o_1=\o^\frac{1}{2}$), the third term in $II_1$ is bounded by
\begin{equation}\label{pt-mu-sqrt-mu-1}
\begin{split}
%& \left| \int_{0}^{t} \iint_{\Omega \times \R^3} \pt_t \Big(\frac{1}{\sqrt{\tilde{\mu}}} %\Big) \sqrt{\tilde{\mu}} \tilde{\G}(\tilde{f},\tilde{f})  \pt_t \tilde{f} \dd v \dd x \dd %s \right| \\
%----------
%\lesssim &
&
\e^{-1}\int_{0}^{t}\normm{ \tilde{\nu}^{\frac{1}{2}} \pt_t \Big(\frac{1}{\sqrt{\tilde{\mu}}} \Big) \sqrt{\tilde{\mu}} \pt_t \tilde{f}}_{L^2_{x,v}}^2
 +\e\int_{0}^{t}
\normm{ \tilde{\nu}^{-\frac{1}{2}} \tilde{\G}(\tilde{f},\tilde{f})}_{L^2_{x,v}}^2 \\
%----------
%\lesssim & \e \int_{0}^{t} %\normm{{\nu}^{-\frac{1}{2}}\tilde{\G}(\tilde{f},\tilde{f})}_{L^2_{x,v}}
%\Big(\frac{\a}{\e}\mathfrak{h}_{1} + \frac{\a}{\e} \norm{\tilde{f}}_{L^2_{\g_{+}}} %\Big)\\
%-----------------
%& \times \Big( \frac{1}{\e} \normm{\ipt \pt_t \tilde{f}}_{L^2_{x,v}(\tilde{\nu})} %+\e^{\frac{3}{2}} \normm{ w^{\frac{1}{2}} \pt_t \tilde{f} }_{L^\infty_{x,v}} + %\normm{\tilde{\P} \pt_t \tilde{f}}_{L^2_{x,v}} \Big)\dd s \\
%----------
\lesssim & \e  \int_{0}^{t} \Big(\frac{\a}{\e}\mathfrak{h}_{1} + \frac{\a}{\e} \norm{\tilde{f}}_{L^2_{\g_{+}}} \Big)^2\Big( \frac{1}{\e} \normm{\ipt \pt_t \tilde{f}}_{L^2_{x,v}(\tilde{\nu})} +\e^{\frac{3}{2}} \normm{ \o_1 \pt_t \tilde{f} }_{L^\infty_{x,v}} + \normm{\tilde{\P} \pt_t \tilde{f}}_{L^2_{x,v}} \Big)^2 \\
& + \e \int_{0}^{t} \normm{{\nu}^{-\frac{1}{2}}\tilde{\G}(\tilde{f},\tilde{f})}_{L^2_{x,v}}^2 \\
%----------------------
\lesssim & \sup_{0 \le s \le t}\Big(
\e^{\frac{3}{2}} \normm{\o_1 \pt_t \tilde{f} }_{L^\infty_{x,v}}
+ \normm{\pt_t \tilde{f}}_{L^2_{x,v}}
\Big)^2
\int_{0}^{t} \Big(\frac{\a}{\e}\mathfrak{h}_{1} + \frac{\a}{\e} \norm{\tilde{f}}_{L^2_{\g_{+}}} \Big)^2  \\
%----------------------
&+\sup_{0\le s \le t} \Big(\frac{\a}{\e}\mathfrak{h}_{1} + \frac{\a}{\e} \norm{\tilde{f}}_{L^2_{\g_{+}}} \Big)^2
\frac{1}{\e^2}\int_{0}^{t} \normm{ \ipt \pt_t \tilde{f}}_{L^2_{x,v}(\tilde{\nu})} ^2  + \e \int_{0}^{t} \normm{{\nu}^{-\frac{1}{2}}\tilde{\G}(\tilde{f},\tilde{f})}_{L^2_{x,v}}^2  \\
%----------------------
\lesssim & \e \normmm{\tilde{f}(t)}_{2}^{2}\Big(\normmm{\tilde{f}(t)}_{2}^{2} + \norm{\tilde{f}(0)}^2_{L^2_{\g_{+}}} \Big)+ \e \int_{0}^{t} \normm{{\nu}^{-\frac{1}{2}}\tilde{\G}(\tilde{f},\tilde{f})}_{L^2_{x,v}}^2.
\end{split}
\end{equation}
where we used \eqref{sup-f-gamma-L2}, Lemma \ref{tilde-f-bounded-by-f} and the definition of $\normmm{\tilde{f}(t)}_{2}$.
Combining these estimates yields
\begin{align*}
\norm{ II_{1}} \lesssim& \e\int_{0}^{t} \Big(
\normm{ \tilde{\nu}^{-\frac{1}{2}} \tilde{\G}(\pt_t \tilde{f},\tilde{f})}_{L^2_{x,v}}^2
+ \normm{ \tilde{\nu}^{-\frac{1}{2}} \tilde{\G}(\tilde{f},\pt_t \tilde{f})}_{L^2_{x,v}}^2
+ \normm{ \tilde{\nu}^{-\frac{1}{2}} \tilde{\G}(\tilde{f},\tilde{f})}_{L^2_{x,v}}^2 \Big) \\
&+ \e \normmm{\tilde{f}(t)}_{2}^{2}\Big(1+\normmm{\tilde{f}(t)}_{2}^{2}+ \left[\!\left[\tilde{f}_{0}\right]\!\right]_2^2\Big).
\end{align*}

Next, we estimate $II_2$. Similar to \eqref{pt-mu-sqrt-mu-1}, by Corollary \ref{tildeL tildeGammat tildeLt - L2 estimate}, the first term in $II_2$ is bounded by
\begin{align*}
%&  \left| \int_{0}^{t} \iint_{\Omega \times \R^3}\e^{-1} %\pt_t\Big(\frac{1}{\sqrt{\tilde{\mu}}} \Big) \sqrt{\tilde{\mu}} \tilde{L}\tilde{f} \pt_t %\tilde{f} \dd v \dd x \dd s \right| \\
%-------------
%\lesssim
 &  \e^{-1} \int_{0}^{t}  \normm{\tilde{\nu}^{\frac{1}{2}}\pt_t\Big(\frac{1}{\sqrt{\tilde{\mu}}} \Big) \sqrt{\tilde{\mu}} \pt_t \tilde{f}}_{L^2_{x,v}}^2
 +  \e^{-1} \int_{0}^{t} \normm{ \ipt f}_{L^2_{x,v}(\tilde{\nu})}^2 \\
%-------------
%\lesssim &  \e \int_{0}^{t} \normm{\e^{-1} \ipt f}_{L^2_{x,v}(\tilde{\nu})} %\Big(\frac{\a}{\e}\mathfrak{h}_{1} + \frac{\a}{\e} \norm{\tilde{f}}_{L^2_{\g_{+}}} \Big) %\\
%-------------
%&\quad \times   \Big( \normm{\ipt \pt_t \tilde{f}}_{L^2_{x,v}(\tilde{\nu})} %+\e^{\frac{3}{2}} \normm{ w^{\frac{1}{2}} \pt_t \tilde{f} }_{L^\infty_{x,v}} + %\normm{\tilde{\P} \pt_t \tilde{f}}_{L^2_{x,v}} \Big) \dd s.  \\
%-------------
\lesssim  & \e  \normmm{\tilde{f}(t)}_{2}^{2}\Big(1+\normmm{\tilde{f}(t)}_{2}^{2} + \norm{\tilde{f}(0)}^2_{L^2_{\g_{+}}}  \Big).
\end{align*}
Since $\tilde{L}^{t}$ is orthogonal to $\tilde{\P}$, by \eqref {Lt-L2-estimate}  (with $\o_1=\o^\frac{1}{2}$) and Lemma \ref{tilde-f-bounded-by-f}, the second term in $II_2$ is bounded by
\begin{align*}
&   \e^{-1}\int_{0}^{t} \normm{\tilde{\nu}^{-\frac{1}{2}} \tilde{L}_t\tilde{f}}_{L^2_{x,v}}^2
+\e^{-1}\int_{0}^{t}\normm{ \ipt \pt_t \tilde{f}}_{L^2_{x,v}(\tilde{\nu})}^2 \\
%-------------
\lesssim &  \e \int_{0}^{t} \Big(\frac{\a}{\e}\mathfrak{h}_{1} + \frac{\a}{\e} \norm{\tilde{f}}_{L^2_{\g_{+}}} \Big)^2 \Big( \normm{\e^{-1} \ipt \tilde{f}}_{L^2_{x,v}}  +\normm{\omega_1 \tilde{f}}_{L^\infty_{x,v}} + \normm{\tilde{\P}\tilde{f}}_{L^2_{x,v}}  \Big)^2 \\
&+ \e \int_{0}^{t} \normm{\e^{-1} \ipt \pt_t \tilde{f}}_{L^2_{x,v}(\tilde{\nu})}^2  \\
%-------------
\lesssim  & \e^{\frac{1}{2}} \normmm{\tilde{f}(t)}_{2}^{2}\Big(1+\normmm{\tilde{f}(t)}_{2}^{2} + \norm{\tilde{f}(0)}^2_{L^2_{\g_{+}}} \Big).
\end{align*}
By \eqref{Gt-L2-estimate} (with $\o_1=\o^\frac{1}{2}$) and Lemma \ref{tilde-f-bounded-by-f}, the third term in $II_2$ is controlled as
 \begin{align*}
%&  \left| \int_{0}^{t} \iint_{\Omega \times \R^3}  \tilde{\G}^{t}(f,f) \pt_t \tilde{f}\dd %v \dd x \dd s \right| \\
%-------------
 &  \e \int_{0}^{t} \normm{\e^{-1} \ipt \pt_t \tilde{f}}_{L^2_{x,v}(\tilde{\nu})} \Big(\frac{\a}{\e}\mathfrak{h}_{1} + \frac{\a}{\e} \norm{\tilde{f}}_{L^2_{\g_{+}}} \Big)  \normm{\omega \tilde{f}}_{L^\infty_{x,v}} \\
%-------------
&\qquad \times \Big(  \normm{\e^{-1} \ipt \tilde{f}}_{L^2_{x,v}} + \e^{\frac{3}{2}} \normm{ \o_1  \tilde{f}}_{L^\infty_{x,v}} + \normm{\tilde{\P} \tilde{f}}_{L^2_{x,v}} \Big) \dd s.  \\
%-------------
\lesssim  &  \normmm{\tilde{f}(t)}_{2}\Big(\normmm{\tilde{f}(t)}_{2} + \norm{\tilde{f}(0)}_{L^2_{\g_{+}}} \Big)   \e^{-\frac{3}{2}} \int_{0}^{t} \Big(  \normm{\ipt \pt_t \tilde{f}}^2_{L^2_{x,v}(\tilde{\nu})}
+ \normm{\ip \tilde{f}}^2_{L^2_{x,v}(\tilde{\nu})} \Big)
\\
& + \e^{\frac{1}{2}}  \normmm{\tilde{f}(t)}_{2}^{2}\normmm{\tilde{f}(t)}_{2}^{2}\\
%-------------
\lesssim  & \e^{\frac{1}{2}}  \normmm{\tilde{f}(t)}_{2}^{2}\Big(\normmm{\tilde{f}(t)}_{2}^{2} + \norm{\tilde{f}(0)}^{2}_{L^2_{\g_{+}}} \Big).
\end{align*}
 Collecting the above estimates yields
\begin{align*}
\norm{ II_{2}} \lesssim &  \e^{\frac{1}{2}}  \normmm{\tilde{f}(t)}_{2}^{2}\Big(1+\normmm{\tilde{f}(t)}_{2}^{2} + \left[\!\left[\tilde{f}_{0}\right]\!\right]_2^2 \Big).
\end{align*}

Finally, we estimate $II_3$. For the first term, note that
\begin{align*}
&\pt_t \left( \frac{\pt_t \tilde{\mu}}{\sqrt{\tilde{\mu}}} \right) = \pt_t \Big[\Big( \frac{\norm{v-\c}^2-3T}{2}\Big)\frac{\pt_t \th}{(1+\th)^2} \sqrt{\tilde{\mu}} + \frac{(v-\c)\cdot \pt_t \c}{1+\th} \sqrt{\tilde{\mu}} + \pt_t \rho \sqrt{\tilde{\mu}} \Big],\\
%-------------
&\pt_t \tilde{f} = \frac{1}{\sqrt{\tilde{\mu}}}\Big[ \pt_t (\sqrt{\tilde{\mu}} \tilde{f}) - \pt_t (\sqrt{\tilde{\mu}}) \tilde{f} \Big],
\end{align*}
and $\pt_t(\sqrt{\tilde{\mu}} \tilde{f})$ also satisfies the conservation laws \eqref{abc - average zero}.
%
%\begin{equation}\label{pt-tilde-f-conservation-law}
%\begin{split}
%\iint_{\Omega \times \R^{3}} \pt_t(\sqrt{\tilde{\mu}} \tilde{f}) \dd v \dd x
%&=\pt_t \iint_{\Omega \times \R^{3}} \tilde{\chi}_{0} \tilde{f} \dd v \dd x    = 0,\\
%-------------
%\iint_{\Omega \times \R^{3}} Ax \cdot (v_{1}, v_{2}, v_{3}) \pt_t(\sqrt{\tilde{\mu}} %\tilde{f}) \dd v \dd x
%&= \pt_t \iint_{\Omega \times \R^{3}} Ax \cdot (\tilde{\chi}_{1}, \tilde{\chi}_{2}, %\tilde{\chi}_{3}) \tilde{f} \dd v \dd x  =   0,\\
%-------------
%\iint_{\Omega \times \R^{3}} \frac{\norm{v}^2-3}{\sqrt{6}} \pt_t(\sqrt{\tilde{\mu}} %\tilde{f}) \dd v \dd x
%&= \pt_t \iint_{\Omega \times \R^{3}} \tilde{\chi}_{4} \tilde{f} \dd v \dd x   = 0.
%\end{split}
%\end{equation}
Consequently, all linear terms in $\pt_t \left( \frac{\pt_t \tilde{\mu}}{\sqrt{\tilde{\mu}}} \right)$ (e.g., $ \pt_t \pt_t \th \norm{v}^2  \sqrt{\tilde{\mu}}$, $\pt_t \pt_t \c \cdot  v\sqrt{\tilde{\mu}}$ and $\pt_t \pt_t \rho \sqrt{\tilde{\mu}} $) are orthogonal to $\pt_t(\sqrt{\tilde{\mu}} \tilde{f})$. Therefore, only the remaining nonlinear terms contribute, giving
\begin{align*}
& \Big| \int_{0}^{t} \iint_{\Omega \times \R^3} \pt_t \Big( \frac{\pt_t \tilde{\mu}}{\sqrt{\tilde{\mu}}} \Big)  \pt_t \tilde{f} \Big|
%-------------
\lesssim  \e \int_{0}^{t} \mathfrak{h}_{1}  \Big ( \a \mathfrak{h}_{1,t} + \a \norm{\pt_t \tilde{f}}_{L^2_{\g_{+}}} \Big )
\Big( \normm{ \pt_t \tilde{f} }_{L^2_{x,v}} + \normm{ \tilde{f}}_{L^2_{x,v} } \Big )   %-------------
\lesssim  \e^{2}\normmm{\tilde{f}(t)}_{2}^{3}.
\end{align*}
Next, similar to the proof of Lemma \ref{v polynomial - L2 estimate}, the second (cubic) term in $II_3$ is bounded by
\begin{align*}
%&  \Big|  \int_{0}^{t} \iint_{\Omega \times \R^3} \e \pt_t \Big ( \frac{\pt_t %\sqrt{\tilde{\mu}}}{\sqrt{\tilde{\mu}}}\Big ) \tilde{f} \pt_t \tilde{f} \Big| \\
%\lesssim
& \e \int_{0}^{t} \iint_{\Omega \times \R^3} \big|\tilde{f} \pt_t \tilde{f}\big| \langle v\rangle^4  \Big( \norm{\pt_t \rho}^2 + \norm{\pt_t \c}^2 + \norm{\pt_t \th}^2 + \norm{\pt_t \pt_t \rho} + \norm{\pt_t \pt_t \c} + \norm{\pt_t \pt_t \th}\Big)  \\
\lesssim & \e \int_{0}^{t}  \normm{\pt_t \tilde{f}}_{L^2_{x,v}} \Big (\a \mathfrak{h}_{1,t} + \a \norm{\pt_t \tilde{f}}_{L^2_{\g_{+}}} \Big )
 \Big(\e^{-1} \normm{\ipt \tilde{f}}_{L^2_{x,v}} + \e^{\frac{3}{2}} \normm{\o_1 \tilde{f}}_{L^\infty_{x,v}} + \normm{\tilde{\P} \tilde{f}}_{L^2_{x,v}} \Big) \\
\lesssim& \e^{2} \normmm{\tilde{f}(t)}_{2}^{3},
\end{align*}
where we used Proposition \ref{th tht w wt - ODE} and Lemma \ref{tilde-f-bounded-by-f}.
Finally, the last term in $II_3$ is controlled as
\begin{align*}
%\Big|  \int_{0}^{t} \iint_{\Omega \times \R^3} \e \frac{\pt_t %\sqrt{\tilde{\mu}}}{\sqrt{\tilde{\mu}}} \pt_t \tilde{f} \pt_t \tilde{f}\Big|
%\lesssim &
&\e \int_{0}^{t} \iint_{\Omega \times \R^3} |\pt_t \tilde{f}|^{2} \langle v\rangle^2 \left( \norm{\pt_t \rho} + \norm{\pt_t \c} + \norm{\pt_t \th}\right)  \\
\lesssim & \e^2 \int_{0}^{t}   \normm{\pt_t \tilde{f}}_{L^2_{x,v}} \Big(\frac{\a}{\e}\mathfrak{h}_{1} + \frac{\a}{\e} \norm{\tilde{f}}_{L^2_{\g_{+}}} \Big)
 \Big(  \e^{-1} \normm{ \ipt \pt_t\tilde{f}}_{L^2_{x,v}} + \e^{\frac{3}{2}} \normm{ \o_1  \pt_t\tilde{f}}_{L^\infty_{x,v}} + \normm{\tilde{\P} \pt_t\tilde{f}}_{L^2_{x,v}} \Big) \\
\lesssim& \e^2 \normmm{\tilde{f}(t)}_{2}^{3}.
\end{align*}
Collecting these estimates gives $\norm{ II_{3} } \lesssim \e^{2} \normmm{\tilde{f}(t)}_{2}^{3}$.

 Combining the bounds for  $II_{1}, II_{2}$ and $II_{3}$ establishes \eqref{gt-ft-nonlinear-estimate}.
\medskip

\noindent\textbf{Step 3.  Estimate for \eqref{g-1over-w-estimate} and \eqref{gt-1over-w-estimate}.}

For $\normm{\tilde{g}\omega^{-1}}_{L^2_{x,v}}^2$ and $\normm{\tilde{g}^{t} \omega^{-1}}_{L^2_{x,v}}^2$, the algebraic growth in $v$ is absorbed by the exponential decay of $\omega^{-1}$. Therefore, using Lemma \ref{r ptr s - smallness} and arguing as in  Steps 1 and 2, we obtain \eqref{g-1over-w-estimate} and  \eqref{gt-1over-w-estimate}. The details are omitted for brevity.

This completes the proof of Proposition \ref{tildeg tildegt - L2 estimate}.
\end{proof}
\bigskip

%\begin{align*}
%\int_{0}^{t} \normm{\tilde{g}\omega^{-1}}_{L^2_{x,v}}^2 \dd s \lesssim \int_{0}^{t} %\normm{ \tilde{\nu}^{-\frac{1}{2}} \tilde{\G}(\tilde{f},\tilde{f})}_{L^2_{x,v}}^2\dd s + %\e^2 \Big(\normmm{\tilde{f}(t)}_{2}^2+\normmm{\tilde{f}(t)}_{2}^4\Big).
%\end{align*}
%\begin{align*}
%\int_{0}^{t} \normm{\tilde{g}^{t}\omega^{-1}}_{L^2_{x,v}}^2 \dd t \lesssim
%& \int_{0}^{t} \Big( \normm{ \tilde{\nu}^{-\frac{1}{2}}  \tilde{\G}(\pt_t %\tilde{f},\tilde{f})}_{L^2_{x,v}}^2 + \normm{\tilde{\nu}^{-\frac{1}{2}} %\tilde{\G}(\tilde{f},\pt_t \tilde{f})}_{L^2_{x,v}}^2 + \normm{\tilde{\nu}^{-\frac{1}{2}}  %\tilde{\G}(\tilde{f},\tilde{f})}_{L^2_{x,v}}^2\Big)\dd s\\
%& + \normmm{\tilde{f}(t)}_{2}^2 \Big(\e^2 + \normmm{\tilde{f}(t)}_{2}^2 + %\left[\!\left[\tilde{f}_{0}\right]\!\right]_2^2\Big),
%\end{align*}
%where we also used
%\begin{equation}\label{pt-mu-sqrt-mu-2}
%\begin{split}
%&\int_{0}^{t} \Big(\frac{\a}{\e}\mathfrak{h}_{1} + \frac{\a}{\e} %\norm{\tilde{f}}_{L^2_{\g_{+}}} \Big)^2
%\frac{1}{\e^2} \normm{\ipt \tilde{f}}^2_{L^2_{x,v}(\tilde{\nu})}
%  \dd s\\
%----------------------
%\lesssim & \sup_{0\le s \le t} \Big(\frac{\a}{\e}\mathfrak{h}_{1} + \frac{\a}{\e} %\norm{\tilde{f}}_{L^2_{\g_{+}}} \Big)^2
%\frac{1}{\e^2}\int_{0}^{t} \normm{ \ipt  \tilde{f}}_{L^2_{x,v}(\tilde{\nu})} ^2  \dd s \\
%----------------------
%\lesssim & \normmm{\tilde{f}(t)}_{2}^{2}\Big(\normmm{\tilde{f}(t)}_{2}^{2} + %\norm{\tilde{f}(0)}^2_{L^2_{\g_{+}}} \Big).
%\end{split}
%\end{equation}

\subsection{Proof of Main Result for the Case $0 \leq \a \ll \e$} \label{proof-main-2}\
\medskip

In this subsection, we give the proof of Theorem \ref{main-th-2}.
\medskip

\begin{proof}[\textbf{Proof of Theorem \ref{main-th-2}}]  \
In the regime $0 \leq \a \ll \e$, we work with the perturbation equation \eqref{tildef - Boltzmanneq ch1} around the rotating Maxwellian $\tilde{\mu}$.
The argument follows the same pattern as that of Theorem \ref{main-th-1}. For conciseness,  we only point out the main differences and omit most of the repetitive details.

\medskip

\noindent\textbf{Step 1. Global existence and uniform $\e$-independent estimates.}

To obtain the global a priori estimate \eqref{uniform-bound-tilde}, we follow the argument from Step 1 in  the proof of Theorem \ref{main-th-1}.

First, applying Corollary \ref{tildeG - L2 Linfty L3 L6 estimate} and Proposition \ref{f - L2L3 estimate} with source terms $g=-\e^{-1}\tilde{L}\tilde{f}+\tilde{g}$ (for $S_1\tilde{f}$) and $g=-\e^{-1}\tilde{L}\pt_t\tilde{f}+\tilde{g}_t$ (for $S_1\pt_t \tilde{f}$), and then using Proposition \ref{tildeg tildegt - L2 estimate}, we obtain
\begin{equation}\label{Gamma-L2-estimate}
\begin{split}
 &\normm{ \tilde{\nu}^{-\frac{1}{2}} \tilde{\G}(\tilde{f}, \tilde{f})}^2_{L^2_{t,x,v}}
+\normm{ \tilde{\nu}^{-\frac{1}{2}} \tilde{\G}(\tilde{f},\pt_t \tilde{f})}^2_{L^2_{t,x,v}}
+\normm{ \tilde{\nu}^{-\frac{1}{2}} \tilde{\G}(\pt_t \tilde{f},\tilde{f})}^2_{L^2_{t,x,v}} %------------------
\lesssim  \left[\!\left[\tilde{f}_{0}\right]\!\right]_{2}^2 \normmm{\tilde{f}}_{2}^{2}(t) + \normmm{\tilde{f}}_{2}^{4}(t).
\end{split}
\end{equation}

Second, multiplying the estimate \eqref{P-tilde-f-macro-L2} in Proposition \ref{Psi - L2 and L6 estimate} and the estimate \eqref{P-pt-tilde-f-macro-L2} in Corollary \ref{Psi - L2 and L6 estimate-pt} by a small constant, and adding the result to the estimates \eqref
{tilde-f-energy-estimate} and \eqref
{tilde-ft-energy-estimate} in Proposition \ref{tildef tildeft - Energy estimate}, we deduce
\begin{equation}\label{EfDf - bound1}
\begin{split}
\mathscr{E}_{2}[\tilde{f}](t)  + \mathscr{D}_{2}[\tilde{f}](t)
\lesssim
\left[\!\left[\tilde{f}_{0}\right]\!\right]_{2}^2
+ \left[\!\left[\tilde{f}_{0}\right]\!\right]_{2}^2  \normmm{\tilde{f}}_{2}^{2}(t)
+ \normmm{\tilde{f}}_{2}^{3}(t)
+ \normmm{\tilde{f}}_{2}^{4}(t).
\end{split}
\end{equation}

Third, applying Proposition \ref{w-tildef-tildeft-infty-estimate} and Lemma \ref{tildeGamma - L2 and Linfty estimate} gives
\begin{equation} \label{wfinfty - bound2}
\begin{split}
\e\normm{\omega {f}}_{L^{\infty}_{t,x,v}}^2 + \e^{3}\normm{\omega \pt_t {f}}_{L^{\infty}_{t,x,v}}^{2}
\lesssim & \left[\!\left[\tilde{f}_{0}\right]\!\right]_{2}^2
+ \mathscr{E}_{2}[\tilde{f}](t)
+ \mathscr{D}_{2}[\tilde{f}](t)
%------------
+ \normmm{\tilde{f}(t)}_{2}^{4}
+ \normm{\tilde{\P} \tilde{f}}_{L^{\infty}_{t}L^6_{x,v}}^2.
\end{split}
\end{equation}
Using Proposition \ref{Psi - L2 and L6 estimate}, Corollary \ref{tilde-f-L6-estimate-final} and Proposition \ref{w-tildef-tildeft-infty-estimate}, we derive the bound
\begin{equation} \label{pf6 - bound3}
\begin{split}
\normm{\tilde{\P} \tilde{f}}_{L^{\infty}_{t}L^6_{x,v}}^2
\lesssim &\left[\!\left[\tilde{f}_{0}\right]\!\right]_{2}^2
+\left[\left[\tilde{f}_{0}\right]\right]_{2}^4
+ \mathscr{E}_{2}[\tilde{f}](t)
+ \mathscr{D}_{2}[\tilde{f}](t)
+ \normmm{\tilde{f}}_{2}^{4}(t)
%------------
+ \normmm{\tilde{f}}_{2}^{6}(t)
+ \d\e\normm{\omega {f}}_{L^{\infty}_{t,x,v}}^2,
\end{split}
\end{equation}
where $\d>0$ is a sufficiently small constant arising from Corollary \ref{tilde-f-L6-estimate-final}.
Combining \eqref{wfinfty - bound2} and \eqref{pf6 - bound3} and absorbing the
terms  $\d\e\normm{\omega {f}}_{L^{\infty}_{t,x,v}}^2$ and $\normm{\tilde{\P} \tilde{f}}_{L^{\infty}_{t}L^6_{x,v}}^2$ on the right-hand side, we obtain
\begin{equation}\label{infty-macro-bound}
\begin{split}
&\e\normm{\omega {f}}_{L^{\infty}_{t,x,v}}^2 + \e^{3}\normm{\omega \pt_t {f}}_{L^{\infty}_{t,x,v}}^{2}  + \normm{\tilde{\P} \tilde{f}}_{L^{\infty}_{t}L^6_{x,v}}^2
\\
\lesssim &
\left[\!\left[\tilde{f}_{0}\right]\!\right]_{2}^2
+\left[\!\left[\tilde{f}_{0}\right]\!\right]_{2}^4
+ \mathscr{E}_{2}[\tilde{f}](t)
+ \mathscr{D}_{2}[\tilde{f}](t)
+ \normmm{\tilde{f}}_{2}^{4}(t)
+ \normmm{\tilde{f}}_{2}^{6}(t).
\end{split}
\end{equation}

Finally, multiplying \eqref{infty-macro-bound} by a small constant, adding the result to \eqref{EfDf - bound1} and absorbing small contributions on the right, we find that
\begin{equation}
\normmm{\tilde{f}}_{2}^2(t)
\lesssim
\left[\!\left[\tilde{f}_{0}\right]\!\right]_{2}^2
+ \normmm{\tilde{f}}_{2}^{3}(t)
+ \normmm{\tilde{f}}_{2}^{4}(t)
+ \normmm{\tilde{f}}_{2}^{6}(t)
\end{equation}
holds for any $0\leq t\leq T$, provided  $\left[\!\left[\tilde{f}_{0}\right]\!\right]_{2}^2 \leq \d_0$ is sufficiently small. Consequently, the a priori assumption \eqref{theta-u-smallness-assumption} is verified if $\d_0$ is chosen further small such that $\d_0\ll \d_1$. The global a priori estimate \eqref{uniform-bound-tilde} on $[0,\infty)$ is then establish via standard continuity argument.
\medskip

\noindent \textbf{Step 2. Derivation of strong convergence \eqref{tilde-f-strong-convergence}--\eqref{momentums-strong-convergence}  and INSF system \eqref{INSF-unst}.}

The uniform bound on $\normmm{\tilde{f}}_{2}(\infty)$ given by \eqref{uniform-bound-tilde} implies:
\begin{align}
&\sup_{0\leq s\leq \infty}\Big(
 \normm{\tilde{f}(s)}_{L^2_{x,v}}
+\normm{\pt_t \tilde{f}(s)}_{L^2_{x,v}}
+\normm{\tilde{\P}\tilde{f}(s)}_{L^6_{x,v}}
\Big)\leq C\delta_0, \label{tilde-f-uniform-bound} \\
%-------------
&\sup_{0\leq s\leq \infty}\Big(
 \norm{\frac{\th(s)}{\e}}
+\norm{\frac{w_i(s)}{\e}}
+\norm{\frac{\pt_t \th(s)}{\e}}
+\norm{\frac{\pt_t{w_i}(s)}{\e}}
\Big)\leq C\delta_0, \label{th-w-uniform-bound}\\
%----------------------
& \int_{0}^{\infty} \Big( \normm{\tilde{\P} \tilde{f}(s)}_{L^2_{x,v}}^{2}
+  \int_{0}^{t} \normm{\pt_t \tilde{\P} \tilde{f}(s)}_{L^2_{x,v}}^{2} \Big)\dd s
\leq C\delta_0,\label{tilde-f-L2-bound} \\
%----------------
&\int_{0}^{\infty}
\Big(
\normm{\ipt \tilde{f}(s)}_{L^2_{x,v}(\tilde{\nu})}^{2}
+\normm{\ipt \pt_t \tilde{f}(s)}_{L^2_{x,v}(\tilde{\nu})}^{2}
\Big) \dd s  \to  0\;\;\; \hbox{ as } \; \e\to  0.\label{ipt-to-0}
\end{align}
Hence, there exist $f^{*}\in L^\infty\left( \mathbb{R}^+; L^2(\O\times \mathbb{R}^3)\right)$ and $\theta^*, w_i^* \in L^\infty(\mathbb{R}^+)$ such that, up to a subsequence,
\begin{align}
\tilde{f} \to  f^{*} \;\; &\text{~~weakly}\!-\!* ~\text{in}~ L^\infty\left( \mathbb{R}^+; L^2(\O\times \mathbb{R}^3)\right), \label{f-es-L2-lim}\\
\frac{\th}{\e}  \to {\th^{*}}, \;\;
\frac{w_{i}}{\e}\to  {w_{i}^{*}}  \;\;& \text{~~weakly}\!-\!* ~\text{in} ~L^\infty(\mathbb{R}^+), \label{th-lim}\\
\frac{\pt_t \th}{\e} \to {\pt_t \th^{*}},\;\;
\frac{\pt_t w_{i}}{\e} \to {\pt_t w_{i}^{*}}  \;\; &\text{~~weakly}\!-\!* ~\text{in} ~L^\infty(\mathbb{R}^+), \label{th-pt-lim}\\
\frac{\c }{\e} \to {\c^{*}},\;\;
\frac{\pt_t \c}{\e} \to {\pt_t \c^{*}} \;\; &\text{~~weakly}\!-\!* ~\text{in} ~L^\infty\left( \mathbb{R}^+; L^\infty(\O)\right)  \label{u-lim}
\end{align}
as $\varepsilon \to 0$, where $\c^{*}=\sum w_{i}^{*} A_{i}x$. For notational simplicity, $w_{i}(t)$ denotes either $w_{i}(t)$ or $w(t)$, and similarly for $w^{*}_{i}(t)$.
Owing to the initial conditions \eqref{parameters-initials}, we have
\begin{equation}\label{th0-initial}
\th^{*}(0)=0,\;\;  w_{i}^{*}(0) =0,\;\; \c^{*}(0,x)=0 \;\;\forall x\in\O.
\end{equation}

Using a Taylor expansion of $\tilde{\mu}$ together with \eqref{th-w-uniform-bound} and \eqref{u-lim} yields
\begin{align}
\frac{\o( \tilde{\mu} - \mu)}{\e} \to \o\Big( \c^{*} \cdot v + \th^{*}\frac{\norm{v}^2 -3}{2}\Big)\mu \;\;
&\hbox{ strongly in } L^\infty\left(\R_{+}; L^1\cap L^\infty(\Omega \times \R^3)\right), \label{tildemu - mu-over-e}\\
%----------------
\o \tilde{\mu} \to \o\mu,\qquad  \o \sqrt{\tilde{\mu}} \to \o\sqrt{\mu}\;\;
&\hbox{ strongly in } L^\infty\left(\R_{+}; L^1\cap L^\infty(\Omega \times \R^3)\right) \label{tildemu - mu}
\end{align}
as $\e \to 0$, where $\o$ is the weight function defined in
\eqref{weight-w}.
The convergence \eqref{ipt-to-0} and \eqref{L-L2-estimate} in Lemma \ref{tildeL tildeGammat tildeLt - L2 estimate} imply $\tilde{L}\tilde{f} \to 0$ strongly in $L^2( \mathbb{R}^+\times \O\times \mathbb{R}^3)$.
Moreover, \eqref{f-es-L2-lim} and \eqref{tildemu - mu} indicate
$\tilde{L}\tilde{f} \to Lf^{*}$ in the sense of distributions. By uniqueness of distribution limits, we obtain $Lf^{*}=0$. Hence, there exist functions $\varrho_{f^{*}},u_{f^{*}},\vartheta_{f^{*}}\in L^\infty ( \mathbb{R}^+; L^2( \O))$ such that
\begin{equation} \label{Ltildef and tildef}
 f^{*} = \Big(\varrho_{f^{*}} + u_{f^{*}} \cdot v + \vartheta_{f^{*}} \frac{\norm{v}^2-3}{2}\Big)\sqrt{\mu}.
\end{equation}
%Furthermore, the uniform bound on $\normmm{\tilde{f}}_{2}(\infty)$ implies that  %$\tilde{\G}(\tilde{f},\tilde{f})$, $\e^{-1}\tilde{L}\tilde{f}$ and $\pt_t \tilde{f}$ are %bounded in $L^2\left( \mathbb{R}^+; L^2(\O\times \mathbb{R}^3)\right)$. Consequently, $v %\cdot \nabla_x \tilde{f}$ is also bounded in $L^2\left( \mathbb{R}^+; L^2(\O\times %\mathbb{R}^3)\right)$ and admits a weak limit.
%On the other hand, \eqref{tilde-f-uniform-bound} implies that
%$v \cdot \nabla_x \tilde{f} \to v \cdot \nabla_x f^{*}$ in the sense of distributions.
%By uniqueness of distribution limit,
Proceeding as in the proof of Theorem \ref{main-th-1}, we also have $\varrho_{f^{*}},u_{f^{*}},\vartheta_{f^{*}} \in L^2\left( \mathbb{R}^+;  H^{1}(\O)\right)$ and
\begin{equation}\label{nabla-f-convergence}
\tilde{\nu}^{-\frac{1}{2}}v \cdot \nabla_x \tilde{f} \to \nu^{-\frac{1}{2}}v \cdot \nabla_x f^{*} \;\; \text{ weakly} ~\text{in}\; L^2\left( \mathbb{R}^+; L^2(\O\times \mathbb{R}^3)\right) \text{ as } \e\to 0.
\end{equation}

We now claim that
\begin{equation} \label{th00-u00}
\begin{split}
  \vartheta^{*}(t) \equiv 0,  \quad w_{i}^{*}(t) \equiv 0 \; (i=1,2,3), \quad \c^{*}(t,x) \equiv 0.
\end{split}
\end{equation}
To this end, observe the identity
\begin{equation*}
\begin{split}
\iint_{\partial\Omega\times\mathbb{R}^3}\tilde{\nu}^{-\frac{1}{2}} \phi \tilde{f}|_{\partial\Omega}\dd \g
&=\iint_{\Omega\times\mathbb{R}^3}\tilde{\nu}^{-\frac{1}{2}}(v\cdot\nabla_x\phi) \tilde{f}
+\iint_{\Omega\times\mathbb{R}^3}\tilde{\nu}^{-\frac{1}{2}}(v\cdot\nabla_x \tilde{f})\phi,
\end{split}
\end{equation*}
where $\phi(x, v)$ is test function satisfying $\phi(\cdot, v)\in C^\infty(\bar{\Omega})$ and $\phi(x, \cdot)\in
C^\infty_0(\mathbb{R}^3)$.  Combining this with \eqref{f-es-L2-lim} and \eqref{nabla-f-convergence} implies
\begin{equation}\label{g-gamma-distr-lim}
\begin{split}
\tilde{\nu}^{-\frac{1}{2}} \tilde{f}|_{\partial\Omega} \to \nu^{-\frac{1}{2}} f^{*} |_{\partial\Omega} \quad  \hbox{ in the sense of distributions} \text{ as } \e\to 0.
\end{split}
\end{equation}
The uniform bound of $\normmm{\tilde{f}}_{2}(\infty)$ implies that
$\displaystyle \sqrt{\frac{\a}{\e}} \norm{\tilde{f}}_{L^2_tL^2_{\gamma,+}}$ is uniformly bounded and hence, up to a subsequence, has a weak limit in $L^2(\mathbb{R}_+\times\dd \gamma)$.
 By \eqref{a-e-limit-0}, \eqref{g-gamma-distr-lim} and the uniqueness of distribution limits, we conclude
\begin{equation}\label{boundary-converge-zero}
\begin{split}
 \sqrt{\frac{\a}{\e}}  \tilde{\nu}^{-\frac{1}{2}} \tilde{f}\big|_{\pt\O}  \to
 0 \quad\hbox{weakly in }  L^2(\mathbb{R}_+\times\dd \gamma)  \text{ as } \e\to 0.
\end{split}
\end{equation}
Recall the ODEs for $\th^{\e}$ and $w_{i}^{\e}$ in Proposition \ref{th tht w wt - ODE}.
Passing to the limit in  \eqref{th - axi 1 estimate} and \eqref{w - axi 1 estimate} and using \eqref{th-lim}, we derive
\begin{equation}\label{theta0-equation}
\begin{split}
&3 \frac{\dd }{\dd t} \int_{\Omega} \theta^{*} \dd x  + \l    \int_{\pt \Omega} 4\theta^{*} \dd S_x =0, \\
&\frac{\dd }{\dd t} \int_{\Omega} w_{i}^{*} \norm{A_{i}x}^2 \dd x  + \l   \int_{\pt \Omega} w_{i}^{*} \norm{A_{i}x}^2 \dd S_x = 0,\;\;i=1,2,3.
\end{split}
\end{equation}
Owing to the initial conditions in \eqref{th0-initial} and the fact $\l=0$ in \eqref{a-e-limit-0}, the ODEs in
\eqref{theta0-equation} admit trivial solutions $\th^{*}(t) \equiv 0$ and $w_{i}^{*}(t) \equiv 0$ for all $i=1,2,3$. This proves the claim \eqref{th00-u00}.

%Elementary computation of $\pt_t \tilde{\mu}$ and $\pt_t \sqrt{\tilde{\mu}}$ analogous to %\eqref{pt-sqrt-mu},
%\begin{equation*}
%\begin{split}
%&\pt_t \tilde{\mu} = \Big( \frac{\norm{v-\c}^2-3(1+\th^\e)}{2}\Big)\frac{\pt_t %\th}{(1+\th)^2} \tilde{\mu}^{\e} + \frac{(v-\c)\cdot \pt_t u}{1+\th^\e} \tilde{\mu}^{\e} + %\pt_t \rho \tilde{\mu}^{\e},\\
%&\pt_t \sqrt{\tilde{\mu}} = \Big( \frac{\norm{v-\c}^2-3(1+\th^\e)}{2}\Big)\frac{\pt_t %\th^\e}{2(1+\th^\e)^2} \sqrt{\tilde{\mu}^{\e}} + \frac{(v-\c)\cdot \pt_t u}{2(1+\th^\e)} %\sqrt{\tilde{\mu}^{\e}} + \frac{\pt_t\rho}{2}  \sqrt{\tilde{\mu}^{\e}},
%\end{split}\end{equation*}
% together with \eqref{th-pt-lim}, \eqref{u-lim}, \eqref{tildemu - mu} and \eqref{th00-u00}, %imply
%\begin{equation} \label{remainder- goes to zero}
%\e^{-1} \o  \pt_t \tilde{\mu} \to 0, \;\;\;\;\; \o  \pt_t \sqrt{\tilde{\mu}} \to 0  \;\; %\text{ strongly in } \; L^\infty\big(\mathbb{R}^{+},L^{1}\cap L^{\infty}(\Omega \times %\R^3)\big).
%\end{equation}

We now prove the strong convergence stated in  \eqref{tilde-f-strong-convergence}--\eqref{momentums-strong-convergence}.
The uniform bound on $\normmm{\tilde{f}}_{2}(t)$ from Step 1, combined with \eqref{Gamma-L2-estimate}, \eqref{rho th u - smallness} from Proposition \ref{th tht w wt - ODE} and Lemma \ref{v polynomial - L2 estimate}, implies
\begin{equation*}%\label{L2-bound-average}
\begin{split}
\pt_t\tilde{f}, \;   \e^{-1}\tilde{\nu}^{-\frac{1}{2}} \tilde{L}\tilde{f}, \;   \tilde{\nu}^{-\frac{1}{2}} \tilde{\G}(\tilde{f}, \tilde{f}), \;
\frac{\pt_t \tilde{\mu}}{\sqrt{\tilde{\mu}}}, \;
\e\frac{\pt_t \sqrt{\tilde{\mu}}}{\sqrt{\tilde{\mu}}} \tilde{f} \in L^2(\mathbb{R}^{+}\times\O\times \mathbb{R}^3).
\end{split}
\end{equation*}
Arguing as in the proof of \eqref{tilde-f-epsilon-strong-convergence} and using velocity averaging lemma, we obtain
\begin{equation*}%\label{tilde-f-epsilon-strong-convergence-2}
\tilde{f} \to  f^* \quad \hbox{ strongly in } L^2\big(\mathbb{R}^{+};L^2(\O\times \mathbb{R}^3)\big)  \text{ as } \e\to 0.
\end{equation*}
Combined this with \eqref{tildemu - mu} implies
\begin{equation*}%\label{Pf-strong-convergence}
\begin{split}
&\int_{\mathbb{R}^3} \tilde{f} \sqrt{\tilde{\mu}}\Big[1,v, \frac{|v|^2-3}{2}\Big] \dd v  \to  \big( \varrho_{f^{*}},  u_{f^{*}}, \vartheta_{f^{*}}\big)
\;\;  \hbox{ strongly in }L^2\big(\mathbb{R}^{+};L^2(\O)\big) \text { as } \e\to  0.
\end{split}
\end{equation*}
In view of \eqref{f-es-L2-lim} and \eqref{tildemu - mu},
the strong convergence \eqref{tilde-f-strong-convergence}--\eqref{momentums-strong-convergence} follow
readily.

Finally, the convergence of \eqref{tildef - Boltzmanneq ch1} to the fluid system \eqref{INSF-unst} can be treated analogously to the case $\e\lesssim \a \leq  1$. We omit the details for brevity.
\medskip

\noindent\textbf{Step 3. Derivation of the perfect Navier slip boundary  \eqref{Navier-bdy-unst-perfect}.}

Define  the weighted boundary average
\begin{equation*}
  \begin{split}
         \langle g\rangle^{\tilde{\mu}}_{\partial\Omega}:=&\sqrt{2\pi}\int_{v\cdot n>0}g\big|_{\partial\Omega}\sqrt{\tilde{\mu}}[n\cdot v]\dd v.
  \end{split}
\end{equation*}
%By \eqref{boundary-converge-zero} and the fact that $\langle {f}\rangle_{\partial\Omega}$ %is independent of $v$, we have
%\begin{equation}\label{ges-bd-con}
%\begin{split}
%\langle \tilde{f}\rangle^{\tilde{\mu}}_{\partial\Omega}\to  \langle %f^{*}\rangle_{\pt\O}\;\;
%\hbox{ weakly in } L^2(\mathbb{R}^+\times\pt\O)  \text{ as } \e\to 0.
%\end{split}
%\end{equation}
Combining this with \eqref{boundary-converge-zero} and \eqref{tildemu - mu}, we obtain
\begin{equation}\label{ges-<ges>-cov-0}
\begin{split}
\sqrt{\frac{\a}{\e}} \tilde{\nu}^{-\frac{1}{2}} \Big (  \tilde{f}\big|_{\pt\O}-\frac{\mu}{\sqrt{\tilde{\mu}}} \langle \tilde{f}\rangle^{\tilde{\mu}}_{\pt\O}  \Big)
  \to
0 \hbox{ weakly in } L^2(\mathbb{R}^{+}\times \dd \g)  \text { as } \e\to  0.
\end{split}
\end{equation}
By \eqref{r-expression}, \eqref{th-lim}, \eqref{u-lim}, \eqref{tildemu - mu} and \eqref{th00-u00}, we have
\begin{equation}\label{bdy-remainder- goes to zero}
\begin{split}
r=\frac{\mathscr{P}\tilde{\mu} - \tilde{\mu}} {\e\sqrt{\tilde{\mu}}}
%&=\,  \frac{\rho (1+\th)^{1/2} \mu - \tilde{\mu}} {\e\sqrt{\tilde{\mu}}}
&= \sqrt{\tilde{\mu}}\Big[\Big(2-\frac{\norm{v}^2}{2}\Big)\frac{\theta}{\e} - v\cdot \frac{\c}{\e} + \e O\Big( \norm{\frac{\theta}{\e}}^2, \norm{\frac{\c}{\e}}^2\Big)p(v)\Big]
 \to 0
\end{split}\end{equation}
 strongly in $ L^\infty\big(\mathbb{R}^{+},L^{1}\cap L^{\infty}(\Omega \times \R^3)\big)$ as $\e\to  0$.

 Following the same pattern as in Section 3.4, we now derive the weak formulations \eqref{theta-star-weak-equation} and \eqref{u-star-weak-equation} with $\l=0$ for the INSF system subject to the perfect Navier slip boundary \eqref{Navier-bdy-unst-perfect}. The details are omitted here for brevity.

This completes the proof of Theorem \ref{main-th-2}.
\end{proof}

\bigskip

%%%%%%%%%%%%%%%%%%%%%%%%%%%%%%%%%%
%%%%%%%%%%%%%%%%%%%%%%%%%%%%%%%%%%

\appendix

\section{$L^2_{t}L^{3}_{x}$ Estimate}\
\medskip

The main goal of this section is to establish the following $L^{2}_{t} L^{3}_{x}$ estimate.
\medskip

\begin{proposition} \label{f - L2L3 estimate} \
Let $g \in L^{2} (\mathbb{R}^{+} \times \Omega \times \mathbb{R}^{3})$ and $r\in L^{2} (  \mathbb{R}^{+}\times \g_{-})$, and let $\tilde{f}, f \in L^\infty(\mathbb{R}^{+}\times \Omega \times \R^3)\cap L^2(\mathbb{R}^{+}\times \Omega \times \R^3)$ be distributional solutions of the transport equation with Maxwell boundary condition
\begin{align}
  \e \partial_t \tilde{f} +  v \cdot \nabla_x \tilde{f} = g  \ \ \ \    & \text{ in } \mathbb{R}^{+}\times \O\times \mathbb{R}^3, \label{f-linear-eq-unst}\\
\tilde{f} |_{ \g_{-} } = ( 1 - \a ) \mathscr{R} \tilde{f} + \a \mathscr{P}_{\g} \tilde{f} + \a r  \ \ \ \    & \text{ on } \mathbb{R}^{+}\times \pt\O\times \mathbb{R}^3. \label{f-linear-bdy-unst}
\end{align}
Denote by $\bar{a},  \bar{b}_{i}$, $\bar{c}$ the coefficients of $ \tilde{\P} \tilde{f}$ with respect to the basis $\{ \bar{\chi}_i \}$, and by ${a},  b_{i}$, $c$ the coefficients of $ \P f$ with respect to the basis $\{ {\chi}_i \}$,

\noindent $(1)$ For $0\leq \a \ll \e$, under the a priori assumption \eqref{theta-u-smallness-assumption},
there exist $\mathbf{S}_{1}\tilde{f}(t,x)$ and $\mathbf{S}_{2}\tilde{f}(t,x)$
such that
\begin{equation}\label{S1S2-esti}
\begin{split}
&\norm{ \bar{a} (t,x) } + \norm{ \bar{b}(t,x) } + \norm{ \bar{c}(t,x) } \le  \mathbf{S}_{1} \tilde{f}(t,x) + \mathbf{S}_{2} \tilde{f}(t,x),\\
%--------------------
&\normm{ \mathbf{S}_{2}\tilde{f}}_{L^{2}_{t,x}} \lesssim    \normm{ \ipt f}_{L^{2}_{t,x,v}},\\
&\normm{ \mathbf{S}_{1} \tilde{f} }_{L^{2}_{t} L^{3}_{x} }
\lesssim    \normm{  \tilde{\nu}^{-\frac{1}{2}} g }_{L^{2}_{t,x,v}}+ \normm{\tilde{\nu}^{\frac{1}{2}} \tilde{f}}_{L^{2}_{t,x,v}} +   \a \norm{\tilde{f}}_{L^2_t L^2_{\g_+}}
+  \a \norm{r}_{L^2_t L^2_{\g_-}}
   +  \normm{ \tilde{f}_0 }_{L^{2}_{x,v}} +
\normm{v\cdot \nabla_x \tilde{f}_0 }_{L^{2}_{x,v}}.
 \end{split}
\end{equation}

\noindent $(2)$ For $\e \lesssim \a \leq 1$,  there exist $\mathbf{S}_{1}f(t,x)$ and $\mathbf{S}_{2}f(t,x)$ such that
\begin{equation}\label{S1S2-esti-0}
\begin{split}
&\norm{ {a} (t,x) } + \norm{ {b}(t,x) } + \norm{ {c}(t,x) } \le  \mathbf{S}_{1} f(t,x) + \mathbf{S}_{2} f(t,x),\\
%----------------
&\normm{ \mathbf{S}_{2} f}_{L^{2}_{t,x}} \lesssim    \normm{ \ip f}_{L^{2}_{t,x,v}},\\
&\normm{ \mathbf{S}_{1} f }_{L^{2}_{t} L^{3}_{x} }
\lesssim    \normm{  {\nu}^{-\frac{1}{2}} g }_{L^{2}_{t,x,v}}+ \normm{{\nu}^{\frac{1}{2}} f}_{L^{2}_{t,x,v}} +   \a \norm{f}_{L^2_t L^2_{\g_+}}
+  \a \norm{r}_{L^2_t L^2_{\g_-}}
   +  \normm{ f_0 }_{L^{2}_{x,v}} +
\normm{v\cdot \nabla_x f_0 }_{L^{2}_{x,v}}.
 \end{split}
\end{equation}
\end{proposition}

\begin{proof}[\textbf{Proof.}]\
The argument follows that of Proposition 3.4 in \cite{Esposito2017}. We provide details only for case (1), as case (2) is analogous.

To isolate the interior and non-grazing part of $\tilde{f}$ near the boundary, we introduce a truncation $\tilde{f}_{\d}$. For $(t,x,v) \in \R \times \bar{\Omega} \times \R^3$ and a small parameter $0 < \d \ll 1$, define
\begin{equation} \label{fdelta- Defintion}
\begin{split}
\tilde{f}_{\d} (t,x,v) := &\big[ 1- \chi \big(\frac{n(x) \cdot v}{\d} \big) \chi \big( \frac{\xi(x)}{\d} \big)  \big]  \big[ 1- \chi \big(\frac{|v|}{2\d} \big)  \big]
 \chi(\d \norm{v})  \big[ \1_{t \in [0, \infty)} \tilde{f}(t,x,v) + \1_{t \in (-\infty, 0]}\tilde{f}_{0}(x,v)\big].
\end{split}
\end{equation}
Here the cutoff function $\chi \in C^{\infty}_{c}(\mathbb{R})$ satisfies
\begin{equation*}
\begin{split}
  0 \leq \chi \leq 1, \ \ \
\chi^{\prime} (x) \geq -4 \times \mathbf{1}_{\frac{1}{2} \leq |x| \leq 1} \ \ \
 \text{and} \ \   \chi(x)=
    \begin{cases}
    1 & \text{if } |x| \leq \frac{1}{2},\\
      0& \text{if } |x| \geq 1 .
    \end{cases}
 \end{split}
 \end{equation*}
Consequently, $\tilde{f}_{\delta}(t,x,v)$ vanishes on the near-grazing set:
\begin{equation}\label{f-delta-support}
  \tilde{f}_{\delta}(t,x,v) = 0  \   \text { for }   \ (x,v)\in \g\backslash \g_{\pm}^{\delta},
\end{equation}
with the non-grazing sets $\g_{\pm}^{\delta}$ defined in \eqref{non-grazing-set}.
Moreover, the following estimates hold:
\begin{equation} \label{f-delta-f}
\begin{split}
& \| \tilde{f}_{\d} \|_{L^{2} (\mathbb{R} \times \O \times \mathbb{R}^{3})}
  \lesssim  \| \tilde{f} \|_{L^{2} (\mathbb{R}_{+} \times \O \times \mathbb{R}^{3})}  + \| \tilde{f}_{0}\|_{L^{2} (\O \times \mathbb{R}^{3})}, \\
& \| \tilde{f}_{\d} \|_{L^{2} ( \mathbb{R} \times \g)}
  \lesssim   \| \tilde{f}{\bf 1}_{\g^{\d}_{\pm}} \|_{L^{2} ( \mathbb{R}_{+} \times \g)}  + \| \tilde{f}_{0} {\bf 1}_{\g^{\d}_{\pm}} \|_{L^{2} (\g)}.
\end{split}
\end{equation}

Under the a priori assumption \eqref{theta-u-smallness-assumption},
there exists a constant $T_M>0$ such that
\begin{equation*}
T_M < T(t)=1+\th(t) <2 T_M \;\;  \text{ for all } \; t\geq 0.
\end{equation*}
Then, for some constants $C_1, C_2>0$ and $p\in (\frac{1}{2}, 1)$, the global Maxwellian
\begin{equation*}
\mu_{M} :=  \frac{1}{(2\pi T_M)^{3/2}} \exp \Big( -\frac{\norm{v}^2}{2T_M}\Big)
\end{equation*}
%\begin{equation*}
%\mu_{\max} := \sum_{\text{All combinations of }\pm} \frac{(1\pm \e)}{(2\pi(1\pm %\e))^{3/2}} \exp \left( -\frac{\sum_{i=1}^{3}\norm{v_i \pm \e}^2}{2(1\pm \e)}\right),
%\end{equation*}
satisfies
\begin{equation}
C_1\mu_{M} \leq  \tilde{\mu} \lesssim C_2\mu_{M}^p,
\end{equation}
as shown in \cite{Guo-Jang-Jiang2009}. Consequently,
\begin{equation}
\norm{\bar{\chi}_{i}(v)} \lesssim  \inn{v}^2 \mu_{M}^{\frac{p}{2}}  \;\; \text{ and } \;\;
\norm{\bar{a}}, \norm{\bar{b}}, \norm{\bar{c}}\lesssim \norm{\int_{\R^3} {\inn{v}^2 \mu_{M}^{\frac{p}{2}}  }  f \dd v}.
\end{equation}

For each $i\in \{0,1,\cdots,4\}$, the truncation $\tilde{f}_\d$ satisfies
\begin{equation}\label{int-f-delta-decomposition}
\begin{split}
&\int_{\mathbb{R}^{3}} \tilde{f}_{\delta} \bar{\chi}_{i} (v) \mathrm{d} v \\
%--------------------------------------
%=& \int_{\mathbb{R}^{3}} \big[ 1- \chi( \frac{n(x) \cdot v}{\delta} ) \chi(
%\frac{\xi(x)}{\delta}) \big]
%--------
%\big[ 1- \chi \big(\frac{|v|}{2\d} \big)  \big]
 %--------
% \chi( \delta |v|) \\
%--------------------------------------
% &\times\Big\{ \mathbf{1}_{t \geq 0}
%\tilde{f}(t,x,v) + \mathbf{1}_{t \leq 0} \chi(t) f_{0} (x,v) \Big\}\bar{\chi}_{i} (v)
%\mathrm{d} v \\
%--------------------------------------
=& \mathbf{1}_{t \geq 0} \int_{\mathbb{R}^{3}} \big[ 1- \chi( \frac{n(x)
\cdot v}{\delta} ) \chi( \frac{\xi(x)}{\delta}) \big]
%--------
 \big[ 1- \chi \big(\frac{|v|}{2\d} \big)  \big]
 %--------
\chi( \delta |v|)%
%--------------------------------------
 \Big\{ \sum_{j=0}^{4} \bar{a}_{j} \bar{\chi}_{j}(v) + \ipt
\tilde{f}\Big\}\bar{ \chi}_{i}(v) \mathrm{d} v \\
%--------------------------------------
&+ \mathbf{1}_{t \leq 0} \int_{\mathbb{R}^{3}} \big[ 1- \chi( \frac{n(x)
\cdot v}{\delta} ) \chi( \frac{\xi(x)}{\delta}) \big]
 %--------
  \big[ 1- \chi \big(\frac{|v|}{2\d} \big)  \big]
 %--------
 \chi( \delta |v|)
\chi(t) f_{0} \bar{\chi}_{i}(v) \mathrm{d} v \\
%--------------------------------------
=& \mathbf{1}_{t\geq 0} \Big\{ \bar{a}_{i}  + O(\delta) \sum_{j=0}^{4} |\bar{a}_{j}
| \Big\} \\
& +\mathbf{1}_{t\geq 0}  \int_{\mathbb{R}^{3}} \big[ 1- \chi( \frac{n(x) \cdot v}{\delta} ) \chi(
\frac{\xi(x)}{\delta}) \big]
%--------
 \big[ 1- \chi \big(\frac{|v|}{2\d} \big)  \big]
 %--------
 \chi( \delta |v|)  \ipt
\tilde{f}  \bar{\chi}_{i}(v) \mathrm{d} v \\
%--------------------------------------
&+ \mathbf{1} _{ t \leq 0} \chi(t) \int_{\mathbb{R}^{3}} \big[ 1- \chi( \frac{n(x) \cdot v}{\delta} ) \chi(
\frac{\xi(x)}{\delta}) \big]
%--------
  \big[ 1- \chi \big(\frac{|v|}{2\d} \big)  \big]
 %--------
 \chi( \delta |v|)  \tilde{f}_{0}
\bar{\chi}_{i}(v) \mathrm{d} v,
\end{split}
\end{equation}
where temporary notations $\bar{a}_{0}=\bar{a}$, $\bar{a}_{i}=\bar{b}_{i}$ $(i=1,2,3)$ and $\bar{a}_{4}=\bar{c}$ are used (see \eqref{barabc-def}).
% \begin{equation*}
%\begin{split}
%\bar{a}_{0}(t,x):=\bar{a}(t,x), \;\;\;\; \bar{a}_{i}(t,x):=\bar{b}_i(t,x) \text{ for } %i=1,2,3, \;\;\;\; \bar{a}_{4}(t,x):=\bar{c}(t,x),
%\end{split}
%\end{equation*}
%so that $\tilde{\P} \tilde{f}= \sum_{j=0}^{4} \bar{a}_{i}\bar{ \chi}_{i} (v)$.
Therefore,
\begin{equation*}
\begin{split}
\sum_{i=0}^{4}\mathbf{1}_{t \geq 0} |\bar{a}_{i} | \leq& \sum_{i=0}^{4}\Big|%
\int_{\mathbb{R}^{3}} \tilde{f}_{\delta}  \bar{\chi}_{i} (v) \mathrm{d} v\Big| +
\mathbf{1} _{ t \leq 0} \chi(t) \int_{\mathbb{R}^{3}} | \tilde{f}_{0} |
\sum_{i=0}^{4}|\bar{\chi}_{i} (v)| \mathrm{d} v \\
&+ \mathbf{1}_{t \geq 0} \Big\{5 O(\delta) \sum_{j=0}^{4} |\bar{a}_{j}| +
O_{\delta} (1) \int_{\mathbb{R}^{3}} | \ipt \tilde{f} |
\sum_{i=0}^{4} |\bar{\chi}_{i} (v)| \mathrm{d} v \Big\}.
\end{split}
\end{equation*}
Hence, for sufficiently small $\d$, we obtain for each $i=0,1,2,3,4$:
\begin{equation}  \label{bound_a_i}
\begin{split}
|\bar{a}_{i} (t,x)| \leq& \ 10 \int_{\mathbb{R}^{3}} |\tilde{f}_{\delta} | \langle
v\rangle^{2} \mu_{M}^{\frac{p}{2}}  \mathrm{d} v + 10\chi(t)\mathbf{1} _{ t \leq 0}
\int_{\mathbb{R}^{3}} |\tilde{f}_{0} | \langle v \rangle^{2} \mu_{M}^{\frac{p}{2}}
\mathrm{d} v+ 10 \int_{\mathbb{R}^{3}} | \ipt \tilde{f}| \langle
v\rangle^{2} \mu_{M}^{\frac{p}{2}}  \mathrm{d} v.
\end{split}%
\end{equation}

We now focus on the term involving $\tilde{f}_{\delta}$ in \eqref{bound_a_i}. By Lemma 3.6 and Lemma 3.7 in \cite{Esposito2017}, there exists an extension $\overline{\tilde{f}_{\delta}}\in L^{2} (\mathbb{R} \times \Omega \times \mathbb{R}^{3})$ of $\tilde{f}_{\delta}$ such that
%\begin{equation}
%\int_{\mathbb{R}^{3}} |\tilde{f}_{\delta} (t,x,v) | \langle v\rangle^{2} %\sqrt{\mu_{M}^p(v)} \mathrm{d} v \leq \int_{\mathbb{R}^{3}} | %\overline{\tilde{f}_{\delta}}(t,x,v) |
%\langle v\rangle^{2} \sqrt{\mu_{M}^p(v)} \mathrm{d} v.
%\end{equation}
%and satisfies
%\begin{equation}
%\e \pt_t \overline{\tilde{f}_{\delta}}+v\cdot \nabla_x \overline{\tilde{f}_{\delta}} + %\overline{\tilde{f}_{\delta}}=h+
%\overline{\tilde{f}_{\delta}}  \label{bastransp0}
%\end{equation}
%in the sense of distributions, where $h=\sum_{i=1}^4 h_i$ and $h_i$ are
%defined in Lemma \ref{extension_dyn}. Applying Lemma \ref{f- velocity averaging lemma} to %(\ref{bastransp0}) with $f= \overline{\tilde{f}_{\delta}}$, $q= h+ %\overline{\tilde{f}_{\delta}}$ and $\psi=\langle v\rangle^{2}\sqrt{\mu_{M}^p(v)}$, we %obtain
%\begin{equation}  \label{average_fdelta}
%\begin{split}
%\Big\|\int_{\mathbb{R}^{3}}\langle %v\rangle^{2}\sqrt{\mu_{M}^p(v)}\overline{\tilde{f}_{\delta}} \mathrm{%
%d} v\Big\|_{L^{2}_{t} L^{3}_{x}} &\lesssim \| \omega^{-1} h \|_{L^{2}_{t,x,v}} +
%\| \omega^{-1} \overline{\tilde{f}_{\delta}} \|_{L^{2}_{t,x,v}}.
%\end{split}%
%\end{equation}
%From Lemma \ref{extension_dyn}, we have the bound
\begin{equation}
\begin{split}  \label{forcing_fdelta}
 \| \omega^{-1} \overline{\tilde{f}_{\delta}}
\|_{L^{2}_{t,x,v}} \lesssim & \| \omega^{-1} g \|_{L^{2}_{t,x,v}} +\|  \tilde{f } \|_{L^{2}_{t,x,v}}+ \| \tilde{f}_{0} \|_{L^{2}_{x,v}} + \| v\cdot
\nabla_{x} \tilde{f}_{0} \|_{L^{2}_{x,v}}  +|\tilde{f}{\bf 1}_{\g^{\d}_{\pm}} |_{L^{2}_tL^{2}_\gamma}
 + |\tilde{f}_{0}{\bf 1}_{\g^{\d}_{\pm}}|_{L^{2}_\gamma}.
\end{split}%
\end{equation}
Note that the boundary term $|\tilde{f}{\bf 1}_{\g^{\d}_{\pm}} |_{L^{2}_tL^{2}_\gamma}$ arises from the definition
of $\tilde{f}_{\delta}$, \eqref{f-delta-support} and \eqref{f-delta-f}.
%In the following, we further estimate
% $|\tilde{f}_{0}{\bf 1}_{\g^{\d}_{\pm}}|_{L^{2}_\gamma}$  and $|\tilde{f}{\bf %1}_{\g^{\d}_{\pm}} |_{L^{2}_tL^{2}_\gamma}$.

To bound $|\tilde{f}_{0}{\bf 1}_{\g^{\d}_{\pm}}|_{L^{2}_\gamma}$, we apply Ukai's trace Lemma in \cite{Ukai1986} or Lemma 2.3 in \cite{Esposito2017} on $\g^{\d}_{\pm}$, yielding
\begin{equation}\label{f0-non-grazing}
\begin{split}
  |\tilde{f}_{0}{\bf 1}_{\g^{\d}_{\pm}}|_{L^{2}_\gamma}
& \lesssim_{\d} \| \tilde{f}_0^2 \|_{L^{1}_{x,v}} +  \|v \cdot \nabla_x  ( \tilde{f}_0^2 ) \|_{L^{1}_{x,v}}  \lesssim_{\d} \| \tilde{f}_0 \|_{L^{2}_{x,v}}^2  + \|  v \cdot \nabla_x  \tilde{f}_0\|_{L^{2}_{x,v}}^2.
\end{split}
\end{equation}
To estimate $|\tilde{f}{\bf 1}_{\g^{\d}_{+}} |_{L^{2}_tL^{2}_\gamma}$, we apply Lemma 3.2 in \cite{Esposito2017}
on the out-going non-grazing set $\g^{\d}_{+}$:
\begin{equation}\label{f-non-grazing+}
\begin{split}
 \int_0^t \int_{\g_{+}} |\tilde{f} \mathbf{1}_{\g^\delta_+}|^2 \dd \g
& \lesssim_{\d} \e \| \tilde{f}_0^2 \|_{L^{1}_{x,v}} +  \int_0^t \| \tilde{f}^2 \|_{L^{1}_{x,v}} + \int_0^t \| \big ( \e \partial_t  + v \cdot \nabla_x
  \big ) ( \tilde{f}^2 ) \|_{L^{1}_{x,v}}  \\
&\lesssim_{\d}  \e \| \tilde{f}_0 \|_{L^{2}_{x,v}}^2 +  \int_0^t \|  \tilde{f} \|_{L^{2}_{x,v}}^2 + \int_0^t \|\tilde{f}\|_{L^{2}_{x,v}(\tilde{\nu})}^2 +  \int_0^t \| \tilde{\nu}^{-\frac{1}{2}} g \|_{L^{2}_{x,v}}^2.
\end{split}
\end{equation}
For $|\tilde{f}{\bf 1}_{\g^{\d}_{-}} |_{L^{2}_tL^{2}_\gamma}$, where trace lemma does not apply on $\g^{\d}_{-}$, we use boundary condition \eqref{f-linear-bdy-unst} and the change of variable $v\mapsto R_x v$ on $\g^{\d}_{-}$:
\begin{equation}\label{f-non-grazing-}
\begin{split}
  \int_0^t \int_{\g_{-}} | \tilde{f} {\bf 1}_{\g^{\d}_{-}}|^2 \dd \g
%  & = \int_0^t \int_{\g^{\d}_{-}}  \big  [ ( 1 - \a ) \mathcal{L} f + \a P_\g f + \a r \big ]^2  \dd \g  \\
  \lesssim   & \int_0^t \int_{\g^{\d}_{-}}   |\mathscr{R} (\tilde{f}) |^2  \dd \g  + \a^2  \int_0^t \int_{\g^{\d}_{+}} |\mathscr{P}_\g \tilde{f}|^2 \dd \g
     + \a^2  \int_0^t \int_{\g^{\d}_{-}}  |r|^2  \dd \g  \\
%   & =   \int_0^t \int_{\g^{\d}_{+}}   | f |^2  \dd \g   + \a^2  \int_0^t \int_{\g^{\d}_{+}} |P_\g f|^2 \dd \g
%     + \a^2  \int_0^t \int_{\g^{\d}_{-}}  |r|^2  \dd \g  \\
      \lesssim  &  \int_0^t \int_{\g^{\d}_{+}}   | \tilde{f} |^2  \dd \g
      + \a^2  \int_0^t \int_{\g_{+}} |\tilde{f}|^2 \dd \g
     + \a^2  \int_0^t \int_{\g^{\d}_{-}}  |r|^2  \dd \g  \\
\lesssim  & \e \| \tilde{f}_0 \|_{L^{2}_{x,v}}^2  + \int_0^t \|\tilde{f}\|_{L^{2}_{x,v}(\tilde{\nu})}^2 +  \int_0^t \| \tilde{\nu}^{-\frac{1}{2}} g \|_{L^{2}_{x,v}}^2
+ \a^2\int_0^t \big [  |\tilde{f}|_{L^2_{\g_{+}}}^2 + |r|_{L^2_{\g_{-}}}^2  \big ],
\end{split}
\end{equation}
where we used \eqref{f-non-grazing+} in the last inequality.

Finally, we define
\begin{equation}  \label{def_S2}
\begin{split}
\mathbf{S}_{1} \tilde{f}(t,x) &: = \int_{\mathbb{R}^{3}} | \overline{\tilde{f}_{\delta}} | \langle v\rangle^{2} \mu_{M}^{\frac{p}{2}} \mathrm{d} v, \quad
\mathbf{S}_{2} \tilde{f}(t,x) : = 4 \int_{\mathbb{R}^{3}} | \ipt \tilde{f}| \langle v\rangle^{2} \mu_{M}^{\frac{p}{2}} \mathrm{d} v.
\end{split}%
\end{equation}
Combining \eqref{bound_a_i}--\eqref{def_S2}, we obtain \eqref{S1S2-esti}. This completes the proof of Proposition \ref{f - L2L3 estimate}.
\end{proof}
\bigskip

\section{Uniqueness of Weak Solutions to INSF}

In the following, we establish the uniqueness of weak solutions to the INSF system in the setting of Theorem \ref{main-th-1} and Theorem  \ref{main-th-2}.

\begin{lemma}[Uniqueness of weak solutions to the INSF system]\label{A-lemma-NS} \
Under the assumptions of Theorem 1.1 and Theorem  \ref{main-th-2}, the weak solution $(u,\vartheta)$ to the INSF system \eqref{INSF-unst} --- subject to either the Dirichlet boundary condition \eqref{Diri-bdy-unst} or the Navier boundary condition \eqref{Navier-bdy-unst} (which reduces \eqref{Navier-bdy-unst-perfect} when $\l=0$) with initial data $(u_0, \vartheta_0)\in \mathbb{H}_{u}\times\mathbb{H}_{\vartheta}$ (defined in \eqref{initial-space-Navier}) --- is unique.
\end{lemma}

\begin{proof}[\textbf{Proof}.]\
We prove uniqueness of weak solution $(u,\vartheta)$ to the INSF system \eqref{INSF-unst} only for the Navier boundary condition \eqref{Navier-bdy-unst}, as the proof for the Dirichlet case \eqref{Diri-bdy-unst} follows analogously and is simpler.

% Recall the INSF system with the Navier boundary condition
% \begin{equation}\label{A-INSF-unst}
% \begin{split}
% \pt_t {u} + {u} \cdot\nabla_x {u}  +\nabla_x {p}
% = \s \Delta {u}, \ \ \  \nabla_x\cdot {u} = 0 \ \ \ \    & \text{ in } % % \mathbb{R}^{+}\times \O, \\
% \pt_t {\vartheta} +  u\cdot\nabla_x {\vartheta}
% = \kappa  \Delta_x {\vartheta} \ \ \ \
%   &\text{ in } \mathbb{R}^{+}\times \O, \\
%---------
%  \left [ \sigma\big(\nabla_x u + {(\nabla_x u)}^{\mathrm{T}}\big) \cdot n
%   +\l u \right ]^{\mathrm{tan}} = 0, \ \ \  u\cdot n =0 \ \ \ \
%     &\text{ on }\mathbb{R}^{+}\times\pt\O, \\
%  \kappa \partial_n \vartheta + \frac{4}{5}\l \vartheta = 0 \ \ \ \  &\text{ on % }\mathbb{R}^{+}\times\pt\O.\\
%-------------
% {u}|_{t=0}={u}_0, \ \ \   {\vartheta}|_{t=0}={\vartheta}_0  \ \ \ \     &\text{ on }  % \O.
% \end{split}
% \end{equation}
As a limit point of solutions to the Boltzmann equation when $\displaystyle \lim_{\e\rightarrow 0}\frac{\a}{\e}= \sqrt{2\pi}\l \in [0,\infty)$, the pair $(u,\vartheta)$ inherits the smallness of $\normmm{f}_{1}$ or $\normmm{\tilde{f}}_{2}$.
More precisely, from the uniform bound \eqref{0-uniform-bound} or \eqref{uniform-bound-tilde} and the uniqueness of distribution limit, up to a subsequence,
\begin{equation*}\label{A-u1-u2-L6-small}
\begin{split}
\tilde{\P}\tilde{f}, \P {f}  \to  \P f^{*}
\text{~~weakly}\!-\!* \text{ in } L^\infty\left( \mathbb{R}^+; L^6(\O\times \mathbb{R}^3)\right).
\end{split}
\end{equation*}
By the lower semi-continuity of the norm under weak$-*$ convergence, the limit  $(u,\vartheta)$ inherits the smallness in $L^\infty_{t}L^6_{x}$:
\begin{equation}\label{A-u-theta-L6-small}
\begin{split}
&\normm{u}_{L^\infty_{t}L^6_{x}}
\lesssim  \normm{\P f^{*}}_{L^\infty_{t}L^6_{x,v}}
\lesssim  \normm{\P f}_{L^\infty_{t}L^6_{x,v}}
\lesssim  \normmm{f}_{1}\ll 1 \text{ when using norm } \normmm{\cdot}_{1},\\
&\normm{u}_{L^\infty_{t}L^6_{x}}
\lesssim  \normm{\P f^{*}}_{L^\infty_{t}L^6_{x,v}}
\lesssim  \normm{\tilde{\P}\tilde{f}}_{L^\infty_{t}L^6_{x,v}}
\lesssim  \normmm{\tilde{f}}_{2}\ll 1 \text{ when using norm } \normmm{\cdot}_{2}.
\end{split}
\end{equation}

For uniqueness, let $(u_1, \vartheta_1)$ and $(u_2, \vartheta_2)$ be two solutions of \eqref{INSF-unst} and  \eqref{Navier-bdy-unst} (which reduces \eqref{Navier-bdy-unst-perfect} when $\l=0$)  with the same initial data $(u_0, \vartheta_0)\in \mathbb{H}_{u}\times\mathbb{H}_{\vartheta}$.
Then it follows from \eqref{A-u-theta-L6-small} that
\begin{equation}\label{u1-u2-L6-small}
\begin{split}
\normm{u_1}_{L^\infty_{t}L^6_{x}}\ll 1, \quad \normm{u_2}_{L^\infty_{t}L^6_{x}}\ll 1.
\end{split}
\end{equation}
Write $w=u_1-u_2, \chi=\vartheta_1-\vartheta_2$. Then $(w, \chi)$ satisfies
\begin{equation}\label{A-INSF-unst-w-chi}
\begin{split}
\pt_t w + u_1 \cdot\nabla_{x} w+ w \cdot\nabla_{x} u_2 +\nabla_{x} (p_1-p_2)
   = \s \Delta_{x} w, \ \ \  \nabla_x\cdot {w} = 0 \ \ \ \    & \text{ in } \mathbb{R}^{+}\times \O, \\
\pt_t {\chi} +  u_1\cdot\nabla_x {\chi} + w\cdot \nabla_x\vartheta_2
= \kappa  \Delta {\chi}  \ \ \ \
  &\text{ in } \mathbb{R}^{+}\times \O, \\
{w}|_{t=0}=0, \ \ \   {\chi}|_{t=0}=0  \ \ \ \     &\text{ on }  \O,\\
%---------------
 \Big [ \sigma\big(\nabla_x w + {(\nabla_x w)}^{\mathrm{T}}\big) \cdot n
  +\l w \Big ]^{\mathrm{tan}} = 0, \ \ \  w\cdot n =0 \ \ \ \
    &\text{ on }\mathbb{R}^{+}\times\pt\O, \\
 \kappa \partial_n \chi + \frac{4}{5}\l \chi = 0 \ \ \ \  &\text{ on }\mathbb{R}^{+}\times\pt\O.
\end{split}
\end{equation}
Standard $L^2$ energy estimate on \eqref{A-INSF-unst-w-chi}
 leads to the energy equality
 \begin{equation}\label{A-w-L2-energy}
\begin{split}
&\frac{1}{2} \|w(t)\|^2_{L^2(\O)} + 2\sigma \int_{0}^{t}  \|\nabla_x^{\text{s}} w\|^2_{L^2(\O)} + \lambda \int_{0}^{t} |w_\tau|^2_{L^2(\pt\O)}\\
= & - \int_0^t \int_{\O}  (u_1\cdot \nabla_{x} w)\cdot  w\dd x \dd s-\int_0^t \int_{\O} (w\cdot \nabla_{x} u_2)\cdot  w \dd x \dd s,
\end{split}
\end{equation}
where $w_\tau$ denotes the tangential component of $w$ on $\pt\O$ (in fact, $w_\tau=w$ because $w\cdot n|_{\pt\O}=0$). Here we used $\nabla_x\cdot {w}=0$, $\Delta_x w= 2 \text{div}(\nabla_x^{\text{s}} w)- \text{grad} (\nabla_x\cdot  w) $ and the Navier boundary condition in \eqref{A-INSF-unst-w-chi}.
The first integral on the right-hand side vanishes because $\nabla_{x} \cdot u_1=0$
and $n\cdot u_1|_{\pt\O}=0$.
Using $\nabla_{x} \cdot w=0$, $n\cdot w|_{\pt\O}=0$ and integrating by parts, we have
\begin{equation}\label{A-term-1}
\begin{split}
 \Big | \int_0^t \int_{\O} (w\cdot \nabla_{x} u_2)\cdot  w \dd x \dd s\Big |
    =&\Big | - \int_0^t \int_{\O}   (w\cdot \nabla_{x} w)\cdot u_2 \dd x \dd s\Big |\\
%----------
    \lesssim &\|u_2\|_{L^{\infty}_t L^6_x}   \|w\|_{L^2_t L^3_x}\|\nabla_{x} w\|_{L^2_t L^2_x} \\
%----------
    \lesssim & \|u_2\|_{L^{\infty}_t L^6_x}    \Big \| \|w\|_{L^2_x}^{\frac{1}{2}}   \|w\|_{H^1_x}^{\frac{1}{2}} \Big\|_{L^2_t}     \|\nabla_{x} w\|_{L^2_t L^2_x}\\
%----------
    \lesssim & \|u_2\|_{L^{\infty}_t L^6_x}\big(  \|w\|_{L^2_t L^2_x}^2   +  \|\nabla_{x} w\|_{L^2_t L^2_x}^2\big),
\end{split}
\end{equation}
where we used the Gagliardo-Nirenberg inequality. Substituting \eqref{A-term-1} into \eqref{A-w-L2-energy} yields
\begin{equation}\label{A-w-L2-energy-inequality}
\begin{split}
&\frac{1}{2} \|w(t)\|^2_{L^2(\O)} + 2\sigma \int_{0}^{t}  \|\nabla_x^{\text{s}} w\|^2_{L^2(\O)} + \lambda \int_{0}^{t} |w_\tau|^2_{L^2(\pt\O)}
\lesssim  \|u_2\|_{L^{\infty}_t L^6_x}\int_{0}^{t} \|w\|_{H^1(\O)}^2.
\end{split}
\end{equation}
Similarly,
\begin{equation}\label{A-chi-L2-energy-inequality}
\begin{split}
&\frac{1}{2} \|\chi(t)\|^2_{L^2(\O)} + \kappa \int_{0}^{t}  \|\nabla\chi\|^2_{L^2(\O)} + \frac{4}{5}\lambda \int_{0}^{t} |\chi|^2_{L^2(\pt\O)}
\lesssim   \|\vartheta_2\|_{L^{\infty}_t L^6_x}\int_{0}^{t} \big( \|w\|_{H^1(\O)}^2+  \|\nabla \chi\|_{L^2(\O)}^2\big).
\end{split}
\end{equation}

Because the coefficient  $\l$ influences the boundary dissipation, we treat the cases  $\l> 0$ and $\l=0$ separately. The geometry of $\O$  also affects the solution when $\l=0$ (perfect Navier slip boundary).
\medskip

\noindent\textbf{Step 1. Case $\l> 0$.}

\noindent\textbf{Step 1.1. Estimate for $u$.}

To close \eqref{A-w-L2-energy-inequality}, we use the following Korn-type inequality (see Proposition 3.13 in \cite{Tapia-Amrouche-Conca-Ghosh-2021}): for any $g\in H^1(\O)$ with $g\cdot n|_{\pt\O}=0$,
\begin{equation}\label{Korn-type-inequality}
\begin{split}
&\| g\|_{H^1(\O)} \simeq \left\{
                           \begin{array}{ll}
                             \|\nabla_x^{\text{s}}  g\|_{L^2(\O)}, & \text{ if } \; \O \text{ is non-axisymmetric}; \\
                             \|\nabla_x^{\text{s}}  g\|_{L^2(\O)}+ |g_{\tau}|_{L^2(\pt\O)}, & \text{ if } \; \O \text{ is axisymmetric or spherical}.
                           \end{array}
                         \right.
\end{split}
\end{equation}
From \eqref{u1-u2-L6-small} and \eqref{Korn-type-inequality}, we obtain
\begin{equation}\label{w-L2-energy-close-non0}
\begin{split}
&\frac{1}{2} \|w(t)\|^2_{L^2(\O)} + \sigma \int_{0}^{t}  \|\nabla_x^{\text{s}} w\|^2_{L^2(\O)} + \frac{\l}{2}\int_{0}^{t} |w_\tau|^2_{L^2(\pt\O)}\leq 0\;\; \;\; \text{ for all } t\geq 0.
\end{split}
\end{equation}
regardless of whether $\O$ is axisymmetric, spherical or non-axisymmetric.
Together with \eqref{Korn-type-inequality}, this gives $w\equiv 0$; hence $u$ is unique.
\medskip

\noindent\textbf{Step 1.2. Estimate for $\vartheta$.}

Recall the Friedrich inequality
\begin{equation}\label{Friedrich-inequality}
\begin{split}
\|\chi\|_{H^1(\O)} \lesssim   \|\nabla_x\chi\|_{L^2(\O)}+ |\chi|_{L^2(\pt\O)}.
\end{split}
\end{equation}
Using \eqref{u1-u2-L6-small}, \eqref{Korn-type-inequality}--\eqref{Friedrich-inequality}, we can close the energy equality \eqref{A-chi-L2-energy-inequality}
% as
%\begin{equation}\label{A-chi-L2-energy-close}
%\begin{split}
%& \|w(t)\|^2_{L^2(\O)} + \int_{0}^{t}  \|w\|_{H^1(\O)}^2 +\|\chi(t)\|^2_{L^2(\O)}+ %\int_{0}^{t}  \|\chi\|^2_{H^1(\O)}
%\leq 0
%\end{split}
%\end{equation}
and deduce uniqueness of $\vartheta$.
\medskip

\noindent\textbf{Step 2. Case $\l= 0$.}

For $\l= 0$, if $\O$ is axisymmetric or spherical, the incompressible Navier-Stokes equation with perfect Navier slip boundary  admits nontrivial kernels $u=R(x)$ (see \cite{Amrouche-Rejaiba-2014}), where $R(x)=Ax$ is a basis element of $\mathcal{R}_\O$ defined in \eqref{RO}. The heat equation with homogeneous Neumann boundary also has constants as kernels.  Therefore, to ensure uniqueness of $(u,\vartheta)$ when $\l= 0$, we must require $(u_0, \vartheta_0)\in \mathbb{H}_{u}\times\mathbb{H}_{\vartheta}$ (cf. \eqref{initial-space-Navier}).
\medskip

\noindent\textbf{Step 2.1. Estimate for $u$.}

\noindent\textbf{Step 2.1.1. $\O$ axisymmetric or spherical.}

First note that the Navier-Stokes equation with perfect Navier slip boundary $\l=0$
satisfies conservation law of angular momentum:
  \begin{equation}\label{conservation-law-angular-momentum-u}
\begin{split}
 \pt_t \int_{\O}u(t,x)\cdot R(x)\dd x=0 \;\;  \text{ for all } R(x) \in\mathcal{R}_{\O}\text{ and all } t>0.
\end{split}
\end{equation}
Indeed, for any $R \in\mathcal{R}_{\O}$,
  \begin{equation}\label{Delta-u-R}
\begin{split}
  \s \int_{\O}\Delta_x u\cdot R\dd x
  = & 2\s  \int_{\O}\text{div}(\nabla_x^{\text{s}} u) \cdot R\dd x\;\;  \;\;(\text{by } \Delta_x u= 2 \text{div}(\nabla_x^{\text{s}} u)- \text{grad} (\nabla_x\cdot  u) ) \\
  = & \s \int_{\O}\pt_i(\pt_iu_k+\pt_k u_i)R_k\dd x \\
   = & \s \int_{\O}\pt_i[(\pt_iu_k+\pt_k u_i)R_k]\dd x-\s \int_{\O}(\pt_iu_k+\pt_k u_i)\pt_iR_k\dd x \\
    = & \s \int_{\pt\O}\underbrace{ n_i(\pt_iu_k+\pt_k u_i)}R_k\dd S_x-2\s \int_{\O}\nabla_x^{\text{s}} u: \nabla_x^{\text{s}} R \dd x \;\;  \;\;(\text{by } \nabla_x^{\text{s}} R=0) \\
    = & \s \int_{\pt\O}\underbrace{[(\pt_iu_j+\pt_j u_i)n_in_j]n_k } R_k\dd S_x \;\;\;\;
    (\text{by } n\cdot R|_{\pt \O}=0) \\
    =&0,
\end{split}
\end{equation}
where for the under braced term we used
\begin{equation}
\begin{split}
0=&\Big  [ \sigma\big(\nabla_x u + {(\nabla_x u)}^{\mathrm{T}}\big) \cdot n
  \Big  ]^{\mathrm{tan}} \\
  %--------------
  =& \sigma\big(\nabla_x u + {(\nabla_x u)}^{\mathrm{T}}\big) \cdot n
   - \sigma \Big [ n\cdot \big(\nabla_x u + {(\nabla_x u)}^{\mathrm{T}}\big) \cdot n
  \Big ] n\\
   %--------------
  =&\sigma   \big(\nabla_x u + {(\nabla_x u)}^{\mathrm{T}}\big) \cdot n -  \sigma \Big[\big(\nabla_x u + {(\nabla_x u)}^{\mathrm{T}}\big): \big ( n\otimes n \big)  \Big ] n .
\end{split}\end{equation}
Moreover, for the nonlinear term and pressure term
  \begin{equation}\label{nabla-p-R}
\begin{split}
  \int_{\O} (u\cdot \nabla_x) u \cdot R\dd x
   = & \int_{\O}u_i\pt_i(u_j R_j)\dd x -\int_{\O}u_i u_j \pt_i R_j\dd x \\
   = & \int_{\O}\pt_i(u_iu_j R_j)\dd x - \int_{\O}\pt_i u_i u_j R_j\dd x  -\int_{\O}u\otimes u: \nabla_x^{\text{s}} R\dd x \\
   = & \int_{\pt\O}n_iu_iu_j R_j\dd S_x  =0,\\
%---------------
  \int_{\O} \nabla_x p \cdot R\dd x
   = & \int_{\O} \pt_i(p R_i)\dd x - \int_{\O} p \nabla_x\cdot R  \dd x
   =  \int_{\pt\O}p n\cdot R\dd S_x
    =0.\end{split}
\end{equation}
Here we used the facts $\nabla_x^{\text{s}} R=0$,  $\nabla_x\cdot R=0$ and $n\cdot R|_{\pt\O}=0=n\cdot u|_{\pt\O}$. Combining \eqref{Delta-u-R}--\eqref{nabla-p-R}, we prove the claim \eqref{conservation-law-angular-momentum-u}.

 Therefore, if $u_0\in \mathbb{H}_{u}$, then $u\in \mathbb{H}_{u}$ for all $t>0$.
%\begin{equation}\label{compatible-condition-u}
% \begin{split}
% \int_{\O}u\cdot R\dd x=0 \;\;  \text{ for } R\in\mathcal{R}_{\O}, \forall t>0.
%\end{split}
%\end{equation}
Thus, $w=u_1-u_2$ satisfies
\begin{equation}\label{compatible-condition-u}
 \begin{split}
 \int_{\O}w\cdot R\dd x=0\;\;  \text{ for all } R \in\mathcal{R}_{\O} \text{ and all } t>0.
\end{split}
\end{equation}
For axisymmetric or spherical domains, Proposition 3.15 in \cite{Tapia-Amrouche-Conca-Ghosh-2021} gives the Poincar\'{e} type inequality:
\begin{equation}
\begin{split}
\|g\|_{L^2(\O)} \lesssim &  \|\nabla_x^{\text{s}}  g\|_{L^2(\O)}
+\Big(\sum_{j=1}^{\dim{\mathcal{R}_{\O}}}\Big|\int_{\O} g\cdot R_j\dd x\Big|^2\Big)^{\frac{1}{2}}
\end{split}
\end{equation}
for $g\in H^1(\O)$ with $g\cdot n|_{\pt\O}=0$.
Combined with  \eqref{compatible-condition-u}, this implies
\begin{equation}\label{uL2-nabla-s-uL2}
\begin{split}
\|w\|_{L^2(\O)} \lesssim  \|\nabla_x^{\text{s}}  w\|_{L^2(\O)} \;\; \text{  for all } t\geq 0.
\end{split}
\end{equation}
Combining \eqref{uL2-nabla-s-uL2} with the standard Korn-type inequality (Theorem 2.1 in \cite{Ciarlet-Ciarlet-2005})
\begin{equation}\label{Korn-type-inequality-H1}
\begin{split}
\|g\|_{H^1(\O)} \lesssim &  \|\nabla_x^{\text{s}}  g\|_{L^2(\O)}
 + \| g \|_{L^2(\O)}, \;\; \forall g\in H^1(\O),
\end{split}
\end{equation}
 we  obtain
\begin{equation}\label{korn-uH1-nabla-s-uL2}
\begin{split}
\|w\|_{H^1(\O)} \lesssim  \|\nabla_x^{\text{s}}  w\|_{L^2(\O)}.
\end{split}
\end{equation}
Inserting \eqref{u1-u2-L6-small} and \eqref{korn-uH1-nabla-s-uL2} into \eqref{A-w-L2-energy-inequality} gives
\begin{equation}\label{w-L2-energy-close}
\begin{split}
&\frac{1}{2} \|w(t)\|^2_{L^2(\O)} + \sigma \int_{0}^{t}  \|w\|_{H^1(\O)}^2\lesssim \frac{1}{2} \|w(t)\|^2_{L^2(\O)} + \sigma \int_{0}^{t}  \|\nabla_x^{\text{s}} w\|^2_{L^2(\O)}\leq 0,
\end{split}
\end{equation}
provided $u_0\in \mathbb{H}_{u}$. Hence $w\equiv 0$ and uniqueness follows.
\medskip

\noindent\textbf{Step 2.1.2. $\O$ non-axisymmetric.}

Here $ \mathcal{R}_\Omega=\{0\}$. Using \eqref{u1-u2-L6-small} and the first case of  \eqref{Korn-type-inequality} directly closes \eqref{A-w-L2-energy-inequality}.
\medskip

\noindent\textbf{Step 2.2. Estimate for $\vartheta$.}

 From \eqref{u1-u2-L6-small}, \eqref{A-chi-L2-energy-inequality} and \eqref{w-L2-energy-close}, we have
\begin{equation}\label{A-chi-L2-energy-close}
\begin{split}
&\|w(t)\|^2_{L^2(\O)} + \int_{0}^{t}  \|w\|_{H^1(\O)}^2+ \|\chi(t)\|^2_{L^2(\O)} + \int_{0}^{t}  \|\nabla\chi\|^2_{L^2(\O)}
\leq 0.
\end{split}
\end{equation}
Moreover, the heat equation with homogeneous Neumann condition $\partial_n \vartheta|_{\pt\O}=0$  satisfies conservation law:
  \begin{equation*}
\begin{split}
 \pt_t \int_{\O}\vartheta\dd x=0 \;\; \text{  for all } t> 0,
\end{split}
\end{equation*}
where we have used $\nabla_x \cdot u=0$ and $n\cdot u|_{\pt\O}=0$. Thus, if $\vartheta_0\in \mathbb{H}_{\th}$, then
$\int_{\O}\chi\dd x=0$ for all $t\geq 0$.
The Poincar\'{e} inequality therefore gives
\begin{equation*}
\begin{split}
\|\chi\|_{H^1(\O)} \lesssim \|\nabla_x\chi\|_{L^2(\O)}.
\end{split}
\end{equation*}
Combined with \eqref{A-chi-L2-energy-close}, this yields uniqueness of $\vartheta$.
 This completes the proof.
\end{proof}

\begin{lemma}[Global energy estimates and existence]
\label{lem:global-estimate-NS}
Under the assumptions of Lemma~\ref{A-lemma-NS}, the INSF system
admits a unique global weak solution $(u,\vartheta)$. Moreover, for
every $t\geq 0$,
\begin{equation}\label{eq:global-estimate-NS}
\begin{aligned}
 &\sup_{0\leq s\leq t}\|u(s)\|_{L^2(\O)}^2
 +\sigma\int_0^t\|u(s)\|_{H^1(\O)}^2\,ds
 \lesssim \|u_0\|_{L^2(\O)}^2,\\
 &\sup_{0\leq s\leq t}\|\vartheta(s)\|_{L^2(\O)}^2
 +\kappa\int_0^t\|\vartheta(s)\|_{H^1(\O)}^2\,ds
 \lesssim \|\vartheta_0\|_{L^2(\O)}^2 .
\end{aligned}
\end{equation}
In particular,
\begin{equation}\label{eq:global-existence-NS}
 (u,\vartheta)
 \in C\bigl([0,\infty);L^2(\O)\bigr)
 \cap L^\infty\bigl(\mathbb{R}^+;L^2(\O)\bigr)
 \cap L^2\bigl(\mathbb{R}^+;H^1(\O)\bigr).
\end{equation}
The implicit constants may depend on $\O$, $\sigma$, $\kappa$, and
$\lambda$, but are independent of $t$.
\end{lemma}

\begin{proof}
The estimates in \eqref{eq:global-estimate-NS} follow from the standard
$L^2$ energy inequalities, together with the Korn and
Poincar\'e--Friedrichs inequalities used in the proof of
Lemma~\ref{A-lemma-NS}. When $\lambda=0$, one also uses the conservation
of the constraints defining $\mathbb{H}_u$ and $\mathbb{H}_\vartheta$, as established
there.

A standard Faedo--Galerkin approximation, followed by the
Aubin--Lions compactness argument, yields a global Leray--Hopf weak
solution satisfying \eqref{eq:global-estimate-NS}. Letting
$t\to\infty$ in this estimate gives the $L^\infty_tL^2_x$ and
$L^2_tH^1_x$ assertions in \eqref{eq:global-existence-NS}.
Furthermore, the equations and the
$L^\infty(\mathbb{R}^+;L^6(\O))$ bound established in the proof of
Lemma~\ref{A-lemma-NS} give the corresponding local-in-time
time-derivative estimates in the dual energy spaces. The
Lions--Magenes lemma therefore yields the strong $L^2$-continuity in
\eqref{eq:global-existence-NS}. Uniqueness follows from
Lemma~\ref{A-lemma-NS}.
\end{proof}

\bigskip

%%%%%%%%%%%%%%%%%%%%%%%%%%%%%%%%%%
%%%%%%%%%%%%%%%%%%%%%%%%%%%%%%%%%%

\section{Gaussian Integration and Elliptic Estimates}

\begin{lemma}[Gaussian integrals on the half-line]
The following integrals hold:
\begin{align*}
&\int_{\R_+} \exp \Big ( -\frac{{v_1}^2}{2T} \Big) \dd v_1 = \sqrt{\frac{\pi}{2}} T^{1/2}, \quad\;\;\;\;\;
\int_{\R_+} v_1 \exp \Big( -\frac{v_1^2}{2T} \Big ) \dd v_1 = {T}, \\
&\int_{\R_+} v_1^2 \exp \Big (-\frac{v_1^2}{2T} \Big ) \dd v_1 = \sqrt{\frac{\pi}{2}} T^{3/2}, \quad\;
\int_{\R_+} v_1^3 \exp \Big (-\frac{v_1^2}{2T} \Big ) \dd v_1 = 2 {T}^2,\\
&\int_{\R_+} v_1^4 \exp \Big (-\frac{v_1^2}{2T} \Big ) \dd v_1 = 3\sqrt{\frac{\pi}{2}} T^{5/2}, \quad
\int_{\R_+} v_1^5 \exp \Big (-\frac{v_1^2}{2T} \Big ) \dd v_1 = 8 {T}^3.
\end{align*}
\end{lemma}
\begin{proof}[\textbf{Proof}.]\
These follow directly from standard Gaussian integral formulas.
\end{proof}
\medskip

\begin{lemma}\label{wholeinteg}
Let $\mu$ be the global Maxwellian defined in \eqref{def-mu} and $\tilde{\mu}$ the rotating Maxwellian defined in \eqref{tilde-mu-def}. Then the following integrals hold:
\begin{align*}
&\int_{\R^3} \tilde{\mu} \dd v = \rho, \quad
\int_{\R^3} (v-\c)\tilde{\mu} \dd v = 0, \quad \int_{\R^3} v\tilde{\mu} \dd v = \rho \c,\\
&\int_{\R^3} \norm{v-\c}^2\tilde{\mu} \dd v = 3\rho T, \quad \quad\quad\int_{\R^3} \norm{v}^2\tilde{\mu} \dd v = 3\rho T +\rho \norm{\c}^2, \\
&\int_{\R^3} (v-\c)\norm{v-\c}^2\tilde{\mu} \dd v = 0, \quad \int_{\R^3} v\norm{v}^2\tilde{\mu} \dd v = 5\rho T \c  + \rho \c \norm{\c}^2,\\
&\int_{\R^3} \norm{v-\c}^4\tilde{\mu} \dd v = 15\rho T^2, \quad\quad  \int_{\R^3} \norm{v}^4\tilde{\mu} \dd v = 15\rho T^2 + 10 \rho T \norm{\c}^2 + \rho \norm{\c}^4.
\end{align*}
\end{lemma}
\begin{proof}[\textbf{Proof}.]\  These follow from direct computation using the definition of $\tilde{\mu}$ and Gaussian integration.
\end{proof}
\medskip

%\begin{align*}
%\int_{\R^3} \tilde{\mu} \dd v =& \int_{\R^3} \frac{\rho}{(2 \pi T)^{3/2}} \exp %\Big(-\frac{\norm{v-\c}^2}{2T}\Big) \dd v = \rho,
%\end{align*}
%\begin{align*}
%\int_{\R^3} (v-\c) \tilde{\mu} \dd v =& \int_{\R^3} (v-\c)\frac{\rho}{(2 \pi T)^{3/2}} \exp %\Big(-\frac{\norm{v-\c}^2}{2T}\Big) \dd v =0,
%\end{align*}
%\begin{align*}
%\int_{\R^3} \norm{v-\c}^2 \tilde{\mu} \dd v =& \int_{\R^3} \norm{v-\c}^2\frac{\rho}{(2 \pi %T)^{3/2}} \exp \Big(-\frac{\norm{v-\c}^2}{2T}\Big) \dd v \\
%=& \frac{3 \rho}{(2 \pi T)^{3/2}}\Big[\int_{\R^3} v_1^2 \exp \Big(-\frac{v_1^2}{2T}\Big) \dd %v_1\Big]\Big[\int_{\R^3} \exp \Big(-\frac{v_1^2}{2T}\Big) \dd v_1\Big]^2\\
%=& 3\rho T.
%\end{align*}

\begin{lemma} \label{whole-integ-error-terms}
Let $\tilde{\mu}$ be the rotating Maxwellian defined in \eqref{tilde-mu-def}.
Assume that $|\c|+|\th|\ll 1$. Then the following almost orthogonality relations hold:
%\noindent $(1)$ Orthogonality relations for $\tilde{\chi}_4$:
\begin{align}
&\begin{cases}\label{tilde-chi-4-orthogonal}
%&\displaystyle
%\int_{\R^{3}}   \tilde{\chi}_{4} \tilde{\chi}_{0} \dd v
%=   \int_{\R^{3}}   \frac{\norm{v}^2-3}{\sqrt{6}}\tilde{\mu}
%=      O(|\c|+|\th|),  \\
&\displaystyle
\int_{\R^{3}}   \tilde{\chi}_{4}  \tilde{\chi}_{k}\dd v
%=   \int_{\R^{3}}   \frac{\norm{v}^2-3}{\sqrt{6}} v_k \tilde{\mu}
=      O(|\c|+|\th|), \;\; k=0, 1,2,3,   \\
&\displaystyle
\int_{\R^{3}}   \tilde{\chi}_{4} \tilde{\chi}_{4}\dd v
%=   \int_{\R^{3}}   \frac{(\norm{v}^2-3)^2}{\sqrt{6}} \tilde{\mu}
=     1+ O(|\c|+|\th|);
\end{cases}\\
%\end{align}
%\noindent $(2)$ Orthogonality relations for $v_iv_j\sqrt{\tilde{\mu}}$:
%\begin{equation}
& \begin{cases}\label{vivj-chi-k-orthogonal}
&\displaystyle
\int_{\R^{3}}   v_iv_j\sqrt{\tilde{\mu}} \tilde{\chi}_{0} \dd v
%=   \int_{\R^{3}}   \frac{\norm{v}^2-3}{\sqrt{6}}\tilde{\mu}
=    1+  O(|\c|+|\th|),   \quad   i,j=1,2,3,  \\
&\displaystyle
\int_{\R^{3}}   v_iv_j\sqrt{\tilde{\mu}} \tilde{\chi}_{k}\dd v
%=   \int_{\R^{3}}   \frac{\norm{v}^2-3}{\sqrt{6}} v_k \tilde{\mu}
=      O(|\c|+|\th|),    \quad   i,j,k=1,2,3,   \\
&\displaystyle
\int_{\R^{3}}   v_iv_j\sqrt{\tilde{\mu}} \tilde{\chi}_{4}\dd v
%=   \int_{\R^{3}}   \frac{(\norm{v}^2-3)^2}{\sqrt{6}} \tilde{\mu}
=      \frac{2}{\sqrt{6}} \d_{ij}+ O(|\c|+|\th|),   \quad   i,j=1,2,3;
\end{cases}\\
%\noindent $(3)$ Mixed moments with $v_i$ and $\tilde{\chi}_4$:
%\begin{equation}
&\begin{cases}\label{vi-tilde-chi4-orthogonal}
&\displaystyle \int_{\R^{3}}
 v_i \tilde{\chi}_4 \tilde{\chi}_j\dd v = O(|\c|+|\th|),\quad i=1,2,3,j=0,4, \\
%&\displaystyle
% \int_{\R^{3}}
%v_i \tilde{\chi}_4 \tilde{\chi}_4\dd v = O(|\c|+|\th|),\quad
%  i=1,2,3, \\
%---------------
&\displaystyle\int_{\R^{3}}
 v_i \tilde{\chi}_4 \tilde{\chi}_j\dd v = \frac{2}{\sqrt{6}} \d_{ij}+ O(|\c|+|\th|), \quad
 i,j=1,2,3;
\end{cases}\\
%\noindent $(4)$ Third and fourth moments:
&\begin{cases}\label{b2-estimate-Gauss-theta3}
&\displaystyle\int_{\R^{3}}
 v_i v_jv_k \tilde{\mu}\dd v = O(|\c|+|\th|), \quad    i,j,k=1,2,3,
 \\
 &\displaystyle\int_{\R^{3}}
 v_i v_jv_k (|v|^2-3) \tilde{\mu}\dd v = O(|\c|+|\th|),  \quad    i,j,k=1,2,3, \\
%---------------
&\displaystyle\int_{\mathbb{R}^3}v_i^2 v_j^2\tilde{\mu} \dd v=
 \left\{
   \begin{array}{ll}
     3+O(|\c|+|\th|), & \hbox{if $i=j$,} \\
     1+O(|\c|+|\th|), & \hbox{if $i\neq j$;}
   \end{array}
 \right.
\end{cases}\\
%\end{align}
%\noindent $(5)$ Orthogonality with $v_i(|v|^2 - 10)$:
%\begin{align}
&\begin{cases}\label{vi(|v|^2-10)-orthogonality}
&\displaystyle \int_{\R^{3}}   v_i(|v|^2 - 10)\sqrt{\tilde{\mu}}  \tilde{\chi}_{j}\dd v
%=   \int_{\R^{3}}   v_i  (|v|^2 - 10)\tilde{\mu}\dd v
=       O(|\c|+|\theta|),\quad   i=1,2,3, j=0,4,   \\
%&\displaystyle\int_{\R^{3}}   v_i(|v|^2 - 10)\sqrt{\tilde{\mu}}  \tilde{\chi}_{4}\dd v
%=   \frac{1}{\sqrt{6}}\int_{\R^{3}}   v_i  (|v|^2 - 10)(|v|^2 - 3)\tilde{\mu}\dd v
%=       O(|\c|+|\theta|),\quad   i=1,2,3,\\
&\displaystyle\int_{\R^{3}}   v_i(|v|^2 - 10)\sqrt{\tilde{\mu}}  \tilde{\chi}_{j}\dd v
%=   \int_{\R^{3}}   v_i v_j (|v|^2 - 10)\tilde{\mu}\dd v
=     -5\delta_{ij} +   O(|\c|+|\theta|),\quad    i,j=1,2,3;
\end{cases}\\
%\end{equation}
%\noindent $(6)$ Orthogonality with $v_i v_j(|v|^2 - 10)$:
%\begin{equation}
&\begin{cases}\label{v_iv_j(|v|^2 - 10)-orthogonal 1}
&\displaystyle\int_{\R^{3}}   v_iv_j(|v|^2 - 10)\sqrt{\tilde{\mu}}  \tilde{\chi}_{k}
%---------------
%=   \int_{\R^{3}}   v_iv_jv_k  (|v|^2 - 10)\tilde{\mu}
%---------------
=      O(|\c|+|\theta|),   \quad    i,j,k=1,2,3,   \\
%&\displaystyle
%\int_{\R^{3}}   v_iv_j(|v|^2 - 10)\sqrt{\tilde{\mu}}  \tilde{\chi}_{k}
%---------------
%=   \frac{1}{\sqrt{6}}\int_{\R^{3}}   v_i v_j (|v|^2 - 10)(|v|^2 - 3)\tilde{\mu}
%---------------
%=      O(|\c|+|\theta|),   \quad   i,j=1,2,3,k=0,4\\
&\displaystyle
\int_{\R^{3}}   v_iv_j(|v|^2 - 10)\sqrt{\tilde{\mu}}  \tilde{\chi}_{k}
%---------------
%=   \int_{\R^{3}}   v_i v_j (|v|^2 - 10)\tilde{\mu}
%---------------
=     -5\delta_{ij} +  O(|\c|+|\theta|),   \quad   i,j=1,2,3, k=0,4;
\end{cases}\\
%\end{equation}
%\noindent $(7)$ Orthogonality with $v_i v_j(|v|^2 - 5)$:
%\begin{equation}
&\begin{cases}\label{v_iv_j(|v|^2 - 5)-orthogonal}
%&\displaystyle
%\int_{\R^{3}}   v_iv_j(|v|^2 - 5)\sqrt{\tilde{\mu}}  \tilde{\chi}_{0}\dd v
%=   \int_{\R^{3}}   v_i v_j (|v|^2 - 5)\tilde{\mu}
%=      O(|\c|^2+|\theta|^2),   \quad    i,j=1,2,3, \\
&\displaystyle
\int_{\R^{3}}   v_iv_j(|v|^2 - 5)\sqrt{\tilde{\mu}}  \tilde{\chi}_{k}\dd v
%=   \int_{\R^{3}}   v_iv_jv_k  (|v|^2 - 5)\tilde{\mu}
=      O(|\c|+|\theta|),   \quad    i,j=1,2,3, k=0,1,2,3, \\
&\displaystyle
\int_{\R^{3}}   v_iv_j(|v|^2 - 5)\sqrt{\tilde{\mu}}  \tilde{\chi}_{4}\dd v
%=   \frac{1}{\sqrt{6}}\int_{\R^{3}}   v_i v_j (|v|^2 - 5)(|v|^2 - 3)\tilde{\mu}
=    \frac{10}{\sqrt{6}}\d_{ij}+  O(|\c|+|\theta|) \quad    i,j=1,2,3.
\end{cases}
\end{align}
\end{lemma}
\begin{proof}[\textbf{Proof}.]\
All relations follow from Lemma \ref{wholeinteg} together with the definition of $\tilde{\mu}$.
\end{proof}
\medskip

\begin{lemma}[Boundary integrals] \label{boundaryinteg}
Let $\tilde{\mu}$ be the rotating Maxwellian defined in \eqref{tilde-mu-def}.
Then
\begin{align}
&\int_{n \cdot v >0} \tilde{\mu} [n \cdot v] \dd v = \frac{\rho T^{1/2}}{(2 \pi)^{1/2}}, \quad
\int_{n \cdot v >0} (v-\c) \tilde{\mu} [n \cdot v]  \dd v = {\frac{\rho nT}{2}}, \label{v-c-boundary-integral}\\
&\int_{n \cdot v >0} \norm{v-\c}^2 \tilde{\mu} [n \cdot v]  \dd v = \frac{4\rho T^{3/2}}{(2 \pi)^{1/2}},\quad
\int_{n \cdot v >0} (v-\c)\norm{v-\c}^2 \tilde{\mu} [n \cdot v]  \dd v = {\frac{5\rho nT^2}{2}},\label{v-c-v-c-2-boundary-integral}\\
&\int_{n \cdot v >0} \norm{v-\c}^4 \tilde{\mu} [n \cdot v]  \dd v = \frac{24 \rho T^{5/2}}{(2 \pi)^{1/2}},\label{v-c-4-boundary-integral}\\
%---------------------------
&\int_{n \cdot v >0} v \tilde{\mu} [n \cdot v]  \dd v = \frac{\rho  \c T^{1/2}}{(2 \pi)^{1/2}}  {+ \frac{\rho nT}{2}},\label{v1-boundary-integral}\\
&\int_{n \cdot v >0} \norm{v}^2 \tilde{\mu} [n \cdot v]  \dd v = \frac{4\rho T^{3/2}}{(2 \pi)^{1/2}} + \frac{\rho \norm{ \c }^2T^{1/2}}{(2 \pi)^{1/2}},\label{v2-boundary-integral}\\
&\int_{n \cdot v >0} v\norm{v}^2 \tilde{\mu} [n \cdot v]  \dd v = \frac{6\rho \c T^{3/2}}{(2 \pi)^{1/2}} + \frac{\rho  \c \norm{ \c }^2T^{1/2}}{(2 \pi)^{1/2}} {+ \frac{5\rho nT^2}{2} + \frac{\rho nT\norm{\c}^2}{2}},\label{v3-boundary-integral}\\
&\int_{n \cdot v >0} \norm{v}^4 \tilde{\mu} [n \cdot v]  \dd v = \frac{24\rho T^{5/2}}{(2 \pi)^{1/2}}  + \frac{12\rho \norm{ \c }^2T^{3/2}}{(2 \pi)^{1/2}}
+\frac{\rho \norm{ \c }^4T^{1/2}}{(2 \pi)^{1/2}}. \label{v4-boundary-integral}
\end{align}
\end{lemma}

\begin{proof}[\textbf{Proof}.]\
Decompose $v= v_{\perp} + {v_{\|}}n$ with $v_\| \in \R$, $v_{{\perp}} \in \R^2$, where   $v_{\|} n \,\| \, n$ and $v_{\perp} \perp n$. By the definition of $\c$ in \eqref{u-definition}, we have $\c \perp n$.

Direct computation using Lemma \ref{wholeinteg} gives
\begin{align*}
&\int_{n \cdot v >0} \tilde{\mu} [n \cdot v] \dd v = \frac{\rho }{(2 \pi T)^{3/2}}
 \int_{\R_{+} \times \R \times \R} v_{\|} \exp \Big(-\frac{v_{\|}^2 + \norm{v_{\perp}- \c }^2 }{2T}\Big) \dd v_{\|} \dd v_{\perp}
% =& \frac{\rho }{(2 \pi T)^{1/2}} \int_{\R_{+}} v_{\|} \exp \Big(-\frac{v_{\|}^2}{2T}\Big) \dd v_{\|}
=\frac{\rho T^{1/2}}{(2 \pi)^{1/2}}, \\
%----------------------
&\int_{n \cdot v >0}(v- \c ) \tilde{\mu} [n \cdot v] \dd v = \frac{\rho }{(2 \pi T)^{3/2}}
\int_{\R_{+} \times \R \times \R} \big(v_{\perp}- \c  + {v_{\|}}n\big) {v_{\|}} \exp \Big(-\frac{v_{\|}^2 + \norm{v_{\perp}- \c }^2}{2T}\Big) \dd v_{\|} \dd v_{\perp}
 = \frac{\rho nT}{2}.
\end{align*}
This establishes \eqref{v-c-boundary-integral}. Proceeding similarly, we obtain
\begin{align*}
&\int_{n \cdot v >0}\norm{v- \c }^2 \tilde{\mu} [n \cdot v] \dd v \\
=& \frac{\rho }{(2 \pi T)^{3/2}}
\int_{\R_{+} \times \R \times \R} v_{\|}(v_{\|}^2 + \norm{v_{\perp}- \c }^2)\exp \Big(-\frac{v_{\|}^2 + (v_{\perp}- \c )^2}{2T}\Big) \dd v_{\|} \dd v_{\perp}
%=& \frac{\rho T^{3/2}}{(2 \pi)^{1/2}}(2+1+1)
=\frac{4\rho T^{3/2}}{(2 \pi)^{1/2}},\\
%---------------
& \int_{n \cdot v >0}(v- \c )\norm{v- \c }^2 \tilde{\mu} [n \cdot v] \dd v \\
=& \frac{1}{(2 \pi T)^{3/2}}
\int_{\R_{+} \times \R \times \R} (v_{\perp}- \c  + v_{\|}n) {v_{\|}} (v_{\|}^2 + \norm{v_{\perp}- \c }^2)\exp \Big(-\frac{v_{\|}^2 + \norm{v_{\perp}- \c }^2}{2T}\Big) \dd v_{\|} \dd v_{\perp}\\
=& \frac{1}{(2 \pi T)^{3/2}}
\int_{\R_{+} \times \R \times \R} n {v_{\|}}^2 (v_{\|}^2 + \norm{v_{\perp}- \c }^2)\exp \Big(-\frac{v_{\|}^2 + \norm{v_{\perp}- \c }^2}{2T}\Big) \dd v_{\|} \dd v_{\perp}
 = \frac{5\rho nT^2}{2},
\\
%---------------
&\int_{n \cdot v >0}\norm{v- \c }^4 \tilde{\mu} [n \cdot v] \dd v \\
=& \frac{\rho }{(2 \pi T)^{3/2}}
\int_{\R_{+} \times \R \times \R} v_{\|}(v_{\|}^2 + \norm{v_{\perp}- \c }^2)^2 \exp \Big(-\frac{v_{\|}^2 +\norm{v_{\perp}- \c }^2}{2T}\Big) \dd v_{\|} \dd v_{1} \dd v_{2} \\
=& \frac{\rho }{(2 \pi T)^{3/2}}
\int_{\R_{+} \times \R \times \R} v_{\|}\Big(v_{\|}^4 + v_{1}^4 +v_{2}^4 + 2v_{\|}^2v_{1}^2 + 2 v_{1}^2v_{2}^2 + 2v_{2}^2v_{\|}^2 \Big)
\exp \Big(-\frac{v_{\|}^2 + v_{1}^2 +v_{2}^2}{2T}\Big) \dd v_{\|} \dd v_{1} \dd v_{2} \\
%=& \frac{\rho T^{5/2}}{(2 \pi)^{1/2}}(8 + 3+3 +2\times 2 + 2\times 1+ 2\times 2)
 =&\frac{24 \rho T^{5/2}}{(2 \pi)^{1/2}}.
\end{align*}
This proves \eqref{v-c-v-c-2-boundary-integral}--\eqref{v-c-4-boundary-integral}.

%\begin{align*}
%\int_{n \cdot v >0} v \tilde{\mu} [n \cdot v] \dd v =& \int_{n \cdot v >0}  \c  %\tilde{\mu} [n \cdot v] \dd v + \frac{\rho nT}{2}= \frac{\rho  \c T^{1/2}}{(2 \pi)^{1/2}} %+ \frac{\rho nT}{2}.
%\end{align*}

The relation \eqref{v1-boundary-integral} follows from \eqref{v-c-boundary-integral}.
For \eqref{v2-boundary-integral}, using $|v|^2 = |v - \c|^2 + 2(v - \c) \cdot \c + |\c|^2$ and noting $\c \cdot n = 0$, we have
\begin{align*}
\int_{n \cdot v >0}\norm{v}^2 \tilde{\mu} [n \cdot v] \dd v
= &\int_{n \cdot v >0}\Big(\norm{v- \c }^2 + 2 (v- \c ) \cdot  \c  + \norm{ \c }^2 \Big) \tilde{\mu} [n \cdot v] \dd v\\
=& \frac{4\rho T^{3/2}}{(2 \pi)^{1/2}} + \frac{\rho \norm{ \c }^2T^{1/2}}{(2 \pi)^{1/2}}.
\end{align*}

For \eqref{v3-boundary-integral}, we decompose
\begin{align*}
v\norm{v}^2 = (v- \c )\norm{v- \c }^2 +  \c  \norm{v- \c }^2 + 2 (v- \c )  \c  \cdot (v- \c ) + 2  \c   \c  \cdot (v- \c ) + (v- \c )\norm{ \c }^2 +  \c  \norm{ \c }^2.
\end{align*}
Splitting $v_{\perp}$ into components parallel $v_{\c}$ and perpendicular $v_{\c^{\perp}}$ to $\c$, direct computation yields:
\begin{align*}
& \int_{n \cdot v >0}(v- \c )\norm{v- \c }^2 \tilde{\mu} [n \cdot v] \dd v = \frac{5\rho nT^2}{2},\qquad  \int_{n \cdot v >0} \c  \norm{v- \c }^2 \tilde{\mu} [n \cdot v] \dd v = \frac{4\rho  \c T^{3/2}}{(2 \pi)^{1/2}},\\
%--------------
& \int_{n \cdot v >0} \c   \c  \cdot (v- \c ) \tilde{\mu} [n \cdot v] \dd v
=  \c   \c \cdot\frac{\rho nT}{2}
= 0,\quad
%------------
 \int_{n \cdot v >0}(v-u)\norm{ \c }^2 \tilde{\mu} [n \cdot v] \dd v
= \frac{\rho nT\norm{ \c }^2}{2},\\
%--------------
& \int_{n \cdot v >0} \c  \norm{ \c }^2 \tilde{\mu} [n \cdot v] \dd v
= \frac{\rho  \c \norm{ \c }^2T^{1/2}}{(2 \pi)^{1/2}},\qquad
 \int_{n \cdot v >0}(v- \c )  \c  \cdot (v- \c ) \tilde{\mu} [n \cdot v] \dd v
%= \int_{\R_{+}\times \R \times \R} \c  v_{\|} \norm{v_{ \c } - \c }^2 \tilde{\mu} \dd v_{\|} \dd v_{ \c } \dd v_{ \c ^{\perp}}
= \frac{\rho  \c T^{3/2}}{(2 \pi)^{1/2}}.
\end{align*}
Combining these results proves \eqref{v3-boundary-integral}.
%\begin{align*}
%&\int_{n \cdot v >0} v\norm{v}^2 \tilde{\mu} [n \cdot v] \dd v
%=\frac{6\rho  \c T^{3/2}}{(2 \pi)^{1/2}}  + \frac{\rho  \c \norm{ \c }^2T^{1/2}}{(2 %\pi)^{1/2}} + \frac{5\rho nT^2}{2} + \frac{\rho nT\norm{ \c }^2}{2},
%\end{align*}
%which proves \eqref{v3-boundary-integral}.

For \eqref{v4-boundary-integral}, we use the decomposition
\begin{align*}
\norm{v}^4 = \norm{v- \c }^4 + \norm{ \c }^4 + 4((v- \c ) \cdot  \c )^2 + 2 \norm{v- \c }^2 \norm{ \c }^2 + \big(\text{odd order of } (v- \c ) \cdot  \c\big).
\end{align*}
Then, the above calculations indicate
\begin{align*}
\int_{n \cdot v >0} \norm{v}^4 \tilde{\mu} [n \cdot v] \dd v
=&\int_{n \cdot v >0} (\norm{v- \c }^4 + \norm{ \c }^4 + 4((v- \c ) \cdot  \c )^2 + 2 \norm{v- \c }^2 \norm{ \c }^2)\tilde{\mu} [n \cdot v] \dd v \\
=&\frac{24\rho T^{5/2}}{(2 \pi)^{1/2}}  + \frac{\rho \norm{ \c }^4T^{1/2}}{(2 \pi)^{1/2}} + \frac{4\rho \norm{ \c }^2T^{3/2}}{(2 \pi)^{1/2}} + \frac{8\rho \norm{ \c }^2T^{3/2}}{(2 \pi)^{1/2}},
%=&\frac{24\rho T^{5/2}}{(2 \pi)^{1/2}}  + \frac{12\rho \norm{ \c }^2T^{3/2}}{(2 %\pi)^{1/2}}
%+\frac{\rho \norm{ \c }^4T^{1/2}}{(2 \pi)^{1/2}},
\end{align*}
which further leads to \eqref{v4-boundary-integral}.
This complete the proof.
\end{proof}
\medskip

The next result is standard in elliptic theory  (see, e.g., \cite{gilbarg1977elliptic}).

\begin{lemma}\label{Poisson-equation-theory}
Let $p\in \{2, \frac{6}{5}\}$, and let $\xi \in L^p(\Omega)$ and satisfy the compatible condition $\int_{\Omega} \xi \dd x =0$. Then the elliptic equation
\begin{align}\label{Poisson-equation}
- \Delta_x \phi = \xi \;\text{ in }\O, \quad& {\pt_n} \phi = 0  \;\text{ on }\pt\O,  \quad \int_{\Omega} \phi  \dd x =0.
\end{align}
  admits a unique solution $\phi\in W^{2,p}(\Omega)$ satisfying
\begin{align}
\normm{\nabla_x^2 \phi}_{L^2_x} + \normm{\nabla_x \phi}_{L^2_x} + \normm{\phi}_{L^2_x} &\lesssim \normm{\xi}_{L^2_x}, \quad \text{ if } \xi \in L^{2}(\Omega), \label{elliptic-equation-estimate-L2} \\
\normm{\nabla_x^2 \phi}_{L^{\frac{6}{5}}_x} + \normm{\nabla_x \phi}_{L^2_x} + \normm{\phi}_{L^{6}_x} & \lesssim \normm{\xi}_{L^{\frac{6}{5}}_x}, \quad \text{ if } \xi \in L^{\frac{6}{5}}(\Omega).\label{elliptic-equation-estimate-L6}
\end{align}
\end{lemma}
\medskip

The following lemma is adapted from Theorem 2.11 in \cite{Bernou2022}
and Lemma 3 in \cite{Chen2024}.

\begin{lemma} \label{elliptic-system-theory}
Let $\xi : \Omega \to \R^3$, and let $\phi$ satisfy the elliptic system
\begin{equation}\label{elliptic-system}
\begin{split}
-\text{div}  (\nabla^{\text{s}}_x \phi) = \xi \quad & \text{in } \Omega,\\
\phi \cdot n =0 \quad& \text{on } \pt \Omega, \\
(\nabla^{\text{s}}_x \phi) n = (\nabla^{\text{s}}_x\phi: n \otimes n)n \quad& \text{on } \pt \Omega.
\end{split}
\end{equation}

\noindent (1)  If $\xi \in L^2(\Omega)$, then the variational formulation
\begin{equation}\label{elliptic-system-variational-formulation}
\begin{split}
\int_{\Omega} \nabla^{\text{s}}_x \phi: \nabla^{\text{s}}_x \sigma \dd x
= \int_{\Omega} \xi\cdot \sigma \dd x   \;\; \text{ for all }  \sigma \in \mathscr{H}(\O)
\end{split}
\end{equation}
admits a unique weak solution $\phi\in \mathscr{H}(\O)$. Here
 \begin{equation}\label{elliptic-system-variational-formulation-space}
\begin{split}
\mathscr{H}(\O):=\Big\{ \sigma:\O \to \R^3 :\; \sigma\in H^1_x(\O),\;  \sigma\cdot n\big|_{\pt \O}=0, \; P_{\Omega}\Big(\int_{\Omega} \nabla^{\text{a}}_x\sigma \dd x\Big)=0 \Big\}
\end{split}
\end{equation}
and $P_{\O}$ denotes the orthogonal projection onto the set
$A_\O:=\big \{ A\in \mathfrak{so}(3,\mathbb{R}): \; Ax\in \mathcal{R}_{\O} \big \}$.

\noindent (2) Let $p\in \{2, \frac{6}{5}\}$ and assume $\xi \in L^p(\Omega)$  satisfies the compatible condition
\begin{equation}\label{elliptic-system-compatible-condition}
\begin{split}
\int_{\Omega} Ax \cdot \xi(x) \dd x =0 \quad \text{ for any } Ax \in \mathcal{R}_{\Omega}.
\end{split}
\end{equation}
 Then \eqref{elliptic-system} admits a unique strong solution $\phi\in W^{2,\frac{6}{5}}_x(\O)\cap \mathscr{H}(\O)$ with
\begin{align}
\normm{\nabla^2_x \phi}_{L^2_x} + \normm{\nabla_x \phi}_{L^2_x} + \normm{\phi}_{L^2_x} \lesssim &\normm{\xi}_{L^2_x}, \quad \text{ if } \xi \in L^2(\Omega), \label{elliptic-system-estimate-L2} \\
\normm{\nabla^2_x \phi}_{L^{\frac{6}{5}}_x} + \normm{\nabla_x \phi}_{L^{2}_x} + \normm{\phi}_{L^{6}_x} \lesssim &\normm{\xi}_{L^{\frac{6}{5}}_x}, \quad \text{ if } \xi \in L^{\frac{6}{5}}(\Omega). \label{elliptic-system-estimate-L6}
\end{align}
\end{lemma}
\bigskip

%%%%%%%%%%%%%%%%%%%%%%%%%%%%%%%%%%
%%%%%%%%%%%%%%%%%%%%%%%%%%%%%%%%%%

\noindent\textbf{Acknowledgements.}
F. Zhou is supported by NSFC grant 12271179.

\medskip

\noindent\textbf{Conflict of Interest Statement}

The authors declare that they have no conflicts of interest.
\medskip

\noindent\textbf{Data Availability Statement}

No data was used for the research described in the article.
\medskip

\bibliographystyle{plain}

%%%%%%%%%%%%%%%%%%%%%%%%%%%%%%%%%%
%%%%%%%%%%%%%%%%%%%%%%%%%%%%%%%%%%
%%%%%%%%%%%%%%%%%%%%%%%%%%%%%%%%%%
%%%%%%%%%%%%%%%%%%%%%%%%%%%%%%%%%%

\end{document}